\documentclass[12pt]{amsbook}
\usepackage{amssymb}
\usepackage[all]{xy}

\usepackage[colorlinks=true,linkcolor=magenta,citecolor=magenta]{hyperref}
\makeindex

\newtheorem{theorem}{Theorem}[chapter]
\newtheorem{answer}[theorem]{Answer}
\newtheorem{cat}[theorem]{Cat}
\newtheorem{claim}[theorem]{Claim}

\newtheorem{convention}[theorem]{Convention}
\newtheorem{criticism}[theorem]{Criticism}
\newtheorem{definition}[theorem]{Definition}
\newtheorem{exercise}[theorem]{Exercise}
\newtheorem{fact}[theorem]{Fact}
\newtheorem{goal}[theorem]{Goal}
\newtheorem{plan}[theorem]{Plan}

\newtheorem{problem}[theorem]{Problem}
\newtheorem{proposition}[theorem]{Proposition}
\newtheorem{question}[theorem]{Question}
\newtheorem{questions}[theorem]{Questions}

\newtheorem{thought}[theorem]{Thought}

\begin{document}

\title{Easy quantum groups}

\author{Teo Banica}
\address{Department of Mathematics, University of Cergy-Pontoise, F-95000 Cergy-Pontoise, France. {\tt teo.banica@gmail.com}}


\begin{abstract}
A closed subgroup $G\subset_uU_N^+$ is called easy when its associated Tannakian category $C_{kl}=Hom(u^{\otimes k},u^{\otimes l})$ appears from a category of partitions, $C=span(D)$ with $D=(D_{kl})\subset P$, via the standard implementation of partitions as linear maps. The examples abound, and the main known subgroups $G\subset U_N^+$ are either easy, or not far from being easy. We discuss here the basic theory, examples and known classification results for the easy quantum groups $G\subset U_N^+$, as well as various generalizations of the formalism, known as super-easiness theories, and the unification problem for them.
\end{abstract}

\maketitle

\chapter*{Preface}

One problem in quantum mechanics is that the particles there are of ``slippery'' nature, having no clear positions and speeds. The trick, going back to Heisenberg, Schr\"odinger and others, is to relax, and look at this with a mix of algebraic and probabilistic thoughts. Following Heisenberg, the point is that of accepting that, at that tiny scales, our usual $\mathbb R^3$ might get replaced by something more complicated, of ``matrix'' type, where functions do not necessarily commute, $fg\neq gf$. Also, following Schr\"odinger, what about lowering our ambitions too, and instead of looking for exact positions and speeds, rather look for the probability of finding the particle here or there, and with this or that speed. 

\bigskip

The two points of view, of Heisenberg and Schr\"odinger, are in fact equivalent, and the credit for the unification, establishing quantum mechanics as a sort of peculiar mix of matrix type algebra and probability, goes to Dirac. Later on Feynman and others came with something even better, quantum electrodynamics (QED), but the badge ``mixture of bizarre algebra plus probability'' still remained attached to the theory.

\bigskip

To us, mathematical physicists, this suggests to get involved into building abstract theories based on noncommutativity, $fg\neq gf$, and with some probabilistic ideas in mind. With a bit of luck, maybe one day such a theory will prove to be useful in physics, helping in substantially improving what Feynman and others are saying.

\bigskip

The quantum groups come from this philosophy. To be more precise, a quantum group $G$ is something similar to a group, except for the fact that the functions on it do not necessarily commute, $fg\neq gf$. Here ``function'' can mean many things, such as usual complex function $f:G\to\mathbb C$, or function taking values in an arbitrary field $f:G\to F$, or even function defined infinitesimally around the origin, $f:\mathfrak g\to F$, or why not element of the enveloping Lie algebra $f\in U\mathfrak g$, and so on. In short, many choices here, but you get the point, up to some algebra, nocommutativity somewhere, $fg\neq gf$.

\bigskip

What are then the simplest, the ``easiest'' quantum groups?

\bigskip

This is a tricky question, whose answer depends on your knowledge of mathematics, and perhaps on your physics background too. At the white belt level, meaning decent knowledge of group theory, and skipped physics classes when in school, many interesting theories can be developed. At a more advanced level, however, you will learn that no matter what you want to do with your quantum group $G$, be that algebra, geometry, analysis, probability or physics, your favorite pastime will be that of using the Tannakian category of $G$, which is the collection $C=\{C_{kl}\}$ of linear spaces given by:
$$C_{kl}=Hom(u^{\otimes k},u^{\otimes l})$$

Here $u$ is a fundamental representation of $G$, and the simplest way of having such a fundamental representation is by assuming that $G$ appears as a closed subgroup of the free unitary group, $G\subset_uU_N^+$. That is, $G$ must be a ``compact quantum Lie group''. But with this picture in mind, easiness of $G$ should mean that $C$ is as easy as possible, and this technically means that $C$ must appear from a category of partitions $D=(D_{kl})\subset P$, via the standard implementation of partitions as linear maps, as follows:
$$C=span(D)$$

Thanks to theorems of Brauer and others, the examples of easy quantum groups abound, and in fact all the main known subgroups $G\subset U_N^+$ are either easy, or not far from being easy. Which is something remarkable, and technically, very useful.

\bigskip

The present book is about this, introduction and basic theory, examples and classification results for the easy quantum groups $G\subset U_N^+$, and various generalizations of the formalism, known as super-easiness theories. The book is organized in 4 parts, with Part I discussing the basics of easiness, Parts II and III getting more in detail into the continuous case, and the discrete case, and Part IV dealing with super-easiness.

\bigskip

This is my third quantum group book, coming as a continuation of \cite{ba1}, general introduction to quantum groups, and \cite{ba2}, technical book on quantum permutations. There will be a fourth and final book as well \cite{ba3}, discussing arithmetic aspects.

\bigskip

It is a pleasure to thank everyone, young or old, involved in the theory discussed below. Starting of course with Hermann Weyl. Many thanks go as well to my cats. They say that electrons at rest purr too, but this remains to be mathematically confirmed.

\bigskip

\

{\em Cergy, July 2025}

\smallskip

{\em Teo Banica}

\baselineskip=15.95pt
\tableofcontents
\baselineskip=14pt

\part{Easy quantum groups}

\ \vskip50mm

\begin{center}
{\em Love is a stranger

In an open car

To tempt you in

And drive you far away}
\end{center}

\chapter{Brauer theorems}

\section*{1a. Abelian groups}

Welcome to easiness. As the name tends to indicate, this is something rather algebraic, of yoga nature. And think here for instance, for a vague comparison, to the modern concepts of algebraic geometry, aren't these something easy too, perhaps even trivial.

\bigskip

However, there are no reasons to be scared about this. The idea indeed is that, once you are at a reasonably advanced level in representation theory, easiness is something truly simple, deserving its name, and that you should perfectly master, and build your theories upon. Our goal in this first chapter will be to reach to the ``reasonably advanced level'' that is needed, in order to talk about easiness. On the menu, compact Lie groups $G\subset U_N$, with lots of computations for them, then Brauer theorems and easiness.

\bigskip

For simplicity, let us  first restrict the attention to the real case, $G\subset O_N$. Normally the first thing to be done with such a beast is to consider its Lie algebra $\mathfrak g=T_1G$. Indeed, this Lie algebra $\mathfrak g$ is a vector space, so what we have here is a potentially fruitful ``linearization'' idea. However, and here comes our point, the construction $G\to\mathfrak g$ is not the only possible linearization of $G$. As a rival construction, we have:

\index{Lie group}
\index{Lie algebra}
\index{Tannakian category}
\index{tensor category}

\begin{definition}
Associated to any closed subgroup $G\subset O_N$ are the vector spaces
$$C_{kl}=\left\{T\in\mathcal L(H^{\otimes k},H^{\otimes l})\Big|Tg^{\otimes k}=g^{\otimes l}T,\forall g\in G\right\}$$
where $H=\mathbb C^N$. We call Tannakian category of $G$ the collection of spaces $C=(C_{kl})$.
\end{definition}

Observe that, due to $g\in G\subset O_N\subset\mathcal L(H)$, we have $g^{\otimes k}\in\mathcal L(H^{\otimes k})$ for any $k$, so the equality $Tg^{\otimes k}=g^{\otimes l}T$ makes indeed sense, as an equality of maps as follows: 
$$Tg^{\otimes k},g^{\otimes l}T\in \mathcal L(H^{\otimes k},H^{\otimes l})$$

It is also clear by definition that each $C_{kl}$ is a complex vector space. Moreover, it is also clear by definition that $C=(C_{kl})$ is indeed a category, in the sense that:
$$T\in C_{kl}\ ,\ S\in C_{lm}\implies ST\in C_{km}$$

Quite remarkably, the closed subgroup $G\subset O_N$ can be reconstructed from its Tannakian category $C=(C_{kl})$, and in a very simple way. More precisely, we have:

\begin{claim}
Given a closed subgroup $G\subset O_N$, we have
$$G=\left\{g\in O_N\Big|Tg^{\otimes k}=g^{\otimes l}T,\forall k,l,\forall T\in C_{kl}\right\}$$
where $C=(C_{kl})$ is the associated Tannakian category.
\end{claim}

Summarizing, we have some sort of idea here, and the question that you surely have in mind is probably, okay, but why not formulating our claim as a theorem, then proving our theorem, like mathematicians do, and then going ahead, saying what we have to say. Good point, but the problem is that, what we have to say is precisely this claim.

\bigskip

In short, and please don't misunderstand me, but Claim 1.2 is something extremely deep, and in any case far deeper than what you can possibly imagine. So, we will first attempt to understand what Claim 1.2 really says, with some illustrations, examples, and so on. And then, once we'll reach to the point where we are stunned by the beauty of this claim, and eager of looking for applications, and developing whole theories based on it, we will get of course motivated, do the math, and pull out a proof.

\bigskip

Let us begin with some simple observations. We first have:

\begin{proposition}
Given a closed subgroup $G\subset O_N$, set as before
$$C_{kl}=\left\{T\in\mathcal L(H^{\otimes k},H^{\otimes l})\Big|Tg^{\otimes k}=g^{\otimes l}T,\forall g\in G\right\}$$
where $H=\mathbb C^N$, and then set as in Claim 1.2:
$$\widetilde{G}=\left\{g\in O_N\Big|Tg^{\otimes k}=g^{\otimes l}T,\forall k,l,\forall T\in C_{kl}\right\}$$
Then $\widetilde{G}$ is closed subgroup of $O_N$, and we have inclusions $G\subset\widetilde{G}\subset O_N$.
\end{proposition}

\begin{proof}
Let us first prove that $\widetilde{G}$ is a group. Assuming $g,h\in\widetilde{G}$, we have $gh\in\widetilde{G}$, due to the following computation, valid for any $k,l$ and any $T\in C_{kl}$:
\begin{eqnarray*}
T(gh)^{\otimes k}
&=&Tg^{\otimes k}h^{\otimes k}\\
&=&g^{\otimes l}Th^{\otimes k}\\
&=&g^{\otimes l}h^{\otimes l}T\\
&=&(gh)^{\otimes l}T
\end{eqnarray*}

Also, we have $1\in\widetilde{G}$, trivially. Finally, assuming $g\in\widetilde{G}$, we have:
\begin{eqnarray*}
T(g^{-1})^{\otimes k}
&=&(g^{-1})^{\otimes l}[g^{\otimes l}T](g^{-1})^{\otimes k}\\
&=&(g^{-1})^{\otimes l}[Tg^{\otimes k}](g^{-1})^{\otimes k}\\
&=&(g^{-1})^{\otimes l}T
\end{eqnarray*}

Thus we have $g^{-1}\in\widetilde{G}$, and so $\widetilde{G}$ is a group, as claimed. Finally, the fact that we have an inclusion $G\subset\widetilde{G}$, and that $\widetilde{G}\subset O_N$ is closed, are both clear from definitions.
\end{proof}

Getting back now to Claim 1.2, what is going on there is that given a closed subgroup $G\subset O_N$, we can associated to it a sort of ``completion'' $\widetilde{G}$, which is an intermediate subgroup $G\subset\widetilde{G}\subset O_N$. And we want to prove that $G\subset\widetilde{G}$ is an equality.

\bigskip

This certainly looks a bit technical, but at least we are now into some familiar ground, that of abstract mathematics, so based on our knowledge here, we can start thinking about a strategy for proving our claim. And here, things are quite clear. Indeed, the relations $Tg^{\otimes k}=g^{\otimes l}T$ can be thought of as being some sort of ``commutation relations'', with the remark that at $k=l$ they are indeed usual commutation relations, namely:
$$[T,g^{\otimes k}]=0$$

Thus, what our claim basically says is that ``the commutant of the commutant must be the group itself''. Which leads right away, assuming some knowledge of pure mathematics, either advanced algebra, or some functional analysis, to the following thought:

\index{bicommutant theorem}
\index{von Neumann}

\begin{thought}
Our Claim 1.2 is some sort of bicommutant theorem for groups, and it is based on this that we should look for both its proof, and beauty and applications.
\end{thought}

We will see later that this thought is indeed correct, with one proof of Claim 1.2 coming from the bicommutant theorem of von Neumann, and with its beauty and applications coming in relation with the operator algebras introduced by the same von Neumann. Which is something very nice, for the story, John von Neumann being one of the last serious mathematicians, knowing well both math and physics, as we should all do.

\bigskip

This being said, what we have so far as theory, namely Definition 1.1, Claim 1.2, Proposition 1.3 and Thought 1.4, still remains something quite abstract. So, let us work out some examples. The orthogonal diagonal matrices form a subgroup $\mathbb Z_2^N\subset O_N$, and for the various subgroups $G\subset\mathbb Z_2^N$ our theory is quite exciting, and we have:

\index{abelian group}
\index{diagonal group}

\begin{theorem}
For the abelian groups of diagonal matrices, $G\subset\mathbb Z_2^N$, we have
$$C_{kl}=\left\{T\in\mathcal L(H^{\otimes k},H^{\otimes l})\Big|\exists g\in G,g_{i_1}\ldots g_{i_k}\neq g_{j_1}\ldots g_{j_l}\implies T_{j_1\ldots j_l,i_1\ldots i_k}=0\right\}$$
with the notation $g=diag(g_1,\ldots,g_N)$, and Claim 1.2 holds when $|G|=1,2,2^{N-1},2^N$.
\end{theorem}

\begin{proof}
We have several things to be proved, the idea being as follows:

\medskip

(1) Case $G=\{1\}$. Here we obviously have, for any two integers $k,l$, the following formula, which confirms the general formula in the statement:
$$C_{kl}=\mathcal L(H^{\otimes k},H^{\otimes l})$$

Regarding now Claim 1.2, consider the intermediate subgroup $G\subset\widetilde{G}\subset O_N$, constructed in Proposition 1.3, that we must prove to be equal to $G$ itself. Since any element $g\in\widetilde{G}$ must commute with the algebra $C_{11}=M_N(\mathbb C)$, we must have:
$$g=\pm1$$

But from the relation $T=gT$, which must hold for any $T\in C_{01}=H$, we conclude that we must have $g=1$, so we obtain $\widetilde{G}=\{1\}$, as desired.

\medskip

(2) Case $G=\mathbb Z_2$, with this meaning $G=\{1,-1\}$. This is something just a bit more complicated. Let us look at the relations defining $C_{kl}$, namely:
$$Tg^{\otimes k}=g^{\otimes l}T$$

These relations are automatic for $g=1$. As for the other group element, namely $g=-1$, here the relations hold either when $k+l$ is even, or when $T=0$. Thus, we have the following formula, which confirms again the general formula in the statement:
$$C_{kl}=
\begin{cases}
\mathcal L(H^{\otimes k},H^{\otimes l})&(k+l\in 2\mathbb N)\\
\{0\}&(k+l\notin 2\mathbb N)
\end{cases}$$

As for Claim 1.2 for our group, this follows from the computation done in (1) above, the point being that $g\in\widetilde{G}$ commutes with $C_{11}=M_N(\mathbb C)$ precisely when $g=\pm1$.

\medskip

(3) General case $G\subset\mathbb Z_2^N$. Let us look at the relations defining $C_{kl}$. We have:
\begin{eqnarray*}
T\in C_{kl}
&\iff&Tg^{\otimes k}=g^{\otimes l}T,\forall g\in G\\
&\iff&(Tg^{\otimes k})_{ji}=(g^{\otimes l}T)_{ji},\forall i,j,\forall g\in G\\
&\iff&T_{j_1\ldots,j_l,i_1\ldots i_k}g_{i_1}\ldots g_{i_k}=g_{j_1}\ldots g_{j_l}T_{j_1\ldots,j_k,i_1\ldots i_l},\forall i,j,\forall g\in G\\
&\iff&(g_{j_1}\ldots g_{i_k}-g_{j_1}\ldots g_{j_l})T_{j_1\ldots,j_l,i_1\ldots i_k},\forall i,j,\forall g\in G
\end{eqnarray*}

Thus, we are led to the formula in the statement, namely:
$$C_{kl}=\left\{T\in\mathcal L(H^{\otimes k},H^{\otimes l})\Big|\exists g\in G,g_{i_1}\ldots g_{i_k}\neq g_{j_1}\ldots g_{j_l}\implies T_{j_1\ldots j_l,i_1\ldots i_k}=0\right\}$$

(4) Case $G=\mathbb Z_2^N$. Here the formula from (3) can be turned into something better, because due to the fact that the entries $g_1,\ldots,g_N\in\{-1,1\}$ of a group element $g\in G$ can take all possible values, we have the following equivalence, with the symbol $\{\ \}_2$ standing for set with repetitions, with the pairs of elements of type $\{x,x\}$ removed:
$$g_{i_1}\ldots g_{i_k}=g_{j_1}\ldots g_{j_l},\forall g\in G\iff\{i_1,\ldots,i_k\}_2=\{j_1,\ldots,j_l\}_2$$

Thus, in this case we obtain the following formula, with $\{\ \}_2$ being as above:
$$C_{kl}=\left\{T\in\mathcal L(H^{\otimes k},H^{\otimes l})\Big|\{i_1,\ldots,i_k\}_2\neq\{j_1,\ldots,j_l\}_2\implies T_{j_1\ldots j_l,i_1\ldots i_k}=0\right\}$$

Regarding now Claim 1.2, the idea is that, a bit as for $G=\mathbb Z_2$, we can get away with the commutation with $C_{11}$. Indeed, according to the above formulae, we have:
$$C_{11}=\left\{T\in M_N(\mathbb C)\Big|i\neq j\implies T_{ij}=0\right\}$$

Thus we have $C_{11}=\Delta$, with $\Delta\subset M_N(\mathbb C)$ being the algebra of diagonal matrices. Now if we construct $G\subset\widetilde{G}\subset O_N$ as before, we have, as desired:
\begin{eqnarray*}
g\in\widetilde{G}
&\implies&g\in C_{11}'=\Delta'=\Delta\\
&\implies&g\in\Delta\cap O_N=G
\end{eqnarray*}

(5) Before getting into more examples, let us go back to the case where $G\subset\mathbb Z_2^N$ is arbitrary, and have a look at Claim 1.2 in this case. We know that we have $\{1\}\subset G\subset\mathbb Z_2^N$, and by functoriality, at the level of the associated $C_{11}$ spaces, we have:
$$\Delta\subset C_{11}\subset M_N(\mathbb C)$$

Now construct the intermediate group $G\subset\widetilde{G}\subset O_N$ as before. For $g\in\widetilde{G}$ we have:
$$g\in C_{11}'\cap O_N\subset\Delta'\cap O_N=\Delta\cap O_N=\mathbb Z_2^N$$

Thus, we have $G\subset\widetilde{G}\subset\mathbb Z_2^N$. This looks encouraging, because our Claim 1.2 becomes now something regarding the abelian groups, that can be normally solved with group theory. However, as we will soon discover, the combinatorics can be quite complicated.

\medskip

(6) General case $|G|=2$. This is the same as saying that $G\simeq\mathbb Z_2$, or equivalently, that $G=\{1,g\}$ with $g\in\mathbb Z_2$, $g\neq1$. By permuting the basis of $\mathbb R^N$ we can assume that our non-trivial group element $g\in G$ is as follows, for a certain integer $M<N$:
$$g=\begin{pmatrix}1_M&0\\ 0&-1_{N-M}\end{pmatrix}$$

By using the general formula found in (3), we obtain the following formula:
$$C_{11}=\left\{T\in M_N(\mathbb C)\Big|T_{ij}=0\ {\rm when}\ i\leq M,j>M\ {\rm or}\ i>M,j\leq M\right\}$$

But this means that, in this case, the algebra $C_{11}$ is block-diagonal, as follows:
$$C_{11}=\left\{ \begin{pmatrix}A&0\\0&B\end{pmatrix}\Big|A\in M_M(\mathbb C), B\in M_{N-M}(\mathbb C)\right\}$$

Now since any element $h\in\widetilde{G}$ must commute with this algebra, we must have:
$$\widetilde{G}\subset\left\{
\begin{pmatrix}1&0\\0&1\end{pmatrix}, 
\begin{pmatrix}1&0\\0&-1\end{pmatrix}, 
\begin{pmatrix}-1&0\\0&1\end{pmatrix}, 
\begin{pmatrix}-1&0\\0&-1\end{pmatrix} 
\right\}$$

Summarizing, well done, but we are still not there. In order to finish we must use, as in (1), the relations $T=hT$ with $T\in C_{01}$. In order to do so, by using again the general formula from (3), this time with $k=0,l=1$, we obtain the following formula:
$$C_{01}=\left\{T\in \mathbb C^N\Big|j>M\implies T_j=0\right\}$$

But this formula tells us that the space $C_{01}$ appears as follows:
$$C_{01}=\left\{ \binom{\xi}{0}\Big|\,\xi\in\mathbb C^M\right\}$$

Now since any element $h\in\widetilde{G}$ must satisfy $T=hT$, for any $T\in C_{01}$, this rules out half of the 4 solutions found above, and we end up with $\widetilde{G}=\{1,g\}$, as desired. 

\medskip

(7) A next step would be to investigate the case $|G|=4$. Here we have $G=\{1,g,h,gh\}$ with $g,h\in\mathbb Z^2-\{1\}$ distinct, and by permuting the basis, we can assume that:
$$g=\begin{pmatrix}1\\&1\\&&-1\\&&&-1\end{pmatrix}\quad,\quad
h=\begin{pmatrix}1\\&-1\\&&1\\&&&-1\end{pmatrix}\quad,\quad
gh=\begin{pmatrix}1\\&-1\\&&-1\\&&&1\end{pmatrix}$$

However, the computations as in the proof of (6) become quite complicated, and in addition we won't get away in this case with $C_{11},C_{01}$ only, so all this becomes too technically involved, and we will stop here, in the lack of a better idea.

\medskip

(8) Case $|G|=2^{N-1}$. This is the last situation, announced in the statement, still having a reasonably simple direct proof, and we will discuss this now. At the level of examples, given a non-empty subset $I\subset\{1,\ldots,N\}$, we have an example, as follows:
$$G_I=\left\{g\in\mathbb Z_2^N\Big|\prod_{i\in I}g_i=1\right\}$$

Indeed, this set $G_I\subset\mathbb Z_2^N$ is clearly a group, and since it is obtained by using one binary relation, namely $\prod_ig_i=\pm1$ being assumed to be $1$, the number of elements is:
$$|G_I|=\frac{|\mathbb Z_2^N|}{2}=\frac{2^N}{2}=2^{N-1}$$

Our claim now is that all the index 2 subgroups $G\subset\mathbb Z_2^N$ appear in this way. Indeed, by taking duals these subgroups correspond to the order 2 subgroups $H\subset\mathbb Z_2^N$, and since we must have $H=\{1,g\}$ with $g\neq1$, we have $2^N-1$ choices for such subgroups. But this equals the number of choices for a non-empty subset $I\subset\{1,\ldots,N\}$, as desired.

\medskip

(9) Case $|G|=2^{N-1}$, continuation. We know from the above that we have $G=G_I$, for a certain non-empty subset $I\subset\{1,\ldots,N\}$, and we must prove Claim 1.2 for this group. In order to do so, let us go back to the formula of $C_{kl}$ found in (4) for the group $\mathbb Z_2^N$. In the case of the subgroup $G_I\subset\mathbb Z_2^N$, which appears via the relation $\prod_ig_i=1$, that formula adapts as follows, with the symbol $\{\ \}_{2I}$ standing for set with repetitions, with the pairs of elements of type $\{x,x\}$ removed, and with the subsets equal to $I$ being removed too:
$$C_{kl}=\left\{T\in\mathcal L(H^{\otimes k},H^{\otimes l})\Big|\{i_1,\ldots,i_k\}_{2I}\neq\{j_1,\ldots,j_l\}_{2I}\implies T_{j_1\ldots j_l,i_1\ldots i_k}=0\right\}$$

In order to prove now Claim 1.2 for our group, we already know from (5) that we have $\widetilde{G}\subset\mathbb Z_2^N$. It is also clear that, given $h\in\widetilde{G}$, when using $T=hT$ with $T\in C_{01}$, or more generally $T=h^{\otimes l}T$ with $T\in C_{0l}$ at small values of $l\in\mathbb N$, we won't obtain anything new. However, at $l=|I|$ we do obtain a constraint, and since this constaint must cut the target group $\mathbb Z_2^N$ by at least half, we end up with $G=\widetilde{G}$, as desired.
\end{proof}

The proof of Theorem 1.5 contains many interesting computations, that are useful in everyday life, and among the many things that can be highlighted, we have:

\index{planar algebra}

\begin{fact}
The diagonal part of $C=(C_{kl})$, formed by the algebras
$$C_{kk}=\left\{T\in\mathcal L(H^{\otimes k})\Big|Tg^{\otimes k}=g^{\otimes k}T,\forall g\in G\right\}$$
does not determine $G$. For instance $G=\{1\},\mathbb Z_2$ are not distinguished by it.
\end{fact}

Obviously, this is something quite annoying, because there are countless temptations to use $\Delta C=(C_{kk})$ instead of $C$, for instance because the spaces $C_{kk}$ are algebras, and also, at a more advanced level, because $\Delta C$ is a planar algebra in the sense of Jones \cite{jo3}. But, we are not allowed to do this, at least in general. More on this later.

\bigskip

What we have so far is quite interesting, and suggests further working on our problem. Unfortunately, at the other end, where $G\subset O_N$ is big, things become fairly complicated, and the only result that we can state and prove with bare hands is:

\begin{proposition}
Our Claim 1.2 holds for $G=O_N$ itself, trivially.
\end{proposition}

\begin{proof}
For the orthogonal group $G=O_N$ itself we have indeed $\widetilde{G}=G$, due to the inclusions $G\subset \widetilde{G}\subset O_N$. Observe however that some mystery remains for this group $G=O_N$, because the spaces $C_{kl}$ do not look easy to compute. We will be back to this.
\end{proof}

As a conclusion now, we are definitely into interesting mathematics, and Claim 1.2 is definitely worth some attention, and a proof. So, time for a theorem about it:

\index{Tannakian duality}
\index{Peter-Weyl theory}

\begin{theorem}
Given a closed subgroup $G\subset O_N$, we have
$$G=\left\{g\in O_N\Big|Tg^{\otimes k}=g^{\otimes l}T,\forall k,l,\forall T\in C_{kl}\right\}$$
where $C=(C_{kl})$ is the associated Tannakian category.
\end{theorem}

\begin{proof}
We already know that this is something non-trivial. However, this can be proved by using either Peter-Weyl theory, or Tannakian duality, as follows:

\medskip

(1) Consider, as before in Proposition 1.3 and afterwards, the following set:
$$\widetilde{G}=\left\{g\in O_N\Big|Tg^{\otimes k}=g^{\otimes l}T,\forall k,l,\forall T\in C_{kl}\right\}$$

We know that $\widetilde{G}\subset O_N$ is a closed subgroup, and that  $G\subset\widetilde{G}$. Thus, we have an intermediate subgroup as follows, that we want to prove to be equal to $G$ itself:
$$G\subset\widetilde{G}\subset O_N$$

(2) In order to prove this, consider the Tannakian category of $\widetilde{G}$, namely:
$$\widetilde{C}_{kl}=\left\{T\in\mathcal L(H^{\otimes k},H^{\otimes l})\Big|Tg^{\otimes k}=g^{\otimes l}T,\forall g\in\widetilde{G}\right\}$$

By functoriality, from $G\subset\widetilde{G}$ we obtain $\widetilde{C}\subset C$. On the other hand, according to the definition of $\widetilde{G}$, we have $C\subset\widetilde{C}$. Thus, we have the following equality:
$$C=\widetilde{C}$$

(3) Assume now by contradiction that $G\subset\widetilde{G}$ is not an equality. Then, at the level of algebras of functions, the following quotient map is not an isomorphism either:
$$C(\widetilde{G})\to C(G)$$ 

On the other hand, we know from Peter-Weyl that we have decompositions as follows, with the sums being over all the irreducible unitary representations:
$$C(\widetilde{G})=\overline{\bigoplus}_{\pi\in Irr(\widetilde{G})}M_{\dim\pi}(\mathbb C)
\quad,\quad 
C(G)=\overline{\bigoplus}_{\nu\in Irr(G)}M_{\dim\nu}(\mathbb C)$$

Now observe that each unitary representation $\pi:\widetilde{G}\to U_K$ restricts into a certain representation $\pi':G\to U_K$. Since the quotient map $C(\widetilde{G})\to C(G)$ is not an isomorphism, we conclude that there is at least one representation $\pi$ satisfying:
$$\pi\in Irr(\widetilde{G})\quad,\quad \pi'\notin Irr(G)$$

(4) We are now in position to conclude. By using Peter-Weyl theory again, the above representation $\pi\in Irr(\widetilde{G})$ appears in a certain tensor power of the fundamental representation $u:\widetilde{G}\subset U_N$. Thus, we have inclusions of representations, as follows:
$$\pi\in u^{\otimes k}\quad,\quad\pi'\in u'^{\otimes k}$$

Now since we know that $\pi$ is irreducible, and that $\pi'$ is not, by using one more time Peter-Weyl theory, we conclude that we have a strict inequality, as follows:
$$\dim(\widetilde{C}_{kk})=dim(End(u^{\otimes k}))<dim(End(u'^{\otimes k}))=
\dim(C_{kk})$$ 

But this contradicts the equality $C=\widetilde{C}$ found in (2), which finishes the proof.

\medskip

(5) Alternatively, we can use Tannakian duality. This duality states that any compact group $G$ appears as the group of endomorphisms of the canonical inclusion functor $Rep(G)\subset\mathcal H$, where $Rep(G)$ is the category of final dimensional continuous unitary representations of $G$, and $\mathcal H$ is the category of finite dimensional Hilbert spaces. 

\medskip

(6) Now in the case of a closed subgroup $G\subset_uO_N$, we know from Peter-Weyl theory that any $r\in Rep(G)$ appears as a subrepresentation $r\in u^{\otimes k}$. In categorical terms, this means that, with suitable definitions, $Rep(G)$ appears as a ``completion'' of the category $C=(C_{kl})$. Thus $C$ uniquely determines $G$, and we obtain the result.
\end{proof}

All the above was of course quite brief, but we will be back to this topic, and to Tannakian duality in general, on numerous occasions, in what follows.

\section*{1b. Reflection groups}

What we have so far is quite interesting, and a first way of extending our knowledge is by looking at the unitary case. We have the following extension of Definition 1.1:

\index{Tannakian category}
\index{colored integer}
\index{tensor category}

\begin{definition}
Associated to any closed subgroup $G\subset U_N$ are the vector spaces
$$C_{kl}=\left\{T\in\mathcal L(H^{\otimes k},H^{\otimes l})\Big|Tg^{\otimes k}=g^{\otimes l}T,\forall g\in G\right\}$$
where $H=\mathbb C^N$, and with $k,l=\circ\bullet\bullet\circ\ldots$ being colored integers, with the conventions
$$g^{\otimes\circ}=g\quad,\quad g^{\otimes\bullet}=\bar{g}\quad,\quad g^{\otimes\emptyset}=1$$
and multiplicativity. We call Tannakian category of $G$ the collection of spaces $C=(C_{kl})$.
\end{definition}

As before in the real case, any compact Lie group $G\subset U_N$ can be reconstructed from its Tannakian category $C=(C_{kl})$, in a very simple way, as follows:

\index{Tannakian duality}

\begin{theorem}
Given a closed subgroup $G\subset U_N$, we have
$$G=\left\{g\in U_N\Big|Tg^{\otimes k}=g^{\otimes l}T,\forall k,l,\forall T\in C_{kl}\right\}$$
where $C=(C_{kl})$ is the associated Tannakian category.
\end{theorem}

\begin{proof}
This follows indeed as in the real case, by using either Peter-Weyl theory, or Tannakian duality, as explained in the proof of Theorem 1.8.
\end{proof}

Regarding the basic examples, we have here the following result:

\index{abelian group}
\index{diagonal group}

\begin{proposition}
For the abelian groups of diagonal matrices, $G\subset\mathbb T^N$, we have
$$C_{kl}=\left\{T\in\mathcal L(H^{\otimes k},H^{\otimes l})\Big|\exists g\in G,g_{i_1}^{k_1}\ldots g_{i_r}^{k_r}\neq g_{j_1}^{l_1}\ldots g_{j_s}^{l_s}\implies T_{j_1\ldots j_l,i_1\ldots i_k}=0\right\}$$
where $r=|k|$, $s=|l|$, and with $g=diag(g_1,\ldots,g_N)$. Also, Claim 1.2 holds.
\end{proposition}

\begin{proof}
The formula in the statement is something that we know from Theorem 1.5 in the real case, and the proof in the complex case is similar. As for the last assertion, this is something that we know to hold, from Theorem 1.10, but some explicit computations for special subgroups $G\subset\mathbb T^N$, along the lines of those from the proof of Theorem 1.5, can be quite instructive, and we will leave this as an interesting exercise.
\end{proof}

The question is now, can we get beyond this, with some further explicit verifications of Claim 1.2, in the orthogonal case, or in the general unitary case. And the problem is that this is something quite tricky, because as already mentioned, for $G=O_N$ itself, and in fact for $G=U_N$ itself too, the computation of the spaces $C_{kl}$ is something quite complicated, which cannot be done with bare hands, and we will leave this for later.

\bigskip

Fortunately the symmetric group $S_N\subset O_N$ comes to the rescue, and we have:

\index{symmetric group}
\index{partitions}
\index{Kronecker symbols}
\index{fitting indices}

\begin{theorem}
For the symmetric group $S_N\subset O_N$ we have the formula
$$C_{kl}=span\left(T_\pi\Big|\pi\in P(k,l)\right)$$
with $P(k,l)$ being the set of partitions of $k$ upper points, and $l$ lower points, and where
$$T_\pi(e_{i_1}\otimes\ldots\otimes e_{i_k})=\sum_{j_1\ldots j_l}\delta_\pi
\begin{pmatrix}i_1&\ldots&i_k\\ j_1&\ldots&j_l\end{pmatrix}
e_{j_1}\otimes\ldots\otimes e_{j_l}$$
where $\delta_\pi\in\{0,1\}$ is $1$ when the indices fit, and is $0$ otherwise. Also, Claim 1.2 holds.
\end{theorem}

\begin{proof}
This is something nice and elementary, and, importantly for us, fundamental for what is to follow in this book, the idea being as follows:

\medskip

(1) First, at the level of notations, we denote as mentioned by $P(k,l)$ the set of partitions of $k$ upper points, and $l$ lower points, and with these partitions being represented as pictures, with the blocks being represented by strings. As an example, here is a partition in $P(3,3)$, with two blocks, represented by strings, in the obvious way:
$$
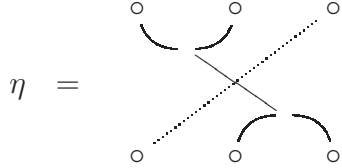
\xymatrix@R=2mm@C=3mm{\\ \\ \eta\ \ =\\ \\}\ \ \ 
\xymatrix@R=2mm@C=3mm{
\circ\ar@/_/@{-}[dr]&&\circ&&\circ\ar@{.}[ddddllll]\\
&\ar@/_/@{-}[ur]\ar@{-}[ddrr]\\
\\
&&&\ar@/^/@{-}[dr]\\
\circ&&\circ\ar@/^/@{-}[ur]&&\circ}$$

(2) Now given $\pi\in P(k,l)$ and multi-indices $i=(i_1,\ldots,i_k)$ and $j=(j_1,\ldots,j_l)$, we can put $i,j$ on the legs of $\pi$, in the obvious way. Then, if the ``indices fit'', meaning that all the strings of $\pi$ join equal indices of $i,j$, we set $\delta_\pi\binom{i}{j}=1$. Otherwise, we set $\delta_\pi\binom{i}{j}=0$. As an example, for the above partition $\eta\in P(3,3)$, we have the following formula:
$$\delta_\eta\begin{pmatrix}a&b&c\\d&e&f\end{pmatrix}=\delta_{abef}\delta_{cd}$$

(3) In order to prove now the result, let us first work out the case $k=0$, that we will regularly need in what follows. It is traditional here, and convenient, to change a bit notations. So, let us associate to any $\pi\in P(k)$ a vector, as follows:
$$\xi_\pi=\sum_{i_1\ldots i_k}\delta_\pi(i_1 \ldots i_k)e_{i_1}\otimes\ldots\otimes e_{i_k}$$

With this notation, we must prove that we have the following equality:
$$C_{0k}=span\left(\xi_\pi\Big|\pi\in P(k)\right)$$

(4) Let us first prove that we have $\supset$. Given $\sigma\in S_N$, we have indeed:
\begin{eqnarray*}
\sigma^{\otimes k}\xi_\pi
&=&\sum_{i_1\ldots i_k}\delta_\pi(i_1 \ldots i_k)\sigma(e_{i_1})\otimes\ldots\otimes\sigma(e_{i_k})\\
&=&\sum_{i_1\ldots i_k}\delta_\pi(i_1 \ldots i_k)e_{\sigma(i_1)}\otimes\ldots\otimes e_{\sigma(i_k)}\\
&=&\sum_{i_1\ldots i_k}\delta_\pi(\sigma^{-1}(j_1)\ldots\sigma^{-1}(j_k))e_{j_1}\otimes\ldots\otimes e_{j_k}\\
&=&\sum_{i_1\ldots i_k}\delta_\pi(j_1\ldots j_k)e_{j_1}\otimes\ldots\otimes e_{j_k}\\
&=&\xi_\pi
\end{eqnarray*}

(5) In order to prove now $\subset$, consider an arbitrary vector of $\mathbb C^N$, as follows:
$$\xi=\sum_{i_1\ldots i_k}\lambda_{i_1\ldots i_k}e_{i_1}\otimes\ldots\otimes e_{i_k}$$

Given $\sigma\in S_N$, by reasoning as before, we have the following formula:
$$\sigma^{\otimes k}\xi=\sum_{i_1\ldots i_k}\lambda_{\sigma^{-1}(i_1)\ldots\sigma^{-1}(i_k)}e_{i_1}\otimes\ldots\otimes e_{i_k}$$

Thus the condition $\sigma^{\otimes k}\xi=\xi$ for any $\sigma\in S_N$ is equivalent to:
$$\lambda_{i_1\ldots i_k}=\lambda_{\sigma(i_1)\ldots\sigma(i_k)}\quad,\quad\forall i,\sigma$$

But this latter condition is equivalent to the following condition:
$$\ker i=\ker j\implies\lambda_i=\lambda_j$$

Thus, we are led to the conclusion that $\lambda:\{1,\ldots,N\}^k\to\mathbb C$ must come from a function $\varphi:P(k)\to\mathbb C$, via a formula of type $\lambda_i=\varphi(\ker i)$, and it follows that the inclusion $\supset$ that we established in (4) is indeed an equality, as desired.

\medskip

(6) Summarizing, and getting back now to our theorem as stated, we have proved the formula of $C_{kl}$ there, in the case $k=0$. In order to pass now to the general case, two methods are available. We can either fine-tune the above computation, and we will leave this as an instructive exercise, or we can argue that the result at $k=0$ gives the result in general, via Frobenius duality, and we will leave this as an instructive exercise too.

\medskip

(7) Regarding now Claim 1.2, that we definitely know to hold from Theorem 1.8, but that we would like to prove now explicitely, consider the intermediate subgroup $S_N\subset\widetilde{S}_N\subset O_N$, constructed as in Proposition 1.3, that we must prove to be equal to $S_N$. In order to prove this equality, we use the following one-block ``fork'' partition:
$$\xymatrix@R=1mm@C=2mm{\\ \\ \mu\ \ =\\ \\ }\ \ \ 
\xymatrix@R=2mm@C=3mm{
\circ\ar@/_/@{-}[dr]&&\circ\\
&\ar@/_/@{-}[ur]\ar@{-}[dd]\\
&&&\\
&\circ}$$

The linear map associated to $\mu$ is then given by the following formula:
$$T_\mu(e_i\otimes e_j)=\delta_{ij}e_i$$

We therefore have the following formula, valid for any $g\in O_N$:
$$(T_\mu g^{\otimes 2})_{i,jk}
=\sum_{lm}(T_\mu)_{i,lm}(g^{\otimes 2})_{lm,jk}
=g_{ij}g_{ik}$$

On the other hand, we have as well the following formula:
$$(gT_\mu)_{i,jk}
=\sum_lg_{il}(T_\mu)_{l,jk}
=\delta_{jk}g_{ij}$$

Thus, we have the following equivalence, valid for any $g\in O_N$:
$$T_\mu g^{\otimes 2}=gT_\mu\iff g_{ij}g_{ik}=\delta_{jk}g_{ij},\forall i,j,k$$

Now by assuming $g\in\widetilde{S}_N$, the formula on the left holds, so the formula on the right must hold too. But this shows that we must have $g_{ij}\in\{0,1\}$, with exactly one 1 entry on each of the rows of $g$. Thus we must have $g\in S_N$, which finishes the proof. 
\end{proof}

The above result is quite encouraging, and suggests looking into other reflection groups. For $H_N$, which is the group of symmetries of the unit cube in $\mathbb R^N$, we have:

\index{hyperoctahedral group}
\index{hypercube}
\index{blocks of even size}

\begin{theorem}
For the hyperoctadedral group $H_N=\mathbb Z_2\wr S_N\subset O_N$ we have
$$C_{kl}=span\left(T_\pi\Big|\pi\in P_{even}(k,l)\right)$$
where $P_{even}$ means partitions all whose blocks have even size. Also, Claim 1.2 holds.
\end{theorem}

\begin{proof}
This follows a bit as for $S_N$. Consider indeed a vector of $\mathbb C^N$, as follows:
$$\xi=\sum_{i_1\ldots i_k}\lambda_{i_1\ldots i_k}e_{i_1}\otimes\ldots\otimes e_{i_k}$$

Then the condition $g^{\otimes k}\xi=\xi$ for any $g=\sigma^w\in H_N$ is equivalent to:
$$\lambda_{i_1\ldots i_k}=w_{i_1}\ldots w_{i_k}\lambda_{\sigma(i_1)\ldots\sigma(i_k)}\quad,\quad\forall i,\sigma$$

But this latter condition is equivalent to the following condition, along with the fact that we must have $w_{i_1}\ldots w_{i_k}=1$, which amounts in saying that $\ker i\in P_{even}$:
$$\ker i=\ker j\implies\lambda_i=\lambda_j$$

Thus, we are led to the conclusion in the statement. Finally, regarding the explicit verification of Claim 1.2, this follows again as for $S_N$, by using the following partition:
$$\xymatrix@R=0.5mm@C=2mm{\\ \\ \\ \chi\ \ =\\ \\ }\ \ \ 
\xymatrix@R=2mm@C=3mm{
\circ\ar@/_/@{-}[dr]&&\circ\\
&\ar@/_/@{-}[ur]\ar@{-}[dd]\\
&&&\\
&\ar@/^/@{-}[dr]\ar@/_/@{-}[dl]\\
\circ&&\circ}$$

And we will leave the verifications here, which are elementary, as an exercice.
\end{proof}

Quite remarkably, the results for $S_N,H_N$ can be generalized as follows:

\index{wreath product}
\index{complex reflection group}

\begin{theorem}
For the complex reflection group $H_N^s=\mathbb Z_s\wr S_N\subset U_N$ we have
$$C_{kl}=span\left(T_\pi\Big|\pi\in P^s(k,l)\right)$$
where $P^s(k,l)$ is the set of partitions of $k$ upper points and $l$ lower points, satisfying
$$\#\circ=\#\bullet(s)$$
as a weighted equality, in each block. Also, Claim 1.2 holds.
\end{theorem}

\begin{proof}
This follows a bit as for $S_N,H_N$. Consider indeed a vector, as follows:
$$\xi=\sum_{i_1\ldots i_k}\lambda_{i_1\ldots i_k}e_{i_1}\otimes\ldots\otimes e_{i_k}$$

Then the condition $g^{\otimes k}\xi=\xi$ for any $g=\sigma^w\in H_N^s$ is equivalent to:
$$\lambda_{i_1\ldots i_k}=w_{i_1}\ldots w_{i_k}\lambda_{\sigma(i_1)\ldots\sigma(i_k)}\quad,\quad\forall i,\sigma$$

But this latter condition is equivalent to the following condition, along with the fact that we must have $w_{i_1}\ldots w_{i_k}=1$, which amounts in saying that $\ker i\in P^s(k)$:
$$\ker i=\ker j\implies\lambda_i=\lambda_j$$

Thus, we are led to the conclusion in the statement.
\end{proof}

Finally, let us record the $s=\infty$ particular case of Theorem 1.14, as follows:

\index{full reflection group}

\begin{theorem}
For the full complex reflection group $K_N=\mathbb T\wr S_N\subset U_N$ we have
$$C_{kl}=span\left(T_\pi\Big|\pi\in\mathcal P_{even}(k,l)\right)$$
where $\mathcal P_{even}(k,l)$ is the set of partitions of $k$ upper points and $l$ lower points, satisfying
$$\#\circ=\#\bullet$$
as a weighted equality, in each block. Also, Claim 1.2 holds.
\end{theorem}

\begin{proof}
This appears indeed as the $s=\infty$ particular case of Theorem 1.14.
\end{proof}

\section*{1c. Rotation groups}

We have so far a beginning of theory based on Claim 1.2 and its philosophy, with the main examples being the reflection groups $G=\mathbb Z_s\wr S_N$. We would like now to look into the continuous groups $G\subset U_N$. Let us start with a basic investigation of the simplest such group, namely $O_N\subset U_N$. We can say the following, about this group:

\begin{proposition}
For the orthogonal group $O_N$ we have inclusions as follows,
$$span\left(T_\pi\Big|\pi\in P_2(k,l)\right)\subset C_{kl}\subset span\left(T_\pi\Big|\pi\in P(k,l)\right)$$
where $P_2(k,l)$ is the set of pairings of $k$ upper points, and $l$ lower points.
\end{proposition}

\begin{proof}
Since we have $S_N\subset O_N$, by functoriality and Theorem 1.12 we have:
$$C_{kl}\subset span\left(T_\pi\Big|\pi\in P(k,l)\right)$$

For the other inclusion, the one on the left, let us first work out, as usual, the case $k=0$, and with the change $k\leftrightarrow l$. For a pairing $\pi\in P_2(k)$ we set, as before:
$$\xi_\pi=\sum_{i_1\ldots i_k}\delta_\pi(i_1 \ldots i_k)e_{i_1}\otimes\ldots\otimes e_{i_k}$$

We must prove $\xi_\pi\in C_{0k}$. For this purpose, let us pick $g\in O_N$, and write:
$$g(e_i)=\sum_jg_{ji}e_j$$

We have then the following computation:
\begin{eqnarray*}
g^{\otimes k}\xi_\pi
&=&\sum_{i_1\ldots i_k}\delta_\pi(i_1 \ldots i_k)ge_{i_1}\otimes\ldots\otimes ge_{i_k}\\
&=&\sum_{i_1\ldots i_k}\sum_{j_1\ldots j_k}\delta_\pi(i_1 \ldots i_k)g_{j_1i_1}\ldots g_{j_ki_k}e_{j_1}\otimes\ldots\otimes e_{j_k}
\end{eqnarray*}

As an illustration now, let us see what happens for a simple pairing, such as $\pi=\cap\,\cap$. Here the above computation can be continued as follows:
\begin{eqnarray*}
g^{\otimes k}\xi_{\cap\,\cap}
&=&\sum_{i_1\ldots i_4}\sum_{j_1\ldots j_4}\delta_{\cap\,\cap}(i_1 \ldots i_4)g_{j_1i_1}\ldots g_{j_4i_4}e_{j_1}\otimes\ldots\otimes e_{j_4}\\
&=&\sum_{ab}\sum_{j_1\ldots j_4}g_{j_1a}g_{j_2a}g_{j_3b}g_{j_4b}e_{j_1}\otimes\ldots\otimes e_{j_4}\\
&=&\sum_{j_1\ldots j_4}\delta_{j_1j_2}\delta_{j_3j_4}e_{j_1}\otimes\ldots\otimes e_{j_4}\\
&=&\xi_{\cap\,\cap}
\end{eqnarray*}

The same computation works in general, and by using $gg^t=1$, we obtain $g^{\otimes k}\xi_\pi=\xi_\pi$, for any $\pi\in P_2(k)$. Thus, we have indeed inclusions as in the statement.
\end{proof}

The above result, perhaps coupled with a few more computations, that you can do by yourself, in order to evaluate the situation, suggests that we should have:
$$C_{kl}=span\left(T_\pi\Big|\pi\in P_2(k,l)\right)$$

However, this is hard to prove with bare hands, and we will have to trick. Our trick will be something quite natural, the idea being first to prove that the above spaces $span(T_\pi)$ ``qualify'' for what is expected from a Tannakian category, and then, conclude that we have equality, because $span(T_\pi)$ can only correspond to $O_N$. Let us start with:

\index{tensor category}
\index{flip operator}
\index{colored integer}

\begin{definition}
A tensor category over $H=\mathbb C^N$ is a collection $C=(C_{kl})$ of linear spaces $C_{kl}\subset\mathcal L(H^{\otimes k},H^{\otimes l})$ satisfying the following conditions:
\begin{enumerate}
\item $S,T\in C$ implies $S\otimes T\in C$.

\item If $S,T\in C$ are composable, then $ST\in C$.

\item $T\in C$ implies $T^*\in C$.

\item Each $C_{kk}$ contains the identity operator.

\item $C_{\emptyset k}$ with $k=\circ\bullet,\bullet\circ$ contain the operator $R:1\to\sum_ie_i\otimes e_i$.

\item $C_{kl,lk}$ with $k,l=\circ,\bullet$ contain the flip operator $\Sigma:a\otimes b\to b\otimes a$.
\end{enumerate}
\end{definition}

Here, as usual, the tensor powers $H^{\otimes k}$, which are Hilbert spaces depending on a colored integer $k=\circ\bullet\bullet\circ\ldots\,$, are defined by the following formulae, and multiplicativity:
$$H^{\otimes\emptyset}=\mathbb C\quad,\quad 
H^{\otimes\circ}=H\quad,\quad
H^{\otimes\bullet}=\bar{H}\simeq H$$

We have already met such categories, when dealing with the Tannakian categories of the closed subgroups $G\subset U_N$, and our knowledge can be summarized as follows:

\index{tensor category}
\index{Tannakian category}

\begin{proposition}
Given a closed subgroup $G\subset U_N$, its Tannakian category
$$C_{kl}=\left\{T\in\mathcal L(H^{\otimes k},H^{\otimes l})\Big|Tg^{\otimes k}=g^{\otimes l}T,\forall g\in G\right\}$$
is a tensor category over $H=\mathbb C^N$. Conversely, given a tensor category $C$ over $\mathbb C^N$,
$$G=\left\{g\in U_N\Big|Tg^{\otimes k}=g^{\otimes l}T,\forall k,l,\forall T\in C_{kl}\right\}$$
is a closed subgroup of $U_N$.
\end{proposition}

\begin{proof}
This is something that we basically know, the idea being as follows:

\medskip

(1) Regarding the first assertion, we have to check here the axioms (1-6) in Definition 1.17. The axioms (1-4) being all clear from definitions, let us establish (5). But this follows from the fact that each element $g\in G$ is a unitary, which can be reformulated as follows, with $R:1\to\sum_ie_i\otimes e_i$ being the map in Definition 1.17:
$$R\in Hom(1,g\otimes\bar{g})\quad,\quad 
R\in Hom(1,\bar{g}\otimes g)$$

Regarding now the condition in Definition 1.17 (6), this comes from the fact that the matrix coefficients $g\to g_{ij}$ and their conjugates $g\to\bar{g}_{ij}$ commute with each other.

\medskip

(2) Finally, the last assertion is clear from definitions, with the verifications being almost identical to those in the proof of Proposition 1.3.
\end{proof}

Summarizing, we have so far precise axioms for the tensor categories $C=(C_{kl})$, given in Definition 1.17, as well as correspondences as follows:
$$G\to C_G\quad,\quad 
C\to G_C$$

We will prove in what follows that these correspondences are inverse to each other. In order to get started, we first have the following technical result:

\begin{proposition}
Consider the following conditions:
\begin{enumerate}
\item $C=C_{G_C}$, for any tensor category $C$.

\item $G=G_{C_G}$, for any closed subgroup $G\subset U_N$.
\end{enumerate}
We have then $(1)\implies(2)$. Also, $C\subset C_{G_C}$ is automatic.
\end{proposition}

\begin{proof}
Given $G\subset U_N$, we have $G\subset G_{C_G}$. On the other hand, by using (1) we have $C_G=C_{G_{C_G}}$. Thus, we have an inclusion of closed subgroups of $U_N$, which becomes an isomorphism at the level of the associated Tannakian categories, so $G=G_{C_G}$. Finally, the fact that we have an inclusion $C\subset C_{G_C}$ is clear from definitions.
\end{proof}

In order to establish Tannakian duality, we will need some abstract constructions. Following Malacarne \cite{mal}, let us start with the following elementary fact:

\begin{proposition}
Given a tensor category $C=(C_{kl})$ over a Hilbert space $H$,
$$E_C
=\bigoplus_{k,l}C_{kl}
\subset\bigoplus_{k,l}B(H^{\otimes k},H^{\otimes l})
\subset B\left(\bigoplus_kH^{\otimes k}\right)$$
is a closed $*$-subalgebra. Also, inside this algebra,
$$E_C^{(s)}
=\bigoplus_{|k|,|l|\leq s}C_{kl}
\subset\bigoplus_{|k|,|l|\leq s}B(H^{\otimes k},H^{\otimes l})
=B\left(\bigoplus_{|k|\leq s}H^{\otimes k}\right)$$
is a finite dimensional $*$-subalgebra.
\end{proposition}

\begin{proof}
This is clear indeed from the categorical axioms from Definition 1.17.
\end{proof}

By using now the bicommutant theorem of von Neumann, we have:

\begin{proposition}
Given a Tannakian category $C$, the following are equivalent:
\begin{enumerate}
\item $C=C_{G_C}$.

\item $E_C=E_{C_{G_C}}$.

\item $E_C^{(s)}=E_{C_{G_C}}^{(s)}$, for any $s\in\mathbb N$.

\item $E_C^{(s)'}=E_{C_{G_C}}^{(s)'}$, for any $s\in\mathbb N$.
\end{enumerate}
In addition, the inclusions $\subset$, $\subset$, $\subset$, $\supset$ are automatically satisfied.
\end{proposition}

\begin{proof}
The equivalences are clear from definitions, and from the bicommutant theorem. As for the last assertion, we have $C\subset C_{G_C}$, which shows that we have as well: 
$$E_C\subset E_{C_{G_C}}$$

We therefore obtain by truncating $E_C^{(s)}\subset E_{C_{G_C}}^{(s)}$, and by taking the commutants, this gives $E_C^{(s)}\supset E_{C_{G_C}}^{(s)}$. Thus, we are led to the conclusion in the statement.
\end{proof}

Summarizing, we would like to prove that we have $E_C^{(s)'}\subset E_{C_{G_C}}^{(s)'}$. But this can be done by doing some algebra, and we are led to the following conclusion:

\index{Tannakian duality}

\begin{theorem}
The Tannakian duality constructions 
$$C\to G_C\quad,\quad 
G\to C_G$$
are inverse to each other.
\end{theorem}

\begin{proof}
This follows by doing some algebra, in order to prove that we have indeed $E_C^{(s)'}\subset E_{C_{G_C}}^{(s)'}$, as mentioned above, and we refer here to the paper of Malacarne \cite{mal}. Alternatively, this can be proved via standard Tannakian methods, and we refer here to the paper of Woronowicz \cite{wo2}. For more on all this, you have as well my book \cite{ba1}.
\end{proof}

With this piece of general theory in hand, let us go back now to the orthogonal group $O_N$, and to partitions and pairings, as in Proposition 1.16. In order to construct a Tannakian category out of the pairings, via the operation $\pi\to T_\pi$, we will need:

\begin{proposition}
The assignement $\pi\to T_\pi$ is categorical, in the sense that
$$T_\pi\otimes T_\sigma=T_{[\pi\sigma]}\quad,\quad 
T_\pi T_\sigma=N^{c(\pi,\sigma)}T_{[^\sigma_\pi]}\quad,\quad 
T_\pi^*=T_{\pi^*}$$
where $c(\pi,\sigma)$ are certain integers, coming from the erased components in the middle.
\end{proposition}

\begin{proof}
This is something elementary, the computations being as follows:

\medskip

(1) The concatenation axiom follows from the following computation:
\begin{eqnarray*}
&&(T_\pi\otimes T_\sigma)(e_{i_1}\otimes\ldots\otimes e_{i_p}\otimes e_{k_1}\otimes\ldots\otimes e_{k_r})\\
&=&\sum_{j_1\ldots j_q}\sum_{l_1\ldots l_s}\delta_\pi\begin{pmatrix}i_1&\ldots&i_p\\j_1&\ldots&j_q\end{pmatrix}\delta_\sigma\begin{pmatrix}k_1&\ldots&k_r\\l_1&\ldots&l_s\end{pmatrix}e_{j_1}\otimes\ldots\otimes e_{j_q}\otimes e_{l_1}\otimes\ldots\otimes e_{l_s}\\
&=&\sum_{j_1\ldots j_q}\sum_{l_1\ldots l_s}\delta_{[\pi\sigma]}\begin{pmatrix}i_1&\ldots&i_p&k_1&\ldots&k_r\\j_1&\ldots&j_q&l_1&\ldots&l_s\end{pmatrix}e_{j_1}\otimes\ldots\otimes e_{j_q}\otimes e_{l_1}\otimes\ldots\otimes e_{l_s}\\
&=&T_{[\pi\sigma]}(e_{i_1}\otimes\ldots\otimes e_{i_p}\otimes e_{k_1}\otimes\ldots\otimes e_{k_r})
\end{eqnarray*}

(2) The composition axiom follows from the following computation:
\begin{eqnarray*}
&&T_\pi T_\sigma(e_{i_1}\otimes\ldots\otimes e_{i_p})\\
&=&\sum_{j_1\ldots j_q}\delta_\sigma\begin{pmatrix}i_1&\ldots&i_p\\j_1&\ldots&j_q\end{pmatrix}
\sum_{k_1\ldots k_r}\delta_\pi\begin{pmatrix}j_1&\ldots&j_q\\k_1&\ldots&k_r\end{pmatrix}e_{k_1}\otimes\ldots\otimes e_{k_r}\\
&=&\sum_{k_1\ldots k_r}N^{c(\pi,\sigma)}\delta_{[^\sigma_\pi]}\begin{pmatrix}i_1&\ldots&i_p\\k_1&\ldots&k_r\end{pmatrix}e_{k_1}\otimes\ldots\otimes e_{k_r}\\
&=&N^{c(\pi,\sigma)}T_{[^\sigma_\pi]}(e_{i_1}\otimes\ldots\otimes e_{i_p})
\end{eqnarray*}

(3) Finally, the involution axiom follows from the following computation:
\begin{eqnarray*}
&&T_\pi^*(e_{j_1}\otimes\ldots\otimes e_{j_q})\\
&=&\sum_{i_1\ldots i_p}<T_\pi^*(e_{j_1}\otimes\ldots\otimes e_{j_q}),e_{i_1}\otimes\ldots\otimes e_{i_p}>e_{i_1}\otimes\ldots\otimes e_{i_p}\\
&=&\sum_{i_1\ldots i_p}\delta_\pi\begin{pmatrix}i_1&\ldots&i_p\\ j_1&\ldots& j_q\end{pmatrix}e_{i_1}\otimes\ldots\otimes e_{i_p}\\
&=&T_{\pi^*}(e_{j_1}\otimes\ldots\otimes e_{j_q})
\end{eqnarray*}

Summarizing, our correspondence is indeed categorical.
\end{proof}

Good news, we can now prove the Brauer theorem for $O_N$, as follows:

\index{pairings}
\index{orthogonal group}
\index{Brauer theorem}

\begin{theorem}
For the orthogonal group $O_N$ we have
$$C_{kl}=span\left(T_\pi\Big|\pi\in P_2(k,l)\right)$$
where $P_2(k,l)$ is the set of pairings of $k$ upper points, and $l$ lower points.
\end{theorem}

\begin{proof}
We know from Proposition 1.16 that we have inclusions as follows:
$$C_{kl}\supset span\left(T_\pi\Big|\pi\in P_2(k,l)\right)$$

On the other hand, Proposition 1.23 shows that the spaces on the right form a Tannakian category, in the sense of Definition 1.17. Thus the Tannakian duality result from Theorem 1.22 applies, and provides us with a certain subgroup $G\subset U_N$, such that:
$$G=\left\{g\in U_N\Big|T_\pi g^{\otimes k}=g^{\otimes l}T_\pi,\forall k,l,\forall \pi\in P_2(k,l)\right\}$$

Moreover, by functoriality of Tannakian duality, we have $O_N\subset G$. But the relation $g^{\otimes 2}T_\cap=T_\cap$ with $g\in U_N$ implies $g\in O_N$, as explained in the proof of Proposition 1.18, so we have as well $G\subset O_N$. Thus, we have $G=O_N$, which gives the result.
\end{proof}

For the unitary group $U_N$ now, the result is similar, as follows:

\index{unitary group}
\index{matching pairings}
\index{Brauer theorem}

\begin{theorem}
For the unitary group $U_N$ we have
$$C_{kl}=span\left(T_\pi\Big|\pi\in\mathcal P_2(k,l)\right)$$
where $\mathcal P_2(k,l)$ is the set of matching pairings of $k$ upper points, and $l$ lower points.
\end{theorem}

\begin{proof}
The proof here is very similar to the proof for $O_N$, and in fact, even a bit simpler, after a close examination, and with the convention, in the statement, that matching means joining $\circ-\circ$ or $\bullet-\bullet$ on the horizontal, and $\circ-\bullet$ on the vertical.
\end{proof}

As a comment here, the above proofs might seem quite wizarding, and you may wonder where the computations, which are usually needed for proving such things, have dissapeared. In answer, we have used Theorem 1.22, and Proposition 1.23 as well.

\section*{1d. Easy groups}

We have now a substantial number of results based on Claim 1.2 and its philosophy, so time for some axiomatics. We can formulate a key definition, as follows:

\index{easy group}

\begin{definition}
A closed subgroup $G\subset U_N$ is called easy when
$$C_{kl}=span\left(T_\pi\Big|\pi\in D(k,l)\right)$$
with $D(k,l)\subset P(k,l)$ being certain sets of partitions.
\end{definition}

We already know, from our various computations performed above, that many interesting closed subgroups $G\subset U_N$ are easy. In fact, at the level of examples, we have the following result, which summarizes our main results so far, in this book:

\begin{theorem}
The following closed subgroups $G\subset U_N$ are easy,
\begin{enumerate}
\item $U_N$ itself, coming from $D=\mathcal P_2$,

\item $O_N$, coming from $D=P_2$,

\item $H_N^s$, coming from $D=P^s$,
\end{enumerate}
with the last result covering $S_N,H_N,K_N$, which appear from $D=P,P_{even},\mathcal P_{even}$.
\end{theorem}

\begin{proof}
This is indeed a reformulation of our main results so far:

\medskip

-- The results regarding $O_N,U_N$ are reformulations of Theorems 1.24 and 1.25.

\medskip

-- The result regarding $H_N^s=\mathbb Z_s\wr S_N$ is a reformulation of Theorem 1.14.

\medskip

-- At $s=1$ we have $H_N^1=S_N$, and $P^1=P$, and we recover Theorem 1.12.

\medskip

-- At $s=2$ we have $H_N^2=H_N$, and $P^2=P_{even}$, and we recover Theorem 1.13.

\medskip

-- At $s=\infty$ we have $H_N^\infty=K_N$, and $P^\infty=\mathcal P_{even}$, and we recover Theorem 1.15.
\end{proof}

We can further improve our formalism, with the following definition:

\index{category of partitions}
\index{horizontal concatenation}
\index{vertical concatenation}
\index{upside-down turning}
\index{semicircle partition}
\index{crossing partition}

\begin{definition}
A category of partitions is a collection of sets $D=\bigsqcup_{k,l}D(k,l)$, with $D(k,l)\subset P(k,l)$, having the following properties:
\begin{enumerate}
\item Stability under the horizontal concatenation, $(\pi,\sigma)\to[\pi\sigma]$.

\item Stability under vertical concatenation $(\pi,\sigma)\to[^\sigma_\pi]$, with matching middle symbols.

\item Stability under the upside-down turning $*$, with switching of colors, $\circ\leftrightarrow\bullet$.

\item Each set $P(k,k)$ contains the identity partition $||\ldots||$.

\item The sets $P(\emptyset,\circ\bullet)$ and $P(\emptyset,\bullet\circ)$ both contain the semicircle $\cap$.

\item The sets $P(k,\bar{k})$ with $|k|=2$ contain the crossing partition $\slash\hskip-2.0mm\backslash$.
\end{enumerate}
\end{definition} 

Observe that all the sets of partitions that we used so far in this book, and notably those appearing in Theorem 1.27, are categories of partitions. There are many other examples of such categories, and we will explore this later in this book. Now back to theory, we have the following result, improving our easiness formalism:

\index{category of partitions}
\index{easy group}
\index{Tannakian duality}

\begin{theorem}
Any category of partitions $D\subset P$ produces a series of easy groups $G=(G_N)$, with $G_N\subset U_N$ for any $N\in\mathbb N$, via the formula
$$C_{kl}=span\left(T_\pi\Big|\pi\in D(k,l)\right)$$
for any $k,l$, and Tannakian duality. Any easy group appears in this way.
\end{theorem}

\begin{proof}
This follows indeed from Tannakian duality. To be more precise:

\medskip

(1) In what regards the first assertion, once we fix an integer $N\in\mathbb N$, the various axioms in Definition 1.28 show that the following spaces form a Tannakian category:
$$span\left(T_\pi\Big|\pi\in D(k,l)\right)$$

Thus, Tannakian duality applies, and provides us with a closed subgroup $G_N\subset U_N$ such that the following equalities are satisfied, for any colored integers $k,l$:
$$C_{kl}=span\left(T_\pi\Big|\pi\in D(k,l)\right)$$

But this closed subgroup $G_N\subset U_N$ is easy by definition, and we get the result.

\medskip

(2) Conversely now, consider an easy quantum group $G\subset U_N$, in the sense of Definition 1.26, with this meaning that the corresponding Tannakian category appears as follows, with $D(k,l)\subset P(k,l)$ being certain sets of partitions:
$$C_{kl}=span\left(T_\pi\Big|\pi\in D(k,l)\right)$$

We can then ``saturate'' the collection of sets $D=D(k,l)$, by setting:
$$\widetilde{D}(k,l)=\left\{\pi\in P(k,l)\Big|T_\pi\in C_{kl}\right\}$$

To be more precise, with this definition made, we have inclusions as follows, with the collection $\widetilde{D}$ formed by the sets on the right being a category of partitions:
$$D(k,l)\subset\widetilde{D}(k,l)$$

Also, it follows from definitions that we have equalities as follows:
$$C_{kl}=span\left(T_\pi\Big|\pi\in\widetilde{D}(k,l)\right)$$

Thus, $G$ appears from a category of partitions, namely $\widetilde{D}$, as desired.
\end{proof}

In relation with the above results, which seem to close axiomatic discussions, observe however that the correspondence $D\leftrightarrow G$ coming from the above theorem and its proof is not exactly bijective, for instance because at $N=1$ the group $G=\{1\}$, which is easy, appears from any $D$. This is a subtle issue, and we will be back to this.

\bigskip

As already mentioned, at the practical level the question of fully classifying the categories of partitions, and the easy groups, appears. We will be back to this, later.

\bigskip

Finally, at the level of main examples, let us record the following result:

\begin{theorem}
The following groups, with $H_N=\mathbb Z_2\wr S_N$ and $K_N=\mathbb T\wr S_N$,
$$\xymatrix@R=50pt@C=50pt{
K_N\ar[r]&U_N\\
H_N\ar[u]\ar[r]&O_N\ar[u]}$$
are all easy, the corresponding categories of partitions being as follows,
$$\xymatrix@R=17mm@C=18mm{
\mathcal P_{even}\ar[d]&\mathcal P_2\ar[l]\ar[d]\\
P_{even}&P_2\ar[l]}$$
with $2$ standing for pairings, and ``even'' standing for partitions with even blocks. 
\end{theorem}

\begin{proof}
This is indeed something that we know from Theorem 1.27.
\end{proof}

And good news, that is all. As mentioned in the beginning of this chapter, our aim here was to get a bit familiar with Brauer theorems and easiness, and job done, I hope. As for the deeper understanding of all this, we have chapter 2 below, where we will review all this, from a quantum group perspective, and then the rest of the book too.

\section*{1e. Exercises}

What we did in this chapter was basically Brauer theorems for the closed subgroups $G\subset U_N$, assuming some familiarity with these, and before getting to exercises about this, we can only recommend getting more familiar with the subgroups $G\subset U_N$, via:

\begin{exercise}
Get more familiar with the subgroups $G\subset U_N$, by learning some:
\begin{enumerate}
\item Finite group theory: abelian groups, permutation groups, reflection groups. 

\item Compact group theory: Haar integration and Peter-Weyl theory for them.

\item Compact Lie group theory: Lie theory and Tannakian duality for them.
\end{enumerate}
\end{exercise}

In relation now with what we did in this chapter, there have been many computations in relation with the direct verification of Claim 1.2, and doing some more computations here is the best exercise that we can recommend, for better understanding all this:

\begin{exercise}
Prove directly Claim 1.2 for the following groups:
\begin{enumerate}
\item Subgroups $G\subset \mathbb Z_2^N$: cases $|G|=4$, $|G|=2^{N-2}$.

\item Subgroups $G\subset\mathbb T^N$: cases $|G|$ small, and $|\mathbb T^N/G|$ small.

\item Complex reflection groups $G=\mathbb Z_s\wr S_N$, with full details.

\item The basic continuous groups, namely $O_N$ and $U_N$.
\end{enumerate}
\end{exercise}

At the theoretical level, passed what can be done with bare hands, our main tool was Tannakian duality, so coming as a continuation of Exercise 1.31, we have:

\begin{exercise}
Learn more about Tannakian duality, and various proofs of it, and more about Brauer theorems as well, and various proofs of them too.
\end{exercise}

In relation now with easiness, we have the following exercise:

\begin{exercise}
Prove that the real and complex bistochastic groups, $B_N\subset O_N$ and $C_N\subset U_N$, consisting of matrices having sum $1$ on each row and column, are easy.
\end{exercise}

Finally, getting back again to theory, we have a key exercise, as follows:

\begin{exercise}
You might know from Peter-Weyl theory that for $G\subset U_N$ we have
$$\dim(C_{0k})=\int_G\chi^k$$
where $\chi(g)=Tr(g)$. Can we reformulate what we did in this chapter, in analytic terms?
\end{exercise}

Some of the above exercises might be quite difficult. But do not worry, we will come back to most of them, later in this book.

\chapter{Quantum groups}

\section*{2a. Operator algebras}

Welcome to easiness, take two. What we learned in the previous chapter was in fact just half of the story, and the other half, involving quantum analogues of the compact groups $G\subset U_N$ that we considered, is to be discussed here. Among others, the square formed by the main examples of easy groups will evolve into a cube, as follows:
$$
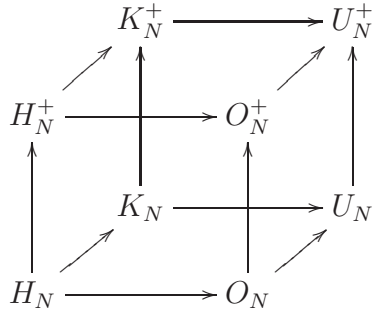
\xymatrix@R=16pt@C=16pt{
&K_N^+\ar[rr]&&U_N^+\\
H_N^+\ar[rr]\ar[ur]&&O_N^+\ar[ur]\\
&K_N\ar[rr]\ar[uu]&&U_N\ar[uu]\\
H_N\ar[uu]\ar[ur]\ar[rr]&&O_N\ar[uu]\ar[ur]
}$$

In order to get started, we need to know what a ``quantum space'' is, and for quantum groups we will see afterwards. However, this is not an easy question, and as the term ``quantum'' tends to indicate, we are a bit into physics here. So, we must either know some quantum physics, or trust mathematical physicists, and their findings. We will opt for this latter way. Following, as usual, John von Neumann, we have:

\index{quantum space}
\index{operator algebra}
\index{von Neumann}

\begin{fact}
A quantum space is the dual of an operator algebra.
\end{fact}

So, our plan will be that of explaining what an operator algebra is, then what a quantum space is, and then what a quantum group is. Afterwards, we will axiomatize the easy quantum groups, and work out some basic examples, including those appearing in the above cube. Getting started now, what we need is the following definition:

\index{Hilbert space}
\index{scalar product}

\begin{definition}
A Hilbert space is a complex vector space $H$ given with a scalar product $<x,y>$, satisfying the following conditions:
\begin{enumerate}
\item $<x,y>$ is linear in $x$, and antilinear in $y$.

\item $\overline{<x,y>}=<y,x>$, for any $x,y$.

\item $<x,x>>0$, for any $x\neq0$.

\item $H$ is complete with respect to the norm $||x||=\sqrt{<x,x>}$.
\end{enumerate}
\end{definition}

Here the fact that $||.||$ is indeed a norm comes from the Cauchy-Schwarz inequality, which states that if the conditions (1,2,3) above are satisfied, then we have:
$$|<x,y>|\leq||x||\cdot||y||$$

Indeed, this inequality comes from the fact that the following degree 2 polynomial, with $t\in\mathbb R$ and $w\in\mathbb T$, being positive, its discriminant must be negative:
$$f(t)=||x+twy||^2$$

At the level of basic examples, we first have the Hilbert space $H=\mathbb C^N$, with its usual scalar product, taken by definition linear at left, namely:
$$<x,y>=\sum_{i=1}^Nx_i\bar{y}_i$$

More generally, given an index set $I$, we can form the Hilbert space $H=l^2(I)$ of the square-summable sequences $(x_i)_{i\in I}$, with similar scalar product, namely:
$$<x,y>=\sum_{i\in I}x_i\bar{y}_i$$

Even more generally, given a measured space $X$, we can form the space $H=L^2(X)$ of square-summable functions $f:X\to\mathbb C$, with similar scalar product, namely:
$$<f,g>=\int_X f(x)\overline{g(x)}\,dx$$

Quite remarkably, this latter extension is, at least at the level of the very abstract theory, not really needed, and this due to the following result:

\index{orthonormal basis}
\index{Gram-Schmidt}

\begin{theorem}
Each Hilbert space $H$ has a basis, meaning a set $\{e_i\}_{i\in I}$ which spans a dense subspace of $H$, whose vectors are pairwise orthogonal, and of norm one:
$$<e_i,e_j>=\delta_{ij}$$
Moreover, the cardinality of the indexing set $|I|=\dim H$ is uniquely determined by $H$, and we have an isomorphism $H\simeq l^2(I)$.  
\end{theorem}

\begin{proof}
Here the first assertion follows from Gram-Schmidt, the idea being that any algebraic basis $\{f_i\}_{i\in I}$ can be turned into an orthonormal basis $\{e_i\}_{i\in I}$, by using a recurrence method. As for the second assertion, this is clear from the first one.
\end{proof}

The above is something quite tricky. For instance with $H=L^2(\mathbb R)$, or even better, with $H=L^2(\mathbb R^3)$, which is the space of wave functions of the electron, by the Weierstrass theorem $I$ can be taken countable, so we obtain an isomorphism as follows:
$$L^2(\mathbb R^3)\simeq l^2(\mathbb N)$$

However, this isomorphism is not very explicit, and struggling with it is a main occupation when doing basic quantum mechanics, such as solving the hydrogen atom.

\bigskip

Speaking physics, observe that in Definition 2.2 our scalar products are linear at left. This is the so-called mathematicians' convention, as opposed to Dirac's convention, but the point is that in recent times many fine theoretical physicists, genuinely interested in physics, and not knowing much mathematics, including of course myself and my collaborators, have opted for this mathematicians' convention, for various technical reasons.

\bigskip

Getting now to operators and operator algebras, we first have:

\index{operator}
\index{bounded operator}
\index{operator norm}

\begin{proposition}
For a linear operator $T:H\to H$, the following are equivalent:
\begin{enumerate}
\item $T$ is continuous.

\item $T$ is continuous at $0$.

\item $T$ maps the unit ball of $H$ into something bounded.

\item $T$ is bounded, in the sense that $||T||=\sup_{||x||=1}||Tx||$ is finite.
\end{enumerate}
\end{proposition}

\begin{proof}
Here the equivalences $(1)\iff(2)\iff(3)\iff(4)$ all follow from definitions, by using the linearity of $T$, and performing various rescalings, and with the number $||T||$ needed in (4) being the bound coming from (3).
\end{proof}

With the above result in hand, we can now formulate:

\index{operator algebra}
\index{infinite matrix}

\begin{theorem}
The bounded operators $T:H\to H$ form an algebra $B(H)$, on which 
$$||T||=\sup_{||x||=1}||Tx||$$
is a norm, and which is complete with respect to this norm. In the case where the Hilbert space $H$ comes with a basis $\{e_i\}_{i\in I}$, we have an embedding
$$B(H)\subset M_I(\mathbb C)$$
which is $M_N(\mathbb C)\subset M_N(\mathbb C)$ for $H=\mathbb C^N$, but which is not an isomorphism in general.
\end{theorem}

\begin{proof}
Again, all this is standard, with the algebra property of $B(H)$ being clear, with the norm property of $||.||$ being clear too, and with the norm closedness of $B(H)$ coming by constructing the limit of a Cauchy sequence $\{T_n\}$ as follows:
$$Tx=\lim_{n\to\infty}T_nx\quad,\quad\forall x\in H$$

Finally, in what regards the embedding $B(H)\subset M_I(\mathbb C)$, this can be constructed by using the same formula as in usual linear algebra, namely:
$$T_{ij}=<Te_j,e_i>$$

As for the fact that this embedding is not an isomorphism, when $\dim H=\infty$, the point here is that with $I=\mathbb N$ the infinite matrix $T=diag(0,1,2,3,\ldots)$ does not come from a bounded operator, providing us with the desired counterexample.
\end{proof}

Summarizing, the correct infinite analogue of the algebra $M_N(\mathbb C)$ is not the infinite matrix algebra $M_I(\mathbb C)$, which is actually not even an algebra, when $|I|=\infty$, but rather the algebra $B(H)$ of bounded linear operators $T:H\to H$ on a Hilbert space $H$.

\bigskip

With this understood, we can go now into truly intreresting material. Everything advanced that you know about $M_N(\mathbb C)$, be that projections, rotations, other special matrices, or spectral theorem and so on, uses adjoint matrices. So, let us talk first about adjoint operators, in our framework. The result here is as follows:

\index{adjoint operator}

\begin{proposition}
Any bounded operator $T\in B(H)$ has an adjoint $T^*\in B(H)$, given by the following formula, valid for any two vectors $x,y\in H$:
$$<Tx,y>=<x,T^*y>$$
The operation $T\to T^*$ is then an isometric involution of $B(H)$, and we have: 
$$||TT^*||=||T||^2$$
When $H$ comes with an orthonormal basis $\{e_i\}_{i\in I}$, we have $(T^*)_{ij}=\overline{T}_{ji}$.
\end{proposition}

\begin{proof}
As before, all this is standard material. Given an operator $T\in B(H)$, let us pick a vector $y\in H$, and consider the following linear form: 
$$x\to<Tx,y>$$

This linear form must then come from a scalar product with a vector $T^*y$, as in the statement, and we obtain in this way a definition for $T^*$, namely $y\to T^*y$. It is then routine to check that we have indeed $T^*\in B(H)$, with this coming from:
$$||T^*||=||T||$$

The fact that $T\to T^*$ is then an involution of $B(H)$ is routine too. Regarding now the formula $||TT^*||=||T||^2$, in one sense we have the following estimate:
$$||TT^*||\leq||T||\cdot||T^*||=||T||^2$$

In the other sense, we have the following estimate:
\begin{eqnarray*}
||T||^2
&=&\sup_{||x||=1}|<Tx,Tx>|\\
&=&\sup_{||x||=1}|<x,T^*Tx>|\\
&\leq&||T^*T||
\end{eqnarray*}

Now by replacing in this formula $T\to T^*$ we obtain $||T||^2\leq||TT^*||$, as desired. Finally, $(T^*)_{ij}=\overline{T}_{ji}$ is clear from the formula $T_{ij}=<Te_j,e_i>$, applied to $T,T^*$.
\end{proof}

Good news, we can now talk about operator algebras, as follows:

\index{operator algebra}

\begin{definition}
An operator algebra is an algebra of bounded operators $A\subset B(H)$ which contains the unit, is closed under taking adjoints, 
$$T\in A\implies T^*\in A$$
and is closed as well under the norm.
\end{definition}

This definition is in fact one of the many possible ones, with the choice here being a matter of knowledge of mathematics, and physics, and taste. But more on this later. Getting now to where we wanted to get, with this, we can formulate some tough results, inspired by the usual linear algebra, that of the algebra $M_N(\mathbb C)$, as follows:

\index{self-adjoint operator}
\index{normal operator}
\index{spectral theorem}
\index{commutative algebra}
\index{diagonalization}

\begin{theorem}
The following happen:
\begin{enumerate}
\item Any self-adjoint operator, $T=T^*$, is diagonalizable.

\item More generally, any normal operator, $TT^*=T^*T$, is diagonalizable.

\item In fact, any family $\{T_i\}$ of commuting normal operators is diagonalizable.
\end{enumerate}
Thus, any commutative operator algebra is of the form $A=C(X)$, with $X$ compact space.
\end{theorem}

\begin{proof}
This is certainly a tough theorem, with (1,2,3) coming by generalizing the Spectral Theorem, in its various incarnations, for the usual matrices $M\in M_N(\mathbb C)$. As for the final conclusion, this follows from (3), because if we write $A=span(T_i)$, then the family $\{T_i\}$ consists of commuting normal operators, and this leads to the conclusion $A=C(X)$, with $X$ being a certain compact space associated to the family $\{T_i\}$.
\end{proof}

In relation with the above result, there are some good news and some bad news. The good news first, we can, eventually, talk about quantum spaces, as follows:

\index{quantum space}
\index{compact quantum space}

\begin{definition}
We can think of any operator algebra $A\subset B(H)$ as being of the form
$$A=C(X)$$
with $X$ compact quantum space. When $A$ is commutative, $X$ is a usual compact space.
\end{definition}

As for the bad news, all this is based on Theorem 2.8, which remains something terribly complicated, and that we would rather like to avoid, when building foundations. Also, there is a problem with functoriality, because a morphism a quantum spaces $X\to Y$ should normally come from a morphism of algebras $C(Y)\to C(X)$, but shall we ask or not something in relation with the embeddings $C(X)\subset B(H)$ and $C(Y)\subset B(K)$. And finally, we have a philosophical problem too, the Hilbert spaces are certainly nice objects, but do we really need them for talking about basic things like quantum spaces. 

\bigskip

But above everything, we have the following dumb argument:

\index{Hilbert space}

\begin{criticism}
Needing a Hilbert space for talking about the circle is ridiculous.
\end{criticism}

To be more precise, the circle $\mathbb T$ is the simplest compact space that we know, and this since childhood. However, in order to view it as a compact quantum space, as in Definition 2.9, we need something of type $C(\mathbb T)\subset B(L^2(\mathbb T))$, which is ridiculous.

\bigskip

Summarizing, Definition 2.7, Theorem 2.8 and Definition 2.9 are not good, and we must invent something else. And here is the magic trick, due to Gelfand:

\index{operator algebra}
\index{abstract algebra}
\index{normed algebra}
\index{Banach algebra}

\begin{definition}
An abstract operator algebra, or $C^*$-algebra, is a complex algebra $A$ having a norm $||.||$ and an involution $*$, subject to the following conditions:
\begin{enumerate}
\item $A$ is closed with respect to the norm.

\item We have $||aa^*||=||a||^2$, for any $a\in A$.
\end{enumerate}
\end{definition}

In other words, what we did here is to axiomatize the abstract properties of the operator algebras $A\subset B(H)$, without any reference to the Hilbert space $H$. We will see in a moment that our axiomatization is indeed complete, in the sense that any $C^*$-algebra appears as an operator algebra, $A\subset B(H)$. Thus, getting back now to our quantum space questions, we will be able to recycle Defintion 2.9, simply by replacing there ``operator algebra'' by $C^*$-algebra, and everything, or almost, will be fine. In particular, no one in this world will ever be able to come with something like Criticism 2.10.

\bigskip

Getting to work now, let us develop the theory of $C^*$-algebras. We first have:

\index{spectrum}
\index{spectral radius}

\begin{proposition}
Given an element $a\in A$ of a $C^*$-algebra, define its spectrum as:
$$\sigma(a)=\left\{\lambda\in\mathbb C\Big|a-\lambda\notin A^{-1}\right\}$$
The following spectral theory results hold, exactly as in the $A=B(H)$ case:
\begin{enumerate}
\item We have $\sigma(ab)\cup\{0\}=\sigma(ba)\cup\{0\}$.

\item We have $\sigma(f(a))=f(\sigma(a))$, for any $f\in\mathbb C(X)$ having poles outside $\sigma(a)$.

\item The spectrum $\sigma(a)$ is compact, non-empty, and contained in $D_0(||a||)$.

\item The spectra of unitaries $(u^*=u^{-1})$ and self-adjoints $(a=a^*)$ are on $\mathbb T,\mathbb R$.

\item The spectral radius of normal elements $(aa^*=a^*a)$ is given by $\rho(a)=||a||$.
\end{enumerate}
In addition, assuming $a\in A\subset B$, the spectra of $a$ with respect to $A$ and to $B$ coincide.
\end{proposition}

\begin{proof}
Here the assertions (1-5), which are formulated a bit informally, are well-known for the full operator algebra $A=B(H)$, and the proof in general is similar:

\medskip

(1) Assuming that $1-ab$ is invertible, with inverse $c$, we have $abc=cab=c-1$, and it follows that $1-ba$ is invertible too, with inverse $1+bca$. Thus $\sigma(ab),\sigma(ba)$ agree on $1\in\mathbb C$, and by linearity, it follows that $\sigma(ab),\sigma(ba)$ agree on any point $\lambda\in\mathbb C^*$.

\medskip

(2) The formula $\sigma(f(a))=f(\sigma(a))$ is clear for polynomials, $f\in\mathbb C[X]$, by factorizing $f-\lambda$, with $\lambda\in\mathbb C$. Then, the extension to the rational functions is straightforward, because $P(a)/Q(a)-\lambda$ is invertible precisely when $P(a)-\lambda Q(a)$ is.

\medskip

(3) By using $1/(1-b)=1+b+b^2+\ldots$ for $||b||<1$ we obtain that $a-\lambda$ is invertible for $|\lambda|>||a||$, and so $\sigma(a)\subset D_0(||a||)$. It is also clear that $\sigma(a)$ is closed, so what we have is a compact set. Finally, assuming $\sigma(a)=\emptyset$ the function $f(\lambda)=\varphi((a-\lambda)^{-1})$ is well-defined, for any $\varphi\in A^*$, and by Liouville we get $f=0$, contradiction.

\medskip

(4) Assuming $u^*=u^{-1}$ we have $||u||=1$, and so $\sigma(u)\subset D_0(1)$. But with $f(z)=z^{-1}$ we obtain via (2) that we have as well $\sigma(u)\subset f(D_0(1))$, and this gives $\sigma(u)\subset\mathbb T$. As for the result regarding the self-adjoints, this can be obtained from the result for the unitaries, by using (2) with functions of type $f(z)=(z+it)/(z-it)$, with $t\in\mathbb R$.

\medskip

(5) It is routine to check, by integrating quantities of type $z^n/(z-a)$ over circles centered at the origin, and estimating, that the spectral radius is given by $\rho(a)=\lim||a^n||^{1/n}$. But in the self-adjoint case, $a=a^*$, this gives $\rho(a)=||a||$, by using exponents of type $n=2^k$, and then the extension to the general normal case is straightforward.

\medskip 

(6) Regarding now the last assertion, the inclusion $\sigma_B(a)\subset\sigma_A(a)$ is clear. For the converse, assume $a-\lambda\in B^{-1}$, and set $b=(a-\lambda )^*(a-\lambda )$. We have then:
$$\sigma_A(b)-\sigma_B(b)=\left\{\mu\in\mathbb C-\sigma_B(b)\Big|(b-\mu)^{-1}\in B-A\right\}$$

Thus this difference in an open subset of $\mathbb C$. On the other hand $b$ being self-adjoint, its two spectra are both real, and so is their difference. Thus the two spectra of $b$ are equal, and in particular $b$ is invertible in $A$, and so $a-\lambda\in A^{-1}$, as desired.
\end{proof}

With these ingredients, we can now a prove a key result, as follows:

\index{Gelfand theorem}
\index{commutative algebra}
\index{spectrum of algebra}
\index{algebra character}
\index{compact space}

\begin{theorem}[Gelfand]
If $X$ is a compact space,  the algebra $C(X)$ of continuous functions on it $f:X\to\mathbb C$ is a $C^*$-algebra, with usual norm and involution, namely:
$$||f||=\sup_{x\in X}|f(x)|\quad,\quad 
f^*(x)=\overline{f(x)}$$
Conversely, any commutative $C^*$-algebra is of this form, $A=C(X)$, with 
$$X=\Big\{\chi:A\to\mathbb C\ ,\ {\rm normed\ algebra\ character}\Big\}$$
with topology making continuous the evaluation maps $ev_a:\chi\to\chi(a)$.
\end{theorem}

\begin{proof}
There are several things going on here, the idea being as follows:

\medskip

(1) The first assertion is clear from definitions. Observe that we have indeed:
$$||ff^*||
=\sup_{x\in X}|f(x)|^2
=||f||^2$$

Observe also that the algebra $C(X)$ is commutative, because $fg=gf$.

\medskip

(2) Conversely, given a commutative $C^*$-algebra $A$, let us define $X$ as in the statement. Then $X$ is compact, and $a\to ev_a$ is a morphism of algebras, as follows:
$$ev:A\to C(X)$$

(3) We first prove that $ev$ is involutive. We use the following formula, which is similar to the $z=Re(z)+iIm(z)$ decomposition formula for usual complex numbers:
$$a=\frac{a+a^*}{2}+i\cdot\frac{a-a^*}{2i}$$

Thus it is enough to prove $ev_{a^*}=ev_a^*$ for the self-adjoint elements $a$. But this is the same as proving that $a=a^*$ implies that $ev_a$ is a real function, which is in turn true, by Proposition 2.12, because $ev_a(\chi)=\chi(a)$ is an element of $\sigma(a)$, contained in $\mathbb R$.

\medskip

(4) Since $A$ is commutative, each element is normal, so $ev$ is isometric:
$$||ev_a||
=\rho(a)
=||a||$$

It remains to prove that $ev$ is surjective. But this follows from the Stone-Weierstrass theorem, because $ev(A)$ is a closed subalgebra of $C(X)$, which separates the points.
\end{proof}

The above result is something truly remarkable, and we can now formulate:

\index{quantum space}
\index{compact quantum space}

\begin{definition}
Given an arbitrary $C^*$-algebra $A$, we write it as
$$A=C(X)$$
with $X$ compact quantum space. When $A$ is commutative, $X$ is a usual compact space.
\end{definition}

Observe the similarity with Definition 2.9, which is now to be forgotten. Indeed, our theory based on $C^*$-algebras is much better, not using Hilbert spaces, and free as a bird, and all issues mentioned after Definition 2.9, including Criticism 2.10, now dissapear.

\bigskip

Of course, what we have is still just a beginning of something, and we will soon see, once we will be more advanced, that there are in fact some bugs with Definition 2.14 too. To be more precise, there are certain natural quantum spaces $X$, such as the duals $X=\widehat{\Gamma}$ of the non-amenable groups $\Gamma$, corresponding to several $C^*$-algebras $A$. Thus, the correspondence $A\to X$ from Definition 2.14 is not bijective, and needs a fix.

\bigskip

But more on this later, for the moment let us enjoy what we have. A quick comparison between Theorem 2.8 and Theorem 2.13 suggests that operator algebra and $C^*$-algebra might be actually the same thing. And this is indeed the case, the result being:

\index{representation theorem}
\index{GNS theorem}
\index{operator algebra}

\begin{theorem}
Any $C^*$-algebra appears as an operator algebra:
$$A\subset B(H)$$
Moreover, when $A$ is separable, which is usually the case, $H$ can be taken separable.
\end{theorem}

\begin{proof}
This result, called GNS representation theorem after Gelfand-Naimark-Segal, comes as a continuation of the Gelfand theorem, the idea being as follows:

\medskip

(1) Let us first prove that the result holds in the commutative case, $A=C(X)$. Here, we can pick a positive measure on $X$, and construct our embedding as follows:
$$C(X)\subset B(L^2(X))\quad,\quad f\to[g\to fg]$$

(2) In general the proof is similar, the idea being that given a $C^*$-algebra $A$ we can construct a Hilbert space $H=L^2(A)$, and then an embedding as above:
$$A\subset B(L^2(A))\quad,\quad a\to[b\to ab]$$

(3) Finally, the last assertion is clear, because when $A$ is separable, meaning that it has a countable algebraic basis, so does the associated Hilbert space $H=L^2(A)$.
\end{proof}

All this is very nice, and getting back to our original motivations, we have now a beautiful notion of compact quantum space, coming from Definition 2.11, Theorem 2.13 and Definition 2.14. Also, as a bonus, we have as well some spectral theory tools for the study of such spaces, coming from Proposition 2.12, and even a theorem allowing us to pull out of a hat a Hilbert space, in case we ever get lost, namely Theorem 2.15.

\section*{2b. Quantum groups}

We can now go ahead and develop our quantum group, and general easiness program. As a starting point, we have the following key definition, due to Woronowicz \cite{wo1}:

\index{Woronowicz algebra}
\index{comultiplication}
\index{counit}
\index{antipode}
\index{tensor product}
\index{opposite algebra}

\begin{definition}
A Woronowicz algebra is a $C^*$-algebra $A$, given with a unitary matrix $u\in M_N(A)$ whose coefficients generate $A$, such that the formulae
$$\Delta(u_{ij})=\sum_ku_{ik}\otimes u_{kj}\quad,\quad
\varepsilon(u_{ij})=\delta_{ij}\quad,\quad 
S(u_{ij})=u_{ji}^*$$
define morphisms of $C^*$-algebras $\Delta:A\to A\otimes A$, $\varepsilon:A\to\mathbb C$, $S:A\to A^{opp}$, called comultiplication, counit and antipode.
\end{definition}

In this definition $\otimes$ can be any $C^*$-algebraic completion of the usual algebraic tensor product $\otimes_{alg}$, and the symbol $A^{opp}$ denotes the opposite algebra. More on this later.

\bigskip

We say that $A$ is cocommutative when $\Sigma\Delta=\Delta$, where $\Sigma(a\otimes b)=b\otimes a$ is the flip. We have the following result, which justifies the terminology and axioms:

\index{cocommutative algebra}
\index{flip map}
\index{compact Lie group}
\index{commutative algebra}
\index{discrete group}

\begin{proposition}
The following are Woronowicz algebras:
\begin{enumerate}
\item $C(G)$, with $G\subset U_N$ compact Lie group. Here the structural maps are:
$$\Delta(\varphi)=[(g,h)\to \varphi(gh)]\quad,\quad
\varepsilon(\varphi)=\varphi(1)\quad,\quad
S(\varphi)=[g\to\varphi(g^{-1})]$$

\item $C^*(\Gamma)$, with $F_N\to\Gamma$ finitely generated group. Here the structural maps are:
$$\Delta(g)=g\otimes g\quad,\quad
\varepsilon(g)=1\quad,\quad 
S(g)=g^{-1}$$
\end{enumerate}
Moreover, we obtain in this way all the commutative/cocommutative algebras.
\end{proposition}

\begin{proof}
In both cases, we have to indicate a certain matrix $u$. For the first assertion, we can use the matrix $u=(u_{ij})$ formed by matrix coordinates of $G$, given by:
$$g=\begin{pmatrix}
u_{11}(g)&\ldots&u_{1N}(g)\\
\vdots&&\vdots\\
u_{N1}(g)&\ldots&u_{NN}(g)
\end{pmatrix}$$

As for the second assertion, we can use here the diagonal matrix formed by generators:
$$u=\begin{pmatrix}
g_1&&0\\
&\ddots&\\
0&&g_N
\end{pmatrix}$$

Finally, the last assertion follows from the Gelfand theorem, in the commutative case. In the cocommutative case this follows from the Peter-Weyl theory, explained below.
\end{proof}

In view of Proposition 2.17, we can formulate the following definition:

\index{compact quantum group}
\index{discrete quantum group}
\index{abelian group}
\index{Pontrjagin duality}

\begin{definition}
Given a Woronowicz algebra $A$, we formally write
$$A=C(G)=C^*(\Gamma)$$
and call $G$ compact quantum group, and $\Gamma$ discrete quantum group.
\end{definition}

When $A$ is both commutative and cocommutative, $G$ is a compact abelian group, $\Gamma$ is a discrete abelian group, and these groups are dual to each other:
$$G=\widehat{\Gamma}\quad,\quad\Gamma=\widehat{G}$$

In general, we still agree to write the formulae $G=\widehat{\Gamma},\Gamma=\widehat{G}$, but in a formal sense. Finally, let us make as well the following key convention:

\index{equivalence of algebras}

\begin{definition}
We identify two Woronowicz algebras $(A,u)$ and $(B,v)$, as well as the corresponding quantum groups, when we have an isomorphism of $*$-algebras 
$$<u_{ij}>\simeq<v_{ij}>$$
mapping standard coordinates to standard coordinates.
\end{definition}

This convention is here for avoiding amenability issues, as for any quantum group to correspond to a unique Woronowicz algebra, and more on this later. Moving ahead now, we need tools, for the study of our quantum groups. In the classical case, the main tool for the study of the groups $G$ are the group axioms, namely:
$$m(m\times id)=m(id\times m)$$
$$m(u\times id)=m(id\times u)=id$$
$$m(i\times id)\delta=m(id\times i)\delta=1$$

Here $\delta(g)=(g,g)$. The point now is that all these formulae hold as well for our quantum groups, in algebra formulation of course, the result being as follows:

\index{group axioms}
\index{Hopf algebra}
\index{square of antipode}

\begin{proposition}
The maps $\Delta,\varepsilon,S$ satisfy the Hopf algebra axioms, namely:
$$(\Delta\otimes id)\Delta=(id\otimes \Delta)\Delta$$
$$(\varepsilon\otimes id)\Delta=(id\otimes\varepsilon)\Delta=id$$
$$m(S\otimes id)\Delta=m(id\otimes S)\Delta=\varepsilon(.)1$$
In addition, the square of the antipode is the identity, $S^2=id$.
\end{proposition}

\begin{proof}
As a first observation, the result holds in the commutative case, $A=C(G)$ with $G\subset U_N$. Indeed, here we know from Proposition 2.17 that we have:
$$\Delta=m^t\quad,\quad 
\varepsilon=u^t\quad,\quad 
S=i^t$$

Thus, in this case, the various conditions in the statement on $\Delta,\varepsilon,S$ simply come by transposition from the group axioms satisfied by $m,u,i$. In general now, we have:
$$(\Delta\otimes id)\Delta(u_{ij})=(id\otimes \Delta)\Delta(u_{ij})=\sum_{kl}u_{ik}\otimes u_{kl}\otimes u_{lj}$$

As for the other axioms, their verification is similar, with the technical remark that the first two formulae hold on $A$, while the third formula only holds on $<u_{ij}>$.
\end{proof}

All this is very nice. In order to reach now to more advanced results, following again Woronowicz \cite{wo1}, let us call corepresentation of $A$ any unitary matrix $v\in M_n(\mathcal A)$, where $\mathcal A=<u_{ij}>$, satisfying the same conditions as those satisfied by $u$, namely:
$$\Delta(v_{ij})=\sum_kv_{ik}\otimes v_{kj}\quad,\quad 
\varepsilon(v_{ij})=\delta_{ij}\quad,\quad 
S(v_{ij})=v_{ji}^*$$

We have the following key result, due to Woronowicz \cite{wo1}:

\index{Haar integration}
\index{Ces\`aro limit}
\index{convolution}
\index{corepresentation}

\begin{theorem}
Any Woronowicz algebra has a unique Haar integration functional, 
$$\left(\int_G\otimes id\right)\Delta=\left(id\otimes\int_G\right)\Delta=\int_G(.)1$$
which can be constructed by starting with any faithful positive form $\varphi\in A^*$, and setting
$$\int_G=\lim_{n\to\infty}\frac{1}{n}\sum_{k=1}^n\varphi^{*k}$$
where $\phi*\psi=(\phi\otimes\psi)\Delta$. Moreover, for any corepresentation $v\in M_n(\mathbb C)\otimes A$ we have
$$\left(id\otimes\int_G\right)v=P$$
where $P$ is the orthogonal projection onto $Fix(v)=\{\xi\in\mathbb C^n|v\xi=\xi\}$.
\end{theorem}

\begin{proof}
Following \cite{wo1}, this can be done in 3 steps, as follows:

\medskip

(1) Given $\varphi\in A^*$, our claim is that the following limit converges, for any $a\in A$:
$$\int_\varphi a=\lim_{n\to\infty}\frac{1}{n}\sum_{k=1}^n\varphi^{*k}(a)$$

Indeed, by linearity we can assume that $a\in A$ is the coefficient of certain corepresentation, $a=(\tau\otimes id)v$. But in this case, an elementary computation gives the following formula, with $P_\varphi$ being the orthogonal projection onto the $1$-eigenspace of $(id\otimes\varphi)v$:
$$\left(id\otimes\int_\varphi\right)v=P_\varphi$$

(2) Since $v\xi=\xi$ implies $[(id\otimes\varphi)v]\xi=\xi$, we have $P_\varphi\geq P$, where $P$ is the orthogonal projection onto the fixed point space in the statement, namely:
$$Fix(v)=\left\{\xi\in\mathbb C^n\Big|v\xi=\xi\right\}$$

The point now is that when $\varphi\in A^*$ is faithful, by using a standard positivity trick, we can prove that we have $P_\varphi=P$, exactly as in the classical case.

\medskip

(3) With the above formula in hand, the left and right invariance of $\int_G=\int_\varphi$ is clear on coefficients, and so in general, and this gives all the assertions. See \cite{wo1}. 
\end{proof}

We can now develop, again following \cite{wo1}, the Peter-Weyl theory for the corepresentations of $A$. Consider the dense subalgebra $\mathcal A\subset A$ generated by the coefficients of the fundamental corepresentation $u$, and endow it with the following scalar product: 
$$<a,b>=\int_Gab^*$$

With this convention, we have the following result, from \cite{wo1}:

\index{Peter-Weyl theory}
\index{corepresentation}
\index{irreducible corepresentation}
\index{character}

\begin{theorem}
We have the following Peter-Weyl type results, with the various operations on corepresentations being defined in a straightforward way:
\begin{enumerate}
\item Any corepresentation decomposes as a sum of irreducible corepresentations.

\item Each irreducible corepresentation appears inside a certain $u^{\otimes k}$.

\item $\mathcal A=\bigoplus_{v\in Irr(A)}M_{\dim(v)}(\mathbb C)$, the summands being pairwise orthogonal.

\item The characters of irreducible corepresentations form an orthonormal system.
\end{enumerate}
\end{theorem}

\begin{proof}
All these results are from \cite{wo1}, the idea being as follows:

\medskip

(1) For a corepresentation $v\in M_n(A)$, the algebra $End(v)=\{T\in M_n(\mathbb C)|Tv=vT\}$ is a finite dimensional $C^*$-algebra, and so decomposes as $End(v)=M_{n_1}(\mathbb C)\oplus\ldots\oplus M_{n_k}(\mathbb C)$. But this gives a decomposition of type $v=v_1+\ldots+v_k$, as desired.

\medskip

(2) Consider the Peter-Weyl corepresentations, $u^{\otimes k}$ with $k$ being a colored integer, defined by $u^{\otimes\emptyset}=1$, $u^{\otimes\circ}=u$, $u^{\otimes\bullet}=\bar{u}$ and multiplicativity. The coefficients of these corepresentations span the dense algebra $\mathcal A$, and by using (1), this gives the result.

\medskip

(3) Here the direct sum decomposition, which is a $*$-coalgebra isomorphism, follows from (2). As for the second assertion, this follows from the fact that $(id\otimes\int_G)v$ is the orthogonal projection $P_v$ onto the space $Fix(v)$, for any corepresentation $v$.

\medskip

(4) Let us define indeed the character of a corepresentation $v\in M_n(A)$ to be the trace, $\chi_v=Tr(v)$. Since this character is a coefficient of $v$, the orthogonality assertion follows from (3). As for the norm 1 claim, this follows once again from $(id\otimes\int_G)v=P_v$. 
\end{proof}

As a first consequence of the above result, we can now clarify the structure of the cocommutative Woronowicz algebras, closing a discussion started in Proposition 2.17:

\index{cocommutative algebra}
\index{discrete group}

\begin{proposition}
For a Woronowicz algebra $A$, the following are equivalent:
\begin{enumerate}
\item $A$ is cocommutative, $\Sigma\Delta=\Delta$.

\item The irreducible corepresentations of $A$ are all $1$-dimensional.

\item $A=C^*(\Gamma)$, for some group $\Gamma=<g_1,\ldots,g_N>$, up to equivalence.
\end{enumerate}
\end{proposition}

\begin{proof}
This follows from the Peter-Weyl theory, as follows:

\medskip

$(1)\implies(2)$ The assumption $\Sigma\Delta=\Delta$ tells us that the inclusion $\mathcal A_{central}=<\chi_v>\subset\mathcal A$ is an isomorphism, and by using Peter-Weyl theory we conclude that any irreducible corepresentation of $A$ must be equal to its character, and so must be 1-dimensional.

\medskip

$(2)\implies(3)$ This follows once again from Peter-Weyl, because if we denote by $\Gamma$ the group formed by the 1-dimensional corepresentations, then we have $\mathcal A=\mathbb C[\Gamma]$, and so $A=C^*(\Gamma)$ up to the standard equivalence relation for Woronowicz algebras.

\medskip

$(3)\implies(1)$ This is something trivial, that we know from Proposition 2.17.
\end{proof}

Still in relation with the discrete groups, but at a more advanced level, following as before Woronowicz \cite{wo1}, we have the following result:

\index{full algebra}
\index{reduced algebra}
\index{amenability}
\index{coamenability}
\index{Kesten criterion}

\begin{theorem}
Let $A_{full}$ be the enveloping $C^*$-algebra of $\mathcal A$, and $A_{red}$ be the quotient of $A$ by the null ideal of the Haar integration. The following are then equivalent:
\begin{enumerate}
\item The Haar functional of $A_{full}$ is faithful.

\item The projection map $A_{full}\to A_{red}$ is an isomorphism.

\item The counit map $\varepsilon:A_{full}\to\mathbb C$ factorizes through $A_{red}$.

\item We have $N\in\sigma(Re(\chi_u))$, the spectrum being taken inside $A_{red}$.
\end{enumerate}
If this is the case, we say that the underlying discrete quantum group $\Gamma$ is amenable.
\end{theorem}

\begin{proof}
This is well-known in the group dual case, $A=C^*(\Gamma)$, with $\Gamma$ being a usual discrete group. In general, the result follows by adapting the group dual case proof:

\medskip

$(1)\iff(2)$ This simply follows from the fact that the GNS construction for the algebra $A_{full}$ with respect to the Haar functional produces the algebra $A_{red}$.

\medskip

$(2)\iff(3)$ Here $\implies$ is trivial, and conversely, a counit $\varepsilon:A_{red}\to\mathbb C$ produces an isomorphism $\Phi:A_{red}\to A_{full}$, by slicing the map $\widetilde{\Delta}:A_{red}\to A_{red}\otimes A_{full}$.

\medskip

$(3)\iff(4)$ Here $\implies$ is clear, coming from $\varepsilon(N-Re(\chi (u)))=0$, and the converse can be proved by doing some functional analysis. See \cite{wo1}.
\end{proof}

\section*{2c. Diagrams, easiness}

In order to reach now to a theory of easiness for the compact quantum groups, we need several extra ingredients. First, we need free analogues of the orthogonal and unitary groups. The constructions here, due to Wang \cite{wa1}, are as follows:

\index{free orthogonal group}
\index{free unitary group}
\index{free rotation group}
\index{free quantum group}

\begin{theorem}
The following universal algebras are Woronowicz algebras,
$$C(O_N^+)=C^*\left((u_{ij})_{i,j=1,\ldots,N}\Big|u=\bar{u},u^t=u^{-1}\right)$$
$$C(U_N^+)=C^*\left((u_{ij})_{i,j=1,\ldots,N}\Big|u^*=u^{-1},u^t=\bar{u}^{-1}\right)$$
so the underlying quantum spaces $O_N^+,U_N^+$ are compact quantum groups.
\end{theorem}

\begin{proof}
This follows from the elementary fact that if a matrix $u=(u_{ij})$ is orthogonal or biunitary, then so must be the following matrices:
$$(u^\Delta)_{ij}=\sum_ku_{ik}\otimes u_{kj}\quad,\quad 
(u^\varepsilon)_{ij}=\delta_{ij}\quad,\quad
(u^S)_{ij}=u_{ji}^*$$

Thus, we can indeed define morphisms $\Delta,\varepsilon,S$ as in Definition 2.16, by using the universal properties of $C(O_N^+)$, $C(U_N^+)$, and this gives the result.
\end{proof}

Next, we need to talk about Tannakian duality.  The result here, which is very similar to the Tannakian duality result from chapter 1, is as follows:

\index{Tannakian duality}
\index{tensor category}
\index{Tannakian category}

\begin{theorem}
The following operations are inverse to each other:
\begin{enumerate}
\item The construction $G\to C$, which associates to a closed subgroup $G\subset_uU_N^+$ the tensor category formed by the intertwiner spaces $C_{kl}=Hom(u^{\otimes k},u^{\otimes l})$.

\item The construction $C\to G$, associating to a tensor category $C$ the closed subgroup $G\subset_uU_N^+$ coming from the relations $T\in Hom(u^{\otimes k},u^{\otimes l})$, with $T\in C_{kl}$.
\end{enumerate}
\end{theorem}

\begin{proof}
The idea is that we have indeed a construction $G\to C_G$, producing a subcategory of the tensor $C^*$-category of finite dimensional Hilbert spaces, as follows:
$$(C_G)_{kl}=Hom(u^{\otimes k},u^{\otimes l})$$

We have as well a construction $C\to G_C$, obtained by setting:
$$C(G_C)=C(U_N^+)\big/\left<T\in Hom(u^{\otimes k},u^{\otimes l})\Big|\forall k,l,\forall T\in C_{kl}\right>$$

Regarding now the bijection claim, some elementary algebra shows that $C=C_{G_C}$ implies $G=G_{C_G}$, and that $C\subset C_{G_C}$ is automatic. Thus we are left with proving:
$$C_{G_C}\subset C$$

But this latter inclusion can be proved indeed, by doing some algebra, and using von Neumann's bicommutant theorem, in finite dimensions. See Malacarne \cite{mal}. 
\end{proof}

Following the material from chapter 1, we can now talk about easiness. Let us first recall from chapter 1 that the partitions produce linear maps, as follows:

\index{partition}
\index{Kronecker symbol}

\begin{definition}
Associated to any partition $\pi\in P(k,l)$ between an upper row of $k$ points and a lower row of $l$ points is the linear map $T_\pi:(\mathbb C^N)^{\otimes k}\to(\mathbb C^N)^{\otimes l}$ given by 
$$T_\pi(e_{i_1}\otimes\ldots\otimes e_{i_k})=\sum_{j_1\ldots j_l}\delta_\pi\begin{pmatrix}i_1&\ldots&i_k\\ j_1&\ldots&j_l\end{pmatrix}e_{j_1}\otimes\ldots\otimes e_{j_l}$$
with the Kronecker type symbols $\delta_\pi\in\{0,1\}$ depending on whether the indices fit or not. 
\end{definition}

To be more precise, we agree to put the two multi-indices on the two rows of points, in the obvious way. The Kronecker symbols are then defined by $\delta_\pi=1$ when all the strings of $\pi$ join equal indices, and by $\delta_\pi=0$ otherwise. This construction is motivated by:

\begin{proposition}
The assignement $\pi\to T_\pi$ is categorical, in the sense that we have
$$T_\pi\otimes T_\sigma=T_{[\pi\sigma]}\quad,\quad
T_\pi T_\sigma=N^{c(\pi,\sigma)}T_{[^\sigma_\pi]}\quad,\quad
T_\pi^*=T_{\pi^*}$$
where $c(\pi,\sigma)$ are certain integers, coming from the erased components in the middle.
\end{proposition}

\begin{proof}
This is something that we know well from chapter 1, coming from some elementary computations for the above compositions, explained there.
\end{proof}

Let us axiomatize now the categories of partitions. The definition here, from \cite{bsp}, \cite{tw1}, which is very similar to the one from the classical case, is as follows:

\index{category of partitions}
\index{semicircle partition}
\index{horizontal concatenation}
\index{vertical concatenation}
\index{upside-down turning}

\begin{definition}
A collection of sets $D=\bigsqcup_{k,l}D(k,l)$ with $D(k,l)\subset P(k,l)$ is called a category of partitions when it has the following properties:
\begin{enumerate}
\item Stability under the horizontal concatenation, $(\pi,\sigma)\to[\pi\sigma]$.

\item Stability under vertical concatenation $(\pi,\sigma)\to[^\sigma_\pi]$, with matching middle symbols.

\item Stability under the upside-down turning $*$, with switching of colors, $\circ\leftrightarrow\bullet$.

\item Each set $P(k,k)$ contains the identity partition $||\ldots||$.

\item The sets $P(\emptyset,\circ\bullet)$ and $P(\emptyset,\bullet\circ)$ both contain the semicircle $\cap$.
\end{enumerate}
\end{definition} 

Generally speaking, the axioms in Definition 2.29 can be thought of as being a ``delinearized version'' of the categorical conditions which are verified by the Tannakian categories. We have in fact the following result, going back to \cite{bsp}:

\index{category of partitions}
\index{easy quantum group}

\begin{theorem}
Each category of partitions $D=(D(k,l))$ produces a family of compact quantum groups $G=(G_N)$, one for each $N\in\mathbb N$, via the formula
$$Hom(u^{\otimes k},u^{\otimes l})=span\left(T_\pi\Big|\pi\in D(k,l)\right)$$
which produces a Tannakian category, and therefore a closed subgroup $G_N\subset U_N^+$. The quantum groups which appear in this way are called easy.
\end{theorem}

\begin{proof}
This follows indeed from Woronowicz's Tannakian duality, in its ``soft'' form from \cite{mal}, as explained in Theorem 2.26. Indeed, let us set:
$$C_{kl}=span\left(T_\pi\Big|\pi\in D(k,l)\right)$$

By using the axioms in Definition 2.29, and the various categorical properties of the operation $\pi\to T_\pi$, from Proposition 2.28, we deduce that $C=(C_{kl})$ is a Tannakian category. Thus the Tannakian duality result applies, and gives the result.
\end{proof}

As a first application, we can formulate a general Brauer theorem, as follows:

\begin{theorem}
The basic classical and quantum rotation groups are all easy,
$$\xymatrix@R=51pt@C=50pt{
O_N^+\ar[r]&U_N^+\\
O_N\ar[u]\ar[r]&U_N\ar[u]}
\qquad\xymatrix@R=25pt@C=50pt{\\ :\\}\qquad
\xymatrix@R=18.5mm@C=16.2mm{
NC_2\ar[d]&\mathcal{NC}_2\ar[l]\ar[d]\\
P_2&\mathcal P_2\ar[l]}$$
with the quantum groups on the left corresponding to the categories on the right.
\end{theorem}

\begin{proof}
This is something that we already know for $O_N,U_N$, but since these results follow easily from those for $O_N^+,U_N^+$, let us just prove everything, as follows:

\medskip

(1) The quantum group $U_N^+$ is defined via the following relations:
$$u^*=u^{-1}\quad,\quad 
u^t=\bar{u}^{-1}$$ 

But, via our correspondence between partitions and maps, these relations tell us that the following two operators must be in the associated Tannakian category $C$:
$$T_\pi\quad,\quad \pi={\ }^{\,\cap}_{\circ\bullet}\ ,\,{\ }^{\,\cap}_{\bullet\circ}$$

Thus the associated Tannakian category is $C=span(T_\pi|\pi\in D)$, with:
$$D
=<{\ }^{\,\cap}_{\circ\bullet}\,\,,{\ }^{\,\cap}_{\bullet\circ}>
={\mathcal NC}_2$$

(2) The quantum group $O_N^+\subset U_N^+$ is defined by imposing the following relations:
$$u_{ij}=\bar{u}_{ij}$$

Thus, the following operators must be in the associated Tannakian category $C$:
$$T_\pi\quad,\quad\pi=|^{\hskip-1.32mm\circ}_{\hskip-1.32mm\bullet}\ ,\,|_{\hskip-1.32mm\circ}^{\hskip-1.32mm\bullet}$$

Thus the associated Tannakian category is $C=span(T_\pi|\pi\in D)$, with:
$$D
=<\mathcal{NC}_2,|^{\hskip-1.32mm\circ}_{\hskip-1.32mm\bullet},|_{\hskip-1.32mm\circ}^{\hskip-1.32mm\bullet}>
=NC_2$$

(3) The group $U_N\subset U_N^+$ is defined via the following relations:
$$[u_{ij},u_{kl}]=0\quad,\quad
[u_{ij},\bar{u}_{kl}]=0$$

Thus, the following operators must be in the associated Tannakian category $C$:
$$T_\pi\quad,\quad \pi={\slash\hskip-2.1mm\backslash}^{\hskip-2.5mm\circ\circ}_{\hskip-2.5mm\circ\circ}\ ,\,{\slash\hskip-2.1mm\backslash}^{\hskip-2.5mm\circ\bullet}_{\hskip-2.5mm\bullet\circ}$$

Thus the associated Tannakian category is $C=span(T_\pi|\pi\in D)$, with:
$$D
=<\mathcal{NC}_2,{\slash\hskip-2.1mm\backslash}^{\hskip-2.5mm\circ\circ}_{\hskip-2.5mm\circ\circ},{\slash\hskip-2.1mm\backslash}^{\hskip-2.5mm\circ\bullet}_{\hskip-2.5mm\bullet\circ}>
=\mathcal P_2$$

(4) In order to deal now with $O_N$, we can simply use the following formula: 
$$O_N=O_N^+\cap U_N$$

Indeed, at the categorical level, this formula tells us that the associated Tannakian category is given by $C=span(T_\pi|\pi\in D)$, with:
$$D
=<NC_2,\mathcal P_2>
=P_2$$

Thus, we are led to the conclusions in the statement.
\end{proof}

\section*{2d. The standard cube} 

Our purpose now is to unify and extend the squares from chapter 1 and from Theorem 2.31, consisting respectively of $H_N,K_N,O_N,U_N$ and of $O_N,U_N,O_N^+,U_N^+$, as to reach to the nice cube pictured at the beginning of the present chapter, namely:
$$
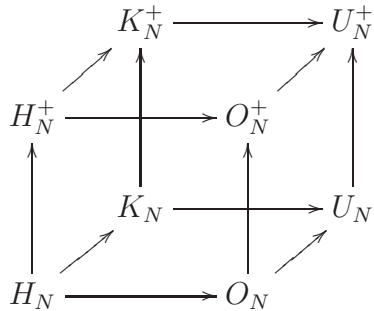
\xymatrix@R=16pt@C=16pt{
&K_N^+\ar[rr]&&U_N^+\\
H_N^+\ar[rr]\ar[ur]&&O_N^+\ar[ur]\\
&K_N\ar[rr]\ar[uu]&&U_N\ar[uu]\\
H_N\ar[uu]\ar[ur]\ar[rr]&&O_N\ar[uu]\ar[ur]
}$$

Thus, we need to talk about $H_N^+,K_N^+$, with this meaning both their definition, and easiness property. But for this, we first need to talk about the quantum permutation group $S_N^+$, again both definition, and easiness property. This will be something which is quite tricky, and will take some time. Following Wang \cite{wa2}, let us start with:

\index{quantum permutation group}
\index{free symmetric group}
\index{magic unitary}

\begin{theorem}
The following universal $C^*$-algebra, with magic meaning formed by projections $(p^2=p^*=p)$, summing up to $1$ on each row and each column,
$$C(S_N^+)=C^*\left((u_{ij})_{i,j=1,\ldots,N}\Big|u={\rm magic}\right)$$
is a Woronowicz algebra, with comultiplication, counit and antipode given by:
$$\Delta(u_{ij})=\sum_ku_{ik}\otimes u_{kj}\quad,\quad 
\varepsilon(u_{ij})=\delta_{ij}\quad,\quad 
S(u_{ij})=u_{ji}$$
Thus $S_N^+$ is a compact quantum group, called quantum permutation group, and the classical version of this quantum group is the usual permutation group $S_N$.
\end{theorem}

\begin{proof}
We have several assertions here, the idea being as follows:

\medskip

(1) As a first observation, the universal $C^*$-algebra in the statement is indeed well-defined, because the conditions $p^2=p^*=p$ satisfied by the coordinates give:
$$||u_{ij}||\leq1$$

In order to prove now that we have a Woronowicz algebra, we must construct maps $\Delta,\varepsilon,S$ given by the formulae in the statement. Consider the following matrices:
$$u^\Delta_{ij}=\sum_ku_{ik}\otimes u_{kj}\quad,\quad 
u^\varepsilon_{ij}=\delta_{ij}\quad,\quad 
u^S_{ij}=u_{ji}$$

Our claim is that, since $u$ is magic, so are these three matrices. Indeed, regarding $u^\Delta$, its entries are idempotents, as shown by the following computation:
$$(u_{ij}^\Delta)^2
=\sum_{kl}u_{ik}u_{il}\otimes u_{kj}u_{lj}
=\sum_{kl}\delta_{kl}u_{ik}\otimes\delta_{kl}u_{kj}
=u_{ij}^\Delta$$

These elements are self-adjoint as well, as shown by the following computation:
$$(u_{ij}^\Delta)^*
=\sum_ku_{ik}^*\otimes u_{kj}^*
=\sum_ku_{ik}\otimes u_{kj}
=u_{ij}^\Delta$$

The row and column sums for the matrix $u^\Delta$ can be computed as follows:
$$\sum_ju_{ij}^\Delta
=\sum_{jk}u_{ik}\otimes u_{kj}
=\sum_ku_{ik}\otimes 1
=1$$
$$\sum_iu_{ij}^\Delta
=\sum_{ik}u_{ik}\otimes u_{kj}
=\sum_k1\otimes u_{kj}
=1$$

Thus, $u^\Delta$ is magic. Regarding now $u^\varepsilon,u^S$, these matrices are magic too, and this for obvious reasons. Thus, all our three matrices $u^\Delta,u^\varepsilon,u^S$ are magic, so we can define $\Delta,\varepsilon,S$ by the formulae in the statement, by using the universality property of $C(S_N^+)$. 

\medskip

(2) Regarding the last assertion, consider the symmetric group $S_N$, viewed as permutation group of the $N$ coordinate axes of $\mathbb R^N$. The action of $S_N$ on the standard basis $e_1,\ldots,e_N\in\mathbb R^N$ being given by $\sigma:e_j\to e_{\sigma(j)}$, the coordinate functions on $S_N$ are:
$$u_{ij}=\chi\left(\sigma\in G\Big|\sigma(j)=i\right)$$

Observe also that the matrix $u=(u_{ij})$ that these functions form is magic, in the sense that its entries are projections, summing up to $1$ on each row and each column. Our claim now, which will prove the last assertion, is that we have the following formula:
$$C(S_N)=C^*_{comm}\left((u_{ij})_{i,j=1,\ldots,N}\Big|u={\rm magic}\right)$$

Indeed, the algebra $A$ on the right being commutative, by the Gelfand theorem it must be of the form $A=C(X)$, with $X$ being a certain compact space. Now since we have coordinates $u_{ij}:X\to\mathbb R$, we have an embedding $X\subset M_N(\mathbb R)$. Also, since we know that these coordinates form a magic matrix, the elements $g\in X$ must be 0-1 matrices, having exactly one 1 entry on each row and each column. Thus $X=S_N$, as desired.
\end{proof}

Still following Wang \cite{wa2}, we have the following surprising result:

\begin{theorem}
We have an embedding $S_N\subset S_N^+$, given at the algebra level by: 
$$u_{ij}\to\chi\left(\sigma\in S_N\Big|\sigma(j)=i\right)$$
This is an isomorphism at $N\leq3$, but not at $N\geq4$, where $S_N^+$ is not classical, nor finite.
\end{theorem} 

\begin{proof}
The fact that we have an embedding as above follows from Theorem 2.32. Regarding now the second assertion, we can prove this in four steps, as follows:

\medskip

\underline{Case $N=2$}. The fact that $S_2^+$ is indeed classical, and hence collapses to $S_2$, is trivial, because the $2\times2$ magic matrices are as follows, with $p$ being a projection:
$$U=\begin{pmatrix}p&1-p\\1-p&p\end{pmatrix}$$

Indeed, this shows that the entries of $U$ commute. Thus $C(S_2^+)$ is commutative, and so equals its biggest commutative quotient, which is $C(S_2)$. Thus, $S_2^+=S_2$.

\medskip

\underline{Case $N=3$}. By using the same argument as in the $N=2$ case, and the symmetries of the problem, it is enough to check that $u_{11},u_{22}$ commute. But this follows from:
\begin{eqnarray*}
u_{11}u_{22}
&=&u_{11}u_{22}(u_{11}+u_{12}+u_{13})\\
&=&u_{11}u_{22}u_{11}+u_{11}u_{22}u_{13}\\
&=&u_{11}u_{22}u_{11}+u_{11}(1-u_{21}-u_{23})u_{13}\\
&=&u_{11}u_{22}u_{11}
\end{eqnarray*}

Indeed, by applying the involution to this formula, we obtain that we have as well $u_{22}u_{11}=u_{11}u_{22}u_{11}$. Thus, we obtain $u_{11}u_{22}=u_{22}u_{11}$, as desired.

\medskip

\underline{Case $N=4$}. Consider the following matrix, with $p,q$ being projections:
$$U=\begin{pmatrix}
p&1-p&0&0\\
1-p&p&0&0\\
0&0&q&1-q\\
0&0&1-q&q
\end{pmatrix}$$ 

This matrix is magic, and we can choose $p,q\in B(H)$ as for the algebra $<p,q>$ to be noncommutative and infinite dimensional. We conclude that $C(S_4^+)$ is noncommutative and infinite dimensional as well, and so $S_4^+$ is non-classical and infinite, as claimed.

\medskip

\underline{Case $N\geq5$}. Here we can use the standard embedding $S_4^+\subset S_N^+$, obtained at the level of the corresponding magic matrices in the following way:
$$u\to\begin{pmatrix}u&0\\ 0&1_{N-4}\end{pmatrix}$$

Indeed, with this in hand, the fact that $S_4^+$ is a non-classical, infinite compact quantum group implies that $S_N^+$ with $N\geq5$ has these two properties as well.
\end{proof}

The above result might seem quite puzzling, but hey, we are doing quantum here, so take it easy. As a matter of doublechecking our findings, we are not wrong with our formalism, because as explained once again in \cite{wa2}, we have as well:

\index{coaction}
\index{counting measure}
\index{quantum permutation group}

\begin{theorem}
The quantum permutation group $S_N^+$ acts on the set $X=\{1,\ldots,N\}$, the corresponding coaction map $\Phi:C(X)\to C(X)\otimes C(S_N^+)$ being given by:
$$\Phi(e_i)=\sum_je_j\otimes u_{ji}$$
In fact, $S_N^+$ is the biggest compact quantum group acting on $X$, by leaving the counting measure invariant, in the sense that $(tr\otimes id)\Phi=tr(.)1$, where $tr(e_i)=\frac{1}{N},\forall i$.
\end{theorem}

\begin{proof}
Let us first determine when $\Phi$ is multiplicative. We have:
$$\Phi(e_i)\Phi(e_k)
=\sum_{jl}e_je_l\otimes u_{ji}u_{lk}
=\sum_je_j\otimes u_{ji}u_{jk}$$
$$\Phi(e_ie_k)
=\delta_{ik}\Phi(e_i)
=\delta_{ik}\sum_je_j\otimes u_{ji}$$

Thus, the multiplicativity of $\Phi$ is equivalent to the following conditions:
$$u_{ji}u_{jk}=\delta_{ik}u_{ji}\quad,\quad\forall i,j,k$$

The other conditions to be satisfied by $\Phi$ can be processed in a similar way, and we reach to the conclusion that $u$ must be magic, which gives the result. See \cite{wa2}.
\end{proof}

Getting now into easiness, we have the following result, which provides a more reasonable explanation for the liberation operation $S_N\to S_N^+$, and its mysteries:

\index{noncrossing partition}
\index{Brauer theorem}
\index{liberation}
\index{fork partition}

\begin{theorem}
The following hold:
\begin{enumerate}
\item The quantum groups $S_N,S_N^+$ are both easy, coming respectively from the categories $P,NC$ of partitions, and noncrossing partitions.  

\item Thus, $S_N\to S_N^+$ is just a regular easy quantum group liberation, coming from $D\to D\cap NC$ at the level of the associated categories of partitions.
\end{enumerate}
\end{theorem}

\begin{proof}
We already know the result for $S_N$, so we just need to prove the result for $S_N^+$. In order to do so, recall that the subgroup $S_N^+\subset O_N^+$ appears as follows:
$$C(S_N^+)=C(O_N^+)\Big\slash\Big<u={\rm magic}\Big>$$

In order to interpret the magic condition, consider the fork partition:
$$\mu\in P(2,1)$$

Given a corepresentation $u$, we have the following formulae:
$$(T_\mu u^{\otimes 2})_{i,jk}
=\sum_{lm}(T_\mu)_{i,lm}(u^{\otimes 2})_{lm,jk}
=u_{ij}u_{ik}$$
$$(uT_\mu)_{i,jk}
=\sum_lu_{il}(T_\mu)_{l,jk}
=\delta_{jk}u_{ij}$$

We conclude that we have the following equivalence:
$$T_\mu\in Hom(u^{\otimes 2},u)\iff u_{ij}u_{ik}=\delta_{jk}u_{ij},\forall i,j,k$$

The condition on the right being equivalent to the magic condition, we obtain:
$$C(S_N^+)=C(O_N^+)\Big\slash\Big<T_\mu\in Hom(u^{\otimes 2},u)\Big>$$

Thus $S_N^+$ is indeed easy, the corresponding category of partitions being:
$$D=<\mu>=NC$$

Finally, observe that this proves the result for $S_N$ too, because from the formula $S_N=S_N^+\cap O_N$ we obtain that the group $S_N$ is easy, coming from the category of partitions $D=<NC,P_2>=P$. Thus, we are led to the conclusions in the statement.
\end{proof}

With this understood, we can get now into quantum reflection groups, and reach our main objective for this chapter, namely having a nice cube. Let us start with:

\index{quantum reflection group}
\index{hyperoctahedral quantum group}
\index{free wreath product}

\begin{theorem}
We have quantum groups $H_N^+,K_N^+$, constructed as follows,
$$C(H_N^+)=C^*\left((u_{ij})_{i,j=1,\ldots,N}\Big|u_{ij}=u_{ij}^*,\,(u_{ij}^2)={\rm magic}\right)$$
$$C(K_N^+)=C^*\left((u_{ij})_{i,j=1,\ldots,N}\Big|[u_{ij},u_{ij}^*]=0,\,(u_{ij}u_{ij}^*)={\rm magic}\right)$$
which are liberations of $H_N,K_N$. Also, we have $H_N^+=\mathbb Z_2\wr_*S_N^+$, $K_N^+=\mathbb T\wr_*S_N^+$.
\end{theorem}

\begin{proof}
There is a long story with this result, following \cite{bb+}, \cite{bbc} and related papers. In the above form, which is something simplified, that we will need in what follows, the first assertion follows in the usual way, namely from the observation that if a matrix $u=(u_{ij})$ satisfies the relations in the statement, then so do the following matrices:
$$(u^\Delta)_{ij}=\sum_ku_{ik}\otimes u_{kj}\quad,\quad 
(u^\varepsilon)_{ij}=\delta_{ij}\quad,\quad 
(u^S)_{ij}=u_{ji}$$

As for the second assertion, the formulae there, with $\wr_*$ being a so-called free wreath product, are similar to the formulae $H_N=\mathbb Z_2\wr S_N$, $K_N=\mathbb T\wr S_N$ from the classical case, and their proof is routine. For more on this, we refer to \cite{bb+}, \cite{bbc}.
\end{proof}

Good news, we can now complete our cube, as follows:

\index{standard cube}
\index{quantum rotation}
\index{quantum reflection}
\index{Brauer theorem}

\begin{theorem}
We have quantum rotation and reflection groups, as follows,
$$\xymatrix@R=18pt@C=18pt{
&K_N^+\ar[rr]&&U_N^+\\
H_N^+\ar[rr]\ar[ur]&&O_N^+\ar[ur]\\
&K_N\ar[rr]\ar[uu]&&U_N\ar[uu]\\
H_N\ar[uu]\ar[ur]\ar[rr]&&O_N\ar[uu]\ar[ur]
}$$
which are all easy, the corresponding categories of partitions being as follows,
$$\xymatrix@R=19pt@C5pt{
&\mathcal{NC}_{even}\ar[dl]\ar[dd]&&\mathcal {NC}_2\ar[dl]\ar[ll]\ar[dd]\\
NC_{even}\ar[dd]&&NC_2\ar[dd]\ar[ll]\\
&\mathcal P_{even}\ar[dl]&&\mathcal P_2\ar[dl]\ar[ll]\\
P_{even}&&P_2\ar[ll]
}$$
with on top, the symbol $NC$ standing everywhere for noncrossing partitions.
\end{theorem}

\begin{proof}
This is something that we know for all quantum groups under consideration, except for $H_N^+,K_N^+$, and for these two quantum groups, the proof goes as follows:

\medskip

(1) We know that $H_N^+\subset O_N^+$ appears via the cubic relations, namely:
$$u_{ij}u_{ik}=u_{ji}u_{ki}=0\quad,\quad\forall j\neq k$$

Our claim now is that, in Tannakian terms, these relations reformulate as follows, with $\chi\in P(2,2)$ being the 1-block partition, joining all 4 points:
$$T_\chi\in End(u^{\otimes 2})$$

In order to prove our claim, observe first that we have, by definition of $T_\chi$:
$$T_\chi(e_i\otimes e_j)=\delta_{ij}e_i\otimes e_i$$

With this formula in hand, we have the following computation:
\begin{eqnarray*}
T_\chi u^{\otimes 2}(e_i\otimes e_j\otimes1)
&=&T_\chi\left(\sum_{abij}e_{ai}\otimes e_{bj}\otimes u_{ai}u_{bj}\right)(e_i\otimes e_j\otimes1)\\
&=&T_\chi\sum_{ab}e_a\otimes e_b\otimes u_{ai}u_{bj}\\
&=&\sum_ae_a\otimes e_a\otimes u_{ai}u_{aj}
\end{eqnarray*}

On the other hand, we have as well the following computation:
\begin{eqnarray*}
u^{\otimes 2}T_\chi(e_i\otimes e_j\otimes1)
&=&\delta_{ij}u^{\otimes 2}(e_i\otimes e_j\otimes1)\\
&=&\delta_{ij}\left(\sum_{abij}e_{ai}\otimes e_{bj}\otimes u_{ai}u_{bj}\right)(e_i\otimes e_j\otimes1)\\
&=&\delta_{ij}\sum_{ab}e_a\otimes e_b\otimes u_{ai}u_{bi}
\end{eqnarray*}

We conclude that $T_\chi u^{\otimes 2}=u^{\otimes 2}T_\chi$ means that $u$ is cubic, as desired. Thus, our claim is proved. But this shows that $H_N^+$ is easy, coming from the following category:
$$D=<\chi>=NC_{even}$$

(2) Regarding now $K_N^+$, the proof here is similar, leading this time to the category $\mathcal{NC}_{even}$ of noncrossing matching partitions. For details, we refer here to \cite{bb+}. 
\end{proof}

All the above is very nice. Perhaps not as satisfying as solving a Rubik's cube, or doing some Lego stuff, as a kid, but not bad. Further dealing with the above cube, with all sorts of enhancements, will keep us busy, for the rest of this book.

\section*{2e. Exercises} 

As usual with our exercises in this book, these will be for the most a mixture of more things to learn, and finishing computations not fully done in the above. To start with, since we are now into quantum, the following exercice is mandatory:

\begin{exercise}
Learn some quantum mechanics, from a true physicist, meaning from a book which is clearly advertised, on the front or back cover, not to be rigorous.
\end{exercise}

This piece of advice is serious. First because doing quantum groups, or mathematics in general, without any physics motivations will lead you nowhere, with the number of talented young people having tried this, and who ended up in depression, alcoholism, suicide or worse (remember Hell) being hard to count. And second, because quantum mechanics is non-trivial, to the point that even Einstein himself did not understand it, so take it easy there, just throw to trash excessive rigor, and use Love instead. Along the same lines now, but talking mathematics, here is a second exercise, equally important:

\begin{exercise}
Learn more about operator algebras and quantum spaces:
\begin{enumerate}
\item Operator algebras: more about $C^*$-algebras, and von Neumann algebras too.

\item Quantum spaces: have a look at free probability, and at $K$-theory too.

\item Quantum groups: read in detail the papers of Woronowicz \cite{wo1}, \cite{wo2}.
\end{enumerate}
\end{exercise}

At a more concrete level now, for truly learning quantum groups, nothing better than spending some time on $S_N^+$, which is the most exciting object around:

\begin{exercise}
Futher advance in your understanding of $S_N\to S_N^+$, as follows:
\begin{enumerate}
\item Prove that $S_3^+=S_3$, by using a new clever method, of your choice.

\item Prove that $S_4^+\neq S_4$, again by using a new method, of your choice.

\item Prove that $S_4^+$ is coamenable, while $S_5^+$ is not coamenable.

\item Can we talk about quantum permutations of finite quantum spaces?

\item If yes, can you prove that for $M_2$, given by $C(M_2)=M_2(\mathbb C)$, we get $SO_3$?

\item Based on this, can we say that $S_4^+$ should be a kind of twist of $SO_3$?
\end{enumerate}
\end{exercise}

Finally, in direct relation with easiness, besides completing of course the few proofs that were not given in the above with full details, we have:

\begin{exercise}
Extend the easiness theory that we have so far:
\begin{enumerate}
\item Towards bistochastic quantum groups, that you will have to define.

\item By looking for an intermediate liberation $O_N^*$, that you will have to find. 
\end{enumerate}
\end{exercise}

As usual, in what regards the last exercises, some of them might be quite difficult, but in case you do not find, no worries, we will be back to this, later in this book.

\chapter{Algebraic theory}

\section*{3a. Basic operations}

Wecome to easiness, again. Now that we learned the basics, time to make a to-do list, for the remainder of this book. There are many ways of proceeding here, with the presentation, and we have divided what is to be said in three parts, as follows:

\bigskip

(1) There is some general theory, of both algebraic and analytic nature, to be developed for the easy quantum groups $G\subset U_N^+$, starting from the axioms. We will do this in this chapter and in the next one, first with some algebra, in the present chapter, and then with some analytic results, mostly of probabilistic nature, in the next chapter.

\bigskip

(2) There are also many further examples to be studied, and classification results that can be obtained for them, and importantly, there is also some further general theory, that does not follow straight from the axioms, and requires verification on a case-by-case basis, by using the classification results. We will discuss this in Parts II and III.

\bigskip

(3) Finally, there are many closed subgroups $G\subset U_N^+$ which are not easy, but are not far from being easy either, with an illustrating example here being the symplectic group $Sp_N\subset U_N$, with $N\in2\mathbb N$. We will discuss all this, ``super-easiness'' theories, with super-easy meaning more general, and so harder than easy, in Part IV.

\bigskip

Getting started now, we will be first interested in the various operations that can be performed on the closed subgroups $G\subset U_N^+$, and how these behave in the easy case. Following Wang \cite{wa1}, the most basic operations on quantum groups are as follows:

\index{tensor product}
\index{free product}
\index{dual free product}

\begin{proposition}
The class of Woronowicz algebras is stable under taking:
\begin{enumerate}
\item Tensor products, $A=A'\otimes A''$, with $u=u'+u''$. At the quantum group level we obtain usual products, $G=G'\times G''$ and $\Gamma=\Gamma'\times\Gamma''$.

\item Free products, $A=A'*A''$, with $u=u'+u''$. At the quantum group level we obtain dual free products $G=G'\,\hat{*}\,G''$ and free products $\Gamma=\Gamma'*\Gamma''$.
\end{enumerate}
\end{proposition}

\begin{proof}
Everything here is clear from definitions. In addition to this, let us mention as well that we have $\int_{A'\otimes A''}=\int_{A'}\otimes\int_{A''}$ and $\int_{A'*A''}=\int_{A'}*\int_{A''}$. Also, the corepresentations of the above products can be explicitly computed. See Wang \cite{wa1}.
\end{proof}

In relation with easiness, we cannot expect much going on here, although an interesting question is the axiomatization of the class of products of easy quantum groups. Moving ahead now, here are some further basic operations, once again from Wang \cite{wa1}:

\index{Woronowicz subalgebra}
\index{quotient algebra}
\index{quantum subgroup}
\index{quotient quantum group}

\begin{proposition}
The class of Woronowicz algebras is stable under taking:
\begin{enumerate}
\item Subalgebras $A'=<u'_{ij}>\subset A$, with $u'$ being a corepresentation of $A$. At the quantum group level we obtain quotients $G\to G'$ and subgroups $\Gamma'\subset\Gamma$.

\item Quotients $A\to A'=A/I$, with $I$ being a Hopf ideal, $\Delta(I)\subset A\otimes I+I\otimes A$. At the quantum group level we obtain subgroups $G'\subset G$ and quotients $\Gamma\to\Gamma'$.
\end{enumerate}
\end{proposition}

\begin{proof}
Once again, everything is clear, and we have as well some straightforward supplementary results, regarding integration and corepresentations. See \cite{wa1}.
\end{proof}

As before with the product operations, all this is a bit too general, for a systematic easiness study. Here are now some further operations, which in contrast to the previous ones, will lead to some interesting theory in relation with easiness:

\index{intersection of subgroups}
\index{generation operation}

\begin{theorem}
The closed subgroups of $U_N^+$ are subject to operations as follows:
\begin{enumerate}
\item Intersection: $G\cap H$ is the biggest quantum subgroup of $G,H$.

\item Generation: $<G,H>$ is the smallest quantum group containing $G,H$.
\end{enumerate}
\end{theorem}

\begin{proof}
We must prove that the universal quantum groups in the statement exist indeed. For this purpose, let us pick writings as follows, with $I,J$ being Hopf ideals:
$$C(G)=C(U_N^+)/I\quad,\quad 
C(H)=C(U_N^+)/J$$

We can then construct our two universal quantum groups, as follows:
$$C(G\cap H)=C(U_N^+)/<I,J>$$
$$C(<G,H>)=C(U_N^+)/(I\cap J)$$

To be more precise, since $I,J$ are Hopf ideals, so are $<I,J>$ and $I\cap J$, so have indeed quantum groups, which have the needed universal properties, as desired.
\end{proof}

In practice now, what we have in Theorem 3.3 is quite theoretical, and in what concerns the operation $\cap$, this can be usually computed by using:

\begin{proposition}
Given subgroups $G,H\subset K$, appearing at the algebra level as follows, with $\mathcal R,\mathcal P$ being certain sets of polynomial $*$-relations between the coordinates $u_{ij}$,
$$C(G)=C(K)/\mathcal R\quad,\quad 
C(H)=C(K)/\mathcal P$$
the intersection $H\cap K$ is given by the formula $C(G\cap H)=C(K)/\{\mathcal R,\mathcal P\}$.
\end{proposition}

\begin{proof}
This follows from Theorem 3.3, or rather from its proof, and from the following trivial fact, regarding relations and ideals:
$$I=<\mathcal R>,J=<\mathcal P>
\quad\implies\quad <I,J>=<\mathcal R,\mathcal P>$$

Thus, we are led to the conclusion in the statement.
\end{proof}

In order to discuss now $<\,,>$, let us call Hopf image of a representation $C(K)\to A$ the smallest Hopf algebra quotient $C(L)$ producing a factorization as follows:
$$C(K)\to C(L)\to A$$

Here the fact that this quotient exists indeed is routine, by dividing by a suitable ideal. More generally, we can talk in the same way about the joint Hopf image of a family of representations $C(K)\to A_i$, and with this notion in hand, we have:

\index{Hopf image}

\begin{proposition}
Assuming $G,H\subset K$, the quantum group $<G,H>$ is such that
$$C(K)\to C(<G,H>)\to C(G),C(H)$$
is the joint Hopf image of the quotient maps $C(K)\to C(G),C(H)$.
\end{proposition}

\begin{proof}
In the particular case from the statement, the joint Hopf image appears as the smallest Hopf algebra quotient $C(L)$ producing factorizations as follows:
$$C(K)\to C(L)\to C(G),C(H)$$

We conclude from this that we have $L=<G,H>$, as desired.
\end{proof}

In order to reformulate the above operations in the easy setting, in terms of the associated categories of partitions, we first need a Tannakian reformulation of the constructions in Theorem 3.3. In the Tannakian setting, we have the following result:

\begin{theorem}
The intersection and generation operations $\cap$ and $<\,,>$ can be constructed via the Tannakian correspondence $G\to C_G$, as follows:
\begin{enumerate}
\item Intersection: defined via $C_{G\cap H}=<C_G,C_H>$.

\item Generation: defined via $C_{<G,H>}=C_G\cap C_H$.
\end{enumerate}
\end{theorem}

\begin{proof}
We know from Tannakian duality for quantum groups, in its soft form explained in chapter 2, that our two quantum groups $G,H\subset U_N^+$ appear as follows:
$$C(G)=C(U_N^+)\Big/\left<T\in Hom(u^{\otimes k},u^{\otimes l})\Big|\forall k,l,\forall T\in (C_G)_{kl}\right>$$
$$C(H)=C(U_N^+)\Big/\left<T\in Hom(u^{\otimes k},u^{\otimes l})\Big|\forall k,l,\forall T\in (C_H)_{kl}\right>$$

Now if we denote by $I,J$ the Hopf ideals on the right, and perform the operations from Theorem 3.3, we obtain quantum groups as follows:
$$C(G\cap H)=C(U_N^+)\Big/\left<T\in Hom(u^{\otimes k},u^{\otimes l})\Big|\forall k,l,\forall T\in<C_G,C_H>_{kl}\right>$$
$$C(<G,H>)=C(U_N^+)\Big/\left<T\in Hom(u^{\otimes k},u^{\otimes l})\Big|\forall k,l,\forall T\in (C_G\cap C_H)_{kl}\right>$$

Thus, once again by using the Tannakian duality for quantum groups, in its soft form explained in chapter 2, we are led to the conclusion in the statement.
\end{proof}

In relation now with our easiness questions, we first have the following result:

\begin{proposition}
Assuming that $G,H$ are easy, so is $G\cap H$, and we have
$$D_{G\cap H}=<D_G,D_H>$$
at the level of the corresponding categories of partitions.
\end{proposition}

\begin{proof}
We have indeed the following computation, with ``span'' standing for the corresponding span of linear maps, via the operation $\pi\to T_\pi$:
\begin{eqnarray*}
C_{G\cap H}
&=&<C_G,C_H>\\
&=&<span(D_G),span(D_H)>\\
&=&span(<D_G,D_H>)
\end{eqnarray*}

Thus, by Tannakian duality we obtain the result.
\end{proof}

Regarding the generation operation, the situation here is more complicated, as follows:

\begin{proposition}
Assuming that $G,H$ are easy, we have an inclusion 
$$<G,H>\subset\{G,H\}$$
coming from an inclusion of Tannakian categories as follows,
$$C_G\cap C_H\supset span(D_G\cap D_H)$$
where $\{G,H\}$ is the easy quantum group having as category of partitions $D_G\cap D_H$.
\end{proposition}

\begin{proof}
This follows from the following computation, with as before ``span'' standing for the corresponding span of linear maps, via the operation $\pi\to T_\pi$:
\begin{eqnarray*}
C_{<G,H>}
&=&C_G\cap C_H\\
&=&span(D_G)\cap span(D_H)\\
&\supset&span(D_G\cap D_H)
\end{eqnarray*}

Indeed, by Tannakian duality we obtain from this all the assertions.
\end{proof}

The problem now is that it is not clear if the inclusions in Proposition 3.8 are isomorphisms or not, and this not even under a supplementary $N>>0$ assumption. Thus, we have some problems here, and we must proceed as follows:

\index{easy generation}

\begin{theorem}
The intersection and easy generation operations $\cap$ and $\{\,,\}$ can be constructed via the Tannakian correspondence $G\to D_G$, as follows:
\begin{enumerate}
\item Intersection: defined via $D_{G\cap H}=<D_G,D_H>$.

\item Easy generation: defined via $D_{\{G,H\}}=D_G\cap D_H$.
\end{enumerate}
\end{theorem}

\begin{proof}
Here the situation is as follows:

\medskip

(1) This is a true result, coming from Proposition 3.7.

\medskip

(2) This is more of an empty statement, coming from Proposition 3.8.
\end{proof}

With the above notions in hand, we can now formulate a nice result, which improves our main result so far, the one from the end of chapter 2, as follows:

\index{intersection/generation diagram}
\index{standard cube}

\begin{theorem}
The basic quantum unitary and reflection groups, namely
$$
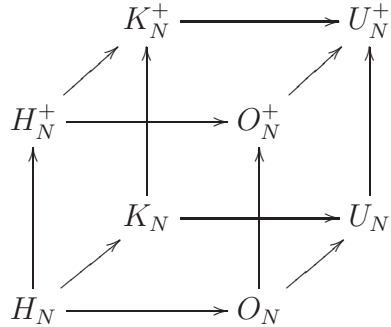
\xymatrix@R=18pt@C=18pt{
&K_N^+\ar[rr]&&U_N^+\\
H_N^+\ar[rr]\ar[ur]&&O_N^+\ar[ur]\\
&K_N\ar[rr]\ar[uu]&&U_N\ar[uu]\\
H_N\ar[uu]\ar[ur]\ar[rr]&&O_N\ar[uu]\ar[ur]
}$$
are all easy, and form an intersection/easy generation diagram, in the sense that any subsquare $P\subset Q,R\subset S$ of this diagram satisfies $Q\cap R=P$, $\{Q,R\}=S$.
\end{theorem}

\begin{proof}
We know from chapter 2 that the above quantum unitary and reflection groups are all easy, the corresponding categories of partitions being as follows:
$$
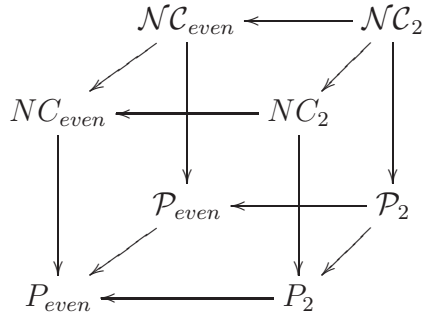
\xymatrix@R=18pt@C5pt{
&\mathcal{NC}_{even}\ar[dl]\ar[dd]&&\mathcal {NC}_2\ar[dl]\ar[ll]\ar[dd]\\
NC_{even}\ar[dd]&&NC_2\ar[dd]\ar[ll]\\
&\mathcal P_{even}\ar[dl]&&\mathcal P_2\ar[dl]\ar[ll]\\
P_{even}&&P_2\ar[ll]
}$$

Now since these categories form an intersection and generation diagram, the quantum groups form an intersection and easy generation diagram, as claimed.
\end{proof}

It is possible to further improve the above result, by proving that the diagram there is actually a plain generation diagram. However, this is something quite technical, requiring advanced quantum group techniques, and we will comment on this later.

\bigskip

By looking at the cube in Theorem 3.10, a natural idea in order to construct new quantum groups would be that of ``cutting it in half'', using $\cap$. Indeed, assume that we managed to find an intermediate easy quantum group $G$, for one of the 3 edges ending at $U_N^+$. Then, we can intersect all the vertices of the cube with $G$, and in practice this will produce 4 new easy quantum groups, including $G$ itself, cutting the cube in half.

\bigskip

This was for the method, and there are many things that can be said here, and all this will be slowly explored, later in this book. However, as a matter of having an illustration for this, and a useful class of new quantum groups, that can be used as examples or counterexamples for many things, let us work out the simplest instance of this method, by using an intermediate quantum group $U_N\subset G\subset U_N^+$. We can use here:

\begin{theorem}
We have an intermediate easy quantum group $U_N\subset U_N^*\subset U_N^+$, given by the following formula, and called half-classical unitary group,
$$C(U_N^*)=C(U_N^+)\Big/\left<abc=cba\Big|\forall a,b,c\in\{u_{ij},u_{ij}^*\}\right>$$
corresponding to the category of matching pairings $\mathcal P_2^*$ having the property that when relabelling clockwise the legs $\circ\bullet\circ\bullet\ldots$, the formula $\#\circ=\#\bullet$ holds in each block.
\end{theorem}

\begin{proof}
Here the fact that $U_N^*$ as constructed above is indeed a quantum group, lying as a proper intermediate subgroup $U_N\subset U_N^*\subset U_N^+$, can be checked via a routine computation, but the best is to view this via Tannakian duality. Indeed, the  half-commutation relations $abc=cba$ come from the map $T_{\slash\hskip-1.6mm\backslash\hskip-1.1mm|\hskip0.5mm}$ associated to the half-classical crossing:
$$\slash\hskip-2.0mm\backslash\hskip-1.7mm|\hskip0.5mm\in P(3,3)$$

Thus, by Tannakian duality, we are led to the first conclusions in the statement, and with the category of partitions associated to $U_N^*$ being as follows, with the convention that the symbol $\slash\hskip-2.0mm\backslash\hskip-1.7mm|\hskip0.5mm$ stands here for all 8 possible colorings of the diagram $\slash\hskip-2.0mm\backslash\hskip-1.7mm|\hskip0.5mm$:
$$D=<\slash\hskip-2.0mm\backslash\hskip-1.7mm|\hskip0.5mm>$$

Regarding now the explicit computation of $D$, observe that no matter how we color the legs of $\slash\hskip-2.0mm\backslash\hskip-1.7mm|\hskip0.5mm$, of course as for strings to join $\circ-\circ$ or $\bullet-\bullet$, we have a matching pairing, having in addition the property that when relabelling clockwise the legs $\circ\bullet\circ\bullet\ldots$, the formula $\#\circ=\#\bullet$ holds in each block. Thus, we have an inclusion as follows:
$$D\subset\mathcal P_2^*$$

On the other hand, by doing some standard combinatorics, we see that any element of $\mathcal P_2^*$ can be written as a composition of diagrams of type $\slash\hskip-2.0mm\backslash\hskip-1.7mm|\hskip0.5mm$, appearing with all its 8 possible colorings, as above. Thus, we have as well an inclusion as follows:
$$\mathcal P_2^*\subset D$$

But this shows that we have $D=\mathcal P_2^*$, which proves the last assertion.
\end{proof}

Now with the above construction in hand, we can perform our ``cutting the cube'' operation, and we are led to the following statement, improving what we have so far:

\begin{theorem}
We have easy quantum groups as follows, obtained via the commutation relations $abc=cba$, applied to the standard coordinates and their adjoints,
$$\xymatrix@R=32pt@C=18pt{
&K_N^*\ar[rr]&&U_N^*\\
H_N^*\ar[rr]\ar[ur]&&O_N^*\ar[ur]
}$$
which fit horizontally, in the middle, into the diagram of basic easy quantum groups,
$$\xymatrix@R=18pt@C=18pt{
&K_N^+\ar[rr]&&U_N^+\\
H_N^+\ar[rr]\ar[ur]&&O_N^+\ar[ur]\\
&K_N\ar[rr]\ar[uu]&&U_N\ar[uu]\\
H_N\ar[uu]\ar[ur]\ar[rr]&&O_N\ar[uu]\ar[ur]
}$$
with the enlarged diagram being an intersection/easy generation diagram.
\end{theorem}

\begin{proof}
There are several things going on here, and we will be quite brief:

\medskip

(1) First, the fact that we have indeed quantum groups as in the statement, which are all easy, follows from Proposition 3.7 and Theorem 3.11. 

\medskip

(2) Once again by using Proposition 3.7 and Theorem 3.11, we conclude as well that the categories of partitions for our new quantum groups are as follows:
$$\xymatrix@R=32pt@C=18pt{
&\ \mathcal P_{even}^*\ \ar[dl]&&\ \mathcal P_2^*\ar[ll]\ar[dl]\\
P_{even}^*\!&&P_2^*\ar[ll]
}$$

(3) The point now is that, when inserting this square diagram into the standard cube of categories of partitions, from the proof of Theorem 3.10, on the horizontal, in the middle, we obtain an intersection and generation diagram. Thus, the diagram formed by quantum groups is an intersection/easy generation diagram, as stated.

\medskip

(4) This was for the idea, and we will be back to this, with full details, in chapter 7 below, which will be dedicated to the half-liberation operation.
\end{proof}

As a conclusion to all this, the class of easy quantum groups behaves well with respect to $\cap$, and there is some interesting theory as well in relation with $<\,,>$. Some further interesting operations include various complexification operations, to be discussed in chapter 11, and the projective version operation, to be discussed in chapter 16.

\section*{3b. Envelopes, tori} 

The easy quantum groups appear as intermediate subgroups $S_N\subset G\subset U_N^+$, and regarding these latter quantum groups, we can say something about them, as follows:

\begin{theorem}
The closed subgroups $G\subset U_N^+$ which are homogeneous, in the sense that they contain the symmetric group $S_N$, and so appear as intermediate subgroups
$$S_N\subset G\subset U_N^+$$
are exacty those having as Tannakian category a collection $C=(C_{kl})$ with
$$span\left(T_\pi\Big|\pi\in\mathcal{NC}_2(k,l)\right)\subset C_{kl}\subset span\left(T_\pi\Big|\pi\in P(k,l)\right)$$
satisfying the general axioms for Tannakian categories. Moreover, $G$ is easy precisely when $C_{kl}=span(D(k,l))$, for a certain category of partitions $D=(D(k,l))\subset P$.
\end{theorem}

\begin{proof}
The inclusions in the statement are clear from the functoriality of the operation $G\to C$, and from our easiness results for $S_N,U_N^+$. As for the rest, the converse comes from Tannakian duality, and the last assertion is clear from definitions.
\end{proof}

The above result is something quite abstract, and ultimately rather trivial, but can be useful for many purposes. In order to further build on it, let us start with:

\begin{proposition}
Given a homogeneous quantum group $S_N\subset G\subset U_N^+$, with associated Tannakian category $C=(C_{kl})$, the sets
$$D(k,l)=\left\{\pi\in P(k,l)\Big|T_\pi\in C_{kl}\right\}$$ 
form a category of partitions $D\subset P$.
\end{proposition}

\begin{proof}
We use the basic properties of the correspondence $\pi\to T_\pi$, namely:
$$T_{[\pi\sigma]}=T_\pi\otimes T_\sigma\quad,\quad T_{[^\sigma_\pi]}\sim T_\pi T_\sigma\quad,\quad T_{\pi^*}=T_\pi^*$$

Together with the fact that $C$ is a tensor category, we deduce from these formulae that we have the following implications:
$$\pi,\sigma\in D\implies T_\pi,T_\sigma\in C\implies T_\pi\otimes T_\sigma\in C\implies T_{[\pi\sigma]}\in C\implies[\pi\sigma]\in D$$
$$\pi,\sigma\in D\implies T_\pi,T_\sigma\in C\implies T_\pi T_\sigma\in C\implies T_{[^\sigma_\pi]}\in C\implies[^\sigma_\pi]\in D$$
$$\pi\in D\implies T_\pi\in C\implies T_\pi^*\in C\implies T_{\pi^*}\in C\implies\pi^*\in D$$

Thus $D$ is indeed a category of partitions, as claimed.
\end{proof}

We can further refine the above observation, in the following way:

\begin{proposition}
Given a quantum group $S_N\subset G\subset U_N^+$, construct $D\subset P$ as above, and let $S_N\subset\bar{G}\subset U_N^+$ be the easy quantum group associated to $D$. Then:
\begin{enumerate}

\item We have $G\subset\bar{G}$, as subgroups of $U_N^+$.

\item $\bar{G}$ is the smallest easy quantum group containing $G$.

\item $G$ is easy precisely when $G\subset\bar{G}$ is an isomorphism.
\end{enumerate}
\end{proposition}

\begin{proof}
All the assertions are elementary, their proofs being as follows:

\medskip

(1) We know that the Tannakian category of $\bar{G}$ is given by:
$$\bar{C}_{kl}=span\left(T_\pi\Big|\pi\in D(k,l)\right)$$

Thus we have $\bar{C}\subset C$, and so $G\subset\bar{G}$, as subgroups of $U_N^+$.

\medskip

(2) Assuming that we have $G\subset G'$, with $G'$ easy, coming from a Tannakian category $C'=span(D')$, we must have $C'\subset\bar{C}$, and so $D'\subset D$. Thus, $\bar{G}\subset G'$, as desired.

\medskip

(3) This is a trivial consequence of (2).
\end{proof}

As an application, we can now introduce a notion of ``easy envelope'', as follows:

\begin{definition}
The easy envelope of a homogeneous quantum group $S_N\subset G\subset U_N^+$ is the easy quantum group $S_N\subset\bar{G}\subset U_N^+$ associated to the category of partitions
$$D(k,l)=\left\{\pi\in P(k,l)\Big|T_\pi\in C_{kl}\right\}$$ 
where $C=(C_{kl})$ is the Tannakian category of $G$.
\end{definition}

At the level of examples, most of the known quantum groups $S_N\subset G\subset U_N^+$ are in fact easy. However, there are many non-easy examples as well, and we will be back to this later. Moving ahead now, we can fine-tune all this, by using an arbitrary parameter $p\in\mathbb N$, which can be thought of as being an ``easiness level'', and we are led to:

\begin{theorem}
Given a quantum group $S_N\subset G\subset U_N^+$, consider the linear spaces
$$E^p_{kl}=\left\{\alpha_1T_{\pi_1}+\ldots+\alpha_pT_{\pi_p}\in C(k,l)\Big|\alpha_i\in\mathbb C,\pi_i\in P(k,l)\right\}$$
let $C^p$ be the smallest tensor category containing the family $E^p=(E^p_{kl})$, and let 
$$S_N\subset\bar{G}^p\subset U_N^+$$
be the quantum group corresponding to this category $C^p$. We have then:
$$G\subset\ldots\subset\bar{G}^3\subset\bar{G}^2\subset\bar{G}^1=\bar{G}\quad,\quad 
G=\bigcap_{p\in\mathbb N}\bar{G}^p$$
We say that $G$ has easiness level $p$ when $G=\bar{G}^p$, with $p\in\mathbb N$ chosen minimal.
\end{theorem}

\begin{proof}
This is something quite self-explanatory, the idea being as follows:

\medskip

(1) As a first observation, at $p=1$ we have $C^1=E^1=span(D)$, where $D$ is the category of partitions constructed in Proposition 3.14. Thus the quantum group $\bar{G}^1$ constructed above coincides with the easy envelope of $G$, from Definition 3.16.

\medskip

(2) In the general case, $p\in\mathbb N$, the family $E^p=(E^p_{kl})$ constructed above is not necessarily a tensor category, but we can consider the tensor category $C^p$ generated by it, as indicated. By definition of $E^p_{kl}$, and by using Proposition 3.15, these linear spaces $E^p_{kl}$ form an increasing filtration of $C_{kl}$. The same remains true when completing into tensor categories, so we have an increasing filtration, as follows:
$$C=\bigcup_{p\in\mathbb N}C^p$$

(3) At the quantum group level now, we obtain the decreasing intersection in the statement. Finally, the last statement is a definition, coming from all this.
\end{proof}

All this is quite abstract, and we will be back to it later, when having more examples of homogeneous quantum groups. Changing topics now, as a last purely algebraic question that we would like to discuss here, and which is of crucial importance, we have:

\begin{question}
What are the group dual subgroups $T\subset G$ of an easy quantum group $G$, and how can these be used in order to get information about $G$? 
\end{question}

To be more precise, recall from the classical Lie group theory that many things about a compact Lie group $G\subset U_N$ can be said once knowing its maximal torus $T\subset G$. So, what we are asking here for is a ``maximal torus theory'' for the closed subgroups $G\subset U_N^+$, and more specifically for the easy ones, that we are interested in.

\bigskip

Getting started now, following \cite{bv2}, we have the following definition:

\index{diagonal torus}
\index{torus}
\index{group dual}

\begin{proposition}
Given a closed subgroup $G\subset U_N^+$, consider its ``diagonal torus'', which is the closed subgroup $T\subset G$ constructed as follows:
$$C(T)=C(G)\Big/\left<u_{ij}=0\Big|\forall i\neq j\right>$$
This torus is then a group dual, $T=\widehat{\Lambda}$, where $\Lambda=<g_1,\ldots,g_N>$ is the discrete group generated by the elements $g_i=u_{ii}$, which are unitaries inside $C(T)$.
\end{proposition}

\begin{proof}
This is something going back to \cite{bv2}, which is elementary. The idea indeed is that since $u$ is unitary, its diagonal entries $g_i=u_{ii}$ are unitaries inside $C(T)$. Moreover, from $\Delta(u_{ij})=\sum_ku_{ik}\otimes u_{kj}$ we obtain, when passing inside the quotient:
$$\Delta(g_i)=g_i\otimes g_i$$

It follows that we have $C(T)=C^*(\Lambda)$, modulo identifying as usual the $C^*$-completions of the various group algebras, and so that we have $T=\widehat{\Lambda}$, as claimed.
\end{proof}

In the general easy case, the diagonal torus can be computed as follows:

\begin{theorem}
For an easy quantum group $G\subset U_N^+$, coming from a category of partitions $D\subset P$, the associated diagonal torus is $T=\widehat{\Gamma}$, with:
$$\Gamma=F_N\Big/\left<g_{i_1}\ldots g_{i_k}=g_{j_1}\ldots g_{j_l}\Big|\forall i,j,k,l,\exists\pi\in D(k,l),\delta_\pi\begin{pmatrix}i\\ j\end{pmatrix}\neq0\right>$$
Moreover, we can just use partitions $\pi$ which generate the category $D$.
\end{theorem}

\begin{proof}
Let $g_i=u_{ii}$ be the standard coordinates on the diagonal torus $T$, and set $g=diag(g_1,\ldots,g_N)$. We have then the following computation:
\begin{eqnarray*}
C(T)
&=&\left[C(U_N^+)\Big/\left<T_\pi\in Hom(u^{\otimes k},u^{\otimes l})\Big|\forall\pi\in D\right>\right]\Big/\left<u_{ij}=0\Big|\forall i\neq j\right>\\
&=&\left[C(U_N^+)\Big/\left<u_{ij}=0\Big|\forall i\neq j\right>\right]\Big/\left<T_\pi\in Hom(u^{\otimes k},u^{\otimes l})\Big|\forall\pi\in D\right>\\
&=&C^*(F_N)\Big/\left<T_\pi\in Hom(g^{\otimes k},g^{\otimes l})\Big|\forall\pi\in D\right>
\end{eqnarray*}

The associated discrete group, $\Gamma=\widehat{T}$, is therefore given by:
$$\Gamma=F_N\Big/\left<T_\pi\in Hom(g^{\otimes k},g^{\otimes l})\Big|\forall\pi\in D\right>$$

Now observe that, with $g=diag(g_1,\ldots,g_N)$ as above, we have:
$$T_\pi g^{\otimes k}(e_{i_1}\otimes\ldots\otimes e_{i_k})=\sum_{j_1\ldots j_l}\delta_\pi\binom{i}{j}e_{j_1}\otimes\ldots\otimes e_{j_l}\cdot g_{i_1}\ldots g_{i_k}$$
$$g^{\otimes l}T_\pi(e_{i_1}\otimes\ldots\otimes e_{i_k})=\sum_{j_1\ldots j_l}\delta_\pi\binom{i}{j}e_{j_1}\otimes\ldots\otimes e_{j_l}\cdot g_{j_1}\ldots g_{j_l}$$

We conclude that the relation $T_\pi\in Hom(g^{\otimes k},g^{\otimes l})$ reformulates as follows:
$$\sum_{j_1\ldots j_l}\delta_\pi\binom{i}{j}e_{j_1}\otimes\ldots\otimes e_{j_l}\cdot g_{i_1}\ldots g_{i_k}
=\sum_{j_1\ldots j_l}\delta_\pi\binom{i}{j}e_{j_1}\otimes\ldots\otimes e_{j_l}\cdot g_{j_1}\ldots g_{j_l}$$

Thus, we obtain the formula in the statement. Finally, the last assertion follows from Tannakian duality, because we can replace everywhere $D$ by a generating subset.
\end{proof}

More generally now, we have the following result, generalizing Proposition 3.19:

\index{spinned diagonal torus}

\begin{proposition}
Given a closed subgroup $G\subset U_N^+$ and a matrix $Q\in U_N$, we let $T_Q\subset G$ be the diagonal torus of $G$, with fundamental representation spinned by $Q$:
$$C(T_Q)=C(G)\Big/\left<(QuQ^*)_{ij}=0\Big|\forall i\neq j\right>$$
This torus is then a group dual, $T_Q=\widehat{\Lambda}_Q$, where $\Lambda_Q=<g_1,\ldots,g_N>$ is the discrete group generated by the elements $g_i=(QuQ^*)_{ii}$, which are unitaries inside $C(T_Q)$.
\end{proposition}

\begin{proof}
This follows from Proposition 3.19, because, as said in the statement, $T_Q$ is by definition a diagonal torus. Equivalently, since $v=QuQ^*$ is a unitary corepresentation, its diagonal entries $g_i=v_{ii}$, when regarded inside $C(T_Q)$, are unitaries, and satisfy:
$$\Delta(g_i)=g_i\otimes g_i$$

Thus $C(T_Q)$ is a group algebra, and more specifically we have $C(T_Q)=C^*(\Lambda_Q)$, where $\Lambda_Q=<g_1,\ldots,g_N>$ is the group in the statement, and this gives the result.
\end{proof}

The interest in the above construction comes from:

\begin{proposition}
Any torus $T\subset G$ appears as follows, for a certain $Q\in U_N$:
$$T\subset T_Q\subset G$$
In other words, any torus appears inside a spinned diagonal torus.
\end{proposition}

\begin{proof}
Given a torus $T\subset G$, we have $T\subset G\subset U_N^+$. On the other hand, we know that each torus $T\subset U_N^+$ has a fundamental corepresentation as follows, with $Q\in U_N$:
$$u=Q
\begin{pmatrix}g_1\\&\ddots\\&&g_N\end{pmatrix}
Q^*$$

But this shows that we have $T\subset T_Q$, and this gives the result.
\end{proof}

Getting back now to our motivations, namely ``maximal torus theory'' for the compact quantum groups, things here exist, but are rather conjectural. We first have:

\begin{theorem}
The following results hold, both for the compact Lie groups, and for the duals of the finitely generated discrete groups:
\begin{enumerate}
\item Generation: any closed quantum subgroup $G\subset U_N^+$ has the generation property $G=<T_Q|Q\in U_N>$. In other words, $G$ is generated by its tori.

\item Characters: if $G$ is connected, for any nonzero $P\in C(G)_{central}$ there exists $Q\in U_N$ such that $P$ becomes nonzero, when mapped into $C(T_Q)$.

\item Amenability: a closed subgroup $G\subset U_N^+$ is coamenable if and only if each of the tori $T_Q$ is coamenable, in the usual discrete group sense.

\item Growth: assuming $G\subset U_N^+$, the discrete quantum group $\widehat{G}$ has polynomial growth if and only if each the discrete groups $\widehat{T_Q}$ has polynomial growth.
\end{enumerate}
\end{theorem}

\begin{proof}
For group duals everything is trivial, and for classical groups, $G\subset U_N$, all this comes from standard facts from linear algebra and Lie theory, as follows:

\medskip

(1) Generation. We use the following formula, established above: 
$$T_Q=G\cap Q^*\mathbb T^NQ$$

Since any group element $U\in G$ is diagonalizable, $U=Q^*DQ$ with $Q\in U_N,D\in\mathbb T^N$, we have $U\in T_Q$ for this value of $Q\in U_N$, and this gives the result.

\medskip

(2) Characters. We can take here $Q\in U_N$ to be such that $QTQ^*\subset\mathbb T^N$, where $T\subset U_N$ is a maximal torus for $G$, and this gives the result.

\medskip

(3) Amenability. This conjecture holds trivially in the classical case, $G\subset U_N$, due to the fact that these latter quantum groups are all coamenable.

\medskip

(4) Growth. This is something non-trivial, well-known from the theory of compact Lie groups, and we refer here for instance to \cite{dpr}.
\end{proof}

The above statements are conjectured to hold for any compact quantum group, and for a number of verifications, we refer to \cite{bpa} and subsequent papers. In the easy case, it is possible to formulate a straightforward analogue of Theorem 3.20, in the spinned diagonal torus setting, and with this in hand, it is conjectured that the relevant standard tori of $G$, in relation with such conjectures, are the Fourier ones. See \cite{ba1}, \cite{bpa}.

\section*{3c. Gram determinants}

Let us discuss now a key algebraic problem, that we already met in the context of Proposition 3.8, namely the linear independence of the vectors $\xi_\pi$. We first have:

\index{order on partitions}
\index{supremum of partitions}
\index{lattice of partitions}

\begin{definition}
Let $P(k)$ be the set of partitions of $\{1,\ldots,k\}$, and $\pi,\sigma\in P(k)$.
\begin{enumerate}
\item We write $\pi\leq\sigma$ if each block of $\pi$ is contained in a block of $\sigma$.

\item We let $\pi\vee\sigma\in P(k)$ be the partition obtained by superposing $\pi,\sigma$.
\end{enumerate}
Also, we denote by $|.|$ the number of blocks of the partitions $\pi\in P(k)$.
\end{definition}

As an illustration here, at $k=2$ we have $P(2)=\{||,\sqcap\}$, and we have:
$$||\leq\sqcap$$

Also, at $k=3$ we have $P(3)=\{|||,\sqcap|,\sqcap\hskip-3.2mm{\ }_|\,,|\sqcap,\sqcap\hskip-0.7mm\sqcap\}$, and the order relation is as follows:
$$|||\ \leq\ \sqcap|\ ,\ \sqcap\hskip-3.2mm{\ }_|\ ,\ |\sqcap\ \leq\ \sqcap\hskip-0.7mm\sqcap$$

In relation with our linear independence questions, the idea will be that of using:

\index{Gram matrix}

\begin{proposition}
The Gram matrix of the vectors $\xi_\pi$ is given by the formula
$$<\xi_\pi,\xi_\sigma>=N^{|\pi\vee\sigma|}$$
where $\vee$ is the superposition operation, and $|.|$ is the number of blocks.
\end{proposition}

\begin{proof}
According to the formula of the vectors $\xi_\pi$, we have:
\begin{eqnarray*}
<\xi_\pi,\xi_\sigma>
&=&\sum_{i_1\ldots i_k}\delta_\pi(i_1,\ldots,i_k)\delta_\sigma(i_1,\ldots,i_k)\\
&=&\sum_{i_1\ldots i_k}\delta_{\pi\vee\sigma}(i_1,\ldots,i_k)\\
&=&N^{|\pi\vee\sigma|}
\end{eqnarray*}

Thus, we have obtained the formula in the statement.
\end{proof}

In order to study the Gram matrix $G_k(\pi,\sigma)=N^{|\pi\vee\sigma|}$, and more specifically to compute its determinant, we will use several standard facts about partitions. We have:

\index{M\"obius function}

\begin{definition}
The M\"obius function of any lattice, and so of $P$, is given by
$$\mu(\pi,\sigma)=\begin{cases}
1&{\rm if}\ \pi=\sigma\\
-\sum_{\pi\leq\tau<\sigma}\mu(\pi,\tau)&{\rm if}\ \pi<\sigma\\
0&{\rm if}\ \pi\not\leq\sigma
\end{cases}$$
with the construction being performed by recurrence.
\end{definition}

As an illustration here, for $P(2)=\{||,\sqcap\}$, we have by definition:
$$\mu(||,||)=\mu(\sqcap,\sqcap)=1$$

Also, $||<\sqcap$, with no intermediate partition in between, so we obtain:
$$\mu(||,\sqcap)=-\mu(||,||)=-1$$

Finally, we have $\sqcap\not\leq||$, and so we have as well the following formula:
$$\mu(\sqcap,||)=0$$

Back to the general case now, the main interest in the M\"obius function comes from the M\"obius inversion formula, which states that the following happens:
$$f(\sigma)=\sum_{\pi\leq\sigma}g(\pi)\quad
\implies\quad g(\sigma)=\sum_{\pi\leq\sigma}\mu(\pi,\sigma)f(\pi)$$

In linear algebra terms, the statement and proof of this formula are as follows:

\index{M\"obius matrix}
\index{M\"obius inversion}

\begin{theorem}
The inverse of the adjacency matrix of $P(k)$, given by
$$A_k(\pi,\sigma)=\begin{cases}
1&{\rm if}\ \pi\leq\sigma\\
0&{\rm if}\ \pi\not\leq\sigma
\end{cases}$$
is the M\"obius matrix of $P$, given by $M_k(\pi,\sigma)=\mu(\pi,\sigma)$.
\end{theorem}

\begin{proof}
This is well-known, coming for instance from the fact that $A_k$ is upper triangular. Indeed, when inverting, we are led into the recurrence from Definition 3.26.
\end{proof}

As an illustration, for $P(2)$ the formula $M_2=A_2^{-1}$ appears as follows:
$$\begin{pmatrix}1&-1\\ 0&1\end{pmatrix}=
\begin{pmatrix}1&1\\ 0&1\end{pmatrix}^{-1}$$

Now back to our Gram matrix considerations, we have the following key result:

\begin{proposition}
The Gram matrix of the vectors $\xi_\pi$ with $\pi\in P(k)$,
$$G_{\pi\sigma}=N^{|\pi\vee\sigma|}$$
decomposes as a product of upper/lower triangular matrices, $G_k=A_kL_k$, where
$$L_k(\pi,\sigma)=
\begin{cases}
N(N-1)\ldots(N-|\pi|+1)&{\rm if}\ \sigma\leq\pi\\
0&{\rm otherwise}
\end{cases}$$
and where $A_k$ is the adjacency matrix of $P(k)$.
\end{proposition}

\begin{proof}
We have the following computation, based on Proposition 3.25:
\begin{eqnarray*}
G_k(\pi,\sigma)
&=&N^{|\pi\vee\sigma|}\\
&=&\#\left\{i_1,\ldots,i_k\in\{1,\ldots,N\}\Big|\ker i\geq\pi\vee\sigma\right\}\\
&=&\sum_{\tau\geq\pi\vee\sigma}\#\left\{i_1,\ldots,i_k\in\{1,\ldots,N\}\Big|\ker i=\tau\right\}\\
&=&\sum_{\tau\geq\pi\vee\sigma}N(N-1)\ldots(N-|\tau|+1)
\end{eqnarray*}

According now to the definition of $A_k,L_k$, this formula reads:
\begin{eqnarray*}
G_k(\pi,\sigma)
&=&\sum_{\tau\geq\pi}L_k(\tau,\sigma)\\
&=&\sum_\tau A_k(\pi,\tau)L_k(\tau,\sigma)\\
&=&(A_kL_k)(\pi,\sigma)
\end{eqnarray*}

Thus, we are led to the formula in the statement.
\end{proof}

As an illustration for the above result, at $k=2$ we have $P(2)=\{||,\sqcap\}$, and the above decomposition $G_2=A_2L_2$ appears as follows:
$$\begin{pmatrix}N^2&N\\ N&N\end{pmatrix}
=\begin{pmatrix}1&1\\ 0&1\end{pmatrix}
\begin{pmatrix}N^2-N&0\\N&N\end{pmatrix}$$

We are led in this way to the following formula, due to Lindst\"om \cite{lin}:

\index{Gram matrix}
\index{Gram determinant}
\index{Lindst\"om formula}

\begin{theorem}
The determinant of the Gram matrix $G_k$ is given by
$$\det(G_k)=\prod_{\pi\in P(k)}\frac{N!}{(N-|\pi|)!}$$
with the convention that in the case $N<k$ we obtain $0$.
\end{theorem}

\begin{proof}
If we order $P(k)$ as usual, with respect to the number of blocks, and then lexicographically, $A_k$ is upper triangular, and $L_k$ is lower triangular. Thus, we have:
\begin{eqnarray*}
\det(G_k)
&=&\det(A_k)\det(L_k)\\
&=&\det(L_k)\\
&=&\prod_\pi L_k(\pi,\pi)\\
&=&\prod_\pi N(N-1)\ldots(N-|\pi|+1)
\end{eqnarray*}

Thus, we are led to the formula in the statement.
\end{proof}

The above computation can be thought of as corresponding to the group $S_N$, but we can do such things for any easy quantum group. As a first illustration, let us discuss the case of the orthogonal group $O_N$. Here the combinatorics is that of the Young diagrams. We denote by $|.|$ the number of boxes, and we use quantity $f^\lambda$, which gives the number of standard Young tableaux of shape $\lambda$. We have then the following result:

\index{Young tableaux}
\index{Gram determinant}

\begin{theorem}
The determinant of the Gram matrix of $O_N$ is given by
$$\det(G_{kN})=\prod_{|\lambda|=k/2}f_N(\lambda)^{f^{2\lambda}}$$
where the quantities on the right are $f_N(\lambda)=\prod_{(i,j)\in\lambda}(N+2j-i-1)$.
\end{theorem}

\begin{proof}
For the group $O_N$ the Gram matrix is diagonalizable, as follows:
$$G_{kN}=\sum_{|\lambda|=k/2}f_N(\lambda)P_{2\lambda}$$

Here $1=\sum P_{2\lambda}$ is the standard partition of unity associated to the Young diagrams having $k/2$ boxes, and the coefficients $f_N(\lambda)$ are those in the statement. Now since we have $Tr(P_{2\lambda})=f^{2\lambda}$, this gives the formula in the statement. For details here, see \cite{bcu}.
\end{proof}

In order to deal now with $O_N^+,S_N^+$, we will need the following well-known fact:

\index{fattening of partitions}
\index{shrinking of partitions}

\begin{proposition}
We have a bijection $NC(k)\simeq NC_2(2k)$, as follows:
\begin{enumerate}
\item The application $NC(k)\to NC_2(2k)$ is the ``fattening'' one, obtained by doubling all the legs, and doubling all the strings as well.

\item Its inverse $NC_2(2k)\to NC(k)$ is the ``shrinking'' application, obtained by collapsing pairs of consecutive neighbors.
\end{enumerate}
\end{proposition}

\begin{proof}
The fact that the above two operations are indeed inverse to each other is clear, by drawing pictures, and computing the corresponding compositions.
\end{proof}

At the level of the associated Gram matrices, the result is as follows:
\begin{proposition}
The Gram matrices of $NC_2(2k)\simeq NC(k)$ are related by
$$G_{2k,n}(\pi,\sigma)=n^k(\Delta_{kn}^{-1}G_{k,n^2}\Delta_{kn}^{-1})(\pi',\sigma')$$
where $\pi\to\pi'$ is the shrinking operation, and $\Delta_{kn}$ is the diagonal of $G_{kn}$.
\end{proposition}

\begin{proof}
In the context of the bijection from Proposition 3.31, we have:
$$|\pi\vee\sigma|=k+2|\pi'\vee\sigma'|-|\pi'|-|\sigma'|$$

We therefore have the following formula, valid for any $n\in\mathbb N$:
$$n^{|\pi\vee\sigma|}=n^{k+2|\pi'\vee\sigma'|-|\pi'|-|\sigma'|}$$

Thus, we are led to the formula in the statement.
\end{proof}

Now back to $O_N^+,S_N^+$, let us begin with some examples. We first have:

\begin{proposition}
The first Gram matrices and determinants for $O_N^+$ are
$$\det\begin{pmatrix}N^2&N\\N&N^2\end{pmatrix}=N^2(N^2-1)$$
$$\det\begin{pmatrix}
N^3&N^2&N^2&N^2&N\\
N^2&N^3&N&N&N^2\\
N^2&N&N^3&N&N^2\\
N^2&N&N&N^3&N^2\\
N&N^2&N^2&N^2&N^3
\end{pmatrix}=N^5(N^2-1)^4(N^2-2)$$
with the matrices being written by using the lexicographic order on $NC_2(2k)$.
\end{proposition}

\begin{proof}
The formula at $k=2$, where $NC_2(4)=\{\sqcap\sqcap,\bigcap\hskip-4.9mm{\ }_\cap\,\}$, is clear from definitions. At $k=3$ however, things are tricky. The partitions here are as follows:
$$NC(3)=\{|||,\sqcap|,\sqcap\hskip-3.2mm{\ }_|\,,|\sqcap,\sqcap\hskip-0.7mm\sqcap\}$$

The Gram matrix and its determinant are, according to Theorem 3.29:
$$\det\begin{pmatrix}
N^3&N^2&N^2&N^2&N\\
N^2&N^2&N&N&N\\
N^2&N&N^2&N&N\\
N^2&N&N&N^2&N\\
N&N&N&N&N
\end{pmatrix}=N^5(N-1)^4(N-2)$$

By using now Proposition 3.32, this gives the formula in the statement.
\end{proof}

In general, such tricks won't work, because $NC(k)$ is strictly smaller than $P(k)$ at $k\geq4$. However, following Di Francesco \cite{dif}, we have the following result:

\index{meander determinant}
\index{Gram determinant}
\index{Chebycheff polynomials}
\index{Di Francesco formula}

\begin{theorem}
The determinant of the Gram matrix for $O_N^+$ is given by
$$\det(G_{kN})=\prod_{r=1}^{[k/2]}P_r(N)^{d_{k/2,r}}$$
where $P_r$ are the Chebycheff polynomials, given by
$$P_0=1\quad,\quad 
P_1=X\quad,\quad 
P_{r+1}=XP_r-P_{r-1}$$
and $d_{kr}=f_{kr}-f_{k,r+1}$, with $f_{kr}$ being the following numbers, depending on $k,r\in\mathbb Z$,
$$f_{kr}=\binom{2k}{k-r}-\binom{2k}{k-r-1}$$
with the convention $f_{kr}=0$ for $k\notin\mathbb Z$. 
\end{theorem}

\begin{proof}
This is something quite technical, obtained by using a decomposition as follows of the Gram matrix $G_{kN}$, with the matrix $T_{kN}$ being lower triangular:
$$G_{kN}=T_{kN}T_{kN}^t$$

Thus, a bit as in the proof of the Lindst\"om formula, we obtain the result, but the problem lies however in the construction of $T_{kN}$, which is non-trivial. See \cite{dif}.
\end{proof}

We refer to \cite{bcu} for further details regarding the above result, including a short proof, based on the bipartite graph planar algebra combinatorics developed by Jones in \cite{jo4}. Moving ahead now, regarding $S_N^+$, we have here the following formula, which is quite similar, obtained via shrinking, also from Di Francesco \cite{dif}:

\index{meander determinant}
\index{Gram determinant}
\index{Chebycheff polynomials}
\index{Di Francesco formula}

\begin{theorem}
The determinant of the Gram matrix for $S_N^+$ is given by
$$\det(G_{kN})=(\sqrt{N})^{a_k}\prod_{r=1}^kP_r(\sqrt{N})^{d_{kr}}$$
where $P_r$ are the Chebycheff polynomials, given by
$$P_0=1\quad,\quad 
P_1=X\quad,\quad 
P_{r+1}=XP_r-P_{r-1}$$
and $d_{kr}=f_{kr}-f_{k,r+1}$, with $f_{kr}$ being the following numbers, depending on $k,r\in\mathbb Z$,
$$f_{kr}=\binom{2k}{k-r}-\binom{2k}{k-r-1}$$
with the convention $f_{kr}=0$ for $k\notin\mathbb Z$, and where $a_k=\sum_{\pi\in \mathcal P(k)}(2|\pi|-k)$.
\end{theorem}

\begin{proof}
This follows indeed from Theorem 3.34, by using Proposition 3.32.
\end{proof}

There are many more things that can be said about the Gram determinants of the easy quantum groups, some being open problems, and we refer here to \cite{bcu}.

\section*{3d. Representation theory} 

Let us discuss now the representation theory of the easy quantum groups $G\subset U_N^+$. In the classical case, $G\subset U_N$, things are quite tricky, and this even for simple groups like $S_N,O_N$, and we had a taste of this in the previous section, when talking about Gram determinants. So, we will basically restrict the attention here to the free case:

\begin{definition}
An easy quantum group $S_N\subset G\subset U_N^+$, coming from a category $\mathcal{NC}_2\subset D\subset P$, is called free when the following equivalent conditions are satisfied:
\begin{enumerate}
\item $G$ appears as an intermediate quantum group $S_N^+\subset G\subset U_N^+$.

\item $D$ appears as an intermediate category of partitions $\mathcal{NC}_2\subset D\subset NC$.
\end{enumerate}
Equivalently, via fattening, $G$ must appear from a category of Temperley-Lieb diagrams.
\end{definition}

Getting now to the representation theory problem, generally speaking, following the modern approach from \cite{fwe}, the idea is to use the following formula:
$$End(u^{\otimes k})=span\left(T_\pi\Big|\pi\in D(k,k)\right)$$

Indeed, assuming as above that we are in the free case, $D(k,k)\subset NC(k,k)$, the partitions $\pi\in D(k,k)$ which are symmetric, $\pi=\pi^*$, will produce linear maps $T_\pi$ which are projections, up to a scalar. Thus, with a bit of care, we will be able to write a formula as follows, with the sum being over certain symmetric partitions $D(k,k)\subset NC(k,k)$, and with $T_\pi'$ being projections, appearing as certain modifications of the linear maps $T_\pi$:
$$1=\sum_\pi T_\pi'$$

But this formula, once properly formulated, is exactly what we need, allowing us to decompose each Peter-Weyl representation $u^{\otimes k}$ into a sum of irreducibles. 

\bigskip

In practice now, all this is quite complicated, and besides assuming that we are in the free case, we must fatten the partitions, in order to simplify a bit the computations in \cite{fwe}, in connection with the middle components which appear when concatenating. So, it is tempting to go instead to the original proofs, from papers written in the 90s and 00s, doing everything in an elementary way, with a cheap Frobenius trick. We will need:

\index{Catalan number}
\index{Clebsch-Gordan rules}
\index{rotation group}
\index{Frobenius trick}

\begin{theorem}
The irreducible representations of $SU_2$ can be labeled by positive integers, $r_k$ with $k\in\mathbb N$, and the fusion rules for these representations are:
$$r_k\otimes r_l=r_{|k-l|}+r_{|k-l|+2}+\ldots+r_{k+l}$$ 
The dimensions of these representations are $\dim r_k=k+1$.
\end{theorem}

\begin{proof}
Our claim is that we can construct, by recurrence on $k\in\mathbb N$, a sequence $r_k$ of irreducible, self-adjoint and distinct representations of $SU_2$, satisfying:
$$r_0=1\quad,\quad
r_1=u\quad,\quad 
r_{k-1}\otimes r_1=r_{k-2}+r_k$$

Indeed, assume that $r_0,\ldots,r_{k-1}$ are constructed, and let us construct $r_k$. We have:
$$r_{k-2}\otimes r_1=r_{k-3}+r_{k-1}$$

Thus $r_{k-1}\subset r_{k-2}\otimes r_1$, and since $r_{k-2}$ is irreducible, by Frobenius we have:
$$r_{k-2}\subset r_{k-1}\otimes r_1$$

We conclude there exists a certain representation $r_k$ such that:
$$r_{k-1}\otimes r_1=r_{k-2}+r_k$$

By recurrence, $r_k$ is self-adjoint. Now observe that according to our recurrence formula, we can split $u^{\otimes k}$ as a sum of the following type, with positive coefficients:  
$$u^{\otimes k}=c_kr_k+c_{k-2}r_{k-2}+c_{k-4}r_{k-4}+\ldots$$

We conclude by Peter-Weyl that we have an inequality as follows, with equality precisely when $r_k$ is irreducible, and non-equivalent to the other summands $r_i$:
$$\sum_ic_i^2\leq\dim(End(u^{\otimes k}))=\int_{SU_2}\chi^{2k}$$

Now recall that we have a well-known isomorphism $SU_2\simeq S^3_\mathbb R$, coming from:
$$SU_2=\left\{\begin{pmatrix}x+iy&z+it\\ -z+it&x-iy\end{pmatrix}\Big|x^2+y^2+z^2+t^2=1\right\}$$

In this picture the moments of the main character $\chi=2x$ can be computed via spherical coordinates and some calculus, and follow to be the Catalan numbers:
$$\int_{SU_2}\chi^{2k}=C_k$$

On the other hand, some straightforward combinatorics, based on the fusion rules, shows that we have $\sum_ic_i^2=C_k$ as well. Thus, the estimate found above reads:
$$C_k=\sum_ic_i^2\leq\dim(End(u^{\otimes k}))=\int_{SU_2}\chi^{2k}=C_k$$

We conclude that we have equality in our estimate, so our representation $r_k$ is irreducible, and non-equivalent to $r_{k-2},r_{k-4},\ldots$ Moreover, this representation $r_k$ is not equivalent to $r_{k-1},r_{k-3},\ldots$ either, with this coming from $r_p\subset u^{\otimes p}$, and from:
$$\dim(Fix(u^{\otimes 2s+1}))=\int_{SU_2}\chi^{2s+1}=0$$

Thus, we proved our claim. Now since each irreducible representation of $SU_2$ must appear in some tensor power $u^{\otimes k}$, and we know how to decompose each $u^{\otimes k}$ into sums of representations $r_k$, these representations $r_k$ are all the irreducible representations of $SU_2$, and we are done. As for the dimension formula, this is clear by recurrence.
\end{proof}

Quite remarkably, a similar result holds for $O_N^+$, as follows:

\index{free orthogonal group}
\index{Catalan number}
\index{Clebsch-Gordan rules}
\index{Frobenius trick}

\begin{theorem}
The irreducible representations of $O_N^+$  with $N\geq2$ can be labeled by positive integers, $r_k$ with $k\in\mathbb N$, the fusion rules for these representations are
$$r_k\otimes r_l=r_{|k-l|}+r_{|k-l|+2}+\ldots+r_{k+l}$$ 
and the dimensions are $\dim r_k=(q^{k+1}-q^{-k-1})/(q-q^{-1})$, with $q+q^{-1}=N$.
\end{theorem}

\begin{proof}
The idea is to skilfully recycle the proof of Theorem 3.37. As before, our claim is that we can construct, by recurrence on $k\in\mathbb N$, a sequence $r_0,r_1,r_2,\ldots$ of irreducible, self-adjoint and distinct representations of $O_N^+$, satisfying:
$$r_0=1\quad,\quad
r_1=u\quad,\quad 
r_{k-1}\otimes r_1=r_{k-2}+r_k$$

In order to do so, we can use as before $r_{k-2}\otimes r_1=r_{k-3}+r_{k-1}$ and Frobenius, and we conclude there exists a certain representation $r_k$ such that:
$$r_{k-1}\otimes r_1=r_{k-2}+r_k$$

As a first observation, $r_k$ is self-adjoint, because its character is a certain polynomial with integer coefficients in $\chi$, which is self-adjoint. In order to prove now that $r_k$ is irreducible, and non-equivalent to $r_0,\ldots,r_{k-1}$, let us split as before $u^{\otimes k}$, as follows:  
$$u^{\otimes k}=c_kr_k+c_{k-2}r_{k-2}+c_{k-4}r_{k-4}+\ldots$$

The point now is that we have the following equalities and inequalities:
$$C_k=\sum_ic_i^2\leq\dim(End(u^{\otimes k}))\leq |NC_2(k,k)|=C_k$$

Indeed, the equality at left is clear as before, then comes a standard inequality, then an inequality coming from easiness, then a standard equality. Thus, we have equality, so $r_k$ is irreducible, and non-equivalent to $r_{k-2},r_{k-4},\ldots$ Moreover, $r_k$ is not equivalent to $r_{k-1},r_{k-3},\ldots$ either, by using the same argument as for $SU_2$, and the end of the proof is exactly as for $SU_2$. As for dimensions, by recurrence we obtain, with $q+q^{-1}=N$:
$$\dim r_k=q^k+q^{k-2}+\ldots+q^{-k+2}+q^{-k}$$

But this gives the dimension formula in the statement, and we are done.
\end{proof}

It is possible to use similar methods for the other free quantum groups, or follow the modern approach from \cite{fwe}, as explained above, and we first obtain in this way:

\begin{theorem}
The irreducible representations of the quantum group $U_N^+$ with $N\geq2$ can be indexed by $\mathbb N*\mathbb N$, with fusion rules as follows,
$$r_k\otimes r_l=\sum_{k=xy,l=\bar{y}z}r_{xz}$$
and the corresponding dimensions $\dim r_k$ can be computed by recurrence.
\end{theorem}

\begin{proof}
This follows indeed by using the same arguments as in the proof of Theorem 3.38, via recurrence and a Frobenius trick, and for details here, we refer to \cite{ba1}. Alternatively, we have a proof based on Theorem 3.38 itself, going as follows:

\medskip

(1) Given a Woronowicz algebra $(A,v)$, we can construct its free complexification $(\tilde{A},\tilde{v})$ as follows, with the fact that we get indeed a Woronowicz algebra being clear:
$$\tilde{A}=<zv_{ij}>\subset C(\mathbb T)*A\quad,\quad\tilde{v}=zv$$

At the quantum group level we obtain in this way free complexifications of the corresponding compact and discrete quantum groups $G$ and $\Gamma$, denoted $\widetilde{G}$ and $\widetilde{\Gamma}$.

\medskip

(2) Now observe that we have embeddings as follows, with the first one coming by using the counit, and the second one coming from the universality property of $U_N^+$:
$$O_N^+
\subset\widetilde{O_N^+}
\subset U_N^+$$

Our claim is that the embedding on the right is an isomorphism, so that we have an isomorphism of compact quantum groups, as follows:
$$U_N^+=\widetilde{O_N^+}$$

(3) In order to prove this claim, which will afterwards lead to our theorem, let us denote by $v,zv,u$ the fundamental representations of the 3 quantum groups from (2). At the level of the associated Hom spaces we obtain reverse inclusions, as follows:
$$Hom(v^{\otimes k},v^{\otimes l})
\supset Hom((zv)^{\otimes k},(zv)^{\otimes l})
\supset Hom(u^{\otimes k},u^{\otimes l})$$

But the spaces on the left and on the right are known from chapter 2, the easiness result there stating that these spaces are as follows:
$$span\left(T_\pi\Big|\pi\in NC_2(k,l)\right)\supset span\left(T_\pi\Big|\pi\in\mathcal{NC}_2(k,l)\right)$$

Regarding the spaces in the middle, these are obtained from those on the left by coloring, and we obtain the same spaces as those on the right. Thus, by Tannakian duality, our embedding $\widetilde{O_N^+}\subset U_N^+$ is an isomorphism, modulo the usual equivalence relation.

\medskip

(4) The point now is that the fusion rules for $U_N^+=\widetilde{O}_N^+$ can be computed  by using those of $O_N^+$, from Theorem 3.38, and we end up with a certain ``free complexification'' of the Clebsch-Gordan rules, corresponding to the formulae in the statement.
\end{proof}

Regarding now quantum reflections, let us start with something very basic, namely:

\index{Catalan number}
\index{Clebsch-Gordan rules}
\index{rotation group}

\begin{theorem}
The irreducible representations of $SO_3$ can be labeled by positive integers, $r_k$ with $k\in\mathbb N$, and the fusion rules for these representations are:
$$r_k\otimes r_l=r_{|k-l|}+r_{|k-l|+1}+\ldots+r_{k+l}$$ 
The dimensions of these representations are $\dim r_k=2k+1$.
\end{theorem}

\begin{proof}
As before with $SU_2$, this is something very classical, and there are many possible proofs. Along the lines of our proof of Theorem 3.37, we can argue that, by using the double cover map $SU_2\to SO_3$, we are led to the following formula:
$$\int_{SO_3}\chi^k=C_k$$

But this is what we need for concluding, by reasoning as before, for $SU_2$.
\end{proof}

Quite remarkably, a similar result holds for $S_N^+$, as follows:

\index{free orthogonal group}
\index{Catalan number}
\index{Clebsch-Gordan rules}
\index{Frobenius trick}

\begin{theorem}
The irreducible representations of $S_N^+$ with $N\geq4$ can be labeled by positive integers, $r_k$ with $k\in\mathbb N$, the fusion rules for these representations are
$$r_k\otimes r_l=r_{|k-l|}+r_{|k-l|+1}+\ldots+r_{k+l}$$ 
and the dimensions are $\dim r_k=(q^{k+1}-q^{-k})/(q-1)$, with $q+q^{-1}=N-2$.
\end{theorem}

\begin{proof}
Same story here as for the passage $SU_2\to O_N^+$, the idea being to skilfully recycle the proof of Theorem 3.40. We refer to \cite{ba1} for details on all this.
\end{proof}

Regarding now the general quantum reflections, the result here is as follows:

\begin{theorem}
The irreducible representations of the quantum group $H_N^{s+}$ can be labeled $r_x$ with $x\in<\mathbb Z_s>$, with the fusion rules being as follows,
$$r_x\otimes r_y=\sum_{x=vz,y=\bar{z}w}r_{vw}+r_{v\cdot w}$$
where the fusion product is given by $(i_1\ldots i_k)\cdot (j_1\ldots j_l)=i_1\ldots i_{k-1}(i_k+j_1)j_2\ldots j_l$.
\end{theorem}

\begin{proof}
Again, this is something of Clebsch-Gordan type, that we know from Theorem 3.41 to hold at $s=1$, and whose proof in general is similar. See \cite{bv1}.
\end{proof}

As a philosophical consequence of all this, regarding the main examples of free quantum groups, we have the following result, coming from the above theorems:

\index{irreducible corepresentation}
\index{fusion rules}
\index{Clebsch-Gordan rules}

\begin{theorem}
The irreducible representations of the main free quantum groups, namely the free rotation and reflection groups
$$\xymatrix@R=50pt@C=50pt{
K_N^+\ar[r]&U_N^+\\
H_N^+\ar[u]\ar[r]&O_N^+\ar[u]}$$
can be classified by using easiness, and their fusion rules are given by simple formulae, which are simpler than those for the corresponding classical counterparts.
\end{theorem}

\begin{proof}
This is something rather philosophical, which follows from our various results above, and from the comparison of these with the known results for $O_N,U_N,H_N,K_N$.
\end{proof}

We refer to \cite{ba1} for more details on all this. Let us also mention that it is possible to further build on the above, for instance with the computation of the growth exponents of the above quantum groups, and again we refer here to \cite{ba1} and the literature. 

\bigskip

Finally, getting back to Theorems 3.37 and 3.38, the following question appears:

\begin{question}
Is it possible to unify the representation theory results for $SU_2$ and for $O_N^+$ with $N\geq2$, in an extended easiness framework?
\end{question}

This is an interesting question, but we will leave it for later, towards the end of the present book. We will see there that $SU_2$ is ``super-easy'', coming exactly as $O_N^+$ from the category of pairings $P_2$, but via a different implementation of the $\pi\to T_\pi$ operation.

\section*{3e. Exercises}

We are now into quite specialized quantum group theory, and most of the things that we have seen in this chapter have open problems right around the corner, so our exercises here will be of the same type, rather difficult, of research flavor. We first have:

\begin{exercise}
Axiomatize the easiness-type property of the products of easy quantum groups. Can you do the same with the dual free products of easy quantum groups?
\end{exercise}

As a second exercise, that I don't know how to solve myself, we have:

\begin{exercise}
Try fixing the formula $D_{<G,H>}=D_G\cap D_H$, with a $N>>0$ idea.
\end{exercise}

As a third exercise, which is certainly more reasonable, we have:

\begin{exercise}
Check all the details for the intersection/easy generation property of the standard cube, and for half-liberations too. What about plain generation?
\end{exercise}

If you like combinatorics, here is an exercise for you:

\begin{exercise}
Compute the Gram determinants and the fusion rules for irreducible representations for all the easy quantum groups that we have, so far.
\end{exercise}

And if you like group theory, here is an exercise for you:

\begin{exercise}
Compute the easy envelopes of all compact Lie groups that you know, finite or continuous.
\end{exercise}

Finally, at a more advanced level, the questions mentioned in the middle of this chapter, regarding the maximal tori, are for the most open, and are all very interesting.

\chapter{Probabilistic aspects}

\section*{4a. Laws of characters}

We have seen so far a lot of interesting mathematics, mostly of algebraic type, regarding the easy quantum groups $G\subset U_N^+$, with everything usually coming from the following formula, where $D\subset P$ is the category of partitions associated to $G$:
$$Fix(u^{\otimes k})=span\left(\xi_\pi\Big|\pi\in D(k)\right)$$

Time to do now some analysis, along the same lines. We will do this slowly, as a continuation of our algebraic work. The point indeed is that a quick look at what we did so far in this book, namely algebra, leads to the following conclusion:

\begin{fact}
For most questions regarding $G\subset U_N^+$, the knowledge of the numbers
$$M_k=\dim\left(Fix(u^{\otimes k})\right)$$
is enough. When $G$ is easy, coming from a category of partitions $D\subset P$, we have
$$M_k=\dim\left(span\left(\xi_\pi\Big|\pi\in D(k)\right)\right)$$
and in the $N>>0$ regime, and more specifically when $N\geq k$, we have $M_k=|D(k)|$. 
\end{fact}

Here the first claim comes for instance from our representation theory computations from chapter 3, with the Frobenius tricks there being based on the knowledge of the numbers $M_k$. As for the second claim, this comes from the Lindst\"om determinant formula for $S_N$. The point now is that we have an analytic approach to this, coming from:

\begin{proposition}
Given a closed subgroup $G\subset U_N^+$, we have
$$\int_G\chi^k=\dim\left(Fix(u^{\otimes k})\right)$$
where $\chi=\sum_iu_{ii}$ is the main character. In the easy case we have
$$\lim_{N\to\infty}\int_G\chi^k=|D(k)|$$
where $D\subset P$ is the associated category of partitions.
\end{proposition}

\begin{proof}
Here the first assertion comes from Peter-Weyl theory, which gives:
$$\int_G\chi^k=\int_G\chi_{u^{\otimes k}}=\dim\left(Fix(u^{\otimes k})\right)$$

As for the second assertion, this comes from easiness, as explained in Fact 4.1.
\end{proof}

Now let us look more in detail at Fact 4.1 and Proposition 4.2, taken altogether. These tell us to try to compute the law of the main character $\chi=\sum_iu_{ii}$, because the moments $M_k$ of this law are what we need, in relation with most algebraic questions. So, we are led to the following question, in general, and then in the easy case:

\begin{question}
Given a closed subgroup $G\subset U_N^+$, what is the law of the character
$$\chi=\sum_iu_{ii}$$
with respect to Haar integration? Also, in the easy case, what is the asymptotic law
$$\lim_{N\to\infty}law(\chi)$$
say regarded as a formal measure, having as moments the numbers $|D(k)|$?
\end{question}

All this looks quite interesting, taking us away from the usual algebraic paths, that we have been following so far in this book. But, shall we really get into this, which looks quite scary? This is not an easy question, and so time to ask the cat. And cat says:

\begin{cat}
Advanced algebra exists, and is called analysis.
\end{cat}

Which sounds interesting, so time to upgrade our mathematical knowledge and wisdom, in the hope that one day we'll be able to survive in the wild, like cats do. The simplest easy group that we know is the permutation group $S_N$, and here we have the following beautiful result, adding substantial support for the analytic philosophy:

\index{main character}
\index{derangement}
\index{Poisson law}
\index{random permutation}

\begin{theorem}
Consider the symmetric group $S_N$, regarded as a compact group of matrices, $S_N\subset O_N$, via the standard permutation matrices.
\begin{enumerate}
\item The main character $\chi\in C(S_N)$, defined as usual as $\chi=\sum_iu_{ii}$, counts the number of fixed points, $\chi(\sigma)=\#\{i|\sigma(i)=i\}$.

\item The probability for a permutation $\sigma\in S_N$ to be a derangement, meaning to have no fixed points at all, becomes, with $N\to\infty$, equal to $1/e$.

\item The law of the main character $\chi\in C(S_N)$ becomes with $N\to\infty$ the Poisson law $p_1=\frac{1}{e}\sum_k\delta_k/k!$, with respect to the counting measure.
\end{enumerate}
\end{theorem}

\begin{proof}
This is something very classical, the proof being as follows:

\medskip

(1) We have indeed the following computation, which gives the result:
$$\chi(\sigma)
=\sum_iu_{ii}(\sigma)
=\sum_i\delta_{\sigma(i)i}
=\#\left\{i\Big|\sigma(i)=i\right\}$$

(2) We use the inclusion-exclusion principle. Consider the following sets:
$$S_N^i=\left\{\sigma\in S_N\Big|\sigma(i)=i\right\}$$

The probability that we are interested in is then given by:
\begin{eqnarray*}
P(\chi=0)
&=&\frac{1}{N!}\left(|S_N|-\sum_i|S_N^i|+\sum_{i<j}|S_N^i\cap S_N^j|-\sum_{i<j<k}|S_N^i\cap S_N^j\cap S_N^k|+\ldots\right)\\
&=&\frac{1}{N!}\sum_{r=0}^N(-1)^r\sum_{i_1<\ldots<i_r}(N-r)!\\
&=&\frac{1}{N!}\sum_{r=0}^N(-1)^r\binom{N}{r}(N-r)!\\
&=&\sum_{r=0}^N\frac{(-1)^r}{r!}
\end{eqnarray*}

Since we have here the expansion of $1/e$, this gives the result.

\medskip

(3) This follows by generalizing the computation in (2). Indeed, we get:
$$\lim_{N\to\infty}P(\chi=k)=\frac{1}{k!e}$$

Thus, we obtain in the limit a Poisson law of parameter 1, as stated.
\end{proof}

The above result is very beautiful, but we can in fact do even better. Indeed, you might remember from probability that the Poisson law $p_1$ is part of a 1-parameter family $\{p_t\}_{t\geq0}$. And, good news, we have the following generalization of Theorem 4.5:

\index{truncated character}
\index{derangement}
\index{Poisson law}
\index{random permutation}

\begin{theorem}
For the symmetric group $S_N\subset O_N$, the truncated character
$$\chi_t=\sum_{i=1}^{[tN]}u_{ii}$$
with $t\in[0,1]$ follows with $N\to\infty$ the Poisson law $p_t=e^{-t}\sum_k\delta_kt^k/k!$.
\end{theorem}

\begin{proof}
This follows by suitably modifying the proof of Theorem 4.5. Indeed, according to our definition of $\chi_t$, we first have the following formula:
$$\chi(\sigma)=\#\left\{i\in\{1,\ldots,[tN]\}\Big|\sigma(i)=i\right\}$$

Then, the inclusion-exclusion principle gives the following formula:
$$\lim_{N\to\infty}P(\chi_t=k)=\frac{t^k}{k!e^t}$$

Thus, we obtain in the limit a Poisson law of parameter $t$, as stated.
\end{proof}

Summarizing, we have some interesting theory going on here, and this is probably just the tip of the iceberg. So, let us update Question 4.3, as follows: 

\begin{question}
Given $G\subset U_N^+$, what is the law of the truncated character
$$\chi_t=\sum_{i=1}^{[tN]}u_{ii}$$
with respect to Haar integration? Also, in the easy case, what is the asymptotic law
$$\lim_{N\to\infty}law(\chi_t)$$
say regarded as a formal measure, determined by its moments?
\end{question}

Our goal in what follows will be that of answering this question, by extending what we have about $S_N$ to all the easy groups and easy quantum groups that we know.

\section*{4b. Weingarten formula}

In order to reach to our goals, we first need to know how to integrate on the easy quantum groups. The formula here, coming from Peter-Weyl theory, is as follows:

\index{Gram matrix}
\index{Weingarten matrix}
\index{Weingarten formula}
\index{Haar functional}

\begin{theorem}
Assuming that a closed subgroup $G\subset U_N^+$ is easy, coming from a category of partitions $D\subset P$, we have the Weingarten formula
$$\int_Gu_{i_1j_1}^{e_1}\ldots u_{i_kj_k}^{e_k}=\sum_{\pi,\sigma\in D(k)}\delta_\pi(i)\delta_\sigma(j)W_{kN}(\pi,\sigma)$$
where $\delta\in\{0,1\}$ are the usual Kronecker type symbols, and where the Weingarten matrix $W_{kN}=G_{kN}^{-1}$ is the inverse of the Gram matrix $G_{kN}(\pi,\sigma)=N^{|\pi\vee\sigma|}$. 
\end{theorem}

\begin{proof}
We know from chapter 2 that the integrals in the statement form altogether the orthogonal projection $P^k$ onto the following space:
$$Fix(u^{\otimes k})=span\left(\xi_\pi\Big|\pi\in D(k)\right)$$

In order to prove the result, consider the following linear map:
$$E(x)=\sum_{\pi\in D(k)}<x,\xi_\pi>\xi_\pi$$

By a standard linear algebra computation, it follows that we have $P=WE$, where $W$ is the inverse on $Fix(u^{\otimes k})$ of the restriction of $E$. But this restriction is the linear map given by $G_{kN}$, and so $W$ is the linear map given by $W_{kN}$, and this gives the result.
\end{proof}

Observe that there is a bit of confusion in Theorem 4.8, because the partitions in $D(k)$ are known to be linearly independent only at $N\geq k$, and so their Gram matrix is invertible only when adding such a $N\geq k$ assumption. In what follows we will only use Theorem 4.8 in the $N>>0$ regime, where the above formulation will do. In general, the point is that $W_{kN}$ should be taken to be the quasi-inverse of $G_{kN}$. See \cite{csn}.

\bigskip

In relation now with truncated characters, we have the following formula:

\index{truncated character}

\begin{proposition}
The moments of truncated characters are given by the formula
$$\int_G(u_{11}+\ldots +u_{ss})^k=Tr(W_{kN}G_{ks})$$
where $G_{kN}$ and $W_{kN}=G_{kN}^{-1}$ are the associated Gram and Weingarten matrices.
\end{proposition}

\begin{proof}
We have indeed the following computation:
\begin{eqnarray*}
\int_G(u_{11}+\ldots +u_{ss})^k
&=&\sum_{i_1=1}^{s}\ldots\sum_{i_k=1}^s\int u_{i_1i_1}\ldots u_{i_ki_k}\\
&=&\sum_{\pi,\sigma\in D(k)}W_{kN}(\pi,\sigma)\sum_{i_1=1}^{s}\ldots\sum_{i_k=1}^s\delta_\pi(i)\delta_\sigma(i)\\
&=&\sum_{\pi,\sigma\in D(k)}W_{kN}(\pi,\sigma)G_{ks}(\sigma,\pi)\\
&=&Tr(W_{kN}G_{ks})
\end{eqnarray*}

Thus, we have obtained the formula in the statement.
\end{proof}

All this is very good, and normally we have here what is needed in order to answer Question 4.7, in the easy case. In practice, however, before doing so, we will need some training in probability theory. Imagine for instance that we are trying to solve the simplest question around, namely computing the asymptotic law of $\chi_1$ for the simplest easy group, namely $S_N$. Here not even need for the Weingarten formula and Proposition 4.9, because we know from easiness that the moments in question are the Bell numbers:
$$B_k=|P(k)|$$

But how on Earth, without knowing any advanced probability, can we find the measure having as moments the Bell numbers? Mission impossible, hope you agree with me.

\bigskip

So, let us do this, learn some probability, and we will come back to quantum groups later. Obviously, we are in need of ``quantum probability'', so let us start with:

\index{random variable}
\index{law of variable}
\index{moments}
\index{colored moments}

\begin{definition}
Let $A$ be a $C^*$-algebra, given with a trace $tr:A\to\mathbb C$.
\begin{enumerate}
\item The elements $a\in A$ are called random variables.

\item The moments of such a variable are the numbers $M_k(a)=tr(a^k)$.

\item The law of such a variable is the functional $\mu:P\to tr(P(a))$.
\end{enumerate}
\end{definition}

Here $k=\circ\bullet\bullet\circ\ldots$ is by definition a colored integer, and the corresponding powers $a^k$ are defined by the following formulae, and multiplicativity: 
$$a^\emptyset=1\quad,\quad
a^\circ=a\quad,\quad
a^\bullet=a^*$$

As for the polynomial $P$, this is a noncommuting $*$-polynomial in one variable:
$$P\in\mathbb C<X,X^*>$$

Observe that the law is uniquely determined by the moments, because we have:
$$P(X)=\sum_k\lambda_kX^k\implies\mu(P)=\sum_k\lambda_kM_k(a)$$

Generally speaking, the above definition is something quite abstract, but there is no other way of doing things, at least at this level of generality. However, in certain special cases, the formalism simplifies, and we recover more familiar objects, as follows:

\index{spectral theorem}
\index{spectral measure}

\begin{proposition}
Assuming that $a\in A$ is normal, $aa^*=a^*a$, its law corresponds to a probability measure on its spectrum $\sigma(a)\subset\mathbb C$, according to the following formula:
$$tr(P(a))=\int_{\sigma(a)}P(x)d\mu(x)$$
When the trace is faithful we have $supp(\mu)=\sigma(a)$. Also, in the particular case where the variable is self-adjoint, $a=a^*$, this law is a real probability measure.
\end{proposition}

\begin{proof}
Since the $C^*$-algebra $<a>$ generated by $a$ is commutative, the Gelfand theorem applies to it, and gives an identification of $C^*$-algebras, as follows:
$$<a>=C(\sigma(a))$$

Now by using the Riesz theorem, the restriction of $tr$ to this algebra must come from a probability measure $\mu$ as in the statement, and this gives all the assertions.
\end{proof}

Getting now where we wanted to get, we want our ``quantum probability'' theory to apply to the two main cases that we have in mind, namely classical and free. So, following Voiculescu \cite{vo1}, let us introduce the following two notions of independence:

\index{independence}
\index{freeness}

\begin{definition}
Two subalgebras $A,B\subset C$ are called independent when
$$tr(a)=tr(b)=0\implies tr(ab)=0$$
holds for any $a\in A$ and $b\in B$, and free when
$$tr(a_i)=tr(b_i)=0\implies tr(a_1b_1a_2b_2\ldots)=0$$
holds for any $a_i\in A$ and $b_i\in B$.
\end{definition}
 
In other words, what we have here is a straightforward noncommutative extension of the usual notion of independence, along with a natural free analogue of it. In order to understand now what is going on, let us discuss some basic models for independence and freeness. We have the following result, from \cite{vo1}, which clarifies things:

\begin{proposition}
Given two algebras $(A,tr)$ and $(B,tr)$, the following hold:
\begin{enumerate}
\item $A,B$ are independent inside their tensor product $A\otimes B$.

\item $A,B$ are free inside their free product $A*B$.
\end{enumerate}
\end{proposition}

\begin{proof}
Both the assertions are clear from definitions, after some standard discussion regarding the tensor product and free product trace. See Voiculescu \cite{vo1}.
\end{proof}

In relation with groups and algebra, we have the following result:

\begin{proposition}
We have the following results, valid for group algebras:
\begin{enumerate}
\item $C^*(\Gamma),C^*(\Lambda)$ are independent inside $C^*(\Gamma\times\Lambda)$.

\item $C^*(\Gamma),C^*(\Lambda)$ are free inside $C^*(\Gamma*\Lambda)$.
\end{enumerate}
\end{proposition}

\begin{proof}
This follows from the general results in Proposition 4.13, along with the following two isomorphisms, which are both standard:
$$C^*(\Gamma\times\Lambda)=C^*(\Lambda)\otimes C^*(\Gamma)\quad,\quad 
C^*(\Gamma*\Lambda)=C^*(\Lambda)*C^*(\Gamma)$$

Alternatively, we can prove this directly, by using the fact that each algebra is spanned by the corresponding group elements, and checking the result on group elements.
\end{proof}

In order to study independence and freeness, our main tool will be: 

\index{Fourier transform}
\index{Cauchy transform}
\index{R-transform}
\index{convolution}
\index{free convolution}

\begin{theorem}
The convolution is linearized by the log of the Fourier transform,
$$F_f(x)=\mathbb E(e^{ixf})$$
and the free convolution is linearized by the $R$-transform, which is given by
$$G_\mu(\xi)=\int_\mathbb R\frac{d\mu(t)}{\xi-t}\implies G_\mu\left(R_\mu(\xi)+\frac{1}{\xi}\right)=\xi$$
and so is the inverse of the Cauchy transform, up to a $\xi^{-1}$ factor.
\end{theorem}

\begin{proof}
For the first assertion, if $f,g$ are independent, we have indeed:
\begin{eqnarray*}
F_{f+g}(x)
&=&\int_\mathbb Re^{ixz}d(\mu_f*\mu_g)(z)\\
&=&\int_{\mathbb R\times\mathbb R}e^{ix(z+t)}d\mu_f(z)d\mu_g(t)\\
&=&F_f(x)F_g(x)
\end{eqnarray*}

For the second assertion, we need a good model for free convolution, and the best is to use the semigroup algebra of the free semigroup on two generators:
$$A=C^*(\mathbb N*\mathbb N)$$

Indeed, we have some freeness in the semigroup setting, a bit in the same way as for the group algebras $C^*(\Gamma*\Lambda)$, from Proposition 4.14. In addition to this fact, and to what happens in the group algebra case, the following things happen:

\medskip

(1) The variables of type $S^*+f(S)$, with $S\in C^*(\mathbb N)$ being the shift, and with $f\in\mathbb C[X]$ being a polynomial, model in moments all the distributions $\mu:\mathbb C[X]\to\mathbb C$. This is indeed something elementary, which can be checked via a direct algebraic computation.

\medskip

(2) Given $f,g\in\mathbb C[X]$, the variables $S^*+f(S)$ and $T^*+g(T)$, where $S,T\in C^*(\mathbb N*\mathbb N)$ are the shifts corresponding to the generators of $\mathbb N*\mathbb N$, are free, and their sum has the same law as $S^*+(f+g)(S)$. This follows indeed by using a $45^\circ$ argument.

\medskip

(3) But with this in hand, we can see that $\mu\to f$ linearizes the free convolution. We are therefore left with a computation inside $C^*(\mathbb N)$, whose conclusion is that $R_\mu=f$ can be recaptured from $\mu$ via the Cauchy transform $G_\mu$, as stated. See \cite{vo1}.
\end{proof}

As a first main result now, we have the following statement, called Central Limit Theorem (CLT), which in the free case is due to Voiculescu \cite{vo1}:

\index{CLT}
\index{Wigner law}
\index{semicircle law}
\index{normal law}
\index{Gaussian law}

\begin{theorem}[CLT]
Given self-adjoint variables $x_1,x_2,x_3,\ldots$ which are i.i.d./f.i.d., centered, with variance $t>0$, we have, with $n\to\infty$, in moments,
$$\frac{1}{\sqrt{n}}\sum_{i=1}^nx_i\sim g_t/\gamma_t$$
where the limiting laws $g_t/\gamma_t$ are the following measures,
$$g_t=\frac{1}{\sqrt{2\pi t}}e^{-x^2/2t}dx\quad,\quad 
\gamma_t=\frac{1}{2\pi t}\sqrt{4t^2-x^2}dx$$
called normal, or Gaussian, and Wigner semicircle law of parameter $t$.
\end{theorem}

\begin{proof}
This is routine, by using the linearization properties of the Fourier transform and the $R$-transform from Theorem 4.15, and for details here, we refer to any classical probability book for the classical result, and to \cite{vdn} for the free result.
\end{proof}

Next, we have the following complex version of the CLT:

\index{CCLT}
\index{Complex CLT}
\index{Voiculescu law}
\index{circular law law}
\index{complex normal law}
\index{complex Gaussian law}

\begin{theorem}[CCLT]
Given variables $x_1,x_2,x_3,\ldots$ which are i.i.d./f.i.d., centered, with variance $t>0$, we have, with $n\to\infty$, in moments,
$$\frac{1}{\sqrt{n}}\sum_{i=1}^nx_i\sim G_t/\Gamma_t$$
where $G_t/\Gamma_t$ are the complex normal and Voiculescu circular law of parameter $t$, given by:
$$G_t=law\left(\frac{1}{\sqrt{2}}(a+ib)\right)\quad,\quad 
\Gamma_t=law\left(\frac{1}{\sqrt{2}}(\alpha+i\beta)\right)$$
where $a,b/\alpha,\beta$ are independent/free, each following the law $g_t/\gamma_t$.
\end{theorem}

\begin{proof}
This follows indeed from the CLT, by taking real and imaginary parts of all the variables involved, and for details and more here, including the combinatorics of the Voiculescu circular law $\Gamma_t$, which is quite subtle, we refer again to \cite{vdn}.
\end{proof}

We denote by $\boxplus$ the free convolution of real probability measures, given by the fact that $\mu_{a+b}=\mu_a\boxplus\mu_b$, when $a,b$ are free. With this convention, we have the following discrete version of the CLT, called Poisson Limit Theorem (PLT):

\index{PLT}
\index{Poisson limit}
\index{Poisson law}
\index{free Poisson law}
\index{Marchenko-Pastur law}

\begin{theorem}[PLT]
The following Poisson limits converge, for any $t>0$,
$$p_t=\lim_{n\to\infty}\left(\left(1-\frac{t}{n}\right)\delta_0+\frac{t}{n}\delta_1\right)^{*n}\quad,\quad 
\pi_t=\lim_{n\to\infty}\left(\left(1-\frac{t}{n}\right)\delta_0+\frac{t}{n}\delta_1\right)^{\boxplus n}$$
the limiting measures being the Poisson law $p_t$, and the Marchenko-Pastur law $\pi_t$,
$$p_t=\frac{1}{e^t}\sum_{k=0}^\infty\frac{t^k\delta_k}{k!}\quad,\quad 
\pi_t=\max(1-t,0)\delta_0+\frac{\sqrt{4t-(x-1-t)^2}}{2\pi x}\,dx$$
with at $t=1$, the Marchenko-Pastur law being $\pi_1=\frac{1}{2\pi}\sqrt{4x^{-1}-1}\,dx$. 
\end{theorem}

\begin{proof}
This is again routine, by using the Fourier and $R$-transform, and as before we refer here to any classical probability book, and to \cite{vdn}.
\end{proof}

Finally, we have the following ``compound'' generalization of the PLT:

\index{CPLT}
\index{Compound PLT}
\index{compound Poisson law}

\begin{theorem}[CPLT]
Given a compactly supported positive measure $\rho$, of mass $c=mass(\rho)$, the following compound Poisson limits converge,
$$p_\rho=\lim_{n\to\infty}\left(\left(1-\frac{c}{n}\right)\delta_0+\frac{1}{n}\rho\right)^{*n}
\quad,\quad 
\pi_\rho=\lim_{n\to\infty}\left(\left(1-\frac{c}{n}\right)\delta_0+\frac{1}{n}\rho\right)^{\boxplus n}$$
and if we write $\rho=\sum_{i=1}^sc_i\delta_{z_i}$ with $c_i>0$ and $z_i\in\mathbb R$ we have the formula
$$p_\rho/\pi_\rho={\rm law}\left(\sum_{i=1}^sz_i\alpha_i\right)$$
where the variables $\alpha_i$ are Poisson/free Poisson$(c_i)$, independent/free.
\end{theorem}

\begin{proof}
As before, this follows by using the Fourier and the $R$-transform, and details can be found in any probability book, and in \cite{vdn}.
\end{proof}

So long for limiting results in classical and free probability. To finish with, for our purposes here, we will need the following notions, coming from Theorem 4.19:

\index{Bessel law}
\index{free Bessel law}

\begin{definition}
The Bessel and free Bessel laws, depending on parameters $s\in\mathbb N\cup\{\infty\}$ and $t>0$, are the following compound Poisson and free Poisson laws,
$$b^s_t=p_{t\varepsilon_s}\quad,\quad 
\beta^s_t=\pi_{t\varepsilon_s}$$
with $\varepsilon_s$ being the uniform measure on the $s$-th roots of unity. In particular:
\begin{enumerate}
\item At $s=1$ we recover the Poisson laws $p_t,\pi_t$.

\item At $s=2$ we have the real Bessel laws $b_t,\beta_t$.

\item At $s=\infty$ we have the complex Bessel laws $B_t,\mathfrak B_t$.
\end{enumerate}
\end{definition}

Here the terminology comes from the fact that the density of the measure $b_t$ from (2) is a Bessel function of the first kind, the formula, from \cite{bb+}, \cite{bbc}, being as follows:
$$b_t=e^{-t}\sum_{r=-\infty}^\infty\delta_r\sum_{p=0}^\infty \frac{(t/2)^{|r|+2p}}{(|r|+p)!p!}$$

Good news, with the above general theory in hand, we can now formulate our truncated character results for the main examples of easy quantum groups, as follows:

\index{standard cube}
\index{truncated character}

\begin{theorem}
For the main quantum rotation and reflection groups,
$$
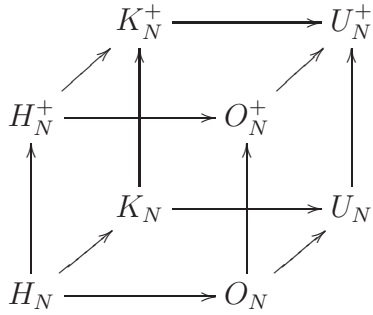
\xymatrix@R=16pt@C=16pt{
&K_N^+\ar[rr]&&U_N^+\\
H_N^+\ar[rr]\ar[ur]&&O_N^+\ar[ur]\\
&K_N\ar[rr]\ar[uu]&&U_N\ar[uu]\\
H_N\ar[uu]\ar[ur]\ar[rr]&&O_N\ar[uu]\ar[ur]
}$$
the corresponding truncated characters follow with $N\to\infty$ the laws
$$
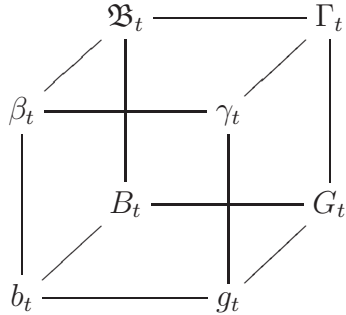
\xymatrix@R=18pt@C=20pt{
&\mathfrak B_t\ar@{-}[rr]\ar@{-}[dd]&&\Gamma_t\ar@{-}[dd]\\
\beta_t\ar@{-}[rr]\ar@{-}[dd]\ar@{-}[ur]&&\gamma_t\ar@{-}[dd]\ar@{-}[ur]\\
&B_t\ar@{-}[rr]\ar@{-}[uu]&&G_t\ar@{.}[uu]\\
b_t\ar@{-}[uu]\ar@{-}[ur]\ar@{-}[rr]&&g_t\ar@{-}[uu]\ar@{-}[ur]
}$$
which are the main limiting laws in classical and free probability.
\end{theorem}

\begin{proof}
We know from chapters 1-2 that the above quantum groups are all easy, coming from the following categories of partitions:
$$
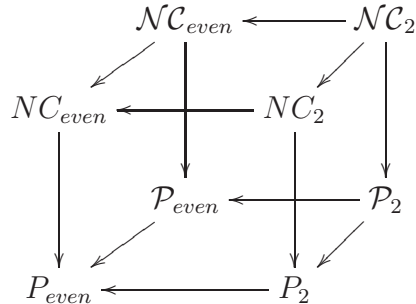
\xymatrix@R=17pt@C4pt{
&\mathcal{NC}_{even}\ar[dl]\ar[dd]&&\mathcal {NC}_2\ar[dl]\ar[ll]\ar[dd]\\
NC_{even}\ar[dd]&&NC_2\ar[dd]\ar[ll]\\
&\mathcal P_{even}\ar[dl]&&\mathcal P_2\ar[dl]\ar[ll]\\
P_{even}&&P_2\ar[ll]
}$$

(1) At $t=1$, we can use the following general formula, from Proposition 4.2:
$$\lim_{N\to\infty}\int_{G_N}\chi^k=|D(k)|$$

But this gives the laws in the statement, via some standard calculus. 

\medskip

(2) In order to compute now the asymptotic laws of truncated characters, at any $t>0$, we can use the general moment formula in Proposition 4.9, namely:
$$\int_G(u_{11}+\ldots +u_{ss})^k=Tr(W_{kN}G_{ks})$$

To be more precise, what happens is that in each of the cases under consideration, the Gram matrix is asymptotically diagonal, and so the Weingarten matrix is asymptotically diagonal too. Thus, in the limit we obtain the following moment formula:
$$\lim_{N\to\infty}\int_{G_N}\chi_t^k=\sum_{\pi\in D(k)}t^{|\pi|}$$

But this gives the laws in the statement, via some standard calculus.
\end{proof}

Summarizing, we have solved the questions raised in the beginning of this chapter. All this was of course quite quick, assuming a bit of familiarity with probability theory. But we will be back to this on several occasions, and notably, in what comes next.

\section*{4c. Cumulant theory}

What we have so far, analytically, is quite exciting, but remains quite wizarding, with the statement of Theorem 4.21 being something quite advanced, and with the proof being something quite advanced too, based on way too many things. So, time to further clarify all this. As a basic question that we would like to solve, we have:

\begin{question}
Given a liberation of easy groups, $G_N\to G_N^+$, can we say, in a simple way, that at the level of asymptotic laws of truncated characters
$$\chi_t=\sum_{i=1}^{[tN]}u_{ii}$$ 
the laws for $G_N^+$ appear as ``liberations'' of the laws for $G_N$?
\end{question}

As a first idea in order to answer the above question, let us focus on the asymptotic moment formula from the end of the proof of Theorem 4.21, namely:
$$\lim_{N\to\infty}\int_{G_N}\chi_t^k=\sum_{\pi\in D(k)}t^{|\pi|}$$

This is something very nice, purely combinatorial, and we are led in this way to:

\begin{answer}
Yes, we can say that according to the Weingarten formula we have
$$\lim_{N\to\infty}\int_{G_N}\chi_t^k=\sum_{\pi\in D(k)}t^{|\pi|}$$
and since $G_N\to G_N^+$ corresponds to $D\to D\cap NC$ at the level of partitions, we have
$$\sum_{\pi\in D(k)}t^{|\pi|}\to\sum_{\pi\in (D\cap NC)(k)}t^{|\pi|}$$
at the level of moments, and we can declare this to mean ``liberation'', probabilistically.
\end{answer}

This is not bad as an answer, because gone in this way all the limiting theorems, classical and free, no need anymore to learn all that material. However, as an obvious downside, what we say at the end does not look very clean, and suggests:

\begin{question}
Sure yes, but how do probability measures liberate in general, with no reference to moment formulae as above, which look very special?
\end{question}

To summarize now, we have an idea, but we must do some more probability, classical and free, in relation with combinatorics and partitions, before formulating our idea. Let us start with the classical case. We have here the following well-known definition:

\index{Fourier transform}
\index{cumulant-generating function}
\index{cumulant}

\begin{definition}
Associated to any real probability measure $\mu=\mu_f$ is the following modification of the logarithm of the Fourier transform $F_\mu(\xi)=\mathbb E(e^{i\xi f})$,
$$K_\mu(\xi)=\log\mathbb E(e^{\xi f})$$
called cumulant-generating function. The Taylor coefficients $k_n(\mu)$ of this series, given by
$$K_\mu(\xi)=\sum_{n=1}^\infty k_n(\mu)\,\frac{\xi^n}{n!}$$
are called cumulants of the measure $\mu$. We also use the notations $k_f,K_f$ for these cumulants and their generating series, where $f$ is a variable following the law $\mu$.
\end{definition}

In other words, the cumulants are more or less the coefficients of the logarithm of the Fourier transform $\log F_\mu$, up to some normalizations. To be more precise, we have $K_\mu(\xi)=\log F_\mu(-i\xi)$, so the formula relating $\log F_\mu$ to the cumulants $k_n(\mu)$ is:
$$\log F_\mu(-i\xi)=\sum_{n=1}^\infty k_n(\mu)\,\frac{\xi^n}{n!}$$

Equivalently, the formula relating $\log F_\mu$ to the cumulants $k_n(\mu)$ is:
$$\log F_\mu(\xi)=\sum_{n=1}^\infty k_n(\mu)\,\frac{(i\xi)^n}{n!}$$

The interest in the cumulants comes from the following result:

\begin{theorem}
The cumulants have the following properties:
\begin{enumerate}
\item $k_n(cf)=c^nk_n(f)$.

\item $k_1(f+d)=k_1(f)+d$, and $k_n(f+d)=k_n(f)$ for $n>1$.

\item $k_n(f+g)=k_n(f)+k_n(g)$, if $f,g$ are independent.
\end{enumerate}
\end{theorem}

\begin{proof}
Here (1) and (2) are both clear from definitions, because we have the following computation, valid for any $c,d\in\mathbb R$, which gives the results:
\begin{eqnarray*}
K_{cf+d}(\xi)
&=&\log\mathbb E(e^{\xi(cf+d)})\\
&=&\log[e^{\xi d}\cdot\mathbb E(e^{\xi cf})]\\
&=&\xi d+K_f(c\xi)
\end{eqnarray*}

As for (3), this follows from the fact that the Fourier transform $F_f(\xi)=\mathbb E(e^{i\xi f})$ satisfies the following formula, whenever $f,g$ are independent random variables:
$$F_{f+g}(\xi)=F_f(\xi)F_g(\xi)$$

Indeed, by applying the logarithm, we obtain the following formula:
$$\log F_{f+g}(\xi)=\log F_f(\xi)+\log F_g(\xi)$$

With the change of variables $\xi\to-i\xi$, we obtain the following formula:
$$K_{f+g}(\xi)=K_f(\xi)+K_g(\xi)$$

Thus, at the level of coefficients, we obtain $k_n(f+g)=k_n(f)+k_n(g)$, as claimed.
\end{proof}

At the level of examples, we have the following result:

\begin{proposition}
The sequence of cumulants $k_1,k_2,k_3,\ldots$ is as follows:
\begin{enumerate}
\item For $\mu=\delta_c$ the cumulants are $c,0,0,\ldots$

\item For $\mu=g_t$ the cumulants are $0,t,0,0,\ldots$

\item For $\mu=p_t$ the cumulants are $t,t,t,\ldots$

\item For $\mu=b_t$ the cumulants are $0,t,0,t,\ldots$
\end{enumerate}
Also, for the compound Poisson laws the cumulants are $k_n(p_\nu)=M_n(\nu)$.
\end{proposition}

\begin{proof}
We have 5 computations to be done, the idea being as follows:

\medskip

(1) For $\mu=\delta_c$ we have the following computation:
$$K_\mu(\xi)
=\log\mathbb E(e^{c\xi})
=\log(e^{c\xi})
=c\xi$$

But the plain coefficients of this series are the numbers $c,0,0,\ldots\,$, and so the Taylor coefficients of this series are these same numbers $c,0,0,\ldots\,$, as claimed.

\medskip

(2) For $\mu=g_t$ we have the following computation:
\begin{eqnarray*}
K_\mu(\xi)
&=&\log F_\mu(-i\xi)\\
&=&\log\exp\left[-t(-i\xi)^2/2\right]\\
&=&t\xi^2/2
\end{eqnarray*}

But the plain coefficients of this series are the numbers $0,t/2,0,0,\ldots\,$, and so the Taylor coefficients of this series are the numbers $0,t,0,0,\ldots\,$, as claimed.

\medskip

(3) For $\mu=p_t$ we have the following computation:
\begin{eqnarray*}
K_\mu(\xi)
&=&\log F_\mu(-i\xi)\\
&=&\log\exp\left[(e^{i(-i\xi)}-1)t\right]\\
&=&(e^\xi-1)t
\end{eqnarray*}

But the plain coefficients of this series are the numbers $t/n!$, and so the Taylor coefficients of this series are the numbers $t,t,t,\ldots\,$, as claimed.

\medskip

(4) For $\mu=b_t$ we have the following computation:
\begin{eqnarray*}
K_\mu(\xi)
&=&\log F_\mu(-i\xi)\\
&=&\log\exp\left[\left(\frac{e^\xi+e^{-\xi}}{2}-1\right)t\right]\\
&=&\left(\frac{e^\xi+e^{-\xi}}{2}-1\right)t
\end{eqnarray*}

But the plain coefficients of this series are the numbers $(1+(-1)^n)t/n!$, so the Taylor coefficients of this series are the numbers $0,t,0,t,\ldots\,$, as claimed.

\medskip

(5) We can assume, by using a continuity argument, that our measure $\nu$ is discrete, as follows, with $t_i>0$ and $z_i\in\mathbb R$, and with the sum being finite:
$$\nu=\sum_i t_i\delta_{z_i}$$

By using now the well-known Fourier transform formula for $p_\nu$, we obtain:
\begin{eqnarray*}
K_{p_\nu}(\xi)
&=&\log F_{p_\nu}(-i\xi)\\
&=&\log\exp\left[\sum_it_i(e^{\xi z_i}-1)\right]\\
&=&\sum_it_i\sum_{n\geq1}\frac{(\xi z_i)^n}{n!}\\
&=&\sum_{n\geq1}\frac{\xi^n}{n!}\sum_it_iz_i^n\\
&=&\sum_{n\geq1}\frac{\xi^n}{n!}\,M_n(\nu)
\end{eqnarray*}

Thus, we are led to the conclusion in the statement.
\end{proof}

All the above is quite nice and elementary, and we can see here emerging some sort of relationship with what we say in Answer 4.23. But in order to really get to that, we still have to introduce partitions, in our discussion. So, let us introduce:

\index{cumulant}
\index{classical cumulant}

\begin{definition}
We define quantities $M_\pi(f),k_\pi(f)$, depending on partitions 
$$\pi\in P(k)$$
by starting with $M_n(f),k_n(f)$, and using multiplicativity over the blocks. 
\end{definition}

To be more precise, the convention here is that for the one-block partition $1_n\in P(n)$, the corresponding moment and cumulant are the usual ones, namely:
$$M_{1_n}(f)=M_n(f)\quad,\quad k_{1_n}(f)=k_n(f)$$

Then, for an arbitrary partition $\pi\in P(k)$, we decompose this partition into blocks, having sizes $b_1,\ldots,b_s$, and we set, by multiplicativity over blocks:
$$M_\pi(f)=M_{b_1}(f)\ldots M_{b_s}(f)\quad,\quad k_\pi(f)=k_{b_1}(f)\ldots k_{b_s}(f)$$

With this convention, following Rota and others, we can now formulate a key result, fully clarifying the relation between moments and cumulants, as follows:

\index{moment-cumulant formula}

\begin{theorem}
We have the moment-cumulant formulae
$$M_n(f)=\sum_{\nu\in P(n)}k_\nu(f)\quad,\quad 
k_n(f)=\sum_{\nu\in P(n)}\mu(\nu,1_n)M_\nu(f)$$
or, equivalently, we have the moment-cumulant formulae
$$M_\pi(f)=\sum_{\nu\leq\pi}k_\nu(f)\quad,\quad 
k_\pi(f)=\sum_{\nu\leq\pi}\mu(\nu,\pi)M_\nu(f)$$
where $\mu$ is the M\"obius function of $P(n)$.
\end{theorem}

\begin{proof}
By using the M\"obius inversion formula explained in chapter 3, the four formulae in the statement are equivalent, so it is enough to prove the first one, namely:
$$M_n(f)=\sum_{\nu\in P(n)}k_\nu(f)$$

In order to do this, we can use the very definition of the cumulants, namely:
$$\log\mathbb E(e^{\xi f})=\sum_{s=1}^\infty k_s(f)\,\frac{\xi^s}{s!}$$

By exponentiating, we obtain from this the following formula:
$$\mathbb E(e^{\xi f})=\exp\left(\sum_{s=1}^\infty k_s(f)\,\frac{\xi^s}{s!}\right)$$

But this leads, via some standard calculus, to the above formula for $M_n(f)$.
\end{proof}

Time now to go back to groups, and answer our questions. We will restrict the attention to the orthogonal easy groups that we know, namely $S_N,H_N,O_N$. We will need the following observation regarding these groups, from \cite{bsp}:

\begin{proposition}
The orthogonal easy groups $G\subset O_N$ that we know, namely
$$S_N\subset H_N\subset O_N$$
come respectively from the following categories of partitions $D\subset P$,
$$P\supset P_{even}\supset P_2$$
which in turn come from certain subsets $L\subset\mathbb N$, as follows,
$$\mathbb N\supset 2\mathbb N\supset\{2\}$$
with $D$ consisting of the partitions $\pi\in P$ whose blocks have lengths belonging to $L\subset\mathbb N$. 
\end{proposition}

\begin{proof}
This is something quite trivial, as stated, but which is conceptual as well, having behind it some interesting mathematics, the idea being as follows:

\medskip

(1) We already know that the first assertion, regarding the easiness correspondence $G\leftrightarrow D$, holds for our 3 groups. As for the second assertion, regarding the correspondence $D\leftrightarrow L$, for our 3 categories of partitions, this is clear from definitions. Thus, everything proved, and with the result itself looking like something quite anecdotical. 

\medskip

(2) However, we will see later, in chapter 6 below, that the easy groups $G$ having the property that their categories of partitions $D$ are stable under removing blocks, and so appear from sets $L\subset\mathbb N$ are in the statement, which are called ``uniform'', are very interesting objects, and that a lot of general theory can be developed for them.
\end{proof}

In relation now with cumulants, we have the following result, also from \cite{bsp}:

\index{truncated character}
\index{easy group}
\index{asymptotic moment}

\begin{theorem}
The cumulants of the asymptotic truncated characters for the groups $S_N,H_N,O_N$ are given by the formula
$$k_n(\chi_t)=t\delta_{n\in L}$$
with $L\subset\mathbb N$ being the associated subset, and at the level of asymptotic moments this gives 
$$M_k(\chi_t)=\sum_{\pi\in D(k)}t^{|\pi|}$$
with $D\subset P$ being the associated category of partitions.
\end{theorem}

\begin{proof}
This is clear indeed from Proposition 4.30, by performing a case-by-case analysis, with the cases $G=S,H,O$ under consideration corresponding to the cumulant computations for the measures $p_t,b_t,g_t$ from Proposition 4.27. See \cite{bsp}. 
\end{proof}

All this is very nice, answering our philosophical questions in the classical case. Now that we killed the problem in the classical case, let us do the same in the free case. Following Speicher \cite{spe} and subsequent work, we have the following definition:

\index{R-transform}
\index{free cumulant}

\begin{definition}
The free cumulants $\kappa_n(a)$ of a variable $a\in A$ are defined by:
$$R_a(\xi)=\sum_{n=1}^\infty\kappa_n(a)\xi^{n-1}$$
That is, the free cumulants are the coefficients of the $R$-transform.
\end{definition}

As before in the classical case, there are many interesting things that can be said about free cumulants. At the level of basic general results, we first have:

\begin{theorem}
The free cumulants have the following properties:
\begin{enumerate}
\item $\kappa_n(\lambda a)=\lambda^n\kappa_n(a)$.

\item $\kappa_n(a+b)=\kappa_n(a)+\kappa_n(b)$, if $a,b$ are free.
\end{enumerate}
\end{theorem}

\begin{proof}
In what regards (1), we have here the following computation:
\begin{eqnarray*}
G_{\lambda a}(\xi)
&=&\int_\mathbb R\frac{d\mu_{\lambda a}(t)}{\xi-t}\\
&=&\frac{1}{\lambda}\int_\mathbb R\frac{d\mu_a(s)}{\xi/\lambda-s}\\
&=&\frac{1}{\lambda}\,G_a\left(\frac{\xi}{\lambda}\right)
\end{eqnarray*}

But, according to the definition of the $R$-transform, this gives the following formula:
$$G_{\lambda a}\left(\lambda R_a(\lambda\xi)+\frac{1}{\xi}\right)
=\frac{1}{\lambda}\,G_a\left(R_a(\lambda\xi)+\frac{1}{\lambda\xi}\right)
=\xi$$

Thus $R_{\lambda a}(\xi)=\lambda R_a(\lambda\xi)$, which gives (1). As for (2), this follows from the fact, that we know well, that the $R$-transform linearizes the free convolution operation. 
\end{proof}

Again in analogy with the classical case, at the level of examples, we have:

\begin{proposition}
The sequence of free cumulants $\kappa_1,\kappa_2,\kappa_3,\ldots$ is as follows:
\begin{enumerate}
\item For $\mu=\delta_c$ the free cumulants are $c,0,0,\ldots$

\item For $\mu=\gamma_t$ the free cumulants are $0,t,0,0,\ldots$

\item For $\mu=\pi_t$ the free cumulants are $t,t,t,\ldots$

\item For $\mu=\beta_t$ the free cumulants are $0,t,0,t,\ldots$
\end{enumerate}
Also, for compound free Poisson laws the free cumulants are $k_n(\pi_\nu)=M_n(\nu)$.
\end{proposition}

\begin{proof}
The proofs are analogous to those from the classical case, as follows:

\medskip

(1) For $\mu=\delta_c$ we have $G_\mu(\xi)=1/(\xi-c)$, and so $R_\mu(\xi)=c$, as desired.

\medskip

(2) For $\mu=\gamma_t$ we have, as computed before, $R_\mu(\xi)=t\xi$, as desired.

\medskip

(3) For $\mu=\pi_t$ we have, also from before, $R_\mu(\xi)=t/(1-\xi)$, as desired.

\medskip

(4) For $\mu=\beta_t$ this follows from the formulae that we have, but the best is to prove directly the last assertion, which generalizes (3,4). With $\nu=\sum_ic_i\delta_{z_i}$ we have:
\begin{eqnarray*}
R_{\pi_\nu}(\xi)
&=&\sum_i\frac{c_iz_i}{1-\xi z_i}\\
&=&\sum_ic_iz_i\sum_{n\geq0}(\xi z_i)^n\\
&=&\sum_{n\geq0}\xi^n\sum_ic_iz_i^{n+1}\\
&=&\sum_{n\geq1}\xi^{n-1}\sum_ic_iz_i^n\\
&=&\sum_{n\geq 1}\xi^{n-1}\,M_n(\nu)
\end{eqnarray*}

Thus, we are led to the conclusion in the statement.
\end{proof}

Also as before in the classical case, we can define generalized free cumulants $\kappa_\pi(a)$ with $\pi\in P(k)$ by starting with the numeric free cumulants $\kappa_n(a)$, as follows:

\begin{definition}
We define free cumulants $\kappa_\pi(a)$, depending on partitions 
$$\pi\in P(k)$$
by starting with $\kappa_n(a)$, and using multiplicativity over the blocks. 
\end{definition}

To be more precise, the convention here is that for the one-block partition $1_n\in P(n)$, the corresponding free cumulant is the usual one, namely:
$$\kappa_{1_n}(a)=\kappa_n(a)$$

Then, for an arbitrary partition $\pi\in P(k)$, we decompose this partition into blocks, having sizes $b_1,\ldots,b_s$, and we set, by multiplicativity over blocks:
$$\kappa_\pi(a)=\kappa_{b_1}(a)\ldots\kappa_{b_s}(a)$$

With this convention, we have the following result, due to Speicher \cite{spe}:

\index{moment-cumulant formula}

\begin{theorem}
We have the moment-cumulant formulae
$$M_n(a)=\sum_{\nu\in NC(n)}\kappa_\nu(a)\quad,\quad 
\kappa_n(a)=\sum_{\nu\in NC(n)}\mu(\nu,1_n)M_\nu(a)$$
or, equivalently, we have the moment-cumulant formulae
$$M_\pi(a)=\sum_{\nu\leq\pi}\kappa_\nu(a)\quad,\quad 
\kappa_\pi(a)=\sum_{\nu\leq\pi}\mu(\nu,\pi)M_\nu(a)$$
where $\mu$ is the M\"obius function of $NC(n)$.
\end{theorem}

\begin{proof}
This follows indeed by doing some combinatorics and calculus, in the spirit of the combinatorics and calculus from the classical case.
\end{proof}

In relation now with easiness, we will need the following result, from \cite{bsp}:

\begin{proposition}
The free orthogonal easy groups $G\subset O_N^+$ that we know,
$$S_N^+\subset H_N^+\subset O_N^+$$
come respectively from the following categories of noncrossing partitions $D\subset NC$,
$$NC\supset NC_{even}\supset NC_2$$
which in turn come from certain subsets $L\subset\mathbb N$, as follows,
$$\mathbb N\supset 2\mathbb N\supset\{2\}$$
with $D$ consisting of the partitions $\pi\in NC$ whose blocks have lengths belonging to $L\subset\mathbb N$. 
\end{proposition}

\begin{proof}
As before with Proposition 4.30, this is something quite trivial, which follows from what we know about the quantum groups in the statement, but which hides behind some interesting mathematics, to be further discussed later, in chapter 6 below.
\end{proof}

In relation with cumulants, we have the following result, also from \cite{bsp}:

\index{truncated character}
\index{free quantum group}
\index{asymptotic moment}

\begin{theorem}
The free cumulants of the asymptotic truncated characters for the quantum groups $S_N^+,H_N^+,O_N^+$ are given by the formula
$$\kappa_n(\chi_t)=t\delta_{n\in L}$$
with $L\subset\mathbb N$ being the associated subset, and at the level of asymptotic moments this gives 
$$M_k(\chi_t)=\sum_{\pi\in D(k)}t^{|\pi|}$$
with $D\subset NC$ being the associated category of partitions.
\end{theorem}

\begin{proof}
This is clear indeed from Proposition 4.37, by performing a case-by-case analysis, with the cases $G=S^+,H^+,O^+$ under consideration corresponding to the cumulant computations for the measures $\pi_t,\beta_t,\gamma_t$ from Proposition 4.34. See \cite{bsp}. 
\end{proof}

Summarizing, and getting back now to Theorem 4.21, we have reached to a perfect understanding of what is going on there, via classical and free cumulants. And, as a bonus, we can answer Question 4.22 as well, in a very conceptual way, as follows:

\begin{answer}
For the main examples of liberations of easy groups, $G_N\to G_N^+$, the asymptotic laws of truncated characters
$$\chi_t=\sum_{i=1}^{[tN]}u_{ii}$$ 
are in Bercovici-Pata bijection, in the sense that the classical cumulants of the laws for $G_N$ coincide with the free cumulants of the laws for $G_N^+$.
\end{answer}

Here the bijection mentioned at the end, which is something purely probabilistic, and very beautiful and conceptual, comes from the paper of Bercovici-Pata \cite{bpa}, and we refer to that paper and the subsequent free probability literature for more on all this.

\bigskip

Finally, as already mentioned on several occasions, all the above is part of some general theory regarding the ``uniform'' easy quantum groups, which remains to be developed. We will be back to this, with full details, in chapter 6 below.

\section*{4d. Invariance questions} 

As a last topic for this chapter, following \cite{ez3}, let us discuss now probabilistic invariance questions with respect to the basic orthogonal quantum groups. We first have:

\begin{definition}
Given a closed subgroup $G\subset O_N^+$, we denote by
$$\alpha:\mathbb C<t_1,\ldots,t_N>\to\mathbb C<t_1,\ldots,t_N>\otimes\,C(G)$$
$$t_i\to\sum_jt_j\otimes v_{ji}$$
the standard coaction of $C(G)$ on the free complex algebra on $N$ variables.
\end{definition}

Observe that the map $\alpha$ constructed above is indeed a coaction, in the sense that it satisfies the following standard coassociativity and counitality conditions:
$$(id\otimes\Delta)\alpha=(\alpha\otimes id)\alpha\quad,\quad 
(id\otimes\varepsilon)\alpha=id$$

With the above notion of coaction in hand, we can now talk about invariant sequences of classical or noncommutative random variables, in the following way:

\begin{definition}
Let $(B,tr)$ be a $C^*$-algebra with a trace, and $x_1,\ldots,x_N\in B$. We say that $x=(x_1,\ldots,x_N)$ is invariant under $G\subset O_N^+$ if the distribution functional 
$$\mu_x:\mathbb C<t_1,\ldots,t_N>\to\mathbb C\quad,\quad 
P\to tr(P(x_1,\ldots,x_N))$$
is invariant under the coaction $\alpha$, in the sense that we have
$$(\mu_x\otimes id)\alpha(P)=\mu_x(P)$$
for any noncommuting polynomial $P\in\mathbb C<t_1,\ldots,t_N>$.
\end{definition}

Observe that in the classical case, where $G\subset O_N$ is a usual group, we recover in this way the usual invariance notion from classical probability. We have the following equivalent formulation of the above invariance condition: 

\index{invariant sequence}

\begin{proposition}
Let $(B,tr)$ be a $C^*$-algebra with a trace, and $x_1,\ldots,x_N\in B$. Then $x=(x_1,\ldots,x_N)$ is invariant under $G\subset O_N^+$ precisely when
$$tr(x_{i_1}\ldots x_{i_k})=\sum_{j_1\ldots j_k}tr(x_{j_1}\ldots x_{j_k})v_{j_1i_1}\ldots v_{j_ki_k}$$
as an equality in $C(G)$, for any $k\in\mathbb N$, and any $i_1,\ldots,i_k\in\{1,\ldots,N\}$.  
\end{proposition}

\begin{proof}
By linearity, in order for a sequence $x=(x_1,\ldots,x_N)$ to be $G$-invariant in the sense of Definition 4.41, the formula there must be satisfied for any noncommuting monomial $P\in\mathbb C<t_1,\ldots,t_N>$. But, with $P=t_{i_1}\ldots t_{i_k}$, we have:
\begin{eqnarray*}
(\mu_x\otimes id)\alpha(P)
&=&(\mu_x\otimes id)\sum_{j_1,\ldots,j_k}t_{j_1}\ldots t_{j_k}\otimes v_{j_1i_1}\ldots v_{j_ki_k}\\
&=&\sum_{j_1,\ldots,j_k}\mu_x(t_{j_1}\ldots t_{j_k})v_{j_1i_1}\ldots v_{j_ki_k}\\
&=&\sum_{j_1\ldots j_k}tr(x_{j_1}\ldots x_{j_k})v_{j_1i_1}\ldots v_{j_ki_k}
\end{eqnarray*}

On the other hand, by definition of the distribution $\mu_x$, we have:
$$\mu_x(P)
=\mu_x(t_{i_1}\ldots t_{i_k})
=tr(x_{i_1}\ldots x_{i_k})$$

Thus, we are led to the conclusion in the statement.
\end{proof}

We can now investigate invariance questions for the sequences of classical or noncommutative random variables, with respect to the main examples of easy quantum groups. We first have a reverse De Finetti theorem, from \cite{ez3}, as follows:

\index{reverse De Finetti}

\begin{theorem}
Let $(x_1,\ldots,x_N)$ be a sequence in $A$.
\begin{enumerate}
\item If $x_1,\ldots,x_N$ are freely independent and identically distributed with amalgamation over $B$, then the sequence is $S_N^+$-invariant.

\item If $x_1,\ldots,x_N$ are freely independent and identically distributed with amalgamation over $B$, and have centered semicircular distributions with respect to $E$, then the sequence is $O_N^+$-invariant.

\item If $<B,x_1,\ldots,x_N>$ is commutative and $x_1,\ldots,x_N$ are conditionally independent and identically distributed given $B$, then the sequence is $S_N$-invariant.

\item If $<x_1,\ldots,x_N>$ is commutative and $x_1,\ldots,x_N$ are conditionally independent and identically distributed given $B$, and have centered Gaussian distributions with respect to $E$, then the sequence is $O_N$-invariant.
\end{enumerate}
\end{theorem}

\begin{proof}
Assume that the joint distribution of $(x_1,\ldots,x_N)$ satisfies one of the conditions in the statement, and let $D$ be the category of partitions associated to the corresponding easy quantum group.  We have then the following computation:
\begin{eqnarray*}
\sum_{j_1\ldots j_k}tr(x_{j_1}\ldots x_{j_k})v_{j_1i_1}\ldots v_{j_ki_k}
&=&\sum_{j_1\ldots j_k}tr(E(x_{j_1}\ldots x_{j_k}))v_{j_1i_1}\ldots v_{j_ki_k}\\
&=&\sum_{j_1\ldots j_k}\sum_{\pi\leq\ker j}tr(\xi^{(\pi)}_E(x_1,\ldots,x_1))v_{j_1i_1}\ldots v_{j_ki_k}\\
&=&\sum_{\pi\in D(k)}tr(\xi^{(\pi)}_E(x_1,\ldots,x_1))\sum_{\ker j\geq\pi}v_{j_1i_1}\ldots v_{j_ki_k}
\end{eqnarray*}

Here $\xi$ denotes the free and classical cumulants in the cases (1,2) and (3,4) respectively. On the other hand, it follows from a direct computation that if $\pi\in D(k)$ then we have the following formula, in each of the 4 cases in the statement:
$$\sum_{\ker j\geq\pi}v_{j_1i_1}\ldots v_{j_ki_k}=
\begin{cases}1&{\rm if}\ \pi\leq\ker i\\ 
0&{\rm otherwise}
\end{cases}$$

By using this formula, we can finish our computation, in the following way:
\begin{eqnarray*}
\sum_{j_1\ldots j_k}tr(x_{j_1}\ldots x_{j_k})v_{j_1i_1}\ldots v_{j_ki_k}
&=&\sum_{\pi\in D(k)}tr(\xi^{(\pi)}_E(x_1,\ldots,x_1))\delta_{\pi\leq\ker i}\\
&=&\sum_{\pi\leq\ker i}tr(\xi_E^{(\pi)}(x_1,\ldots,x_1))\\
&=&tr(x_{i_1}\ldots x_{i_k})
\end{eqnarray*}

Thus, we are led to the conclusions in the statement.
\end{proof}

Still following \cite{ez3}, as main result, we have the following converse to Theorem 4.43:

\index{De Finetti theorem}

\begin{theorem}
Let $(x_i)_{i\in\mathbb N}$ be a $G$-invariant sequence of self-adjoint random variables in $(M,tr)$, and assume that $M=<(x_i)_{i\in\mathbb N}>$. Then there exists a subalgebra $B\subset M$ and a trace-preserving conditional expectation $E:M\to B$ such that:
\begin{enumerate}
\item If $G=(S_N)$, then $(x_i)_{i\in\mathbb N}$ are conditionally independent and identically distributed given $B$.

\item If $G=(S_N^+)$, then $(x_i)_{i\in\mathbb N}$ are freely independent and identically distributed with amalgamation over $B$.

\item If $G=(O_N)$, then $(x_i)_{i\in\mathbb N}$ are conditionally independent, and have Gaussian distributions with mean zero and common variance, given $B$.

\item If $G=(O_N^+)$, then $(x_i)_{i\in\mathbb N}$ form a $B$-valued free semicircular family with mean zero and common variance.
\end{enumerate}
\end{theorem}

\begin{proof}
This is something quite technical, heavily based on the Weingarten formula. Let $j_1,\ldots,j_k \in \mathbb N$ and $b_0,\ldots,b_k\in B$. We have then the following computation:
\begin{eqnarray*}
E(b_0x_{j_1}\ldots x_{j_k}b_k)
&=&\lim_{N\to\infty}E_N(b_0x_{j_1}\ldots x_{j_k}b_k)\\
&=&\lim_{N\to\infty}\sum_{\sigma\leq\ker j}\sum_\pi W_{kN}(\pi,\sigma) \sum_{\pi\leq\ker i}b_0x_{i_1}\ldots x_{i_k}b_k\\
&=&\lim_{N\to\infty}\sum_{\sigma\leq\ker j}\sum_{\pi\leq\sigma}\mu_{D(k)}(\pi,\sigma)N^{-|\pi|}\sum_{\pi\leq\ker i}b_0x_{i_1}\ldots x_{i_k}b_k
\end{eqnarray*}

By using this formula, and various cumulant computations, we have:
$$E(b_0x_{j_1}\ldots x_{j_k}b_k)=\lim_{N\to\infty}\sum_{\sigma\leq\ker j} \sum_{\pi\leq\sigma}\mu_{D(k)}(\pi,\sigma)E_N^{(\pi)}(b_0x_1b_1,\ldots,x_1b_k)$$

We therefore obtain the following formula:
$$E(b_0x_{j_1}\ldots x_{j_k}b_k)=\sum_{\sigma\leq\ker j}\sum_{\pi\leq\sigma} \mu_{D(k)}(\pi,\sigma)E^{(\pi)}(b_0x_1b_1,\ldots,x_1b_k)$$

We can replace the sum of expectation functionals by cumulants, as to obtain:
$$E(b_0x_{j_1}\ldots x_{j_k}b_k)=\sum_{\sigma\leq\ker j}\xi_E^{(\sigma)}(b_0x_1b_1,\ldots,x_1b_k)$$

Here and in what follows $\xi$ denotes as usual the relevant classical or free cumulants, depending on the quantum group that we are dealing with, either classical or free. Now since the classical or free cumulants are determined by the classical or free moment-cumulant formulae, explained before, we conclude that we have the following formula:
$$\xi_E^{(\sigma)}(b_0x_{j_1}b_1,\ldots,x_{j_k}b_k)
=\begin{cases}
\xi_E^{(\sigma)}(b_0x_1b_1,\ldots,x_1b_k)&{\rm if}\ \sigma\in D(k)\ {\rm and}\ \sigma\leq\ker j\\
0&{\rm otherwise}
\end{cases}$$

With this formula in hand, the result then follows from the characterizations of these joint distributions in terms of cumulants, and we are done.
\end{proof}

All the above was of course quite brief, and for more on all this, including full proof, and various versions of the above results, we refer to \cite{ez3} and related papers.

\section*{4e. Exercises} 

We are now into advanced probability theory, and in order to fully understand what has been said above, a lot of extra work is needed. Generally speaking, the best here would be to carefully read a probability book, and then a free probability book. In relation now with what we did in this chapter, as a first instructive exercise, we have:

\begin{exercise}
Recover the Poisson law results for $S_N$ via the formula
$$\int_{S_N}u_{i_1j_1}\ldots u_{i_rj_r}=\begin{cases}
\frac{(N-|\ker i|)!}{N!}&{\rm if}\ \ker i=\ker j\\
0&{\rm otherwise}
\end{cases}$$
that you will have to establish first.
\end{exercise}

Along the same lines, we have as well:

\begin{exercise}
Recover the Bessel law results for $H_N$, then for $H_N^s$, via inclusion-exclusion, or via an explicit integration formula, as in the previous exercise.
\end{exercise}

Regarding Weingarten integration, which is the method that must be mastered, for anything probabilistic in relation with groups and quantum groups, we have:

\begin{exercise}
Write down the Gram and Weingarten matrices of small order, for all the easy groups and quantum groups that you know.
\end{exercise}

At the theoretical level now, you cannot really escape from:

\begin{exercise}
Learn classical and free probability, as much as you can, and in any case enough of them, as to declare the present chapter fully understood.
\end{exercise}

In relation with classical cumulants, we have the following instructive exercise:

\begin{exercise}
Establish the following formulae for the classical cumulants,
$$k_1=M_1$$
$$k_2=-M_1^2+M_2$$
$$k_3=2M_1^3-3M_1M_2+M_3$$
$$k_4=-6M_1^4+12M_1^2M_2-3M_2^2-4M_1M_3+M_4$$
$$\ldots$$
then go on and prove the classical moment-cumulant formula.
\end{exercise}

Similarly, in relation with the free cumulants, we have:

\begin{exercise}
Establish the following formulae for the free cumulants,
$$\kappa_1=M_1$$
$$\kappa_2=-M_1^2+M_2$$
$$\kappa_3=2M_1^3-3M_1M_2+M_3$$
$$\kappa_4=-5M_1^4+10M_1^2M_2-2M_2^2-4M_1M_3+M_4$$
$$\ldots$$
then go on and prove the free moment-cumulant formula.
\end{exercise}

As a comment here, the above two exercises look quite similar, but doing only one of them will not do, because they are in fact quite different. Just try a bit, and you will understand. Finally, in relation with De Finetti, waiting for more here, from you.

\part{The continuous case}

\ \vskip50mm

\begin{center}
{\em We belong to the light

We belong to the thunder

We belong to the sound of the words

We've both fallen under}
\end{center}

\chapter{Bistochastic groups}

\section*{5a. Bistochastic groups}

Welcome to advanced easiness. With the basics understood, we will go, in the present Part II,  and in  Part III as well, into the study of more advanced aspects. We will be interested, as usual, in constructing more examples of easy quantum groups, and in developing some algebraic, geometric, analytic and probabilistic theory for them.

\bigskip

Needless to say, there has been a lot of work on the subject, all over the 2010s, and the things to talk about, both further examples and further theory, abound. So obviously, we need some sort of plan and strategy, for explaining this material. But the answer here comes from the cube formed by the main examples of easy quantum groups, namely: 
$$
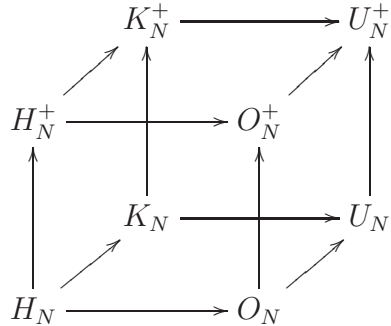
\xymatrix@R=18pt@C=18pt{
&K_N^+\ar[rr]&&U_N^+\\
H_N^+\ar[rr]\ar[ur]&&O_N^+\ar[ur]\\
&K_N\ar[rr]\ar[uu]&&U_N\ar[uu]\\
H_N\ar[uu]\ar[ur]\ar[rr]&&O_N\ar[uu]\ar[ur]
}$$

Indeed, this cube provides us with some sort of ``3D orientation'' into the whole subject, with the 3 coordinate axes corresponding to the dualities real/complex, discrete/continuous, and classical/free. So, speaking presentation, the question is now: shall we talk first in detail about the real case, and deal with the complex case afterwards? Or shall we first settle the questions regarding the discrete case, and then talk about the continuous case? Or, why not vice versa, continuous first, and discrete after? Or, shall we finish the work in the classical case, and discuss liberation afterwards?

\bigskip

These are not easy questions, and a look at the various papers written on the subject, the serious ones I mean, does not provide any clue, with all options being on the table. So, in the lack of any good idea here, we will ask the cat. And cat says:

\begin{cat}
Continuous comes first.
\end{cat}

This sounds reasonable, thanks cat. Actually my cat is named Felix, and thinking well, he might be well a reincarnation of Felix Klein. Indeed, at the time of Klein there was no abstract mathematics, and such boring things as abstract groups, fields of scalars, and so on, and Klein was simply interested in the ``continuous transformations'' of the world surrounding us. And so will be us, to start with: this is definitely reasonable.

\bigskip

In practice, this suggests focusing towards the right face of the cube, namely:
$$\xymatrix@R=51pt@C=50pt{
O_N^+\ar[r]&U_N^+\\
O_N\ar[u]\ar[r]&U_N\ar[u]}$$

For instance, we can say that any intermediate easy quantum group $O_N\subset G\subset U_N^+$ is by definition ``continuous'', and then try to understand the structure of such quantum groups, with a classification result for them, and then with algebraic, geometric, analytic and probabilistic results for the objects that we found. This certainly sounds reasonable, and looks doable as well, because after all we will deal with pairings, the diagram of categories of partitions for our quantum groups above being as follows:
$$\xymatrix@R=18.5mm@C=16.2mm{
NC_2\ar[d]&\mathcal{NC}_2\ar[l]\ar[d]\\
P_2&\mathcal P_2\ar[l]}$$

However, this is a bit superficial, because there might be well interesting easy quantum groups $S_N\subset G\subset U_N^+$ which definitely deserve the name ``continuous'', and which are excluded by our formalism, not containing $O_N$. So, let us look into this first. It is probably safe here to assume that $G$ is classical, at least to start with, because in the classical case we perfectly know what ``continuous'' means. So, we are led to:

\begin{problem}
Besides $O_N,U_N$, what are the simplest easy Lie groups?
\end{problem}

This does not look as a trivial question, and normally we should ask the cat here, but Felix is now away. As for the other cat, Sophus, who surely knows the answer too, he is nowhere to be found, either. So, I guess this is a problem for me and you, reader.

\bigskip

In order to solve this question, let us relax, and looks towards applied linear algebra. A very basic definition there, that you might already know, is as follows:

\index{row-stochastic}
\index{column-stochastic}
\index{bistochastic}

\begin{definition}
A square matrix $M\in M_N(\mathbb C)$ is called bistochastic if each row and each column sum up to the same number:
$$\begin{matrix}
M_{11}&\ldots&M_{1N}&\to&\lambda\\
\vdots&&\vdots\\
M_{N1}&\ldots&M_{NN}&\to&\lambda\\
\downarrow&&\downarrow\\
\lambda&&\lambda
\end{matrix}$$
If this happens only for the rows, or only for the columns, the matrix is called row-stochastic, respectively column-stochastic.
\end{definition}

As the name indicates, these matrices are useful in statistics, and perhaps in other fields like graph theory, computer science and so on, with the case of the matrices having entries in $[0,1]$, which sum up to $\lambda=1$, being the important one. As a basic example of a bistochastic matrix, we have the flat matrix, which is as follows:
$$\mathbb I_N=\begin{pmatrix}
1&\ldots&1\\
\vdots&&\vdots\\
1&\ldots&1
\end{pmatrix}$$

Observe that the rescaling $P_N=\mathbb I_N/N$ has the property mentioned above, namely entries in $[0,1]$, summing up to $1$. In fact, this matrix $P_N=\mathbb I_N/N$ is a very familiar object in linear algebra, being the projection on the all-one vector, namely:
$$\xi=\begin{pmatrix}
1\\
\vdots\\
1
\end{pmatrix}$$ 

Getting back now to the general case, the various notions of stochasticity in Definition 5.3 are closely related to this vector $\xi$, due to the following simple fact:

\begin{proposition}
Let $M\in M_N(\mathbb C)$ be a square matrix.
\begin{enumerate}
\item $M$ is row stochastic, with sums $\lambda$, when $M\xi=\lambda\xi$.

\item $M$ is column stochastic, with sums $\lambda$, when $M^t\xi=\lambda\xi$.

\item $M$ is bistochastic, with sums $\lambda$, when $M\xi=M^t\xi=\lambda\xi$.
\end{enumerate}
\end{proposition}

\begin{proof}
The first assertion is clear from definitions, because when multiplying a matrix by $\xi$, we obtain the vector formed by the row sums:
$$\begin{pmatrix}
M_{11}&\ldots&M_{1N}\\
\vdots&&\vdots\\
M_{N1}&\ldots&M_{NN}
\end{pmatrix}
\begin{pmatrix}
1\\
\vdots\\
1
\end{pmatrix}
=\begin{pmatrix}
M_{11}+\ldots+M_{1N}\\
\vdots\\
M_{N1}+\ldots+M_{NN}
\end{pmatrix}$$

As for the second, and then third assertion, these are both clear from this.
\end{proof}

As an observation here, we can reformulate if we want the above statement in a purely matrix-theoretic form, by using the flat matrix $\mathbb I_N$, as follows:

\index{flat matrix}

\begin{proposition}
Let $M\in M_N(\mathbb C)$ be a square matrix.
\begin{enumerate}
\item $M$ is row stochastic, with sums $\lambda$, when $M\mathbb I_N=\lambda\mathbb I_N$.

\item $M$ is column stochastic, with sums $\lambda$, when $\mathbb I_NM=\lambda\mathbb I_N$.

\item $M$ is bistochastic, with sums $\lambda$, when $M\mathbb I_N=\mathbb I_NM=\lambda\mathbb I_N$.
\end{enumerate}
\end{proposition}

\begin{proof}
This follows from Proposition 5.4, and from the fact that both the rows and columns of the flat matrix $\mathbb I_N$ are copies of the all-one vector $\xi$. Alternatively, we have the following formula, $S_1,\ldots,S_N$ being the row sums of $M$, which gives (1):
$$\begin{pmatrix}
M_{11}&\ldots&M_{1N}\\
\vdots&&\vdots\\
M_{N1}&\ldots&M_{NN}
\end{pmatrix}
\begin{pmatrix}
1&\ldots&1\\
\vdots&&\vdots\\
1&\ldots&1
\end{pmatrix}
=\begin{pmatrix}
S_1&\ldots&S_1\\
\vdots&&\vdots\\
S_N&\ldots&S_N
\end{pmatrix}$$

As for the second, and then third assertion, these are both clear from this.
\end{proof}

In what follows, we will be mainly interested in the bistochastic matrices which are unitary, $M\in U_N$. As the simplest example here, which is a familiar object in quantum physics, we have the following matrix $K_N\in O_N$, obtained by suitably modifying the flat matrix $\mathbb I_N$, as to make the rows pairwise orthogonal, and of norm one:
$$K_N=\frac{1}{N}\begin{pmatrix}
2-N&2&\ldots&2\\
2&\ddots&\ddots&\vdots\\
\vdots&\ddots&\ddots&2\\
2&\ldots&2&2-N
\end{pmatrix}$$

As a first result regarding the unitary bistochastic matrices, we have:

\begin{theorem}
For a unitary matrix $U\in U_N$, the following conditions are equivalent:
\begin{enumerate}
\item $H$ is bistochastic, with sums $\lambda$.

\item $H$ is row stochastic, with sums $\lambda$, and $|\lambda|=1$.

\item $H$ is column stochastic, with sums $\lambda$, and $|\lambda|=1$.
\end{enumerate}
\end{theorem}

\begin{proof}
By using a symmetry argument we just need to prove $(1)\iff(2)$, and both the implications are elementary, as follows:

\medskip

$(1)\implies(2)$ If we denote by $U_1,\ldots,U_N\in\mathbb C^N$ the rows of $U$, we have indeed:
\begin{eqnarray*}
1
&=&<U_1,U_1>\\
&=&\sum_i<U_1,U_i>\\
&=&\sum_i\sum_jU_{1j}\bar{U}_{ij}\\
&=&\sum_jU_{1j}\sum_i\bar{U}_{ij}\\
&=&\sum_jU_{1j}\cdot\bar{\lambda}\\
&=&\lambda\cdot\bar{\lambda}\\
&=&|\lambda|^2
\end{eqnarray*}

$(2)\implies(1)$ Consider the all-one vector $\xi=(1)_i\in\mathbb C^N$. The fact that $U$ is row-stochastic with sums $\lambda$ reads:
\begin{eqnarray*}
\sum_jU_{ij}=\lambda,\forall i
&\iff&\sum_jU_{ij}\xi_j=\lambda\xi_i,\forall i\\
&\iff&U\xi=\lambda\xi
\end{eqnarray*}

Also, the fact that $U$ is column-stochastic with sums $\lambda$ reads:
\begin{eqnarray*}
\sum_iU_{ij}=\lambda,\forall j
&\iff&\sum_jU_{ij}\xi_i=\lambda\xi_j,\forall j\\
&\iff&U^t\xi=\lambda\xi
\end{eqnarray*}

We must prove that the first condition implies the second one, provided that the row sum $\lambda$ satisfies $|\lambda|=1$. But this follows from the following computation:
\begin{eqnarray*}
U\xi=\lambda\xi
&\implies&U^*U\xi=\lambda U^*\xi\\
&\implies&\xi=\lambda U^*\xi\\
&\implies&\xi=\bar{\lambda}U^t\xi\\
&\implies&U^t\xi=\lambda\xi
\end{eqnarray*}

Thus, we have proved both the implications, and we are done.
\end{proof}

Getting now to our questions regarding groups and quantum groups, we would like to talk about the orthogonal and unitary groups of bistochastic matrices, and their free analogues. However, there is a choice to be made here, in connection with what we know from Theorem 5.6, namely shall we let the row and column sums to be arbitrary, $\lambda\in\mathbb T$, or shall we make the normalization $\lambda=1$. For various reasons that will become clear later on, it is better to choose this latter way, so let us formulate:

\begin{convention}
From now on all our bistochastic matrices will be assumed to be normalized, with the sum on all rows and columns being equal to $\lambda=1$.
\end{convention}

To be more precise, this convention is here in order for the resulting groups of bistochastic orthogonal and unitary matrices to be easy. But more on this later. For the moment, let us introduce these groups. We have the following result:

\index{bistochastic group}
\index{complex bistochastic group}

\begin{theorem}
We have closed subgroups as follows:
\begin{enumerate}
\item $B_N\subset O_N$, consisting of the orthogonal matrices which are bistochastic.

\item $C_N\subset U_N$, consisting of the unitary matrices which are bistochastic.
\end{enumerate}
\end{theorem}

\begin{proof}
We know from Theorem 5.6 that the sets of bistochastic matrices $B_N,C_N$ in the statement appear as follows, with $\xi$ being the all-one vector:
$$B_N=\left\{U\in O_N\Big|U\xi=\xi\right\}$$
$$C_N=\left\{U\in U_N\Big|U\xi=\xi\right\}$$

It is then clear that both $B_N,C_N$ are stable under the multiplication, contain the unit, and are stable by inversion. Thus, we have indeed closed subgroups, as stated.
\end{proof}

As already mentioned, the bistochastic matrices and groups are important in statistics, applied linear algebra, quantum physics, and many more. Following Idel-Wolf \cite{iwo} and related papers, we would like to discuss now, as an introduction to this, a non-trivial result regarding the unitary bistochastic group $C_N$. Some advertisement first:

\begin{question}
Is there anything in linear algebra having a one-line statement, that everyone can understand, but a proof so complicated, that no one really understands?
\end{question}

We will see that the answer to this question is yes. Let us begin with some geometric preliminaries. The complex projective space appears by definition as follows:
$$P^{N-1}_\mathbb C=(\mathbb C^N-\{0\})\big/<x=\lambda y>$$

Inside this projective space, we have the Clifford torus, constructed as follows:
$$\mathbb T^{N-1}=\left\{(z_1,\ldots,z_N)\in P^{N-1}_\mathbb C\Big||z_1|=\ldots=|z_N|\right\}$$

With these conventions, we have the following result, from \cite{iwo}:

\index{complex projective space}
\index{Clifford torus}

\begin{proposition}
For a unitary matrix $U\in U_N$, the following are equivalent:
\begin{enumerate}
\item There exist $L,R\in U_N$ diagonal such that $U'=LUR$ is bistochastic.

\item The standard torus $\mathbb T^N\subset\mathbb C^N$ satisfies $\mathbb T^N\cap U\mathbb T^N\neq\emptyset$.

\item The Clifford torus $\mathbb T^{N-1}\subset P^{N-1}_\mathbb C$ satisfies $\mathbb T^{N-1}\cap U\mathbb T^{N-1}\neq\emptyset$.
\end{enumerate}
\end{proposition}

\begin{proof}
These equivalences are all elementary, as follows:

\medskip

$(1)\implies(2)$ Assuming that $U'=LUR$ is bistochastic, which in terms of the all-1 vector $\xi$ means $U'\xi=\xi$, if we set $f=R\xi\in\mathbb T^N$ we have:
$$Uf
=\bar{L}U'\bar{R}f
=\bar{L}U'\xi
=\bar{L}\xi\in\mathbb T^N$$

Thus we have $Uf\in\mathbb T^N\cap U\mathbb T^N$, which gives the conclusion.

\medskip

$(2)\implies(1)$ Given $g\in\mathbb T^N\cap U\mathbb T^N$, we can define $R,L$ as follows:
$$R=\begin{pmatrix}g_1\\&\ddots\\&&g_N\end{pmatrix}\quad,\quad 
\bar{L}=\begin{pmatrix}(Ug)_1\\&\ddots\\&&(Ug)_N\end{pmatrix}$$

With these values for $L,R$, we have then the following formulae:
$$R\xi=g\quad,\quad 
\bar{L}\xi=Ug$$

Thus the matrix $U'=LUR$ is bistochastic, because:
$$U'\xi
=LUR\xi
=LUg
=\xi$$

$(2)\implies(3)$ This is clear, because $\mathbb T^{N-1}\subset P^{N-1}_\mathbb C$ appears as the projective image of $\mathbb T^N\subset\mathbb C^N$, and so $\mathbb T^{N-1}\cap U\mathbb T^{N-1}$ appears as the projective image of $\mathbb T^N\cap U\mathbb T^N$.

\medskip

$(3)\implies(2)$ We have indeed the following equivalence:
$$\mathbb T^{N-1}\cap U\mathbb T^{N-1}\neq\emptyset
\iff\exists\lambda\neq 0,\lambda\mathbb T^N\cap U\mathbb T^N\neq\emptyset$$

But $U\in U_N$ implies $|\lambda|=1$, and this gives the result.
\end{proof}

The point now is that the condition (3) above is something familiar in symplectic geometry, and known to hold for any $U\in U_N$. Thus, following \cite{iwo}, we have:

\index{bistochastic form}
\index{Sinkhorn normal form}
\index{Lagrangian submanifold}
\index{Hamiltonian isotopy}
\index{Arnold conjecture}

\begin{theorem}
Any unitary matrix $U\in U_N$ can be put in bistochastic form,
$$U'=LUR$$
with $L,R\in U_N$ being both diagonal, via a certain non-explicit method.
\end{theorem}

\begin{proof}
As already mentioned, the condition $\mathbb T^{N-1}\cap U\mathbb T^{N-1}\neq\emptyset$ in Proposition 5.10 (3) is something quite natural in symplectic geometry. To be more precise:

\medskip

(1) The Clifford torus $\mathbb T^{N-1}\subset P^{N-1}_\mathbb C$ is a Lagrangian submanifold, and the map $\mathbb T^{N-1}\to U\mathbb T^{N-1}$ is a Hamiltonian isotopy. For more on this, see Arnold \cite{arn}.

\medskip

(2) A non-trivial result of Biran-Entov-Polterovich \cite{bep} and Cho \cite{cho} states that $\mathbb T^{N-1}$ cannot be displaced from itself via a Hamiltonian isotopy. 

\medskip

(3) Thus, the results in \cite{bep}, \cite{cho} tells us that  $\mathbb T^{N-1}\cap U\mathbb T^{N-1}\neq\emptyset$ holds indeed, for any $U\in U_N$. We therefore obtain the result, via Proposition 5.10. See \cite{iwo}.
\end{proof}

There are many further things that can be said here. As explained in \cite{iwo}, the various technical results from \cite{bep}, \cite{cho} show that in the generic, ``transverse'' situation, there are at least $2^{N-1}$ ways of putting a unitary matrix $U\in U_N$ in bistochastic form, and this modulo the obvious transformation $U\to zU$, with $|z|=1$. However, all this still does not clarify things with Theorem 5.11, whose only known proof in general is the one above.

\section*{5b. Singletons and pairings}

Getting back now to quantum groups, we would like to talk about the free analogues $B_N^+,C_N^+$ of the orthogonal and unitary bistochastic groups. Following \cite{bsp}, we have the following result, including as well our previous construction of $B_N,C_N$:

\index{bistochsastic group}
\index{free bistochastic group}
\index{complex bistochastic group}
\index{bistochastic quantum group}

\begin{theorem}
We have the following groups and quantum groups:
\begin{enumerate}
\item $B_N\subset O_N$, consisting of the orthogonal matrices which are bistochastic.

\item $C_N\subset U_N$, consisting of the unitary matrices which are bistochastic.

\item $B_N^+\subset O_N^+$, coming via $u\xi=\xi$, where $\xi$ is the all-one vector.

\item $C_N^+\subset U_N^+$, coming via $u\xi=\xi$, where $\xi$ is the all-one vector.
\end{enumerate}
Also, we have inclusions $B_N\subset B_N^+$ and $C_N\subset C_N^+$, which are both liberations.
\end{theorem}

\begin{proof}
There are several things to be proved, the idea being as follows:

\medskip

(1) We already know from Theorem 5.8 that $B_N,C_N$ are indeed groups, with this coming from the following formulae, with $\xi$ being the all-one vector:
$$B_N=\left\{U\in O_N\Big|U\xi=\xi\right\}$$
$$C_N=\left\{U\in U_N\Big|U\xi=\xi\right\}$$

(2) In what regards now $B_N^+,C_N^+$, these appear by definition as follows:
$$C(B_N^+)=C(O_N^+)\Big/\Big<\xi\in Fix(u)\Big>$$
$$C(C_N^+)=C(U_N^+)\Big/\Big<\xi\in Fix(u)\Big>$$

But since the relation $\xi\in Fix(u)$ is categorical, we have indeed quantum groups.

\medskip

(3) Finally, in what regards the last assertion, since we already know that $O_N\subset O_N^+$ and $U_N\subset U_N^+$ are liberations, we must prove that we have isomorphisms as follows:
$$C(B_N)=C(O_N)\Big/\Big<\xi\in Fix(u)\Big>$$
$$C(C_N)=C(U_N)\Big/\Big<\xi\in Fix(u)\Big>$$

But these isomorphisms are both clear from the formulae of $B_N,C_N$ in (1).
\end{proof}

The above might seem a bit puzzling, because you might have heard from Lie algebras that the ``basic'' compact groups can be of type ABCD, and there was no mention there, in that Lie algebra theory and classification, of bistochastic groups as above. 

\bigskip

Good point, and in answer, indeed, when talking from an abstract algebra perspective, meaning groups and quantum groups $G$ taken up to isomorphism, and ignoring the embeddings $G\subset U_N$, the bistochastic groups and quantum groups are not really ``new'', because, following Raum \cite{rau} and related papers, we have the following result:

\index{Fourier matrix}
\index{discrete Fourier transform}

\begin{theorem}
We have isomorphisms as follows:
\begin{enumerate}
\item $B_N\simeq O_{N-1}$.

\item $B_N^+\simeq O_{N-1}^+$.

\item $C_N\simeq U_{N-1}$.

\item $C_N^+\simeq U_{N-1}^+$.
\end{enumerate}
\end{theorem}

\begin{proof}
Let us pick indeed a matrix $F\in U_N$ satisfying the following condition, where $\xi$ is the all-one vector, and $e_0=(1,0,\ldots,0)$ is the first vector of the standard basis of $\mathbb C^N$, written with indices $0,1,\ldots,N-1$, as usual in discrete Fourier analysis:
$$Fe_0=\frac{1}{\sqrt{N}}\xi$$

Such matrices exist of course, the basic example being the Fourier matrix:
$$F_N=\frac{1}{\sqrt{N}}(w^{ij})_{ij}\quad,\quad w=e^{2\pi i/N}$$

We have then the following computation, for any corepresentation $u$:
\begin{eqnarray*}
u\xi=\xi
&\iff&uFe_0=Fe_0\\
&\iff&F^*uFe_0=e_0\\
&\iff&F^*uF=diag(1,w)
\end{eqnarray*}

Thus we have an isomorphism given by $w_{ij}\to(F^*uF)_{ij}$, as desired.
\end{proof}

Summarizing, we are deviating here a bit from the standard viewpoint on groups and quantum groups, coming from Lie theory, somehow by constructing brand new objects out of the old ones, a bit like a magician pulls out a rabbit from a hat. This being said, the bistochastic quantum groups $B_N,C_N$ and $B_N^+,C_N^+$ will be fundamental objects for us, and will appear on numerous occasions, in the remainder of this book.

\bigskip

More on this in a moment, in connection with easiness questions. As a main motivation, however, we have the fact that the bistochastic matrices, and especially the unitary ones, $U\in C_N$, are cult objects in advanced matrix analysis. So, Lie theory is not everything, we all do mistakes, and ignoring interesting objects like $B_N,C_N$ is one of the flaws of Lie theory, and everything that will follow will be certainly worth developing.

\bigskip

Getting now to the real thing, the point is that our newly constructed bistochastic groups and quantum groups are all easy, and following \cite{bsp}, \cite{tw2}, we have:

\index{singletons and pairings}

\begin{theorem}
The classical and quantum bistochastic groups are all easy, with the quantum groups on the left corresponding to the categories on the right,
$$\xymatrix@R=51pt@C=50pt{
B_N^+\ar[r]&C_N^+\\
B_N\ar[u]\ar[r]&C_N\ar[u]}
\qquad\xymatrix@R=25pt@C=50pt{\\ :\\}\qquad
\xymatrix@R=18.5mm@C=16.2mm{
NC_{12}\ar[d]&\mathcal{NC}_{12}\ar[l]\ar[d]\\
P_{12}&\mathcal P_{12}\ar[l]}$$
where the symbol $12$ stands for ``category of singletons and pairings''.
\end{theorem}

\begin{proof}
This comes from the fact that the all-one vector $\xi$ used in the constructions in Theorem 5.12 is the vector associated to the singleton partition:
$$\xi=T_|$$

Indeed, we obtain that $B_N,C_N,B_N^+,C_N^+$ are inded easy, appearing from the categories of partitions for $O_N,U_N,O_N^+,U_N^+$, by adding singletons. Thus, we get the result.
\end{proof}

Now that we have easiness, we can do many things with it. Let us first discuss laws of characters. In order to formulate our results, we will need:

\begin{definition}
The shifted version of a measure $\mu_t$, depending on a parameter $t>0$, is the law of the variable $t+X$, with $X$ following  the law $\mu_t$. We denote by
$$\xymatrix@R=51pt@C=50pt{
\gamma_t\ar@{-}[d]&\Gamma_t\ar@{-}[l]\ar@{-}[d]\\
g_t&G_t\ar@{-}[l]}
\qquad\xymatrix@R=25pt@C=50pt{\\ \to\\}\qquad
\xymatrix@R=51pt@C=50pt{
\sigma_t\ar@{-}[d]&\Sigma_t\ar@{-}[l]\ar@{-}[d]\\
s_t&S_t\ar@{-}[l]}$$
the shifted versions of the normal, complex normal, semicircle and circular laws.
\end{definition}

With this definition in hand, we can now formulate our character result, as follows:

\begin{theorem}
The asymptotic laws of truncated characters for the bistochastic quantum groups are the shifted versions of the normal and semicircle laws:
$$\xymatrix@R=51pt@C=50pt{
B_N^+\ar[r]&C_N^+\\
B_N\ar[u]\ar[r]&C_N\ar[u]}
\qquad\xymatrix@R=25pt@C=50pt{\\ :\\}\qquad
\xymatrix@R=51pt@C=50pt{
\sigma_t\ar@{-}[d]&\Sigma_t\ar@{-}[l]\ar@{-}[d]\\
s_t&S_t\ar@{-}[l]}$$
Moreover, these laws form convolution semigroups, in Bercovici-Pata bijection.
\end{theorem}

\begin{proof}
This can be done by following the proof for $O_N,U_N,O_N^+,U_N^+$ from chapter 4, and performing modifications where needed, as follows:

\medskip

(1) As before for $O_N,U_N,O_N^+,U_N^+$, for any of our quantum groups $B_N,C_N,B_N^+,C_N^+$ we have the following formula, for the moments of the truncated characters:
$$\int_{G_N}(u_{11}+\ldots +u_{ss})^k
=Tr(W_{kN}G_{ks})$$

(2) Also as before for $O_N,U_N,O_N^+,U_N^+$, for any of our quantum groups $B_N,C_N,B_N^+,C_N^+$, we deduce from this that we have the following asymptotic formula:
$$\lim_{N\to\infty}\int_{G_N}\chi_t^k=\sum_{\pi\in D(k)}t^{|\pi|}$$

(3) In order to finish now the computation, and prove the first assertion, in the real case, assume that we have variables $X,Y$ following the classical/free laws for the groups/quantum groups $O_N,B_N$ or $O_N^+,B_N^+$. By using the above formula, we obtain:
\begin{eqnarray*}
\mathbb E(Y^k)
&=&\sum\left\{t^{|\pi|}\Big|\text{$\pi\in P(k)$ or $NC(k)$, consisting of
singletons and pairings}\right\}\\
&=&\sum_{r=0}^k\binom{k}{r}t^r\sum\left\{t^{|\pi|}\Big|\text{$\pi\in
P(k-r)$ or $NC(k-r)$, consisting of pairings}\right\}\\
&=&\sum_{r=0}^k\binom{k}{r}t^r\mathbb E(X^{k-r})\\
&=&\mathbb E((t+X)^k)
\end{eqnarray*}

Thus the law of $Y$ is the same as the law of $X$ shifted by $t$, as claimed.

\medskip

(4) In what regards now the semigroup assertion and the Bercovici-Pata bijection assertion, these are best proved via cumulants. Indeed, if we denote by $m_t$ the asymptotic law of the variable $\chi_t$, for the group $B_N$, we know from the above that we have:
$$M_k(m_t)=\sum_{\pi\in P_{12}(k)}t^{|\pi|}$$

But this shows that the classical cumulants of our measure $m_t$ are given by:
$$k_n(m_t)=\begin{cases}
t&{\rm if}\ n=1,2\\
0&{\rm if}\ n\geq 3
\end{cases}$$

Similarly, if we denote by $\mu_t$ the asymptotic law of the variable $\chi_t$, for the quantum group $B_N^+$, we know from the above that we have:
$$M_k(\mu_t)=\sum_{\pi\in NC_{12}(k)}t^{|\pi|}$$

But this shows that the free cumulants of our measure $\mu_t$ are given by:
$$\kappa_n(\mu_t)=\begin{cases}
t&{\rm if}\ n=1,2\\
0&{\rm if}\ n\geq 3
\end{cases}$$

(5) Now since the classical and free cumulants computed above are linear in $t$, the measures $\{m_t\}$ form a convolution semigroup, and the measures $\{\mu_t\}$ form a free convolution semigroup. Moreover, since the classical cumulants of $m_t$ equal the free cumulants of $\mu_t$, these semigroups are in Bercovici-Pata bijection, as claimed.

\medskip

(6) In the complex case, involving $U_N,C_N$ and $U_N^+,C_N^+$, the computation is similar. Finally, regarding the semigroup assertions and the Bercovici-Pata bijection, these follow as well from the general moment formula above, via standard cumulant theory.
\end{proof}

Regarding now the Gram determinants, things here are more technical, and the most convenient is to use Theorem 5.13. Let us just record here the results in the orthogonal case, from \cite{bcu}. For the bistochastic group $B_N$, the result is as follows:

\begin{theorem}
For the bistochastic group $B_N$ we have
$$\det(G_{kN})=N^{a_k}\prod_{|\lambda|\leq k/2}f_N(\lambda)^{\binom{k}{2|\lambda|}f^{2\lambda}}$$
where $a_k=\sum_{\pi\in P(k)}(2|\pi|-k)$, and $f_N(\lambda)=\prod_{(i,j)\in\lambda}(N+2j-i-2)$.
\end{theorem}

\begin{proof}
We recall from chapter 3 that the Gram determinant for $O_N$ is given by the following formula, where $f_N'(\lambda)=\prod_{(i,j)\in\lambda}(N+2j-i-1)$:
$${\det}'(G_{kN})=\prod_{|\lambda|=k/2}f_N'(\lambda)^{f^{2\lambda}}$$

On the other hand, by Theorem 5.13 we have an isomorphism $B_N\simeq O_{N-1}$, given by $u=v+1$, where $u,v$ are the fundamental representations of $B_N,O_{N-1}$. But this gives:
$$Fix(u^{\otimes k})
=Fix\left((v+1)^{\otimes k}\right)
=Fix\left(\sum_{r=0}^k\binom{k}{r}v^{\otimes r}\right)$$

Now if we denote by ${\rm det}',f'$ the objects computed for $O_N$, as above, we obtain:
$$\det(G_{kN})
=N^{a_k}\prod_{r=1}^k{\det}'(G_{r,N-1})^{\binom{k}{r}}
=N^{a_k}\prod_{r=1}^k\left(\prod_{|\lambda|=r/2}f'_{N-1}(\lambda)^{f^{2\lambda}}\right)^{\binom{k}{r}}$$

Thus, we are led to the formula in the statement.
\end{proof}

The same method works for $B_N^+$, the result here being as follows:

\begin{theorem}
For the free bistochastic group $B_n^+$ we have
$$\det(G_{kN})=N^{a_k}\prod_{r=1}^{[k/2]}P_r(N-1)^{\sum_{l=1}^{[k/2]}\binom{k}{2l}d_{lr}}$$
where $a_k=\sum_{\pi\in \mathcal P(k)}(2|\pi|-k)$, as before.
\end{theorem}

\begin{proof}
We recall from chapter 3 that the Gram determinant for $O_N^+$ is given by the following formula, where $P_r(N)$ are the Chebycheff polynomials:
$$\det(G_{kN})=\prod_{r=1}^{[k/2]}P_r(N)^{d_{k/2,r}}$$

On the other hand, by Theorem 5.13 we have an isomorphism $B_N^+\simeq O_{N-1}^+$, given by $u=v+1$, where $u,v$ are the fundamental representations of $B_N^+,O_{N-1}^+$. But this gives:
$$Fix(u^{\otimes k})
=Fix\left((v+1)^{\otimes k}\right)
=Fix\left(\sum_{r=0}^k\binom{k}{r}v^{\otimes r}\right)$$

Now if we denote by ${\det}'$ the determinant computed for $O_N^+$, as above, we obtain:
$$\det(G_{kN})
=N^{a_k}\prod_{l=1}^{[k/2]}{\det}'(G_{2l,N-1})^{\binom{k}{2l}}
=N^{a_k}\prod_{l=1}^{[k/2]}\left(\prod_{r=1}^lP_r(N-1)^{d_{lr}}\right)^{\binom{k}{2l}}$$

Thus, we are led to the formula in the statement.
\end{proof}

There are a number of further things that can be said about the bistochastic groups and quantum groups, in relation with advanced representation theory, growth exponents, De Finetti theorems, and more. As a main result here, we have:

\begin{theorem}
Let $(x_i)_{i\in \mathbb N}$ be a $\{G_N\}$-invariant sequence of self-adjoint random variables in a noncommutative probability space $(A,\varphi)$, which generates A. Then there is a subalgebra $C\subset A$ and a $\varphi$-preserving conditional expectation $E:A\to C$ such that:
\begin{enumerate}
\item If $G_N=B_N$, then $(x_i)_{i\in \mathbb N}$ are conditionally independent, and have Gaussian distributions with common mean and variance, given $C$.

\item If $G_N=B_N^+$, then $(x_i)_{i \in \mathbb N}$ form a $C$-valued free semicircular family with common mean and variance.
\end{enumerate}
\end{theorem}

\begin{proof}
This is a result from \cite{ez3}, whose proof is very similar to the results for $G_N=O_N,O_N^+$ discussed in chapter 4, and for full details here, we refer to \cite{ez3}.
\end{proof}

Finally, let us mention the following question, which is open:

\begin{question}
What is the free analogue of the Idel-Wolf theorem?
\end{question}

To be more precise here, let us recall from Theorem 5.11 that any unitary matrix $U\in U_N$ can be decomposed, at least in theory, as follows, with $V\in C_N$ being a unitary bistochastic matrix, and with $L,R\in\mathbb T^N$ being unitary diagonal matrices:
$$U=LVR$$

Generaly speaking, having a free analogue of such things looks quite complicated. However, we can look for some weaker statements. For instance, the Idel-Wolf theorem implies the following generation formula, for closed subgroups of $U_N$, which in fact is something elementary, which can be proved as well starting from definitions:
$$<C_N,\mathbb T_N>=U_N$$

In the quantum case now, as a first step towards stating and proving Idel-Wolf type theorems, we have the following formula to be proved, for closed subgroups of $U_N^+$:
$$<C_N^+,\mathbb T_N^+>=U_N^+$$

As a first observation here, since the presence of either $C_N^+$ or $\mathbb T_N^+$, which are both free, in the generation operation normally guarantees the freeness of the resulting quantum group, we have as well two stronger versions of this formula, as follows:
$$<C_N^+,\mathbb T_N>=U_N^+\quad,\quad 
<C_N,\mathbb T_N^+>=U_N^+$$

In any case, the discussion here is quite elementary, using the known generation results for groups and quantum groups, and the next step is that of going beyond this, with something closer to the Idel-Wolf theorem. Now in order to get started here, the first observation is that the Idel-Wolf theorem tells us that the following map is surjective:
$$\mathbb T_N\times C_N\times\mathbb T_N\to U_N\quad,\quad (L,V,R)\to LVR$$

But with this in hand, we can transpose everything, and we are led to a functional analytic formulation of the Idel-Wolf theorem, involving a map as follows:
$$C(U_N)\to C(\mathbb T_N)\otimes C(C_N)\otimes C(\mathbb T_N)$$

For quantum groups, however, things are quite tricky, due to noncommutativity, and formulating something which looks plausible is not an easy task, requiring some good imagination. We will be back to this in chapter 7 below, when talking half-liberation, with the amount of noncommutativity there being substantially lower.

\section*{5c. The continuous cube} 

The results obtained so far in this chapter are quite interesting, and look rather fundamental, and job for us now to see how this material fits with the general theory developed in chapters 1-4, and notably with the standard cube considered there, namely:
$$
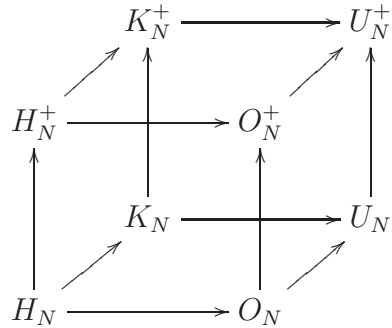
\xymatrix@R=18pt@C=18pt{
&K_N^+\ar[rr]&&U_N^+\\
H_N^+\ar[rr]\ar[ur]&&O_N^+\ar[ur]\\
&K_N\ar[rr]\ar[uu]&&U_N\ar[uu]\\
H_N\ar[uu]\ar[ur]\ar[rr]&&O_N\ar[uu]\ar[ur]
}$$

The first thought goes to an improved, 4D cube, obtained by adding a new dimension to those that we already have, in relation with the notion of bistochasticity. However, there is a problem with all this, which is actually a bit hard to explain and understand, so let us start modestly. We first have the following result, dealing exclusively with the continuous objects that we have, namely those on the right face of the above standard cube, and the bistochastic quantum groups that we introduced in this chapter:

\begin{theorem}
The basic orthogonal and unitary quantum groups and their bistochastic versions are all easy, and they form a diagram as follows,
$$
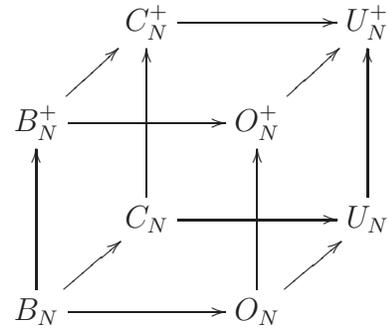
\xymatrix@R=18pt@C=18pt{
&C_N^+\ar[rr]&&U_N^+\\
B_N^+\ar[rr]\ar[ur]&&O_N^+\ar[ur]\\
&C_N\ar[rr]\ar[uu]&&U_N\ar[uu]\\
B_N\ar[uu]\ar[ur]\ar[rr]&&O_N\ar[uu]\ar[ur]
}$$
which is an intersection and easy generation diagram, in the sense that any subsquare $P\subset Q,R\subset S$ of this diagram satisfies $Q\cap R=P$, $\{Q,R\}=S$.
\end{theorem}

\begin{proof}
We know that the quantum groups in the statement are indeed easy, the corresponding categories of partitions being as follows:
$$
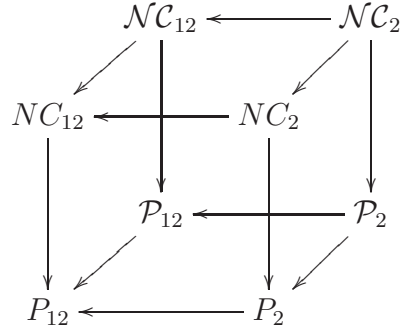
\xymatrix@R=20pt@C8pt{
&\mathcal{NC}_{12}\ar[dl]\ar[dd]&&\mathcal {NC}_2\ar[dl]\ar[ll]\ar[dd]\\
NC_{12}\ar[dd]&&NC_2\ar[dd]\ar[ll]\\
&\mathcal P_{12}\ar[dl]&&\mathcal P_2\ar[dl]\ar[ll]\\
P_{12}&&P_2\ar[ll]
}$$

Now since both this diagram and the one the statement are intersection diagrams, the quantum groups form an intersection and easy generation diagram, as stated.
\end{proof}

The above result is quite nice, and here is a related one, of the same nature:

\begin{theorem}
The basic orthogonal quantum groups are all easy, with
$$
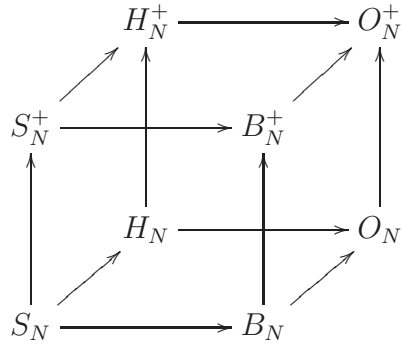
\xymatrix@R=20pt@C=20pt{
&H_N^+\ar[rr]&&O_N^+\\
S_N^+\ar[rr]\ar[ur]&&B_N^+\ar[ur]\\
&H_N\ar[rr]\ar[uu]&&O_N\ar[uu]\\
S_N\ar[uu]\ar[ur]\ar[rr]&&B_N\ar[uu]\ar[ur]
}$$
being an intersection and easy generation diagram.
\end{theorem}

\begin{proof}
We know that the quantum groups in the statement are indeed easy, the corresponding categories of partitions being as follows:
$$
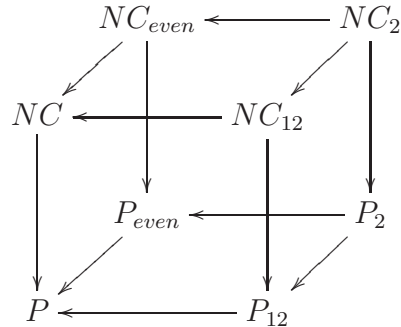
\xymatrix@R=20pt@C6pt{
&NC_{even}\ar[dl]\ar[dd]&&NC_2\ar[dl]\ar[ll]\ar[dd]\\
NC\ar[dd]&&NC_{12}\ar[dd]\ar[ll]\\
&P_{even}\ar[dl]&&P_2\ar[dl]\ar[ll]\\
P&&P_{12}\ar[ll]
}$$

Now since both this diagram and the one the statement are intersection diagrams, the quantum groups form an intersection and easy generation diagram, as stated.
\end{proof}

All this suggests merging everything into a nice 4D cube. Unfortunately, this is not possible, with the problem coming from the following negative result:

\begin{proposition}
The cube from Theorem 5.22, and its unitary analogue
$$
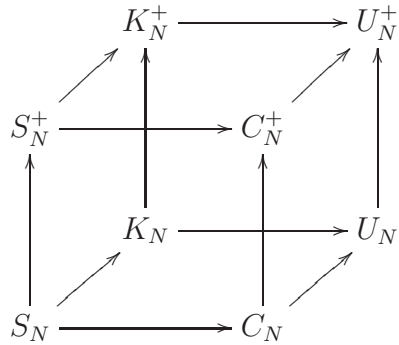
\xymatrix@R=20pt@C=20pt{
&K_N^+\ar[rr]&&U_N^+\\
S_N^+\ar[rr]\ar[ur]&&C_N^+\ar[ur]\\
&K_N\ar[rr]\ar[uu]&&U_N\ar[uu]\\
S_N\ar[uu]\ar[ur]\ar[rr]&&C_N\ar[uu]\ar[ur]
}$$
cannot be merged, without degeneration, into a $4$-dimensional cubic diagram.
\end{proposition}

\begin{proof}
All this is a bit philosophical, with the problem coming from the ``taking the bistochastic version'' operation, and more specifically, from the following equalities:
$$H_N\cap C_N=K_N\cap C_N=S_N$$

Indeed, these equalities do hold, and so the 3D cube obtained by merging the classical faces of the orthogonal and unitary cubes is something degenerate, as follows:
$$
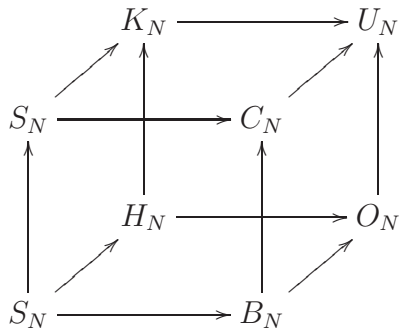
\xymatrix@R=20pt@C=20pt{
&K_N\ar[rr]&&U_N\\
S_N\ar[rr]\ar[ur]&&C_N\ar[ur]\\
&H_N\ar[rr]\ar[uu]&&O_N\ar[uu]\\
S_N\ar[uu]\ar[ur]\ar[rr]&&B_N\ar[uu]\ar[ur]
}$$

Thus, the 4D cube, having this 3D cube as one of its faces, is degenerate too.
\end{proof}

Summarizing, when positioning ourselves at $U_N^+$, we have 4 natural directions to be followed, namely taking the classical, discrete, real and bistochastic versions. And the problem is that, while the first three operations are ``good'', the fourth one is ``bad''. This is good to know, and we will come back to this, later in this book.

\section*{5d. Further constructions}

We would like to discuss now certain versions of the classical and free bistochastic groups $B_N,C_N,B_N^+,C_N^+$, which are quite interesting, for various reasons. First, it is possible to talk about groups and quantum groups $B_N^I,C_N^I,B_N^{I+},C_N^{I+}$ defined by using other vectors $\xi=\xi_I$, which are not easy in general, the result being as follows:

\begin{theorem}
We have a closed subgroup $C_N^{I+}\subset U_N^+$, defined via the formula
$$C(C_N^{I+})=C(U_N^+)\Big/\left<u\xi_I=\xi_I\right>$$
with the vector $\xi_I$ being as follows, with $I\subset\{1,\ldots,N\}$ being a subset:
$$\xi_I=\frac{1}{\sqrt{|I|}}(\delta_{i\in I})_i$$
Moreover, we can talk about quantum groups $B_N^I,C_N^I,B_N^{I+}$ too, in a similar way.
\end{theorem}

\begin{proof}
We must check Woronowicz' axioms, and the proof goes as follows:

\medskip

(1) Let us set $U_{ij}=\sum_ku_{ik}\otimes u_{kj}$. We have then the following computation:
\begin{eqnarray*}
(U\xi_I)_i
&=&\frac{1}{\sqrt{|I|}}\sum_{j\in I}U_{ij}\\
&=&\frac{1}{\sqrt{|I|}}\sum_{j\in I}\sum_ku_{ik}\otimes u_{kj}\\
&=&\sum_ku_{ik}\otimes(u\xi_I)_k\\
&=&\sum_ku_{ik}\otimes(\xi_I)_k\\
&=&\frac{1}{\sqrt{|I|}}\sum_{k\in I}u_{ik}\otimes1\\
&=&(u\xi_I)_i\otimes1\\
&=&(\xi_I)_i\otimes1
\end{eqnarray*}

Thus we can define indeed a comultiplication map, by $\Delta(u_{ij})=U_{ij}$.

\medskip

(2) In order to construct the counit map, $\varepsilon(u_{ij})=\delta_{ij}$, we must prove that the identity matrix $1=(\delta_{ij})_{ij}$ satisfies $1\xi_I=\xi_I$. But this is clear.

\medskip

(3) In order to construct the antipode, $S(u_{ij})=u_{ji}^*$, we must prove that the adjoint matrix $u^*=(u_{ji}^*)_{ij}$ satisfies $u^*\xi_I=\xi_I$. But this is clear from $u\xi_I=\xi_I$.

\medskip

(4) Finally, we can talk about quantum groups $B_N^I,C_N^I,B_N^{I+}$ in a similar way, with these appearing as quantum subgroups of $C_N^{I+}$. 
\end{proof}

Summarizing, for any index set $I\subset\{1,\ldots,N\}$, we have a cube as follows:
$$
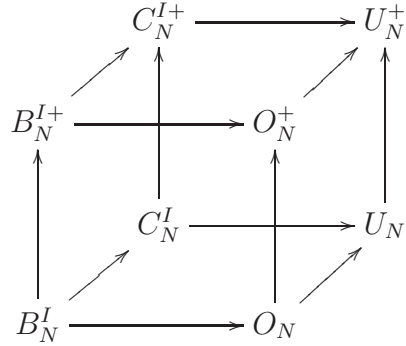
\xymatrix@R=18pt@C=18pt{
&C_N^{I+}\ar[rr]&&U_N^+\\
B_N^{I+}\ar[rr]\ar[ur]&&O_N^+\ar[ur]\\
&C_N^I\ar[rr]\ar[uu]&&U_N\ar[uu]\\
B_N^I\ar[uu]\ar[ur]\ar[rr]&&O_N\ar[uu]\ar[ur]
}$$

We should mention that the above quantum groups are of key importance in connection with the notion of ``affine homogeneous space'', which is something generalizing the usual spheres, and requires an index set $I\subset\{1,\ldots,N\}$, as above. We will be back to this in chapter 16 below, with a brief introduction to the affine homogeneous spaces.

\bigskip

As a final topic for this chapter, let us discuss the half-liberation operation. In connection with the bistochastic groups, we have the following negative result:

\begin{proposition}
The half-classical versions of $B_N^+,C_N^+$ are given by:
$$B_N^+\cap O_N^*=B_N\quad,\quad 
C_N^+\cap U_N^*=C_N$$
In other words, the half-classical versions collapse to the classical versions.
\end{proposition}

\begin{proof}
This follows indeed from Tannakian duality, by using the fact that when capping the half-classical crossing with 2 singletons, we obtain the classical crossing. Equivalently, this follows as well from definitions, via some standard computations.
\end{proof}

However, it is possible to construct intermediate objects for $B_N\subset B_N^+$ and $C_N\subset C_N^+$, by going beyond easiness. Let us recall from Theorem 5.13 that we have an isomorphism as follows, whenever $F\in O_N$ satisfies $Fe_0=\frac{1}{\sqrt{N}}\xi$, where $\xi$ is the all-one vector:
$$C(O_{N-1}^+)\to C(B_N^+)\quad,\quad w_{ij}\to(F^*uF)_{ij}$$

Here, and in what follows, we use indices $i,j=0,1,\ldots,N-1$ for the $N\times N$ compact quantum groups, and indices $i,j=1,\ldots,N-1$ for their $(N-1)\times(N-1)$ subgroups. But with this, we can construct intermediate objects for $B_N\subset B_N^+$, as follows:

\begin{proposition}
Assuming that $F\in O_N$ satisfies $Fe_0=\frac{1}{\sqrt{N}}\xi$, we have inclusions as follows, with the intermediate quantum group $B_F^\circ$ being not easy,
$$B_N\subset B_F^\circ\subset B_N^+$$
obtained by taking the image of the inclusions $O_{N-1}\subset O_{N-1}^*\subset O_{N-1}^+$, via the above isomorphism $O_{N-1}^+\simeq B_N^+$ induced by $F$.
\end{proposition}

\begin{proof}
The fact that we have inclusions as in the statement follows from the above isomorphism, which produces a diagram as follows:
$$\xymatrix@R=15mm@C=20mm{
O_{N-1}\ar[r]&O_{N-1}^*\ar[r]&O_{N-1}^+\\
B_N\ar[r]\ar@{=}[u]&B_F^\circ\ar[r]\ar@{=}[u]&B_N^+\ar@{=}[u]
}$$

To be more precise, the quantum group $B_F^\circ$ from the bottom is by definition the image of the quantum group $O_{N-1}^*$ from the top. Since we know that $B_N\subset B_N^+$ is maximal in the easy setting, this new quantum group $B_F^\circ$ is not easy, as claimed.
\end{proof}

Observe that the relations $abc=cba$ do not hold for all the entries of the modified fundamental corepresentation $v=F^*uF$, due to the fact that we have $v_{00}=1$, and that the relations $ab1=1ba$ corresponds to the commutativity. We have in fact:

\begin{proposition}
The quantum group $B_F^\circ\subset B_N^+$  appears via the relations
$$R^{\otimes 3}T_{\slash\hskip-1.2mm|\hskip-1.2mm\backslash}R^{*\otimes 3}\in End(u^{\otimes 3})$$
where $R=FP$, with $R$ being the projection onto $span(e_1,\ldots,e_{N-1})$.
\end{proposition}

\begin{proof}
With $F^*uF=diag(1,w)$ as above, the relations defining $B_F^\circ$ are:
\begin{eqnarray*}
T_{\slash\hskip-1.2mm|\hskip-1.2mm\backslash}\in End(w^{\otimes 3})
&\iff&P^{\otimes 3}T_{\slash\hskip-1.2mm|\hskip-1.2mm\backslash}P^{\otimes 3}\in End((F^*uF)^{\otimes 3})\\
&\iff&(FP)^{\otimes 3}T_{\slash\hskip-1.2mm|\hskip-1.2mm\backslash}(FP)^{*\otimes 3}\in End(u^{\otimes 3})
\end{eqnarray*}

Thus, we obtain the formula in the statement.
\end{proof}

Now observe that, due to the conditions $F\in O_N$ and $Fe_0=\frac{1}{\sqrt{N}}\xi$, the linear map associated to $R=FP$ maps $e_0\to0\to0$ and $e_i\to e_i\to f_i$, where $\{f_1,\ldots,f_{N-1}\}$ is a certain orthonormal basis of $\xi^\perp$. Thus $R=FP$ must be a partial isometry $e_0^\perp\to\xi^\perp$. In view of this observation, we can further process the above result, as follows:

\begin{proposition}
The quantum groups $B_F^\circ\subset B_N^+$ with $F\in O_N$, $Fe_0=\frac{1}{\sqrt{N}}\xi$ all coincide, and appear via the relations $T\in End(u^{\otimes 3})$, where
\begin{eqnarray*}
T(e_i\otimes e_j\otimes e_k)
&=&e_k\otimes e_j\otimes e_i-(e_k\otimes e_j\otimes\xi'+e_k\otimes\xi'\otimes e_i+\xi'\otimes e_j\otimes e_i)\\
&&+(e_k\otimes\xi'\otimes\xi'+\xi'\otimes e_j\otimes\xi'+\xi'\otimes\xi'\otimes e_i)-\xi'\otimes\xi'\otimes\xi'
\end{eqnarray*}
with $\xi'=\frac{1}{N}\xi$, and with $\xi$ being as usual the all-one vector.
\end{proposition}

\begin{proof}
The linear map $R^{\otimes 3}T_{\slash\hskip-1.2mm|\hskip-1.2mm\backslash}R^{*\otimes 3}$ from Proposition 5.27 acts as follows:
\begin{eqnarray*}
R^{\otimes 3}T_{\slash\hskip-1.2mm|\hskip-1.2mm\backslash}R^{*\otimes 3}(e_i\otimes e_j\otimes e_k)
&=&R^{\otimes 3}\sum_{abc}R_{ia}R_{jb}R_{kc}\ e_c\otimes e_b\otimes e_a\\
&=&\sum_{abcpqr}R_{ia}R_{jb}R_{kc}R_{pc}R_{qb}R_{ca}\ e_p\otimes e_q\otimes e_r\\
&=&\sum_{pqr}(RR^t)_{ir}(RR^t)_{jq}(RR^t)_{kp}\ e_p\otimes e_q\otimes e_r
\end{eqnarray*}

On the other hand, since $R=FP$ must be a partial isometry $e_0^\perp\to\xi^\perp$, we have:
$$RR^*=1-Proj(\xi)\quad,\quad(RR^*)_{ij}=\delta_{ij}-\frac{1}{N}$$

We conclude that the map $R^{\otimes 3}T_{\slash\hskip-1.2mm|\hskip-1.2mm\backslash}R^{*\otimes 3}$ is given by:
$$e_i\otimes e_j\otimes e_k\to\sum_{pqr}\left(\delta_{ir}-\frac{1}{N}\right)\left(\delta_{jq}-\frac{1}{N}\right)\left(\delta_{kp}-\frac{1}{N}\right)\ e_p\otimes e_q\otimes e_r$$

Now by developing, we obtain the formula in the statement.
\end{proof}

An even better statement, which looks more conceptual, is as follows:

\begin{theorem}
The quantum group $B_N^\circ\subset B_N^+$ constructed above, which equals the various quantum groups $B_F^\circ$, appears via the relations $T\in End(u^{\otimes 3})$, where
$$T=T_{\slash\hskip-1.2mm|\hskip-1.2mm\backslash}-
(T_{{\ }^\cdot_\cdot\slash\hskip-1.4mm\backslash}+T_{{\ }^{{\ }^\cdot}_{{\ }_\cdot}\hskip-1.5mm\slash\hskip-1.4mm\backslash}+T_{\slash\hskip-1.4mm\backslash\!\!\!{\ }^\cdot_\cdot})+
(T_{{\ }_{..}\!\!\backslash^{\!\cdot\cdot}}+T_{{\ }^\cdot_\cdot|\!\!\!{\ }^\cdot_\cdot}+T_{{\ }^{..}\!\slash_{\!\cdot\cdot}})-
T_{{\ }^{\cdot\cdot\cdot}_{\cdot\cdot\cdot}}$$
with the convention that the various dots represent singletons.
\end{theorem}

\begin{proof}
This follows indeed from the formula in Proposition 5.28, because the 8 terms there correspond to the 8 partitions in the statement. 
\end{proof}

Observe that we can in fact write an even more compact formula for the linear map $T$ in the above result, in terms of the M\"obius function of $P_{12}$, as follows:
$$T=\sum_{\pi\leq\slash\hskip-1.2mm|\hskip-1.2mm\backslash}\mu(\pi)T_\pi$$

We can perform a similar construction in the unitary case, as follows: 

\begin{theorem}
The image of $U_{N-1}^*$ via the isomorphism $U_{N-1}^+\simeq C_N^+$ is an intermediate quantum group $C_N\subset C_N^\circ\subset C_N^+$, which
appears via the relations
$$T\in End(u^{\otimes k})$$
where $T$ is the linear map from Theorem 5.29, and where $k\in\{\circ\circ\circ,\circ\circ\bullet,\ldots,\bullet\bullet\bullet\}$ ranges over all colored integers of length $3$.
\end{theorem}

\begin{proof}
This follows indeed by proceeding as before in the real case, and replacing where needed the tensor powers $u^{\otimes 3}$ by the colored tensor powers $u^{\otimes k}$, as above.
\end{proof}

We will back to this, with some further examples of bistochastic quantum groups, in chapter 7 below, when systematically discussing the half-liberation operation.

\section*{5e. Exercises} 

Before anything, as already mentioned on several occasions, with bistochastic matrices, groups and quantum groups, we are now into applied linear algebra, of quite exciting type. As a first exercise, recommended, and that you can fully follow or not, we have:

\begin{exercise}
Read and fully understand the proof of the Idel-Wolf theorem, with all the needed symplectic geometry preliminaries, Arnold and so on, and then:
\begin{enumerate}
\item Either work a bit on quantum Idel-Wolf, which is an excellent research topic.

\item Or ditch this book and stay with Idel-Wolf, which is first-class mathematics too.
\end{enumerate}
\end{exercise}

As a second research question now, which is equally interesting, and that we struggled a bit with, in the above, without clear conclusions, we have:

\begin{exercise}
Find a clever way of adding to the continuous cube
$$
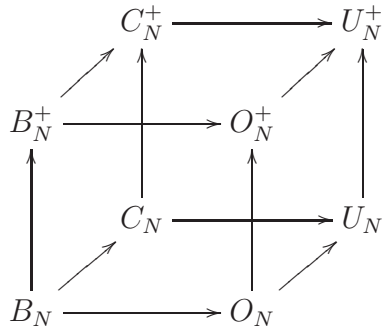
\xymatrix@R=18pt@C=18pt{
&C_N^+\ar[rr]&&U_N^+\\
B_N^+\ar[rr]\ar[ur]&&O_N^+\ar[ur]\\
&C_N\ar[rr]\ar[uu]&&U_N\ar[uu]\\
B_N\ar[uu]\ar[ur]\ar[rr]&&O_N\ar[uu]\ar[ur]
}$$
various reflection groups $H_N^s,H_N^{s+}$. Think about half-liberations too.
\end{exercise}

This is a very good question, because the further examples of easy quantum groups will start to accumulate, in what comes next, and we will soon lose control of everything, geometrically speaking, and the origins of this control losing are right here, in this chapter 5, in relation with the above question. As a more standard exercise now, we have:

\begin{exercise}
Work out laws of characters, Gram determinants and De Finetti theorems for $C_N,C_N^+$, in analogy with what we did in the above for $B_N,B_N^+$.
\end{exercise}

Finally, the material at end of this chapter, going beyond easiness, with various tricky constructions, can certainly lead to many interesting exercises, of research level, and you can try some of these, after inventing them of course. But we will be back to this. 

\chapter{The uniform case}

\section*{6a. Uniform groups}

We already have enough examples of easy quantum groups, and knowledge of the subject, for starting some classification work. However, the general classification question for the easy quantum groups is something quite difficult. In what follows we will cut a bit from complexity, by adding some extra axioms, chosen as ``natural'' as possible. This will lead to some useful classification results, and only after having these results, we will think of relaxing some of our extra axioms, and see what we get.  

\bigskip

So, this will be our plan, and for the philosophy, when dealing with difficult problems, modesty is your main weapon. Getting to work now, a first natural axiom comes from:

\index{uniform quantum group}
\index{removing blocks}

\begin{theorem}
For an easy quantum group $G=(G_N)$, coming from a category of partitions $D\subset P$, the following conditions are equivalent:
\begin{enumerate}
\item $G_{N-1}=G_N\cap U_{N-1}^+$, via the embedding $U_{N-1}^+\subset U_N^+$ given by $u\to diag(u,1)$.

\item $G_{N-1}=G_N\cap U_{N-1}^+$, via the $N$ possible diagonal embeddings $U_{N-1}^+\subset U_N^+$.

\item $D$ is stable under the operation which consists in removing blocks.
\end{enumerate}
\end{theorem}

\begin{proof}
We use the general easiness theory from chapters 1-4:

\medskip

$(1)\iff(2)$ This is something standard, coming from the inclusion $S_N\subset G_N$, which makes everything $S_N$-invariant. The result follows as well from the proof of $(1)\iff(3)$ below, which can be converted into a proof of $(2)\iff(3)$, in the obvious way.

\medskip

$(1)\iff(3)$ Given a subgroup $K\subset U_{N-1}^+$, with fundamental corepresentation $u$, consider the $N\times N$ matrix $v=diag(u,1)$. Our claim is that for any $\pi\in P(k)$ we have:
$$\xi_\pi\in Fix(v^{\otimes k})\iff\xi_{\pi'}\in Fix(v^{\otimes k'}),\,\forall\pi'\in P(k'),\pi'\subset\pi$$

In order to prove this, we must study the condition on the left. We have:
\begin{eqnarray*}
\xi_\pi\in Fix(v^{\otimes k})
&\iff&(v^{\otimes k}\xi_\pi)_{i_1\ldots i_k}=(\xi_\pi)_{i_1\ldots i_k},\forall i\\
&\iff&\sum_j(v^{\otimes k})_{i_1\ldots i_k,j_1\ldots j_k}(\xi_\pi)_{j_1\ldots j_k}=(\xi_\pi)_{i_1\ldots i_k},\forall i\\
&\iff&\sum_j\delta_\pi(j_1,\ldots,j_k)v_{i_1j_1}\ldots v_{i_kj_k}=\delta_\pi(i_1,\ldots,i_k),\forall i
\end{eqnarray*}

Now let us recall that our corepresentation has the special form $v=diag(u,1)$. We conclude from this that for any index $a\in\{1,\ldots,k\}$, we must have:
$$i_a=N\implies j_a=N$$

With this observation in hand, if we denote by $i',j'$ the multi-indices obtained from $i,j$ obtained by erasing all the above $i_a=j_a=N$ values, and by $k'\leq k$ the common length of these new multi-indices, our condition becomes:
$$\sum_{j'}\delta_\pi(j_1,\ldots,j_k)(v^{\otimes k'})_{i'j'}=\delta_\pi(i_1,\ldots,i_k),\forall i$$

Here the index $j$ is by definition obtained from $j'$ by filling with $N$ values. In order to finish now, we have two cases, depending on $i$, as follows:

\medskip

\underline{Case 1}. Assume that the index set $\{a|i_a=N\}$ corresponds to a certain subpartition $\pi'\subset\pi$. In this case, the $N$ values will not matter, and our formula becomes:
$$\sum_{j'}\delta_\pi(j'_1,\ldots,j'_{k'})(v^{\otimes k'})_{i'j'}=\delta_\pi(i'_1,\ldots,i'_{k'})$$

\underline{Case 2}. Assume now the opposite, namely that the set $\{a|i_a=N\}$ does not correspond to a subpartition $\pi'\subset\pi$. In this case the indices mix, and our formula reads:
$$0=0$$

Thus, we are led to $\xi_{\pi'}\in Fix(v^{\otimes k'})$, for any subpartition $\pi'\subset\pi$, as claimed. Thus our claim is proved, and with this in hand, the result follows from Tannakian duality.
\end{proof}

Based on the above result, let us formulate the following definition:

\index{uniform quantum group}
\index{removing blocks}

\begin{definition}
An easy quantum group $G=(G_N)$, coming from a category of partitions $D\subset P$, is called uniform when we have, for any $N\in\mathbb N$:
$$G_{N-1}=G_N\cap U_{N-1}^+$$
Equivalently, $D$ must be stable under the operation which consists in removing blocks.
\end{definition}

We will see later on in this chapter some further motivations for this notion, probabilistic this time, the idea being that, in order for the computations for asymptotic truncated characters to work well, we must assume uniformity. For instance we will show that, in order for a liberation of easy quantum groups $G_N\to G_N^+$ to be compatible with the Bercovici-Pata bijection, we must assume uniformity. But more on this later.

\bigskip

Let us also mention that the notion of uniformity plays as well a key role in noncommutative geometry, the idea there being that, in order for the basic homogeneous spaces over an easy group $G_N$ to behave well, we must assume uniformity. We refer to \cite{ba3} for the story here, and we will be back to this, with a few details, in chapter 16 below.

\bigskip

At the level of the basic examples, the situation is as follows:

\begin{theorem}
The following happen:
\begin{enumerate}
\item The easy groups $O_N,U_N,B_N,C_N$ and $H_N^s$ are all uniform.

\item Their free versions $O_N^+,U_N^+,B_N^+,C_N^+$ and $H_N^{s+}$ are uniform too.

\item However, the various half-liberations are not uniform.
\end{enumerate}
\end{theorem}

\begin{proof}
This follows by using either of the criteria from Definition 6.2, as follows:

\medskip

(1) The fact that $O_N,U_N,B_N,C_N$ are uniform follows either from $G_{N-1}=G_N\cap U_{N-1}$, which is clearly satisfied in all cases, or from the fact that the corresponding categories of partitions, namely $P_2,\mathcal P_2,P_{12},\mathcal P_{12}$, are stable under removing blocks. As for the reflection groups $H_N^s$, these once again satisfy $G_{N-1}=G_N\cap U_{N-1}$, and the corresponding categories of partitions $P^s$ are clearly stable under removing blocks, too.

\medskip

(2) In the free case, we can prove that $O_N^+,U_N^+,B_N^+,C_N^+$ and $H_N^{s+}$ are uniform via $G_{N-1}=G_N\cap U_{N-1}^+$, with this requiring however some playing with the generators and relations defining our quantum groups, and we will leave this as an instructive exercise, or we can simply argue that the corresponding categories of partitions being those from (1) intersected with $NC$, these categories are stable under removing blocks.

\medskip

(3) Finally, in what regards half-liberations, here the result can be seen either with categories of partitions, or with intersections, the point in the half-classical case being that the relations $abc=cba$, when applied to the coefficients of a matrix of type $v=diag(u,1)$, collapse with $c=1$ to the usual commutation relations $ab=ba$. 
\end{proof}

All this is quite nice, and before going forward let us mention that, contrary to what Theorem 6.3 might suggest, the uniformity axiom is in fact something quite strong, which kills of sorts of pathologies that might appear. But, the point is that we are still at the beginning of this book, and we haven't talked about pathologies yet.

\bigskip

And to end with a funny story, when I wrote \cite{bsp} with Speicher we were mainly interested at that time in probability, where uniformity brings results, and we briefly thought of including uniformity in our easiness axioms. But, since we were quite excited at that time by our discovery of half-liberation, we did not do it, and wrote \cite{bsp} without uniformity, matter of talking about half-liberation too. I once told this story to a young researcher, having spent considerable time in his life in upgrading classification results from uniform to non-uniform, as to be fine with \cite{bsp}, and he was not happy at all.

\section*{6b. Classification results} 

For classification purposes the uniformity axiom is something very natural and useful, substantially cutting from complexity, and we have the following result, from \cite{bsp}:

\begin{theorem}
The classical and free uniform orthogonal easy quantum groups, with inclusions between them, are as follows:
$$
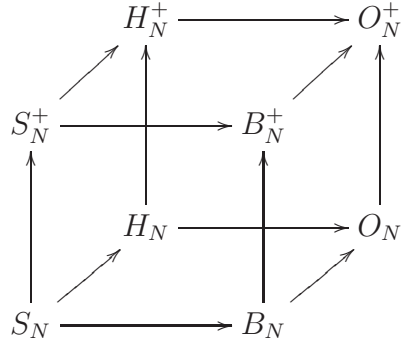
\xymatrix@R=20pt@C=20pt{
&H_N^+\ar[rr]&&O_N^+\\
S_N^+\ar[rr]\ar[ur]&&B_N^+\ar[ur]\\
&H_N\ar[rr]\ar[uu]&&O_N\ar[uu]\\
S_N\ar[uu]\ar[ur]\ar[rr]&&B_N\ar[uu]\ar[ur]
}$$
Moreover, this is an intersection/easy generation diagram, in the sense that for any of its square subdiagrams $P\subset Q,R\subset S$ we have $P=Q\cap R$ and $\{Q,R\}=S$.
\end{theorem}

\begin{proof}
There are several things to be proved, the idea being as follows:

\medskip

(1) We know that the quantum groups in the statement are indeed easy and uniform, the corresponding categories of partitions being as follows:
$$
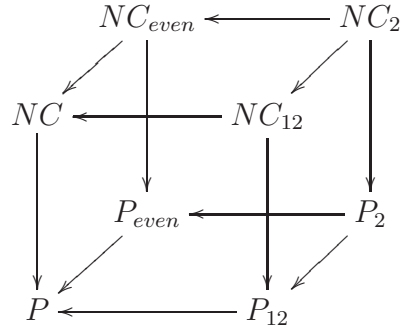
\xymatrix@R=20pt@C6pt{
&NC_{even}\ar[dl]\ar[dd]&&NC_2\ar[dl]\ar[ll]\ar[dd]\\
NC\ar[dd]&&NC_{12}\ar[dd]\ar[ll]\\
&P_{even}\ar[dl]&&P_2\ar[dl]\ar[ll]\\
P&&P_{12}\ar[ll]
}$$

Since this latter diagram is an intersection and generation diagram, we conclude that we have an intersection and easy generation diagram of quantum groups, as stated.

\medskip

(2) Regarding now the classification, consider first an easy group $S_N\subset G_N\subset O_N$. This must come from a certain category $P_2\subset D\subset P$, and if we assume $G=(G_N)$ to be uniform, then $D$ is uniquely determined by the subset $L\subset\mathbb N$ consisting of the sizes of the blocks of the partitions in $D$. Our claim is that the admissible sets are as follows:

\medskip

-- $L=\{2\}$, producing $O_N$.

\medskip

-- $L=\{1,2\}$, producing $B_N$.

\medskip

-- $L=\{2,4,6,\ldots\}$, producing $H_N$.

\medskip

-- $L=\{1,2,3,\ldots\}$, producing $S_N$.

\medskip

(3) Indeed, in one sense, this follows from our easiness results for $O_N,B_N,H_N,S_N$. In the other sense now, assume that $L\subset\mathbb N$ is such that the set $P_L$ consisting of partitions whose sizes of the blocks belong to $L$ is a category of partitions. We know from the axioms of the categories of partitions that the semicircle $\cap$ must be in the category, so we have $2\in L$. We claim that the following conditions must be satisfied as well:
$$k,l\in L,\,k>l\implies k-l\in L$$
$$k\in L,\,k\geq 2\implies 2k-2\in L$$

(4) Indeed, we will prove that both conditions follow from the axioms of the categories of
partitions. Let us denote by $b_k\in P(0,k)$ the one-block partition:
$$b_k=\left\{\begin{matrix}\sqcap\hskip-0.7mm \sqcap&\ldots&\sqcap\\
1\hskip2mm 2&\ldots&k\end{matrix} \right\}$$

For $k>l$, we can write $b_{k-l}$ in the following way:
$$b_{k-l}=\left\{\begin{matrix}\sqcap\hskip-0.7mm
\sqcap&\ldots&\ldots&\ldots&\ldots&\sqcap\\ 1\hskip2mm 2&\ldots&l&l+1&\ldots&k\\
\sqcup\hskip-0.7mm \sqcup&\ldots&\sqcup&|&\ldots&|\\ &&&1&\ldots&k-l\end{matrix}\right\}$$

In other words, we have the following formula:
$$b_{k-l}=(b_l^*\otimes |^{\otimes k-l})b_k$$

Since all the terms of this composition are in $P_L$, we have $b_{k-l}\in P_L$, and this proves our first claim. As for the second claim, this can be proved in a similar way, by capping two adjacent $k$-blocks with a $2$-block, in the middle.

\medskip

(5) With these conditions in hand, we can conclude in the following way:

\medskip

\underline{Case 1}. Assume $1\in L$. By using the first condition with $l=1$ we get:
$$k\in L\implies k-1\in L$$

This condition shows that we must have $L=\{1,2,\ldots,m\}$, for a certain number $m\in\{1,2,\ldots,\infty\}$. On the other hand, by using the second condition we get:
\begin{eqnarray*}
m\in L
&\implies&2m-2\in L\\
&\implies&2m-2\leq m\\
&\implies&m\in\{1,2,\infty\}
\end{eqnarray*}

The case $m=1$ being excluded by the condition $2\in L$, we reach to one of the two sets producing the groups $S_N,B_N$.

\medskip

\underline{Case 2}. Assume $1\notin L$. By using the first condition with $l=2$ we get:
$$k\in L\implies k-2\in L$$

This condition shows that we must have $L=\{2,4,\ldots,2p\}$, for a certain number $p\in\{1,2,\ldots,\infty\}$. On the other hand, by using the second condition we get:
\begin{eqnarray*}
2p\in L
&\implies&4p-2\in L\\
&\implies&4p-2\leq 2p\\
&\implies&p\in\{1,\infty\}
\end{eqnarray*}

Thus $L$ must be one of the two sets producing $O_N,H_N$, and we are done. 

\medskip

(6) In the free case, $S_N^+\subset G_N\subset O_N^+$, the situation is quite similar, the admissible sets being once again the above ones, producing this time $O_N^+,B_N^+,H_N^+,S_N^+$. See \cite{bsp}. 
\end{proof}

When removing the uniformity axiom things become more complicated, and the classification result here, from \cite{bsp}, \cite{rw3}, \cite{web}, is as follows:

\begin{theorem}
The classical and free orthogonal easy quantum groups are
$$
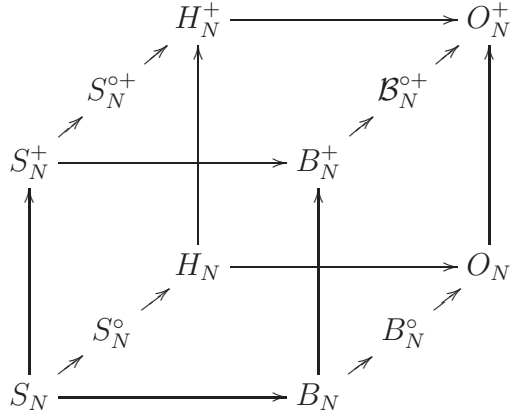
\xymatrix@R=7pt@C=7pt{
&&H_N^+\ar[rrrr]&&&&O_N^+\\
&S_N^{\circ+}\ar[ur]&&&&\mathcal B_N^{\circ+}\ar[ur]\\
S_N^+\ar[rrrr]\ar[ur]&&&&B_N^+\ar[ur]\\
\\
&&H_N\ar[rrrr]\ar[uuuu]&&&&O_N\ar[uuuu]\\
&S_N^\circ\ar[ur]&&&&B_N^\circ\ar[ur]\\
S_N\ar[uuuu]\ar[ur]\ar[rrrr]&&&&B_N\ar[uuuu]\ar[ur]
\\
}$$
with $S_N^\circ=S_N\times\mathbb Z_2$, $B_N^\circ=B_N\times\mathbb Z_2$, and with $S_N^{\circ+},\mathcal B_N^{\circ+}$ being their liberations, where $\mathcal B_N^{\circ+}$ stands for the two possible such liberations, $B_N^{\circ+}\subset B_N^{\circ\circ+}$.
\end{theorem}

\begin{proof}
The idea here is that of jointly classifying the ``classical'' categories of partitions $P_2\subset D\subset P$, and the ``free'' ones $NC_2\subset D\subset NC$. The situation is as follows:

\medskip

(1) At the classical level this leads, via a study which is quite similar to that from the uniform case, to 2 more groups, namely the groups $S_N^\circ,B_N^\circ$. See \cite{bsp}. 

\medskip

(2) At the free level we obtain 3 more quantum groups, $S_N^{\circ+},B_N^{\circ+},B_N^{\circ+}$, with the inclusion $B_N^{\circ+}\subset B_N^{\circ+}$, which is something a bit surprising, being best thought of as coming from an inclusion $B_N^\circ\subset B_N^{\circ\circ}$, which is in fact an isomorphism. See \cite{web}.

\medskip

(3) In short, all this is routine, except for some subtleties in the continuous case, that we will explain now. The quantum groups concerned are as follows, with all being objects that we know, except for $B_N^{\circ+},B_N^{\circ\circ+}$, whose definition will come in a moment:
$$\xymatrix@R=55pt@C=46pt{
B_N^+\ar[r]&B_N^{\circ+}\ar[r]&B_N^{\circ\circ+}\ar[r]&O_N^+\\
B_N\ar[u]\ar[r]&B_N^\circ\ar[u]\ar[rr]&&O_N\ar[u]}$$

As for the corresponding categories of partitions, these are as follows, with again objects that we know, except for $NC_{12}^\circ,NC_{12}^{\circ\circ}$, whose definition will come in a moment:
$$\xymatrix@R=50pt@C=44pt{
NC_{12}\ar[d]&NC_{12}^\circ\ar[l]\ar[d]&NC_{12}^{\circ\circ}\ar[l]&NC_2\ar[l]\ar[d]\\
P_{12}&P_{12}^\circ\ar[l]&&P_2\ar[ll]}$$

(4) Getting now to the core of the problem, we know that $B_N^\circ=B_N\times\mathbb Z_2$ appears from the category $P_{12}^\circ$ of singletons and pairings, having an even total length, and with this coming from the basics of the $\times\mathbb Z_2$ operation. But this suggests to define $NC_{12}^\circ$ as follows, and say that $B_N^{\circ+}$ is the quantum group associated to this category:
$$NC_{12}^\circ=P_{12}^\circ\cap NC$$

(5) However, this is a beginner mistake, done in \cite{bsp}, and no wonder here, because that paper was the beginning of general easiness. Later Weber came in \cite{web} with the correct solution. We have the following formulae for $P_{12}$, which are both clear:
$$P_{12}
=<\slash\hskip-2.1mm\backslash\,,|\cap\hskip-4.2mm{\ }_|>
=<\slash\hskip-2.1mm\backslash\,,||>$$

The point now is that when liberating at the level of these formulae, that is, when removing the crossing, we obtain two distinct categories, as follows:
$$NC_{12}^\circ=<|\cap\hskip-4.2mm{\ }_|>
\quad\supset\quad 
NC_{12}^{\circ\circ}=<||>$$

To be more precise, the category on the left $NC_{12}^\circ$ is the one that we know from (4), noncrossing singletons and pairings, having an even total length. As for the category on the right $NC_{12}^{\circ\circ}$, this is defined as above, and is certainly a subcategory of $NC_{12}^\circ$, because it consists of certain noncrossing singletons and pairings, having an even total length, and there is even a direct, rock-solid proof of this inclusion, as follows:
$$||=[^\pi_\sigma]\in<|\cap\hskip-4.2mm{\ }_|>\quad:\quad\pi=|\cap\hskip-4.2mm{\ }_|\quad,\quad\sigma=||\cup$$

However, we do not have equality, due to the following somewhat bizarre fact:
$$|\cap\hskip-4.2mm{\ }_|\ \notin\ <||>$$

(6) So, this was for the story, we have two noncrossing versions of $P_{12}^\circ$, and so two liberations of $B_N^\circ$, as constructed in (3) above. And for details, regarding all this, and then the fix of the previous classification from \cite{bsp}, we refer to \cite{web}.
\end{proof}

Let us discuss now the unitary case. Here things are considerably more complicated, and even when imposing the uniformity condition, there are too many examples. So, in order to get started, the best is to use uniformity, along with a second axiom. 

\bigskip

In order to formulate our second axiom, which is something very natural too, and has its own interest, consider the cube $T_N=\mathbb Z_2^N$, regarded as diagonal torus of $O_N$. We have then the following result, which is something nice, providing us with our second axiom:

\index{twistable quantum group}

\begin{proposition}
For an easy quantum group $G=(G_N)$, coming from a category of partitions $D\subset P$, the following conditions are equivalent:
\begin{enumerate}
\item $T_N\subset G_N$.

\item $H_N\subset G_N$.

\item $D\subset P_{even}$.
\end{enumerate}
If these conditions are satisfied, we say that $G_N$ is twistable.
\end{proposition}

\begin{proof}
We use the general easiness theory developed in chapters 1-4, and more specifically, the easy envelope operation $G\to\bar{G}$ introduced in chapter 3:

\medskip

$(1)\iff(2)$ Here it is enough to check that the easy envelope $\bar{T}_N$ of the cube equals the hyperoctahedral group $H_N$. But this follows from:
$$\bar{T}_N
=\overline{<T_N,S_N>}
=\bar{H}_N
=H_N$$

$(2)\iff(3)$ This follows by functoriality, from the fact that $H_N$ comes from the category of partitions $P_{even}$, that we know from chapter 1.
\end{proof}

The teminology in the above result comes from the fact that, assuming $D\subset P_{even}$, we can indeed twist $G_N$, into a certain ``quizzy'' quantum group $G_N'$, and vice versa:
$$[H_N\subset G_N\subset U_N]\ \longleftrightarrow\ [H_N\subset G_N'\subset U_N']$$

We refer to chapter 13 for details regarding the operation $G_N\to G_N'$. In what follows we will not need this twisting procedure, and we will just use Proposition 6.6 as it is, as a statement providing us with a simple and natural condition to be imposed on $G_N$. 

\bigskip

In practice now, imposing this second axiom leads to something nice, namely:

\begin{proposition}
Among the easy quantum groups that we know so far, from Theorem 6.3 and its proof, those which are uniform and twistable are as follows,
$$
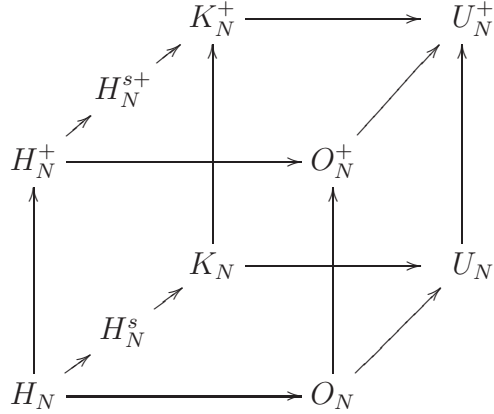
\xymatrix@R=7pt@C=7pt{
&&K_N^+\ar[rrrr]&&&&\ \ U_N^+\\
&H_N^{s+}\ar[ur]&&&&\\
H_N^+\ar[rrrr]\ar[ur]&&&&O_N^+\ar[uurr]\\
\\
&&K_N\ar[rrrr]\ar[uuuu]&&&&\ \ U_N\ar[uuuu]\\
&H_N^s\ar[ur]&&&&\\
H_N\ar[uuuu]\ar[ur]\ar[rrrr]&&&&O_N\ar[uuuu]\ar[uurr]
\\
}$$
where $H_N^s=\mathbb Z_s\wr S_N$ and $H_N^{s+}=\mathbb Z_s\wr_*S_N^+$, as usual, with $s\in\{2,4,\ldots,\infty\}$.
\end{proposition}

\begin{proof}
There are two assertions here, the idea being as follows:

\medskip

(1) We know from Theorem 6.3 that all the quantum groups in the statement are uniform, and since all these quantum groups contain $H_N$, provided that we assume $s\in\{2,4,\ldots,\infty\}$ at the end, as indicated, these quantum groups are twistable too. 

\medskip

(2) In what concerns the uniqueness assertion, this is of course something informal, and with this coming again from Theorem 6.3, which excludes the half-liberations, and from the fact that $B_N,C_N,B_N^+,C_N^+$ and the missing $H_N^s,H_N^{s+}$ are clearly not twistable.
\end{proof}

All this classification business becomes a bit complicated, so time for a pause, and some thinking. Looking at what we have in Proposition 6.7, in comparison with our previous results, Theorem 6.4 and Theorem 6.5, suggests restricting the attention to the upper and lower faces of the cube, with the aim of proving that we have there is everything. 

\bigskip

However, a bit suprisingly, this is not all. Recall that the free complexification $(\tilde{G},\tilde{u})$ of a quantum group $(G,u)$ is obtained by considering the subalgebra $C(\tilde{G})\subset C(\mathbb T)*C(G)$ generated by the entries of $\tilde{u}=zu$, where $z$ is the standard generator of $C(\mathbb T)$. With this convention, we have the following intriguing extra example, from \cite{tw2}:

\begin{proposition}
The free complexification of the full quantum reflection group 
$$K_N^{++}=\widetilde{K_N^+}$$ 
is easy, and appears as an intermediate object, as follows,
$$K_N^+\subset K_N^{++}\subset U_N^+$$
with both inclusions being proper.
\end{proposition}

\begin{proof}
By composing the canonical inclusion $C(K_N^{++})\subset C(\mathbb T)*C(K_N^+)$ with the map $\varepsilon*id$, with $\varepsilon$ being the counit, we obtain a morphism $C(K_N^{++})\to C(K_N^+)$ mapping $\tilde{u}_{ij}\to u_{ij}$, so we have inclusions as in the statement. Now since the following elements are projections, and are not equal, both these inclusions are proper:
$$p_{ij}=\tilde{u}_{ij}\tilde{u}_{ij}^*=zu_{ij}u_{ij}^*z^*$$
$$q_{ij}=\tilde{u}_{ij}^*\tilde{u}_{ij}=u_{ij}^*u_{ij}$$

Regarding the easiness claim, this follows from the general theory of the representations of free complexifications \cite{rau}. To be more precise, as explained in \cite{tw2}, the associated category $\mathcal{NC}_{even}^-$ is that of the even noncrossing partitions which, when rotated on one line, have alternating colors in each block. Observe that the inclusions in the statement correspond then to  the inclusions at the partition level, which are as follows:
$$\mathcal{NC}_{even}\supset\mathcal{NC}_{even}^-\supset\mathcal{NC}_2$$

Thus, we are led to the conclusions in the statement.
\end{proof}

Obviously, $K_N^{++}$ is something quite annoying. However, in connection with our classification questions, the news are good, because we can now turn Proposition 6.7 into a theorem. The result is as follows, where by ``classical/twisted'' and ``free'' we mean $\backslash\hskip-2.1mm/\in D$ and $D\subset NC_{even}$, where $D\subset P_{even}$ is the associated category of partitions:

\begin{theorem}
The classical and free uniform twistable quantum groups are
$$
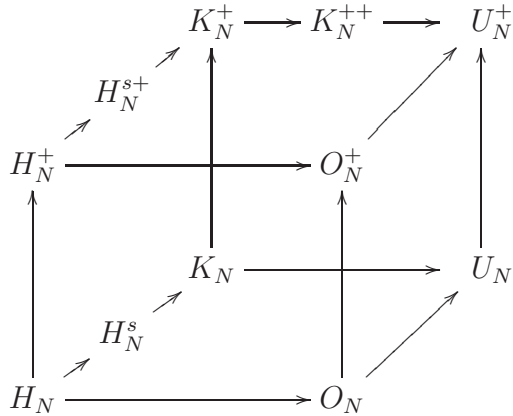
\xymatrix@R=7pt@C=7pt{
&&K_N^+\ar[rr]&&K_N^{++}\ar[rr]&&\ \ U_N^+\\
&H_N^{s+}\ar[ur]&&&&\\
H_N^+\ar[rrrr]\ar[ur]&&&&O_N^+\ar[uurr]\\
\\
&&K_N\ar[rrrr]\ar[uuuu]&&&&\ \ U_N\ar[uuuu]\\
&H_N^s\ar[ur]&&&&\\
H_N\ar[uuuu]\ar[ur]\ar[rrrr]&&&&O_N\ar[uuuu]\ar[uurr]
\\
}$$
where $H_N^s=\mathbb Z_s\wr S_N$, $H_N^{s+}=\mathbb Z_s\wr_*S_N^+$, with $s\in\{2,4,\ldots,\infty\}$, and $K_N^{++}=\widetilde{K_N^+}$. 
\end{theorem}

\begin{proof}
This is a consequence of the classification results in \cite{tw1}, \cite{tw2}. Consider indeed a uniform category of partitions, as follows:
$$\mathcal{NC}_2\subset D\subset P_{even}$$ 

We must prove that in the classical case, where $\backslash\hskip-2.1mm/\in D$, the only solutions are the following categories, corresponding to the lower face of the above cube:
$$D=P_2,\mathcal P_2,P_{even}^s$$

We must prove as well that in the free case, where $D\subset NC_{even}$, the only solutions are the following categories, corresponding to the upper face of the above cube:
$$D=NC_2,\mathcal{NC}_2,\mathcal{NC}_{even}^-,NC_{even}^s$$

We jointly investigate these two problems. Let $B$ be the set of all possible labelled blocks in $D$, having no upper legs. Observe that $B$ is stable under the switching of colors operation, $\circ\leftrightarrow\bullet$. We have two possible situations, as follows:

\medskip

(1) $B$ consists of pairings only. Here the pairings in question can be either all labelled pairings, namely $\circ-\circ$, $\circ-\bullet$, $\bullet-\circ$, $\bullet-\bullet$, or just the matching ones, namely $\circ-\bullet$, $\bullet-\circ$, and we obtain here $P_2,\mathcal P_2$ in the classical case, and $NC_2,\mathcal{NC}_2$ in the free case.

\medskip

(2) $B$ has at least one block of size $\geq 4$. In this case we can let $s\in\{2,4,\ldots,\infty\}$ to be the length of the smallest $\circ\ldots\circ$ block, and we obtain in this way the category $P_{even}^s$ in the classical case, and the categories $\mathcal{NC}_{even}^-,NC_{even}^s$ in the free case. See \cite{tw1}.
\end{proof}

When removing the twistability and the uniformity assumptions, things become more complicated, because we have to deal with both the phenomena appearing from Theorem 6.6 and Theorem 6.9. However, a full classification result in the classical and free cases is available, and for full details on all this, we refer to \cite{tw1}, \cite{tw2}.

\section*{6c. Ground Zero} 

We further discuss now the general uniform and twistable case, $H_N\subset G_N\subset U_N^+$. In this case, we can imagine $G_N$ as sitting inside the standard cube:
$$
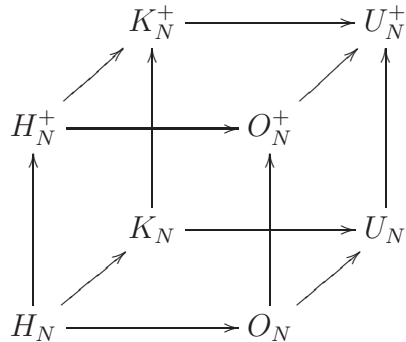
\xymatrix@R=20pt@C=20pt{
&K_N^+\ar[rr]&&U_N^+\\
H_N^+\ar[rr]\ar[ur]&&O_N^+\ar[ur]\\
&K_N\ar[rr]\ar[uu]&&U_N\ar[uu]\\
H_N\ar[uu]\ar[ur]\ar[rr]&&O_N\ar[uu]\ar[ur]
}$$

The point now is that, by using the operations $\cap$ and $\{\,,\}$, we can in principle ``project'' $G_N$ on the faces and edges of the cube, and then use some kind of 3D orientation coming from this, in order to deduce structure and classification results. Let us start with:

\index{classical version}
\index{discrete version}
\index{real version}
\index{free version}
\index{smooth version}
\index{unitary version}

\begin{definition}
Associated to any twistable easy quantum group $H_N\subset G_N\subset U_N^+$ are its classical, discrete and real versions, given by the following formulae,
$$G_N^c=G_N\cap U_N\quad,\quad 
G_N^d=G_N\cap K_N^+\quad,\quad 
G_N^r=G_N\cap O_N^+$$
as well as its free, smooth and unitary versions, given by the following formulae,
$$G_N^f=\{G_N,H_N^+\}\quad,\quad
G_N^s=\{G_N,O_N\}\quad,\quad
G_N^u=\{G_N,K_N\}$$
where $\cap$ and $\{\,,\}$ are respectively the intersection and easy generation operations.
\end{definition}

Here the classical, real and unitary versions are something quite standard. Regarding the discrete and smooth versions, in the classical case, $G_N\subset U_N$, our constructions produce indeed discrete and smooth versions, and this is where our terminology comes from. Of course, it would be nice to have more results on these operations.

\bigskip

Finally, regarding the free version, the various results that we have show that the liberation operation $G_N\to G_N^+$ usually appears via the following formula:
$$G_N^+=\{G_N,S_N^+\}$$

But in the twistable setting, where we have $H_N\subset G_N$, this is the same as setting:
$$G_N^+=\{G_N,H_N^+\}$$

All this is of course a bit theoretical, and this is why we use the symbol $f$ for free versions in the above sense, and keep $+$ for well-known, established liberations. 

\bigskip

In relation now with our questions, and our 3D plan, we can now formulate:

\begin{proposition}
Given an intermediate quantum group $H_N\subset G_N\subset U_N^+$, we have a diagram of closed subgroups of $U_N^+$, obtained by inserting
$$
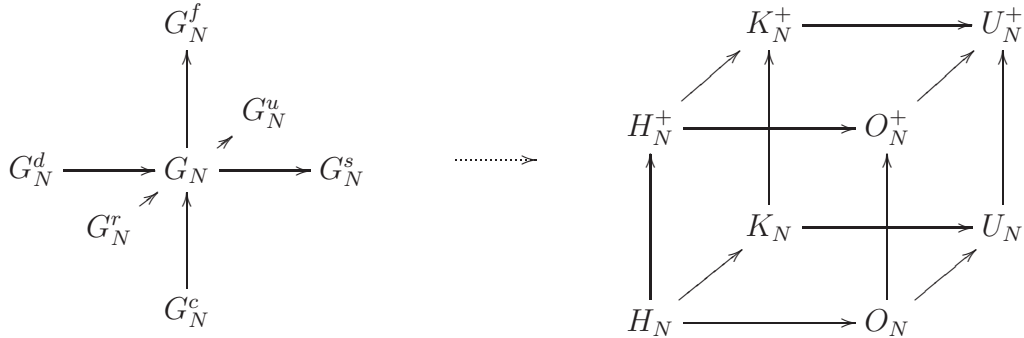
\xymatrix@R=3pt@C=5pt{
&&G_N^f&&\\
\\
&&&G_N^u&\\
G_N^d\ar[rr]&&G_N\ar[rr]\ar[uuu]\ar[ur]&&G_N^s\\
&G_N^r\ar[ur]&&&\\
\\
&&G_N^c\ar[uuu]&&}
\qquad\xymatrix@R=10pt@C=30pt{\\ \\ \\ \ar@.[r]&}\qquad
\xymatrix@R=19pt@C=20pt{
&K_N^+\ar[rr]&&U_N^+\\
H_N^+\ar[rr]\ar[ur]&&O_N^+\ar[ur]\\
&K_N\ar[rr]\ar[uu]&&U_N\ar[uu]\\
H_N\ar[uu]\ar[ur]\ar[rr]&&O_N\ar[uu]\ar[ur]
}$$
in the obvious way, with each $G_N^x$ belonging to the main diagonal of each face.
\end{proposition}

\begin{proof}
The fact that we have indeed the diagram of inclusions on the left is clear from the constructions of the quantum groups involved. Regarding the insertion procedure, consider any of the faces of the cube, denoted as follows:
$$P\subset Q,R\subset S$$

Our claim is that the corresponding quantum group $G=G_N^x$ can be inserted on the corresponding main diagonal $P\subset S$, as follows:
$$
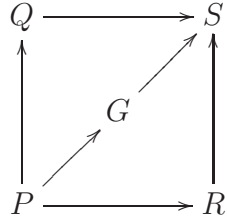
\xymatrix@R=20pt@C=20pt{
Q\ar[rr]&&S\\
&G\ar[ur]\\
P\ar[rr]\ar[uu]\ar[ur]&&R\ar[uu]}$$

We have to check here a total of $6\times 2=12$ inclusions. But, according to Definition 6.10, these inclusions that must be checked are as follows:

\medskip

(1) $H_N\subset G_N^c\subset U_N$, where $G_N^c=G_N\cap U_N$.

\medskip

(2) $H_N\subset G_N^d\subset K_N^+$, where $G_N^d=G_N\cap K_N^+$.

\medskip

(3) $H_N\subset G_N^r\subset O_N^+$, where $G_N^r=G_N\cap O_N^+$.

\medskip

(4) $H_N^+\subset G_N^f\subset U_N^+$, where $G_N^f=\{G_N,H_N^+\}$.

\medskip

(5) $O_N\subset G_N^s\subset U_N^+$, where $G_N^s=\{G_N,O_N\}$.

\medskip

(6) $K_N\subset G_N^u\subset U_N^+$, where $G_N^u=\{G_N,K_N\}$.

\medskip

All these statements being trivial from the definition of the intersection operation $\cap$ and of the easy generation operation $\{\,,\}$, and from our assumption $H_N\subset G_N\subset U_N^+$, our insertion procedure works indeed, and we are done.
\end{proof}

In order now to complete the diagram, we have to project as well our quantum group $G_N$ on the edges of the cube. For this purpose we can basically assume, by replacing $G_N$ with each of its 6 projections on the faces, that $G_N$ actually lies on one of the six faces.

\bigskip

The technical result that we will need here is as follows:

\begin{proposition}
Given an intersection and easy generation diagram $P\subset Q,R\subset S$ and an intermediate easy quantum group $P\subset G\subset S$, as follows,
$$
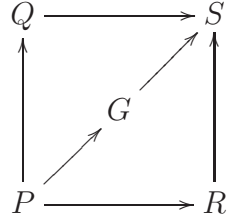
\xymatrix@R=20pt@C=20pt{
Q\ar[rr]&&S\\
&G\ar[ur]\\
P\ar[rr]\ar[uu]\ar[ur]&&R\ar[uu]}$$
we can extend this diagram into a diagram as follows:
$$
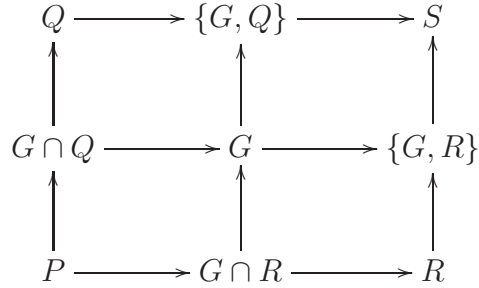
\xymatrix@R=30pt@C=30pt{
Q\ar[r]&\{G,Q\}\ar[r]&S\\
G\cap Q\ar[u]\ar[r]&G\ar[r]\ar[u]&\{G,R\}\ar[u]\\
P\ar[r]\ar[u]&G\cap R\ar[u]\ar[r]&R\ar[u]}$$
In addition, $G$ ``slices the square'', in the sense that this is an intersection and easy generation diagram, precisely when $G=\{G\cap Q,G\cap R\}$ and $G=\{G,Q\}\cap\{G,R\}$.
\end{proposition}

\begin{proof}
This is clear from definitions, because the intersection and easy generation conditions are automatic for the upper left and lower right squares, and so are half of the intersection and easy generation conditions for the lower left and upper right squares.
\end{proof}

Now back to 3 dimensions, and to the cube, we have the following result:

\begin{proposition}
Assuming that $H_N\subset G_N\subset U_N^+$ satisfies the conditions
$$G_N^{cs}=G_N^{sc}\quad,\quad
G_N^{cu}=G_N^{uc}\quad,\quad
G_N^{df}=G_N^{fd}$$
$$G_N^{du}=G_N^{ud}\quad,\quad 
G_N^{rf}=G_N^{fr}\quad,\quad
G_N^{rs}=G_N^{sr}$$
the diagram in Proposition 6.11 can be completed, via the construction in Proposition 6.12, into a diagram dividing the cube along the $3$ coordinates axes, into $8$ small cubes.
\end{proposition}

\begin{proof}
We have to prove that the 12 projections on the edges are well-defined, with the problem coming from the fact that each of these projections can be defined in 2 possible ways, depending on the face that we choose first. The verification goes as follows:

\medskip

(1) Regarding the $3$ edges emanating from $H_N$, the result here follows from:
$$G_N^{cd}=G_N^{dc}=G_N\cap K_N$$ 
$$G_N^{cr}=G_N^{rc}=G_N\cap O_N$$
$$G_N^{dr}=G_N^{rd}=G_N\cap H_N^+$$

These formulae are indeed all trivial, of type:
$$(G\cap Q)\cap R=(G\cap R)\cap Q=G\cap P$$

(2) Regarding the $3$ edges landing into $U_N^+$, the result here follows from:
$$G_N^{fs}=G_N^{sf}=\{G_N,O_N^+\}$$
$$G_N^{fu}=G_N^{uf}=\{G_N,K_N^+\}$$
$$G_N^{su}=G_N^{us}=\{G_N,U_N\}$$

These formulae are once again trivial, of type:
$$\{\{G,Q\},R\}=\{\{G,R\},Q\}=\{G,S\}$$

(3) Finally, regarding the remaining $6$ edges, not emanating from $H_N$ or landing into $U_N^+$, here the result follows from our assumptions in the statement.
\end{proof}

Unfortunately, we are not done yet, because nothing guarantees that we obtain an intersection and easy generation diagram. Thus, we must add more axioms, as follows:

\index{slicing the cube}

\begin{theorem}
Assume that $H_N\subset G_N\subset U_N^+$ satisfies the following conditions, where by ``intermediate'' we mean in each case ``parallel to its neighbors'':
\begin{enumerate}
\item The $6$ compatibility conditions in Proposition 6.13,

\item $G_N^c,G_N,G_N^f$ slice the classical/intermediate/free faces,

\item $G_N^d,G_N,G_N^s$ slice the discrete/intermediate/smooth faces,

\item $G_N^r,G_N,G_N^u$ slice the real/intermediate/unitary faces,
\end{enumerate}
Then $G_N$ ``slices the cube'', in the sense that the diagram obtained in Proposition 6.13 is an intersection and easy generation diagram.
\end{theorem}

\begin{proof}
This follows indeed from Proposition 6.12 and Proposition 6.13.
\end{proof}

Summarizing, we are done now with our geometric program, and we have a whole collection of natural geometric conditions that can be imposed to $G_N$. It is quite clear that $G_N$ can be reconstructed from its edge projections, so in order to do the classification, we first need a ``coordinate system''. Common sense would suggest to use the one emanating from $H_N$, or perhaps the one landing into $U_N^+$. However, technically speaking, the best is to use the coordinate system based at $O_N$, highlighted below:
$$
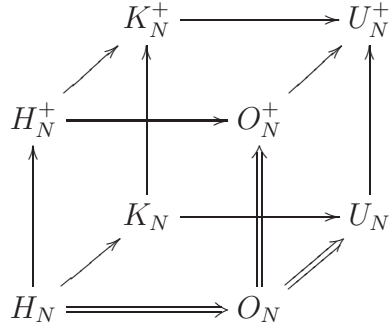
\xymatrix@R=18pt@C=18pt{
&K_N^+\ar[rr]&&U_N^+\\
H_N^+\ar[rr]\ar[ur]&&O_N^+\ar[ur]\\
&K_N\ar[rr]\ar[uu]&&U_N\ar[uu]\\
H_N\ar[uu]\ar[ur]\ar@=[rr]&&O_N\ar@=[uu]\ar@=[ur]
}$$

This choice comes from the fact that the classification result for $O_N\subset O_N^+$, explained below, is something very simple. And this is not the case with the results for $H_N\subset H_N^+$ and for $U_N\subset U_N^+$, from \cite{mw2}, \cite{rw3} which are quite complicated, with uncountably many solutions, in the general non-uniform case. As for the result for $K_N\subset K_N^+$, this is not available yet, but it is known that there are uncountably many solutions here as well.

\bigskip

So, here is now the key result, from \cite{bv2}, dealing with the vertical direction:

\index{orthogonal quantum group}
\index{half-classical orthogonal group}
\index{category of pairings}

\begin{theorem}
There is only one proper intermediate easy quantum group
$$O_N\subset G_N\subset O_N^+$$
namely the quantum group $O_N^*$, which is not uniform.
\end{theorem}

\begin{proof}
We must compute here the categories of pairings, as follows:
$$NC_2\subset D\subset P_2$$

But this can be done via some standard combinatorics, in three steps, as follows:

\medskip

(1) Let $\pi\in P_2-NC_2$, having $s\geq 4$ strings. Our claim is that:

\medskip

-- If $\pi\in P_2-P_2^*$, there exists a semicircle capping $\pi'\in P_2-P_2^*$.

\smallskip

-- If $\pi\in P_2^*-NC_2$, there exists a semicircle capping $\pi'\in P_2^*-NC_2$.

\medskip

Indeed, both these assertions can be easily proved, by drawing pictures.

\medskip

(2) Consider now a partition $\pi\in P_2(k,l)-NC_2(k,l)$. Our claim is that:

\medskip

-- If $\pi\in P_2(k, l)-P_2^*(k,l)$ then $<\pi>=P_2$.

\smallskip

-- If $\pi\in P_2^*(k,l)-NC_2(k,l)$ then $<\pi>=P_2^*$.

\medskip

This can be indeed proved by recurrence on the number of strings, $s=(k+l)/2$, by using (1), which provides us with a descent procedure $s\to s-1$, at any $s\geq4$.

\medskip

(3) Finally, assume that we are given an easy quantum group $O_N\subset G\subset O_N^+$, coming from certain sets of pairings $D(k,l)\subset P_2(k,l)$. We have three cases:

\medskip

-- If $D\not\subset P_2^*$, we obtain $G=O_N$.

\smallskip

-- If $D\subset P_2,D\not\subset NC_2$, we obtain $G=O_N^*$.

\smallskip

-- If $D\subset NC_2$, we obtain $G=O_N^+$.

\medskip

Thus, we have proved the uniquess result. As for the non-uniformity of the unique solution, $O_N^*$, this is something that we already know, from the above.
\end{proof}

Here is now another basic result that we will need, in order to perform our classification work here, dealing this time with the ``discrete vs. continuous'' direction:

\begin{theorem}
There are no proper intermediate easy groups
$$H_N\subset G_N\subset O_N$$
except for $H_N,O_N$ themselves.
\end{theorem}

\begin{proof}
We must prove that there are no proper intermediate categories as follows:
$$P_2\subset D\subset P_{even}$$

But this can done via some combinatorics. For details here, see \cite{bsp}.
\end{proof}

Finally, here is a third and last result that we will need, for our classification work here, regarding the missing direction, namely the  ``real vs. complex'' one:

\begin{theorem}
The proper intermediate easy groups
$$O_N\subset G_N\subset U_N$$
are the groups $\mathbb Z_rO_N$ with $r\in\{2,3,\ldots,\infty\}$, which are not uniform.
\end{theorem}

\begin{proof}
This is standard and well-known, from \cite{tw2}, the proof being as follows:

\medskip

(1) Our first claim, which is elementary, is that $\mathbb TO_N\subset U_N$ is easy, the corresponding category of partitions being the subcategory $\bar{P}_2\subset P_2$ consisting of the pairings having the property that when flatenning, we have the global formula $\#\circ=\#\bullet$. 

\medskip

(2) Our second claim, which is elementary too, is that, more generally, the group $\mathbb Z_rO_N\subset U_N$ is easy, with the corresponding category $P_2^r\subset P_2$ consisting of the pairings having the property that when flatenning, we have the global formula $\#\circ=\#\bullet(r)$. 

\medskip

(3) In what regards now the converse, stating that the above groups $O_N\subset\mathbb Z_rO_N\subset U_N$ are the only ones, we must compute the following categories of pairings:
$$\mathcal P_2\subset D\subset P_2$$

But this can be done, via some standard combinatorics, and we obtain the result. We refer here to Tarrago-Weber \cite{tw2}, and we will be back to this, later in this book.
\end{proof}

We can now formulate a classification result, as follows:

\begin{theorem}[Ground Zero]
There are exactly eight closed subgroups $G_N\subset U_N^+$ having the following properties,
\begin{enumerate}
\item Easiness,

\item Uniformity,

\item Twistability,

\item Slicing property,
\end{enumerate}
namely the quantum groups $O_N,U_N,H_N,K_N$ and $O_N^+,U_N^+,H_N^+,K_N^+$.
\end{theorem}

\begin{proof}
We already know that the 8 quantum groups in the statement have indeed the properties (1-4), and form a cube, as follows:
$$
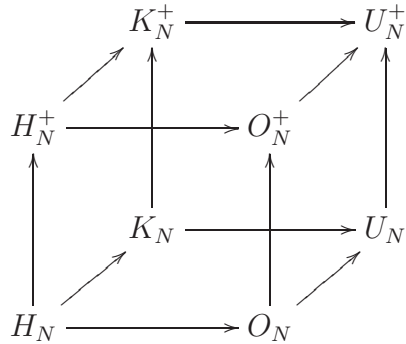
\xymatrix@R=20pt@C=20pt{
&K_N^+\ar[rr]&&U_N^+\\
H_N^+\ar[rr]\ar[ur]&&O_N^+\ar[ur]\\
&K_N\ar[rr]\ar[uu]&&U_N\ar[uu]\\
H_N\ar[uu]\ar[ur]\ar[rr]&&O_N\ar[uu]\ar[ur]
}$$

Conversely now, assuming that an easy quantum group $G=(G_N)$ has the above properties (2-4), the twistability property, (3), tells us that we have:
$$H_N\subset G_N\subset U_N^+$$

Thus $G_N$ sits inside the cube, and the above discussion applies. To be more precise, let us project $G$ on the faces of the cube, as in Proposition 6.11:
$$
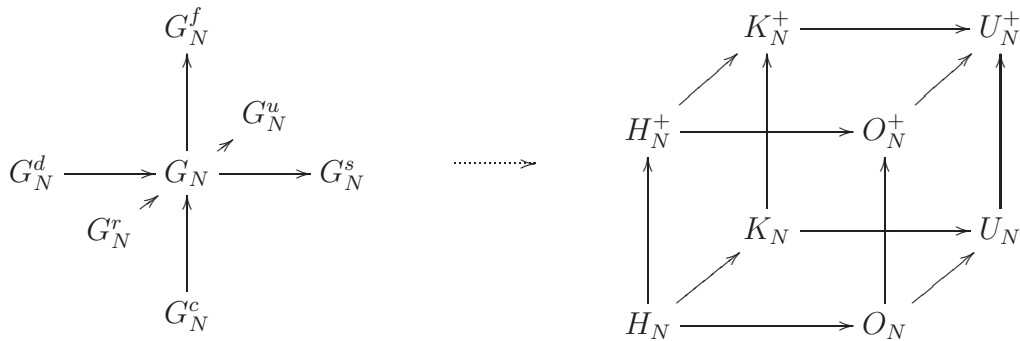
\xymatrix@R=3pt@C=5pt{
&&G_N^f&&\\
\\
&&&G_N^u&\\
G_N^d\ar[rr]&&G_N\ar[rr]\ar[uuu]\ar[ur]&&G_N^s\\
&G_N^r\ar[ur]&&&\\
\\
&&G_N^c\ar[uuu]&&}
\qquad\xymatrix@R=10pt@C=30pt{\\ \\ \\ \ar@.[r]&}\qquad
\xymatrix@R=19pt@C=20pt{
&K_N^+\ar[rr]&&U_N^+\\
H_N^+\ar[rr]\ar[ur]&&O_N^+\ar[ur]\\
&K_N\ar[rr]\ar[uu]&&U_N\ar[uu]\\
H_N\ar[uu]\ar[ur]\ar[rr]&&O_N\ar[uu]\ar[ur]
}$$

In order to compute these projections, and eventually prove that $G_N$ is one of the vertices of the cube, we can use use the coordinate system based at $O_N$:
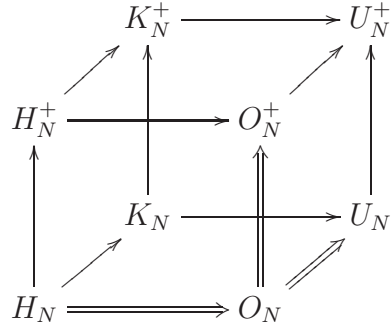
$$\xymatrix@R=18pt@C=18pt{
&K_N^+\ar[rr]&&U_N^+\\
H_N^+\ar[rr]\ar[ur]&&O_N^+\ar[ur]\\
&K_N\ar[rr]\ar[uu]&&U_N\ar[uu]\\
H_N\ar[uu]\ar[ur]\ar@=[rr]&&O_N\ar@=[uu]\ar@=[ur]
}$$

Now by using our classification results, Theorem 6.15, Theorem 6.16 and Theorem 6.17, along with the uniformity condition, (2), we conclude that the edge projections of $G_N$ must be among the vertices of the cube. Moreover, by using the slicing axiom, (4), we deduce from this that $G_N$ itself must be a vertex of the cube. Thus, we have exactly 8 solutions to our problem, namely the vertices of the cube, as claimed.
\end{proof}

Now that we have our Ground Zero, obtained by heavily bombarding the quantum group world, with all the natural axioms available, we can start building. And the point is that when dropping the easiness axiom, some classification results are possible:

\bigskip

(1) In the classical case, we believe that the uniform, half-homogeneous, oriented groups are the obvious ones, with some bistochastic versions excluded. This is of course something quite heavy, well beyond easiness, with the potential tools available for proving such things coming from advanced finite group theory and Lie algebra theory.  Our uniformity axiom could play a key role here, when combined with \cite{sto}, in order to exclude all the exceptional objects which might appear on the way.

\bigskip

(2) In the free case, under similar assumptions, we believe that the solutions should be again the obvious ones, once again with some bistochastic versions excluded. This is something heavy, too, related to the generation conjecture $<G_N,S_N^+>=\{\bar{G}_N,S_N^+\}$. Indeed, assuming that we would have such a formula, and perhaps some more formulae of the same type as well, we can in principle work out our way inside the cube, from the edge and face projections to $G_N$ itself, and in this process $G_N$ would become easy.

\bigskip

(3) In the group dual case, the orientability axiom simplifies, because the group duals are discrete in our sense. We believe that the uniform, twistable, oriented group duals should appear as combinations of certain abelian groups, which appear in the classical case, with duals of varieties of real reflection groups, which appear in the real case. This is probably the easiest question in the present series, and the most reasonable one, to start with. However, there are no concrete results so far, in this direction.

\section*{6d. Bercovici-Pata} 

Getting back to analysis questions, we will prove now that for a liberation of uniform easy groups, $G_N\to G_N^+$, the asymptotic laws of the truncated characters $\chi_t$ are in Bercovici-Pata bijection. We certainly know that this happens, based on our classification results, and explicit computations of laws of truncated characters, but we would like to have an abstract, conceptual proof of this, not based on classification results.

\bigskip

We will need a number of combinatorial preliminaries. Let us start with:

\index{Weingarten function}

\begin{proposition}
For the group $S_N$ the Weingarten function is given by
$$W_{kN}(\pi,\sigma)=\sum_{\tau\leq\pi\wedge\sigma}\mu(\tau,\pi)\mu(\tau,\sigma)\frac{(N-|\tau|)!}{N!}$$
and satisfies the following estimate,
$$W_{kN}(\pi,\sigma)=N^{-|\pi\wedge\sigma|}(
\mu(\pi\wedge\sigma,\pi)\mu(\pi\wedge\sigma,\sigma)+O(N^{-1}))$$
with $\mu$ being the M\"obius function of $P(k)$.
\end{proposition}

\begin{proof}
The first assertion follows from the Weingarten formula. Indeed, in that formula the integrals are known, due to the following well-known, explicit formula:
$$\int_{S_N}g_{i_1j_1}\ldots g_{i_kj_k}=\begin{cases}
\frac{(N-|\ker i|)!}{N!}&{\rm if}\ \ker i=\ker j\\
0&{\rm otherwise}
\end{cases}$$

But this allows the computation of the right term, via the M\"obius inversion formula, and we get the result. As for the second assertion, this follows from the first one.
\end{proof}

In general, things are of course more complicated than this. We will need:

\index{distance on partitions}

\begin{proposition}
The following happen, for the partitions in $P(k)$:
\begin{enumerate}
\item $|\pi|+|\sigma|\leq|\pi\vee\sigma|+|\pi\wedge\sigma|$.

\item $|\pi\vee\tau|+|\tau\vee\sigma|\leq|\pi\vee\sigma|+|\tau|$.

\item $d(\pi,\sigma)=\frac{|\pi|+|\sigma|}{2}-|\pi\vee\sigma|$ is a distance.
\end{enumerate}
\end{proposition}

\begin{proof}
All this is well-known, the idea being as follows:

\medskip

(1) This comes from the fact that $P(k)$ is a semi-modular lattice.

\medskip

(2) This follows from (1), as explained for instance in \cite{ez1}.

\medskip

(3) This follows from (2), which says that the following holds:
$$\frac{|\pi|+|\tau|}{2}-d(\pi,\tau)+\frac{|\tau|+|\sigma|}{2}-d(\tau,\sigma)
\leq\frac{|\pi|+|\sigma|}{2}-d(\pi,\sigma)+|\tau|$$

Thus, we obtain the triangle inequality, and the other axioms are all clear.
\end{proof}

Actually in what follows we will only need (3) in the above statement. As a main result now regarding the Weingarten functions, we have:

\index{geodesicity defect}
\index{Weingarten matrix}

\begin{theorem}
The Weingarten matrix $W_{kN}$ has a series expansion in $N^{-1}$,
$$W_{kN}(\pi,\sigma)=N^{|\pi\vee\sigma|-|\pi|-|\sigma|}\sum_{g=0}^\infty K_g(\pi,\sigma)N^{-g}$$
where the various objects on the right are defined as follows:
\begin{enumerate}
\item A path from $\pi$ to $\sigma$ is a sequence $p=[\pi=\tau_0\neq\tau_1\neq\ldots\neq\tau_r=\sigma]$.

\item The signature of such a path is $+$ when $r$ is even, and $-$ when $r$ is odd.

\item The geodesicity defect of such a path is $g(p)=\sum_{i=1}^rd(\tau_{i-1},\tau_i)-d(\pi,\sigma)$.

\item $K_g$ counts the signed paths from $\pi$ to $\sigma$, with geodesicity defect $g$.
\end{enumerate} 
\end{theorem}

\begin{proof}
The Gram matrix can be written in the following way:
\begin{eqnarray*}
G_{kN}(\pi,\sigma)
&=&N^{|\pi\vee\sigma|}\\
&=&N^{\frac{|\pi|}{2}}N^{|\pi\vee\sigma|-\frac{|\pi|+|\sigma|}{2}}N^{\frac{|\sigma|}{2}}\\
&=&N^{\frac{|\pi|}{2}}N^{-d(\pi,\sigma)}N^{\frac{|\sigma|}{2}}
\end{eqnarray*}

This suggests considering the following diagonal matrix:
$$\Delta=diag(N^{\frac{|\pi|}{2}})$$

So, let us do this, and consider as well the following matrix:
$$H(\pi,\sigma)=\begin{cases}
0&(\pi=\sigma)\\
N^{-d(\pi,\sigma)}&(\pi\neq\sigma)
\end{cases}$$

In terms of these two matrices, the above formula for $G_{kN}$ simply reads:
$$G_{kN}=\Delta(1+H)\Delta$$

Thus, the Weingarten matrix $W_{kN}$ is given by the following formula:
$$W_{kN}=\Delta^{-1}(1+H)^{-1}\Delta^{-1}$$

Consider now the set $P_r(\pi,\sigma)$ of length $r$ paths between $\pi$ and $\sigma$. We have:
\begin{eqnarray*}
H^r(\pi,\sigma)
&=&\sum_{p\in P_r(\pi,\sigma)}H(\tau_0,\tau_1)\ldots H(\tau_{r-1},\tau_r)\\
&=&\sum_{p\in P_r(\pi,\sigma)}N^{-d(\pi,\sigma)-g(p)}
\end{eqnarray*}

Thus by using $(1+H)^{-1}=1-H+H^2-H^3+\ldots$ we obtain:
\begin{eqnarray*}
(1+H)^{-1}(\pi,\sigma)
&=&\sum_{r=0}^\infty(-1)^rH^r(\pi,\sigma)\\
&=&N^{-d(\pi,\sigma)}\sum_{r=0}^\infty\sum_{p\in P_r(\pi,\sigma)}(-1)^rN^{-g(p)}
\end{eqnarray*}

It follows that the Weingarten matrix is given by the following formula:
\begin{eqnarray*}
W_{kN}(\pi,\sigma)
&=&\Delta^{-1}(\pi)(1+H)^{-1}(\pi,\sigma)\Delta^{-1}(\sigma)\\
&=&N^{-\frac{|\pi|}{2}-\frac{|\sigma|}{2}-d(\pi,\sigma)}\sum_{r=0}^\infty\sum_{p\in P_r(\pi,\sigma)}(-1)^rN^{-g(p)}\\
&=&N^{|\pi\vee\sigma|-|\pi|-|\sigma|}\sum_{r=0}^\infty\sum_{p\in P_r(\pi,\sigma)}(-1)^rN^{-g(p)}
\end{eqnarray*}

Now by rearranging the various terms in the above double sum according to their geodesicity defect $g=g(p)$, this gives the following formula:
$$W_{kN}(\pi,\sigma)=N^{|\pi\vee\sigma|-|\pi|-|\sigma|}\sum_{g=0}^\infty K_g(\pi,\sigma)N^{-g}$$

Thus, we are led to the conclusion in the statement.
\end{proof}

As an illustration for all this, we have the following explicit estimates:

\begin{theorem}
Consider an easy quantum group $G=(G_N)$, coming from a category of partitions $D=(D(k))$. For any $\pi\leq\nu$ we have the estimate
$$W_{kN}(\pi,\sigma)=N^{-|\pi|}(\mu(\pi,\sigma)+O(N^{-1}))$$
and for $\pi,\sigma$ arbitrary we have
$$W_{kN}(\pi,\sigma)=O(N^{|\pi\vee\sigma|-|\pi|-|\sigma|})$$
with $\mu$ being the M\"obius function of $D(k)$.
\end{theorem}

\begin{proof}
We have two assertions here, the idea being as follows:

\medskip

(1) The first estimate is clear from the formula in Theorem 6.21, namely: 
$$W_{kN}(\pi,\sigma)=N^{|\pi\vee\sigma|-|\pi|-|\sigma|}\sum_{g=0}^\infty K_g(\pi,\sigma)N^{-g}$$

(2) In the case $\pi\leq\sigma$ it is known that $K_0$ coincides with the M\"obius function of $NC(k)$, as explained for instance in \cite{ez1}, so we obtain once again from Theorem 6.21 the fine estimate in the statement as well, namely:
$$W_{kN}(\pi,\sigma)=N^{-|\pi|}(\mu(\pi,\sigma)+O(N^{-1}))\qquad\forall\pi\leq\sigma$$

Observe that, by symmetry of $W_{kN}$, we obtain as well that we have:
$$W_{kN}(\pi,\sigma)=N^{-|\nu|}(\mu(\nu,\sigma)+O(N^{-1}))\qquad\forall\pi\geq\sigma$$

Thus, we are led to the conclusions in the statement.
\end{proof}

In the case of a category of partitions which is stable under removing blocks, the above estimates can be improved, and lead to the following result:

\index{uniform quantum group}
\index{liberation}
\index{truncated character}
\index{Bercovici-Pata bijection}

\begin{theorem}
For a liberation of uniform easy groups $G_N\to G_N^+$, the laws of truncated characters
$$\chi_t=\sum_{i=1}^{[tN]}u_{ii}$$
are in Bercovici-Pata bijection, in the $N\to\infty$ limit.
\end{theorem} 

\begin{proof}
As a first observation, we already know that the result holds, because we classified in the above all the uniform easy groups $G_N$, and for all the objects found, we have the asymptotic laws of the truncated characters $\chi_t$ computed, for both $G_N,G_N^+$, and with the Bercovici-Pata correspondence verified. However, and here comes our point, we can prove as well this result abstractly, without using the classification, as follows:

\medskip

(1) We have the following computation, to start with, for any $s\in\{1,\ldots,N\}$:
\begin{eqnarray*}
\int_{G_N}(u_{11}+\ldots +u_{ss})^k
&=&\sum_{i_1=1}^{s}\ldots\sum_{i_k=1}^s\int_{G_N}u_{i_1i_1}\ldots u_{i_ki_k}\\
&=&\sum_{\pi,\sigma\in D(k)}W_{kN}(\pi,\sigma)\sum_{i_1=1}^{s}\ldots\sum_{i_k=1}^s\delta_\pi(i)\delta_\sigma(i)\\
&=&\sum_{\pi,\sigma\in D(k)}W_{kN}(\pi,\sigma)G_{ks}(\sigma,\pi)\\
&=&Tr(W_{kN}G_{ks})
\end{eqnarray*}

(2) The point now is that we have the following estimates:
$$G_{kN}(\pi,\sigma):
\begin{cases}
=N^k&(\pi=\sigma)\\
\leq N^{k-1}&(\pi\neq\sigma)
\end{cases}$$

Thus with $N\to\infty$ we have the following estimate:
$$G_{kN}\sim N^k1$$

But this gives the following estimate, for our moment:
\begin{eqnarray*}
\int_{G_N}(u_{11}+\ldots +u_{ss})^k
&=&Tr(G_{kN}^{-1}G_{ks})\\
&\sim&Tr((N^k1)^{-1} G_{ks})\\
&=&N^{-k}Tr(G_{ks})\\
&=&N^{-k}s^k|D(k)|
\end{eqnarray*}

(3) With $s=[tN]$ and $N\to\infty$, the above formula gives:
$$\lim_{N\to\infty}\int_{G_N}\chi_t^k=\sum_{\pi\in D(k)}t^{|\pi|}$$

But this leads to the conclusion in the statement.
\end{proof}

There are many other things that can be said about the uniform easy quantum groups, and their characters and other variables, and we refer here to \cite{ez1}, \cite{ez2}, \cite{ez3} and related papers. As already mentioned, the uniform easy quantum groups play as well a key role in noncommutative geometry, via the associated homogeneous spaces, which again enjoy uniformity properties. We refer here to \cite{ba3} and the noncommutative geometry literature, and we will be back to this in chapter 16 below, with an introduction to the subject.

\section*{6e. Exercises} 

We have been here, in this chapter, into quite recent and specialized quantum group theory, all research matters, and the open questions abound. Basically no matter where you look, there are interesting exercises to start with, and beautiful theorems hidden behind them, if you are truly interested. So, we have no particular exercise to formulate, just look around, enjoy, and try what you like, and that will be an excellent exercise. However, since we had a Ground Zero theorem in this chapter, let us formulate:

\begin{exercise}
Start building on Ground Zero.
\end{exercise}

Instructions above, in the discussion after the Ground Zero theorem, the idea being that all this is mostly likely a beautiful mixture of abelian groups, finite groups, reflection groups, Lie groups and Lie algebras, and quantum extensions of these. However, before getting into this, better finish a first reading of the present book, because we still have to say a lot of interesting things, which can be useful in connection with such questions.

\chapter{Half-liberation}

\section*{7a. Half-liberation}

Let us go back to the standard cube of easy quantum groups, whose continuous face, the one on the right, we are mainly interested in, in this third part of this book:
$$
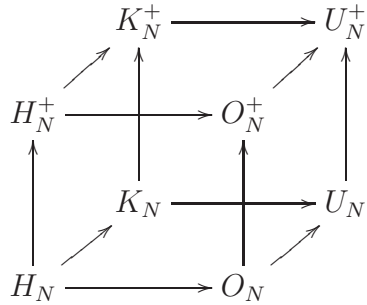
\xymatrix@R=15pt@C=15pt{
&K_N^+\ar[rr]&&U_N^+\\
H_N^+\ar[rr]\ar[ur]&&O_N^+\ar[ur]\\
&K_N\ar[rr]\ar[uu]&&U_N\ar[uu]\\
H_N\ar[uu]\ar[ur]\ar[rr]&&O_N\ar[uu]\ar[ur]
}$$

In order to further understand the right face, which is something quite broad, a natural idea would be that of cutting it in half, both on the horizontal and the vertical, with the help of some suitably chosen intermediate easy quantum groups, as follows: 
$$U_N\subset U_N^\times\subset U_N^+$$
$$O_N^+\subset\dot{O}_N^+\subset U_N^+$$

Indeed, assuming that we have such intermediate quantum groups, by intersecting we would reach to a diagram as follows, that we can further investigate afterwards: 
$$
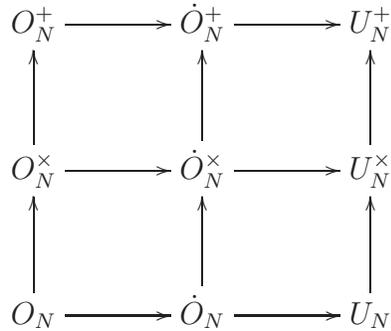
\xymatrix@R=12mm@C=14mm{
O_N^+\ar[r]&\dot{O}_N^+\ar[r]&U_N^+\\
O_N^\times\ar[r]\ar[u]&\dot{O}_N^\times\ar[r]\ar[u]&U_N^\times\ar[u]\\
O_N\ar[r]\ar[u]&\dot{O}_N\ar[r]\ar[u]&U_N\ar[u]}$$

So, this will be our plan. In practice now, finding intermediate quantum groups $U_N\subset U_N^\times\subset U_N^+$ and $O_N^+\subset\dot{O}_N^+\subset U_N^+$ as above can be a quite tricky business, due to the many possible choices for such quantum groups, and in order to do so, the best is to try to solve first the following question, which definitely looks simpler:

\begin{question}
What are the intermediate easy quantum groups
$$O_N\subset O_N^\times\subset O_N^+$$
$$O_N\subset\dot{O}_N\subset U_N$$
and can we really cut the continuous square, using these quantum groups?
\end{question}

In other words, our trick here is to say that any intermediate easy quantum groups $U_N\subset U_N^\times\subset U_N^+$ and $O_N^+\subset\dot{O}_N^+\subset U_N^+$ will produce by intersection intermediate easy quantum groups $O_N\subset O_N^\times\subset O_N^+$ and $O_N\subset\dot{O}_N\subset U_N$, in the obvious way, and so in order to find the former, the best is by trying to find first the latter.

\bigskip

Summarizing, we have here a valuable idea. And in addition, and we kept the good news for the end, Question 7.1 is in fact something that we already met, in chapter 6 when talking about Ground Zero, and we have in fact a full answer to it, as follows:

\begin{theorem}
The following happen:
\begin{enumerate}
\item There is only one intermediate easy quantum group $O_N\subset O_N^\times\subset O_N^+$, namely the half-classical orthogonal group $O_N^*$, appearing via the relations $abc=cba$.

\item The intermediate easy groups $O_N\subset\dot{O}_N\subset U_N$ are those of the form $\mathbb Z_rO_N$, with $r\in\{2,3,\ldots,\infty\}$, all being subgroups of $\mathbb TO_N$, appearing at $r=\infty$.
\end{enumerate}
\end{theorem}

\begin{proof}
These are results that we basically know, from chapter 3 and chapter 6, with the whole story, including some missing details to be given now, being as follows:

\medskip

(1) Regarding the quantum group $O_N^*$, we know from chapter 3 that we have indeed such a quantum group, by performing the following construction:
$$C(O_N^*)=C(O_N^+)\Big/\left<abc=cba\Big|\forall a,b,c\in\{u_{ij}\}\right>$$ 

We also know that $O_N^*$ is easy, because imposing the half-commutation relations $abc=cba$ amounts in imposing the condition $T_\pi\in End(u^{\otimes 3})$, with $\pi$ being as follows:
$$
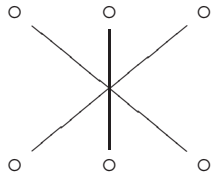
\xymatrix@R=16mm@C=8mm{
\circ\ar@{-}[drr]&\circ\ar@{-}[d]&\circ\ar@{-}[dll]\\
\circ&\circ&\circ}$$

Finally, in order to compute the associated category of pairings, $NC_2\subset P_2^*\subset P_2$, let us label the legs of $\pi$ clockwise $\circ\bullet\circ\bullet\ldots\,$. The diagram that we get is as follows:
$$
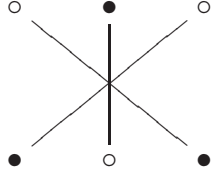
\xymatrix@R=16mm@C=8mm{
\circ\ar@{-}[drr]&\bullet\ar@{-}[d]&\circ\ar@{-}[dll]\\
\bullet&\circ&\bullet}$$

We can see that any string of $\pi$ joins $\circ-\bullet$, and with a bit of combinatorial work, which is routine, we conclude that $P_2^*$ consists of the parings having the property that when relabelling clockwise the legs $\circ\bullet\circ\bullet\ldots\,$, each string joins $\circ-\bullet$.

\medskip

(2) Regarding now the uniqueness of $O_N^*$, as an intermediate easy quantum group $O_N\subset O_N^\times\subset O_N^+$, this is something that we know well from chapter 6. The idea there was that we have to compute the intermediate categories of pairings, as follows:
$$NC_2\subset D\subset P_2$$

But in order to solve this problem, the key ingredient, obtained via semicircle capping, was the fact that given a pairing $\pi\in P_2(k,l)-NC_2(k,l)$, the following happen:

\medskip

-- If $\pi\in P_2(k, l)-P_2^*(k,l)$ then $<\pi>=P_2$.

\medskip

-- If $\pi\in P_2^*(k,l)-NC_2(k,l)$ then $<\pi>=P_2^*$.

\medskip

With this is hand, we conclude that if $D\not\subset P_2^*$, we obtain $G=O_N$, that if $D\subset P_2,D\not\subset NC_2$, we obtain $G=O_N^*$, and that if $D\subset NC_2$, we obtain $G=O_N^+$.

\medskip

(3) Let us investigate now our second question, regarding the intermediate easy groups which are ``hybrid'' between real and complex, as follows:
$$O_N\subset\dot{O}_N\subset U_N$$

This is actually something that we briefly talked about, right before in chapter 6, but without many details, with our discussion there having been simply in order to prove our Ground Zero theorem there. So, time now to review this material, with full proofs.

\medskip

(4) Our first claim is that the group $\mathbb TO_N\subset U_N$ is easy, the corresponding category of partitions being the subcategory $\bar{P}_2\subset P_2$ consisting of the pairings having the property that when flatenning, we have the following global formula:
$$\#\circ=\#\bullet$$

Indeed, if we denote the standard corepresentation by $u=zv$, with $z\in\mathbb T$ and with $v=\bar{v}$, then in order to have $Hom(u^{\otimes k},u^{\otimes l})\neq\emptyset$, the $z$ variabes must cancel, and in the case where they cancel, we obtain the same Hom-space as for $O_N$. Now since the cancelling property for the $z$ variables corresponds precisely to the fact that $k,l$ must have the same numbers of $\circ$ symbols minus $\bullet$ symbols, the associated Tannakian category must come from the category of pairings $\bar{P}_2\subset P_2$, as claimed.

\medskip

(5) Our second claim is that, more generally, the group $\mathbb Z_rO_N\subset U_N$ with $r\in\{2,3,\ldots,\infty\}$ is easy, with the corresponding category $P_2^r\subset P_2$ consisting of the pairings having the property that when flatenning, we have the following global formula:
$$\#\circ=\#\bullet(r)$$

Indeed, this is something that we already know at $r=1,\infty$, where the group in question is $O_N,\mathbb TO_N$. The proof in general is similar, by writing $u=zv$ as above.

\medskip

(6) Summarizing, done with the existence part. For further reference, let us record as well the fact, which is elementary, that we have a presentation result as follows:
$$C(\mathbb TO_N)=C(U_N)\Big/\left<ab^*=a^*b\Big|\forall a,b\in\{u_{ij}\}\right>$$ 

Equivalently, in terms of diagrams, we can say that $\mathbb TO_N\subset U_N$ appears by imposing the condition $T_\pi\in End(u\otimes\bar{u},\bar{u}\otimes u)$, with $\pi$ being as follows:
$$
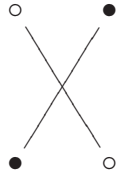
\xymatrix@R=16mm@C=8mm{
\circ\ar@{-}[dr]&\bullet\ar@{-}[dl]\\
\bullet&\circ}$$

We will be actually back to all this in a moment, with the comment that the projective version of $\mathbb TO_N$ is the group $PO_N$, and with $\mathbb TO_N$ being maximal with this property.

\medskip

(7) In order to finish, let us prove now the uniqueness result, stating that the above groups $O_N\subset\mathbb Z_rO_N\subset U_N$ are the only intermediate easy groups $O_N\subset G\subset U_N$. According to our conventions for the easy quantum groups, which apply of course to the classical case, we must compute the following intermediate categories of pairings:
$$\mathcal P_2\subset D\subset P_2$$

So, assume that we have such a category, $D\neq\mathcal P_2$, and pick an element $\pi\in D-\mathcal P_2$, assumed to be flat. We can modify $\pi$, by performing the following operations:

\medskip

-- First, we can compose with the basic crossing, in order to assume that $\pi$ is a partition of type $\cap\ldots\ldots\cap$, consisting of consecutive semicircles. Our assumption $\pi\notin\mathcal P_2$ means that at least one semicircle is colored black, or white.

\medskip

-- Second, we can use the basic mixed-colored semicircles, and cap with them all the mixed-colored semicircles. Thus, we can assume that $\pi$ is a nonzero partition of type $\cap\ldots\ldots\cap$, consisting of consecutive black or white semicircles.

\medskip

-- Third, we can rotate, as to assume that $\pi$ is a partition consisting of an upper row of white semicircles, $\cup\ldots\ldots\cup$, and a lower row of white semicircles, $\cap\ldots\ldots\cap$. Our assumption $\pi\notin\mathcal P_2$ means that this latter partition is nonzero.

\medskip

(8) For any two integers $a,b\in\mathbb N$ consider the partition consisting of an upper row of $a$ white semicircles, and a lower row of $b$ white semicircles, and set:
$$\mathcal C=\left\{\pi_{ab}\Big|a,b\in\mathbb N\right\}\cap D$$

According to the above we have $\pi\in<\mathcal C>$. The point now is that we have:

\medskip

-- There exists $r\in\mathbb N\cup\{\infty\}$ such that $\mathcal C$ equals the following set:
$$\mathcal C_r=\left\{\pi_{ab}\Big|a=b(r)\right\}$$

This is indeed standard, by using the categorical axioms.

\medskip

-- We have the following formula, with $P_2^r$ being as above:
$$<\mathcal C_r>=P_2^r$$

This is standard as well, by doing some diagrammatic work.

\medskip

(9) With these results in hand, the conclusion now follows. Indeed, with $r\in\mathbb N\cup\{\infty\}$ being as above, we know from the beginning of the proof that any $\pi\in D$ satisfies:
$$\pi
\in<\mathcal C>
=<\mathcal C_r>
=P_2^r$$

Thus we have an inclusion $D\subset P_2^r$. Conversely, we have as well:
\begin{eqnarray*}
P_2^r
&=&<\mathcal C_r>\\
&=&<\mathcal C>\\
&\subset&<D>\\
&=&D
\end{eqnarray*}

Thus we have $D=P_2^r$, and this finishes the proof of the uniqueness assertion.
\end{proof}

All this is nice, with Theorem 7.2 providing a full answer to Question 7.1, or at least to the first part of that question. In order to answer now the second part as well, we must suitably cut the right face of the standard cube, by using the quantum groups that we found, and the procedure here is straightforward, leading to the following result:

\begin{theorem}
We have easy quantum groups as follows,
$$\xymatrix@R=13mm@C=13mm{
O_N^+\ar[r]&\mathbb TO_N^+\ar[r]&U_N^+\\
O_N^*\ar[r]\ar[u]&\mathbb TO_N^*\ar[r]\ar[u]&U_N^*\ar[u]\\
O_N\ar[r]\ar[u]&\mathbb TO_N\ar[r]\ar[u]&U_N\ar[u]}$$
with the half-liberations obtained via the relations $abc=cba$, applied to $u_{ij},u_{ij}^*$, and with the hybrids obtained via the relations $ab^*=a^*b$, applied to $u_{ij}$.
\end{theorem}

\begin{proof}
This is standard from what we have in Theorem 7.2 and its proof, with the corresponding categories of partitions being as follows:
$$\xymatrix@R=11mm@C=11mm{
NC_2\ar[d]&\bar{NC}_2\ar[l]\ar[d]&\mathcal{NC}_2\ar[d]\ar[l]\\
P_2^*\ar[d]&\bar{P}_2^*\ar[l]\ar[d]&\mathcal P_2^*\ar[l]\ar[d]\\
P_2&\bar{P}_2\ar[l]&\mathcal P_2\ar[l]}$$

To be more precise, we already have the results, in most of the cases, and the cases left follow by using the categorical properties of the intersection operation $\cap$.
\end{proof}

Quite remarkably, we have as well a discrete version of the above result, as follows:

\begin{theorem}
We have easy quantum groups as follows,
$$\xymatrix@R=12mm@C=12mm{
H_N^+\ar[r]&\mathbb TH_N^+\ar[r]&K_N^+\\
H_N^*\ar[r]\ar[u]&\mathbb TH_N^*\ar[r]\ar[u]&K_N^*\ar[u]\\
H_N\ar[r]\ar[u]&\mathbb TH_N\ar[r]\ar[u]&K_N\ar[u]}$$
with the half-liberations obtained via the relations $abc=cba$, applied to $u_{ij},u_{ij}^*$, and with the hybrids obtained via the relations $ab^*=a^*b$, applied to $u_{ij}$.
\end{theorem}

\begin{proof}
This is again standard, the categories of partitions being as follows:
$$
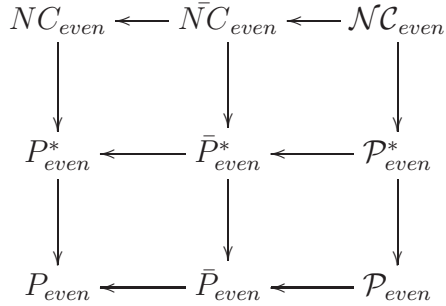
\xymatrix@R=11mm@C=7mm{
NC_{even}\ar[d]&\bar{NC}_{even}\ar[l]\ar[d]&\mathcal{NC}_{even}\ar[d]\ar[l]\\
P_{even}^*\ar[d]&\bar{P}_{even}^*\ar[l]\ar[d]&\mathcal P_{even}^*\ar[l]\ar[d]\\
P_{even}&\bar{P}_{even}\ar[l]&\mathcal P_{even}\ar[l]}$$

As for the proof, again this follows from what we have, via a few computations.
\end{proof}

Finally, let us record as well the following result, dealing with the whole cube:

\begin{theorem}
We can cut the standard cube of easy quantum groups
$$
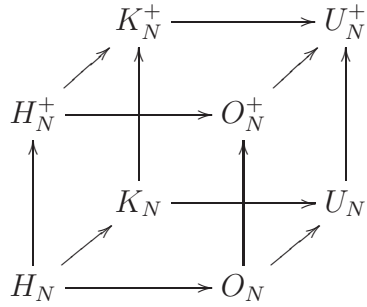
\xymatrix@R=15pt@C=15pt{
&K_N^+\ar[rr]&&U_N^+\\
H_N^+\ar[rr]\ar[ur]&&O_N^+\ar[ur]\\
&K_N\ar[rr]\ar[uu]&&U_N\ar[uu]\\
H_N\ar[uu]\ar[ur]\ar[rr]&&O_N\ar[uu]\ar[ur]
}$$
by using Theorems 7.3 and 7.4, and we obtain an intersection/easy generation diagram.
\end{theorem}

\begin{proof}
We recall from chapter 2 that the categories of partitions for the basic examples of easy quantum groups, which form the standard cube, are as follows:
$$
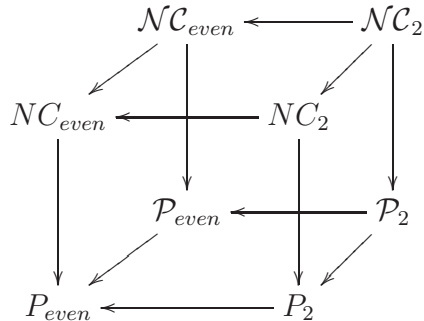
\xymatrix@R=19pt@C5pt{
&\mathcal{NC}_{even}\ar[dl]\ar[dd]&&\mathcal {NC}_2\ar[dl]\ar[ll]\ar[dd]\\
NC_{even}\ar[dd]&&NC_2\ar[dd]\ar[ll]\\
&\mathcal P_{even}\ar[dl]&&\mathcal P_2\ar[dl]\ar[ll]\\
P_{even}&&P_2\ar[ll]
}$$

The point now is that when cutting this cube, using the categories of partitions from the proofs of Theorems 7.3 and 7.4, we obtain intersection/generation diagram. Thus, at the quantum group level, we have an intersection/easy generation diagram, as stated.
\end{proof}

\section*{7b. Matrix models}

Generally speaking, the most powerful tool for the study of the half-liberations are the matrix models. Let us start with some generalities. We first have:

\index{random matrix}
\index{matrix model}

\begin{definition}
A matrix model for $G\subset U_N^+$ is a morphism of $C^*$-algebras
$$\pi:C(G)\to M_K(C(T))$$
where $K\geq1$ is an integer, and $T$ is a compact space.
\end{definition}

The simplest situation is when $\pi$ is faithful in the usual sense. Here $\pi$ obviously reminds $G$. However, this is something quite restrictive, because in this case the algebra $C(G)$ must be quite small, admitting an embedding as follows:
$$\pi:C(G)\subset M_K(C(T))$$

Technically, this means that $C(G)$ must be of type I, as an operator algebra, and this is indeed something quite restrictive. However, there are many interesting examples here, including our half-liberations, and all this is worth a detailed look. We have:

\index{stationary model}
\index{Haar functional}

\begin{definition}
A matrix model $\pi:C(G)\to M_K(C(T))$ is called stationary when
$$\int_G=\left(tr\otimes\int_T\right)\pi$$
where $\int_T$ is the integration with respect to a given probability measure on $T$.
\end{definition}

Here the term ``stationary'' comes from a functional analytic interpretation of all this, with a certain Ces\`aro limit being needed to be stationary, and this will be explained later. Yet another explanation comes from a certain relation with the lattice models, but this relation is rather something folklore, not axiomatized yet. We will be back to this.

\bigskip

The relation between stationarity and faithfulness comes from:

\index{coamenability}
\index{type I algebra}

\begin{proposition}
Assuming that a closed subgroup $G\subset U_N^+$ has a stationary model
$$\pi:C(G)\to M_K(C(T))$$
it follows that $G$ must be coamenable, and that the model is faithful. Moreover, $\pi$ extends into an embedding of von Neumann algebras, as follows,
$$L^\infty(G)\subset M_K(L^\infty(T))$$
which commutes with the canonical integration functionals.
\end{proposition}

\begin{proof}
By performing the GNS construction with respect to $\int_G$, we obtain a factorization as follows, which commutes with the canonical integration functionals:
$$\pi:C(G)\to C(G)_{red}\subset M_K(C(T))$$

Thus, in what regards the coamenability question, we can assume that $\pi$ is faithful. With this assumption made, we have an embedding as follows:
$$C(G)\subset M_K(C(T))$$

Now observe that the GNS construction gives a better embedding, as follows:
$$L^\infty(G)\subset M_K(L^\infty(T))$$

Now since the von Neumann algebra on the right is of type I, so must be its subalgebra $A=L^\infty(G)$. But this means that, when writing the center of this latter algebra as  $Z(A)=L^\infty(X)$, the whole algebra decomposes over $X$, as an integral of type I factors:
$$L^\infty(G)=\int_XM_{K_x}(\mathbb C)\,dx$$

In particular, we can see from this that $C(G)\subset L^\infty(G)$ has a unique $C^*$-norm, and so $G$ is coamenable. Thus we have proved our first assertion, and the second assertion follows as well, because our factorization of $\pi$ consists of the identity, and of an inclusion.
\end{proof}

The above notion of stationary model is very restrictive, and will not apply for instance to the free quantum groups, which are not coamenable. Fortunately, we have as well:

\begin{definition}
Let $\pi:C(G)\to M_K(C(T))$ be a matrix model. 
\begin{enumerate}
\item The Hopf image of $\pi$ is the smallest quotient Hopf $C^*$-algebra $C(G)\to C(H)$ producing a factorization as follows:
$$\pi:C(G)\to C(H)\to M_K(C(T))$$

\item When the inclusion $H\subset G$ is an isomorphism, i.e. when there is no non-trivial factorization as above, we say that $\pi$ is inner faithful.
\end{enumerate}
\end{definition}

Here the existence and uniqueness of the Hopf image come by dividing $C(G)$ by a suitable ideal, with this being standard. Alternatively, in Tannakian terms, we have:

\index{Tannakian duality}

\begin{theorem}
Assuming $G\subset U_N^+$, with fundamental corepresentation $u=(u_{ij})$, the Hopf image of a model $\pi:C(G)\to M_K(C(T))$ comes from the Tannakian category
$$C_{kl}=Hom(U^{\otimes k},U^{\otimes l})$$
where $U_{ij}=\pi(u_{ij})$, and where the spaces on the right are taken in a formal sense.
\end{theorem}

\begin{proof}
Since the morphisms increase the intertwining spaces, when defined either in a representation theory sense, or just formally, we have inclusions as follows:
$$Hom(u^{\otimes k},u^{\otimes l})\subset Hom(U^{\otimes k},U^{\otimes l})$$

More generally, we have such inclusions when replacing $(G,u)$ with any pair producing a factorization of $\pi$. Thus, by Tannakian duality, the Hopf image must be given by the fact that the intertwining spaces must be the biggest, subject to the above inclusions. On the other hand, since $u$ is biunitary, so is $U$, and it follows that the spaces on the right form a Tannakian category. Thus, we have a quantum group $(H,v)$ given by:
$$Hom(v^{\otimes k},v^{\otimes l})=Hom(U^{\otimes k},U^{\otimes l})$$

By the above discussion, $C(H)$ follows to be the Hopf image of $\pi$, as claimed.
\end{proof}

Regarding now the study of the inner faithful models, a key problem is that of computing the Haar integration functional. The result here is as follows:

\index{truncated integrals}

\begin{theorem}
Given an inner faithful model $\pi:C(G)\to M_K(C(T))$, we have
$$\int_G=\lim_{k\to\infty}\frac{1}{k}\sum_{r=1}^k\int_G^r$$
with the truncations of the integration on the right being given by
$$\int_G^r=(\varphi\circ\pi)^{*r}$$
with $\phi*\psi=(\phi\otimes\psi)\Delta$, and with $\varphi=tr\otimes\int_T$ being the random matrix trace.
\end{theorem}

\begin{proof}
This is something quite tricky, the idea being as follows:

\medskip

(1) In order to prove the result, we can proceed as in chapter 2. If we denote by $\int_G'$ the limit in the statement, we must prove that this limit converges, and that:
$$\int_G'=\int_G$$

It is enough to check this on the coefficients of the Peter-Weyl corepresentations, and if we let $v=u^{\otimes k}$ be one of these corepresentations, we must prove that we have:
$$\left(id\otimes\int_G'\right)v=\left(id\otimes\int_G\right)v$$

(2) We know that the matrix on the right is the projection onto $Fix(v)$:
$$\left(id\otimes\int_G\right)v=Proj\Big[Fix(v)\Big]$$

Regarding now the matrix on the left, the trick from \cite{wo1} applies, and gives:
$$\left(id\otimes\int_G'\right)v=Proj\Big[1\in (id\otimes\varphi\pi)v\Big]$$

(3) Now observe that, if we set $V_{ij}=\pi(v_{ij})$, we have the following formula:
$$(id\otimes\varphi\pi)v=(id\otimes\varphi)V$$

Thus, we can apply the trick in \cite{wo1}, and we conclude that the $1$-eigenspace that we are interested in equals $Fix(V)$. But, according to Theorem 7.10, we have:
$$Fix(V)=Fix(v)$$

Thus, we have proved that we have $\int_G'=\int_G$, as desired.
\end{proof}

Now back to the the stationary models, we have the following useful criterion:

\index{stationary on its image}
\index{Hopf image}

\begin{theorem}
For a model $\pi:C(G)\to M_K(C(T))$, the following are equivalent:
\begin{enumerate}
\item $Im(\pi)$ is a Hopf algebra, and the Haar integration on it is:
$$\psi=\left(tr\otimes\int_T\right)\pi$$

\item The linear form $\psi=(tr\otimes\int_T)\pi$ satisfies the idempotent state property:
$$\psi*\psi=\psi$$

\item We have $T_e^2=T_e$, $\forall p\in\mathbb N$, $\forall e\in\{1,*\}^p$, where:
$$(T_e)_{i_1\ldots i_p,j_1\ldots j_p}=\left(tr\otimes\int_T\right)(U_{i_1j_1}^{e_1}\ldots U_{i_pj_p}^{e_p})$$
\end{enumerate}
If these conditions are satisfied, we say that $\pi$ is stationary on its image.
\end{theorem}

\begin{proof}
Given a matrix model $\pi:C(G)\to M_K(C(T))$ as in the statement, we can factorize it via its Hopf image, as follows:
$$\pi:C(G)\to C(H)\to M_K(C(T))$$

Now observe that (1,2,3) above depend only on the factorized representation:
$$\nu:C(H)\to M_K(C(T))$$

Thus, we can assume $G=H$, which means that we can assume that $\pi$ is inner faithful. With this assumption made, the proof of the equivalences goes as follows:

\medskip

$(1)\implies(2)$ This is clear from definitions, because the Haar integration on any compact quantum group satisfies the idempotent state equation:
$$\psi*\psi=\psi$$

$(2)\implies(1)$ Assuming $\psi*\psi=\psi$, we have $\psi^{*r}=\psi$ for any $r\in\mathbb N$, and we obtain by taking a Ces\`aro limit that we have $\int_G=\psi$, which gives the result.

\medskip

In order to establish now $(2)\Longleftrightarrow(3)$, we use the following elementary formula, which comes from the definition of the convolution operation:
$$\psi^{*r}(u_{i_1j_1}^{e_1}\ldots u_{i_pj_p}^{e_p})=(T_e^r)_{i_1\ldots i_p,j_1\ldots j_p}$$

\medskip

$(2)\implies(3)$ Assuming $\psi*\psi=\psi$, by using the above formula at $r=1,2$ we obtain that the matrices $T_e$ and $T_e^2$ have the same coefficients, and so they are equal.

\medskip

$(3)\implies(2)$ Assuming $T_e^2=T_e$, by using the above formula at $r=1,2$ we obtain that the linear forms $\psi$ and $\psi*\psi$ coincide on any product of coefficients $u_{i_1j_1}^{e_1}\ldots u_{i_pj_p}^{e_p}$. Now since these coefficients span a dense subalgebra of $C(G)$, this gives the result.
\end{proof}

Now back to half-liberation, we first have the following result, from \cite{ez2}, \cite{bdu}:

\index{half-liberation}
\index{antidiagonal model}

\begin{proposition}
We have a matrix model as follows, 
$$C(O_N^*)\to M_2(C(U_N))\quad,\quad 
u_{ij}\to\begin{pmatrix}0&v_{ij}\\ \bar{v}_{ij}&0\end{pmatrix}$$
where $v$ is the fundamental corepresentation of $C(U_N)$, as well as a model as follows,
$$C(U_N^*)\to M_2(C(U_N\times U_N))\quad,\quad 
u_{ij}\to\begin{pmatrix}0&v_{ij}\\ w_{ij}&0\end{pmatrix}$$
where $v,w$ are the fundamental corepresentations of the two copies of $C(U_N)$.
\end{proposition}

\begin{proof}
It is routine to check that the matrices on the right are indeed biunitaries, and since the first matrix is also self-adjoint, we obtain in this way models as follows:
$$C(O_N^+)\to M_2(C(U_N))\quad,\quad
C(U_N^+)\to M_2(C(U_N\times U_N))$$

Regarding now the half-commutation relations, this comes from something general, regarding the antidiagonal $2\times2$ matrices. Consider indeed matrices as follows:
$$X_i=\begin{pmatrix}0&x_i\\ y_i&0\end{pmatrix}$$

We have then the following computation:
$$X_iX_jX_k
=\begin{pmatrix}0&x_i\\ y_i&0\end{pmatrix}\begin{pmatrix}0&x_j\\ y_j&0\end{pmatrix}\begin{pmatrix}0&x_k\\ y_k&0\end{pmatrix}
=\begin{pmatrix}0&x_iy_jx_k\\ y_ix_jy_k&0\end{pmatrix}$$

Since this quantity is symmetric in $i,k$, we obtain from this:
$$X_iX_jX_k=X_kX_jX_i$$

Thus, the antidiagonal $2\times2$ matrices half-commute, and we conclude that our models for $C(O_N^+)$ and $C(U_N^+)$ constructed above factorize as in the statement.
\end{proof}

We can now formulate our first concrete modelling theorem, as follows:

\begin{theorem}
The above antidiagonal models, namely
$$C(O_N^*)\to M_2(C(U_N))\quad,\quad 
C(U_N^*)\to M_2(C(U_N\times U_N))$$
are both stationary, and in particular they are faithful.
\end{theorem}

\begin{proof}
Let us first discuss the case of $O_N^*$. We will use Theorem 7.12 (3). Since the fundamental representation is self-adjoint, the various matrices $T_e$ with $e\in\{1,*\}^p$ are all equal. We denote this common matrix by $T_p$. We have, by definition:
$$(T_p)_{i_1\ldots i_p,j_1\ldots j_p}
=\left(tr\otimes\int_H\right)\left[\begin{pmatrix}0&v_{i_1j_1}\\\bar{v}_{i_1j_1}&0\end{pmatrix}\ldots\ldots\begin{pmatrix}0&v_{i_pj_p}\\\bar{v}_{i_pj_p}&0\end{pmatrix}\right]$$

Since when multipliying an odd number of antidiagonal matrices we obtain an atidiagonal matrix, we have $T_p=0$ for $p$ odd. Also, when $p$ is even, we have:
\begin{eqnarray*}
(T_p)_{i_1\ldots i_p,j_1\ldots j_p}
&=&\left(tr\otimes\int_H\right)\begin{pmatrix}v_{i_1j_1}\ldots\bar{v}_{i_pj_p}&0\\0&\bar{v}_{i_1j_1}\ldots v_{i_pj_p}\end{pmatrix}\\
&=&\frac{1}{2}\left(\int_Hv_{i_1j_1}\ldots\bar{v}_{i_pj_p}+\int_H\bar{v}_{i_1j_1}\ldots v_{i_pj_p}\right)\\
&=&\int_HRe(v_{i_1j_1}\ldots\bar{v}_{i_pj_p})
\end{eqnarray*}

We have $T_p^2=T_p=0$ when $p$ is odd, so we are left with proving that for $p$ even we have $T_p^2=T_p$. For this purpose, we use the following formula:
$$Re(x)Re(y)=\frac{1}{2}\left(Re(xy)+Re(x\bar{y})\right)$$

By using this identity for each of the terms which appear in the product, and multi-index notations in order to simplify the writing, we obtain:
\begin{eqnarray*}
(T_p^2)_{ij}
&=&\sum_{k_1\ldots k_p}(T_p)_{i_1\ldots i_p,k_1\ldots k_p}(T_p)_{k_1\ldots k_p,j_1\ldots j_p}\\
&=&\int_H\int_H\sum_{k_1\ldots k_p}Re(v_{i_1k_1}\ldots\bar{v}_{i_pk_p})Re(w_{k_1j_1}\ldots\bar{w}_{k_pj_p})dvdw\\
&=&\frac{1}{2}\int_H\int_H\sum_{k_1\ldots k_p}Re(v_{i_1k_1}w_{k_1j_1}\ldots\bar{v}_{i_pk_p}\bar{w}_{k_pj_p})+Re(v_{i_1k_1}\bar{w}_{k_1j_1}\ldots\bar{v}_{i_pk_p}w_{k_pj_p})dvdw\\
&=&\frac{1}{2}\int_H\int_HRe((vw)_{i_1j_1}\ldots(\bar{v}\bar{w})_{i_pj_p})+Re((v\bar{w})_{i_1j_1}\ldots(\bar{v}w)_{i_pj_p})dvdw
\end{eqnarray*}

Now since $vw\in H$ is uniformly distributed when $v,w\in H$ are uniformly distributed, the quantity on the left integrates up to $(T_p)_{ij}$. Also, since $H$ is conjugation-stable, $\bar{w}\in H$ is uniformly distributed when $w\in H$ is uniformly distributed, so the quantity on the right integrates up to the same quantity, namely $(T_p)_{ij}$. Thus, we have:
$$(T_p^2)_{ij}
=\frac{1}{2}\Big((T_p)_{ij}+(T_p)_{ij}\Big)
=(T_p)_{ij}$$

Summarizing, we have obtained that for any $p$, we have $T_p^2=T_p$. Thus Theorem 7.12 applies, and shows that our model is stationary, as claimed. As for the proof of the stationarity for the model for $U_N^*$, this is similar. 
\end{proof}

As a second illustration, regarding $H_N^*,K_N^*$, we have:

\begin{theorem}
We have a stationary matrix model as follows, 
$$C(H_N^*)\to M_2(C(K_N))\quad,\quad 
u_{ij}\to\begin{pmatrix}0&v_{ij}\\ \bar{v}_{ij}&0\end{pmatrix}$$
where $v$ is the fundamental corepresentation of $C(K_N)$, as well as a stationary model
$$C(K_N^*)\to M_2(C(K_N\times K_N))\quad,\quad 
u_{ij}\to\begin{pmatrix}0&v_{ij}\\ w_{ij}&0\end{pmatrix}$$
where $v,w$ are the fundamental corepresentations of the two copies of $C(K_N)$.
\end{theorem}

\begin{proof}
This follows by adapting the proof of Proposition 7.13 and Theorem 7.14, by adding there the $H_N^+,K_N^+$ relations. All this is in fact part of a more general phenomenon, concerning half-liberation in general, and we refer here to \cite{bdu}. We will be back to this.
\end{proof}

\section*{7c. Representation theory} 

Let us discuss now the modern approach to half-liberation, following Bichon and Dubois-Violette \cite{bdu}, based on crossed products and related $2\times2$ matrix models:

\index{Bichon-Dubois-Violette}
\index{half-liberation}

\begin{theorem}
Given a conjugation-stable closed subgroup $H\subset U_N$, consider the algebra $C([H])\subset M_2(C(H))$ generated by the following variables:
$$u_{ij}=\begin{pmatrix}0&v_{ij}\\ \bar{v}_{ij}&0\end{pmatrix}$$
Then $[H]$ is a compact quantum group, we have $[H]\subset O_N^*$, and any non-classical subgroup $G\subset O_N^*$ appears in this way, with $G=O_N^*$ itself appearing from $H=U_N$.
\end{theorem}

\begin{proof}
We have several things to be proved, the idea being as follows:

\medskip

(1) As a first observation, the matrices in the statement are self-adjoint. Let us prove now that these matrices are orthogonal. We have:
$$\sum_ku_{ik}u_{jk}
=\sum_k\begin{pmatrix}v_{ik}\bar{v}_{jk}&0\\ 0&\bar{v}_{ik}v_{jk}\end{pmatrix}
=\begin{pmatrix}1&0\\0&1\end{pmatrix}$$

In the other sense, the computation is similar, as follows:
$$\sum_ku_{ki}u_{kj}
=\sum_k\begin{pmatrix}v_{ki}\bar{v}_{kj}&0\\ 0&\bar{v}_{ki}v_{kj}\end{pmatrix}
=\begin{pmatrix}1&0\\0&1\end{pmatrix}$$

(2) Our second claim is that the matrices in the statement half-commute. Consider indeed arbitrary antidiagonal $2\times2$ matrices, with commuting entries, as follows:
$$X_i=\begin{pmatrix}0&x_i\\ y_i&0\end{pmatrix}$$

We have then the following computation:
$$X_iX_jX_k
=\begin{pmatrix}0&x_i\\ y_i&0\end{pmatrix}\begin{pmatrix}0&x_j\\ y_j&0\end{pmatrix}\begin{pmatrix}0&x_k\\ y_k&0\end{pmatrix}
=\begin{pmatrix}0&x_iy_jx_k\\ y_ix_jy_k&0\end{pmatrix}$$

Since this quantity is symmetric in $i,k$, we obtain, as desired:
$$X_iX_jX_k=X_kX_jX_i$$

(3) According now to the definition of the quantum group $O_N^*$, we have a representation of algebras, as follows where $w$ is the fundamental corepresentation of $C(O_N^*)$:
$$\pi:C(O_N^*)\to M_2(C(H))\quad,\quad 
w_{ij}\to u_{ij}$$

Thus, with the compact quantum space $[H]$ being constructed as in the statement, we have a representation of algebras, as follows:
$$\rho:C(O_N^*)\to C([H])\quad,\quad 
w_{ij}\to u_{ij}$$

(4) With this in hand, it is routine to check that the compact quantum space $[H]$ constructed in the statement is indeed a compact quantum group, with this being best viewed via an equivalent construction, with a quantum group embedding as follows:
$$C([H])\subset C(H)\rtimes\mathbb Z_2$$

(5) As for the proof of the converse, stating that any non-classical subgroup $G\subset O_N^*$ appears in this way, this is something more tricky, and we refer here to \cite{bdu}. 

\medskip

(6) Finally, for the fact that we have indeed $O_N^*=[U_N]$, we refer here as well to \cite{bdu}, and we will be back to this as well in chapter 8 below.
\end{proof}

In relation with the above, we will need as well the following result, regarding the irreducible corepresentations, also from Bichon-Dubois-Violette \cite{bdu}:

\begin{theorem}
In the context of the correspondence $H\to[H]$ we have a bijection 
$$Irr([H])\simeq Irr_0(H)\coprod Irr_1(H)$$
where the sets on the right are given by
$$Irr_k(H)=\left\{r\in Irr(H)\Big|\exists l\in\mathbb N,r\in u^{\otimes k}\otimes(u\otimes\bar{u})^{\otimes l}\right\}$$
induced by the canonical identification $Irr(H\rtimes\mathbb Z_2)\simeq Irr(H)\coprod Irr(H)$.
\end{theorem}

\begin{proof}
This is something more technical, also from \cite{bdu}. It is easy to see that we have an equality of projective versions $P[H]=PH$, which gives an inclusion as follows:
$$Irr_0(H)=Irr(PH)\subset Irr([H])$$

As for the remaining irreducible representations of $[H]$, these must come from an inclusion $Irr_1(H)\subset Irr([H])$, appearing as above. See \cite{bdu}.
\end{proof}

In relation with the maximal tori, the situation here is very simple, as follows:

\index{group dual}

\begin{theorem}
The group dual subgroups $\widehat{[\Gamma]}_Q\subset[H]$ appear via
$$[\Gamma]_Q=[\Gamma_Q]$$
from the group dual subgroups $\widehat{\Gamma}_Q\subset H$ associated to $H\subset U_N$.
\end{theorem}

\begin{proof}
Let us first discuss the case $Q=1$. Consider the diagonal subgroup $\widehat{\Gamma}_1\subset H$, with the associated quotient map $C(H)\to C(\widehat{\Gamma}_1)$ denoted:
$$v_{ij}\to\delta_{ij}h_i$$

At the level of the algebras of $2\times2$ matrices, this map induces a quotient map:
$$M_2(C(H))\to M_2(C(\widehat{\Gamma}_1))$$

Our claim is that we have a factorization, as follows:
$$\begin{matrix}
C([H])&\subset&M_2(C(H))\\
\\
\downarrow&&\downarrow\\
\\
C([\widehat{\Gamma}_1])&\subset&M_2(C(\widehat{\Gamma}_1))
\end{matrix}$$

Indeed, it is enough to show that the standard generators of $C([H])$ and of $ C([\widehat{\Gamma}_1])$ map to the same elements of $M_2(C(\widehat{\Gamma}_1))$. But these generators map indeed as follows:
$$\begin{matrix}
u_{ij}&\to&\begin{pmatrix}0&v_{ij}\\ \bar{v}_{ij}&0\end{pmatrix}\\
\\
&&\downarrow\\
\\
\delta_{ij}v_{ij}&\to&\begin{pmatrix}0&\delta_{ij}h_i\\ \delta_{ij}h_i^{-1}&0\end{pmatrix}
\end{matrix}$$

Thus we have the above factorization, and since the map on the left is obtained by imposing the relations $u_{ij}=0$ with $i\neq j$, we obtain, as desired:
$$[\Gamma]_1=[\Gamma_1]$$

In the general case now, $Q\in U_N$, the result follows by applying the above $Q=1$ result to the quantum group $[H]$, with fundamental corepresentation $w=QuQ^*$.
\end{proof}

There are many other things that can be said about half-liberations, the idea being that the half-classical geometry is quite close to the classical geometry, due to the above $2\times2$ matrix modelling results. For more on this, we refer to \cite{bdu} and related papers.

\section*{7d. Maximality results} 

Going back to what was said in the beginning of this chapter, we discuss now some interesting related questions, going beyond easiness. Following \cite{bc+}, we first have:

\begin{theorem}
The following inclusion of compact groups is maximal,
$$\mathbb TO_N\subset U_N$$
in the sense that there is no intermediate compact group in between.
\end{theorem}

\begin{proof}
In order to prove this result, consider as well the following group:
$$\mathbb TSO_N=\left\{wU\Big| w\in\mathbb T,U\in SO_N\right\}$$ 

Observe that we have $\mathbb TSO_N=\mathbb TO_N$ if $N$ is odd. If $N$ is even the group $\mathbb TO_N$ has two connected components, with $\mathbb TSO_N$ being the component containing the identity. Also, let us denote by $\mathfrak{so}_N,\mathfrak u_N$ the Lie algebras of $SO_N,U_N$.  It is well-known that $\mathfrak u_N$ consists of the matrices $M\in M_N(\mathbb C)$ satisfying $M^*=-M$, and that:
$$\mathfrak{so}_N=\mathfrak u_N\cap M_N(\mathbb R)$$

Also, it is easy to see that the Lie algebra of $\mathbb TSO_N$ is $\mathfrak{so}_N\oplus i\mathbb R$.

\medskip

\underline{Step 1}. Our first claim is that if $N\geq 2$, the adjoint representation of $SO_N$ on the space of real symmetric matrices of trace zero is irreducible. Let indeed $X \in M_N(\mathbb R)$ be symmetric with trace zero. We must prove that the following space consists of all the real symmetric matrices of trace zero:
$$V=span\left\{UXU^t\Big|U \in SO_N\right\}$$

We first prove that $V$ contains all the diagonal matrices of trace zero. Since we may diagonalize $X$ by conjugating with an element of $SO_N$, our space $V$ contains a nonzero diagonal matrix of trace zero. Consider such a matrix:
$$D=\begin{pmatrix}
d_1\\
&\ddots\\
&&d_N
\end{pmatrix}$$

We can conjugate this matrix by the following matrix:
$$\begin{pmatrix}
0&-1&0\\
1&0&0\\
0&0&I_{N-2}
\end{pmatrix}\in SO_N$$

We conclude that our space $V$ contains as well the following matrix: 
$$D'=\begin{pmatrix}
d_2\\
&d_1\\
&&d_3\\
&&&\ddots\\
&&&&d_N
\end{pmatrix}$$

More generally, we see that for any $1\leq i,j\leq N$ the diagonal matrix obtained from $D$ by interchanging $d_i$ and $d_j$ lies in $V$. Now since $S_N$ is generated by transpositions, it follows that $V$ contains any diagonal matrix obtained by permuting the entries of $D$. But it is well-known that this representation of $S_N$ on the diagonal matrices of trace zero is irreducible, and hence $V$ contains all such diagonal matrices, as claimed.

\medskip

In order to conclude now, assume that $Y$ is an arbitrary real symmetric matrix of trace zero. We can find then an element $U\in SO_N$ such that $UYU^t$ is a diagonal matrix of trace zero.  But we then have $UYU^t \in V$, and hence also $Y\in V$, as desired.

\medskip

\underline{Step 2}. Our claim is that the inclusion $\mathbb TSO_N\subset U_N$ is maximal in the category of connected compact groups. Let indeed $G$ be a connected compact group satisfying:
$$\mathbb TSO_N\subset G\subset U_N$$

Then $G$ is a Lie group. Let $\mathfrak g$ denote its Lie algebra, which satisfies:
$$\mathfrak{so}_N\oplus i\mathbb R\subset\mathfrak g\subset\mathfrak u_N$$

Let $ad_{G}$ be the action of $G$ on $\mathfrak g$ obtained by differentiating the adjoint action of $G$ on itself. This action turns $\mathfrak g$ into a $G$-module. Since $SO_N \subset G$, $\mathfrak g$ is also a $SO_N$-module. Now if $G\neq\mathbb TSO_N$, then since $G$ is connected we must have: $$\mathfrak{so}_N\oplus i\mathbb{R}\neq\mathfrak g$$

It follows from the real vector space structure of the Lie algebras $\mathfrak u_N$ and $\mathfrak{so}_N$ that there exists a nonzero symmetric real matrix of trace zero $X$ such that:
$$iX\in\mathfrak g$$

We know that the space of symmetric real matrices of trace zero is an irreducible representation of $SO_N$ under the adjoint action. Thus $\mathfrak g$ must contain all such $X$, and hence $\mathfrak g=\mathfrak u_N$.  But since $U_N$ is connected, it follows that $G=U_N$.  

\medskip

\underline{Step 3}. Let us compute now the commutant of  $SO_N$ in $ M_N(\mathbb C)$. Our first claim is that at $N=2$, this commutant is as follows: 
$$SO_2'
=\left\{\begin{pmatrix}
\alpha&\beta\\
-\beta&\alpha
\end{pmatrix}\Big|\alpha,\beta\in\mathbb C\right\}$$

As for the case $N\geq3$, our claim here is that this commutant is as follows:
$$SO_N'=\left\{\alpha I_N\Big|\alpha\in\mathbb C\right\}$$

Indeed, at $N=2$, the above formula is clear. At $N\geq 3$ now, an element in $X\in SO_N'$ commutes with any diagonal matrix having exactly $N-2$ entries equal to $1$ and two entries equal to $-1$. Hence $X$ is diagonal. Now since $X$ commutes with any even permutation matrix, and we have assumed $N\geq 3$, it commutes in particular with the permutation matrix associated with the cycle $(i,j,k)$ for any $1<i<j<k$, and hence all the entries of $X$ are the same. We conclude that $X$ is a scalar matrix, as claimed.

\medskip 

\underline{Step 4}. Our claim now is that the set of matrices with nonzero trace is dense in $SO_N$. At $N=2$ this is clear, since the set of elements in $SO_2$ having a given trace is finite.  So assume $N>2$, and consider a matrix as follows:
$$T\in SO_N\simeq SO(\mathbb R^N)\quad,\quad 
Tr(T)=0$$

Let $E\subset\mathbb R^N$ be a 2-dimensional subspace preserved by $T$, such that: 
$$T_{|E} \in SO(E)$$ 
 
Let $\varepsilon>0$ and let $S_\varepsilon \in SO(E)$ satisfying the following condition:
$$||T_{|E}-S_\varepsilon||<\varepsilon$$

Moreover, in the $N=2$ case, we can assume that $T$ satisfies as well:
$$Tr(T_{|E})\neq Tr(S_\varepsilon)$$

Now define $T_\varepsilon\in SO(\mathbb R^N)=SO_N$ by the following formulae:
$$T_{\varepsilon|E}=S_\varepsilon\quad,\quad 
T_{\varepsilon|E^\perp}=T_{|E^\perp}$$

It is clear that we have the following estimate:
$$||T-T_\varepsilon|| \leq ||T_{|E}-S_\varepsilon||<\varepsilon$$

Also, we have the following estimate, which proves our claim:
$$Tr(T_\varepsilon)=Tr(S_\varepsilon)+Tr(T_{|E^\perp})\neq0$$

\underline{Step 5}. Our claim now is that $\mathbb TO_N$ is the normalizer of $\mathbb TSO_N$ in $U_N$, i.e. is the subgroup of $U_N$ consisting of the unitaries $U$ for which, for all $X\in\mathbb TSO_N$:
$$U^{-1}XU \in\mathbb TSO_N$$

Indeed, $\mathbb TO_N$ normalizes $\mathbb TSO_N$, so we must prove that if $U\in U_N$ normalizes $\mathbb TSO_N$ then $U\in\mathbb TO_N$. First note that $U$ normalizes $SO_N$, because if $X \in SO_N$ then:
$$U^{-1}XU \in\mathbb TSO_N$$

Thus we have a formula as follows, for some $\lambda\in\mathbb T$ and $Y\in SO_N$:
$$U^{-1}XU=\lambda Y$$ 

If $Tr(X)\neq0$, we have $\lambda\in\mathbb R$ and hence:
$$\lambda Y=U^{-1}XU \in SO_N$$

The set of matrices having nonzero trace being dense in $SO_N$, we conclude that $U^{-1}XU \in SO_N$ for all $X\in SO_N$. Thus, we have:
\begin{eqnarray*}
X \in SO_N
&\implies&(UXU^{-1})^t(UXU^{-1})=I_N\\
&\implies&X^tU^tUX= U^tU\\
&\implies&U^tU \in SO_N'
\end{eqnarray*}

It follows that at $N\geq 3$ we have $U^tU=\alpha I_N$, with $\alpha \in \mathbb T$, since $U$ is unitary. Hence we have $U=\alpha^{1/2}(\alpha^{-1/2}U)$ with:
$$\alpha^{-1/2}U\in O_N\quad,\quad
U\in\mathbb TO_N$$

If $N=2$, $(U^tU)^t=U^tU$ gives again $U^tU=\alpha I_2$, and we conclude as before.

\medskip

\underline{Step 6}. Our claim is that the inclusion $\mathbb TO_N\subset U_N$ is maximal. Assume indeed that we have $\mathbb TO_N\subset G\subset U_N$, with $G\neq U_N$. From $\mathbb TSO_N\subset G_0 \subset U_N$, we obtain:
$$G_0=\mathbb TSO_N$$

But since $G_0$ is normal in $G$, the group $G$ normalizes $\mathbb TSO_N$, and hence $G\subset\mathbb TO_N$, which finishes the proof.
\end{proof}

Anlong the same lines, still following \cite{bc+}, we have as well the following result:

\index{projective unitary group}

\begin{theorem}
The following inclusion of compact groups is maximal,
$$PO_N\subset PU_N$$
in the sense that there is no intermediate compact group in between.
\end{theorem}

\begin{proof}
This follows from Theorem 7.19. Indeed, assuming $PO_N\subset G \subset PU_N$, the preimage of this subgroup under the quotient map $U_N\to PU_N$ would be then a proper intermediate subgroup of $\mathbb TO_N\subset U_N$, which is a contradiction.
\end{proof}

Still following \cite{bc+}, we have as well the following result:

\index{half-classical orthogonal group}

\begin{theorem}
The following inclusion of compact quantum groups is maximal,
$$O_N\subset O_N^*$$
in the sense that there is no intermediate compact quantum group in between.\end{theorem}

\begin{proof}
Consider indeed a sequence of surjective Hopf $*$-algebra maps as follows, whose  composition is the canonical surjection:
$$C(O_N^*)\overset{f}\longrightarrow A\overset{g}\longrightarrow C(O_N)$$

This produces a diagram of Hopf algebra maps with pre-exact rows, as follows:
$$\xymatrix@R=45pt@C=33pt{
\mathbb C\ar[r]&C(PO_N^*)\ar[d]^{f_|}\ar[r]&C(O_N^*)\ar[d]^f\ar[r]&C(\mathbb Z_2)\ar[r]\ar@{=}[d]&\mathbb C\\
\mathbb C\ar[r]&PA\ar[d]^{g_|}\ar[r]&A\ar[d]^g\ar[r]&C(\mathbb Z_2)\ar[r]\ar@{=}[d]&\mathbb C\\
\mathbb C\ar[r]&PC(O_N)\ar[r]&C(O_N)\ar[r]&C(\mathbb Z_2)\ar[r]&\mathbb C}$$
 
Consider now the following composition, with the isomorphism on the left being something well-known, coming from \cite{bdu}, as explained above: 
$$C(PU_N)\simeq C(PO_N^*)\overset{f_|}\longrightarrow PA\overset{g_|}\longrightarrow PC(O_N)\simeq C(PO_N)$$

Thus $f_|$ or $g_|$ is an isomorphism. If $f_|$ is an isomorphism we get a commutative diagram of Hopf algebra morphisms with pre-exact rows, as follows:
$$\xymatrix@R=45pt@C=33pt{
\mathbb C\ar[r]&C(PO_N^*)\ar@{=}[d]\ar[r]&C(O_N^*)\ar[d]^f\ar[r]&C(\mathbb Z_2)\ar[r]\ar@{=}[d]&\mathbb C\\
\mathbb C\ar[r]&C(PO_N^*)\ar[r]&A\ar[r]&C(\mathbb Z_2)\ar[r]&\mathbb C}$$
 
Then $f$ is an isomorphism. Similarly if $g_|$ is an isomorphism, then $g$ is an isomorphism, and this gives the result. See \cite{bc+}.
\end{proof}

There are many open questions in relation with the above results, but even with these formulated, the discussion is not over here, because we have some similar questions, which are equally interesting, for the quantum permutation groups. Let us start with:

\begin{proposition}
Consider a quantum group $S_N\subset G\subset S_N^+$, with fundamental corepresentation denoted $v$. We have then inclusions as follows, for any $k\in\mathbb N$, 
$$span\left(\xi_\pi\Big|\pi\in P(k)\right)\supset Fix(v^{\otimes k})\supset span\left(\xi_\pi\Big|\pi\in NC(k)\right)$$
and equality on the left or on the right, for any $k\in\mathbb N$, is equivalent to having equality on the left or on the right in the inclusions $S_N\subset G\subset S_N^+$.
\end{proposition}

\begin{proof}
Consider a quantum group $S_N\subset G\subset S_N^+$, and let $w,v,u$ be the fundamental corepresentations of these quantum groups. We have then inclusions as follows:
$$Fix(w^{\otimes k})\supset Fix(v^{\otimes k})\supset Fix(u^{\otimes k})$$

Moreover, by Peter-Weyl, equality on the left or on the right, for any $k\in\mathbb N$, is equivalent to having equality on the left or on the right in the inclusions $S_N\subset G\subset S_N^+$. Now by using the easiness property of $S_N,S_N^+$, this gives the result.
\end{proof}

The above result is good news, because what we have there is a purely combinatorial reformulation of the maximality conjecture, in terms of partitions, noncrossing partitions, and the associated vectors. To be more precise, we have the following statement:

\begin{theorem}
The following conditions are equivalent:
\begin{enumerate}
\item There is no intermediate quantum group, as follows:
$$S_N\subset G\subset S_N^+$$

\item Any linear combination of vectors of type
$$\xi\in span\left(\xi_\pi\Big|\pi\in P(k)\right)-span\left(\xi_\pi\Big|\pi\in NC(k)\right)$$
produces via Tannakian operations the flip map, $\Sigma(a\otimes b)=b\otimes a$.
\end{enumerate}
\end{theorem}

\begin{proof}
According to Proposition 7.22, the non-existence of the quantum groups $S_N\subset G\subset S_N^+$ is equivalent to the non-existence of Tannakian categories as follows:
$$span\left(\xi_\pi\Big|\pi\in P(k)\right)\supset C_k\supset span\left(\xi_\pi\Big|\pi\in NC(k)\right)$$

But this means that whenever we pick an element $\xi$ which is on the left, but not on the right, the Tannakian category that it generates should be the one on the left:
$$<\xi>=span\left(\xi_\pi\Big|\pi\in P(k)\right)$$

Now since the category of all partitions $P=(P(k))$ is generated by the basic crossing $\slash\hskip-2.1mm\backslash$, this amounts in saying that the Tannakian category generated by $\xi$ should contain the vector associated to this basic crossing, which is $\xi_{\slash\hskip-1.5mm\backslash}=\Sigma$, as desired.
\end{proof}

The above result might look quite encouraging, and the first thought goes into inventing some kind of tricky ``averaging operation'', perhaps probability-inspired, made up of Tannakian operations, which in practice means made of basic planar operations, which converts the crossing partitions $\pi\in P(k)-NC(k)$ into the basic crossing $\slash\hskip-2.1mm\backslash$. However, this is something difficult, and in fact such questions are almost always difficult.

\bigskip

Of course, we are not saying here that such things are hopeless, but rather that they require considerable work. In connection with the above-mentioned mysterious ``averaging operation'', our feeling is that this cannot be found with bare hands, and that a heavy use of a computer, in order to understand what is going on, is required. To our knowledge, no one has ever invested much time in all this, and so things here remain open. Getting back to Earth now, here are some concrete results, obtained in this way:

\begin{theorem}
The following happen:
\begin{enumerate}
\item There is no intermediate easy quantum group $S_N\subset G\subset S_N^+$.

\item A generalization of this fact holds, at easiness level $2$, instead of $1$.
\end{enumerate}
\end{theorem}

\begin{proof}
The idea here is that everything follows from Theorem 7.23, with suitable definitions for the various easiness notions involved, and by doing some combinatorics:

\medskip

(1) Here what happens is that any $\pi\in P-NC$ has the following property:
$$<\pi>=P$$

Indeed, the idea is to cap $\pi$ with semicircles, as to preserve one crossing, chosen in advance, and to end up, by a recurrence procedure, with the standard crossing:
$$\slash\hskip-2.1mm\backslash\in<\pi>$$

Now in terms of the notions in Theorem 7.23, the conclusion is that the criterion (2) there holds for the linear combinations $\xi$ having lenght 1, and this gives the result. Indeed, according to \cite{bsp}, the easy quantum groups are by definition those having Tannakian categories as follows, with $D=(D(k))$ being a certain category of partitions:
$$Fix(v^{\otimes k})=span\left(\xi_\pi\Big|\pi\in D(k)\right)$$

Thus, the generation formula $<\pi>=P$ established above does the job, and proves that an intermediate easy quantum group $S_N\subset G\subset S_N^+$ cannot exist. See \cite{ez1}.

\medskip

(2) This is a generalization of (1), the idea being that of looking at the combinations having length 2, of type $\xi=\alpha\xi_\pi+\beta\xi_\sigma$. Our first claim is that, assuming that $G\subset H$ comes from an inclusion of categories $D\subset E$, the maximality at order $2$ is equivalent to the following condition, for any $\pi,\sigma\in E$, not both in $D$, and for any $\alpha,\beta\neq0$:
$$<span(D),\alpha T_\pi+\beta T_\sigma>=span(E)$$

Consider indeed a category $span(D)\subset C\subset span(E)$, corresponding to a quantum group $G\subset K\subset H$ having order 2. The order 2 condition means that we have $C=<C\cap span_2(P)>$, where $span_2$ denotes the space of linear combinations having 2 components. Since we have $span(E)\cap span_2(P)=span_2(E)$, the order 2 formula reads:
$$C=<C\cap span_2(E)>$$

Now observe that the category on the right is generated by the categories $C_{\pi\sigma}^{\alpha\beta}$ constructed in the statement. Thus, the order 2 condition reads:
$$C=\left<C_{\pi\sigma}^{\alpha\beta}\Big|\pi,\sigma\in E,\alpha,\beta\in\mathbb C\right>$$

Now since the maximality at order 2 of the inclusion $G\subset H$ means that we have $C\in\{span(D),span(E)\}$, for any such $C$, we are led to the following condition:
$$C_{\pi\sigma}^{\alpha\beta}\in\{span(D),span(E)\}\quad,\quad\forall\pi,\sigma\in E,\alpha,\beta\in\mathbb C$$

Thus, we have proved our claim. In order to prove now that $S_N\subset S_N^+$ is maximal at order $2$, we can use semicircle capping. The statement that we have to prove is as follows: ``for $\pi\in P-NC,\sigma\in P$ and $\alpha,\beta\neq0$ we have $<\alpha T_\pi+\beta T_\sigma>=span(P)$''. 

\medskip

In order to do this, our claim is that the same method as at level 1 applies, after some suitable modifications. We have indeed two cases, as follows:

\medskip

-- Assuming that $\pi,\sigma$ have at least one different crossing, we can cap the partition $\pi$ as to end up with the basic crossing, and $\sigma$ becomes in this way an element of $P(2,2)$ different from this basic crossing, and so a noncrossing partition, from $NC(2,2)$. Now by substracting this noncrossing partition, which belongs to $C_{S_N^+}=span(NC)$, we obtain that the standard crossing belongs to $<\alpha T_\pi+\beta T_\sigma>$, and we are done.

\medskip

-- In the case where $\pi,\sigma$ have exactly the same crossings, we can start our descent procedure by selecting one common crossing, and then two strings of $\pi,\sigma$ which are different, and then joining the crossing to these two strings. We obtain in this way a certain linear combination $\alpha' T_{\pi'}+\beta'T_{\sigma'}\in <\alpha T_\pi+\beta T_\sigma>$ which satisfies the conditions in the first case discussed above, and we can continuate as indicated there.
\end{proof}

\section*{7e. Exercises} 

We have been once again into a research-flavored chapter, and there are no easy exercises about all this. As a good research exercise, however, we have:

\begin{exercise}
Prove that $S_N\subset S_N^+$ is maximal, at $N=1,2,3,4,5,6,7$.
\end{exercise}

To be more precise, the maximality of $S_N\subset S_N^+$ is trivial at $N=1$, but let us award 1 Dan for that, then you have to think a bit for $S_2\subset S_2^+$, worth 2 Dan, then $S_3\subset S_3^+$ was done by Wang long ago, crucial pioneering work, worth 3 Dan, then for $S_4\subset S_4^+$ I personally asked Bichon and we got awarded 4 Dan for our joint paper, and then for $S_5\subset S_5^+$ that comes from the classification of subfactors of index 5, whose full reading, with all the von Neumann algebra and subfactor needed preliminaries, is worth 5 Dan. Regarding now 6 Dan, this looks certainly possible, again by using the known subfactor results, but this time with the work needed including contributing a bit to that subfactor classification at $N=6$, under suitable transitivity assumptions, and with all this being, potentially, quality research work. As for 7 Dan, this is normally reserved to cats, tigers and other felines, but fearless young humans like you are of course welcome for a try.

\chapter{Unitary groups}

\section*{8a. Basic examples}

In this chapter we discuss a number of more specialized questions, which are however of crucial importance, for the general understanding of the easy quantum groups. Let us go back, as usual, to the standard cube formed by the main quantum groups:
$$
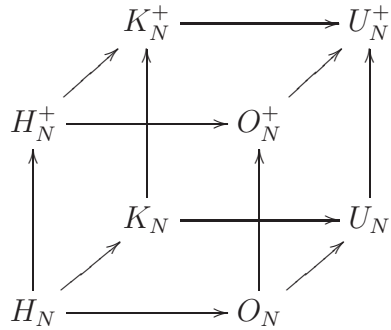
\xymatrix@R=18pt@C=18pt{
&K_N^+\ar[rr]&&U_N^+\\
H_N^+\ar[rr]\ar[ur]&&O_N^+\ar[ur]\\
&K_N\ar[rr]\ar[uu]&&U_N\ar[uu]\\
H_N\ar[uu]\ar[ur]\ar[rr]&&O_N\ar[uu]\ar[ur]
}$$

We will be mainly interested, as before in this Part II of the present book, in the continuous case, corresponding to the right face of the cube, namely:
$$
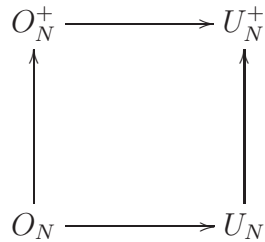
\xymatrix@R=58pt@C=56pt{
O_N^+\ar[r]&U_N^+\\
O_N\ar[u]\ar[r]&U_N\ar[u]}$$

We have seen in chapter 7 that, in what regards the left edge, there is only one intermediate object there, namely the half-classical orthogonal group, $O_N\subset O_N^*\subset O_N^+$. Moreover, this intermediate quantum group is conjecturally unique, even in the general, non-easy setting. Adding to the picture, this quantum group has a unitary counterpart, $U_N\subset U_N^*\subset U_N^+$, which makes our square split in a nice way, as follows:
$$\xymatrix@R=28pt@C=84pt{
O_N^+\ar[r]&U_N^+\\
O_N^*\ar[r]\ar[u]&U_N^*\ar[u]\\
O_N\ar[u]\ar[r]&U_N\ar[u]}$$

In fact, the whole cube can be split in this way, thanks to half-classical reflection groups $H_N^*,K_N^*$, which exist and are non-trivial, obtained by intersecting, and fit well into the diagram, cutting the left face of the cube in the following way:
$$\xymatrix@R=28pt@C=84pt{
H_N^+\ar[r]&K_N^+\\
H_N^*\ar[r]\ar[u]&K_N^*\ar[u]\\
H_N\ar[u]\ar[r]&K_N\ar[u]}$$

In addition, we have also seen in chapter 7 that it is possible to further split the cube, thanks to certain canonical ``hybrid'' objects, lying between real and complex. In the continuous case, the enlarged right square, containing these hybrids, is as follows:
$$\xymatrix@R=14mm@C=14mm{
O_N^+\ar[r]&\mathbb TO_N^+\ar[r]&U_N^+\\
O_N^*\ar[r]\ar[u]&\mathbb TO_N^*\ar[r]\ar[u]&U_N^*\ar[u]\\
O_N\ar[r]\ar[u]&\mathbb TO_N\ar[r]\ar[u]&U_N\ar[u]}$$

As for the enlarged left square, this looks similar, as follows:
$$\xymatrix@R=14mm@C=14mm{
H_N^+\ar[r]&\mathbb TH_N^+\ar[r]&K_N^+\\
H_N^*\ar[r]\ar[u]&\mathbb TH_N^*\ar[r]\ar[u]&K_N^*\ar[u]\\
H_N\ar[r]\ar[u]&\mathbb TH_N\ar[r]\ar[u]&K_N\ar[u]}$$

All this is very nice. However, and here comes our point, things are far from being over, because if we go back to the continuous square, in what regards the right edge of this square, namely $U_N\subset U_N^*\subset U_N^+$, and in contrast with what happens for the left edge, there are in fact many examples of intermediate easy quantum groups, as follows:
$$U_N\subset G\subset U_N^+$$

As a basic example here, whose construction is elementary, going without thinking, consider for instance the free complexification $U_N^\times$ of the unitary group $U_N$. This quantum group appears then as an intermediate quantum group, as above:
$$U_N\subset U_N^\times\subset U_N^+$$

It is quite clear that we have $U_N^\times\neq U_N$, due to noncommutativity, for instance by looking at the diagonal torus. Also, we have $U_N^\times\neq U_N^+$, for instance by looking at the fusion rules. And finally, we have $U_N^\times\neq U_N^*$ as well, once again by looking for instance at the fusion rules. Thus, we have here our new intermediate quantum group $U_N\subset G\subset U_N^+$, and there are probably many more quantum groups, of this type.

\bigskip

Our goal here will be that of understanding these examples, first with the construction and study of some basic classes of such examples, following my early work with Bichon \cite{bb1}, \cite{bb2}, then with the systematic study of the problem, including a full classification result, which is something quite technical, due to Mang and Weber \cite{mw1}, \cite{mw2}, and finally with various algebraic and analytic results about all these new quantum groups. 

\bigskip

Getting started now, let us first discuss some basic examples of unitary quantum groups, following \cite{bb1}, \cite{bb2} and related papers. We already have 3 examples of such quantum groups, namely $U_N\subset U_N^*\subset U_N^+$, that we know well. Our goal will be that of constructing 2 more fundamental examples, appearing as intermediate objects, as follows:
$$U_N\subset U_N^*\subset U_N^{**}\subset U_N^\times\subset U_N^+$$

With respect to the previous square diagram, these quantum groups will fit as follows, making things starting to be a bit crowded, on the right edge:
$$
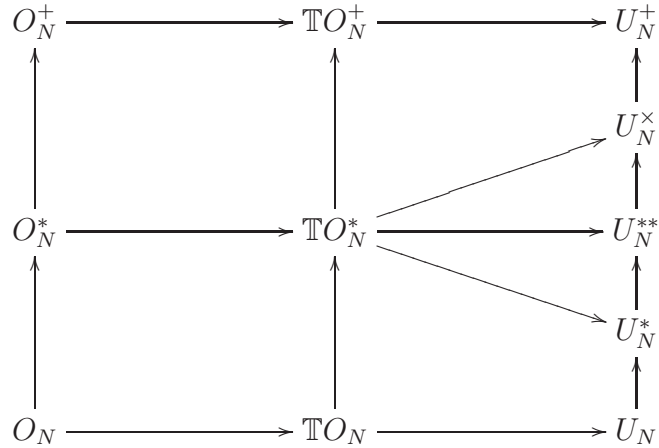
\xymatrix@R=7mm@C=30mm{
O_N^+\ar[r]&\mathbb TO_N^+\ar[r]&U_N^+\\
&&U_N^\times\ar[u]\\
O_N^*\ar[r]\ar[uu]&\mathbb TO_N^*\ar[r]\ar[uu]\ar[ur]\ar[dr]&U_N^{**}\ar[u]\\
&&U_N^*\ar[u]\\
O_N\ar[r]\ar[uu]&\mathbb TO_N\ar[r]\ar[uu]&U_N\ar[u]}$$

Observe in particular that the previous ``canonical'' example of an intermediate quantum group $U_N\subset G\subset U_N^+$, which was the quantum group $U_N^*$, will get in this way downgraded to a secondary status, with the king becoming a certain quantum group $U_N^{**}$, to be introduced in a moment. However, the story will be not over here, because when getting at an even more advanced level, that of the papers \cite{mw1}, \cite{mw2}, the king will change again. In short, be prepared for some exciting story, in what follows next.

\bigskip

In order to get started, a first idea is that of weakening the relations defining $U_N^*$. We know that this quantum group is easy, the result regarding it being as follows:

\begin{proposition}
The quantum group $U_N^*$ is easy, coming from the diagrams
$$
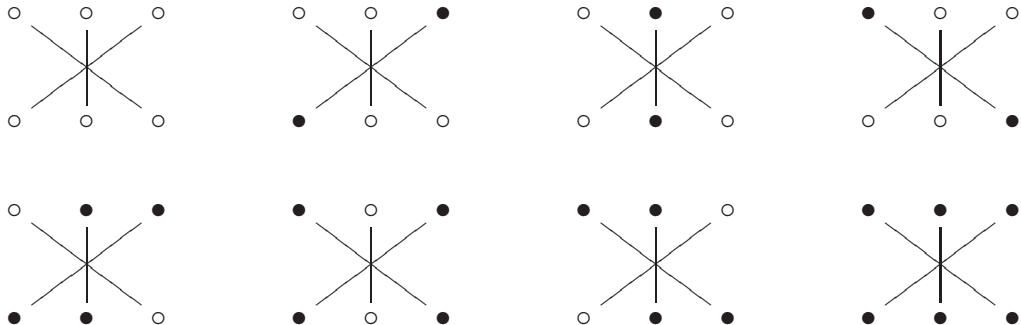
\xymatrix@R=10mm@C=5mm{
\circ\ar@{-}[drr]&\circ\ar@{-}[d]&\circ\ar@{-}[dll]
\\
\circ&\circ&\circ}\quad\ \qquad
\xymatrix@R=10mm@C=5mm{
\circ\ar@{-}[drr]&\circ\ar@{-}[d]&\bullet\ar@{-}[dll]
\\
\bullet&\circ&\circ}\quad\ \qquad
\xymatrix@R=10mm@C=5mm{
\circ\ar@{-}[drr]&\bullet\ar@{-}[d]&\circ\ar@{-}[dll]
\\
\circ&\bullet&\circ}\quad\ \qquad
\xymatrix@R=10mm@C=5mm{
\bullet\ar@{-}[drr]&\circ\ar@{-}[d]&\circ\ar@{-}[dll]
\\
\circ&\circ&\bullet}$$

$$\xymatrix@R=10mm@C=5mm{
\circ\ar@{-}[drr]&\bullet\ar@{-}[d]&\bullet\ar@{-}[dll]
\\
\bullet&\bullet&\circ}\quad\ \qquad
\xymatrix@R=10mm@C=5mm{
\bullet\ar@{-}[drr]&\circ\ar@{-}[d]&\bullet\ar@{-}[dll]
\\
\bullet&\circ&\bullet}\quad\ \qquad
\xymatrix@R=10mm@C=5mm{
\bullet\ar@{-}[drr]&\bullet\ar@{-}[d]&\circ\ar@{-}[dll]
\\
\circ&\bullet&\bullet}\quad\ \qquad
\xymatrix@R=10mm@C=5mm{
\bullet\ar@{-}[drr]&\bullet\ar@{-}[d]&\bullet\ar@{-}[dll]
\\
\bullet&\bullet&\bullet}$$
which must intertwine the various $3$-fold tensor products between $u$ and $\bar{u}$.
\end{proposition}

\begin{proof}
This is something that we know well from the previous chapter, with the diagrams in the statement producing the relations $abc=cba$, $abc^*=c^*ba$, $ab^*c=cb^*a$, and so on up to $a^*b^*c^*=c^*b^*a^*$, between the standard coordinates $u_{ij}$.
\end{proof} 

The above result suggests that there might be $2^8=256$ easy quantum groups $U_N^*\subset G\subset U_N^+$ that can be constructed by using length 3 relations as above, simply by selecting a subset of our set of 8 diagrams, and then constructing $G\subset U_N^+$ by using these diagrams. However, this is far from being true, because the relations produced by our 8 diagrams are often equivalent. To be more precise, we have a number of obvious equivalences, and by erasing the corresponding diagrams, we are led to three diagrams, namely:
$$
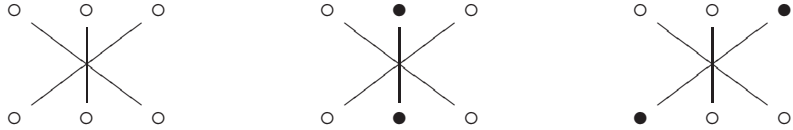
\xymatrix@R=10mm@C=5mm{
\circ\ar@{-}[drr]&\circ\ar@{-}[d]&\circ\ar@{-}[dll]
\\
\circ&\circ&\circ}\qquad\ \qquad
\xymatrix@R=10mm@C=5mm{
\circ\ar@{-}[drr]&\bullet\ar@{-}[d]&\circ\ar@{-}[dll]
\\
\circ&\bullet&\circ}\qquad\ \qquad
\xymatrix@R=10mm@C=5mm{
\circ\ar@{-}[drr]&\circ\ar@{-}[d]&\bullet\ar@{-}[dll]
\\
\bullet&\circ&\circ}$$

Thus, we have in fact at most $2^3=8$ easy quantum groups $U_N^*\subset G\subset U_N^+$ that can be constructed, by picking any of these diagrams, or some combination of these diagrams, and using the corresponding relations. But here, we have the following result:

\begin{proposition}
Consider the following types of relations, between abstract variables $a,b,c\in\{x_i\}$ subject to the relations $\sum_ix_ix_i^*=\sum_ix_i^*x_i=1$:
\begin{enumerate}
\item $abc=cba$.

\item $ab^*c=cb^*a$

\item $abc^*=c^*ba$.
\end{enumerate}
We have then $(1)\iff(3)\implies(2)$.
\end{proposition}

\begin{proof}
The equivalence $(1)\iff(3)$ follows from the following computations, with the semicircle cappings corresponding to the use of $\sum_ix_ix_i^*=\sum_ix_i^*x_i=1$:
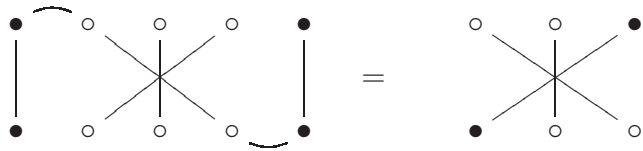
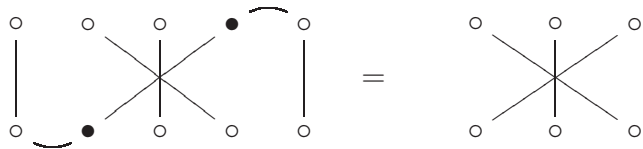
$$\xymatrix@R=10mm@C=5mm{
\bullet\ar@{-}[d]\ar@/^/@{-}[r]&\circ\ar@{-}[drr]&\circ\ar@{-}[d]&\circ\ar@{-}[dll]&\bullet\ar@{-}[d]\\
\bullet&\circ&\circ&\circ\ar@/_/@{-}[r]&\bullet}\quad
\xymatrix@R=4mm@C=6mm{&\\ =\\&\\& }
\xymatrix@R=10mm@C=6mm{
\circ\ar@{-}[drr]&\circ\ar@{-}[d]&\bullet\ar@{-}[dll]\\
\bullet&\circ&\circ}$$

$$\xymatrix@R=10mm@C=5mm{
\circ\ar@{-}[d]&\circ\ar@{-}[drr]&\circ\ar@{-}[d]&\bullet\ar@{-}[dll]\ar@/^/@{-}[r]&\circ\ar@{-}[d]\\
\circ\ar@/_/@{-}[r]&\bullet&\circ&\circ&\circ}\quad
\xymatrix@R=4mm@C=6mm{&\\ =\\&\\& }
\xymatrix@R=10mm@C=6mm{
\circ\ar@{-}[drr]&\circ\ar@{-}[d]&\circ\ar@{-}[dll]\\
\circ&\circ&\circ}$$

As for $(1+3)\implies(2)$, this is best worked out at the algebraic level, again by using the relations $\sum_ix_ix_i^*=\sum_ix_i^*x_i=1$, when summing over $d$, as follows:
\begin{eqnarray*}
ab^*c
&=&\sum_dab^*cdd^*\\
&=&\sum_dadcb^*d^*\\
&=&\sum_dcdab^*d^*\\
&=&\sum_dcb^*add^*\\
&=&cb^*a
\end{eqnarray*}

Thus we have indeed $(1)\iff(3)\implies(2)$, as claimed.
\end{proof}

Now by getting back to our problem, and more specifically, to our above-mentioned idea of using one diagram out of 3 possible ones, we can see, as a consequence of Proposition 8.2, that we have only one good choice, and are led to the following definition:

\begin{definition}
We have an intermediate quantum group as follows,
$$U_N^*\subset U_N^\times\subset U_N^+$$
obtained via the relations $ab^*c=cb^*a$, with $a,b,c\in\{u_{ij}\}$.
\end{definition}

So, this will be our first example of intermediate quantum group $U_N^*\subset G\subset U_N^+$, appearing from some straightforward combinatorial considerations, as explained above. Let us mention right away that $U_N^\times$ appears, as previously indicated, as well as the free complexification of the unitary group $U_N$. But more on this later.

\bigskip

The quantum group $U_N^\times$, in its either incarnations, appearing via Definition 8.3, or as free complexification, has been known for some time, and its theory is well understood. As a first result regarding it, in relation with a diagram drawn before, we have:

\begin{proposition}
We have the formula
$$\mathbb TO_N^+\cap U_N^\times=\mathbb TO_N^*$$
as an equality of quantum subgroups of $U_N^+$.
\end{proposition}

\begin{proof}
According to the definition of $\mathbb TO_N^*$, this quantum group appears as $\mathbb TO_N^*=\mathbb TO_N^+\cap U_N^*$. Thus, we must  prove that we have:
$$\mathbb TO_N^+\cap U_N^\times\subset U_N^*$$

In terms of defining relations, we must prove that, from $ab^*=a^*b$ and $ab^*c=cb^*a$ for any $a,b,c\in\{u_{ij}\}$, we can deduce that we have, for any $a,b,c\in\{u_{ij},u_{ij}^*\}$:
$$abc=cba$$

But this is clear, because by using $ab^*=a^*b$, we can first obtain $a^*bc=cba^*$, and then, by using Proposition 8.2, we can obtain from this the other relations as well. Here we have used the fact that what we know about abstract variables satisfying $\sum_ix_ix_i^*=\sum_ix_i^*x_i=1$ applies to the coordinates to any closed subgroup $G\subset U_N^+$, simply because these coordinates, when rescaled by $\sqrt{N}$, do satisfy these relations.
\end{proof}

As another interesting result about $U_N^\times$, we have:

\begin{proposition}
We have the formula
$$PU_N^\times=PU_N$$
as an equality of quantum subgroups of $PU_N^+$.
\end{proposition}

\begin{proof}
By using the commutation relations $ab^*c=cb^*a$, we obtain:
$$ab^*cd^*=cb^*ad^*=cd^*ab^*$$

Thus the projective coordinates $ab^*$ commute, and this gives the result.
\end{proof}

As already mentioned, the quantum group $U_N^\times$ has been known for some time, and its theory is well understood. As a summary of what can be said about it, we have:

\begin{theorem}
The following happen, regarding the quantum group $U_N^\times$:
\begin{enumerate}
\item The associated category of pairings $D$ can be explicitly described.

\item We have an inclusion of quantum groups $\mathbb TO_N^*\subset U_N^\times$.

\item We have the projective version result $PU_N^\times=PU_N$.

\item The character laws and other probabilistic aspects can be worked out.
\end{enumerate}
\end{theorem}

\begin{proof}
All this is routine, and can be found in the literature, along with the precise statements, which can be quite technical, the idea being as follows:

\medskip

(1) This is something standard, by doing some combinatorics, and we will leave this as an instructive exercise, at this stage of things. We will be back to this question at the end of the present chapter, armed with a more conceptual interpretation of $U_N^\times$, and of some related quantum groups, coming from the recent work of Mang-Weber \cite{mw1}, \cite{mw2}.

\medskip

(2) This is something that we know from Proposition 8.4.

\medskip

(3) This is something that we know too, from Proposition 8.5.

\medskip

(4) This is something standard too, which is in relation with (1), and again we will leave this as an instructive exercise, at this stage of things, and with the promise that we will be back to this, later in this chapter, armed with more powerful technology.
\end{proof}

\section*{8b. Further examples}

Getting back now to our general question of constructing intermediate easy quantum groups $U_N^*\subset G\subset U_N^+$, the approach based on the diagrams from Proposition 8.1 cannot lead to more examples, and we must come up with something more conceptual. The most elegant approach to all this, and in fact to the quantum group $U_N^\times$ itself too, is via projective versions. Let us start with the following standard construction:

\begin{proposition}
Let $G$ be a compact quantum group, and $v=(v_{ij})$ be a corepresentation of $C(G)$. We have then a quotient quantum group $G\to H$, given by:
$$C(H)=<v_{ij}>$$
At the dual level we obtain a discrete quantum subgroup, $\widehat{\Lambda}\subset\widehat{\Gamma}$.
\end{proposition}

\begin{proof}
Here the first assertion follows from the definition of the corepresentations, as being the square matrices satisfying the following conditions:
$$\Delta(v_{ij})=\sum_kv_{ik}\otimes v_{kj}\quad,\quad
\varepsilon(v_{ij})=\delta_{ij}\quad,\quad
S(v_{ij})=v_{ji}^*$$

As for the second assertion, this is just a reformulation of the first assertion, coming from the basic functoriality properties of the Pontrjagin duality.
\end{proof}

We can now talk about projective versions, as follows:

\index{projective version}

\begin{proposition}
Given a compact quantum group $G$, with fundamental corepresentation $u=(u_{ij})$, the $N^2\times N^2$ matrix given in double index notation by
$$v_{ia,jb}=u_{ij}u_{ab}^*$$
is a corepresentation in the above sense, and we have the following results:
\begin{enumerate}
\item The corresponding quotient $G\to PG$ is a compact quantum group.

\item In the classical group case, $G\subset U_N$, we have $PG=G/(G\cap\mathbb T^N)$.

\item In the group dual case, with $\Gamma=<g_i>$, we have $\widehat{P\Gamma}=<g_ig_j^{-1}>$.
\end{enumerate}
\end{proposition}

\begin{proof}
The fact that $v$ is indeed a corepresentation is routine, coming from the definition of the corepresentations. Regarding now other assertions, all these are standard for the classical groups, and in general the proofs are similar, as follows:

\medskip

(1) This follows from Proposition 8.7.

\medskip

(2) This follows from the elementary fact that, via Gelfand duality, $w$ is the matrix of coefficients of the adjoint representation of $G$, whose kernel is the subgroup $G\cap\mathbb T^N$, where $\mathbb T^N\subset U_N$ denotes the subgroup formed by the diagonal matrices.

\medskip

(3) This is something trivial, which follows from definitions.
\end{proof}

As a first interesting result now about projective versions, which is something intimately related to freeness, having no classical counterpart, we have:

\begin{theorem}
We have an identification as follows,
$$PO_N^+=PU_N^+$$
modulo the usual equivalence relation for compact quantum groups.
\end{theorem}

\begin{proof}
We have several proofs for this result, as follows:

\medskip

(1) This follows from the free complexification result $\widetilde{O_N^+}=U_N^+$, that we already met in this book, because by using this, we have right away:
$$PU_N^+=P\widetilde{O_N^+}=PO_N^+$$

(2) We can deduce this as well directly. Indeed, if we denote by $v,u$ the fundamental corepresentations of $O_N^+,U_N^+$, then by easiness we have equalities as follows:
$$Hom\left((v\otimes v)^k,(v\otimes v)^l\right)=span\left(T_\pi\Big|\pi\in NC_2((\circ\bullet)^k,(\circ\bullet)^l)\right)$$
$$Hom\left((u\otimes\bar{u})^k,(u\otimes\bar{u})^l\right)=span\left(T_\pi\Big|\pi\in \mathcal{NC}_2((\circ\bullet)^k,(\circ\bullet)^l)\right)$$

The sets on the right being equal, we conclude that the inclusion $PO_N^+\subset PU_N^+$ preserves the corresponding Tannakian categories, and so must be an isomorphism.
\end{proof}

The above result is really exciting, and suggests:

\begin{thought}
Free geometry, at least in its projective version, might be simpler than classical geometry.
\end{thought}

Which might perhaps sound a bit odd, if this is your first time reading about freeness, and struggling with the details. If this is the case, sorry and I have to admit that the above thought is not mine, but rather comes from my cat. In any case, to be recorded, and we will be back to this at the end of the present chapter, and then more in detail in chapter 16, when discussing noncommutative geometry in general.

\bigskip

Back to work now, in order to further talk about the liberations of $U_N$, we will need in fact  the following more specialized notions:

\index{left projective version}
\index{right projective version}
\index{full projective version}

\begin{definition}
Given a closed subgroup $G\subset U_N^+$, we define quotients as follows, with the variables $x_i$ standing for the standard coordinates $u_{ij}$,
\begin{enumerate}
\item Left projective version: $G\to PG$, with coordinates $p_{ij}=x_ix_j^*$,

\item Right projective version: $G\to P'G$, with coordinates $q_{ij}=x_j^*x_i$,

\item Full projective version: $G\to\mathcal PG$, with coordinates $p_{ij},q_{ij}$,
\end{enumerate}
and we say that $G$ is left/right/full half-classical when these spaces are classical.
\end{definition}

As a first observation, the left projective version $G\to PG$ is the one that we have been using so far in this book, since chapter 2. Thus, what we are doing here is that of fine-tuning our projective version formalism, with 3 notions instead of one.

\bigskip

Observe that in the classical case, where $G\subset U_N$, the three projective versions constructed above obviously coincide, and equal the usual projective version, obtained by dividing under the action of $\mathbb T$. Also, in the real case, $G\subset O_N^+$, the three projective versions coincide as well, and $G$ is left or right half-classical when $G\subset O_N^*$. 

\bigskip

However, in the general complex case $G\subset U_N^+$, the three projective versions from Definition 8.11 do not necessarily coincide, and this can lead to some interesting theory. In order to discuss this, let us start with an elementary result, as follows:

\begin{proposition}
Let $G\subset U_N^+$, with coordinates $x_1,\ldots,x_N$.
\begin{enumerate}
\item $G\subset U_N^\times$ precisely when $\{x_ix_j^*\}$ commute, and $\{x_i^*x_j\}$ commute as well.

\item $G\subset U_N^*$ precisely when the variables $\{x_ix_j,x_ix_j^*,x_i^*x_j,x_i^*x_j^*\}$ all commute.
\end{enumerate}
\end{proposition}

\begin{proof}
Regarding the first assertion, the implication ``$\implies$'' follows from:
$$ab^*cd^*=cb^*ad^*=cd^*ab^*$$
$$a^*bc^*d=c^*ba^*d=c^*da^*b$$

As for the implication ``$\Longleftarrow$'', we can use here the following computation, based on the commutation assumptions in the statement:
\begin{eqnarray*}
ae^*eb^*c
&=&ab^*ce^*e\\
&=&ce^*ab^*e\\
&=&cb^*ee^*a
\end{eqnarray*}

Indeed, by summing over $e=x_i$, we obtain from this, as desired:
$$ab^*c=cb^*a$$

The proof of the second assertion is similar, because we can remove all $*$ signs, except for those concerning $e^*$, and use the above computations with $a,b,c,d\in\{x_i,x_i^*\}$.
\end{proof}

With the above result in hand, we can now formulate:

\begin{theorem}
We have the following results:
\begin{enumerate}
\item $U_N^\times$ is left and right half-classical, and is maximal with this property.

\item $U_N^*$ is fully half-classical, and is maximal with this property.

\item $O_N^*$ is fully half-classical, and is maximal inside $O_N^+$ with this property.
\end{enumerate}
\end{theorem}

\begin{proof}
All these assertions follow indeed from Proposition 8.12.
\end{proof}

We recall that the free complexification of a compact quantum group $G$, with standard coordinates denoted $v_{ij}$, is the compact quantum group $\widetilde{G}$ corresponding to the subalgebra $C(\widetilde{G})\subset C(\mathbb T)*C(G)$ generated by the variables $u_{ij}=zv_{ij}$, where $z$ is the standard generator of $C(\mathbb T)$. Observe that $\widetilde{G}$ is indeed a quantum group, because it appears as a subgroup of $\mathbb T\,\hat{*}\,G$, the quantum group associated to $C(\mathbb T)*C(G)$. We have:

\begin{theorem}
The quantum groups $O_N^*,U_N,U_N^*,U_N^\times$ are as follows:
\begin{enumerate}
\item They have the same left projective version, equal to $PU_N$.

\item They have the same free complexification, equal to $U_N^\times$.
\end{enumerate} 
\end{theorem}

\begin{proof}
This is standard, the idea being as follows:

\medskip

(1) Here $PO_N^*=PU_N$ is a well-known result, that we already know, from chapter 7. It is clear as well that we have inclusions as follows:
$$PU_N\subset PU_N^*\subset PU_N^\times$$

Now by using Proposition 8.12 (1), we conclude that the quantum group $PU_N^\times$ is classical, and so we must have an inclusion as follows:
$$PU_N^\times\subset(PU_N^+)_{class}$$

But this latter quantum group $(PU_N^+)_{class}$ is known to be equal to $PU_N$, and so this inclusion is an isomorphism, and this finishes the proof.

\medskip

(2) If we denote by $v_{ij}$ the standard coordinates on $G=O_N^*,U_N,U_N^*,U_N^\times$, and by $z$ the generator of a copy of $C(\mathbb T)$, free from $C(G)$, then with $a,b,c\in\{v_{ij}\}$ we have:
\begin{eqnarray*}
(za)(zb)^*(zc)
&=&zab^*c\\
&=&zcb^*a\\
&=&(zc)(zb)^*(za)
\end{eqnarray*}

Thus we have $\widetilde{G}\subset\widetilde{U}_N$. Conversely now, it follows from the general theory of the free complexifications of easy quantum groups \cite{rau} that both $K=\widetilde{G},\widetilde{U}_N$ should appear as free complexifications of certain intermediate easy quantum groups, as follows:
$$O_N\subset H\subset O_N^+$$

On the other hand, since we have $PH=P\widetilde{H}=PK=PU_N$, the only choice here is $H=O_N^*$. Thus we have $\widetilde{G}=\widetilde{U}_N=\widetilde{O_N^*}$, and this finishes the proof.
\end{proof}

Going ahead now, our various results above, in relation with projective versions and complexifications, suggest introducing a new quantum group $U_N\subset G\subset U_N^+$, as follows:

\begin{definition}
We have a quantum group as follows,
$$U_N^{**}\subset U_N^+$$
obtained via the relations ``$ab^*,a^*b$ all commute'', with $a,b,c\in\{u_{ij}\}$.
\end{definition}

As a first remark, the real version of $U_N^{**}$, obtained by imposing the conditions $x_i=x_i^*$ to the standard coordinates, is the half-classical orthogonal quantum group $O_N^*$. Also, we have inclusions as follows, coming from the various results in Proposition 8.12:
$$U_N^*\subset U_N^{**}\subset U_N^\times$$

Summarizing, at the level of the examples of new quantum groups that we have, we are led to the diagram announced in the beginning of this chapter, namely:
$$
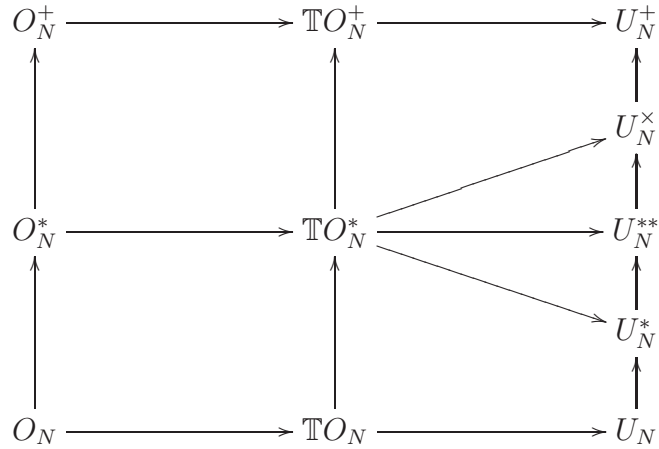
\xymatrix@R=7mm@C=30mm{
O_N^+\ar[r]&\mathbb TO_N^+\ar[r]&U_N^+\\
&&U_N^\times\ar[u]\\
O_N^*\ar[r]\ar[uu]&\mathbb TO_N^*\ar[r]\ar[uu]\ar[ur]\ar[dr]&U_N^{**}\ar[u]\\
&&U_N^*\ar[u]\\
O_N\ar[r]\ar[uu]&\mathbb TO_N\ar[r]\ar[uu]&U_N\ar[u]}$$

To be more precise, the fact that we have this diagram, and also that suitable subdiagrams of it are intersection and generation diagrams, follows from what we have.

\bigskip

Let us develop now some general theory for $U_N^{**}$, according to our usual investigation pattern, for new quantum groups. First, we have the following result:

\begin{proposition}
The quantum group $U_N^{**}$ is easy, coming from the diagrams
$$
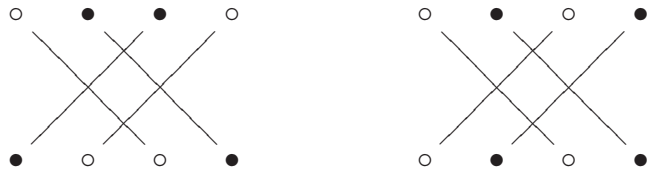
\xymatrix@R=15mm@C=5mm{\circ\ar@{-}[drr]&\bullet\ar@{-}[drr]&\bullet\ar@{-}[dll]&\circ\ar@{-}[dll]\\\bullet&\circ&\circ&\bullet}\qquad \quad\qquad 
\xymatrix@R=15mm@C=5mm{\circ\ar@{-}[drr]&\bullet\ar@{-}[drr]&\circ\ar@{-}[dll]&\bullet\ar@{-}[dll]\\\circ&\bullet&\circ&\bullet}$$
producing intertwiners between the corresponding $4$-fold tensor products between $u,\bar{u}$.
\end{proposition}

\begin{proof}
We know that $U_N^{**}\subset U_N^+$ appears by imposing the conditions stating that the variables $ab^*,a^*b$ all commute, with $a,b,c\in\{u_{ij}\}$. But each such commutation relation corresponds to a certain intertwining relation between certain $4$-fold tensor products between $u,\bar{u}$, and this leads to the conclusion in the statement.
\end{proof}

Next, we have the following key result, further building on the above:

\begin{theorem}
The quantum group $U_N^{**}$ is easy, coming from the category $\mathcal P_2^{**}$ of matching pairings having the property that
$$\#\circ=\#\bullet$$
between the legs of each string, when flattened.
\end{theorem}

\begin{proof}
We know from Proposition 8.16 that the quantum group $U_N^{**}$ is easy. Consider now the diagrams producing this quantum group, namely:
$$
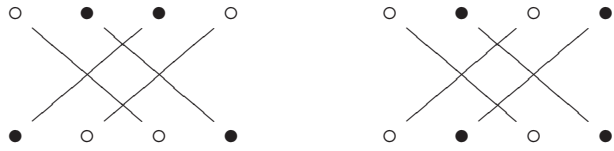
\xymatrix@R=12mm@C=5mm{\circ\ar@{-}[drr]&\bullet\ar@{-}[drr]&\bullet\ar@{-}[dll]&\circ\ar@{-}[dll]\\\bullet&\circ&\circ&\bullet}\qquad \quad \quad 
\xymatrix@R=12mm@C=5mm{\circ\ar@{-}[drr]&\bullet\ar@{-}[drr]&\circ\ar@{-}[dll]&\bullet\ar@{-}[dll]\\\circ&\bullet&\circ&\bullet}$$

By rotating these diagrams, we can see that the condition in the statement, namely $\#\circ=\#\bullet$ between the legs of each string, is satisfied for both. But this leads, via some standard combinatorics, to the conclusion in the statement, and we refer here to \cite{mw1}.
\end{proof}

There are many other things that can be said here, both at the algebraic and the probabilistic level, and for more on all this, we refer to the literature.

\section*{8c. The standard series}

Still following \cite{bb1}, \cite{mw1}, let us further extend now the constructions that we have, of quantum unitary groups. The idea will be that of interpolating between $U_N$ and $U_N^{**}$, with quantum groups $U_N\subset U_N^{(r)}\subset U_N^{**}$, as to have equalities as follows:
$$\left[U_N=U_N^{(1)}\right]\subset\left[U_N^*=U_N^{(2)}\right]\subset\left[U_N^{**}=U_N^{(\infty)}\right]$$

We will do this, construction of $U_N^{(r)}$ and study of this new quantum group, following what we know at $r=1,2,\infty$, in several steps. Following \cite{bb1}, let us start with:

\begin{definition}
The quantum group $U_N^{(r)}\subset U_N^{**}$ is defined according to
$$C(U_N^{(r)})=C(U_N^{**})\Big/\Big<[u_{i_1j_1}\ldots u_{i_rj_r},u_{k_1l_1}\ldots u_{k_rl_r}]=0\Big>$$
with the convention that at $r=\infty$, the relations on the right dissapear.
\end{definition}

As a first observation, the quantum space $U_N^{(r)}\subset U_N^{**}$ constructed above is indeed a quantum group, because the relations on the right are of Tannakian nature. In fact, the above relations are of ``easy'' type, so we can say more, as follows:

\begin{proposition}
The quantum group $U_N^{(r)}\subset U_N^{**}$ is easy, coming from
$$
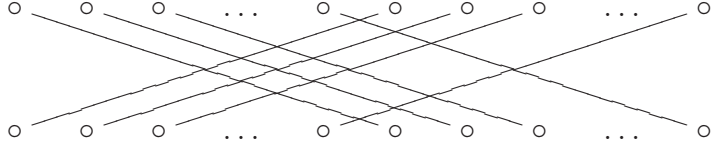
\xymatrix@R=12mm@C=5mm{\circ \ar@{-}[drrrrr]& \circ \ar@{-}[drrrrr] & \circ \ar@{-}[drrrrr]& \ldots & \circ \ar@{-}[drrrrr]& \circ\ar@{-}[dlllll]&\circ\ar@{-}[dlllll]&\circ\ar@{-}[dlllll] & \ldots & \circ\ar@{-}[dlllll] \\ \circ & \circ & \circ& \ldots & \circ & \circ&\circ&\circ& \ldots & \circ } 
$$
regarded as element of $\mathcal P_2(2r,2r)$.
\end{proposition}

\begin{proof}
If we denote by $\pi\in \mathcal P_2(2r,2r)$ the diagram in the statement, an elementary computation  shows that we have the following equivalence:
$$T_\pi\in End(u^{\otimes r})\iff[u_{i_1j_1}\ldots u_{i_rj_r},u_{k_1l_1}\ldots u_{k_rl_r}]=0$$

Thus, we are led to the conclusion in the statement.
\end{proof}

As another basic observation, that we will need as well in what follows, we have:

\begin{proposition}
The standard coordinates of $U_N^{(r)}$ satisfy the relations
$$[u_{i_1j_1}\ldots u_{i_rj_r},u^*_{k_1l_1}\ldots u^*_{k_rl_r}]=0$$
with, as usual, the convention that these relations dissapear at $r=\infty$.
\end{proposition}

\begin{proof}
This follows from the well-known fact that if the coefficients of two unitary corepresentations $(u_{ij})$, $(v_{kl})$ of a quantum group pairwise commute, then the coefficients of $(u_{ij})$, $(v_{kl}^*)$ also pairwise commute. To check this, start with the following relations:
$$u_{ij}v_{kl}=v_{kl}u_{ij}$$

If we multiply on the right by $v_{pl}^*$, and then sum over $l$, we get:
$$\delta_{kp}u_{ij}= \sum_l v_{kl}u_{ij}v_{pl}^*$$

Now multiplying on the left by $v_{kq}^*$ and summing over $k$ gives:
$$v_{pq}^*u_{ij}=u_{ij}v_{pq}^*$$

Thus, we are led to the conclusion in the statement.
\end{proof}

We can now say more about our new quantum groups, as follows:

\begin{theorem}
The quantum groups $U_N^{(r)}$ interpolate between $U_N$ and $U_N^{**}$,
$$U_N\subset U_N^{(r)}\subset U_N^{**}$$
and we have equalities as follows, in relation with our previous quantum groups:
$$\left[U_N=U_N^{(1)}\right]\subset\left[U_N^*=U_N^{(2)}\right]\subset\left[U_N^{**}=U_N^{(\infty)}\right]$$
Moreover, all these quantum groups are easy.
\end{theorem}

\begin{proof}
We have several things to be proved, the idea being as follows:

\medskip

(1) Our first claim is that we have $U_N^{(1)}=U_N$. Indeed, according to Definition 8.18, the quantum group $U_N^{(1)}$ is given by the following formula:
$$C(U_N^{(1)})=C(U_N^{**})\Big/\Big<[u_{ij},u_{kl}]=0\Big>$$

On the other hand, by taking into account Proposition 8.20 as well, we get:
$$C(U_N^{(1)})=C(U_N^{**})\Big/\Big<[u_{ij},u_{kl}]=[u_{ij},u_{kl}^*]=0\Big>$$

But this shows that we have $U_N^{(1)}=U_N$, as claimed.

\medskip

(2) Our second claim is that we have $U_N^{(2)}=U_N^*$. Indeed, according to Definition 8.18, the quantum group $U_N^{(2)}$ is given by the following formula:
$$C(U_N^{(2)})=C(U_N^{**})\Big/\Big<[u_{i_1j_1}u_{i_2j_2},u_{k_1l_1}u_{k_2l_2}]=0\Big>$$

As before, by using as well Proposition 8.20, we get:
$$C(U_N^{(2)})=C(U_N^{**})\Big/\Big<[u_{i_1j_1}u_{i_2j_2},u_{k_1l_1}u_{k_2l_2}]=[u_{i_1j_1}u_{i_2j_2},u^*_{k_1l_1}u^*_{k_2l_2}]=0\Big>$$

But this shows that we have $U_N^{(2)}=U_N^*$, as claimed.

\medskip

(3) But with the above claims in hand, everything from the statement follows, and with the easiness assertion being known from Proposition 8.19.
\end{proof}

Getting now to a more detailed study of $U_N^{(r)}$, the general idea here, from what we have from Theorem 8.21, will be that, perhaps leaving sometimes the limiting cases $r=1,\infty$ aside, our quantum group $U_N^{(r)}$ appears as some sort of technical version of $U_N^*$. So, we should use the same methods as for $U_N^*$, namely matrix models, or easiness.

\bigskip

Following \cite{bb1}, let us start with matrix model techniques. We have here:

\begin{theorem}
We have an embedding of quantum groups, as follows,
$$U_N^{(r)}\subset U_N^r\rtimes\mathbb Z_r$$
and a related cyclic matrix model, as follows,
$$C(U_N^{(r)})\subset M_r(C(U_N^r))$$
and in this latter model, $\int_{U_N^{(r)}}$ appears as the restriction of $tr_r\otimes\int_{U_N^r}$.
\end{theorem}

\begin{proof}
All this is quite routine, following our study of half-liberation from chapter 7, which corresponds to the case $r=2$, and with the remark that at $r=1$ everything is trivial. For details on all this, we refer to \cite{bb1}. But, we will be back to this, right next.
\end{proof}

More generally now, again in analogy with the results from chapter 7, we have:

\begin{proposition}
If $L$ is a compact group, having a $N$-dimensional unitary corepresentation $v$, and an order $r$ automorphism $\sigma:L\to L$, we have a matrix model
$$\pi:C(U_N^{**})\to M_K(C(L))\quad,\quad u_{ij}\to\tau[v_{ij}^{(1)},\ldots,v_{ij}^{(r)}]$$
where $v^{(i)}(g)=v(\sigma^i(g))$, and where $\tau[x_1,\ldots,x_r]$ is obtained by filling the standard $r$-cycle $\tau\in M_r(0,1)$ with the elements $x_1,\ldots,x_r$. We call such models ``cyclic''.
\end{proposition}

\begin{proof}
The matrices $U_{ij}=\tau[v_{ij}^{(1)},\ldots,v_{ij}^{(r)}]$ in the statement appear by definition as follows, with the convention that all the blank spaces denote 0 entries:
$$U_{ij}=\begin{pmatrix}
&&&v_{ij}^{(1)}\\
v_{ij}^{(2)}\\
&\ddots\\
&&v_{ij}^{(r)}
\end{pmatrix}$$

The matrix $U=(U_{ij})$ is then unitary, and so is $\bar{U}=(U_{ij}^*)$. Thus, if we denote by $w=(w_{ij})$ the fundamental corepresentation of $C(U_N^+)$, we have a model as follows:
$$\rho:C(U_N^+)\to M_r(C(L))\quad,\quad w_{ij}\to U_{ij}$$

Now observe that the matrices $U_{ij}U_{kl}^*,U_{ij}^*U_{kl}$ are all diagonal, so in particular, they commute. Thus the above morphism $\rho$ factorizes through $C(U_N^{**})$, as claimed.
\end{proof}

In relation to the above models, we have the following result:

\begin{theorem}
Any cyclic model in the above sense,
$$\pi:C(U_N^{**})\to M_r(C(L))$$
is stationary on its image, with the corresponding closed subgroup $[L]\subset U_N^{**}$, given by 
$$Im(\pi)=C([L])$$
being the quotient $L\rtimes\mathbb Z_r\to[L]$ having as coordinates the variables $u_{ij}=v_{ij}\otimes\tau$.
\end{theorem}

\begin{proof}
Assuming that $(L,\sigma)$ are as in Proposition 8.23, we have an action $\mathbb Z_r\curvearrowright L$, and we can therefore consider the following short exact sequence:
$$1\to\mathbb Z_r\to L\rtimes\mathbb Z_r\to L\to1$$

By doing some standard algebra, we obtain from this a model as follows, where $x^{(i)}=\tilde{\sigma}^i(x)$, with $\tilde{\sigma}:C(L)\to C(L)$ being the automorphism induced by $\sigma:L\to L$:
$$\rho:C(L\rtimes\mathbb Z_r)\subset M_r(C(L))\quad,\quad x\otimes\tau^i\to\tau^i[x^{(1)},\ldots,x^{(r)}]$$

Consider now the quotient quantum group $L\rtimes\mathbb Z_r\to[L]$ having as coordinates the variables $u_{ij}=v_{ij}\otimes\tau$. We have then a injective morphism, as follows:
$$\nu:C([L])\subset C(L\rtimes\mathbb Z_r)\quad,\quad u_{ij}\to v_{ij}\otimes\tau$$

By composing the above two embeddings, we obtain an embedding as follows:
$$\rho\nu:C([L])\subset M_r(C(L))\quad,\quad u_{ij}\to\tau[v_{ij}^{(1)},\ldots,v_{ij}^{(K)}]$$

Now since $\rho$ is stationary, and since $\nu$ commutes with the Haar funtionals as well, it follows that this morphism $\rho\nu$ is stationary, and this finishes the proof.
\end{proof}

As an illustration, we can now recover the following result, from \cite{bdu}:

\begin{proposition}
For any non-classical $G\subset O_N^*$ we have a stationary model
$$\pi:C(G)\to M_2(C(L))\quad,\quad u_{ij}=\begin{pmatrix}0&v_{ij}\\ \bar{v}_{ij}&0\end{pmatrix}$$
where $L\subset U_N$, with coordinates denoted $v_{ij}$, is the lift of $PG\subset PO_N^*=PU_N$.
\end{proposition}

\begin{proof}
Assume first that $L\subset U_N$ is self-conjugate, in the sense that $g\in L\implies\bar{g}\in L$. If we consider the order 2 automorphism of $C(L)$ induced by $g_{ij}\to\bar{g}_{ij}$, we can apply Theorem 8.24, and we obtain a stationary model, as follows:
$$\pi:C([L])\subset M_2(C(L))\quad,\quad u_{ij}\otimes1=\begin{pmatrix}0&v_{ij}\\ \bar{v}_{ij}&0\end{pmatrix}$$

The point now is that, as explained in \cite{bdu}, any non-classical subgroup $G\subset O_N^*$ must appear as $G=[L]$, for a certain self-conjugate subgroup $L\subset U_N$. Moreover, since we have $PG=P[L]$, it follows that $L\subset U_N$ is the lift of $PG\subset PO_N^*=PU_N$, as claimed. 
\end{proof}

In the unitary case now, we have the following result:

\begin{theorem}
For any subgroup $G\subset U_N^{(r)}$ which is $r$-symmetric, in the sense that $u_{ij}\to e^{2\pi i/r}u_{ij}$ defines an automorphism of $C(G)$, we have a stationary model
$$\pi:C(G)\to M_r(C(L))\quad,\quad u_{ij}\to\tau[v_{ij}^{(1)},\ldots,v_{ij}^{(r)}]$$
with $L\subset U_N^r$ being a closed subgroup which is symmetric, in the sense that it is stable under the cyclic action $\mathbb Z_r\curvearrowright U_N^r$.
\end{theorem}

\begin{proof}
This follows from what we have, as follows:

\medskip

(1) Assuming that $L\subset U_N^r$ is symmetric in the above sense, we have representations $v^{(i)}:L\subset U_N^r\to U_N^{(i)}$ for any $i$, and the cyclic action $\mathbb Z_r\curvearrowright U_N^r$ restricts into an order $r$ automorphism $\sigma:L\to L$. Thus we can apply Theorem 8.24, and we obtain a certain closed subgroup $[L]\subset U_N^{(r)}$, having a stationary model as in the statement.

\medskip

(2) Conversely now, assuming that we have a subgroup $G\subset U_N^{(r)}$ which is $r$-symmetric, we must have $C(G)\subset C(L)\rtimes\mathbb Z_r$, for a certain closed subgroup $L\subset U_N^r$ which is symmetric. But this shows that we have $G=[L]$, and we are done.
\end{proof}

The above results can be used in order to say many things about the quantum groups $U_N^{(r)}$, and in particular to compute, via matrix models, the asymptotic laws of characters. We obtain in this way extensions of our previous results regarding $U_N^*$:

\begin{theorem}
The asymptotic laws of truncated characters for $U_N^{(r)}$ can be computed by using matrix models, and we obtain generalizations of the results regarding $U_N^*$.
\end{theorem}

\begin{proof}
We know from Theorem 8.22 that we have a matrix model, as follows:
$$C(U_N^{(r)})\subset M_r(C(U_N^r))\quad,\quad \int_{U_N^{(r)}}=tr_r\otimes\int_{U_N^r}$$

Thus, everything can be computed in the model, and we get the results.
\end{proof}

As already mentioned, all the above, using matrix models, is only half of the story, because we have the easiness methods available as well. Following \cite{mw1}, we have:

\begin{theorem}
The quantum group $U_N^{(r)}$ is easy, coming from the category $\mathcal P_2^{(r)}$ of matching pairings having the property that
$$\#\circ=\#\bullet(r)$$
between the legs of each string, when flattened.
\end{theorem}

\begin{proof}
This can be done in several steps, as follows:

\medskip

(1) At $r=1$ there is nothing to prove, because we know from Theorem 8.21 that we have $U_N^{(1)}=U_N$, corresponding to the category $\mathcal P_2^{(1)}=\mathcal P_2$.

\medskip

(2) At $r=\infty$ now, we know from Theorem 8.21 that we have $U_N^{(\infty)}=U_N^{**}$. On the other hand, we know from Theorem 8.17 that the quantum group $U_N^{**}$ is easy, coming from the category $\mathcal P_2^{(\infty)}=\mathcal P_2^{**}$, so our problem is solved at $r=\infty$ as well.

\medskip

(3) In the general case now, assuming $r<\infty$, we know from Proposition 8.19 that $U_N^{(r)}\subset U_N^{**}$ comes from the following partition in $\mathcal P_2(2r,2r)$:
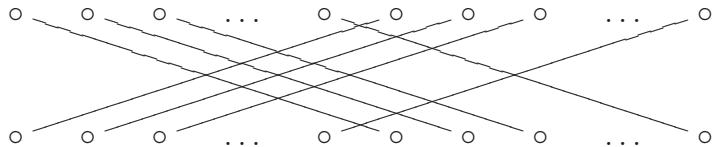
$$\xymatrix@R=12mm@C=5mm{\circ \ar@{-}[drrrrr]& \circ \ar@{-}[drrrrr] & \circ \ar@{-}[drrrrr]& \ldots & \circ \ar@{-}[drrrrr]& \circ\ar@{-}[dlllll]&\circ\ar@{-}[dlllll]&\circ\ar@{-}[dlllll] & \ldots & \circ\ar@{-}[dlllll] \\ \circ & \circ & \circ& \ldots & \circ & \circ&\circ&\circ& \ldots & \circ } 
$$

By rotating this diagram, we can see that the condition in the statement, namely $\#\circ=\#\bullet(r)$ between the legs of each string, is satisfied. But this leads, via some standard combinatorics, to the conclusion in the statement, and we refer here to \cite{mw1}.

\medskip

(4) As an illustration for all this, let us work out the case $r=2$. Here we know that the standard partitions producing $U_N^{(2)}=U_N^*$ are the various colorings of the half-classical crossing, which can be generically denoted as follows, with $a,b,c\in\{\circ,\bullet\}$:
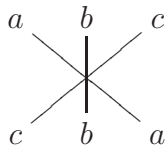
$$\xymatrix@R=10mm@C=5mm{
a\ar@{-}[drr]&b\ar@{-}[d]&c\ar@{-}[dll]
\\
c&b&a}$$

Now by rotating this diagram to the right, as to flatten it, we obtain three interlacing semicircles, with legs labelled as follows:
$$c\ \ b\ \ a\ \ \bar{c}\ \ \bar{b}\ \ \bar{a}$$

Now since each of our three interlacing semicircles has exactly 2 points between its legs, we have $\#\circ+\#\bullet=0(2)$, and so $\#\circ=\#\bullet(2)$, between the legs of each semicircle, regardless of the labels $a,b,c\in\{\circ,\bullet\}$. Thus, we can see that the partitions producing $U_N^{(2)}=U_N^*$ belong to $\mathcal P_2^{(2)}$, and with a bit more work, consisting in showing on pictures that these partitions actually generate $\mathcal P_2^{(2)}$, we are led to the result.
\end{proof}

The above result is quite powerful, and can stand as an alternative to Theorem 8.22, for instance for proving the various probabilistic results from Theorem 8.27, via pure combinatorics. However, in practice, nothing beats matrix models and some calculus, so our proof above of Theorem 8.27, using Theorem 8.22, is the simplest one.

\bigskip

Finally, let us mention that our presentation above of the quantum group $U_N^{(r)}$ was one among others. It is possible to have as well an ``easiness first'' viewpoint on all this, and with this idea in mind, what we have about $U_N^{(r)}$ can be summarized as follows:

\index{unitary quantum group}
\index{cyclic matrix model}

\begin{theorem}
Associated to any $r\in\mathbb N$ is the quantum group 
$$U_N\subset U_N^{(r)}\subset U_N^+$$
coming from the category $\mathcal P_2^{(r)}$ of matching pairings having the property that $$\#\circ=\#\bullet(r)$$
holds between the legs of each string. These quantum groups are as follows:
\begin{enumerate}
\item At $r=1$ we obtain the usual unitary group, $U_N^{(1)}=U_N$.

\item At $r=2$ we obtain the half-classical unitary group, $U_N^{(2)}=U_N^*$.

\item At $r=\infty$ we obtain the quantum group $U_N^{(\infty)}=U_N^{**}$.

\item For any $r|s$ we have an embedding $U_N^{(r)}\subset U_N^{(s)}$.

\item In general, we have an embedding $U_N^{(r)}\subset U_N^r\rtimes\mathbb Z_r$.

\item We have as well a cyclic matrix model $C(U_N^{(r)})\subset M_r(C(U_N^r))$.

\item In this latter model, $\int_{U_N^{(r)}}$ appears as the restriction of $tr_r\otimes\int_{U_N^r}$.
\end{enumerate}
\end{theorem}

\begin{proof}
This is something quite compact, summarizing what we have about $U_N^{(r)}$, with some theorems transformed into definitions, and vice versa. To be more precise:

\medskip

(1) This follows from Theorem 8.21 and Theorem 8.28.

\medskip

(2) This follows again from Theorem 8.21 and Theorem 8.28.

\medskip

(3) Again, this follows from Theorem 8.21 and Theorem 8.28.

\medskip

(4) This follows from the functoriality of the Tannakian correspondence, because when assuming $r|s$ we have $\mathcal P_2^{(s)}\subset\mathcal P_2^{(r)}$, and so $U_N^{(r)}\subset U_N^{(s)}$, as claimed.

\medskip

(5) This follows from Theorem 8.22, using Theorem 8.28.

\medskip

(6) This follows again from Theorem 8.22 and Theorem 8.28.

\medskip

(7) Again, this follows from Theorem 8.22 and Theorem 8.28.
\end{proof}

There are many other interesting things that can be said about the quantum groups $U_N^{(r)}$, and we refer here to \cite{bb2}, \cite{bb3}, \cite{mw1}, \cite{mw2} and the subsequent literature. In what concerns us, we will be back to these quantum groups in chapter 16 below, with some results in relation with their noncommutative gometry meaning.

\section*{8d. The standard family} 

Let us discuss now the second known construction of unitary quantum groups, from \cite{mw2}. This construction uses an additive semigroup $D\subset\mathbb N$, but as pointed out there, using instead the complementary set $C=\mathbb N-D$ leads to several simplifications. So, let us call ``cosemigroup'' any subset $C\subset\mathbb N$ which is complementary to an additive semigroup, $x,y\notin C\implies x+y\notin C$. The construction from \cite{mw2} is then as follows:

\index{cosemigroup}
\index{unitary quantum group}
\index{intermediate liberation}

\begin{theorem}
Associated to any cosemigroup $C\subset\mathbb N$ is the easy quantum group 
$$U_N^{(\infty)}\subset U_N^C\subset U_N^+$$
coming from the category $\mathcal P_2^C\subset\mathcal P_2^{(\infty)}$ of pairings having the property 
$$\#\circ-\#\bullet\in C$$
between each two legs colored $\circ,\bullet$ of two strings which cross. We have:
\begin{enumerate}
\item For $C=\mathbb N$ we obtain the quantum group $U_N^{(\infty)}$.

\item For $C=\emptyset$ we obtain the quantum group $U_N^+$.

\item For $C=\{0\}$ we obtain the quantum group $U_N^\times$.

\item For $C\subset C'$ we have an inclusion $U_N^{C'}\subset U_N^C$.
\end{enumerate}
\end{theorem}

\begin{proof}
Once again this is something very compact, coming from work in \cite{mw2}, with our convention that the semigroup $D\subset\mathbb N$ which is used there is replaced here by its complement $C=\mathbb N-D$. Here are a few explanations on all this:

\medskip

(1) The assumption $C=\mathbb N$ simply tells us that the condition $\#\circ-\#\bullet\in C$ in the statement is irrelevant. Thus, we have $\mathcal P_2^\mathbb N=\mathcal P_2^{(\infty)}$, and so $U_N^\mathbb N=U_N^{(\infty)}$.

\medskip

(2) The assumption $C=\emptyset$ means that the condition $\#\circ-\#\bullet\in C$ can never be applied. Thus, the strings cannot cross, we have $\mathcal P_2^\emptyset=\mathcal{NC}_2$, and so $U_N^\emptyset=U_N^+$.

\medskip

(3) The assumption $C=\{0\}$ means that the pairings in $\mathcal P_2^C\subset\mathcal P_2^{(\infty)}$ must satisfy the condition $\#\circ=\#\bullet$, between each two legs colored $\circ,\bullet$ of two strings which cross. Now consider the standard diagram producing the quantum group $U_N^\times$, namely:
$$
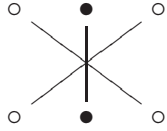
\xymatrix@R=10mm@C=5mm{
\circ\ar@{-}[drr]&\bullet\ar@{-}[d]&\circ\ar@{-}[dll]
\\
\circ&\bullet&\circ}$$

By rotating to the right, as to have this diagram flattened, we obtain three interlacing semicircles, with legs labelled $\circ\bullet\circ\bullet\circ\,\bullet\,$. Thus, we can see that the condition $\#\circ=\#\bullet$, between each two legs colored $\circ,\bullet$ of two strings which cross, holds indeed for this diagram, and with a bit more work, as explained in \cite{mw2}, we obtain from this that we have in fact $\mathcal P_2^{\{0\}}=\mathcal P_2^\times$, and so that we have $U_N^{\{0\}}=U_N^\times$, as claimed.

\medskip

(4) This is clear by functoriality, because $C\subset C'$ implies $\mathcal P_2^{C}\subset\mathcal P_2^{C'}$.
\end{proof}

We have the following key result, from \cite{mw2}:

\index{unitary quantum group}
\index{intermediate liberation}

\begin{theorem}
The easy quantum groups $U_N\subset G\subset U_N^+$ are as follows,
$$U_N\subset\{U_N^{(r)}\}\subset\{U_N^C\}\subset U_N^+$$
with the series covering $U_N$, and the family covering $U_N^+$.
\end{theorem}

\begin{proof}
This is something non-trivial, and we refer here to \cite{mw2}. The general idea is that $U_N^{(\infty)}$ produces a dichotomy for the quantum groups in the statement, and this leads, via some combinatorial computations, to the series and the family. See \cite{mw1}, \cite{mw2}.
\end{proof}

All the above is quite exciting, because we have now a complete point of view on intermediate liberations, at least in the unitary case. We will be back to this on several occasions, first in Part III after discussing some similar problems for the reflection groups too, and then at the end of Part IV, with a discussion of the key problem of constructing full noncommutative geometry theories, based on the liberations of $U_N$ that we have.

\bigskip

As a last topic that we would like to discuss here, we have the notion of projective easiness, which is something general and of independent interest, related to all the above. Let us start with the following straightforward definition:

\begin{definition}
A projective category of pairings is a collection of subsets 
$$NC_2(2k,2l)\subset E(k,l)\subset P_2(2k,2l)$$
stable under the usual categorical operations, and satisfying 
$$\sigma\in E\implies |\sigma|\in E$$
with the vertical bars standing for vertical strings.
\end{definition}

As basic examples here, we have the following projective categories of pairings, where $P_2^*$ is the category of matching pairings:
$$NC_2\subset P_2^*\subset P_2$$

This follows indeed from definitions. Now with the above notion in hand, we can formulate the following projective analogue of the notion of easiness:

\begin{definition}
An intermediate compact quantum group 
$$PO_N\subset H\subset PO_N^+$$
is called projectively easy when its Tannakian category
$$span(NC_2(2k,2l))\subset Hom(v^{\otimes k},v^{\otimes l})\subset span(P_2(2k,2l))$$
comes via via the following formula, using the standard $\pi\to T_\pi$ construction,
$$Hom(v^{\otimes k},v^{\otimes l})=span(E(k,l))$$
for a certain projective category of pairings $E=(E(k,l))$.
\end{definition}

Thus, we have a projective notion of easiness. Observe that, given an easy quantum group $O_N\subset G\subset O_N^+$, its projective version $PO_N\subset PG\subset PO_N^+$ is projectively easy in our sense. In particular the basic projective quantum groups $PO_N\subset PU_N\subset PO_N^+$ are all projectively easy in our sense, coming from the categories $NC_2\subset P_2^*\subset P_2$. 

\bigskip

We have in fact the following general result:

\begin{theorem}
We have a bijective correspondence between the affine and projective categories of partitions, given by the operation
$$G\to PG$$ 
at the level of the corresponding affine and projective easy quantum groups.
\end{theorem}

\begin{proof}
The construction of correspondence $D\to E$ is clear, simply by setting:
$$E(k,l)=D(2k,2l)$$

Indeed, due to the easiness axioms, the conditions in Definition 8.32 are satisfied. Conversely, given $E=(E(k,l))$ as in Definition 8.32, we can set:
$$D(k,l)=\begin{cases}
E(k,l)&(k,l\ {\rm even})\\
\{\sigma:|\sigma\in E(k+1,l+1)\}&(k,l\ {\rm odd})
\end{cases}$$

Our claim is that $D=(D(k,l))$ is a category of partitions. Indeed:

\medskip

(1) The composition action is clear. Indeed, when looking at the numbers of legs involved, in the even case this is clear, and in the odd case, this follows from:
\begin{eqnarray*}
|\sigma,|\sigma'\in E
&\implies&|^\sigma_\tau\in E\\
&\implies&{\ }^\sigma_\tau\in D
\end{eqnarray*}

(2) For the tensor product axiom, we have 4 cases to be investigated, depending on the parity of the number of legs of $\sigma,\tau$, as follows:

\medskip

-- The even/even case is clear. 

\medskip

-- The odd/even case follows from the following computation:
\begin{eqnarray*}
|\sigma,\tau\in E
&\implies&|\sigma\tau\in E\\
&\implies&\sigma\tau\in D
\end{eqnarray*}

-- Regarding now the even/odd case, this can be solved as follows:
\begin{eqnarray*}
\sigma,|\tau\in E
&\implies&|\sigma|,|\tau\in E\\
&\implies&|\sigma||\tau\in E\\
&\implies&|\sigma\tau\in E\\
&\implies&\sigma\tau\in D
\end{eqnarray*}

-- As for the remaining odd/odd case, here the computation is as follows:
\begin{eqnarray*}
|\sigma,|\tau\in E
&\implies&||\sigma|,|\tau\in E\\
&\implies&||\sigma||\tau\in E\\
&\implies&\sigma\tau\in E\\
&\implies&\sigma\tau\in D
\end{eqnarray*}

(3) Finally, the conjugation axiom is clear from definitions. It is also clear that both compositions $D\to E\to D$ and $E\to D\to E$ are the identities, as claimed. As for the quantum group assertion, this is clear as well from definitions.
\end{proof}

Now back to uniqueness issues, we have here the following result:

\begin{theorem}
The following happen:
\begin{enumerate}
\item $O_N^*$ is the only easy quantum group $O_N\subset G\subset O_N^+$.

\item $PU_N$ is the only projectively easy quantum group $PO_N\subset G\subset PO_N^+$.
\end{enumerate}
\end{theorem}

\begin{proof}
The idea here is as follows:

\medskip

(1) The assertion regarding $O_N\subset O_N^*\subset O_N^+$ is from \cite{bv2}, and this is something that we already know, explained in chapter 7.

\medskip

(2) The assertion regarding $PO_N\subset PU_N\subset PO_N^+$ follows from the classification result in (1), and from the duality in Theorem 8.34.
\end{proof}

Summarizing, we have a nice notion of projective easiness, which is in relation with the liberations of $U_N$ too, via the isomorphism $PO_N^+=PU_N^+$. We will be back to this in chapter 16, when discussing noncommutative geometry.

\section*{8e. Exercises}

As before with the last few chapters, we are rather into research matters here, and as a good exercise on all this, which is definitely of research type, we have:

\begin{exercise}
Work out the asymptotic laws of truncated characters, and other probabilistic aspects, for the quantum groups $U_N^C$.
\end{exercise}

Obvioiusly, this is something quite complicated, because in contrast for instance with what happens for the quantum groups $U_N^{(r)}$, where we have explicit matrix models for doing our computations, for the quantum groups $U_N^C$ there is nothing to rely upon, apart from the quite heavy combinatorics producing their definition. However, we have 3 examples that we are familiar with, from Theorem 8.30 (1,2,3), so the first step towards solving the exercise would be that of conjecturing something, based on these 3 examples.

\part{The discrete case}

\ \vskip50mm

\begin{center}
{\em It's in your eyes

I can tell what you're thinking

My heart is sinking too

It's no surprise}
\end{center}

\chapter{Real reflections}

\section*{9a. Basic series}

In this third part of the present book we investigate the discrete, or reflection easy quantum group case. This is the most mysterious case, with many interesting and unexpected examples involved, and reminding the magic world of the complex reflection groups. The study here being the most difficult, a large part of our results will spill into results regarding the general case too. In order to explain our strategy, let us go back to the standard cube formed by the main easy quantum groups, namely: 
$$
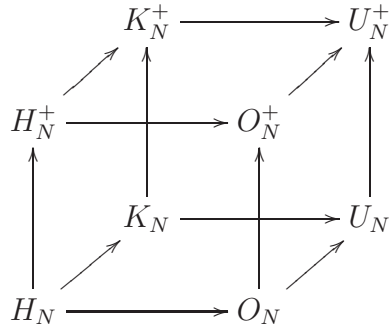
\xymatrix@R=18pt@C=18pt{
&K_N^+\ar[rr]&&U_N^+\\
H_N^+\ar[rr]\ar[ur]&&O_N^+\ar[ur]\\
&K_N\ar[rr]\ar[uu]&&U_N\ar[uu]\\
H_N\ar[uu]\ar[ur]\ar[rr]&&O_N\ar[uu]\ar[ur]
}$$

We will be mainly interested in the left face of the cube, which is of ``discrete'' nature, as opposed to the right face, which is of obvious continuous nature. This left face contains the main examples of quantum reflection groups that we have, namely:
$$
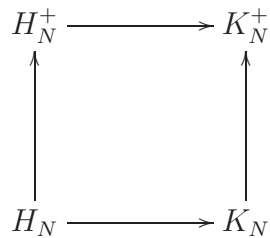
\xymatrix@R=56pt@C=55pt{
H_N^+\ar[r]&K_N^+\\
H_N\ar[u]\ar[r]&K_N\ar[u]}$$

As a first observation, we know in fact far many more quantum reflection groups that this. Indeed, just by looking at the main examples of classical and free quantum groups that we know, we have in fact the following rectangular diagram, with $s\in\{1,2,\ldots,\infty\}$, and with the values $s=1,\infty$ covering the vertices of the rectangle:
$$\xymatrix@R=56pt@C=55pt{
S_N^+\ar[r]&H_N^{s+}\ar[r]&K_N^+\\
S_N\ar[u]\ar[r]&H_N^s\ar[u]\ar[r]&K_N\ar[u]}$$

In addition to this, we can talk about half-liberations, with $K_N^*,H_N^*$ being quantum groups that we already met, and we can talk as well about products with $\mathbb Z_2$, with the quantum group $S_N^\circ=S_N\times \mathbb Z_2$ being again a quantum group that we already met.

\bigskip

In short, we have many examples of quantum reflection groups, and are a bit in trouble with starting something. A good idea would be that of restricting the attention to the real case. Here, with respect to the above rectangular diagram, only the cases $s=1,2$ qualify, and we are left with a rather simple square diagram, as follows:
$$\xymatrix@R=56pt@C=55pt{
S_N^+\ar[r]&H_N^+\\
S_N\ar[u]\ar[r]&H_N\ar[u]}$$

However, as mentioned above, this is not all, because we can talk for instance about the half-liberation $H_N^*$, and with the remark that this will not fit well into our diagram, due to the collapsing result $S_N^*=S_N$. Also, we can talk as well about the quantum groups $G_N^\circ=G_N\times\mathbb Z_2$, and once again here with a bad functoriality remark, namely that we have $H_N^\circ=H_N$. And finally, even worse on this topic, as we will soon discover, there are in fact uncountably many examples of easy quantum groups  $S_N\subset G\subset H_N^+$.

\bigskip

Looks like we are completely lost, so time to ask the cat. And cat says:

\begin{cat}
Go for the beauty, the easy kill is $H_N^\times$.
\end{cat}

Thanks cat, and although there is some confusion here between beauty, easiness, predator and prey and so on, at least from my peaceful human viewpoint, I must confess that, forgetting all the math that I know, $H_N$ looks to be indeed the most beautiful reflection group of them all. So, we should go for its liberations $H_N^\times$, and whether that will be an easy or difficult task, and then what to do afterwards, remains to be seen.

\bigskip

Anticipating now a bit, in order to make a plan out of this, for the present Part III of the present book, let us go back to the left face of the standard cube, containing the main examples of quantum reflection groups that we have, namely:
$$\xymatrix@R=56pt@C=55pt{
H_N^+\ar[r]&K_N^+\\
H_N\ar[u]\ar[r]&K_N\ar[u]}$$

We will investigate in this chapter, following the early paper \cite{ez1}, and then the papers of Raum-Weber \cite{rw1}, \cite{rw2}, \cite{rw3}, the intermediate subgroup question for the left edge, $H_N\subset G\subset H_N^+$. Then, in chapter 10 we will extend this study to the case of the intermediate subgroups $S_N\subset G\subset H_N^+$, and with a bit more work involved, we will obtain from this a classification result in the real case, $S_N\subset G\subset O_N^+$. 

\bigskip

Afterwards, in chapter 11 we will do a similar work for the right edge, $K_N\subset G\subset K_N^+$, and more generally for the intermediate quantum groups $H_N\subset G\subset K_N^+$, and even more generally, for the intermediate quantum groups $S_N\subset G\subset K_N^+$. And finally, in chapter 12 we will discuss, still following \cite{ez1} and Raum-Weber \cite{rw1}, \cite{rw2}, \cite{rw3}, and the more recent work of Mang-Weber \cite{mw3}, \cite{mw4}, various structure and classification results for the whole left face of the cube, and then for the general easy quantum groups $S_N\subset G\subset U_N^+$.

\bigskip

Getting to work now, for the reasons explained above, we will be first interested in this chapter in the various easy liberations of the hyperoctahedral group $H_N$:

\begin{question}
What are the intermediate easy quantum groups
$$H_N\subset G\subset H_N^+$$
lying between $H_N=\mathbb Z_2\wr S_N$, and its free version $H_N^+=\mathbb Z_2\wr_*S_N^+$?
\end{question}

We have so far only three examples of such quantum groups, namely the endpoints $H_N,H_N^+$ themselves, and the half-classical quantum group $H_N^*$, sitting in the middle:
$$H_N\subset H_N^*\subset H_N^+$$

In order to construct more examples, let us first look for intermediate objects for the inclusion on the left, $H_N\subset G\subset H_N^*$. Following \cite{ez1}, we can first introduce a new series of quantum groups, $H_N^{(s)}$ with $s\in\{2,3,\ldots,\infty\}$, which ``interpolates'' between the endpoints $H_N^{(2)}=H_N$ and $H_N^{(\infty)}=H_N^*$, the inclusions being as follows:
$$H_N=H_N^{(2)}\subset\ldots\subset H_N^{(s)}\subset\ldots\subset H_N^{(\infty)}=H_N^*$$

To be more precise, let us define $H_N^{(s)}$ as being obtained from $H_N^*$ by imposing the ``$s$-commutation'' condition $abab\ldots=baba\ldots$ (length $s$ words) to the basic coordinates $u_{ij}$. It is convenient to write down the complete definition of $H_N^{(s)}$, as follows:

\index{cubic relations}
\index{s-mixing relation}

\begin{definition}
$C(H_N^{(s)})$ is the universal $C^*$-algebra generated by $N^2$ self-adjoint variables $u_{ij}$, subject to the following relations:
\begin{enumerate}
\item Orthogonality: $uu^t=u^tu=1$, where $u=(u_{ij})$ and $u^t=(u_{ji})$.

\item Cubic relations: $u_{ij}u_{ik}=u_{ji}u_{ki}=0$, for any $i$ and any $j\neq k$.

\item Half-commutation: $abc=cba$, for any $a,b,c\in \{u_{ij}\}$.

\item $s$-mixing relation: $abab\ldots=baba\ldots$ (length $s$ words), for any $a,b\in\{u_{ij}\}$.
\end{enumerate}
\end{definition}

Observe that at $s=2$ the $s$-mixing relation is the usual commutation $ab=ba$. This relation being stronger than the half-commutation $abc=cba$, we are led to the algebra generated by $N^2$ commuting self-adjoint variables satisfying (1,2), which is $C(H_N)$:
$$H_N^{(2)}=H_N$$

As for the case $s=\infty$, here according to our usual conventions regarding relations of infinite length, used throughout this book, the $s$-mixing relation disappears by definition. Thus we are led to the algebra defined by the relations (1,2,3), which is $C(H_N^*)$:
$$H_N^{(\infty)}=H_N^*$$
 
Summarizing, Definition 9.3 provides us indeed with a new series of hyperoctahedral quantum groups, which are, as previously claimed, as follows:
$$H_N=H_N^{(2)}\subset\ldots\subset H_N^{(s)}\subset\ldots\subset H_N^{(\infty)}=H_N^*$$

All this is quite interesting, so let us present now a detailed study of $H_N^{(s)}$, from an algebraic and probabilistic viewpoint. Our first technical result is as follows:

\begin{proposition}
For a closed subgroup $G\subset H_N^*$, the following are equivalent:
\begin{enumerate}

\item The basic coordinates $u_{ij}$ satisfy $abab\ldots=baba\ldots$ (length $s$ words).

\item We have $T_\pi\in End(u^{\otimes s})$, where $\pi=(135\ldots 2'4'6'\ldots)(246\ldots 1'3'5'\ldots)$.
\end{enumerate}
\end{proposition}

\begin{proof}
According to the definition of the operators $T_\pi$, the operator associated to the partition in the statement is given by the following formula:
$$T_\pi(e_{a_1}\otimes e_{b_1}\otimes e_{a_2}\otimes e_{b_2}\otimes\ldots)=\delta(a)\delta(b)e_b\otimes e_a\otimes e_b\otimes e_a\otimes\ldots$$

Here we use the convention that $\delta(a)=1$ if all the indices $a_i$ are equal, and $\delta(a)=0$ otherwise, along with a similar convention for $\delta(b)$. As for the indices $a,b$ appearing on the right, these are the common values of the $a$ indices and $b$ indices, respectively, in the case $\delta(a)=\delta(b)=1$, and are irrelevant quantities in the remaining cases. Now with the above formula of $T_\pi$ in hand, we have the following computation, for any $u=(u_{ij})$:
$$T_\pi u^{\otimes s}(e_{a_1}\otimes e_{b_1}\otimes e_{a_2}\otimes\ldots)=\sum_{ij}e_i\otimes e_j\otimes e_i\otimes\ldots\otimes u_{ia_1}u_{jb_1}u_{ia_2}\ldots$$

Here the sum is over all indices $i,j$. Similarly, we have the following computation, with the sum being this time over all multi-indices $i=(i_1,\ldots,i_s)$, $j=(j_1,\ldots,j_s)$:
$$u^{\otimes s}T_\pi(e_{a_1}\otimes e_{b_1}\otimes e_{a_2}\otimes\ldots)=\delta(a)\delta(b)\sum_{ij}e_{i_1}\otimes e_{j_1}\otimes e_{i_2}\otimes\ldots\otimes u_{i_1b}u_{j_1a}u_{i_2b}\ldots$$

But with these formulae in hand, the identification of the right terms, after a suitable relabeling of indices, gives the equivalence in the statement.
\end{proof}

Still following \cite{ez1}, we have the following result:

\index{s-balanced partitions}

\begin{theorem}
$H_N^{(s)}$ is an easy quantum group, and its associated category $P_{even}^s$ is that of the ``$s$-balanced'' partitions, i.e. partitions satisfying the following conditions:
\begin{enumerate}
\item The total number of legs is even.

\item In each block, the number of odd legs equals the number of even legs, modulo $s$.
\end{enumerate}
\end{theorem}

\begin{proof}
This is something standard, which can be proved as follows:

\medskip

(1) As a first remark, at $s=2$ the first condition implies the second one, so here we simply get the partitions having an even number of legs, corresponding to $H_N$. Observe also that at $s=\infty$ we get the partitions which are balanced, which correspond to the quantum group $H_N^*$. Thus, we have indeed the result at the endpoints, $s=2,\infty$.

\medskip

(2) Our first claim is that $P_{even}^s$ is indeed a category. But this follows by adapting the $s=\infty$ argument in the proof for $H_N^*$, just by adding ``modulo $s$'' everywhere.

\medskip

(3) It remains to prove that this category corresponds indeed to $H_N^{(s)}$. But this follows from the fact that the partition $\pi$ appearing in Proposition 9.4 generates the category of $s$-balanced partitions, as one can check by a routine computation.
\end{proof}

Observe in particular, coming as a consequence of Theorem 9.5, the fact that the quantum groups $H_N^{(s)}$ are indeed distinct. We will see in a moment, in Theorem 9.6 below, another proof of this fact, which is even nicer, and more intuitive.

\bigskip

As another result now, making us exit the real world, consider the complex reflection group $H_N^s=\mathbb Z_s\wr S_N$, consisting of the monomial matrices having the $s$-roots of unity as nonzero entries. Observe that we have $PH_N^{(s)}=H_N^s/\mathbb T$. We have the following result:

\begin{theorem}
We have an isomorphism of projective quantum groups
$$PH_N^{(s)}=PH_N^s$$
walid for any $s\in\{2,3,\ldots,\infty\}$.
\end{theorem}

\begin{proof}
Observe first that this statement holds indeed at $s=2$, because here we have $H_N^{(2)}=H_N^2=H_N$. This statement holds as well at $s=\infty$. In the general case, observe first that from $H_N^{(s)}\subset H_N^*$ we get $PH_N^{(s)}\subset PH_N^*=PK_N$, so $PH_N^{(s)}$ is indeed a classical group. In order to compute this group, consider the following diagram:
$$\begin{matrix}
H_N^s&\subset&U_N^+\\
\\
\cup&&\cup\\
\\
S_N&\subset&H_N^{(s)}
\end{matrix}$$

The corresponding sets of partitions are then as follows:
$$\begin{matrix}
P_{even}^s(2k,2l)&\supset&\mathcal P_2(2k,2l)\\ 
\\ 
\cup&&\cup\\
\\
P(2k,2l)&\supset&P_{even}^s(2k,2l)
\end{matrix}$$

Let us look now at the projective versions of the above quantum groups:
$$\begin{matrix}
PH_N^s&\subset&PU_N^+\\
\\
\cup&&\cup\\
\\
PH_N&\subset&PH_N^{(s)}
\end{matrix}$$

As before for $H_N^*$, we are in the situation where we have two quantum subgroups having the same diagrams, and we conclude that we have $PH_N^{(s)}=PH_N^s$.
\end{proof}

There are many other things that can be said about $H_N^{(s)}$, of more specialized nature, and we will be back to this, once we will have more tools for studying such quantum groups. The idea indeed is that the above results, which are quite straightforward, exclusively based on easiness, can be complemented by some very concrete results as well, for instance in connection with semidirect products, following \cite{rw1}, \cite{rw2}. But more on this later. 

\bigskip

Moving ahead now, still following the old paper \cite{ez1}, we can introduce as well a second one-parameter series of hyperoctahedral quantum groups, $H_N^{[s]}$ with $s\in\{2,3,\ldots,\infty\}$, again having as main particular case the group $H_N^{[2]}=H_N$, as follows:

\index{ultracubic relations}

\begin{definition}
$C(H_N^{[s]})$ is the universal $C^*$-algebra generated by $N^2$ self-adjoint variables $u_{ij}$, subject to the following relations:
\begin{enumerate}
\item Orthogonality: $uu^t=u^tu=1$, where $u=(u_{ij})$ and $u^t=(u_{ji})$.

\item Ultracubic relations: $acb=0$, for any $a\neq b$ on the same row or column of $u$.

\item $s$-mixing relation: $abab\ldots=baba\ldots$ (length $s$ words), for any $a,b\in\{u_{ij}\}$.
\end{enumerate}
\end{definition}

Our first task is to compare the defining relations for $H_N^{[s]}$ with those for $H_N^{(s)}$. In order to deal at the same time with the cubic and ultracubic relations, it is convenient to use a statement regarding a certain unifying notion, of ``$k$-cubic'' relations:

\begin{proposition}
For a closed subgroup $G\subset O_N^+$, the following are equivalent:
\begin{enumerate}
\item The basic coordinates $u_{ij}$ satisfy the $k$-cubic relations, namely 
$$ac_1\ldots c_kb=0$$
for any $a\neq b$ on the same row or column of $u$, and for any $c_1,\ldots,c_k$.

\item We have $T_\pi\in End(u^{\otimes k+2})$, where $\pi$ is the following partition,
$$\pi=(1,1',k+2,k+2')(2,2')\ldots (k+1,k+1')$$
and where $u=(u_{ij})$ is as usual the fundamental corepresentation.
\end{enumerate}
\end{proposition}

\begin{proof}
According to the definition of the operators $T_\pi$, the operator associated to the partition in the statement is given by the following formula:
$$T_\pi(e_a\otimes e_{c_1}\otimes\ldots\otimes e_{c_k}\otimes e_b)=\delta_{ab}e_a\otimes e_{c_1}\otimes\ldots\otimes e_{c_k}\otimes e_a$$

But this gives the following formula, for any $u=(u_{ij})$, with sum over all indices $i,l$:
\begin{eqnarray*}
&&T_\pi u^{\otimes k+2}(e_a\otimes e_{c_1}\otimes\ldots\otimes e_{c_k}\otimes e_b)\\
&=&\sum_{ij}e_i\otimes e_{j_1}\otimes...\otimes e_{j_k}\otimes e_i\otimes u_{ia}u_{j_1c_1}\ldots u_{j_kc_k}u_{ib}
\end{eqnarray*}

Similarly, we have the following formula, again for any $u=(u_{ij})$, and with the sum being this time over all multi-indices $j=(j_1,\ldots,j_k)$: 
\begin{eqnarray*}
&&u^{\otimes k+2}T_\pi(e_a\otimes e_{c_1}\otimes\ldots\otimes e_{c_k}\otimes e_b)\\
&=&\delta_{ab}\sum_{ijl}e_i\otimes e_{j_1}\otimes...\otimes e_{j_k}\otimes e_l\otimes u_{ia}u_{j_1c_1}\ldots u_{j_kc_k}u_{la}
\end{eqnarray*}

Now the identification of the right terms gives the equivalence in the statement.
\end{proof}

We can now establish the precise relationship between $H_N^{[s]}$ and $H_N^{(s)}$, and also show that no further series can appear in this way, the result being as follows:

\begin{proposition}
For $k\geq 1$ the $k$-cubic relations are all equivalent to the ultracubic relations, and they imply the cubic relations.
\end{proposition}

\begin{proof}
This follows indeed from the following two observations:

\medskip

(1) The $k$-cubic relations imply the $2k$-cubic relations. Indeed, one can connect two copies of the partition $\pi$ in Proposition 9.8, by gluing them with two semicircles in the middle, and the resulting partition is the one implementing the $2k$-cubic relations.

\medskip

(2) The $k$-cubic relations imply the $(k-1)$-cubic relations. Indeed, by capping the partition $\pi$ in Proposition 9.8 with a semicircle at bottom right, we get a certain partition $\pi'\in P(k+2,k)$, and by rotating the upper right leg of this partition we get the partition $\pi''\in P(k+1,k+1)$ implementing the $(k-1)$-cubic relations.
\end{proof}

The above statement shows that replacing in Definition 9.7 the ultracubic condition by any of the $k$-cubic conditions, with $k\geq 2$, won't change the resulting quantum group. The other consequences of Proposition 9.9 can be summarized as follows:

\begin{proposition}
The quantum groups $H_N^{[s]}$ have the following properties:
\begin{enumerate}
\item We have $H_N^{(s)}\subset H_N^{[s]}\subset H_N^+$.

\item At $s=2$ we have $H_N^{[2]}=H_N^{(s)}=H_N$.

\item At $s\geq 3$ we have $H_N^{(s)}\neq H_N^{[s]}$.
\end{enumerate}
\end{proposition}

\begin{proof}
All the assertions basically follow from Proposition 9.9, as follows:

\medskip

(1) For the first inclusion, we need to show that half-commutation + cubic implies ultracubic, and this can be done by placing the half-commutation partition next to the cubic partition, then using 2 semicircle cappings in the middle. The second inclusion follows from Proposition 9.9, because the ultracubic relations (1-cubic relations) imply the cubic relations (0-cubic relations). Thus, we have both inclusions.

\medskip

(2) Observe first that at $s=2$ the $s$-commutation relation is the usual commutation relation $ab=ba$. Thus we are led here to the algebra generated by $N^2$ commuting self-adjoint variables satisfying the cubic condition, which is $C(H_N)$.

\medskip

(3) Finally, $H_N^{(s)}\neq H_N^{[s]}$ will be a consequence of the results below, because at $s\geq 3$ the half-commutation partition $p=(14)(25)(36)$ is $s$-balanced but not locally $s$-balanced.
\end{proof}

Still following \cite{ez1}, we have the following result:

\index{locally s-balanced}

\begin{theorem}
$H_N^{[s]}$ is an easy quantum group, and its associated category is that of the ``locally $s$-balanced'' partitions, i.e. partitions having the property that each of their subpartitions (i.e. partitions obtained by removing certain blocks) are $s$-balanced.
\end{theorem}

\begin{proof}
This is routine from what we have, the idea being as follows:

\medskip

(1) As a first remark, at $s=2$ the locally $s$-balancing condition is automatic for a partition having blocks of even size, so we get indeed the category corresponding to $H_N$. 

\medskip

(2) In the general case now, our first claim is that the locally $s$-balanced partitions from indeed a category of partitions. But this follows simply by adapting the argument in the proof for $H_N^*$, just by adding ``locally'' everywhere. 

\medskip

(3) It remains to prove that this category corresponds indeed to the quantum group $H_N^{[s]}$. But this follows from the fact that the partition generating the category of locally balanced partitions, namely $\pi=(1346)(25)$, is nothing but the one implementing the ultracubic relations, as one can check by a routine computation.
\end{proof}

There are many other things that can be said about the quantum groups $H_N^{(s)}$ and $H_N^{[s]}$, or algebraic and analytic nature. We will be back to all this later, following \cite{rw1}, \cite{rw2}, once we will have more tools for the study of these quantum groups.

\section*{9b. Limiting objects}

The quantum groups $H_N^{(s)}$ and $H_N^{[s]}$ constructed above are all contained in $H_N^{[\infty]}$, and before going further with the construction of more examples, we would like to study in detail this quantum group $H_N^{[\infty]}$. Following as before \cite{rw3}, let us begin with:

\begin{definition}
We let $P_{even}^{[\infty]}$ be the category generated by the partition
$$
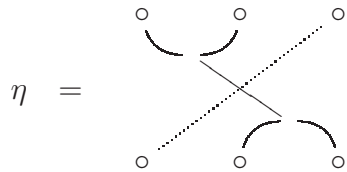
\xymatrix@R=2mm@C=3mm{\\ \\ \eta\ \ =\\ \\}\ \ \ 
\xymatrix@R=2mm@C=3mm{
\circ\ar@/_/@{-}[dr]&&\circ&&\circ\ar@{.}[ddddllll]\\
&\ar@/_/@{-}[ur]\ar@{-}[ddrr]\\
\\
&&&\ar@/^/@{-}[dr]\\
\circ&&\circ\ar@/^/@{-}[ur]&&\circ}$$
and we denote by $H_N^{[\infty]}$ the corresponding easy quantum group $H_N\subset G\subset H_N^+$.
\end{definition}

The partitions $\pi\in P_{even}^{[\infty]}$ can be characterized by the fact that all their subpartitions $\sigma\subset\pi$ belong $P_{even}^*$. That is, the following condition must satisfied:
$$\sigma\subset\pi\implies\sigma\in P_{even}^*$$

As an illustration, let us verify that we have indeed $\eta\in P_{even}^{[\infty]}$.  The standard coloring of $\eta$, with alternating colors, as per the usual $P_{even}^*$ requirements, is as follows:
$$
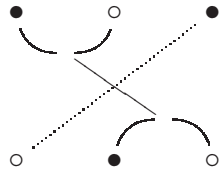
\xymatrix@R=2mm@C=3mm{
\bullet\ar@/_/@{-}[dr]&&\circ&&\bullet\ar@{.}[ddddllll]\\
&\ar@/_/@{-}[ur]\ar@{-}[ddrr]\\
\\
&&&\ar@/^/@{-}[dr]\\
\circ&&\bullet\ar@/^/@{-}[ur]&&\circ}$$

We can see that this partition has then the same number of $\circ,\bullet$ legs. As for the subpartitions, these are as follows, again having the same number of $\circ,\bullet$ legs:
$$
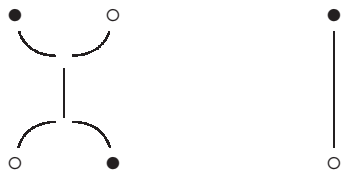
\xymatrix@R=2mm@C=3mm{
\bullet\ar@/_/@{-}[dr]&&\circ\\
&\ar@/_/@{-}[ur]\ar@{-}[dd]\\
\\
&\ar@/^/@{-}[dr]\\
\circ\ar@/^/@{-}[ur]&&\bullet}\qquad\qquad\qquad
\xymatrix@R=15.4mm@C=7mm{\bullet\ar@{-}[d]\\ \circ}
$$

Regarding now the quantum group $H_N^{[\infty]}$, it is known that this contains $H_N^*$, and also that $H_N^{[\infty]}\subset O_N^+$ appears by assuming that the standard coordinates $u_{ij}$ satisfy the relations $abc=0$, for any $a\neq c$ on the same row or column of $u$. In fact, we have:

\begin{theorem}
Let $H_N^{[\infty]}\subset O_N^+$ be the compact quantum group obtained via the relations $abc=0$, whenever $a\neq c$ are on the same row or column of $u$. 
\begin{enumerate}
\item We have inclusions $H_N^*\subset H_N^{[\infty]}\subset H_N^+$.

\item We have $ab_1\ldots b_rc=0$, whenever $a\neq c$ are on the same row or column of $u$.

\item We have $ab^2=b^2a$, for any two entries $a,b$ of $u$.
\end{enumerate}
\end{theorem}

\begin{proof}
We briefly recall the proof in \cite{rw3}, for future use in what follows. Our first claim is that $H_N^{[\infty]}$ comes, as an easy quantum group, from the following diagram:
$$\xymatrix@R=5mm@C=0.1mm{&\\\pi\ \ =\\&}\xymatrix@R=5mm@C=5mm{
\circ\ar@{-}[dd]&\circ\ar@{.}[dd]&\circ\ar@{-}[dd]\\
\ar@{-}[rr]&&\\
\circ&\circ&\circ}$$

Indeed, this diagram acts via the following linear map:
$$T_\pi(e_{ijk})=\delta_{ik}e_{ijk}$$

We therefore have the following formula:
\begin{eqnarray*}
T_\pi u^{\otimes 3}e_{abc}
&=&T_\pi\sum_{ijk}e_{ijk}\otimes u_{ia}u_{jb}u_{kc}\\
&=&\sum_{ijk}e_{ijk}\otimes\delta_{ik}u_{ia}u_{jb}u_{kc}
\end{eqnarray*}

On the other hand, we have as well the following formula:
\begin{eqnarray*}
u^{\otimes 3}T_\pi e_{abc}
&=&u^{\otimes 3}\delta_{ac}e_{abc}\\
&=&\sum_{ijk}e_{ijk}\otimes\delta_{ac}u_{ia}u_{jb}u_{kc}
\end{eqnarray*}

Thus the condition $T_\pi\in End(u^{\otimes 3})$ is equivalent to the following relations:
$$(\delta_{ik}-\delta_{ac})u_{ia}u_{jb}u_{kc}=0$$

The non-trivial cases are $i=k,a\neq c$ and $i\neq k,a=c$, and these produce the following relations between the standard coordinates of our quantum group:

\medskip

-- $u_{ia}u_{jb}u_{ic}=0$ for any $a\neq c$.

\medskip

-- $u_{ia}u_{jb}u_{ka}=0$, for any $i\neq k$. 

\medskip

Thus, we have reached to the standard relations for the quantum group $H_N^{[\infty]}$.

\medskip

(1) We have the following formula:
$$
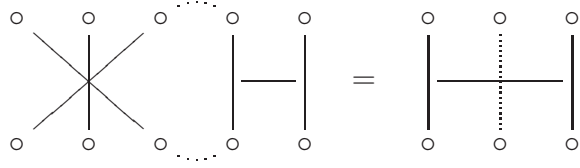
\xymatrix@R=5mm@C=5mm{
\circ\ar@{-}[ddrr]&\circ\ar@{-}[dd]&\circ\ar@{-}[ddll]\ar@/^/@{.}[r]&\circ\ar@{-}[dd]&\circ\ar@{-}[dd]\\
&&&\ar@{-}[r]&\\
\circ&\circ&\circ\ar@/_/@{.}[r]&\circ&\circ}\ \ \xymatrix@R=5mm@C=1mm{&\\=\\&}\xymatrix@R=5mm@C=5mm{
\circ\ar@{-}[dd]&\circ\ar@{.}[dd]&\circ\ar@{-}[dd]\\
\ar@{-}[rr]&&\\
\circ&\circ&\circ}$$

We have as well the following formula:
$$
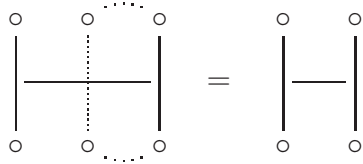
\xymatrix@R=5mm@C=5mm{
\circ\ar@{-}[dd]&\circ\ar@{.}[dd]\ar@/^/@{.}[r]&\circ\ar@{-}[dd]\\
\ar@{-}[rr]&&\\
\circ&\circ\ar@/_/@{.}[r]&\circ}\ \ \xymatrix@R=5mm@C=1mm{&\\=\\&}\xymatrix@R=5mm@C=5mm{
\circ\ar@{-}[dd]&\circ\ar@{-}[dd]\\
\ar@{-}[r]&&\\
\circ&\circ}$$

Thus, we obtain inclusions as desired, namely:
$$H_N^*\subset H_N^{[\infty]}\subset H_N^+$$

(2) At $r=2$, the relations $ab_1b_2c=0$ come indeed from the following diagram:
$$
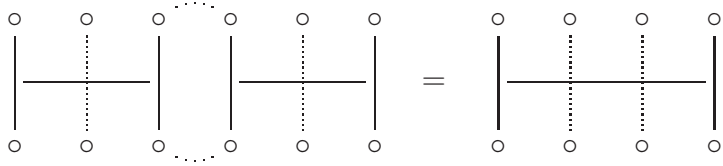
\xymatrix@R=5mm@C=5mm{
\circ\ar@{-}[dd]&\circ\ar@{.}[dd]&\circ\ar@{-}[dd]\ar@/^/@{.}[r]&\circ\ar@{-}[dd]&\circ\ar@{.}[dd]&\circ\ar@{-}[dd]\\
\ar@{-}[rr]&&&\ar@{-}[rr]&&\\
\circ&\circ&\circ\ar@/_/@{.}[r]&\circ&\circ&\circ}
\ \ \xymatrix@R=5mm@C=1mm{&\\=\\&}
\xymatrix@R=5mm@C=5mm{
\circ\ar@{-}[dd]&\circ\ar@{.}[dd]&\circ\ar@{.}[dd]&\circ\ar@{-}[dd]\\
\ar@{-}[rrr]&&&\\
\circ&\circ&\circ&\circ}$$

In the general case $r\geq2$ the proof is similar, see \cite{ez1} for details.

\medskip

(3) We use here an idea from \cite{rw3}. By rotating $\pi$, we obtain:
$$
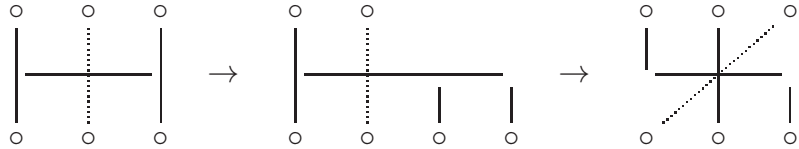
\xymatrix@R=5mm@C=5mm{
\circ\ar@{-}[dd]&\circ\ar@{.}[dd]&\circ\ar@{-}[dd]\\
\ar@{-}[rr]&&\\
\circ&\circ&\circ}
\ \ \xymatrix@R=5mm@C=0.1mm{&\\ \to\\&}\ 
\xymatrix@R=5mm@C=5mm{
\circ\ar@{-}[dd]&\circ\ar@{.}[dd]&\\
\ar@{-}[rrr]&&\ar@{-}[d]&\ar@{-}[d]\\
\circ&\circ&\circ&\circ}
\ \ \xymatrix@R=5mm@C=0.1mm{&\\ \to\\&}\ 
\xymatrix@R=5mm@C=5mm{
\circ\ar@{-}[d]&\circ\ar@{-}[dd]&\circ\ar@{.}[ddll]\\
\ar@{-}[rr]&&\ar@{-}[d]\\
\circ&\circ&\circ}$$

Let us denote by $\sigma$ the partition on the right. Since $T_\sigma(e_{ijk})=\delta_{ij}e_{kji}$, we obtain:
\begin{eqnarray*}
T_\sigma u^{\otimes 3}e_{abc}
&=&T_\sigma\sum_{ijk}e_{ijk}\otimes u_{ia}u_{jb}u_{kc}\\
&=&\sum_{ijk}e_{kji}\otimes\delta_{ij}u_{ia}u_{jb}u_{kc}
\end{eqnarray*}

On the other hand, we obtain as well the following formula:
\begin{eqnarray*}
u^{\otimes 3}T_\sigma e_{abc}
&=&u^{\otimes 3}\delta_{ab}e_{cba}\\
&=&\sum_{ijk}e_{kji}\otimes\delta_{ab}u_{kc}u_{jb}u_{ia}
\end{eqnarray*}

Thus the condition $T_\sigma\in End(u^{\otimes 3})$ is equivalent to the following relations:
$$\delta_{ij}u_{ia}u_{jb}u_{kc}=\delta_{ab}u_{kc}u_{jb}u_{ia}$$

Now by setting $j=i,b=a$ in this formula we obtain the following formula:
$$u_{ia}^2u_{kc}=u_{kc}u_{ia}^2$$

But these are exactly the commutation relations in the statement, as desired.
\end{proof}

In order to discuss some further features of $H_N^{[\infty]}$, we will need some basic twisting theory. We will systematically discuss the twisting in chapter 13, and in what concerns the present chapter, we will only need here the following standard fact:

\begin{proposition}
There is a signature map $\varepsilon:P_{even}\to\{-1,1\}$, given by 
$$\varepsilon(\tau)=(-1)^c$$
where $c$ is the number of switches needed to make $\tau$ noncrossing. In addition:
\begin{enumerate}
\item For $\tau\in S_k$, this is the usual signature.

\item For $\tau\in P_2$ we have $(-1)^c$, where $c$ is the number of crossings.

\item For $\tau\leq\pi\in NC_{even}$, the signature is $1$.
\end{enumerate}
\end{proposition}

\begin{proof}
In order to show that the signature map $\varepsilon:P_{even}\to\{-1,1\}$ in the statement, given by $\varepsilon(\tau)=(-1)^c$, is well-defined, we must prove that the number $c$ in the statement is well-defined modulo 2. It is enough to perform the verification for the noncrossing partitions. More precisely, given $\tau,\tau'\in NC_{even}$ having the same block structure, we must prove that the number of switches $c$ required for the passage $\tau\to\tau'$ is even.

\medskip

In order to do so, observe that any partition $\tau\in P(k,l)$ can be put in ``standard form'', by ordering its blocks according to the appearence of the first leg in each block, counting clockwise from top left, and then by performing the switches as for block 1 to be at left, then for block 2 to be at left, and so on.

\medskip

The point now is that, under the assumption $\tau\in NC_{even}(k,l)$, each of the moves required for putting a leg at left, and hence for putting a whole block at left, requires an even number of switches. Thus, putting $\tau$ is standard form requires an even number of switches. Now given $\tau,\tau'\in NC_{even}$ having the same block structure, the standard form coincides, so the number of switches $c$ required for the passage $\tau\to\tau'$ is indeed even.

\medskip

Regarding now the remaining assertions, these are all elementary:

\medskip

(1) For $\tau\in S_k$ the standard form is $\tau'=id$, and the passage $\tau\to id$ comes by composing with a number of transpositions, which gives the signature. 

\medskip

(2) For a general $\tau\in P_2$, the standard form is of type $\tau'=|\ldots|^{\cup\ldots\cup}_{\cap\ldots\cap}$, and the passage $\tau\to\tau'$ requires $c$ mod 2 switches, where $c$ is the number of crossings. 

\medskip

(3) Assuming that $\tau\in P_{even}$ comes from $\pi\in NC_{even}$ by merging a certain number of blocks, we can prove that the signature is 1 by proceeding by recurrence.
\end{proof}

Getting back now to $H_N^{[\infty]}$, we have the following useful result, regarding it:

\begin{theorem}
We have the following equalities,
\begin{eqnarray*}
P_{even}^*&=&\left\{\pi\in P_{even}\Big|\varepsilon(\tau)=1,\forall\tau\leq\pi,|\tau|=2\right\}
\\
P_{even}^{[\infty]}&=&\left\{\pi\in P_{even}\Big|\sigma\in P_{even}^*,\forall\sigma\subset\pi\right\}\\
P_{even}^{[\infty]}&=&\left\{\pi\in P_{even}\Big|\varepsilon(\tau)=1,\forall\tau\leq\pi\right\}
\end{eqnarray*}
where $\varepsilon:P_{even}\to\{\pm1\}$ is the signature of even permutations.
\end{theorem}

\begin{proof}
This is routine combinatorics, from \cite{ba1}, \cite{rw3}, the idea being as follows:

\medskip

(1) Given $\pi\in P_{even}$, we have $\tau\leq\pi,|\tau|=2$ precisely when $\tau=\pi^\beta$ is the partition obtained from $\pi$ by merging all the legs of a certain subpartition $\beta\subset\pi$, and by merging as well all the other blocks. Now observe that $\pi^\beta$ does not depend on $\pi$, but only on $\beta$, and that the number of switches required for making $\pi^\beta$ noncrossing is $c=N_\bullet-N_\circ$ modulo 2, where $N_\bullet/N_\circ$ is the number of black/white legs of $\beta$, when labelling the legs of $\pi$ counterclockwise $\circ\bullet\circ\bullet\ldots$ Thus $\varepsilon(\pi^\beta)=1$ holds precisely when $\beta\in\pi$ has the same number of black and white legs, and this gives the result.

\medskip

(2) This simply follows from the equality $P_{even}^{[\infty]}=<\eta>$ coming from Theorem 9.13, by computing $<\eta>$, and for the complete proof here we refer to \cite{rw3}.

\medskip

(3) We use the fact, also from \cite{rw3}, that the relations $g_ig_ig_j=g_jg_ig_i$ are trivially satisfied for real reflections. We conclude from this that we have:
$$P_{even}^{[\infty]}(k,l)=\left\{\ker\begin{pmatrix}i_1&\ldots&i_k\\ j_1&\ldots&j_l\end{pmatrix}\Big|g_{i_1}\ldots g_{i_k}=g_{j_1}\ldots g_{j_l}\ {\rm inside}\ \mathbb Z_2^{*N}\right\}$$

Thus, the partitions in $P_{even}^{[\infty]}$ are those describing the relations between variables subject to the conditions $g_i^2=1$. We conclude that $P_{even}^{[\infty]}$ appears from $NC_{even}$ by ``inflating blocks'', in the sense that each $\pi\in P_{even}^{[\infty]}$ can be transformed into a partition $\pi'\in NC_{even}$ by deleting pairs of consecutive legs, belonging to the same block. Now since this operation leaves invariant modulo 2 the number $c\in\mathbb N$ of switches in the definition of the signature, it leaves invariant the signature $\varepsilon=(-1)^c$ itself, and we obtain the inclusion ``$\subset$''. Conversely, given $\pi\in P_{even}$ satisfying $\varepsilon(\tau)=1$, $\forall\tau\leq\pi$, our claim is that:
$$\rho\leq\sigma\subset\pi,|\rho|=2\implies\varepsilon(\rho)=1$$

Indeed, let us denote by $\alpha,\beta$ the two blocks of $\rho$, and by $\gamma$ the remaining blocks of $\pi$, merged altogether. We know that the partitions $\tau_1=(\alpha\wedge\gamma,\beta)$, $\tau_2=(\beta\wedge\gamma,\alpha)$, $\tau_3=(\alpha,\beta,\gamma)$ are all even. On the other hand, putting these partitions in noncrossing form requires respectively $s+t,s'+t,s+s'+t$ switches, where $t$ is the number of switches needed for putting $\rho=(\alpha,\beta)$ in noncrossing form. Thus $t$ is even, and we are done. With the above claim in hand, we conclude, by using the second equality in the statement, that we have $\sigma\in P_{even}^*$. Thus we have $\pi\in P_{even}^{[\infty]}$, which ends the proof of ``$\supset$''.
\end{proof}

There are many other things that can be said about the quantum groups $H_N^{[\infty]}$ introduced above, both at the algebraic and probabilistic level. We will be back to this.

\section*{9c. Varieties of groups}

Following the papers of Raum-Weber \cite{rw1}, \cite{rw2}, \cite{rw3}, we will extend now the construction of the series $H_N^{(s)}$ and $H_N^{[s]}$, that we have so far, into something very general, which will cover in fact all the intermediate easy quantum groups, as follows:
$$H_N\subset G\subset H_N^{[\infty]}$$ 

This will be something quite tricky, with some delicate group theory and algebra involved. Following \cite{rw1}, \cite{rw2}, \cite{rw3}, let us start with the following definition:

\begin{definition}
We call real reflection group any finitely generated discrete group, whose generators are real reflections, in the sense that they square up to $1$:
$$\Gamma=<g_1,\ldots,g_N>\quad,\quad g_i^2=1$$
Such a real reflection group is called uniform if each permutation $\sigma\in S_N$ produces a group automorphism, given on generators by the following formula:
$$g_i\to g_{\sigma(i)}$$
Also, we say that $\Gamma$ is non-degenerate when its biggest abelian quotient, obtained by imposing commutation relations to all the generators $g_i$, equals the group $\mathbb Z_2^N$.
\end{definition}

There are many things that can be said, about these conditions. As a first observation, having a real reflection group as above is the same as having a quotient as follows:
$$\mathbb Z_2^{*N}\to\Gamma$$

When assuming in addition that $\Gamma$ is non-degenerate, in the above sense, we can see that $\Gamma$ must appear as an intermediate discrete group, as follows:
$$\mathbb Z_2^{*N}\to\Gamma\to\mathbb Z_2^N$$

It is quite useful at this point to look as well at the dual of $\Gamma$, from this perspective. If we denote as usual by $T_N^+$ the free real torus, appearing as dual of $\mathbb Z_2^{*N}$, having a real reflection group as above is the same as having a compact quantum group as follows:
$$\widehat{\Gamma}\subset T_N^+$$

Moreover, assuming in addition that $\Gamma$ is non-degenerate, in the above sense, we can see that $\widehat{\Gamma}$ must appear as an intermediate compact quantum group, as follows:
$$T_N\subset\widehat{\Gamma}\subset T_N^+$$

At the level of examples of such groups, we have 3 examples to be always kept in mind, coming from our usual classical / half-classical / free philosophy, namely:

\bigskip

(1) The group $\mathbb Z_2^N$, usual product of $\mathbb Z_2$ with itself.

\bigskip

(2) The group $\mathbb Z_2^{\circ N}$, with $\circ$ standing for free product up to $abc=cba$.

\bigskip

(3) The group $\mathbb Z_2^{*N}$, usual free product of $\mathbb Z_2$ with itself.

\bigskip

Following Raum and Weber \cite{rw1}, \cite{rw2}, \cite{rw3}, we will prove now that the easy quantum groups $H_N\subset G\subset H_N^{[\infty]}$ are in correspondence with the uniform real reflection groups $\mathbb Z_2^{*\infty}\to\Gamma\to\mathbb Z_2^\infty$, with the main instances of this correspondence being as follows:
$$\xymatrix@R=17mm@C=17mm{
\mathbb Z_2^N\ar@{~}[d]&\mathbb Z_2^{\circ N}\ar[l]\ar@{~}[d]&\mathbb Z_2^{*N}\ar[l]\ar@{~}[d]\\
H_N\ar[r]&H_N^*\ar[r]&H_N^{[\infty]}}$$

This will be something quite tricky, with the correspondence involving as well the corresponding categories of partitions, which will be as follows:
$$P_{even}^{[\infty]}\subset D\subset P_{even}$$

Getting started now, as a first result, which is something of purely group-theoretical nature, and without many assumptions on the groups involved, we have:

\begin{proposition}
Given a real reflection group $\mathbb Z_2^{*N}\to\Gamma$, the following family of subsets $D(k,l)\subset P(k,l)$ is a category of partitions
$$D(k,l)=\left\{\pi\in P(k,l)\Big|\ker\binom{i}{j}\leq\pi\implies g_{i_1}\ldots g_{i_k}=g_{j_1}\ldots g_{j_l}\right\}$$
satisfying $P_{even}^{[\infty]}\subset D\subset P$. Moreover, this category appears as
$$P_{even}^{[\infty]}\subset D\subset P_{even}$$
when assuming that our real reflection group $\Gamma$ is non-degenerate.
\end{proposition}

\begin{proof}
There are many things to be checked here, namely the 5 axioms for the categories of partitions, and then the inclusions regarding $D$, the idea being as follows:

\medskip

(1) Composition. We must prove here that the following happens:
$$\pi,\sigma\in D\implies\begin{bmatrix}\pi\\ \sigma\end{bmatrix}\in D$$

So, assume $\pi,\sigma\in D$, which amounts in saying that we have:
$$\ker\binom{i}{j}\leq\pi\implies g_{i_1}\ldots g_{i_k}=g_{j_1}\ldots g_{j_l}$$
$$\ker\binom{j}{s}\leq\pi\implies g_{j_1}\ldots g_{j_l}=g_{s_1}\ldots g_{s_m}$$

In order to prove our result, we must prove that the following happens:
$$\ker\binom{i}{s}\leq\begin{bmatrix}\pi\\ \sigma\end{bmatrix}
\implies g_{i_1}\ldots g_{i_k}=g_{s_1}\ldots g_{s_m}$$

But this is clear from our assumptions $\pi,\sigma\in D$ above. Indeed, since the condition on the left, namely  $\ker(^i_s)\leq[^\pi_\sigma]$, tells us that the indices fit, we can come up with a middle index $j$ which fits with both $i,s$, and we get our result, coming from:
$$g_{i_1}\ldots g_{i_k}=g_{j_1}\ldots g_{j_l}=g_{s_1}\ldots g_{s_m}$$

(2) Tensor products. We must prove here that the following happens:
$$\pi,\sigma\in D\implies[\pi\ \sigma]\in D$$

So, assume $\pi,\sigma\in D$, which amounts in saying that we have:
$$\ker\binom{i}{j}\leq\pi\implies g_{i_1}\ldots g_{i_k}=g_{j_1}\ldots g_{j_l}$$
$$\ker\binom{s}{t}\leq\pi\implies g_{s_1}\ldots g_{s_m}=g_{t_1}\ldots g_{s_n}$$

In order to prove our result, we must prove that the following happens:
$$\ker\binom{p}{q}\leq[\pi\ \sigma]
\implies g_{p_1}\ldots g_{p_{k+m}}=g_{q_1}\ldots g_{q_{l+n}}$$

But this is clear from our assumptions $\pi,\sigma\in D$ above. Indeed, since the condition on the left, namely  $\ker(^p_q)\leq[\pi \ \sigma]$, tells us that the indices fit, we can split each our indices as $p=(is)$ and $q=(jt)$, and we get our result, coming from:
\begin{eqnarray*}
g_{p_1}\ldots g_{p_{k+m}}
&=&g_{i_1}\ldots g_{i_k}\cdot g_{s_1}\ldots g_{s_m}\\
&=&g_{j_1}\ldots g_{j_l}\cdot g_{t_1}\ldots g_{s_n}\\
&=&g_{q_1}\ldots g_{q_{l+n}}
\end{eqnarray*}

(3) Conjugation. We must prove here that the following happens:
$$\pi,\sigma\in D\implies\pi^*\in D$$

So, assume $\pi,\sigma\in D$, which amounts in saying that we have:
$$\ker\binom{i}{j}\leq\pi\implies g_{i_1}\ldots g_{i_k}=g_{j_1}\ldots g_{j_l}$$

In order to prove our result, we must prove that the following happens:
$$\ker\binom{p}{q}\leq\pi^*
\implies g_{p_1}\ldots g_{p_l}=g_{q_1}\ldots g_{q_k}$$

But this is clear from our assumption $\pi\in D$ above. Indeed, the condition on the left, namely $\ker(^p_q)\leq\pi^*$, tells us that the indices fit, and we deduce from this by upside-down turning that we must have $\ker(^q_p)\leq\pi$, and we get our result, coming from:
\begin{eqnarray*}
g_{p_1}\ldots g_{p_l}
&=&g_{j_1}\ldots g_{j_l}\\
&=&g_{i_1}\ldots g_{i_k}\\
&=&g_{q_1}\ldots g_{q_k}
\end{eqnarray*}

(4) Unit. We must prove here that the following happens:
$$\ker\binom{i}{j}\leq 1_k\implies g_{i_1}\ldots g_{i_k}=g_{j_1}\ldots g_{j_k}$$

But this is clear, because the condition on the left, namely $\ker(^i_j)\leq 1_k$, tells us that our indices must be equal, $i=j$, and so the equality to be proved is trivial.

\medskip

(5) Semicircle. Here we must prove that the following happens:
$$\ker(ab)\leq\cap\implies g_ag_b=1$$

But this is clear, because the condition on the left, namely $\ker(ab)\leq\cap$, tells us that our indices must be equal, $a=b$, and the formula to be proved becomes $g_ag_a=1$, which is true, due to our assumption that $\Gamma$ is a real reflection group.

\medskip

(6) Inclusion $P_{even}^{[\infty]}\subset D$. In order to prove this inclusion, consider the following partition, that we alredy met in the above, when investigating the category $P_{even}^{[\infty]}$:
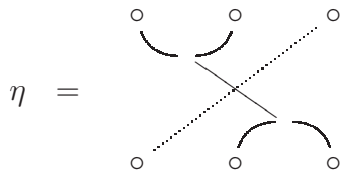
$$\xymatrix@R=2mm@C=3mm{\\ \\ \eta\ \ =\\ \\}\ \ \ 
\xymatrix@R=2mm@C=3mm{
\circ\ar@/_/@{-}[dr]&&\circ&&\circ\ar@{.}[ddddllll]\\
&\ar@/_/@{-}[ur]\ar@{-}[ddrr]\\
\\
&&&\ar@/^/@{-}[dr]\\
\circ&&\circ\ar@/^/@{-}[ur]&&\circ}$$

We have then $\eta\in D$, coming from the following formula for group elements:
$$g_ag_ag_b=g_bg_ag_a$$

Thus we obtain our inclusion of categories, in the following way:
$$<\eta>=P_{even}^{[\infty]}\subset D$$

(7) Inclusion $D\subset P_{even}$, assuming $\Gamma\to\mathbb Z_2^N$. In order to prove this, assume by contradiction that we have $D\not\subset P_{even}$. But this means that we have at least one partition $\pi\in D$ having some blocks of odd length, and by capping with semicircles we conclude that either the singleton, or the double singleton, must be in $D$:
$$|\in D\quad,\quad ||\in D$$

But in both cases the double singleton must be in $D$:
$$||\in D$$

Our claim now is that this is contradictory. Indeed, this condition tells us that we must have $g_ag_b=1$, for any indices $a,b$, and by using now our non-degeneracy assumption $\Gamma\to\mathbb Z_2^N$, or rather its consequence $\Gamma\neq\{1\}$, we obtain our contradiction, as desired.
\end{proof}

The above result is very nice. At the level of main examples, the basic groups, taken in an $N>>0$ sense, produce the following categories of partitions:
$$\xymatrix@R=17mm@C=24mm{
\mathbb Z_2^N\ar@{~}[d]&\mathbb Z_2^{\circ N}\ar[l]\ar@{~}[d]&\mathbb Z_2^{*N}\ar[l]\ar@{~}[d]\\
P_{even}&P_{even}^*\ar[l]&P_{even}^{[\infty]}\ar[l]}$$

More generally, for any $s\in\{2,4,\ldots,\infty\}$, the categories of partitions for the quantum groups $H_N^{(s)}\subset H_N^{[s]}$ come from the quotients of $\mathbb Z_2^{\circ N}\leftarrow\mathbb Z_2^{*N}$ by the relations $(ab)^s=1$:
$$\xymatrix@R=17mm@C=15mm{
\mathbb Z_2^N\ar@{~}[d]&\mathbb Z_2^{\circ N}/<(ab)^s=1>\ar[l]\ar@{~}[d]&\mathbb Z_2^{*N}/<(ab)^s=1>\ar[l]\ar@{~}[d]\\
P_{even}\ar[r]&P_{even}^{(s)}\ar[r]&P_{even}^{[s]}}$$

Conversely now, we have the following result, also from \cite{rw1}, \cite{rw2}, \cite{rw3}:

\begin{proposition}
Given an intermediate category of partitions 
$$P_{even}^{[\infty]}\subset D\subset P_{even}$$
we can associate to it a discrete group, as follows,
$$\Gamma=\left\langle g_1,\ldots g_N\Big|g_{i_1}\ldots g_{i_k}=g_{j_1}\ldots g_{j_l},\forall i,j,k,l,\ker\binom{i}{j}\in D(k,l)\right\rangle$$
which is a uniform reflection group $\mathbb Z_2^{*N}\to\Gamma\to\mathbb Z_2^N$.
\end{proposition}

\begin{proof}
Again, many things to be checked here, the idea being as follows:

\medskip

(1) First of all, the construction in the statement produces indeed a group, with the verification of the group axioms coming from the same computations as those in the beginning of the proof of Proposition 9.17, read in the opposite sense.

\medskip

(2) The fact that our group is indeed a real reflection group, $\mathbb Z_2^{*N}\to\Gamma$, comes from $\cap\in D$. Indeed, this condition tells us that we must have $g_ag_a=1$, as desired.

\medskip

(3) The fact that we group $\Gamma$ that we obtain is indeed uniform is standard too, coming from the fact that our axioms for categories of partitions are permutation-invariant.

\medskip

(4) Finally, the fact the we have $\Gamma\to\mathbb Z_2^N$ can be seen by functoriality. Indeed, for $D=P_{even}$ we obtain $\Gamma=\mathbb Z_2^N$, and so from $D\subset P_{even}$ we get $\Gamma\to\mathbb Z_2^N$, as desired.
\end{proof}

Our claim now is that the correspondences in Proposition 9.17 and Proposition 9.18 are inverse to each other. To be more precise, as explained in \cite{rw2}, the correspondences $\Gamma\to D$ and $D\to\Gamma$ are bijective, and inverse to each other, at $N=\infty$:

\begin{proposition}
The above correspondences are one-to-one between:
\begin{enumerate}
\item Uniform reflection groups $\mathbb Z_2^{*\infty}\to\Gamma\to\mathbb Z_2^\infty$.

\item Categories of partitions $P_{even}^{[\infty]}\subset D\subset P_{even}$.
\end{enumerate}
\end{proposition}

\begin{proof}
This is something quite routine from what we have, and for details here, and for further interpretations of all this, using some abstract algebra, we refer to \cite{rw3}. 
\end{proof}

Let us recall now, from the discussion following Proposition 9.17, that all our examples of easy quantum groups $H_N\subset G\subset H_N^+$, expect for the quantum group $H_N^+$ itself, come from certain categories of partitions which are covered by the correspondence in Proposition 9.19. This suggests to fine-tune the correspondence in Proposition 9.19, by adding to the picture the corresponding easy quantum groups $H_N\subset G\subset H_N^+$ as well, and we have here the following remarkable result, from \cite{rw1}, \cite{rw2}, \cite{rw3}:

\index{variety of groups}

\begin{theorem}
We have correspondences between:
\begin{enumerate}
\item Uniform reflection groups $\mathbb Z_2^{*\infty}\to\Gamma\to\mathbb Z_2^\infty$.

\item Categories of partitions $P_{even}^{[\infty]}\subset D\subset P_{even}$.

\item Easy quantum groups $G=(G_N)$, with $H_N^{[\infty]}\supset G_N\supset H_N$.
\end{enumerate}
\end{theorem}

\begin{proof}
This is something which is quite clear from what we have, and for details here, and for some further interpretations of all this, using abstract algebra, we refer to \cite{rw3}. As an illustration, as mentioned above, we have the following correspondences:
$$\xymatrix@R=17mm@C=27mm{
\mathbb Z_2^N\ar@{~}[d]&\mathbb Z_2^{\circ N}\ar[l]\ar@{~}[d]&\mathbb Z_2^{*N}\ar[l]\ar@{~}[d]\\
H_N\ar[r]&H_N^*\ar[r]&H_N^{[\infty]}}$$

More generally, for any $s\in\{2,4,\ldots,\infty\}$, the quantum groups $H_N^{(s)}\subset H_N^{[s]}$ constructed in \cite{ez1} come from the quotients of $\mathbb Z_2^{\circ N}\leftarrow\mathbb Z_2^{*N}$ by the relations $(ab)^s=1$:
$$\xymatrix@R=17mm@C=17mm{
\mathbb Z_2^N\ar@{~}[d]&\mathbb Z_2^{\circ N}/<(ab)^s=1>\ar[l]\ar@{~}[d]&\mathbb Z_2^{*N}/<(ab)^s=1>\ar[l]\ar@{~}[d]\\
H_N\ar[r]&H_N^{(s)}\ar[r]&H_N^{[s]}}$$

For details on all this, and more, we refer to \cite{rw3}.
\end{proof}

As before with other examples of new quantum reflection groups, there are many other things that can be said, both of algebraic and probabilistic nature.

\section*{9d. Crossed products}

With the results that we have so far, the classification of the easy quantum groups $H_N\subset G\subset H_N^+$ is not over yet, for instance because the correspondence in Theorem 9.20 does not cover the quantum group $H_N^+$ itself. Thus, there is still construction and classification work to be done, and we will be of course back to this, in due time.

\bigskip

Before this, however, let us further examine what we have, from a purely algebraic viewpoint. Since we have decomposition results $H_N=\mathbb Z_2\wr S_N$ and $H_N^+=\mathbb Z_2\wr_*S_N^+$, our classification problem for the intermediate easy quantum groups $H_N\subset G\subset H_N^+$ amounts in classifying the intermediate easy quntum groups, as follows:
$$\mathbb Z_2\wr S_N\subset G\subset \mathbb Z_2\wr_*S_N^+$$

Thus, we can expect the solutions $G$ to appear as some kind of crossed products. Generally speaking, however, this is not exactly true, and we will further comment on this later on, once having the complete list of such intermediate easy quantum groups. 

\bigskip

However, one thing that we can do is to try to work out such decomposition results, for the solutions $G$ that we already have, which are those of the form $G=H_N^\Gamma$, found in the previous section. We are led in this way to the following question:

\begin{question}
Do the quantum groups $G=H_N^\Gamma$ that we found decompose as crossed products, in analogy with the decomposition result $H_N=\mathbb Z_2\wr S_N$?
\end{question}

In order to answer this question, let us go back to the main examples that we have. In what regards the very basic examples, the study here is quite elementary, and we conclude that we do have crossed product decomposition results, as follows:
$$\xymatrix@R=17mm@C=24mm{
\mathbb Z_2^N\ar@{~}[d]&\mathbb Z_2^{\circ N}\ar[l]\ar@{~}[d]&\mathbb Z_2^{*N}\ar[l]\ar@{~}[d]\\
\mathbb Z_2^N\rtimes S_N\ar[r]&\mathbb Z_2^{\circ N}\rtimes S_N\ar[r]&\mathbb Z_2^{*N}\rtimes S_N}$$

To be more precise, the decomposition on the left, $H_N=\mathbb Z_2^N\rtimes S_N$, is the same thing as the usual writing $H_N=\mathbb Z_2\wr S_N$. As for the decompositions in the middle and on the right, $H_N^*=\mathbb Z_2^{\circ N}\rtimes S_N$ and $H_N^{[\infty]}=\mathbb Z_2^{*N}\rtimes S_N$, these can be worked out by using the numerous explicit descriptions of the quantum groups $H_N^*$ and $H_N^{[\infty]}$, found above.

\bigskip

More generally now, a similar study shows that in fact for any $s\in\{2,4,\ldots,\infty\}$, the quantum groups $H_N^{(s)}\subset H_N^{[s]}$ constructed in \cite{ez1} come from the quotients of $\mathbb Z_2^{\circ N}\leftarrow\mathbb Z_2^{*N}$ by the relations $(ab)^s=1$, via a crossed product operation, as follows:
$$\xymatrix@R=17mm@C=10mm{
\mathbb Z_2^N\ar@{~}[d]&\mathbb Z_2^{\circ N}/<(ab)^s=1>\ar[l]\ar@{~}[d]&\mathbb Z_2^{*N}/<(ab)^s=1>\ar[l]\ar@{~}[d]\\
\mathbb Z_2^N\rtimes S_N\ar[r]&\mathbb Z_2^{\circ N}/<(ab)^s=1>\rtimes S_N\ar[r]&\mathbb Z_2^{*N}/<(ab)^s=1>\rtimes S_N}$$

Summarizing, we have here some good evidence towards a ``yes'' answer to Question 9.21, and with the crossed product decomposition being every time something very simple, involving the diagonal torus of $H_N^\Gamma$. And the point now is that we have indeed such a ``yes'' answer, in general, as shown by the following result, from \cite{rw1}, \cite{rw2}, \cite{rw3}:

\begin{theorem}
We have a decomposition result of type
$$H_N^\Gamma=\Gamma\rtimes S_N$$
with $\Gamma$, or rather its dual, being viewed as diagonal torus of $H_N^\Gamma$.
\end{theorem}

\begin{proof}
This is something that we know to hold, from the above, for all the main examples of easy quantum groups of type $H_N^\Gamma$, and the proof in general is similar, basically coming from what we have, and with some care in defining the relevant $\rtimes$ operation for the quantum groups involved. For full details on all this, we refer to \cite{rw1}, \cite{rw2}, \cite{rw3}.
\end{proof}

As a conclusion to all this, with respect to our objectives formulated in the beginning of this chapter, we have solved half of the classification problem for the easy quantum groups $H_N\subset G\subset H_N^+$. We will be back to the other half, in the next chapter.

\section*{9e. Exercises} 

The material in this chapter has been quite exciting, and suprising too, at the first reading, and we have several interesting exercises about all this. First, we have:

\begin{exercise}
Work out in detail the further general theory for the quantum groups $H_N^\Gamma$, including toral subgroups, and crossed product decomposition results for them.
\end{exercise}

All this is normally quite standard, and if getting lost, you can always take a look at the papers of Raum-Weber, where these questions are solved, and report on what you learned. As a second exercise now, which is more of an open question, we have:

\begin{exercise}
Work out the basic probabilistic aspects of the quantum groups $H_N^\Gamma$, in terms of the associated discrete groups $\Gamma$.
\end{exercise}

This is certainly something quite interesting, in waiting to be solved, for some time already, and the work here is most likely quite routine. 

\chapter{The real case}

\section*{10a. General strategy}

Good news, we have now all the needed ingredients for doing some exciting classification work. We will discuss here, still following \cite{ez1}, and then the papers of Raum-Weber \cite{rw1}, \cite{rw2}, \cite{rw3}, the classification problem for the orthogonal easy quantum groups:
$$S_N\subset G\subset O_N^+$$

We will see that a full classification result is available in this case, and with this being one of the main achievements of the classification work for the easy quantum groups. The discussion will use pretty much everything that we learned so far in this book, at the examples and classification level, along with a number of supplementary ingredients.

\bigskip

In order to explain what is to be done, let us go back to the standard cube of easy quantum groups, that we are very familiar with, namely:
$$
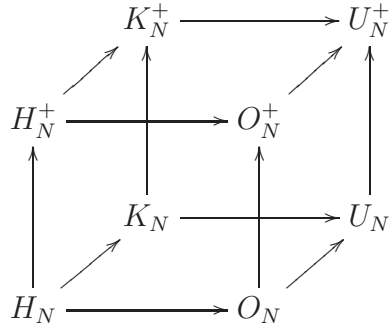
\xymatrix@R=18pt@C=18pt{
&K_N^+\ar[rr]&&U_N^+\\
H_N^+\ar[rr]\ar[ur]&&O_N^+\ar[ur]\\
&K_N\ar[rr]\ar[uu]&&U_N\ar[uu]\\
H_N\ar[uu]\ar[ur]\ar[rr]&&O_N\ar[uu]\ar[ur]
}$$

Assuming that we are in the twistable case, $H_N\subset G\subset O_N^+$, we are interested in understanding the inner objects for the face of the cube facing us, namely:
$$
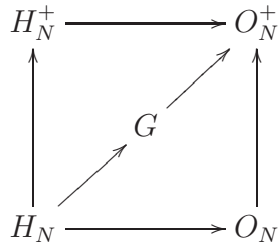
\xymatrix@R=22pt@C=22pt{
H_N^+\ar[rr]&&O_N^+\\
&G\ar[ur]\\
H_N\ar[uu]\ar[rr]\ar[ur]&&O_N\ar[uu]}$$

But here, the idea is very simple, namely that of ``projecting'' $G$ on the upper and lower edges, as to reach to a diagram as follows:
$$\xymatrix@R=24pt@C=44pt{
H_N^+\ar[r]&G_{free}\ar[r]&O_N^+\\
&G\ar[u]\\
H_N\ar[uu]\ar[r]&G_{class}\ar[r]\ar[u]&O_N\ar[uu]}$$

Indeed, we know from chapter 6 what the values of both the classical and free quantum groups $G_{class},G_{free}$ can be, so we will be left in this way with a quite routine study of the intermediate objects for the liberation operation $G_{class}\subset G_{free}$. 

\bigskip

This was for the idea, and I can feel the following question coming right away:

\begin{question}
Yes, but why not projecting $G$ on the left and right edges? Or, even better, projecting $G$ on all $4$ edges?
\end{question}

Good point, so an alternative technique would be indeed to project $G$ on the left and right edges, as to reach to a diagram as follows:
$$\xymatrix@R=30pt@C=42pt{
H_N^+\ar[rr]&&O_N^+\\
G_{disc}\ar[r]\ar[u]&G\ar[r]&G_{cont}\ar[u]\\
H_N\ar[u]\ar[rr]&&O_N\ar[u]}$$

In practice, however, this won't really work, because the correspondence $G_{disc}\leftrightarrow G_{cont}$ is very far from being bijective. To be more precise, according to our various results for the left and right edges, on the left we have an uncountable family, probably followed by some more objects, including $H_N^+$, while on the right we only have 3 objects. That is, the data that we will get from the above diagram will be as follows, not very usable:
$$G_{disc}\in\left\{H_N^\Gamma,\ldots, H_N^+\right\}\quad,\quad G_{cont}\in\left\{O_N,O_N^*,O_N^+\right\}$$

In contrast, our first idea, the one above, using the correspondence $G_{class}\leftrightarrow G_{free}$, looks quite valuable, because this correspondence, that we are already familiar with, is not exactly bijective, but is not far from being bijective either. So, good idea, and as a piece of homework for us, lying ahead, let us record the following question:

\begin{problem}
Clarify the bijectivity aspects of the correspondence
$$G_{class}\leftrightarrow G_{free}$$
before seriously getting to classification work.
\end{problem}

As a funny comment here, the first systematic classification paper for the easy quantum groups was \cite{ez1}, and this problem was of course duly investigated prior to that work, with full scientific integrity, and with a quick ``yes'' answer, based on a previous mistake from \cite{bsp}, leading to the conclusion that $G_{class}\leftrightarrow G_{free}$ is trivially bijective. Later Weber came in \cite{web} with a counterexample/fix for \cite{bsp}, which was of course taken into account afterwards, in the classification work of Raum-Weber \cite{rw1}, \cite{rw2}, \cite{rw3}. More on this later, but believe me, ``no mistakes'' usually means in mathematics mediocrity, so please stay away from that, study interesting problems, and do mistakes from time to time. 

\bigskip

Back to work now, everything that has been said above regards the twistable case, $H_N\subset G\subset O_N^+$, which is the simplest. But normally the same technique should apply as well to the general case, $S_N\subset G\subset O_N^+$, via a diagram as follows:
$$\xymatrix@R=24pt@C=44pt{
S_N^+\ar[r]&G_{free}\ar[r]&O_N^+\\
&G\ar[u]\\
S_N\ar[uu]\ar[r]&G_{class}\ar[r]\ar[u]&O_N\ar[uu]}$$

To be more precise, again we know from chapter 6 what the classical and free orthogonal quantum groups are, and this should normally allow us to recover $G$, by performing a case-by-case classification work, for each liberation operation $G_{class}\subset G_{free}$.

\bigskip

This was for the idea, and in practice, our starting point will be the result from chapter 6 regarding the classical and free orthogonal quantum groups. The statement, along with some further details regarding the proof, that we will need here, is as follows:

\begin{theorem}
The classical and free orthogonal easy quantum groups are
$$
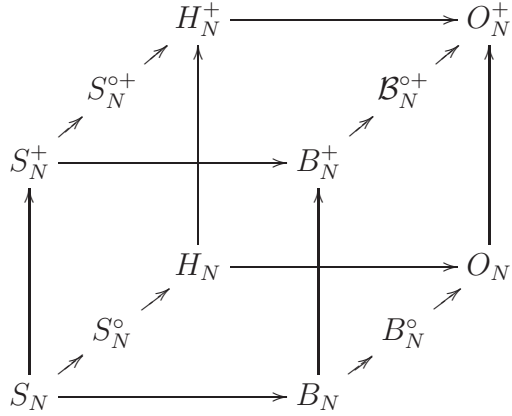
\xymatrix@R=7pt@C=7pt{
&&H_N^+\ar[rrrr]&&&&O_N^+\\
&S_N^{\circ+}\ar[ur]&&&&\mathcal B_N^{\circ+}\ar[ur]\\
S_N^+\ar[rrrr]\ar[ur]&&&&B_N^+\ar[ur]\\
\\
&&H_N\ar[rrrr]\ar[uuuu]&&&&O_N\ar[uuuu]\\
&S_N^\circ\ar[ur]&&&&B_N^\circ\ar[ur]\\
S_N\ar[uuuu]\ar[ur]\ar[rrrr]&&&&B_N\ar[uuuu]\ar[ur]
\\
}$$
with $S_N^\circ=S_N\times\mathbb Z_2$, $B_N^\circ=B_N\times\mathbb Z_2$, and with $S_N^{\circ+},\mathcal B_N^{\circ+}$ being their liberations, where $\mathcal B_N^{\circ+}$ stands for the two possible such liberations, $B_N^{\circ+}\subset B_N^{\circ\circ+}$.
\end{theorem}

\begin{proof}
This is something that we know from chapter 6, with everything being quite standard, except for the ramification question for the liberations of $B_N^\circ$, that we will recall now. The continuous face of the cube, on the right, looks in detail as follows:
$$
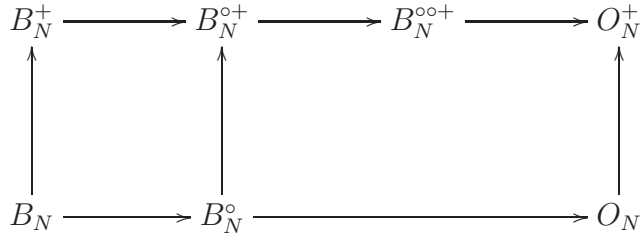
\xymatrix@R=55pt@C=46pt{
B_N^+\ar[r]&B_N^{\circ+}\ar[r]&B_N^{\circ\circ+}\ar[r]&O_N^+\\
B_N\ar[u]\ar[r]&B_N^\circ\ar[u]\ar[rr]&&O_N\ar[u]}$$

As for the corresponding categories of partitions, these are as follows, containing objects that we know, except for $NC_{12}^\circ,NC_{12}^{\circ\circ}$, whose definition will come in a moment:
$$
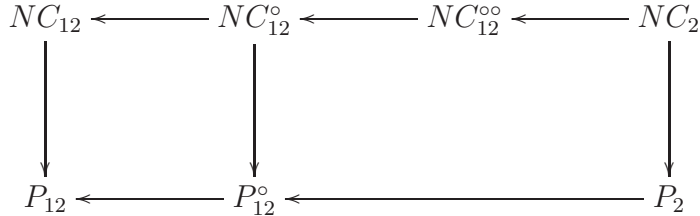
\xymatrix@R=50pt@C=44pt{
NC_{12}\ar[d]&NC_{12}^\circ\ar[l]\ar[d]&NC_{12}^{\circ\circ}\ar[l]&NC_2\ar[l]\ar[d]\\
P_{12}&P_{12}^\circ\ar[l]&&P_2\ar[ll]}$$

Getting now to the core of the problem, we know that $B_N^\circ=B_N\times\mathbb Z_2$ appears from the category $P_{12}^\circ$ of singletons and pairings, having an even total length. The point now is that we have the following formulae for $P_{12}^\circ$, which are both clear:
$$P_{12}^\circ
=<\slash\hskip-2.1mm\backslash\,,|\cap\hskip-4.2mm{\ }_|>
=<\slash\hskip-2.1mm\backslash\,,||>$$

Now when liberating at the level of these formulae, that is, when removing the crossing, we obtain two possibly distinct categories, as follows:
$$NC_{12}^\circ=<|\cap\hskip-4.2mm{\ }_|>
\quad\supset\quad 
NC_{12}^{\circ\circ}=<||>$$

Observe that we have indeed an inclusion as above, due to the following formula:
$$||=[^\pi_\sigma]\in<|\cap\hskip-4.2mm{\ }_|>\quad:\quad\pi=|\cap\hskip-4.2mm{\ }_|\quad,\quad\sigma=||\cup$$

However, we do not have equality, due to the following somewhat bizarre fact:
$$|\cap\hskip-4.2mm{\ }_|\ \notin\ <||>$$

So, this was for the story, the idea being that we have two noncrossing versions of $P_{12}^\circ$, and so two liberations of $B_N^\circ$, as constructed above. And for details, regarding all this, and then the fix of the previous classification from \cite{bsp}, we refer to \cite{web}.
\end{proof}

Getting back now to our classification program for the orthogonal easy quantum groups $S_N\subset G\subset O_N^+$, we recall that our idea is very simple, namely that of regarding $G$ as being inside the cube from Theorem 10.3, and then ``projecting'' it on the upper and lower faces, as to reach to a diagram as follows, where for convenience we have collapsed to 2D:
$$\xymatrix@R=22pt@C=40pt{
S_N^+\ar[r]&G_{free}\ar[r]&O_N^+\\
&G\ar[u]\\
S_N\ar[uu]\ar[r]&G_{class}\ar[r]\ar[u]&O_N\ar[uu]}$$

Then, we will be left with a most likely quite routine classification problem for the intermediate objects for the liberation $G_{class}\subset G_{free}$, and with this latter liberation taking 7 possible values, according to our classification result from Theorem 10.3.

\bigskip

However, in view of the issues with $B_N^\circ$ from Theorem 10.3, we must be very careful with all this, especially when talking about $G_{class}$ and $G_{free}$. So, following \cite{ez1}, and in mind with the update coming from \cite{web}, let us start with the following definition:

\begin{definition}
Consider an easy quantum group $S_N\subset G\subset O_N^+$, coming from a category of partitions $NC_2\subset D\subset P$.
\begin{enumerate}
\item The classical version $S_N\subset G_{class}\subset O_N$ is obtained by setting $G_{class}=G\cap O_N$. Equivalently, $G$ is the easy group coming from $<D,\slash\hskip-2.1mm\backslash>$.

\item The free version $S_N^+\subset G_{free}\subset O_N^+$ is obtained by setting $G_{free}=\{G,S_N^+\}$. That is, $G_{free}$ is the easy quantum group coming from $D\cap NC$.
\end{enumerate}
\end{definition}

Here we have used the general material regarding the operations $\cap$ and $\{\,,\}$ from chapter 3, and we refer to that chapter for more on all this. 

\bigskip

In practice now, at the quantum group level we obtain the rectangular diagram given above, right before Definition 10.4, and at the level of categories of partitions we obtain the following rectangular diagram, which is dual to the quantum group diagram:
$$\xymatrix@R=22pt@C=40pt{
NC\ar[dd]&D\cap NC\ar[l]\ar[d]&NC_2\ar[dd]\ar[l]\\
&D\ar[d]\\
P&<D,\slash\hskip-2.1mm\backslash>\ar[l]&P_2\ar[l]}$$

In relation now with Theorem 10.3, and with the questions that we are interested in, in relation with our classification program, we have the following result:

\begin{proposition}
The following happen:
\begin{enumerate}
\item We have $(G_{free})_{class}=G_{class}$.

\item However, we can have $(G_{class})_{free}\neq G_{free}$.
\end{enumerate}
\end{proposition}

\begin{proof}
At the first glance, this might look as chapter 3 grade material, obtained by playing with the categories in Definition 10.4, but since we have $\neq$ in (2), this is certainly more subtle than that. To be more precise, we must use Theorem 10.3, as follows:

\medskip

(1) With respect to the cube in Theorem 10.3, computing $(G_{free})_{class}$ means going up and down, while computing $G_{class}$ means simply going down. But these operations lead to the same outcome, namely one of the 6 objects on the lower face.

\medskip

(2) We have a similar picture here, with computing $(G_{class})_{free}$ meaning going down and up, and computing $G_{free}$ meaning going up. But these operations lead to the same outcome, except in the case $G_{free}=B_N^{\circ\circ+}$, which cannot appear as $(G_{class})_{free}$, because due to our definitions of $B_N^{\circ+},B_N^{\circ\circ+}$, explained in the proof of Theorem 10.3, we have:
$$(B_N^\circ)_{free}=B_N^{\circ+}\neq B_N^{\circ\circ+}$$ 

Thus, we have the result, the simplest counterexample being $G=B_N^{\circ\circ+}$.
\end{proof}

Following \cite{ez1}, we begin our classification work with a technical result, valid in the general case. Given an easy quantum group $S_N\subset G\subset O_N^+$, let us set, as before:
$$G_{class}=G\cap O_N$$

We know from the general properties of $\cap$ that, if we denote by $D$ the category of partitions associated to $G$, then the category of partitions associated to $G_{class}$ is:
$$D_{class}=<D,\slash\hskip-2.1mm\backslash>$$

Consider as well the free version of $G$, defined as before as follows:
$$G_{free}=\{G,S_N^+\}$$

According to our definition for the easy generation operation $\{\,,\}$, the category of partitions for this latter quantum group is then given by:
$$D_{free}=D\cap NC$$

Finally, let us call a category of partitions ``even'' when it consists of partitions having an even number of legs. That is, $D$ is even when the following happens:
$$k+l\in 2\mathbb N+1\implies D(k,l)=\emptyset$$

Observe that this is the case for the category of partitions $P^\circ$ associated to the group $S_N^\circ=S_N\times\mathbb Z_2$, which consists precisely of the partitions $\pi\in P$ having an even number of legs. In fact, by functoriality, the fact that a category of partitions $D$ is even means precisely that the corresponding easy quantum group $G$ appears as follows:
$$S_N^\circ\subset G\subset O_N^+$$ 

With these conventions, we have the following result, from \cite{ez1}, further building on what we know from the above, regarding the classical and free version:

\begin{proposition}
Given an easy quantum group $S_N\subset G\subset O_N^+$ as above, denote by $\Lambda,\Lambda_{free}\subset\mathbb N$ the sets of possible sizes of blocks of the partitions of $D,D_{free}$.
\begin{enumerate}
\item $\Lambda_{free}\subset\Lambda\subset \Lambda_{free}\cup(\Lambda_{free}-1)$.

\item $1\in\Lambda$ implies $1\in\Lambda_{free}$.

\item If $D_{free}$ is even, so is $D$.
\end{enumerate}
\end{proposition}

\begin{proof}
We will heavily use the various abstract notions and results in \cite{bsp}:

\medskip

(1) The first inclusion in the statement, namely $\Lambda_{free}\subset\Lambda$, follows from the following inclusion of categories of partitions, which itself comes from definitions:
$$D_{free}\subset D$$

(2) As for the second inclusion, namely $\Lambda\subset \Lambda_{free}\cup(\Lambda_{free}-1)$, this is equivalent to the following statement: ``If $\beta$ is a block of a partition $\pi\in D$, then there exists a certain block $\beta'$ of a certain partition $\pi'\in D_{free}$, having size $|\beta|$ or $|\beta|-1$''. 

\medskip

(3) But this latter statement follows by using the ``capping'' method in \cite{bsp}. Indeed, we can cap $\pi$ with semicircles, as for $\beta$ to remain unchanged, and we end up with a certain partition $\pi'$ consisting of $\beta$ and of some extra points, at most one point between any two legs of $\beta$, which might be connected or not. Note that the semicircle capping being a categorical operation, this partition $\pi'$ remains in $D_{free}$.

\medskip

(4) Now by further capping $\pi'$ with semicircles, as to get rid of the extra points, the size of $\beta$ can only increase, and we end up with a one-block partition having size at least that of $\beta$. This one-block partition is obviously noncrossing, and by capping it again with semicircles we can reduce the number of legs up to $|\beta|$ or $|\beta|-1$, and we are done.

\medskip

(5) Getting now to the second assertion, the condition $1\in\Lambda$ in the statement means that there exists $\pi\in D$ having a singleton. By capping $\pi$ with semicircles outside this singleton, we can obtain a singleton, or a double singleton. Since both these partitions are noncrossing, and have a singleton, we obtain $1\in\Lambda_{free}$, and we are done. 

\medskip

(6) Finally, regarding the last assertion, assume by contradiction that $D$ is not even, and consider a partition $\pi\in D$ having an odd number of legs. By capping $\pi$ with enough semicircles we can arrange for ending up with a singleton, and since this singleton is by definition in $D_{free}\cap NC$, we obtain our contradiction, and are done.
\end{proof}

We are now in position of splitting the classification. Recall from Theorem 10.3 that the classical and free orthogonal easy quantum groups are as follows:
$$
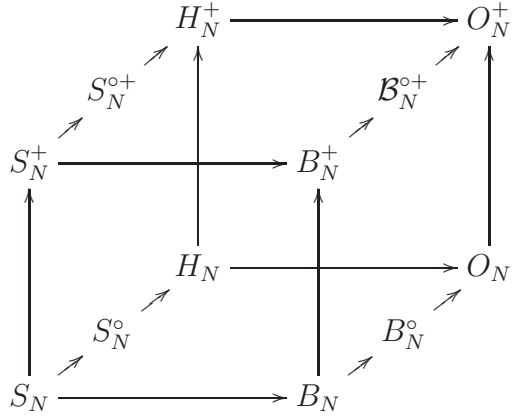
\xymatrix@R=7pt@C=7pt{
&&H_N^+\ar[rrrr]&&&&O_N^+\\
&S_N^{\circ+}\ar[ur]&&&&\mathcal B_N^{\circ+}\ar[ur]\\
S_N^+\ar[rrrr]\ar[ur]&&&&B_N^+\ar[ur]\\
\\
&&H_N\ar[rrrr]\ar[uuuu]&&&&O_N\ar[uuuu]\\
&S_N^\circ\ar[ur]&&&&B_N^\circ\ar[ur]\\
S_N\ar[uuuu]\ar[ur]\ar[rrrr]&&&&B_N\ar[uuuu]\ar[ur]
\\
}$$

With this in mind, we can further build on what we know from Proposition 10.6, on a case-by-case basis, and we are led to the following key result, from \cite{ez1}:

\begin{theorem}
Given an easy quantum group $S_N\subset G\subset O_N^+$ as before, construct its classical version $S_N\subset G_{class}\subset O_N$.
\begin{enumerate}
\item If $G_{class}\neq H_N$ then $G_{class}\subset G\subset G_{free}$.

\item If $G_{class}=H_N$ then $S_N^\circ\subset G\subset H_N^+$.
\end{enumerate}
\end{theorem}

\begin{proof}
We recall that the inclusion $G\subset G_{free}$ follows from definitions. For the other inclusion, we have 7 cases, depending on the exact value of the easy group $G_{class}$, and of $G_{free}$, and we can solve each of these cases by using Proposition 10.6, as follows:

\medskip

(1) $G_{class}=O_N$. Here we have $\Lambda_{free}=\{2\}$, so we get $\{2\}\subset\Lambda\subset\{1,2\}$. Moreover, again by Proposition 10.6, we get $\Lambda=\{2\}$. Thus $D\subset P_2$, which gives $O_N\subset G$.

\medskip

(2) $G_{class}=S_N$. Here there is nothing to prove, because we have an inclusion $S_N\subset G$, by definition of our easy quantum group $G$.

\medskip

(3) $G_{class}=B_N$. Here we have $\Lambda_{free}=\{1,2\}$, so we get $\Lambda=\{1,2\}$. Thus we have $D\subset P_{12}$, which gives an inclusion $B_N\subset G$.

\medskip

(4) $G_{class}=S_N^\circ$. Here we have an inclusion $D\subset P$ by definition, and we deduce that we have $D\subset P^\circ$, which gives an inclusion $S_N^\circ\subset G$.

\medskip

(5) $G_{class}=B_N^\circ$ and $G_{free}=B_N^{\circ+}$. Here we have $\Lambda_{free}=\{1,2\}$, so we get $\Lambda=\{1,2\}$. This gives $D\subset P_{12}$, and we get $D\subset P_{12}^\circ$, which gives an inclusion $B_N^\circ\subset G$.

\medskip

(6) $G_{class}=B_N^\circ$ and $G_{free}=B_N^{\circ\circ+}$. Here we have again $\Lambda_{free}=\{1,2\}$, so we get $\Lambda=\{1,2\}$. This gives $D\subset P_{12}$, and we get $D\subset P_{12}^\circ$, which gives $B_N^\circ\subset G$.

\medskip

(7) $G_{class}=H_N$. Here we have $D\subset P$ by definition, and we deduce that we have $D\subset P^\circ$, which gives an inclusion $S_N^\circ\subset G$.
\end{proof}

We can see from the above result that the case $G_{class}=H_N$ is quite special, and this is in tune with our findings from chapter 9, where we have seen that $H_N$ has uncountably many liberations, and with the classification of such liberations being actually not over yet. This ramification phenomenon will play a key role, in what follows.

\section*{10b. Liberation study}

According to what we have so far, namely Theorem 10.3 and Theorem 10.7, we are left with a case-by-case study, of the easy intermediate objects $G_N$ for various liberation operations as follows, with the endpoints $G_N$ and $G_N^+$ being known:
$$G_N\subset G_N^\times\subset G_N^+$$ 

We have already seen such questions in this book, notably in chapter 6 for the group $S_N$, and then in chapter 7 for the group $O_N$, with the conclusion that in these two cases, the complete lists of such liberations are very short, as follows:
$$S_N\subset S_N^+\quad,\quad O_N\subset O_N^*\subset O_N^+$$

Moreover, in both these cases the proofs were quite similar, basically based on the method of semicircle capping. We refer to chapters 6 and 7 for the whole story, and for our purposes here, let us record these findings a bit informally, as follows:

\begin{fact}
The easy liberations of $G_N=S_N,O_N$ are as follows:
\begin{enumerate}
\item They can be classified via semicircle capping.

\item They consist of the half-liberations $G_N^*$ and of the free versions $G_N^+$. 

\item With the remark that for $G_N=S_N$, we have $S_N^*=S_N$.
\end{enumerate}
\end{fact}

Getting back now to our liberation problem, in general, as formulated above, involving an arbitary inclusion $G_N\subset G_N^+$, our goal will be that of extending this type of finding to all the orthogonal easy groups that we know, by using similar methods, and with of course the remark that the case of $H_N$ is special, needing more study.

\bigskip

In practice now, in order to start now the classification, based on what we have in Theorem 10.7, we will need the following notions:

\index{semicircle capping}
\index{singleton capping}
\index{doubleton capping}

\begin{definition}
Let $\pi\in P(k,l)$ be a partition, with the points counted modulo $k+l$, counterclockwise starting from bottom left.
\begin{enumerate}
\item We call semicircle capping of $\pi$ any partition obtained from $\pi$ by connecting with a semicircle a pair of consecutive neighbors.

\item We call singleton capping of $\pi$ any partition obtained from $\pi$ by capping one of its legs with a singleton.

\item We call doubleton capping of $\pi$ any partition obtained from $\pi$ by capping two of its legs with singletons.
\end{enumerate}
\end{definition}

In other words, the semicircle, singleton and doubleton cappings are elementary operations on partitions, which lower the total number of legs by $2,1,2$ respectively. 

\bigskip

Observe that there are $k+l$ possibilities for placing the semicircle or the singleton, and $(k+l)(k+l-1)/2$ possibilities for placing the double singleton. Observe also that in the case of 2 particular ``semicircle cappings'', namely those at left or at right, the semicircle in question is rather a vertical bar, but we will still call it semicircle.

\bigskip

With these conventions, we have the following technical result, extending some technical results from chapters 6 and 7, that we will heavily use in what follows:

\begin{proposition}
Let $\pi$ be a partition, having $j$ legs.
\begin{enumerate}
\item If $\pi\in P_2-P_2^*$ and $j>4$, there exists a semicircle capping $\pi'\in P_2-P_2^*$.

\item If $\pi\in P_2^*-NC_2$ and $j>6$, there exists a semicircle capping $\pi'\in P_2^*-NC_2$.

\item If $\pi\in P-NC$ and $j>4$, there exists a singleton capping $\pi'\in P-NC$.

\item If $\pi\in P_{12}-NC_{12}$ and $j>4$, there exists a singleton capping $\pi'\in P_{12}-NC_{12}$.

\item If $\pi\in P^\circ-NC^\circ$ and $j>4$, there exists a doubleton capping $\pi'\in P^\circ-NC^\circ$.

\item If $\pi\in P_{12}^\circ-NC_{12}^\circ$ and $j>4$, there exists a doubleton capping $\pi'\in P_{12}^\circ-NC_{12}^\circ$.
\end{enumerate}
\end{proposition}

\begin{proof}
We write $\pi\in P(k,l)$, so that the number of legs is $j=k+l$. In the cases where our partition is a pairing, we use as well the number of strings, $s=j/2$. Let us agree that all partitions are drawn as to have a minimal number of crossings. 

\medskip

We will use the same idea for all the proofs, namely to ``isolate'' a block of $\pi$ having a crossing, or an odd number of crossings, then to ``cap'' $\pi$ as in the statement, as for this block to remain crossing, or with an odd number of crossings. 

\medskip

Here we use of course the observation that the ``balancing'' condition which defines the categories of partitions $P_2^*,P_{even}$ can be interpreted as saying that each block of the partition has an even number of crossings, when the picture of the partition is drawn such that this number of crossings is minimal.

\medskip

(1) The assumption $\pi\notin P_2^*$ means that $\pi$ has certain strings having an odd number of crossings. We fix such an ``odd'' string, and we try to cap $\pi$, as for this string to remain odd in the resulting partition $\pi'$. An examination of all the possible pictures shows that this is possible, provided that our partition has $s>2$ strings, and we are done.

\medskip

(2) The assumption $\pi\notin NC_2$ means that $\pi$ has certain crossing strings. We fix such a pair of crossing strings, and we try to cap $\pi$, as for these strings to remain crossing in $\pi'$. Once again, an examination of all the possible pictures shows that this is possible, provided that our partition has $s>3$ strings, and we are done.

\medskip

(3) Indeed, since $\pi$ is crossing, we can choose two of its blocks which are intersecting. If there are some other blocks left, we can cap one of their legs with a singleton, and we are done. If not, this means that our two blocks have a total of $j'\geq j>4$ legs, so at least one of them has $j''>2$ legs. One of these $j''$ legs can always be capped with a singleton, as for the capped partition to remain crossing, and we are done.

\medskip

(4) Here we can simply cap with a singleton, as in (3).

\medskip

(5) Here we can cap with a doubleton, by proceeding twice as in (3).

\medskip

(6) Here we can cap again with a doubleton, by proceeding twice as in (3).
\end{proof}

As before with what we knew from chapters 6 and 7, involved in the proof of Fact 10.8, we can apply several times what we found in Proposition 10.10, by recurrence, and we are led in this way to the following result, also from \cite{ez1}, which is finer:

\begin{proposition}
Let $\pi$ be a partition.
\begin{enumerate}
\item If $\pi\in P_2-P_2^*$ then $<\pi,NC_2>=P_2$.

\item If $\pi\in P_2^*-NC_2$ then $<\pi,NC_2>=P_2^*$.

\item If $\pi\in P-NC$ then $<\pi,NC>=P$.

\item If $\pi\in P_{12}-NC_{12}$ then $<\pi,NC_{12}>=P_{12}$.

\item If $\pi\in P^\circ-NC^\circ$ then $<\pi,NC^\circ>=P^\circ$.

\item If $\pi\in P_{12}^\circ-NC_{12}^\circ$ then $<\pi,NC_{12}^\circ>=P_{12}^\circ$.
\end{enumerate}
\end{proposition}

\begin{proof}
We use what we have in Proposition 10.10, with the observation that the ``capping partition'' appearing there is always in the good category. That is, we use the following facts, which are all clear from the definition of the categories involved:

\medskip

--  The semicircle is in $NC_2,NC^\circ$.

\medskip

-- The singleton is in $NC,NC_{12}$.

\medskip

-- The doubleton is in $NC_{12}^\circ$. 

\medskip

The point now is that, in the context of the capping operations in Proposition 10.10, these observations tell us that, in each of the cases under consideration, the category to be computed can only decrease when replacing $\pi$ by one of its cappings $\pi'$. 

\medskip

Indeed, for the singleton and doubleton cappings this is clear from definitions, and for the semicircle capping this is clear as well from definitions, unless in the case where the ``capping semicircle'' is actually a ``bar'' added at left or at right, where we can use a categorical rotation operation as in \cite{bsp}.

\medskip

(1) This assertion can be proved by recurrence on the number of strings, $s=(k+l)/2$. Indeed, by using Proposition 10.10 (1), for $s>3$ we have a descent procedure $s\to s-1$, and this leads to the situation $s\in\{1,2,3\}$, where the statement is clear.

\medskip

(2) Again, this can be proved by recurrence on the number of strings, $s=(k+l)/2$. Indeed, by using Proposition 10.10 (2), for $s>3$ we have a descent procedure $s\to s-1$, and this leads to the situation $s\in\{1,2,3\}$, where the statement is clear.

\medskip

(3) We can proceed by recurrence on the number of legs of $\pi$. If the number of legs is $j=4$, then $\pi$ is a basic crossing, and we have $<\pi>=P$. If the number of legs is $j>4$ we can apply Proposition 10.10 (3), and the result follows from:
$$<\pi>\supset <\pi'>=P$$

(4) This is similar to the proof of (1), by using Proposition 10.10 (4).

\medskip

(5) This is again similar to the proof of (1), by using Proposition 10.10 (5).

\medskip

(6) This is again similar to the proof of (1), by using Proposition 10.10 (6).
\end{proof}

All this might seem quite technical, but good news, we are almost there, with what we need in practice. As usual by building on what we already knew from chapters 6 and 7, involved in the proof of Fact 10.8, we can reformulate what we found in Proposition 10.11 in a more convenient way, the result here, still from \cite{ez1}, being as follows:

\begin{proposition}
Let $\pi$ be a partition.
\begin{enumerate}
\item If $\pi\in P_2$ then $<\pi,NC_2>\in\{P_2,P_2^*,NC_2\}$.

\item If $\pi\in P$ then $<\pi,NC>\in\{P,NC\}$.

\item If $\pi\in P_{12}$ then $<\pi,NC_{12}>\in\{P_{12},NC_{12}\}$.

\item If $\pi\in P^\circ$ then $<\pi,NC'>\in\{P^\circ,NC^\circ\}$.

\item If $\pi\in P_{12}^\circ$ then $<\pi,NC_{12}^\circ>\in\{P_{12}^\circ,NC_{12}^\circ\}$.
\end{enumerate}
\end{proposition}

\begin{proof}
This follows indeed by rearranging the various technical results above, and more specifically by suitably interpreting what we found in Proposition 10.11.
\end{proof}

We are now in position of stating a main result. Recall from Theorem 10.3 that the classical and free orthogonal easy quantum groups are as follows:
$$
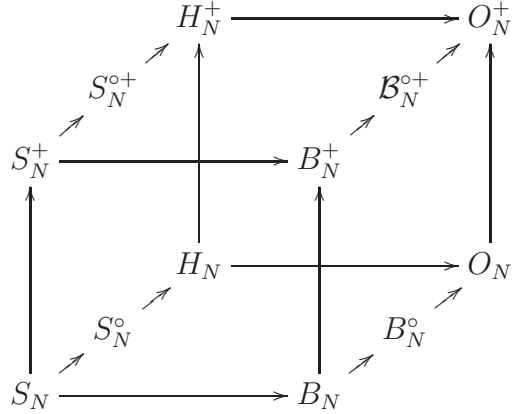
\xymatrix@R=7pt@C=7pt{
&&H_N^+\ar[rrrr]&&&&O_N^+\\
&S_N^{\circ+}\ar[ur]&&&&\mathcal B_N^{\circ+}\ar[ur]\\
S_N^+\ar[rrrr]\ar[ur]&&&&B_N^+\ar[ur]\\
\\
&&H_N\ar[rrrr]\ar[uuuu]&&&&O_N\ar[uuuu]\\
&S_N^\circ\ar[ur]&&&&B_N^\circ\ar[ur]\\
S_N\ar[uuuu]\ar[ur]\ar[rrrr]&&&&B_N\ar[uuuu]\ar[ur]
\\
}$$

With this cube in mind, and by taking as well into account the various issues in the special case of the hyperoctahedral group $H_N$, coming from the previous chapter, and from Theorem 10.7 as well, let us formulate the following definition:

\begin{definition}
We call ``non-hyperoctahedral'' any easy quantum group 
$$S_N\subset G\subset O_N^+$$
such that $G_{class}\neq H_N$. 
\end{definition}

We refer to the above for various interpretations of this condition. Now with this convention made, we have the following classification result, for such quantum groups:

\index{non-hyperoctahedral}
\index{modified symmetric group}
\index{modified bistochastic group}

\begin{theorem}
There are exactly $13$ non-hyperoctahedral orthogonal easy quantum groups, namely:
\begin{enumerate}
\item $O_N,O_N^*,O_N^+$: the orthogonal quantum groups.

\item $S_N,S_N^+$: the symmetric quantum groups.

\item $B_N,B_N^+$: the bistochastic quantum groups.

\item $S_N^\circ,S_N^{\circ+}$: the modified symmetric quantum groups.

\item $B_N^\circ,B_N^{\circ+}$: the modified bistochastic quantum groups.

\item $B_N^{\circ\circ*},B_N^{\circ\circ+}$: the extra modified bistochastic quantum groups.
\end{enumerate}
\end{theorem}

\begin{proof}
This basically follows from what we have, the idea being as follows:

\medskip

(1) We know from Proposition 10.10 that what we have to do is to classify the easy quantum groups satisfying $G_{class}\subset G\subset G_{free}$. 

\medskip

(2) More precisely, leaving the issues with the liberations of $B_N^\circ$ aside, we have to prove that for $G_{class}=S_N,B_N,S_N^\circ,B_N^\circ$ there is no such partial liberation, and that for $G_{class}=O_N$ there is only one partial liberation, namely the quantum group $G_{class}^*$.

\medskip

(3) But this follows from the various results from Proposition 10.12, via the Tannakian results in \cite{bsp}, which provide us with the list of 11 objects in the statement.

\medskip

(4) However, as before with other results from the old papers \cite{bsp} and \cite{ez1}, there are some problems here coming from the exact cube from Theorem 10.3. We refer to \cite{web} and \cite{rw3} for the updates and fixes of all this, the idea being that the quantum group $B_N^{\circ\circ+}$, as well as its half-classical version $B_N^{\circ\circ*}=B_N^{\circ\circ+}\cap O_N^*$, must be added to the list in the statement, and with these changes made, the result holds as stated.
\end{proof}

There are many things that can be said about the above result. As a first observation, our classification so far can be reformulated in the following more intuitive form:

\begin{theorem}
The orthogonal easy quantum groups $S_N\subset G\subset O_N^+$ are
$$
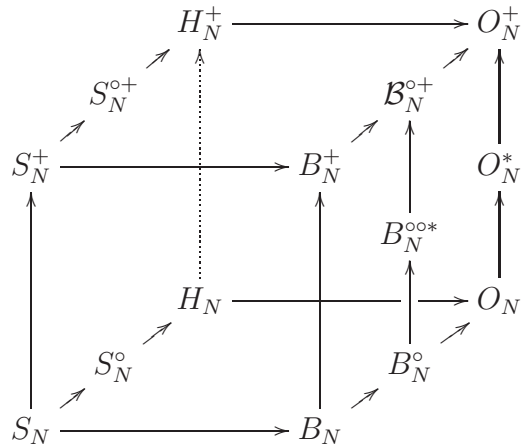
\xymatrix@R=7pt@C=7pt{
&&H_N^+\ar[rrrr]&&&&O_N^+\\
&S_N^{\circ+}\ar[ur]&&&&\mathcal B_N^{\circ+}\ar[ur]\\
S_N^+\ar[rrrr]\ar[ur]&&&&B_N^+\ar[ur]&&O_N^*\ar[uu]\\
&&&&&B_N^{\circ\circ*}\ar[uu]\\
&&H_N\ar@{-}[rrr]\ar@.[uuuu]&&&\ar[r]&O_N\ar[uu]\\
&S_N^\circ\ar[ur]&&&&B_N^\circ\ar[ur]\ar[uu]\\
S_N\ar[uuuu]\ar[ur]\ar[rrrr]&&&&B_N\ar[uuuu]\ar[ur]
\\
}$$
with $\mathcal B_N^{\circ+}$ standing for $B_N^{\circ+},B_N^{\circ\circ+}$, and with the dotted arrow still to be investigated.
\end{theorem}

\begin{proof}
This follows indeed from Theorem 10.14, and with the remark that the vertical arrow $B_N^{\circ\circ*}\to\mathcal B_N^{\circ+}$ lands by definition into $B_N^{\circ\circ+}$, as to have indeed an inclusion.
\end{proof}

In regards with the above result, the right face of the cube, which is the continuous one, still deserves some more work, pictorially speaking, and this due to various non-functoriality phenomena which appear. Here is a better result, regarding that face:

\begin{theorem}
The continuous orthogonal easy quantum groups, that is, the intermediate easy quantum groups $B_N\subset G\subset O_N^+$, are as follows,
$$
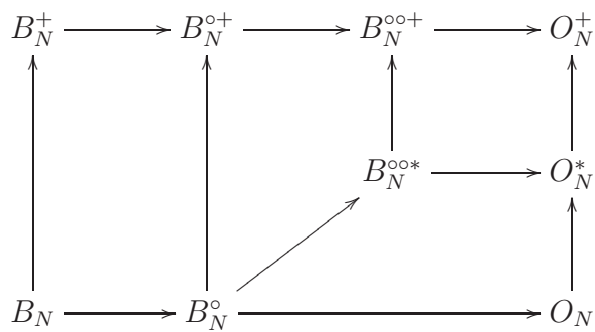
\xymatrix@R=35pt@C=40pt{
B_N^+\ar[r]&B_N^{\circ+}\ar[r]&B_N^{\circ\circ+}\ar[r]&O_N^+\\
&&B_N^{\circ\circ*}\ar[u]\ar[r]&O_N^*\ar[u]\\
B_N\ar[uu]\ar[r]&B_N^\circ\ar[uu]\ar[rr]\ar[ur]&&O_N\ar[u]}$$
with the half-classical versions of $B_N^+,B_N^{\circ+}$ collapsing to $B_N,B_N^\circ$.
\end{theorem}

\begin{proof}
This follows indeed from Theorem 10.15, and with the whole diagram being in fact the correct version of the right face of the cube in Theorem 10.15. For the sake of completness, and for further reference, let us record as well the diagram of the categories of partitions for the quantum groups in the statement. This is as follows:
$$
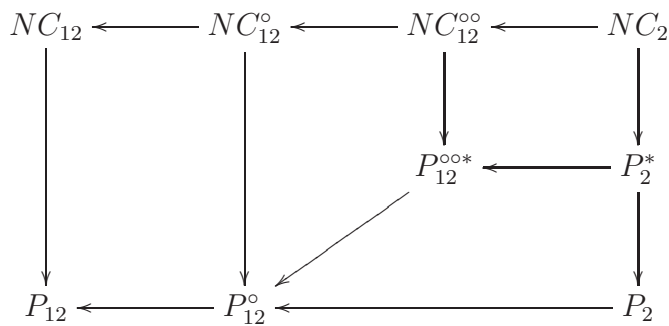
\xymatrix@R=35pt@C=40pt{
NC_{12}\ar[dd]&NC_{12}^\circ\ar[l]\ar[dd]&NC_{12}^{\circ\circ}\ar[l]\ar[d]&NC_2\ar[l]\ar[d]\\
&&P_{12}^{\circ\circ*}\ar[dl]&P_2^*\ar[d]\ar[l]\\
P_{12}&P_{12}^\circ\ar[l]&&P_2\ar[ll]}$$

For more on these categories, and their meaning, we refer to Theorem 10.3.
\end{proof}

It is of course possible to come up with some more diagrams for the above classification results, and their various particular cases of interest. However, we will defer the discussion here to the end of the present chapter, after solving the hyperoctahedral case as well.

\section*{10c. Higher reflections}

We are now in position of finishing the classification. The idea, from \cite{rw3}, is that, with Theorem 10.14 taking care of the non-hyperoctahedral case, we are left with a study in the hyperoctahedral case. But here there is a dichotomy coming from $H_N^{[\infty]}$, which in the combinatorial language of \cite{rw3} corresponds to a dichotomy coming from group-theoretical categories, as opposed to non-group-theoretical categories. 

\bigskip

We recall from chapter 9 that given a uniform reflection group $\mathbb Z_2^{*N}\to\Gamma\to\mathbb Z_2^N$, we can associate to it subsets $D(k,l)\subset P(k,l)$, which form a category of partitions, as follows:
$$D(k,l)=\left\{\pi\in P(k,l)\Big|\ker\binom{i}{j}\leq\pi\implies g_{i_1}\ldots g_{i_k}=g_{j_1}\ldots g_{j_l}\right\}$$

Observe that we have $P_{even}^{[\infty]}\subset D\subset P_{even}$, with the inclusions coming respectively from $\eta\in D$, and from $\Gamma\to\mathbb Z_2^N$. Conversely, given a category of partitions $P_{even}^{[\infty]}\subset D\subset P_{even}$, we can associate to it a uniform reflection group $\mathbb Z_2^{*N}\to\Gamma\to\mathbb Z_2^N$, as follows:
$$\Gamma=\left\langle g_1,\ldots g_N\Big|g_{i_1}\ldots g_{i_k}=g_{j_1}\ldots g_{j_l},\forall i,j,k,l,\ker\binom{i}{j}\in D(k,l)\right\rangle$$

As explained in \cite{rw2}, the correspondences $\Gamma\to D$ and $D\to\Gamma$ are bijective, and inverse to each other, at $N=\infty$. We have in fact the following result, from \cite{rw1}, \cite{rw2}, \cite{rw3}:

\begin{theorem}
We have correspondences between:
\begin{enumerate}
\item Uniform reflection groups $\mathbb Z_2^{*\infty}\to\Gamma\to\mathbb Z_2^\infty$.

\item Categories of partitions $P_{even}^{[\infty]}\subset D\subset P_{even}$.

\item Easy quantum groups $G=(G_N)$, with $H_N^{[\infty]}\supset G_N\supset H_N$.
\end{enumerate}
\end{theorem}

\begin{proof}
This is something quite tricky, and we refer here to \cite{rw3}. As an illustration, as mentioned above, we have the following correspondences:
$$\xymatrix@R=17mm@C=27mm{
\mathbb Z_2^N\ar@{~}[d]&\mathbb Z_2^{\circ N}\ar[l]\ar@{~}[d]&\mathbb Z_2^{*N}\ar[l]\ar@{~}[d]\\
H_N\ar[r]&H_N^*\ar[r]&H_N^{[\infty]}}$$

More generally, for any $s\in\{2,4,\ldots,\infty\}$, the quantum groups $H_N^{(s)}\subset H_N^{[s]}$ constructed in \cite{ez1} come from the quotients of $\mathbb Z_2^{\circ N}\leftarrow\mathbb Z_2^{*N}$ by the relations $(ab)^s=1$:
$$\xymatrix@R=17mm@C=17mm{
\mathbb Z_2^N\ar@{~}[d]&\mathbb Z_2^{\circ N}/<(ab)^s=1>\ar[l]\ar@{~}[d]&\mathbb Z_2^{*N}/<(ab)^s=1>\ar[l]\ar@{~}[d]\\
H_N\ar[r]&H_N^{(s)}\ar[r]&H_N^{[s]}}$$

For details on all this, and more, we refer to \cite{rw3}.
\end{proof}

The structure and classification results discussed above, concerning the intermediate easy quantum groups $H_N\subset G\subset H_N^{[\infty]}$, do not close the classification problem in general, for the easy quantum groups $H_N\subset G\subset H_N^+$. The point indeed is that we have for instance intermediate objects for the following inclusion:
$$H_N^{[\infty]}\subset G\subset H_N^+$$

In order to discuss this question, which does have a non-trivial answer, let us start with the following construction, from \cite{rw3}, which is something quite tricky:

\begin{proposition}
Let $H_N^{\diamond r}\subset H_N^+$ be the easy quantum group coming from the following partition:
$$\pi_r=\ker\begin{pmatrix}1&\ldots&r&r&\ldots&1\\1&\ldots&r&r&\ldots&1\end{pmatrix}$$
We have then inclusions between these quantum groups, as follows,
$$H_N^{[\infty]}\subset\ldots\subset H_N^{\diamond 3}\subset H_N^{\diamond 2}\subset H_N^{\diamond 1}=H_N^+$$
and all these inclusions are proper.
\end{proposition}

\begin{proof}
We have several things to be proved, the idea being as follows:

\medskip

(1) Consider indeed the quantum group $H_N^{\diamond r}\subset H_N^+$ coming from the partition $\pi_r$ in the statement, which is by definition easy. As a first illustration for this construction, let us examine the case $r=1$. Here our partition $\pi_1$ is something familiar, namely:
$$
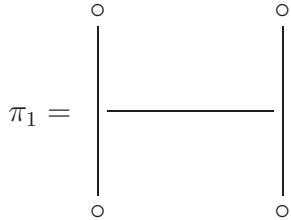
\xymatrix@R=10mm@C=20mm{\\ \pi_1=}
\xymatrix@R=10mm@C=20mm{
\circ\ar@{-}[dd]&\circ\ar@{-}[dd]\\
\ar@{-}[r]&\\
\circ&\circ
}$$

Now since we have $\pi_1\in NC_{even}$, we obtain $H_N^{\diamond 1}=H_N^+$, as claimed.

\medskip

(2) Let us discuss we well the case $r=2$. Here the partition $\pi_2$ in the statement, producing the subgroup $H_N^{\diamond 2}\subset H_N^+$ is as follows:
$$
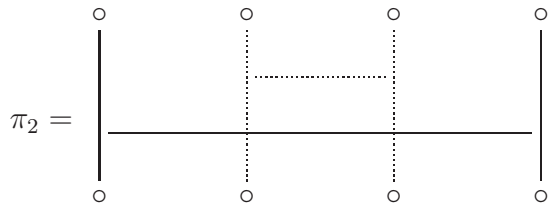
\xymatrix@R=4mm@C=10mm{\\ \\ \pi_2=}
\xymatrix@R=5mm@C=15mm{
\circ\ar@{-}[ddd]&\circ\ar@{.}[ddd]&\circ\ar@{.}[ddd]&\circ\ar@{-}[ddd]\\
&\ar@{.}[r]&\\
\ar@{-}[rrr]&&&\\
\circ&\circ&\circ&\circ
}$$

In order to prove our results regarding $H_N^{\diamond 2}$, our first claim is that we have $H_N^{[\infty]}\subset H_N^{\diamond 2}$. By functoriality, this amounts in checking that we have:
$$<\pi_2>\subset P_{even}^{[\infty]}$$

Thus we must check that we have $\pi_2\subset P_{even}^{[\infty]}$, and this is clear from either of the various explicit descriptions of the category $P_{even}^{[\infty]}$, obtained before.

\medskip

(3) In order to finish now our study of $H_N^{\diamond 2}$, consider the inclusions of quantum groups that we established in the above, namely:
$$H_N^{[\infty]}\subset H_N^{\diamond 2}\subset H_N^{\diamond 1}=H_N^+$$

We must prove that these inclusions are proper, which amounts in proving that the reverse inclusions for the corresponding categories, which are as follows, are proper:
$$P_{even}^{[\infty]}\supset <\pi_2>\supset NC_{even}$$

But this follows by carefully examining the partition $\pi_2$, and the category of partitions that it generates, with the conclusion that this category is indeed as above. To be more precise, since $\pi_2$ is crossing we certainly have a proper embedding on the right, and the fact that the embedding on the left is proper too is standard. See \cite{rw3}.

\medskip

(4) In the general case now, $r\in\mathbb N$, our first claim is that we have $H_N^{[\infty]}\subset H_N^{\diamond r}$. By functoriality, this amounts in checking that we have:
$$<\pi_r>\subset P_{even}^{[\infty]}$$

Thus we must check that, with $\pi_r$ being as in the statement, the following happens:
$$\pi_r\subset P_{even}^{[\infty]}$$

But this is clear from either of the various explicit descriptions of the category $P_{even}^{[\infty]}$, obtained before, and we obtain in this way the result.

\medskip

(5) Let us prove now that we have inclusions $H_N^{\diamond(r+1)}\subset H_N^{\diamond r}$ as in the statement, for any $r\in\mathbb N$. At $r=2$, to start with, the partition $\pi_3$ is as follows:
$$
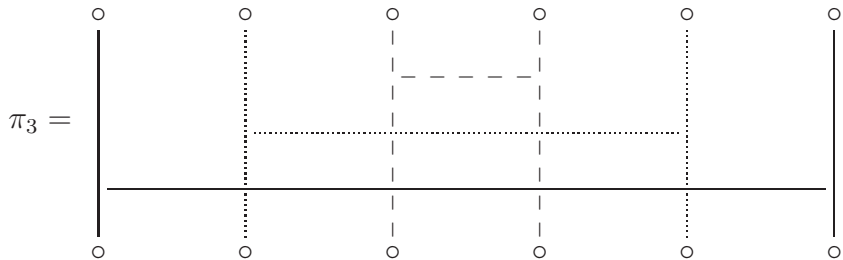
\xymatrix@R=4mm@C=10mm{\\ \ \\ \pi_3=}
\xymatrix@R=5mm@C=15mm{
\circ\ar@{-}[dddd]&\circ\ar@{.}[dddd]&\circ\ar@{--}[dddd]&\circ\ar@{--}[dddd]&\circ\ar@{.}[dddd]&\circ\ar@{-}[dddd]\\
&&\ar@{--}[r]&\\
&\ar@{.}[rrr]&&&\\
\ar@{-}[rrrrr]&&&&&\\
\circ&\circ&\circ&\circ&\circ&\circ
}$$

But, it is clear that by capping with semicircles, in the obvious way, we can obtain the partition $\pi_2$ for this partition. Thus, we have indeed $H_N^{\diamond3}\subset H_N^{\diamond 2}$, and the proof of $H_N^{\diamond(r+1)}\subset H_N^{\diamond r}$ in general is similar, by suitably capping $\pi_{r+1}$ with semicircles.

\medskip

(6) Finally, the fact that the inclusions $H_N^{\diamond(r+1)}\subset H_N^{\diamond r}$ that we obtained are indeed proper is best seen at the categorical level, coming from the fact that we have proper inclusions of categories of partitions, and for details here we refer here to \cite{rw3}.
\end{proof}

Quite remarkably, we have the following uniqueness result, also from \cite{rw3}:

\begin{theorem}
Let $H_N^{\diamond r}\subset H_N^+$ be the easy quantum group coming from:
$$\pi_k=\ker\begin{pmatrix}1&\ldots&r&r&\ldots&1\\1&\ldots&r&r&\ldots&1\end{pmatrix}$$
We have then inclusions between these quantum groups, as follows,
$$H_N^{[\infty]}\subset\ldots\subset H_N^{\diamond 3}\subset H_N^{\diamond 2}\subset H_N^{\diamond 1}=H_N^+$$
and these are all the easy quantum groups $H_N^{[\infty]}\subset G\subset H_N^+$, satisfying $G\neq H_N^{[\infty]}$.
\end{theorem}

\begin{proof}
Here the first part of the statement is something that we already know, from Proposition 10.18, reproduced here for convenience. As for the last assertion, regarding uniqueness, this is something quite technical, and we refer here to \cite{rw3}.
\end{proof}

There are many other things that can be said about the quantum groups $H_N^{\diamond r}$ introduced above, both at the algebraic and probabilistic level. Let us start with:

\begin{proposition}
The quantum group $H_N^{\diamond r}\subset H_N^+$ appears via the relations
$$\delta_{ab}u_{a_1i_1}\ldots u_{a_ri_r}u_{a_rj_r}\ldots u_{a_1j_1}=\delta_{ij}u_{a_1i_1}\ldots u_{a_ri_r}u_{b_ri_r}\ldots u_{b_1i_1}$$
applied to the standard coordinates $u_{ij}$.
\end{proposition}

\begin{proof}
We know that quantum group $H_N^{\diamond r}\subset H_N^+$ appears by definition by imposing the following relations to the standard coordinates $u_{ij}$:
$$T_{\pi_r}\in End(u^{\otimes 2r})$$

In order to interpret these relations, let us first compute the operator $T_{\pi_r}$. We know that the partition $\pi_r$ is, pictorially speaking, as follows:
$$
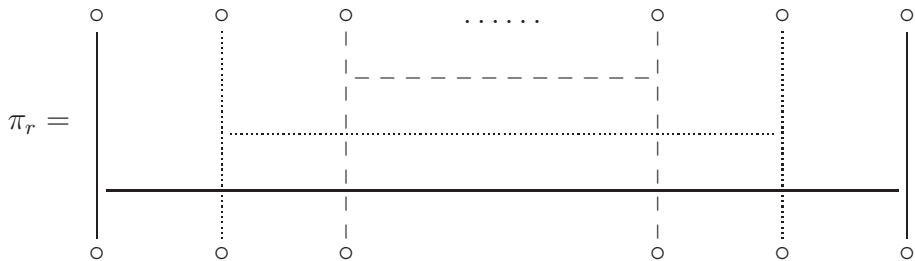
\xymatrix@R=4mm@C=8mm{\\ \ \\ \pi_r=}
\xymatrix@R=5mm@C=12mm{
\circ\ar@{-}[dddd]&\circ\ar@{.}[dddd]&\circ\ar@{--}[dddd]&\ldots\ldots&\circ\ar@{--}[dddd]&\circ\ar@{.}[dddd]&\circ\ar@{-}[dddd]\\
&&\ar@{--}[rr]&&\\
&\ar@{.}[rrrr]&&&&\\
\ar@{-}[rrrrrr]&&&&&&\\
\circ&\circ&\circ&&\circ&\circ&\circ
}$$

Thus, the operator associated to this partition is as follows:
$$T_{\pi_r}(e_{i_1}\otimes\ldots\otimes e_{i_r}\otimes e_{j_r}\otimes\ldots\otimes e_{j_1})=\delta_{ij}e_{i_1}\otimes\ldots\otimes e_{i_r}\otimes e_{i_r}\otimes\ldots\otimes e_{i_1}$$

With this formula in hand, we have the following computation:
\begin{eqnarray*}
&&T_{\pi_r}u^{\otimes 2r}(e_{i_1}\otimes\ldots\otimes e_{i_r}\otimes e_{j_r}\otimes\ldots\otimes e_{j_1}\otimes1)\\
&=&T_{\pi_r}\left(\sum_{abij}e_{a_1i_1}\otimes\ldots\otimes e_{a_ri_r}\otimes e_{b_rj_r}\otimes\ldots\otimes e_{b_1j_1}\otimes u_{a_1i_1}\ldots u_{a_ri_r}u_{b_rj_r}\ldots u_{b_1j_1}\right)\\
&&(e_{i_1}\otimes\ldots\otimes e_{i_r}\otimes e_{j_r}\otimes\ldots\otimes e_{j_1}\otimes1)\\
&=&T_{\pi_r}\sum_{ab}e_{a_1}\otimes\ldots\otimes e_{a_r}\otimes e_{b_r}\otimes\ldots\otimes e_{b_1}\otimes u_{a_1i_1}\ldots u_{a_ri_r}u_{b_rj_r}\ldots u_{b_1j_1}\\
&=&\sum_ae_{a_1}\otimes\ldots\otimes e_{a_r}\otimes e_{a_r}\otimes\ldots\otimes e_{a_1}\otimes u_{a_1i_1}\ldots u_{a_ri_r}u_{a_rj_r}\ldots u_{a_1j_1}
\end{eqnarray*}

On the other hand, we have as well the following computation:
\begin{eqnarray*}
&&u^{\otimes 2r}T_{\pi_r}(e_{i_1}\otimes\ldots\otimes e_{i_r}\otimes e_{j_r}\otimes\ldots\otimes e_{j_1}\otimes1)\\
&=&\delta_{ij}u^{\otimes 2}(e_{i_1}\otimes\ldots\otimes e_{i_r}\otimes e_{i_r}\otimes\ldots\otimes e_{i_1}\otimes1)\\
&=&\delta_{ij}\left(\sum_{abij}e_{a_1i_1}\otimes\ldots\otimes e_{a_ri_r}\otimes e_{b_rj_r}\otimes\ldots\otimes e_{b_1j_1}\otimes u_{a_1i_1}\ldots u_{a_ri_r}u_{b_rj_r}\ldots u_{b_1j_1}\right)\\
&&(e_{i_1}\otimes\ldots\otimes e_{i_r}\otimes e_{i_r}\otimes\ldots\otimes e_{i_1}\otimes1)\\
&=&\delta_{ij}\sum_{ab}e_{a_1}\otimes\ldots\otimes e_{a_r}\otimes e_{b_r}\otimes\ldots\otimes e_{b_1}\otimes u_{a_1i_1}\ldots u_{a_ri_r}u_{b_ri_r}\ldots u_{b_1i_1}
\end{eqnarray*}

We conclude that $T_{\pi_r}u^{\otimes 2r}=u^{\otimes 2r}T_{\pi_r}$ is equivalent to the following condition, which must be satisfied for all the indices involved, namely $i,j$:
\begin{eqnarray*}
&&\sum_ae_{a_1}\otimes\ldots\otimes e_{a_r}\otimes e_{a_r}\otimes\ldots\otimes e_{a_1}\otimes u_{a_1i_1}\ldots u_{a_ri_r}u_{a_rj_r}\ldots u_{a_1j_1}\\
&=&\delta_{ij}\sum_{ab}e_{a_1}\otimes\ldots\otimes e_{a_r}\otimes e_{b_r}\otimes\ldots\otimes e_{b_1}\otimes u_{a_1i_1}\ldots u_{a_ri_r}u_{b_ri_r}\ldots u_{b_1i_1}
\end{eqnarray*}

By looking at the summands, the following must happen, for any $a,b,i,j$:
$$\delta_{ab}u_{a_1i_1}\ldots u_{a_ri_r}u_{a_rj_r}\ldots u_{a_1j_1}=\delta_{ij}u_{a_1i_1}\ldots u_{a_ri_r}u_{b_ri_r}\ldots u_{b_1i_1}$$

Thus, we are led to the conclusion in the statement.
\end{proof}

Here are a number of supplementary results regarding the quantum groups $H_N^{\diamond r}$, which can be useful in practice, when dealing with these quantum groups:

\begin{theorem}
The quantum groups $H_N^{\diamond r}$ have the following properties:
\begin{enumerate}
\item Their diagonal torus is $\mathbb Z_2^{*N}$, independently on $r$.

\item These quantum groups are not coamenable.

\item Their intersection is the quantum group $H_N^{[\infty]}$.
\end{enumerate}
\end{theorem}

\begin{proof}
All this is routine from what we have, the idea being as follows:

\medskip

(1) This is best seen by functoriality. Indeed, we know from Proposition 10.18 that we have inclusions of quantum groups as follows:
$$H_N^{[\infty]}\subset H_N^{\diamond r}\subset H_N^+$$

Thus, at the level of diagonal tori, we obtain inclusions as follows:
$$\widehat{\mathbb Z_2^{*N}}\subset T\subset\widehat{\mathbb Z_2^{*N}}$$

We therefore conclude that we have $T=\widehat{\mathbb Z_2^{*N}}$, independently of $r$, as stated.

\medskip

(2) This follows from (1), the diagonal torus being non-coamenable.
 
\medskip

(3) We know from  Proposition 10.18 that we have inclusions as follows:
$$H_N^{[\infty]}\subset\ldots\subset H_N^{\diamond 3}\subset H_N^{\diamond 2}\subset H_N^{\diamond 1}=H_N^+$$

Now consider the following intersection, which is a decreasing intersection:
$$G=\bigcap_rH_N^{\diamond r}$$

This intersection is then an easy quantum group, appearing as follows:
$$H_N^{[\infty]}\subset G\subset H_N^+$$

Now by using the classification result from Theorem 10.19, along with the fact that the inclusions between the quantum groups $H_N^{\diamond r}$ are proper, that we know to hold from Proposition 10.18, we conclude that we have $G=H_N^{[\infty]}$, as stated.
\end{proof}

Getting back now to classification matters, what we have in Theorem 10.17 and Theorem 10.19 is still not enough. Fortunately, $H_N^{[\infty]}$ produces a dichotomy, and there are no further examples, the final classification result, from \cite{rw3}, being as follows:

\begin{theorem}
The easy quantum groups $H_N\subset G\subset H_N^+$ are as follows,
$$H_N\subset H_N^\Gamma\subset H_N^{[\infty]}\subset H_N^{\diamond r}\subset H_N^+$$
with the family $H_N^\Gamma$ covering $H_N,H_N^{[\infty]}$, and with the series $H_N^{\diamond r}$ covering $H_N^+$.
\end{theorem}

\begin{proof}
This follows from the various classification results above, with a bit more work, the idea being as follows:

\medskip

(1) The easy quantum groups $H_N\subset G\subset H_N^+$ can be shown to be either of the form $H_N\subset G\subset H_N^{[\infty]}$, or of the form $H_N^{[\infty]}\subset G\subset H_N^+$. 

\medskip

(2) But with these two latter classification problems being solved by our various classification results, we obtain the result. We refer here to the paper of Raum-Weber \cite{rw3}. 
\end{proof}

\section*{10d. Classification results}

All this is very nice, and is exactly what we need, in order to finish our classification work. As a first result, in the twistable case, where our orthogonal easy quantum group $S_N\subset G\subset O_N^+$ contains $H_N$, the classification result, from \cite{rw3}, is as follows:

\begin{theorem}
The easy quantum groups $H_N\subset G\subset O_N^+$ are as follows,
$$
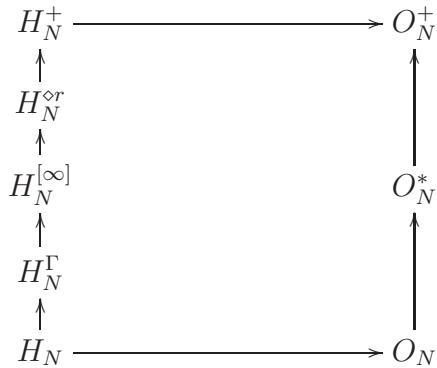
\xymatrix@R=4mm@C=40mm{
H_N^+\ar[r]&O_N^+\\
H_N^{\diamond r}\ar[u]\\
H_N^{[\infty]}\ar[u]&O_N^*\ar[uu]\\
H_N^\Gamma\ar[u]\\
H_N\ar[r]\ar[u]&O_N\ar[uu]}$$
with the family $H_N^\Gamma$ covering $H_N,H_N^{[\infty]}$, and with the series $H_N^{\diamond r}$ covering $H_N^+$.
\end{theorem}

\begin{proof}
This follows indeed from the various results above, and from those in chapter 9. For further details, we refer to the paper of Raum and Weber \cite{rw3}.
\end{proof}

Regarding now the general orthogonal easy quantum group case, $S_N\subset G\subset O_N^+$, we can formulate things here as follows:

\begin{theorem}
The orthogonal easy quantum groups are as follows:
\begin{enumerate}
\item $O_N,O_N^*,O_N^+$: the orthogonal quantum groups.

\item $S_N,S_N^+$: the symmetric quantum groups.

\item $B_N,B_N^+$: the bistochastic quantum groups.

\item $S_N^\circ,S_N^{\circ+}$: the modified symmetric quantum groups.

\item $B_N^\circ,B_N^{\circ+}$: the modified bistochastic quantum groups.

\item $B_N^{\circ\circ*},B_N^{\circ\circ+}$: the extra modified bistochastic quantum groups.

\item $H_N,H_N^\Gamma,H_N^{\diamond r},H_N^+$: the hyperoctahedral quantum groups.
\end{enumerate}
\end{theorem}

\begin{proof}
This follows indeed from what we have, by combining Theorem 10.15, which deals with the non-hyperoctahedral case, and Theorem 10.23 and its versions, dealing with the hyperoctahedral case. For more on all this, we refer as usual to \cite{rw3}.
\end{proof}

There are many things that can be said about the above result. As a first observation, our classification can be reformulated in the following more intuitive form:

\begin{theorem}
The orthogonal easy quantum groups $S_N\subset G\subset O_N^+$ are
$$
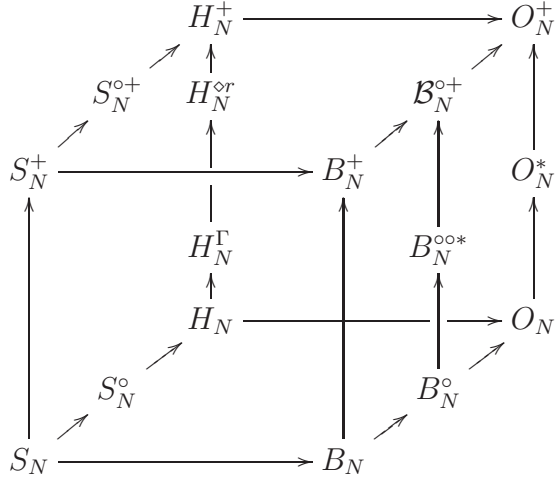
\xymatrix@R=9pt@C=9pt{
&&H_N^+\ar[rrrr]&&&&O_N^+\\
&S_N^{\circ+}\ar[ur]&H_N^{\diamond r}\ar[u]&&&\mathcal B_N^{\circ+}\ar[ur]\\
S_N^+\ar[rrrr]\ar[ur]&&\ar[u]&&B_N^+\ar[ur]&&O_N^*\ar[uu]\\
&&H_N^\Gamma\ar@{-}[u]&&&B_N^{\circ\circ*}\ar[uu]\\
&&H_N\ar@{-}[rrr]\ar[u]&&&\ar[r]&O_N\ar[uu]\\
&S_N^\circ\ar[ur]&&&&B_N^\circ\ar[ur]\ar[uu]\\
S_N\ar[uuuu]\ar[ur]\ar[rrrr]&&&&B_N\ar[uuuu]\ar[ur]
}$$
with $\mathcal B_N^{\circ+}$ standing for the quantum groups $B_N^{\circ+}\subset B_N^{\circ\circ+}$.
\end{theorem}

\begin{proof}
This follows indeed from Theorem 10.24, and with the remark that the vertical arrow $B_N^{\circ\circ*}\to\mathcal B_N^{\circ+}$ lands by definition into $B_N^{\circ\circ+}$, as to have indeed an inclusion.
\end{proof}

All this is very nice, we have reached to our objectives, formulated in the beginning of this chapter. Let us record as well the result in the half-classical case, as follows:

\begin{theorem}
The half-classsical orthogonal easy quantum groups are as follows:
\begin{enumerate}
\item $O_N,O_N^*$: the orthogonal quantum groups.

\item $S_N$: the symmetric group.

\item $B_N$: the bistochastic group.

\item $S_N^\circ$: the modified symmetric group.

\item $B_N^\circ$: the modified bistochastic group.

\item $B_N^{\circ\circ*}$: the half-classical modified bistochastic quantum group.

\item $H_N^\Gamma$, with $\Gamma$ half-classical: the hyperoctahedral quantum groups.
\end{enumerate}
\end{theorem}

\begin{proof}
This follows indeed from Theorem 10.25, by removing from there the free versions, and the quantum groups from (7) which are not half-classical. Alternatively, as explained in \cite{web}, it is possible to obtain this result directly, by classifying the half-classical categories of partitions. For more on all this, we refer to \cite{web}.
\end{proof}

There are of course many other interesting particular cases of Theorem 10.25, such as those concerning the uniform case, and so on. In addition, it is of course possible to draw some nice diagrams for the quantum groups involved, and for the corresponding categories of partitions, and other objects such as the diagonal subgroups, as well.

\section*{10e. Exercises} 

The material in this chapter has been quite exciting, and we have several interesting exercises about all this, which are rather research questions. First, we have:

\begin{exercise}
Work out the basic probabilistic aspects of the quantum groups $H_N^{\diamond r}$, and some further algebraic aspects as well.
\end{exercise}

This is certainly something quite interesting, in waiting to be solved, for some time already. As a second exercise, again as difficult and interesting as they get, we have:

\begin{exercise}
Find an abstract contravariant duality between the intermediate easy objects for the inclusion $U_N\subset U_N^+$, and the inclusion $H_N\subset H_N^+$.
\end{exercise}

This might seem quite puzzling, but the thing is that, we know from chapter 8 that the intermediate objects for $U_N\subset U_N^+$ consist of a series followed by a family, and we also know from this chapter that the intermediate objects for $H_N\subset H_N^+$ consist of a family followed by a series. So, our question makes sense. We will be actually back to this, later in this book, but just with some further comments, in the lack of a solution.

\chapter{Complex reflections}

\section*{11a. Reflection groups}

In this chapter we keep building on the theory developed in chapters 9-10, with complex versions of the constructions performed there, and generalizations of some of the classification results obtained there. In order to explain our strategy, let us go back to the standard cube formed by the main easy quantum groups, namely:
$$
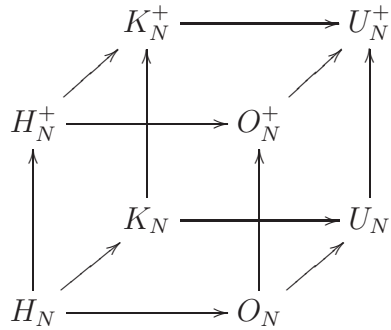
\xymatrix@R=18pt@C=18pt{
&K_N^+\ar[rr]&&U_N^+\\
H_N^+\ar[rr]\ar[ur]&&O_N^+\ar[ur]\\
&K_N\ar[rr]\ar[uu]&&U_N\ar[uu]\\
H_N\ar[uu]\ar[ur]\ar[rr]&&O_N\ar[uu]\ar[ur]
}$$

We have seen that the intermediate easy quantum groups $H_N\subset G\subset H_N^+$ can be fully classified, and that with a bit more work, this leads to a full classification of the easy quantum groups $S_N\subset G\subset H_N^+$, which can be thought of as being the easy ``real quantum reflection groups''. Moreover, with a bit more work, in the continuous case, this even leads to a classification of the orthogonal easy quantum groups, $S_N\subset G\subset O_N^+$. 

\bigskip

Our aim here is to do a similar work in the unitary case, first for the intermediate easy quantum groups $K_N\subset G\subset K_N^+$, which can be thought of as being the easy ``purely complex quantum reflection groups'', then for the intermediate easy quantum groups $H_N\subset G\subset K_N^+$, corresponding to the left face of the cube, and finally for the intermediate easy quantum groups $S_N\subset G\subset K_N^+$, which can be thought of as being the easy ``quantum reflection groups''. We will comment as well on the consequences of this to the classification of the general easy quantum groups, $S_N\subset G\subset U_N^+$.

\bigskip

Getting started now, let us first formulate the following broad definition, which covers all the examples that we have in mind, and perhaps some more examples too:

\begin{definition}
A quantum reflection group is an intermediate subgroup
$$S_N\subset G\subset K_N^+$$
between the symmetric group $S_N$, and the quantum reflection group $K_N^+=\mathbb T\wr_*S_N^+$.
\end{definition}

We are of course mostly interested in the easy case, but it is instructive to start with a study in general, without easiness assumption. Indeed, in the classical case already, the situation is very interesting, and we have here the following celebrated result of  Shephard and Todd, which is arguably on par with the ABCDEFG classification of Lie groups:

\index{complex reflection group}
\index{Shephard-Todd}

\begin{theorem}
The irreducible complex reflection groups are
$$H_N^{sd}=\left\{U\in H_N^s\Big|(\det U)^d=1\right\}$$
along with $34$ exceptional examples.
\end{theorem}

\begin{proof}
This is something quite advanced, the idea being as follows:

\medskip

(1) First of all, we already know that $H_N^s=\mathbb Z_s\wr S_N$ is a subgroup of the unitary group $U_N$, that we are actually very familiar with. The point now is that, the determinant $\det:U_N\to\mathbb T$ being a group morphism, imposing the condition $\det U=1$, or more generally imposing the condition $(\det U)^d=1$, for some $d\in\mathbb N$, still leaves us with a subgroup of $U_N$, that we can denote $H_N^{sd}$, as in the statement.

\medskip

(2) As basic examples of this construction, in the case $d=s$ we have of course $H_N^{sd}=H_N^s$, with this coming from the fact that we have $\det U\in\mathbb Z_s$, for any matrix $U\in H_N^s$. Observe also that this latter observation tells us to assume $d|s$ in our construction, as for the resulting group $H_N^{sd}$ not to degenerate. In fact, with this assumption made, it is easy to see that the resulting quantum groups $H_N^{sd}$ are distinct.

\medskip

(3) At the level of new examples now, of particular interest is the alternating group $A_N$, which appears at the parameter values $s=d=1$. Indeed, we know that we have $H_N^1=S_N$, and since the determinant function $\det:U_N\to\mathbb T$ produces by restriction to the permutation matrices $S_N\subset U_N$ the signature of the permutations, $\varepsilon:S_N\to\{\pm1\}$, by imposing the condition $\det U=1$ we obtain the alternating group $A_N$.

\medskip

(4) This was for the basic theory of the subgroups $H_N^{sd}\subset U_N$ in the statement. The point now is that all these subgroups are complex reflection groups, which are in addition ``irreducible'', in some intuitive sense. Moreover, and here comes the point, any irreducible complex reflection group $G\subset U_N$ can be shown to be of this form, up to some 34 exceptional examples, which can be explicitely classified.

\medskip

(5) So, this is what the statement is about. Regarding now the proof, this is something quite complicated, especially if you wish to have a complete classification, with the 34 exceptional examples involved fully classified and listed, and we refer here to the paper of Shephard and Todd \cite{sto}, and to the subsequent literature on the subject.
\end{proof}

Regarding now easiness, we know from chapter 3 that at $d=s$ the group under consideration, namely $H_s^s$ itself, is easy, the precise result being as follows:

\begin{theorem}
The group $H_N^s=\mathbb Z_s\wr S_N$ is easy, with the corresponding category of partitions $P^s$ consisting of the partitions having the property that each block, when weighted according to the rules $\circ\to +,\bullet\to -$, has as size a multiple of $s$. 
\end{theorem}

\begin{proof}
This is something that we know well, extending some well-known results at $s=1,2$, where $H_N^s$ is respectively the symmetric group $S_N$, and the hyperoctahedral group $H_N$. For full details here, we refer to chapter 1 of the present book.
\end{proof}

The above results raise the interesting question of examining the easiness features of the group $H_N^{sd}$, in general. To be more precise, we would like to know if this group is easy or not, and if not, what is its ``easy envelope'', in the sense of chapter 3.

\bigskip

As a first observation here, the case $N=2,s=4,d=2$ is special, as follows:

\begin{theorem}
The complex reflection group $H_2^{42}$ is easy, the corresponding category of partitions being as follows,
$$D(k,l)=\begin{cases}
P^2(k,l)&{\rm when}\ \underline{k}=\underline{l}(4)\\
\emptyset&{\rm otherwise}
\end{cases}$$
where $\underline{k}$ is the number $\#\circ-\#\bullet$, over the symbols of $k$.
\end{theorem}

\begin{proof}
According to the definition of $H_N^{sd}$, we have:
\begin{eqnarray*}
H_2^{42}
&=&\left\{g\in H_2^4\Big|\det g\in\mathbb Z_2\right\}\\
&=&\left\{\begin{pmatrix}a&0\\0&b\end{pmatrix},\begin{pmatrix}0&a\\b&0\end{pmatrix}\Big|a,b\in\mathbb Z_4,ab\in\mathbb Z_2\right\}\\
&=&\left\{\begin{pmatrix}a&0\\0&b\end{pmatrix},\begin{pmatrix}0&a\\b&0\end{pmatrix}\Big|a,b=\pm1\ {\rm or}\ a,b=\pm i\right\}\\
&=&H_2\cup iH_2
\end{eqnarray*}

Now observe that by functoriality, the associated Tannakian category $C$ satisfies:
$$C\subset C_{H_2^2}=span(P^2)$$

In order to compute $C$, we use the trivial fact that the fixed point relations $g^{\otimes l}\xi=\xi$, $(tg)^{\otimes l}\xi=\xi$ with $t\in\mathbb T$ imply $t^l=1$, with the usual conventions $t^\circ=t,t^\bullet=\bar{t}$ for the colored exponents. In our case, with $t=i$ we obtain that we have:
$$C(0,l)\neq\emptyset\implies i^l=1\implies\underline{l}=0(4)$$

More generally, the same method gives in fact the following implications:
$$C(k,l)\neq\emptyset\implies i^k=i^l\implies\underline{k}=\underline{l}(4)$$

We conclude from this that, with $D=(D(k,l))$ being the collection of sets in the statement, we have an inclusion as follows: 
$$C\subset span(D)$$

But this collection of sets $D$ forms a category of partitions, and by comparing with the classification results in \cite{tw1}, we obtain $C=span(D)$, as stated.
\end{proof}

Our claim now is that, provided that $H_N^{sd}$ is not one of the groups in Theorem 11.3 or Theorem 11.4, which are easy, and that we understand well, this group is not easy, and its easy envelope, in the sense of chapter 3, can be explicitly computed. 

\bigskip

In order to discuss this, let us start with a study in the continuous case. In analogy with the general construction of the complex reflection groups from Theorem 11.2, using the determinant function $\det:U_N\to\mathbb T$, we can formulate the following definition:

\begin{proposition}
Given a number $d\in\mathbb N\cup\{\infty\}$, consider the group
$$U_N^d=\left\{g\in U_N\Big|\det g\in\mathbb Z_d\right\}$$
where $\mathbb Z_d$ is the group of $d$-th roots of unity. This group is homogeneous,
$$S_N\subset U_N^d\subset U_N$$
when the parameter $d$ is even, $d\in 2\mathbb N\cup\{\infty\}$.
\end{proposition}

\begin{proof}
We recall from chapter 1 that the embedding $S_N\subset U_N$ that we use is the one given by the usual permutation matrices, namely:
$$\sigma(e_i)=e_{\sigma(i)}$$

Thus the determinant of a permutation $\sigma\in S_N$ is its signature, $\varepsilon(\sigma)\in\mathbb Z_2$, and this gives both the group property of $U_N^d$, and the last assertion.
\end{proof}

In what follows we will be mostly interested in the case $2|d$. However, the value $d=1$ is interesting and useful as well, because we have inclusions, as follows:
$$SU_N=U_N^1\subset U_N^d\subset U_N^\infty=U_N$$

By functoriality, we therefore obtain inclusions of categories, as follows:
$$C_{U_N}\subset C_{U_N^d}\subset C_{SU_N}$$

The group $U_N$ is well-known to be easy, its category being given by  $C_{U_N}=span(\mathcal P_2)$, where $\mathcal P_2$ is the category of the matching pairings. The representation theory of $SU_N$ is well-known as well, in diagrammatic terms, as explained for instance in \cite{wo2}. 

\bigskip

Regarding now $U_N^d$, with $d\in\mathbb N\cup\{\infty\}$ being arbitrary, we have here:

\begin{theorem}
The Tannakian category of $U_N^d$ appears as a part of the Tannakian category of $SU_N$, obtained by restricting the attention to the spaces $C(k,l)$ with 
$$\underline{k}=\underline{l}(d)$$
where $\underline{k}$ is the number $\#\circ-\#\bullet$, computed over all the symbols of $k$.
\end{theorem}

\begin{proof}
Our first claim is that in the finite case, $d<\infty$, we have a disjoint union decomposition as follows, where $w=e^{2\pi i/Nd}$:
$$U_N^d=SU_N\ \sqcup\ wSU_N\ \sqcup\ w^2SU_N\ \sqcup\ldots\sqcup\ w^{d-1}SU_N$$

Indeed, we have $w^N=e^{2\pi i/d}$, and so the condition $\det g\in\mathbb Z_d$ from Proposition 11.5 means $\det g=w^{Nk}$, for some $k\in\{0,1,\ldots,d-1\}$, and our claim follows from:
\begin{eqnarray*}
\det g=w^{Nk}
&\iff&\det\left(\frac{g}{w^k}\right)=1\\
&\iff&\frac{g}{w^k}\in SU_N\\
&\iff&g\in w^kSU_N
\end{eqnarray*}

Now given $g\in U_N$, $\xi\in(\mathbb C^N)^{\otimes k}$ and $\lambda\in\mathbb C$, consider the following conditions:
$$g^{\otimes k}\xi=\xi\quad,\quad (\lambda g)^{\otimes k}\xi=\xi\quad,\ \ldots\ , 
\quad (\lambda^{d-1}g)^{\otimes k}\xi=\xi$$

These conditions are then equivalent to the following conditions:
$$g^{\otimes k}\xi=\xi\quad,\quad\lambda^k=1$$ 

Now by taking $g\in SU_N$ and $\lambda=w^N$, with $w=e^{2\pi i/Nd}$ being as above, this gives the result. Finally, the assertion at $d=\infty$ can be proved in a similar way.
\end{proof}

Summarizing, the Tannakian category of $U_N^d$ appears as a part of the category computed in \cite{wo2}, and the value $d=\infty$, corresponding to $U_N$ itself, which is easy, is special. It is of course possible to go beyond this remark, but we will not need this here.

\bigskip

Let us discuss now the computation of easy envelopes. We recall from chapter 3 that we have the following definition, in the general easy quantum group case:

\begin{definition}
The easy envelope of a homogeneous quantum group $S_N\subset G\subset U_N^+$ is the easy quantum group $S_N\subset\bar{G}\subset U_N^+$ associated to the category of partitions
$$D(k,l)=\left\{\pi\in P(k,l)\Big|T_\pi\in C_{kl}\right\}$$ 
where $C=(C_{kl})$ is the Tannakian category of $G$.
\end{definition}

As a technical observation, we can in fact generalize the above construction to any closed subgroup $G\subset U_N^+$, and we have the following result:

\begin{proposition}
Given a closed subgroup $G\subset U_N^+$, construct $D\subset P$ as above, and let $S_N\subset\bar{G}\subset U_N^+$ be the easy quantum group associated to $D$. We have then
$$\bar{G}=\overline{<G,S_N>}$$
where $<G,S_N>\subset U_N^+$ is the smallest closed subgroup containing $G,S_N$.
\end{proposition}

\begin{proof}
It is well-known, and elementary to show, using Woronowicz's Tannakian duality results in \cite{wo2}, that the smallest subgroup $<G,S_N>\subset U_N^+$ from the statement exists indeed, and can be obtained by intersecting the Tannakian categories of $G,S_N$:
$$C_{<G,S_N>}=C_G\cap C_{S_N}$$

We conclude from this that for any $\pi\in P(k,l)$ we have:
$$T_\pi\in C_{<G,S_N>}(k,l)\iff T_\pi\in C_G(k,l)$$

It follows that the $D$ categories for the quantum groups $<G,S_N>$ and $G$ coincide, and so the easy envelopes $\overline{<G,S_N>}$ and $\bar{G}$ coincide as well, as stated.
\end{proof}

With these notions in hand, we can say more about the groups $U_N^d$ and $H_N^{sd}$. To start with, the easy envelope of $U_N^d$ can be computed as follows:

\begin{theorem}
The easy envelope of the group $U_N^d$ is given by
$$\overline{U_N^d}=U_N$$
for any $d\geq1$.
\end{theorem}

\begin{proof}
By functoriality, we can restrict the attention to the case $d=1$, where our group is the special unitary group:
$$U_N^1=SU_N$$

We have to prove that the following implication holds:
$$\pi\in P(k),\xi_\pi\in Fix(g^{\otimes k}),\forall g\in SU_N\implies\pi\in\mathcal P_2(k)$$

For this purpose, we will use the following isomorphism of projective versions:
$$PSU_N=PU_N$$

To be more precise, let us start with the following simple fact:
$$g^{\otimes k}\xi_\pi=\xi_\pi\implies(wg)^{\otimes k}\xi_\pi=w^k\xi_\pi,\forall w\in\mathbb T$$

In relation with the above implication, we have two cases, as follows:

\medskip

\underline{Case $\underline{k}=0$}. Here the condition $\underline{k}=0$ means by definition that $k$ has the same number of black and white legs. Thus in the above formula we have $w^k=1$, and we obtain:
$$g^{\otimes k}\xi_\pi=\xi_\pi,\forall g\in SU_N\implies h^{\otimes k}\xi_\pi=\xi_\pi,\forall h\in U_N$$

We can therefore conclude by using the Brauer result for $U_N$, which states that the vectors $\xi_\pi$ on the right are those appearing from the partitions $\pi\in\mathcal P_2(k)$.

\medskip

\underline{Case $\underline{k}\neq0$}. Here we must prove that a partition $\pi\in P(k)$ as above does not exist. In order to do so, observe first that, since $w^{\underline{k}}=\bar{w}^k$, we obtain:
$$g^{\otimes k}\xi_\pi=\xi_\pi,\forall g\in SU_N\implies h^{\otimes k\bar{k}}(\xi_\pi\otimes\xi_\pi)=(\xi_\pi\otimes\xi_\pi),\forall h\in U_N$$

But this shows that $\xi_\pi\otimes\xi_\pi$ must come from a pairing, and so $\xi_\pi$ itself must come from a pairing. Thus, as a first conclusion, we must have $\pi\in P_2(k)$.

\medskip

Since the standard coordinates $u_{ij}$ of our group $SU_N$ commute, we can permute if we want the legs of this pairing, and we are left with a pairing of the following type:
$$\pi=\cap\cap\ldots\cap$$

Now if we take into account the labels, by further permuting the legs we can assume that we are in the case $\pi=[\alpha\beta\gamma]$, where $\alpha,\beta,\gamma$ are all pairings of type $\cap\cap\ldots\cap$, with $\alpha$ being white, $\beta$ being black, and $\gamma$ being matching. Moreover, by using the Brauer result for $U_N$, the invariance condition is trivially satisfied for $\gamma$, so we can assume $\gamma=\emptyset$.

\medskip

Summarizing, we are now in the case $\pi=[\alpha\beta]$, with $\alpha,\beta$ being both of type $\cap\cap\ldots\cap$, and with $\alpha$ being white, and $\beta$ being black. With $\alpha=2r$ and $\beta=2s$, we have:
$$\xi_\pi=\sum_{i_1\ldots i_r}\sum_{j_1\ldots j_s}e_{i_1}\otimes e_{i_1}\otimes\ldots\otimes e_{i_r}\otimes e_{i_r}\otimes e_{j_1}\otimes e_{j_1}\otimes\ldots\otimes e_{j_s}\otimes e_{j_s}$$

An arbitrary matrix $g\in SU_N$ acts in the following way on this vector:
\begin{eqnarray*}
g^{\otimes k}\xi_\pi
&=&\sum_{i_1\ldots i_r}\sum_{j_1\ldots j_s}(gg^t)_{a_1b_1}\ldots(gg^t)_{a_rb_r}(\bar{g}g^*)_{c_1d_1}\ldots(\bar{g}g^*)_{c_sd_s}\\
&&e_{a_1}\otimes e_{b_1}\otimes\ldots\otimes e_{a_r}\otimes e_{b_r}\otimes e_{c_1}\otimes e_{d_1}\otimes\ldots\otimes e_{c_s}\otimes e_{d_s}
\end{eqnarray*}

Thus, in order to have $g^{\otimes k}\xi_\pi=\xi_\pi$, the matrix $gg^t$ must be a scalar multiple of the identity. Now since this latter condition is not satisfied by any $g\in SU_N$, the formula $g^{\otimes k}\xi_\pi=\xi_\pi$ does not hold in general, and so our partition $\pi$ does not exist, as desired.
\end{proof}

Getting back now to our questions regarding the reflection groups, in what follows, the most convenient for the study of $H_N^s$ and its subgroups $H_N^{sd}$ is to use the wreath product decomposition $H_N^s=\mathbb Z_s\wr S_N$. According to this formula, we have:

\begin{proposition}
Assuming that $d\in\mathbb N\cup\{\infty\}$ satisfies $2|d|[2,s]$, we have
$$H_N^{sd}=\left\{\sigma(\rho_1,\ldots,\rho_N)\Big|\sigma\in S_N,\rho_i\in\mathbb Z_s,\rho_1\ldots\rho_N\in\mathbb Z_d\right\}$$
where the group elements are given by the formula
$$\sigma(\rho_1,\ldots,\rho_N)=\sum_i\rho_ie_{\sigma(i)i}$$
and this group is homogeneous, $S_N\subset H_N^{sd}\subset U_N$.
\end{proposition}

\begin{proof}
With the convention in the statement for $\sigma(\rho_1,\ldots,\rho_N)$, we have:
$$H_N^s=\left\{\sigma(\rho_1,\ldots,\rho_N)\Big|\sigma\in S_N,\rho_i\in\mathbb Z_s\right\}$$

Consider now an arbitrary number $d\in\mathbb N\cup\{\infty\}$. According to the definition of $H_N^{sd}$, this group has the following description, where $\varepsilon:S_N\to\{\pm1\}$ is the signature:
$$H_N^{sd}=\left\{\sigma(\rho_1,\ldots,\rho_N)\Big|\sigma\in S_N,\rho_i\in\mathbb Z_s,\varepsilon(\sigma)\rho_1\ldots\rho_N\in\mathbb Z_d\right\}$$

Now when assuming $2|d$ we have $-1\in\mathbb Z_d$, and so $\varepsilon(\sigma)=\pm1\in\mathbb Z_d$, and we obtain the formula in the statement. As for the homogeneity claim, this is clear as well.
\end{proof}

Regarding now the easy envelope of $H_N^{sd}$, we have the following result:

\begin{theorem}
We have the easy envelope formula
$$\overline{H_N^{sd}}=H_N^s$$
unless we are in the case $\overline{H_2^{42}}=H_2^{42}$, which is exceptional.
\end{theorem}

\begin{proof}
We have an inclusion $H_N^{sd}\subset H_N^s$, and by functoriality, and by using as well the easiness result for $H_N^s$, we succesively obtain:
$$H_N^{sd}\subset H_N^s\implies span(P^s)\subset C\implies P^s\subset D$$

In order to prove the reverse inclusion $D\subset P^s$, we must compute the category $D$. For this purpose, it is enough to discuss the fixed points. For a partition $\pi\in P(k)$, the associated vector $T_\pi$, that we will denote here by $\xi_\pi$, is given by:
$$\xi_\pi=\sum_{i_1\ldots i_k}\delta_\pi(i_1,\ldots,i_k)e_{i_1}\otimes\ldots\otimes e_{i_k}$$

Now with $g=\sigma(\rho_1,\ldots,\rho_N)\in H_N^{sd}$, as in Proposition 11.10, we have:
$$g^{\otimes k}\xi_\pi=\sum_{i_1\ldots i_k}\delta_\pi(i_1,\ldots,i_k)\rho_{i_1}\ldots\rho_{i_k}\ e_{i_{\sigma(1)}}\otimes\ldots\otimes e_{i_{\sigma(k)}}$$

On the other hand, by replacing $i_r\to i_{\sigma(r)}$, we have as well:
\begin{eqnarray*}
\xi_\pi
&=&\sum_{i_1\ldots i_k}\delta_\pi(i_{\sigma(1)},\ldots,i_{\sigma(k)})\ e_{i_{\sigma(1)}}\otimes\ldots\otimes e_{i_{\sigma(k)}}\\
&=&\sum_{i_1\ldots i_k}\delta_\pi(i_1,\ldots,i_k)\ e_{i_{\sigma(1)}}\otimes\ldots\otimes e_{i_{\sigma(k)}}
\end{eqnarray*}

We conclude from this that the formula $g^{\otimes k}\xi_\pi=\xi_\pi$ is equivalent to:
$$\delta_\pi(i_1,\ldots,i_k)=1\implies\rho_{i_1}\ldots\rho_{i_k}=1$$

To be more precise, in order for the equality $g^{\otimes k}\xi_\pi=\xi_\pi$ to hold, this formula must hold for any numbers $\rho_1,\ldots,\rho_N\in\mathbb Z_s$ satisfying the following condition:
$$\rho_1\ldots\rho_N\in\mathbb Z_d$$

Observe that in the case $d=s$ the condition $\rho_1\ldots\rho_N\in\mathbb Z_d$ dissapears, and the condition $\delta_\pi(i_1,\ldots,i_k)=1\implies\rho_{i_1}\ldots\rho_{i_k}=1$, for any $\rho_1,\ldots,\rho_N\in\mathbb Z_s$, tells us that all the blocks of $\pi$, when weighted according to the rules $\circ\to +,\bullet\to -$, must have as size a multiple of $s$. Thus $\pi\in P^s$. Now back to our question, so far we have obtained:
$$D(k)=\left\{\pi\Big|\delta_\pi(i_1,\ldots,i_k)=1\implies\rho_{i_1}\ldots\rho_{i_k}=1,\forall\rho_1,\ldots,\rho_N\in\mathbb Z_s,\rho_1\ldots\rho_N\in\mathbb Z_d\right\}$$

In order to compute this set, let $\pi$ and $i_1,\ldots,i_k$ be as above, and consider the partition $\nu=\ker i$. We have then $\nu\leq\pi$, and since $i_1,\ldots,i_k\in\{1,\ldots,N\}$, we have $r\leq N$. Depending now on the value of $r=|\nu|$, we have two cases, as follows:

\medskip

(1) In the case $N>r$ we have a free variable among $\{\rho_1,\ldots,\rho_N\}$, that we can adjust as to have $\rho_1\ldots\rho_N\in\mathbb Z_d$. Thus, the condition $\rho_1\ldots\rho_N\in\mathbb Z_d$ dissapears, and we are left with the $H_N^s$ problem, which gives, as explained above, $\nu\in P_s$.

\medskip

(2) In the case $N=r$, let us denote by $a_1+b_1,\ldots,a_N+b_N$ the lengths of the blocks of $\nu$, with $a_i$ standing for the white legs, and $b_i$ standing for the black legs. We have:
$$\rho_1^{a_1-b_1}\ldots\rho_N^{a_N-b_N}=1,\forall\rho_1,\ldots,\rho_N\in\mathbb Z_s,\rho_1\ldots\rho_N\in\mathbb Z_d$$

With $c_i=a_i-b_i$, and with $\eta_N=\rho_1\ldots\rho_N$, we must have:
$$\rho_1^{c_1-c_N}\ldots\rho_{N-1}^{c_{N-1}-c_N}\eta_N^{c_N}=1,\forall\rho_1,\ldots,\rho_{N-1}\in\mathbb Z_s,\forall \eta_N\in\mathbb Z_d$$

Thus we must have $c_1=\ldots=c_N(s)$, and this common value must be a number $c=0(d)$. Now let us introduce the following sets:
$$P_c^{sd}=\left\{\pi\Big||\pi|=N,a_i-b_i=c(s)\right\}$$

In terms of these sets, and of their union $P^{sd}=\cup_cP_c^{sd}$, we have obtained that $\pi\in D$ happens if and only if any subpartition $\nu\leq\pi$ has the following property:

\medskip

(1) If $|\nu|<N$, then $\nu\in P^s$.

\medskip

(2) If $|\nu|=N$, then $\nu\in P^{s,d}$.

\medskip

(3) If $|\nu|>N$, no condition.

\medskip

But this shows that we must have $\pi\in P^s$, unless we are in the exceptional case, $N=2,s=4,d=2$. Thus we have $\overline{H_N^{sd}}=H_N^s$, as stated.
\end{proof}

Observe in particular that Theorem 11.11 tells us that the group $H_N^{sd}$ is not easy, unless we are in the special cases of the groups $H_N^s$ or $H_2^{42}$. Indeed, this follows from the definition of the easy envelope, from our computation of easy envelope, and from the fact that the inclusion $H_N^{sd}\subset H_N^s$ is proper, unless we are in the case $d=s$.

\section*{11b. Quantum reflections}

With the classical case reasonably understood, let us discuss now the free case. To start with, we will review, with full details, the theory of the quantum reflection groups $H_N^{s+}$. The free analogues of the reflection groups $H_N^s$ can be constructed as follows:

\begin{definition}
The algebra $C(H_N^{s+})$ is the universal $C^*$-algebra generated by $N^2$ normal elements $u_{ij}$, subject to the following relations,
\begin{enumerate}
\item $u=(u_{ij})$ is unitary,

\item $u^t=(u_{ji})$ is unitary,

\item $p_{ij}=u_{ij}u_{ij}^*$ is a projection,

\item $u_{ij}^s=p_{ij}$,
\end{enumerate}
with Woronowicz algebra maps $\Delta,\varepsilon,S$ constructed by universality.
\end{definition}

Here we allow the value $s=\infty$, with the convention that the last axiom simply disappears in this case. Observe that at $s<\infty$ the normality condition is actually redundant. This is because a partial isometry $a$ subject to the relation $aa^*=a^s$ is normal. As a first result now, making the connection with $H_N^s$, we have:

\begin{proposition}
We have an inclusion of quantum groups
$$H_N^s\subset H_N^{s+}$$
which is a liberation, in the sense that the classical version of $H_N^{s+}$, obtained by dividing by the commutator ideal, is the group $H_N^s$.
\end{proposition}

\begin{proof}
This follows as for $O_N\subset O_N^+$ or for $S_N\subset S_N^+$, by using the Gelfand theorem, applied to the quotient of $C(H_N^{s+})$ by its commutator ideal.
\end{proof}

In analogy with the results from the real case, we have the following result:

\begin{proposition}
The algebras $C(H_N^{s+})$ with $s=1,2,\infty$, and their presentation relations in terms of the entries of the matrix $u=(u_{ij})$, are as follows:
\begin{enumerate}
\item For $C(H_N^{1+})=C(S_N^+)$, the matrix $u$ is magic: all its entries are projections, summing up to $1$ on each row and column.

\item For $C(H_N^{2+})=C(H_N^+)$ the matrix $u$ is cubic: it is orthogonal, and the products of pairs of distinct entries on the same row or the same column vanish.

\item For $C(H_N^{\infty+})=C(K_N^+)$ the matrix $u$ is unitary, its transpose is unitary, and all its entries are normal partial isometries.
\end{enumerate}
\end{proposition}

\begin{proof}
This is something elementary, from \cite{bb+}, \cite{bv1}, the idea being as follows:

\medskip

(1) This follows from definitions and from standard operator algebra tricks.

\medskip

(2) This follows as well from definitions and standard operator algebra tricks.

\medskip

(3) This is just a translation of the definition of $C(H_N^{s+})$, at $s=\infty$.
\end{proof}

Let us prove now that $H_N^{s+}$ with $s<\infty$ is a quantum permutation group. For this purpose, we must change the fundamental representation. Let us start with:

\begin{definition}
A $(s,N)$-sudoku matrix is a magic unitary of size $sN$, of the form
$$m=\begin{pmatrix}
a^0&a^1&\ldots&a^{s-1}\\
a^{s-1}&a^0&\ldots&a^{s-2}\\
\vdots&\vdots&&\vdots\\
a^1&a^2&\ldots&a^0
\end{pmatrix}$$
where $a^0,\ldots,a^{s-1}$ are $N\times N$ matrices.
\end{definition}

The basic examples of such matrices come from the group $H_n^s$. Indeed, with $w=e^{2\pi i/s}$, each of the $N^2$ matrix coordinates $u_{ij}:H_N^s\to\mathbb C$ takes values in the following set:
$$S=\{0\}\cup\{1,w,\ldots,w^{s-1}\}$$

Thus, this coordinate function $u_{ij}:H_N^s\to\mathbb C$ decomposes as follows:
$$u_{ij}=\sum_{r=0}^{s-1}w^ra^r_{ij}$$

Here each $a^r_{ij}$ is a function taking values in $\{0,1\}$, and so a projection in the $C^*$-algebra sense, and it follows from definitions that these projections form a sudoku matrix. With this notion in hand, we have the following result, from \cite{bv1}:

\begin{theorem}
The following happen:
\begin{enumerate}
\item The algebra $C(H_N^s)$ is isomorphic to the universal commutative $C^*$-algebra generated by the entries of a $(s,N)$-sudoku matrix.

\item The algebra $C(H_N^{s+})$ is isomorphic to the universal $C^*$-algebra generated by the entries of a $(s,N)$-sudoku matrix.
\end{enumerate}
\end{theorem}

\begin{proof}
The first assertion follows from the second one, via Proposition 11.13. In order to prove the second assertion, consider the universal algebra in the statement:
$$A=C^*\left(a_{ij}^p\ \Big\vert \left(a^{q-p}_{ij}\right)_{pi,qj}=(s,N)-\mbox{sudoku }\right)$$

Consider also the algebra $C(H_N^{s+})$. According to Definition 11.12, this is presented by certain relations $R$, that we will call here level $s$ cubic conditions:
$$C(H_N^{s+})=C^*\left(u_{ij}\ \Big\vert\  u=N\times N\mbox{ level $s$ cubic }\right)$$

We will construct a pair of inverse morphisms between these algebras.

\medskip

(1) Our first claim is that $U_{ij}=\sum_pw^{-p}a^p_{ij}$ is a level $s$ cubic unitary. Indeed, by using the sudoku condition, the verification of (1-4) in Definition 11.12 is routine.

\medskip

(2) Our second claim is that the elements $A^p_{ij}=\frac{1}{s}\sum_rw^{rp}u^r_{ij}$, with the convention $u_{ij}^0=p_{ij}$, form a level $s$ sudoku unitary. Once again, the proof here is routine.

\medskip

(3) According to the above, we can define a morphism $\Phi:C(H_N^{s+})\to A$ by the formula $\Phi(u_{ij})=U_{ij}$, and a morphism $\Psi:A\to C(H_N^{s+})$ by the formula $\Psi(a^p_{ij})=A^p_{ij}$.

\medskip

(4) We check now the fact that $\Phi,\Psi$ are indeed inverse morphisms:
\begin{eqnarray*}
\Psi\Phi(u_{ij})
&=&\sum_pw^{-p}A^p_{ij}\\
&=&\frac{1}{s}\sum_pw^{-p}\sum_rw^{rp}u_{ij}^r\\
&=&\frac{1}{s}\sum_{pr}w^{(r-1)p}u_{ij}^r\\
&=&u_{ij}
\end{eqnarray*}

As for the other composition, we have the following computation:
\begin{eqnarray*}
\Phi\Psi(a^p_{ij})
&=&\frac{1}{s}\sum_rw^{rp}U_{ij}^r\\
&=&\frac{1}{s}\sum_rw^{rp}\sum_qw^{-rq}a_{ij}^q\\
&=&\frac{1}{s}\sum_qa_{ij}^q\sum_rw^{r(p-q)}\\
&=&a^p_{ij}
\end{eqnarray*}

Thus we have an isomorphism $C(H_N^{s+})=A$, as claimed.
\end{proof}

We will need the following simple fact: 

\begin{proposition}
A $sN\times sN$ magic unitary commutes with the matrix
$$\Sigma=
\begin{pmatrix}
0&I_N&0&\ldots&0\\
0&0&I_N&\ldots&0\\
\vdots&\vdots&&\ddots&\\
0&0&0&\ldots&I_N\\
I_N&0&0&\ldots&0
\end{pmatrix}$$
if and only if it is a sudoku matrix in the sense of Definition 11.15.
\end{proposition}

\begin{proof}
This follows from the fact that commutation with $\Sigma$ means that the matrix is circulant. Thus, we obtain the sudoku relations from Definition 11.15.
\end{proof}

Now let $Z_s$ be the oriented cycle with $s$ vertices, and consider the graph $NZ_s$ consisting of $N$ disjoint copies of it. Observe that, with a suitable labeling of the vertices, the adjacency matrix of this graph is the above matrix $\Sigma$. We obtain from this:

\begin{theorem}
We have the following results:
\begin{enumerate}
\item $H_N^s$ is the symmetry group of $NZ_s$.

\item $H_N^{s+}$ is the quantum symmetry group of $NZ_s$.
\end{enumerate}
\end{theorem}

\begin{proof}
This is something elementary, the idea being as follows:

\medskip

(1) This follows from definitions.

\medskip

(2) This follows from Theorem 11.16 and Proposition 11.17, because the algebra $C(H_N^{s+})$ is the quotient of the algebra $C(S_{sN}^+)$ by the relations making the fundamental corepresentation commute with the adjacency matrix of $NZ_s$.
\end{proof}

Next in line, we must talk about wreath products. We have here:

\begin{theorem}
We have isomorphisms as follows,
$$H_N^s=\mathbb Z_s\wr S_N\quad,\quad H_N^{s+}=\mathbb Z_s\wr_*S_N^+$$
with $\wr$ being a wreath product, and $\wr_*$ being a free wreath product.
\end{theorem}

\begin{proof}
This follows from the following formulae, valid for any connected graph $X$, and explained before in this book, applied to the graph $Z_s$:
$$G(NX)=G(X)\wr S_N\quad,\quad 
G^+(NX)=G^+(X)\wr_*S_N^+$$

Alternatively, (1) follows from definitions, and (2) can be proved directly, by constructing a pair of inverse morphisms. For details here, we refer to \cite{bv1}.
\end{proof}

Regarding now the easiness property of $H_N^s,H_N^{s+}$, we already know that this happens at $s=1,2$. The point is that this happens at $s=\infty$ too, the result being as follows:

\begin{theorem}
The quantum groups $K_N,K_N^+$ are easy, the corresponding categories
$$\mathcal P_{even}\subset P\quad,\quad 
\mathcal{NC}_{even}\subset NC$$
consisting of the partitions satisfying $\#\circ=\#\bullet$, as a weighted equality, in each block.
\end{theorem}

\begin{proof}
This is something which is routine, and we refer to \cite{bb+}.
\end{proof}

More generally now, we have the following result, from \cite{bb+}:

\begin{theorem}
The quantum groups $H_N^s,H_N^{s+}$ are easy, the corresponding categories
$$P^s\subset P\quad,\quad 
NC^s\subset NC$$
consisting of partitions satisfying $\#\circ=\#\bullet(s)$, as a weighted sum, in each block.
\end{theorem}

\begin{proof}
Observe that the result holds at $s=1$, trivially, then at $s=2$ as well, where our condition is equivalent to $\#\circ=\#\bullet(2)$ in each block, as found before, and finally at $s=\infty$ too, as explained in Theorem 11.20. In general, this follows as in the case of $H_N,H_N^+$, by using the one-block partition in $P(s,s)$. See \cite{bb+}.
\end{proof}

\section*{11c. Representation theory}

Let us discuss now the representation theory of $H_N^{s+}$. For this purpose, let us go back to the elements $u_{ij},p_{ij}$ in Definition 11.12. We recall from Proposition 11.14 that the matrix $p=(p_{ij})$ is a magic unitary. We first have the following result:

\begin{proposition}
The elements $u_{ij}$ and $p_{ij}$ satisfy:
\begin{enumerate}
\item $p_{ij}u_{ij}=u_{ij}$.

\item $u_{ij}^*=u_{ij}^{s-1}$.

\item $u_{ij}u_{ik}=0$ for $j\neq k$.
\end{enumerate}
\end{proposition}

\begin{proof}
We use the fact that in a $C^*$-algebra, $aa^*=0$ implies $a=0$.

\medskip

(1) This follows from the following computation, with $a=(p_{ij}-1)u_{ij}$:
$$aa^*
=(p_{ij}-1)p_{ij}(p_{ij}-1)
=0$$

(2) With $a=u_{ij}^*-u_{ij}^{s-1}$ we have $aa^*=0$, which gives the result.

\medskip

(3) With $a=u_{ij}u_{ik}$ we have $aa^*=0$, which gives the result.
\end{proof}

In what follows, we make the convention $u_{ij}^0=p_{ij}$. We have then:

\begin{proposition}
The algebra $C(H_N^{s+})$ has a family of $N$-dimensional corepresentations $\{u_k|k\in\mathbb Z\}$, satisfying the following conditions:
\begin{enumerate}
\item $u_k=(u_{ij}^k)$ for any $k\geq 0$.

\item $u_k=u_{k+s}$ for any $k\in\mathbb Z$.

\item $\bar{u}_k=u_{-k}$ for any $k\in\mathbb Z$.
\end{enumerate}
\end{proposition}

\begin{proof}
This is something elementary, the idea being as follows:

\medskip

(1) Let us set $u_k=(u_{ij}^k)$. By using Proposition 11.22 (3), we have:
$$\Delta(u_{ij}^k)
=\sum_{l_1\ldots l_k}u_{il_1}\ldots u_{il_k}\otimes u_{l_1j}\ldots u_{l_kj}
=\sum_lu_{il}^k\otimes u_{lj}^k$$

We have as well, trivially, the following two formulae:
$$\varepsilon(u_{ij}^k)=\delta_{ij}\quad,\quad 
S(u_{ij}^k)=u_{ji}^{*k}$$

(2) This follows once again from Proposition 11.22 (3), as follows:
$$u_{ij}^{k+s}
=u_{ij}^ku_{ij}^s
=u_{ij}^kp_{ij}
=u_{ij}^k$$

(3) This follows from Proposition 11.22 (2), and we are done.
\end{proof}

Let us compute now the intertwiners between the various tensor products between the above corepresentations $u_i$. For this purpose, we make the assumption $N\geq 4$, which brings linear independence. In order to simplify the notations, we will use:

\begin{definition}
For $i_1,\ldots,i_k\in\mathbb Z$ we use the notation
$$u_{i_1\ldots i_k}=u_{i_1}\otimes\ldots\otimes u_{i_k}$$
where $\{u_i|i\in\mathbb Z\}$ are the corepresentations in Proposition 11.23.
\end{definition}

Observe that in the particular case $i_1,\ldots,i_k\in\{\pm 1\}$, we obtain in this way all the possible tensor products between $u=u_1$ and $\bar{u}=u_{-1}$, known by \cite{wo1} to contain any irreducible corepresentation of $C(H_N^{s+})$. Here is now our main result:

\begin{theorem}
We have the following equality of linear spaces,
$$Hom(u_{i_1\ldots i_k},u_{j_1\ldots j_l})=span\left\{T_p\Big|p\in NC_s(i_1\ldots i_k,j_1\ldots j_l)\right\}$$
where the set on the right consists of elements of $NC(k,l)$ having the property that in each block, the sum of $i$ indices equals the sum of $j$ indices, modulo $s$.
\end{theorem}

\begin{proof}
This result is from \cite{bv1}, the idea of the proof being as follows:

\medskip

(1) Our first claim is that, in order to prove $\supset$, we may restrict attention to the case $k=0$. This follow indeed from the Frobenius duality isomorphism.

\medskip

(2) Our second claim is that, in order to prove $\supset$ in the case $k=0$, we may restrict attention to the one-block partitions. Indeed, this follows once again from a standard trick. Consider the following disjoint union:
$$NC_s=\bigcup_{k=0}^\infty\bigcup_{i_1\ldots i_k} NC_s(0,i_1\ldots i_k)$$

This is a set of labeled partitions, having property that each $p\in NC_s$ is noncrossing, and that for $p\in NC_s$, any block of $p$ is in $NC_s$. But it is well-known that under these assumptions, the global algebraic properties of $NC_s$ can be checked on blocks.

\medskip

(3) Proof of $\supset$. According to the above considerations, we just have to prove that the vector associated to the one-block partition in $NC(l)$ is fixed by $u_{j_1\ldots j_l}$, when:
$$s|j_1+\ldots+j_l$$

Consider the standard generators $e_{ab}\in M_N(\mathbb C)$, acting on the basis vectors by:
$$e_{ab}(e_c)=\delta_{bc}e_a$$

The corepresentation $u_{j_1\ldots j_l}$ is given by the following formula:
$$u_{j_1\ldots j_l}=\sum_{a_1\ldots a_l}\sum_{b_1\ldots b_l}u_{a_1b_1}^{j_1}\ldots u_{a_lb_l}^{j_l}\otimes e_{a_1b_1}\otimes\ldots\otimes e_{a_lb_l}$$

As for the vector associated to the one-block partition, this is:
$$\xi_l=\sum_be_b^{\otimes l}$$

By using now several times the relations in Proposition 11.22, we obtain, as claimed: 
\begin{eqnarray*}
u_{j_1\ldots j_l}(1\otimes\xi_l)
&=&\sum_{a_1\ldots a_l}\sum_bu_{a_1b}^{j_1}\ldots u_{a_lb}^{j_l}\otimes e_{a_1}\otimes\ldots\otimes e_{a_l}\\
&=&\sum_{ab}u_{ab}^{j_1+\ldots+j_l}\otimes e_a^{\otimes l}\\
&=&1\otimes\xi_l
\end{eqnarray*}

(4) Proof of $\subset$. The spaces in the statement form a Tannakian category, so they correspond to a Woronowicz algebra $A$, coming with corepresentations $\{v_i\}$, such that:
$$Hom(v_{i_1\ldots i_k},v_{j_1\ldots j_l})={\rm span}\left\{T_p\Big|p\in NC_s(i_1\ldots i_k,j_1\ldots j_l)\right\}$$

On the other hand, the inclusion $\supset$ that we just proved shows that $C(H_N^{s+})$ is a model for the category. Thus we have a quotient map as follows:
$$A\to C(H_N^{s+})\quad,\quad 
v_i\to u_i$$

But this latter map can be shown to be an isomorphism, by suitably adapting the proof from the $s=1$ case, for the quantum permutation group $S_N^+$. See \cite{bb+}, \cite{bv1}.
\end{proof}

As an illustration for the above result, we have the following statement:

\begin{proposition}
The basic corepresentations $u_0,\ldots,u_{s-1}$ are as follows:
\begin{enumerate}
\item $u_1,\ldots,u_{s-1}$ are irreducible.

\item $u_0=1+r_0$, with $r_0$ irreducible.

\item $r_0,u_1,\ldots,u_{s-1}$ are distinct.
\end{enumerate}
\end{proposition}

\begin{proof}
We apply Theorem 11.25 with $k=l=1$ and $i_1=i,j_1=j$. This gives:
$$\dim (Hom(u_i,u_j))=\# NC_s(i,j)$$

We have two candidates for the elements of $NC_s(i,j)$, namely the two partitions in $NC(1,1)$. So, consider these two partitions, with the points labeled by $i,j$:
$$p=\left\{\begin{matrix}i\cr\Big\vert\cr j\end{matrix}\right\}\qquad\ \qquad
q=\left\{\begin{matrix}i\cr|\cr\cr |\cr j\end{matrix}\right\}$$

We have to check for each of these partitions if the sum of $i$ indices equals or not the sum of $j$ indices, modulo $s$, in each block. The answer is as follows:
\begin{eqnarray*}
p\in NC_s(i,j)&\iff&i=j\\
q\in NC_s(i,j)&\iff&i=j=0
\end{eqnarray*}

By collecting together these two answers, we obtain:
$$\# NC_s(i,j)=
\begin{cases}
0&{\rm if\ }i\neq j\\
1&{\rm if\ }i=j\neq 0\\
2&{\rm if\ }i=j=0
\end{cases}$$

We can now prove the various assertions, as follows:

\medskip

(1) This follows from the second equality.

\medskip

(2) This follows from the third equality and from the fact that we have $1\in u_s$.

\medskip

(3) This follows from the first equality.
\end{proof}

We can now compute the fusion rules for $H_N^{s+}$. The result, from \cite{bv1}, is as follows:

\begin{theorem}
Let $F=<\mathbb Z_s>$ be the set of words over $\mathbb Z_s$, with involution given by $(i_1\ldots i_k)^-=(-i_k)\ldots(-i_1)$, and with fusion product given by:
$$(i_1\ldots i_k)\cdot (j_1\ldots j_l)=i_1\ldots i_{k-1}(i_k+j_1)j_2\ldots j_l$$
The irreducible representations of $H_N^{s+}$ can then be labeled $r_x$ with $x\in F$, such that the involution and fusion rules are $\bar{r}_x=r_{\bar{x}}$ and
$$r_x\otimes r_y=\sum_{x=vz,y=\bar{z}w}r_{vw}+r_{v\cdot w}$$
and such that we have $r_i=u_i-\delta_{i0}1$ for any $i\in\mathbb Z_s$.
\end{theorem}

\begin{proof}
This basically follows from Theorem 11.25, the idea being as follows:

\medskip

(1) Consider the monoid $A=\{a_x|x\in F\}$, with multiplication $a_xa_y=a_{xy}$. We denote by $\mathbb NA$ the set of linear combinations of elements in $A$, with coefficients in $\mathbb N$, and we endow it with fusion rules as in the statement:
$$a_x\otimes a_y=\sum_{x=vz,y=\bar{z}w}a_{vw}+a_{v\cdot w}$$

With these notations, $(\mathbb NA,+,\otimes)$ is a semiring. We will use as well the set $\mathbb ZA$, formed by the linear combinations of elements of $A$, with coefficients in $\mathbb Z$. The above tensor product operation extends to $\mathbb ZA$, and $(\mathbb ZA,+,\otimes)$ is a ring.

\medskip

(2) Our claim is that the fusion rules on $\mathbb ZA$ can be uniquely described by conversion formulae as follows, with $C$ being positive integers, and $D$ being integers:
$$a_{i_1}\otimes\ldots\otimes a_{i_k}=\sum_l\sum_{j_1\ldots j_l}C_{i_1\ldots i_k}^{j_1\ldots j_l}a_{j_1\ldots j_l}$$
$$a_{i_1\ldots i_k}=\sum_l\sum_{j_1\ldots j_l}D_{i_1\ldots i_k}^{j_1\ldots j_l}a_{j_1}\otimes\ldots\otimes a_{j_l}$$

Indeed, the existence and uniqueness of such decompositions follow from the definition of the tensor product operation, and by recurrence over $k$ for the $D$ coefficients.

\medskip

(3) Our claim is that there is a unique morphism of rings $\Phi:\mathbb ZA\to R$, such that $\Phi(a_i)=r_i$ for any $i$. Indeed, consider the following elements of $R$:
$$r_{i_1\ldots i_k}=\sum_l\sum_{j_1\ldots j_l}D_{i_1\ldots i_k}^{j_1\ldots j_l}r_{j_1}\otimes\ldots\otimes r_{j_l}$$

In case we have a morphism as claimed, we must have $\Phi(a_x)=r_x$ for any $x\in F$. Thus our morphism is uniquely determined on  $A$, so it is uniquely determined on $\mathbb ZA$. In order to prove now the existence, we can set $\Phi(a_x)=r_x$ for any $x\in F$, then extend $\Phi$ by linearity to the whole $\mathbb ZA$. Since $\Phi$ commutes with the above conversion formulae, which describe the fusion rules, it is indeed a morphism. 

\medskip

(4) Our claim is that $\Phi$ commutes with the linear forms $x\to\#(1\in x)$. Indeed, by linearity we just have to check the following equality:
$$\#(1\in a_{i_1}\otimes\ldots\otimes a_{i_k})=\#(1\in r_{i_1}\otimes\ldots\otimes r_{i_k})$$

Now remember that the elements $r_i$ are defined as $r_i=u_i-\delta_{i0}1$. So, consider the elements $c_i=a_i+\delta_{i0}1$. Since the operations $r_i\to u_i$ and $a_i\to c_i$ are of the same nature, by linearity the above formula is equivalent to:
$$\#(1\in c_{i_1}\otimes\ldots\otimes c_{i_k})=\#(1\in u_{i_1}\otimes\ldots\otimes u_{i_k})$$

Now by using Theorem 6.19, what we have to prove is:
$$\#(1\in c_{i_1}\otimes\ldots\otimes c_{i_k})=\#NC_s(i_1\ldots i_k)$$

In order to prove this formula, consider the product on the left:
$$P=(a_{i_1}+\delta_{i_10}1)\otimes(a_{i_2}+\delta_{i_20}1)\otimes\ldots\otimes (a_{i_k}+\delta_{i_k0}1)$$

This quantity can be computed by using the fusion rules on $A$. A recurrence on $k$ shows that the final components of type $a_x$ will come from the different ways of grouping and summing the consecutive terms of the sequence $(i_1,\ldots,i_k)$, and removing some of the sums which vanish modulo $s$, as to obtain the sequence $x$. But this can be encoded by families of noncrossing partitions, and in particular the 1 components will come from the partitions in $NC_s(i_1\ldots i_k)$. Thus $\#(1\in P)=\# NC_s(i_1\ldots i_k)$, as claimed.

\medskip

(5) Our claim now is that $\Phi$ is injective. Indeed, this follows from the result in the previous step, by using a standard positivity argument, namely:
\begin{eqnarray*}
\Phi(\alpha)=0
&\implies&\Phi(\alpha\alpha^*)=0\\
&\implies&\#(1\in \Phi(\alpha\alpha^*))=0\\
&\implies&\#(1\in \alpha\alpha^*)=0\\
&\implies&\alpha=0
\end{eqnarray*}

Here $\alpha$ is arbitrary in the domain of $\Phi$, we use the notation $a_x^*=a_{\bar{x}}$, where $a\to\#(1,a)$ is the unique linear extension of the operation consisting of counting the number of 1's.  Observe that this latter linear form is indeed positive definite, according to the identity $\#(1,a_xa_y^*)=\delta_{xy}$, which is clear from the definition of the product of $\mathbb ZA$.

\medskip

(6) Our claim is that $\Phi(A)\subset R_{irr}$. This is the same as saying that $r_x\in R_{irr}$ for any $x\in F$, and we will prove it by recurrence. Assume that the assertion is true for all the words of length $<k$, and consider an arbitrary length $k$ word, $x=i_1\ldots i_k$. We have:
$$a_{i_1}\otimes a_{i_2\ldots i_k}=a_x+a_{i_1+i_2,i_3\ldots i_k}+\delta_{i_1+i_2,0}a_{i_3\ldots i_k}$$

By applying $\Phi$ to this decomposition, we obtain:
$$r_{i_1}\otimes r_{i_2\ldots i_k}=r_x+r_{i_1+i_2,i_3\ldots i_k}+\delta_{i_1+i_2,0}r_{i_3\ldots i_k}$$

We have the following computation, which is valid for $y=i_1+i_2,i_3\ldots i_k$, as well as for $y=i_3\ldots i_k$ in the case $i_1+i_2=0$:
\begin{eqnarray*}
\#(r_y\in r_{i_1}\otimes r_{i_2\ldots i_k})
&=&\#(1,r_{\bar{y}}\otimes r_{i_1}\otimes r_{i_2\ldots i_k})\\
&=&\#(1,a_{\bar{y}}\otimes a_{i_1}\otimes a_{i_2\ldots i_k})\\
&=&\#(a_y\in a_{i_1}\otimes a_{i_2\ldots i_k})\\
&=&1  
\end{eqnarray*}

Moreover, we know from the previous step that we have $r_{i_1+i_2,i_3\ldots i_k}\neq r_{i_3\ldots i_k}$, so we conclude that the following formula defines an element of $R^+$:
$$\alpha=r_{i_1}\otimes r_{i_2\ldots i_k}-r_{i_1+i_2,i_3\ldots i_k}-\delta_{i_1+i_2,0}r_{i_3\ldots i_k}$$

On the other hand, we have $\alpha=r_x$, so we conclude that we have $r_x\in R^+$. Finally, the irreducibility of $r_x$ follows from the following computation:
\begin{eqnarray*}
\#(1\in r_x\otimes\bar{r}_x)
&=&\#(1\in r_x\otimes r_{\bar{x}})\\
&=&\#(1\in a_x\otimes a_{\bar{x}})\\
&=&\#(1\in a_x\otimes\bar{a}_x)\\
&=&1
\end{eqnarray*}

(7) Summarizing, we have constructed an injective ring morphism, as follows:
$$\Phi:\mathbb ZA\to R\quad,\quad 
\Phi(A)\subset R_{irr}$$

The remaining fact to be proved, namely that we have $\Phi(A)=R_{irr}$, is clear from the general results in \cite{wo1}. Indeed, since each element of $\mathbb NA$ is a sum of elements in $A$, by applying $\Phi$ we get that each element in $\Phi(\mathbb NA)$ is a sum of irreducible corepresentations in $\Phi(A)$. But since $\Phi(\mathbb NA)$ contains all the tensor powers between the fundamental corepresentation and its conjugate, we get $\Phi(A)=R_{irr}$, and we are done.
\end{proof}

Still following \cite{bv1}, let us present now a useful formulation of Theorem 11.27. We begin with a slight modification of Theorem 11.27, as follows:

\begin{theorem}
Consider the free monoid $A=<a_i|i\in\mathbb Z_s>$ with the involution $a_i^*=a_{-i}$, and define inductively the following fusion rules on it:
$$pa_i\otimes a_jq=pa_ia_jq+pa_{i+j}q+\delta_{i+j,0}p\otimes q$$
Then the irreducible representations of $H_N^{s+}$ can be indexed by the elements of $A$, and the fusion rules and involution are the above ones.
\end{theorem}

\begin{proof}
Our claim is that this follows from Theorem 11.27, by performing the following relabeling of the irreducible corepresentations:
$$r_{i_1\ldots i_k}\to a_{i_1}\ldots a_{i_k}$$

Indeed, with the notations in Theorem 11.27 we have the following computation, valid for any two elements $i,j\in\mathbb Z_s$ and any two words $x,y\in F$:
\begin{eqnarray*}
r_{xi}\otimes r_{jy}
&=&\sum_{xi=vz,jy=\bar{z}w}r_{vw}+r_{v\cdot w}\\
&=&r_{xijy}+r_{x,i+j,y}+\delta_{i+j,0}\sum_{x=vz,y=\bar{z}w}r_{vw}+r_{v\cdot w}\\
&=&r_{xijy}+r_{x,i+j,y}+\delta_{i+j,0}r_x\otimes r_y
\end{eqnarray*}

With the above relabeling $r_{i_1\ldots i_k}\to a_{i_1}\ldots a_{i_k}$, this gives the formula in the statement (with $r_x\to p$ and $r_y\to q$), and we are done.
\end{proof}

Based on the above, we have a second reformulation as well, as follows:

\begin{theorem}
Consider the monoid $M=<a,z|z^s=1>$ with the involution $a^*=a,z^*=z^{-1}$, and define inductively the following fusion rules on it:
$$vaz^i\otimes z^jaw=vaz^{i+j}aw+\delta_{s|i+j}v\otimes w$$
Then the irreducible representations of $H_N^{s+}$ can be indexed by the elements of the monoid $N=<aza>$, and the fusion rules and involution are the above ones.
\end{theorem}

\begin{proof}
It is routine to check that the elements $az^ia$ with $i=1,\ldots,s$ are free inside $M$. In other words, the submonoid $N'=<az^ia>$ is free on $s$ generators, so it can be identified with the free monoid $A$ in Theorem 11.28, via $a_i=az^ia$. We have $(az^ia)^*=az^{-i}a$, so this identification is involution-preserving. Consider now two arbitrary elements $p,q\in N'$. By using twice the formula in the statement, we obtain:
\begin{eqnarray*}
pa_i\otimes a_jq
&=&paz^ia\otimes az^jaq\\
&=&paz^iaaz^jaq+paz^i\otimes z^jaq\\
&=&paz^iaaz^jaq+paz^{i+j}aq+\delta_{i+j,0}p\otimes q\\
&=&pa_ia_jq+pa_{i+j}q+\delta_{i+j,0}p\otimes q
\end{eqnarray*}

Thus our identification $N'\simeq A$ is fusion rule-preserving. In order to conclude, it remains to prove that the inclusion $N\subset N'$ is actually an equality. But this follows from the fact that $A$ is generated as a fusion monoid by $a_1$. Indeed, by using the identification $N'\simeq A$ this shows that $N'$ is generated as a fusion monoid by $aza$, and we are done. 
\end{proof}

We refer to \cite{bv1} and related papers, including \cite{fwe}, for more on the above, including for some further useful technical reformulations of Theorem 11.27.

\section*{11d. Complexification}

With the above done, we have to face now the main problem that we have, the one formulated in the beginning of the present chapter, namely that of coming up with some further examples of intermediate easy quantum groups, as follows:
$$S_N\subset G\subset K_N^+$$

However, this does not look obvious at all, because the world of such quantum groups is quite wild, a bit in analogy with the world of the complex reflection groups. An idea here would be to first discuss the simplest case, which is the ``purely complex'' one:
$$K_N\subset G\subset K_N^+$$

But this does not look obvious either, and in short it seems like we are stuck with some difficult mathematics, and time for a tactical retreat. This being said, let us ask the cat, who is a world-class expert in tactical retreats. And cat says:

\begin{cat}
Yes don't punch above your weight, but have at least some $K_N^\Gamma$ and $K_N^{\diamond r}$ beasts constructed, by whatever complexification method of your choice.
\end{cat}

Thanks cat, this looks like some wise advice, so let us have at least some complex versions of the constructions of $H_N^\Gamma$ and $H_N^{\diamond r}$ from chapters 9-10 done, always good to have this, and leave the tricky further examples, and classification results, to future generations. With the remark that, who knows, maybe when looking for applications and everything, the beasts of type $K_N^\Gamma$ and $K_N^{\diamond r}$ might be enough. Or at least when looking for applications at our weight class. It's all about weight, in life, isn't it, sweet kit-kat.

\bigskip

Getting started now, we already know what ``real version'' and ``complexification'' should mean, in the quantum group context. First, we have following definition:

\begin{definition}
We can talk about real versions of quantum groups, as follows:
\begin{enumerate}
\item The real version of an easy quantum group $S_N\subset G\subset U_N^+$ is the easy quantum group $S_N\subset G_{real}\subset O_N^+$ given by $G_{real}=G\cap O_N^+$. 

\item Equivalently, if $G$ comes from a category of partitions $D\subset P$, then $G_{real}$ comes from the category of partitions $<D,NC_2>$.
\end{enumerate}
\end{definition}

Observe that the operation in (1) is well-defined for any closed subgroup $G\subset U_N^+$, producing a certain closed subgroup $G_{real}\subset O_N^+$, but in what follows we will only need this in the easy case. As for the equivalence with (2), in the easy case, this comes from our general results from chapter 3 regarding the intersection operation $\cap$.

\bigskip

Getting now to complexification, we have here a similar definition, as follows:

\begin{definition}
We can talk about quantum group complexification, as follows:
\begin{enumerate}
\item The complexification of an easy quantum group $S_N\subset G\subset O_N^+$ is the easy quantum group $K_N\subset G_{comp}\subset U_N^+$ given by $G_{comp}=
\{G,K_N\}$. 

\item Equivalently, if $G$ comes from a category of partitions $D\subset P$, then $G_{comp}$ comes from the category of partitions $D\cap\mathcal P_{even}$.
\end{enumerate}
\end{definition}

As before with Definition 11.31, there are several comments to be made here, some being trivial, and some more being subtle, the idea being as follows:

\bigskip

-- First, the operation in (1) can be performed in fact for any easy quantum group $S_N\subset G\subset U_N^+$, but this extension is without much interest in the non-real case, $G\not\subset O_N^+$, because our main examples here tend to contain $K_N$, anyway. 

\bigskip

-- As a more subtle remark now, we have a version of the operation in (1) obtained by using the plain generation operation, $G_{comp}'=<G,K_N>$, and which works for any closed subgroup $G\subset U_N^+$. However, as explained in chapter 3, we have $\{\,,\}\neq<\,,>$ in general.

\bigskip

-- Finally, regarding the equivalence between our operations in (1) and (2), this comes from our results in chapter 3 regarding the easy generation operation $\{\,,\}$. In fact, this operation $\{\,,\}$ was defined there precisely via $\cap$ at the level of categories of partitions.

\bigskip

All this is nice, and as a first task for us, we would like to know to which extent the operations in Definition 11.31 and Definition 11.32 are inverse to each other. However, this does not look exactly obvious, due to a variety of technical reasons. So, stuck again, and time again to ask the cat, who fortunately is still around. And cat says:

\begin{cat}
Dude I told you, don't punch above your weight. Just have some $K_N^\Gamma$ and $K_N^{\diamond r}$ beasts quickly constructed, and then go ahead with chapter 12.
\end{cat}

Thanks cat. So if I understand well I should look for an alternative conceptual way of complexifying the compact quantum groups, with a clear definition, and some nice mathematical theory about this, then define $K_N^\Gamma$ and $K_N^{\diamond r}$ as you say, and done.

\bigskip

And alternative methods, fortunately, do exist. We have for instance the free complexification operation, which works well in a number of important cases, as we know from the previous chapters, and whose definition in general is as follows:

\begin{definition}
The free complexification of a closed subgroup $G\subset U_N^+$, with fundamental corepresentation $u$, is the closed subgroup $\widetilde{G}\subset U_N^+$ given by
$$C(\widetilde{G})=<zu_{ij}>\subset C(\mathbb T)*C(G)$$
with fundamental corepresentation $\widetilde{u}=zu$, where $z$ is the standard generator of $C(\mathbb T)$.
\end{definition}

This sounds very nice, and we already know from the previous chapters that this works well in a number of cases. Indeed, we first have the following key result:
$$\widetilde{O_N^+}=U_N^+$$

At the level of the main intermediate liberations, again in the continuous case, things are nice too, because we know from chapter 8 that we have the following equalities:
$$\widetilde{O_N}=\widetilde{U_N}=\widetilde{O_N^*}=\widetilde{U_N^*}=U_N^\times$$

At the discrete level now, which is the one that we are interested in, in this chapter, again things work fine in the free case, where we have the following result:
$$\widetilde{H_N^+}=K_N^+$$

As for the intermediate liberations, in the discrete case, this remains to be worked out. However, before doing that, we have several questions to be solved, namely:

\begin{questions}
Regarding the free complexification, in the easy case:
\begin{enumerate}
\item Is it true that if $G$ is easy, then so is $\widetilde{G}$?

\item If $G$ is real, is it the real version of $\widetilde{G}$?

\item In fact, do we have $\widetilde{G}=G_{comp}$?
\end{enumerate}
\end{questions}

But, in what regards these questions, although all these look doable, with some work involved, none is trivial, and at the level of what is known, the situation is as follows:

\bigskip

(1) Here the answer is most likely yes, as a consequence of the results of Raum in \cite{rau}, who computed there the representation theory of $\widetilde{G}$, in terms of that of $G$, in general. But this is non-trivial, and still remains to be applied to the easy case. 

\bigskip

(2) This is most likely a rather delicate question, which seems to require a case-by-case analysis, and which perhaps comes after (3). In any case we can't expect a plain yes answer here, for instance because of $\widetilde{O}_N=U_N^\times$, which gives $(\widetilde{O}_N)_{real}=O_N^*$.

\bigskip

(3) Here there is indication from the above-mentioned work of Raum in \cite{rau} that the answer should be yes, at least under some suitable assumptions on our easy quantum group $G$, but all this is not exactly trivial, and still remains to be worked out.

\bigskip

Summarizing, all this not very good news, and we are again stuck, and I am afraid that I will have to ask again the cat. And cat says:

\begin{cat}
Define your objects first, and study them afterwards.
\end{cat}

Thanks cat. I think I eventually got your point, so how can we surround $K_N^\Gamma$ and $K_N^{\diamond r}$ by some nice general theories, without knowing what these objects are. So, time to fix this, definitely. Based on the above, and in the lack of something better, we have:

\begin{theorem}
We have easy quantum groups $K_N^\times$ as follows,
$$\xymatrix@R=50pt@C=50pt{
K_N\ar[r]&K_N^\Gamma\ar[r]&K_N^{[\infty]}\ar[r]&K_N^{\diamond r}\ar[r]&K_N^+\\
H_N\ar[u]\ar[r]&H_N^\Gamma\ar[r]\ar[u]&H_N^{[\infty]}\ar[r]\ar[u]&H_N^{\diamond r}\ar[r]\ar[u]&H_N^+\ar[u]
}$$
obtained by categorical complexification, $G\to G_{comp}$.
\end{theorem}

\begin{proof}
This is more of an empty statement, the idea being that we can perform to the quantum groups on the bottom the complexification construction $G\to G_{comp}$ from Definition 11.32, and up to a few functoriality checks, and some checks at the endpoints too, which are all elementary, we are led to the diagram in the statement.
\end{proof}

All this is quite nice, job done, at least we know one thing. There are of course many questions left, and we will be back to this in the next chapter.

\section*{11e. Exercises} 

This was a difficult chapter, and as an exercise here, of course difficult, we have:

\begin{exercise}
Classify the quantum reflection groups.
\end{exercise}

We will actually comment a bit more on this exercise, in the next chapter.

\chapter{The complex case}

\section*{12a. Liberation theory}

Welcome to advanced easiness. What we did so far in this book was rather standard material, known for some time, and relatively well understood. So, time now to get into really difficult questions, of research flavor. In order to explain the problems, let us go back to the standard cube formed by the main easy quantum groups, namely: 
$$
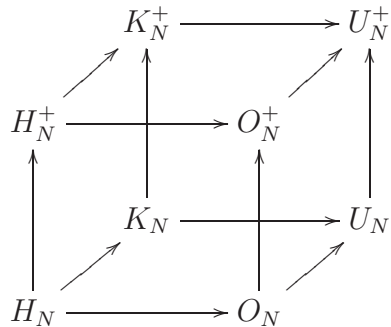
\xymatrix@R=18pt@C=18pt{
&K_N^+\ar[rr]&&U_N^+\\
H_N^+\ar[rr]\ar[ur]&&O_N^+\ar[ur]\\
&K_N\ar[rr]\ar[uu]&&U_N\ar[uu]\\
H_N\ar[uu]\ar[ur]\ar[rr]&&O_N\ar[uu]\ar[ur]
}$$

The question that we would mostly like to solve is the classification problem for the easy quantum groups inside the cube, $H_N\subset G\subset U_N^+$, which are called twistable. More generally, we would like to solve classification problem in the general easy case, $S_N\subset G\subset U_N^+$. And these questions are not trivial, the situation being as follows:

\bigskip

(1) A natural idea would be that of following the strategy from the real case, from chapters 9-10, which was successful, with a classification in the non-hyperoctahedral case, coupled with a classification in the hyperoctahedral case. However, this is something quite difficult, and for recent advances on this program, we refer to Mang-Weber \cite{mw3}, \cite{mw4}.

\bigskip

(2) A second idea, that we already met in the real case, and in other situations, and which appears as a modification of the Mang-Weber program, would be that of imposing, at least to start with, some extra conditions on our easy quantum groups. For instance, having the twistable, uniform case fully solved would be certainly a good thing.

\bigskip

(3) Finally, as a variation of what has been said above, we have the natural question of better understanding, to start with, the ``face to face'' correspondences in the above cube, and more specifically the ``liberation'', ``complexification'' and ``discretization'' procedures. And with this being something that we already met, on several occasions.

\bigskip

So, what to choose? So many things to talk about, and we are actually running out of time and space, because the present chapter 12 will be the end of our standard discussion on easiness, with Part IV coming afterwards being about something else.

\bigskip

Cat is gone hunting, but before leaving, a bit worried about me writing this book, at this point, he told me to be modest. So, we will choose something very modest for the present chapter 12, namely one-third of (3) above, the ``liberation'' question.

\bigskip

Getting started now, our results so far about the easy liberations of $O_N,U_N,H_N,K_N$ are an excellent input for the study of the general liberation problem, that we will study here, in the easy case, and in general, with the idea in mind of talking afterwards about the classification of the easy quantum groups $S_N\subset G\subset U_N^+$, and of more general such quantum groups. Let us start with something very general, as follows:

\begin{definition}
A liberation of a compact Lie group $G\subset U_N$ is a quantum group 
$$G\subset G^\times\subset U_N^+$$
whose classical version, $G^\times\cap U_N$, equals the group $G$ itself. 
\end{definition}

This is obviously a very general definition, which is of course something very natural. However, at this level of generality, nothing much can be said, or at least it is not known yet how to do this. Here are however a few basic remarks on the subject:

\begin{proposition}
The set of liberations of a given compact Lie group $G\subset U_N$ has the following properties:
\begin{enumerate}
\item It is stable under the intersection operation $\cap$.

\item It is not necessarily stable under the generation operation $<\,,>$.
\end{enumerate}
\end{proposition}

\begin{proof}
This is something elementary, the idea being as follows:

\medskip

(1) This is something trivial.

\medskip

(2) The result for the generation operation fails indeed, for instance for the hyperoctahedral group $H_N=\mathbb Z_2\wr S_N$. Indeed, $H_N$ has at least two main liberations, namely the twisted orthogonal group $O_N'$, which appears as quantum symmetry group of the hypercube in $\mathbb R^N$, and the quantum group $H_N^+=\mathbb Z_2\wr_*S_N^+$, which is the quantum symmetry group of the coordinate axes of $\mathbb R^N$. And the point is that we have:
$$<O_N',H_N^+>=U_N^+$$

Thus, (2) fails, and in a particularly bad way, for $H_N$.
\end{proof}

At a more constructive level now, one idea is that the liberations of a compact Lie group $G\subset U_N$ should appear via operations of type $G^\times=<G,I^\times>$, with $I^\times\subset U_N^+$ being a ``basic'' quantum group. In order to discuss this, let us start with:

\index{diagonal subgroup}
\index{soft liberation}
\index{hard liberation}

\begin{definition}
Given $H_N\subset G\subset U_N^+$, the diagonal tori $T=G\cap\mathbb T_N^+$ and reflection subgroups $K=G\cap K_N^+$ for $G$ and for $G_{class}=G\cap U_N$ form a diagram as follows:
$$\xymatrix@R=17mm@C=16mm{
T\ar[r]&K\ar[r]&G\\
T_{class}\ar[r]\ar[u]&K_{class}\ar[r]\ar[u]&G_{class}\ar[u]}$$
We say that $G$ appears as a soft/hard liberation when it is generated by $G_{class}$ and by $K/T$, which means that the right square/whole rectangle should be generation diagrams.
\end{definition}

Observe that hard liberation implies soft liberation, because the diagonal torus being included in the reflection group, $T\subset K$, we have the following implication:
$$T\subset K\implies <G_{class},T>\subset<G_{class},K>$$

Also, it is in fact possible to further complicate the above picture, by adding free versions as well, with these free versions being given by the following formula:
$$G_{free}=<G,S_N^+>$$ 

All this is quite technical, and as a concrete result in connection with the hard liberation notion, we have the following statement, regarding the basic unitary groups:

\begin{theorem}
The diagonal tori of the basic unitary quantum groups
$$\xymatrix@R=17mm@C=16mm{
U_N\ar[r]&U_N^*\ar[r]&U_N^+\\
O_N\ar[r]\ar[u]&O_N^*\ar[r]\ar[u]&O_N^+\ar[u]}$$
are as follows,
$$\xymatrix@R=17mm@C=16mm{
\mathbb T_N\ar[r]&\mathbb T_N^*\ar[r]&\mathbb T_N^+\\
T_N\ar[r]\ar[u]&T_N^*\ar[r]\ar[u]&T_N^+\ar[u]}$$
and these unitary quantum groups all appear via hard liberation.
\end{theorem}

\begin{proof}
The first assertion is something that we already know. As for the second assertion, this is something which is quite routine as well. We will be back to this.
\end{proof}

As an interesting remark now, our notion of hard liberation has its limitations, and some subtleties appear at the level of the quantum reflection groups, as follows:

\begin{theorem}
The diagonal tori of the basic quantum reflection groups
$$\xymatrix@R=17mm@C=16mm{
K_N\ar[r]&K_N^*\ar[r]&K_N^+\\
H_N\ar[r]\ar[u]&H_N^*\ar[r]\ar[u]&H_N^+\ar[u]}$$
are as follows,
$$\xymatrix@R=17mm@C=16mm{
\mathbb T_N\ar[r]&\mathbb T_N^*\ar[r]&\mathbb T_N^+\\
T_N\ar[r]\ar[u]&T_N^*\ar[r]\ar[u]&T_N^+\ar[u]}$$
and these quantum reflection groups do not all appear via hard liberation.
\end{theorem}

\begin{proof}
The first assertion is clear, as a consequence of Theorem 12.4, because the diagonal torus is the same for a quantum group, and for its reflection subgroup:
$$G\cap\mathbb T_N^+=(G\cap K_N^+)\cap\mathbb T_N^+$$

Regarding the second assertion, things are quite tricky here, as follows:

\medskip

(1) In the classical case the hard liberation property definitely holds, because any classical group is by definition a hard liberation of itself.

\medskip

(2) In the half-classical case the answer is again positive, and this can be proved by using the technology developed by Bichon and Dubois-Violette in \cite{bdu}.

\medskip

(3) In the free case the hard liberation property fails, due to the intermediate quantum groups $H_N^{[\infty]}$, $K_N^{[\infty]}$, where ``hard liberation stops''. We will be back to this.
\end{proof}

Summarizing, the notions of soft and hard liberation provide us with some answers, to the questions that we have. However, there are still many open questions regarding these operations, quite often in relation with the generation operation $<\,,>$.

\section*{12b. Generation results}

In order to further comment on the above questions, let us recall that following notion, that we studied in detail in chapter 6, in connection with various algebraic and analytic questions, and which plays a key role in connection with the notion of easiness:

\begin{definition}
A family $G=(G_N)$ with $G_N\subset U_N^+$ is called uniform when 
$$G_{N-1}=G_N\cap U_{N-1}^+$$
for any $N\geq 2$, with the embeddings $U_{N-1}^+\subset U_N^+$ being given by $u\to diag(u,1)$. 
\end{definition}

As a first remark, under this uniformity assumption, when assuming that $G_{N-1}$ is not classical, $G_N$ is not classical either. Thus, there is an integer $n\in\{2,3,\ldots,\infty\}$ such that $G_1,\ldots,G_{n-1}$ are all classical, and then $G_n,G_{n+1},\ldots$ are all non-classical. We have:

\begin{proposition}
Assume that $G=(G_N)$ is uniform, let $n\in\{2,3,\ldots,\infty\}$ be minimal such that $G_n$ is not classical, and consider the following generation conditions:
\begin{enumerate}
\item Strong generation: $G_N=<G_N^c,G_n>$, for any $N>n$.

\item Usual generation: $G_N=<G_N^c,G_{N-1}>$, for any $N>n$.

\item Initial step generation: $G_{n+1}=<G_{n+1}^c,G_n>$.
\end{enumerate}
We have then $(1)\iff(2)\implies(3)$, and $(3)$ is in general strictly weaker.
\end{proposition}

\begin{proof}
All the implications and non-implications are elementary, as follows:

\medskip

$(1)\implies(2)$ This follows from $G_n\subset G_{N-1}$ for $N>n$, coming from uniformity.

\medskip

$(2)\implies(1)$ By using twice the usual generation, and then the uniformity, we have:
\begin{eqnarray*}
G_N
&=&<G_N^c,G_{N-1}>\\
&=&<G_N^c,G_{N-1}^c,G_{N-2}>\\
&=&<G_N^c,G_{N-2}>
\end{eqnarray*}

Thus we have a descent method, and we end up with the strong generation condition.

\medskip

$(2)\implies(3)$ This is clear, because (2) at $N=n+1$ is precisely (3).

\medskip

$(3)\hskip2.3mm\not\hskip-2.3mm\implies(2)$ In order to construct counterexamples here, simplest is to use group duals. Indeed, with $G_N=\widehat{\Gamma_N}$ and $\Gamma_N=<g_1,\ldots,g_N>$, the uniformity condition from Definition 12.6 tells us that we must be in a projective limit situation, as follows:
$$\Gamma_1\leftarrow\Gamma_2\leftarrow\Gamma_3\leftarrow\Gamma_4\leftarrow\ldots\quad,\quad \Gamma_{N-1}=\Gamma_N/<g_N=1>$$

But with this picture in hand, the result is clear. Indeed, assuming for instance that $\Gamma_2$ is given and not abelian, there are many ways of completing the sequence, and so the uniqueness coming from the generation condition in (2) can only fail.
\end{proof}

Let us introduce as well the following more technical notions:

\begin{proposition}
Assume that $G=(G_N)$ is uniform, let $n\in\{2,3,\ldots,\infty\}$ be as above, and consider the following conditions, where $I_N\subset G_N$ is the diagonal torus:
\begin{enumerate}
\item Strong hard liberation: $G_N=<G_N^c,I_n>$, for any $N\geq n$.

\item Technical condition: $G_N=<G_N^c,I_{N-1}>$ for any $N>n$, and $G_n=<G_n^c,I_n>$.

\item Hard liberation: $G_N=<G_N^c,I_N>$, for any $N$.

\item Initial step hard liberation: $G_n=<G_n^c,I_n>$.
\end{enumerate}
We have then $(1)\implies(2)\implies(3)\implies(4)$.
\end{proposition}

\begin{proof}
Our first claim is that when assuming that $G=(G_N)$ is uniform, the family of diagonal tori $I=(I_N)$ follows to be uniform as well. In order to prove this claim, observe first that the definition of the diagonal torus can be reformulated as follows:
$$I_N=G_N\cap\widehat{F_N}$$

WIth this picture in hand, the uniformity claim for $I=(I_N)$ comes from that of $G=(G_N)$, and from that of $\widehat{F}=(\widehat{F_N})$, which is trivial, as follows:
\begin{eqnarray*}
I_N\cap U_{N-1}^+
&=&(G_N\cap\widehat{F_N})\cap U_{N-1}^+\\
&=&(G_N\cap U_{N-1}^+)\cap(\widehat{F_N}\cap U_{N-1}^+)\\
&=&G_{N-1}\cap\widehat{F_{N-1}}\\
&=&I_{N-1}
\end{eqnarray*}

Thus our claim is proved, and this gives the various implications in the statement. 
\end{proof}

Let us discuss now to understand the relationship between the above conditions. In the group dual case, the simplest example to look at is the free real torus:
$$G=(T_N^+)$$

Here, with respect to the $3+4=7$ conditions that we have, the last 2 conditions trivially hold, and the first 5 conditions all require $T_3^+=<T_3,T_2^+>$, which is wrong. Indeed, in order to see this latter fact, consider the following discrete group:
$$\Gamma=\left<a,b,c\Big|a^2=b^2=c^2=1,[a,b]=[a,c]=1\right>$$

We have then $T_3\subset\widehat{\Gamma}$ and $T_2^+\subset\widehat{\Gamma}$ as well, and so we have:
$$<T_3,T_2^+>\subset\widehat{\Gamma}$$

On the other hand we have $\Gamma\neq\mathbb Z_2^{*3}$, and so $\widehat{\Gamma}\neq T_3^+$, and we conclude that we have:
$$<T_3,T_2^+>\neq T_3^+$$

With these preliminaries in hand, we can now formulate our main theoretical observation on the subject, which is something quite useful in practice, as follows:

\begin{theorem}
Assuming that $G=(G_N)$ is uniform, and with $n\in\{2,3,\ldots,\infty\}$ as above, minimal such that $G_n$ is not classical, the following conditions are equivalent,
\begin{enumerate}
\item Generation: $G_N=<G_N^c,G_{N-1}>$, for any $N>n$.

\item Strong generation: $G_N=<G_N^c,G_n>$, for any $N>n$.

\item Hard liberation: $G_N=<G_N^c,I_N>$, for any $N\geq n$.

\item Strong hard liberation: $G_N=<G_N^c,I_n>$, for any $N\geq n$.
\end{enumerate}
modulo their initial steps.
\end{theorem}

\begin{proof}
Our first claim is that generation plus initial step hard liberation imply the technical hard liberation condition. Indeed, the recurrence step goes as follows:
\begin{eqnarray*}
G_N
&=&<G_N^c,G_{N-1}>\\
&=&<G_N^c,G_{N-1}^c,I_{N-1}>\\
&=&<G_N^c,I_{N-1}>
\end{eqnarray*}

In order to pass now from the technical hard liberation condition to the strong hard liberation condition itself, observe that we have:
\begin{eqnarray*}
G_N
&=&<G_N^c,G_{N-1}>\\
&=&<G_N^c,G_{N-1}^c,I_{N-1}>\\
&=&<G_N^c,I_{N-1}>
\end{eqnarray*}

With this condition in hand, we have then as well:
\begin{eqnarray*}
G_N
&=&<G_N^c,G_{N-1}>\\
&=&<G_N^c,G_{N-1}^c,I_{N-2}>\\
&=&<G_N^c,I_{N-2}>
\end{eqnarray*}

This procedure can be of course be continued. Thus we have a descent method, and we end up with the strong hard liberation condition. In the other sense now, we want to prove that we have $G_N=<G_N^c,G_{N-1}>$ at $N\geq n$. At $N=n+1$ this is something that we already have. At $N=n+2$ now, we have:
\begin{eqnarray*}
G_{n+2}
&=&<G_{n+2}^c,I_n>\\
&=&<G_{n+2}^c,G_{n+1}^c,I_n>\\
&=&<G_{n+2}^c,G_{n+1}>
\end{eqnarray*}

This procedure can be of course be continued. Thus, we have a descent method, and we end up with the strong generation condition.
\end{proof}

The above results remain of course quite theoretical. Still at the theoretical level, we believe that the uniformity condition and generation condition are best viewed together. The idea indeed is that given a family of compact quantum groups $G=(G_N)$ with $G_N\subset U_N^+$, we have a ``ladder of cubes'', formed by cubes as follows:
$$
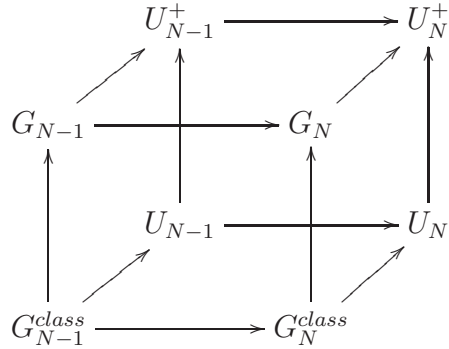
\xymatrix@R=20pt@C=15pt{
&U_{N-1}^+\ar[rr]&&U_N^+\\
G_{N-1}\ar[rr]\ar[ur]&&G_N\ar[ur]\\
&U_{N-1}\ar[rr]\ar[uu]&&U_N\ar[uu]\\
G_{N-1}^{class}\ar[uu]\ar[ur]\ar[rr]&&G_N^{class}\ar[uu]\ar[ur]
}$$

Thus, we have the question of investigating the $2\times 6=12$ intersection and generation properties, for the faces of such cubes, either with $N\in\mathbb N$ arbitrary, or with $N\geq n$.  These questions are quite interesting, and nothing much is known on all this, at least so far.

\section*{12c. Examples, duality}

Moving now forward, in order to avoid the above difficulties with the generation operation $<\,,>$, we can formulate a new definition, in the easy case, as follows:

\begin{definition}
We say that an easy quantum group $S_N\subset G\subset U_N^+$ appears as an easy soft liberation when we have the formula
$$G=\{G_{class},K\}$$
with $K=G\cap K_N^+$ being as usual its reflection subgroup, and with $\{\,,\}$ being the easy generation operation, obtained by $\cap$ at the level of categories of partitions.
\end{definition}

With this notion in hand, let us first go back to the quantum reflection groups, which were in need of liberation results. Let us recall from chapter 9 that we have:

\begin{theorem}
The easy quantum groups $H_N\subset G_N\subset H_N^+$, and the corresponding diagonal tori, are as follows,
$$
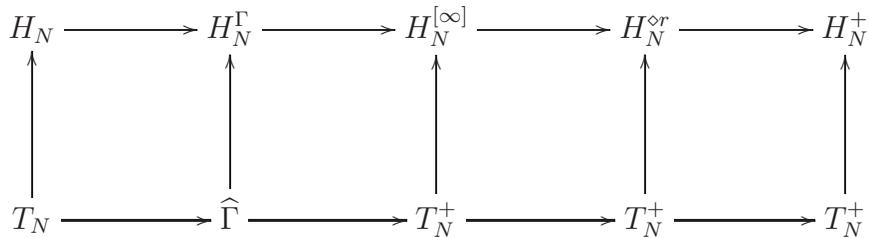
\xymatrix@R=50pt@C=50pt{
H_N\ar[r]&H_N^\Gamma\ar[r]&H_N^{[\infty]}\ar[r]&H_N^{\diamond r}\ar[r]&H_N^+\\
T_N\ar[u]\ar[r]&\widehat{\Gamma}\ar[r]\ar[u]&T_N^+\ar[r]\ar[u]&T_N^+\ar[r]\ar[u]&T_N^+\ar[u]
}$$
with the family $H_N^\Gamma$ and the series $H_N^{\diamond r}$ covering the endpoints.
\end{theorem}

\begin{proof}
The classification result is something that we know well, from chapters 9-10, and the assertion about the diagonal tori is clear as well from definitions. See \cite{rw3}.
\end{proof}

In relation now with our liberation questions, we can see that our hard liberation theory, based on blowing up the diagonal torus, cannot get beyond $H_N^{[\infty]}$. Thus, we have to focus on the quantum groups of type $H_N^\Gamma$. And here we have the following result:

\begin{theorem}
The quantum groups $H_N^\Gamma$ appear via hard liberation, as follows:
$$H_N^\Gamma=<H_N,\widehat{\Gamma}>$$
In particular, we have the ``master formula'' $H_N^{[\infty]}=<H_N,T_N^+>$.
\end{theorem}

\begin{proof}
We use the basic fact, from \cite{rw2}, and which is complementary to the easiness considerations above, that we have a crossed product decomposition as follows:
$$H_N^\Gamma=\widehat{\Gamma}\rtimes S_N$$

With this result in hand, we obtain that we have the missing inclusion, namely:
$$H_N^\Gamma\ =\ <S_N,\widehat{\Gamma}>\ \subset\ <H_N,\widehat{\Gamma}>$$

Finally, the last assertion is clear, by taking $\Gamma=\mathbb Z_2^{*N}$. Indeed, this group produces $H_N^{[\infty]}$, and the corresponding group dual is the free real torus $T_N^+$.
\end{proof}

As an interesting consequence of Theorem 12.12, let us record the following result:

\begin{proposition}
We have the following formula,
$$span(P_{even}^\Gamma)=span(P_{even})\cap C_{\widehat{\Gamma}}$$
where $C_{\widehat{\Gamma}}$ is the Tannakian category associated to $\widehat{\Gamma}$.
\end{proposition}

\begin{proof}
We use the Tannakian approach to the intersection and generation operations $\cap$ and $<\,,>$, which is summarized in the following well-known formulae: 
$$C_{G\cap H}=<C_G,C_H>\quad,\quad C_{<G,H>}=C_G\cap C_H$$

With these general formulae in hand, the generation formula in Theorem 12.12, namely $H_N^\Gamma=<H_N,\widehat{\Gamma}>$, reformulates in terms of Tannakian categories as follows:
$$C_{H_N^\Gamma}=C_{H_N}\cap C_{\widehat{\Gamma}}$$

But this is precisely the equality in the statement.
\end{proof}

In practice now, the category $C_{\widehat{\Gamma}}$ appearing in Proposition 12.13 is given by the following well-known formula, that we know well since chapter 1: 
$$C_{\widehat{\Gamma}}(k,l)=\left\{ T\in M_{N^l\times N^k}(\mathbb C)\Big|g_{i_1}\ldots g_{i_k}\neq g_{j_1}\ldots g_{j_l}\implies T_{j_1\ldots j_l,i_1\ldots i_k}=0\right\}$$

With this formula in hand, it is clear that the $\subset$ inclusion in Proposition 12.13 holds indeed, and that $\supset$ holds as well on $P_{even}$. However, having $\supset$ extended to the span of $P_{even}$ looks like a difficult combinatorial question. Thus, as a philosophical conclusion, the crossed product results in \cite{rw2} solve a difficult combinatorial question.

\bigskip

Let us discuss now the complex reflections. We first have here the following result:

\begin{theorem}
The easy quantum groups $K_N^\times$, which are as follows,
$$\xymatrix@R=50pt@C=50pt{
K_N\ar[r]&K_N^\Gamma\ar[r]&K_N^{[\infty]}\ar[r]&K_N^{\diamond r}\ar[r]&K_N^+\\
H_N\ar[u]\ar[r]&H_N^\Gamma\ar[r]\ar[u]&H_N^{[\infty]}\ar[r]\ar[u]&H_N^{\diamond r}\ar[r]\ar[u]&H_N^+\ar[u]
}$$
appear by easy soft liberation, $K_N^\times=<K_N,H_N^\times>$.
\end{theorem}

\begin{proof}
The point here is that the quantum groups $K_N^\times$, which are already known, appear indeed via easy soft liberation. But this latter fact follows from \cite{ba1}, \cite{tw1}.
\end{proof}

In order to discuss now hard liberation issues, let us start with:

\index{complexification}
\index{free complexification}
\index{soft liberation}
\index{hard liberation}

\begin{proposition}
The diagonal tori of the quantum groups $K_N^\times$ are as follows,
$$\xymatrix@R=50pt@C=50pt{
K_N\ar[r]&K_N^\Gamma\ar[r]&K_N^{[\infty]}\ar[r]&K_N^{\diamond r}\ar[r]&K_N^+\\
\mathbb T_N\ar[u]\ar[r]&\widehat{\Gamma_c}\ar[r]\ar[u]&\mathbb T_N^+\ar[r]\ar[u]&\mathbb T_N^+\ar[r]\ar[u]&\mathbb T_N^+\ar[u]
}$$
with $\Gamma\to\Gamma_c$ being a certain complexification operation, satisfying $<\mathbb T_N,\widehat{\Gamma}>\subset\widehat{\Gamma_c}$.
\end{proposition}

\begin{proof}
As a first observation, the results are clear and well-known for the endpoints $K_N,K_N^+$ and for the middle point $K_N^{[\infty]}$ as well. By functoriality it follows that the diagonal torus of $K_N^{[r]}$ must be the free complex torus $\mathbb T_N^+$, for any $r\in\mathbb N$, so we are done with the right part of the diagram. Regarding now the left part of the diagram, concerning the quantum groups $K_N^\Gamma$, if we denote by $T_1(.)$ the diagonal torus, we have:
\begin{eqnarray*}
T_1(K_N^\Gamma)
&=&T_1(<K_N,H_N^\Gamma>)\\
&\supset&<T_1(K_N),T_1(H_N^\Gamma)>\\
&=&<\mathbb T_N,\widehat{\Gamma}>
\end{eqnarray*}

Thus, we are led to the conclusion in the statement.
\end{proof}

Observe that the above inclusion $<\mathbb T_N,\widehat{\Gamma}>\subset\widehat{\Gamma_c}$ fails to be an isomorphism, and this for instance for $\Gamma=\mathbb Z_2^{*N}$. However, the construction $\Gamma\to\Gamma_c$ can be in principle explicitely computed, for instance by using Tannakian methods. Indeed, our soft liberation formula $K_N^\Gamma=<K_N,H_N^\Gamma>$ translates into a Tannakian formula, as follows:
$$\mathcal P_{even}^\Gamma=\mathcal P_{even}\cap P_{even}^\Gamma$$

The problem is that of explicitely computing the category on the left, corresponding to $K_N^\Gamma$, and then of deducing from this a presentation formula for the associated diagonal torus $\widehat{\Gamma_c}$. Now back to the hard liberation question, we have the following result:

\begin{theorem}
The quantum groups $K_N^\Gamma$ appear via hard liberation, and this even in a stronger form, as follows:
$$K_N^\Gamma=<K_N,\widehat{\Gamma}>$$
In particular, we have the formula $K_N^{[\infty]}=<K_N,T_N^+>$.
\end{theorem}

\begin{proof}
This follows from the above results. Indeed, we have:
\begin{eqnarray*}
K_N^\Gamma
&=&<K_N,H_N^\Gamma>\\
&=&<K_N,H_N,\widehat{\Gamma}>\\
&=&<K_N,\widehat{\Gamma}>
\end{eqnarray*}

Thus we have the formula in the statement, and the fact that this implies the fact that $K_N^\Gamma$ appears indeed via hard liberation follows from the above results as well. Finally, with $\Gamma=\mathbb Z_2^{*N}$ we obtain from this the formula $K_N^{[\infty]}=<K_N,T_N^+>$.
\end{proof}

In relation now with the orthogonal groups, the situation is much simpler, because the quantum groups $O_N\subset O_N^*\subset O_N^+$ are the only easy liberations of $O_N$. In addition, it is known that the inclusion $O_N\subset O_N^*$ is maximal, in the sense that it has no intermediate object at all. Also, as explained in \cite{bbc}, the conjecture is that  $O_N\subset O_N^*\subset O_N^+$ are the only liberations of $O_N$, not necessarily easy. In order to discuss this, we will need:

\begin{proposition}
We have the generation formula
$$O_N^+=<O_N,H_N^{[\infty]}>$$
where $H_N^{[\infty]}$ is the liberation of $H_N$ introduced before.
\end{proposition}

\begin{proof}
We use the Tannakian approach to $\cap$ and $<\,,>$. According to the general formula $C_{<G,H>}=C_G\cap C_H$, the formula in the statement is equivalent to:
$$C_{O_N^+}=C_{O_N}\cap C_{H_N^{[\infty]}}$$

By easiness, we are led into the following combinatorial statement:
$$NC_2=P_2\cap P_{even}^{[\infty]}$$

In order to establish this latter formula, we use one of the explicit descriptions of the category $P_{even}^{[\infty]}$ that we found before in chapter 9, which is as follows:
$$P_{even}^{[\infty]}=\left\{\pi\in P_{even}\Big|\sigma\in P_{even}^*,\forall\sigma\subset\pi\right\}$$

With this formula in hand, the fact that we have $NC_2\subset P_2\cap P_{even}^{[\infty]}$  is of course clear. This is in fact something that we already know, coming from:
$$O_N^+\supset<O_N,H_N^{[\infty]}>$$

Regarding the reverse inclusion, let $\pi\in P_2\cap P_{even}^{[\infty]}$. If we assume that $\pi$ has a crossing, then we have a basic crossing $\sigma\subset\pi$, and since we have $\sigma\not\in P_{even}^*$, we obtain in this way a contradition. Thus our reverse inclusion is proved, and we are done.
\end{proof}

As a comment here, the above result can be deduced as well from the standard easy classification results, by using the fact that the quantum group $O_N^\times=<O_N,H_N^{[\infty]}>$ is easy, and is not classical, nor half-classical. However, all this is ultimately too complicated, and having a direct and clear proof as above is probably something quite useful.

\bigskip

In relation now with our hard liberation questions, we have:

\begin{proposition}
The quantum groups $O_N,O_N^*,O_N^+$ all appear via hard liberation,
$$O_N^\times=<O_N,T_N^\times>$$
where $T_N^\times\subset O_N^\times$ is the diagonal torus, equal respectively to $T_N,T_N^*,T_N^+$.
\end{proposition}

\begin{proof}
This is trivial for $O_N$, and routine for $O_N^*$. In the case of $O_N^+$ the problem looks more difficult, but we have in fact the following proof:
\begin{eqnarray*}
O_N^+
&=&<O_N,H_N^{[\infty]}>\\
&=&<O_N,H_N,T_N^+>\\
&=&<O_N,T_N^+>
\end{eqnarray*}

Thus, we are led to the conclusion in the statement.
\end{proof}

Let us go back now to the conjecture regarding $O_N\subset O_N^*\subset O_N^+$, which is the most interesting statement around. It is known that $O_N\subset O_N^*\subset O_N^+$ are the unique easy liberations of $O_N$. In terms of our present formalism, this means that $O_N\subset O_N^*\subset O_N^+$ are the unique soft liberations of $O_N$. Here is a related result:

\begin{theorem}
The basic orthogonal quantum groups, namely
$$O_N\subset O_N^*\subset O_N^+$$
are the unique hard liberations of $O_N$.
\end{theorem}

\begin{proof}
A hard liberation of $O_N$ must appear by definition as follows, for a certain real reflection group $\mathbb Z_2^{*N}\to\Gamma\to\mathbb Z_2^N$, whose dual is the diagonal torus of the liberation:
$$O_N^\Gamma=<O_N,\widehat{\Gamma}>$$

On the other hand, we have the following computation, based on the fact that the class of easy quantum groups is stable under $<,>$:
\begin{eqnarray*}
O_N^\Gamma
&=&<O_N,\widehat{\Gamma}>\\
&=&<O_N,H_N,\widehat{\Gamma}>\\
&=&<O_N,H_N^\Gamma>\\
&\in&\{O_N,O_N^*,O_N^+\}
\end{eqnarray*}

Thus, we are led to the conclusion in the statement.
\end{proof}

The above is quite nice, and we believe that Theorem 12.19 can be further extended, by using the notion of spinned tori. In fact, all this leads us into the notion of Fourier liberation. We will discuss all this later in this chapter, at the end.

\bigskip

In relation now with the unitary quantum groups, we have:

\begin{theorem}
The basic unitary quantum groups, $U_N,U_N^*,U_N^+$, appear via real and complex soft liberation, and via hard liberation as well, as follows:
\begin{enumerate}
\item If we set $K_N^\times=U_N^\times\cap K_N^+$, we have $U_N^\times=<U_N,K_N^\times>$.

\item In fact, if we set  $H_N^\times=U_N^\times\cap H_N^+$, we have $U_N^\times=<U_N,H_N^\times>$.

\item In the free case, we have as well the formula $U_N^+=<U_N,H_N^{[\infty]}>$.

\item We have $U_N^\times=<U_N,I_N^\times>$, with $I_N^\times\subset U_N^\times$ being the diagonal torus.
\end{enumerate}
\end{theorem}

\begin{proof}
These results are trivial for $U_N$, and for $U_N^*,U_N^+$ the proofs are as follows:

\medskip

(1) This is well-known, coming from the following standard formulae:
$$\mathcal P_2^*=\mathcal P_2\cap\mathcal P_{even}^*\quad,\quad \mathcal{NC}_2=\mathcal P_2\cap\mathcal{NC}_{even}$$

(2) This enhances (1), by using the following standard formulae:
$$\mathcal P_2^*=\mathcal P_2\cap P_{even}^*\quad,\quad \mathcal{NC}_2=\mathcal P_2\cap NC_{even}$$

(3) This enhances (2) in the free case, and can be proved as follows:
\begin{eqnarray*}
U_N^+
&=&<U_N,O_N^+>\\
&=&<U_N,O_N,H_N^{[\infty]}>\\
&=&<U_N,H_N^{[\infty]}>
\end{eqnarray*}

(4) For $U_N^*$ we have indeed the following computation, based on (2):
\begin{eqnarray*}
U_N^*
&=&<U_N,H_N^*>\\
&=&<U_N,H_N,T_N^*>\\
&=&<U_N,T_N^*>\\
&\subset&<U_N,\mathbb T_N^*>
\end{eqnarray*}

For $U_N^+$ we can use a similar method, based on (3), as follows:
\begin{eqnarray*}
U_N^+
&=&<U_N,H_N^{[\infty]}>\\
&=&<U_N,H_N,T_N^+>\\
&=&<U_N,T_N^+>\\
&\subset&<U_N,\mathbb T_N^+>
\end{eqnarray*}

Since the reverse inclusions are clear, this finishes the proof.
\end{proof}

For the quantum groups $U_N^{(r)}$ the corresponding reflection groups $K_N^{(r)}=U_N^{(r)}\cap K_N^+$ can be explicitly computed, because we have a diagram as follows:
$$\xymatrix@R=50pt@C=50pt{
U_N^r\rtimes\mathbb Z_r\ar[r]&U_N^{(r)}\\
S_N^r\rtimes\mathbb Z_r\ar[r]\ar[u]&K_N^{(r)}\ar[u]
}$$

For the quantum groups $U_N^C$, however, the situation is considerably more complicated, because the corresponding reflection groups $K_N^C=U_N^C\cap K_N^+$ seem to collapse to $K_N,K_N^*,K_N^+$. Thus, we are in need of a new method here. 

\bigskip

The classification results for the liberations of $H_N,U_N$ have some obvious similarity between them. We have indeed a family followed by a series, and a series followed by a family, and this suggests the existence of a ``contravariant duality'', as follows:
$$\xymatrix@R=50pt@C=50pt{
U_N\ar[r]&U_N^{(r)}\ar[r]&U_N^C\ar[r]&U_N^+\\
H_N^+\ar@{.}[u]&H_N^{\diamond r}\ar[l]\ar@{.}[u]&H_N^\Gamma\ar[l]\ar@{.}[u]&H_N\ar[l]\ar@{.}[u]
}$$

In what follows we will attempt to axiomatize this duality. However, as we will soon discover, this is something quite complicated. Let us begin with:

\begin{definition}
We have ``covariant'' correspondences $H_N^\times\leftrightarrow U_N^\times$ between the liberations of $H_N$ and the liberations of $U_N$, constructed as follows:
\begin{enumerate}
\item To any $U_N^\times$ we can associate the quantum group $H_N^\times=U_N^\times\cap H_N^+$.

\item To any $H_N^\times$ we can associate the quantum group $U_N^\times=<H_N^\times,U_N>$.
\end{enumerate}
\end{definition}

Observe that both the above correspondences are indeed covariant. In practice now, in the easy case, we have the following result:

\begin{proposition}
The operations $U_N^\times\to U_N^\times\cap H_N^+$ and  $H_N^\times\to<H_N^\times,U_N>$ are both ``controlled'', in the easy case, by the corresponding quantum groups 
$$O_N^\times\in\{O_N,O_N^*,O_N^+\}$$
appearing via $U_N^\times\to U_N^\times\cap O_N^+$ and $H_N^\times\to<O_N,H_N^\times>$ respectively, and their images collapse to $\{H_N,H_N^*,H_N^+\}$ and $\{U_N,U_N^*,U_N^+\}$ respectively.
\end{proposition}

\begin{proof}
With $O_N^\times=U_N^\times\cap O_N^+$, we have the following computation:
\begin{eqnarray*}
H_N^\times
&=&U_N^\times\cap H_N^+\\
&=&U_N^\times\cap O_N^+\cap H_N^+\\
&=&O_N^\times\cap H_N^+\\
&\in&\{H_N,H_N^*,H_N^+\}
\end{eqnarray*}

Also, with $O_N^\times=<O_N,H_N^\times>$ this time, we have the following computation:
\begin{eqnarray*}
U_N^\times
&=&<U_N,H_N^\times>\\
&=&<U_N,O_N,H_N^\times>\\
&=&<U_N,O_N^\times>\\
&\in&\{U_N,U_N^*,U_N^+\}
\end{eqnarray*}

Thus, we are led to the conclusions in the statement. 
\end{proof}

Moving ahead, let us begin with an elementary statement, as follows:

\begin{proposition}
We have quantum groups $H_N\subset G_N\subset U_N^+$ as follows,
$$\xymatrix@R=27pt@C=40pt{
H_N\ar[r]&K_N\ar[r]&U_N\\
H_N^*\ar[u]\ar[r]&K_N^*\ar[r]\ar[u]&U_N^*\ar[u]\\
H_N\ar[r]\ar[u]&K_N\ar[u]\ar[r]&U_N\ar[u]}$$
and this is an intersection and generation diagram.
\end{proposition}

\begin{proof}
The fact that we have a diagram as above is clear from definitions, and the intersection and generation properties follow from easiness.
\end{proof}

In general now, any intermediate quantum group $H_N\subset G_N\subset U_N^+$ will appear inside the square, and we can therefore use some 2D orientation methods in order to deal with it. To be more precise, we can use the following observation:

\begin{proposition}
Given an intersection and generation diagram $P\subset Q,R\subset S$ and an intermediate quantum group $P\subset G\subset S$, we have a diagram as follows:
$$\xymatrix@R=25pt@C=30pt{
Q\ar[r]&<G,Q>\ar[r]&S\\
G\cap Q\ar[u]\ar[r]&G\ar[r]\ar[u]&<G,R>\ar[u]\\
P\ar[r]\ar[u]&G\cap R\ar[u]\ar[r]&R\ar[u]}$$
In addition, $G$ slices the square, in the sense that this is an intersection and generation diagram, precisely when $G=<G\cap Q,G\cap R>$ and $G=<G,Q>\cap<G,R>$.
\end{proposition}

\begin{proof}
This is indeed clear from definitions, because the intersection and generation diagram conditions are automatic for the upper left and lower right squares, as well as half of the generation diagram conditions for the lower left and upper right squares.
\end{proof}

Now back to our classification problem, we have the following result:

\begin{theorem}
The intermediate easy quantum groups $H_N\subset G_N\subset U_N^+$ which slice the square $H_N\subset H_N^+,U_N\subset U_N^+$, in the sense of Proposition 12.24, are as follows,
$$\xymatrix@R=27pt@C=40pt{
H_N^+\ar[r]&E_N^+\ar[r]&U_N^+\\
H_N^*\ar[u]\ar[r]&E_N^*\ar[r]\ar[u]&U_N^*\ar[u]\\
H_N\ar[r]\ar[u]&E_N\ar[u]\ar[r]&U_N\ar[u]}$$
with $H_N\subset E_N\subset U_N$ being an easy quantum group, and with $E_N^*,E_N^+$ being obtained via soft liberation, $E_N^*=<E_N,H_N^*>$ and $E_N^+=<E_N,H_N^+>$.
\end{theorem}

\begin{proof}
Assuming that $H_N\subset G_N\subset U_N^+$ is easy, and slices the square, its unitary version $G_N^u=<G_N,U_N>$ must be easy, and so is one of the easy quantum groups $U_N^\times$. Now observe that the slicing condition tells us in particular that $U_N^\times$ appears via the duality in Proposition 12.22 from its real discrete version $H_N^\times=U_N^\times\cap H_N^+$. Thus by duality we must have $U_N^\times\in\{U_N,U_N^*,U_N^+\}$, and this gives the result.
\end{proof}

As a remark here, when further imposing the uniformity condition the half-liberations dissapear, and we are left with the classical and free solutions, from \cite{tw2}.

\bigskip

Let us go back now to duality considerations, with the idea of ``fixing'' what we have. The classification results for $H_N,U_N$ have some obvious similarity between them. We have indeed a family followed by a series, and a series followed by a family, and this suggests the existence of a ``contravariant duality'', as follows:
$$\xymatrix@R=50pt@C=50pt{
U_N\ar[r]&U_N^{(r)}\ar[r]&U_N^C\ar[r]&U_N^+\\
H_N^+\ar@{.}[u]&H_N^{\diamond r}\ar[l]\ar@{.}[u]&H_N^\Gamma\ar[l]\ar@{.}[u]&H_N\ar[l]\ar@{.}[u]
}$$

As a first, naive attempt here, we could try to construct such a duality $H_N^\times\leftrightarrow U_N^\times$ by using a kind of ``complementation formula'', of the following type:
$$<H_N^\times,U_N^\times>=U_N^+$$

To be more precise, given a quantum group $H_N^\times$, we would like to define its dual $U_N^\times$ to be the ``minimal'' quantum group having the above property, and vice versa. Observe that such a correspondence $H_N^\times\leftrightarrow U_N^\times$ would be indeed contravariant. In practice now, however, the main problem comes from the following formula: 
$$U_N^+=<U_N,H_N^{[\infty]}>$$

Indeed, this formula shows that our naive attempt presented above simply fails, because the dual of $U_N$ would be $H_N^{[\infty]}$, instead of being $H_N^+$, as desired. However, our duality idea above still makes sense, and establishing it is a good open problem.

\section*{12d. Beyond easiness} 

Our aim here is to present some classification results, beyond easiness. Let us first discuss the half-classical case. We have the following definition, to start with:

\begin{definition}
The half-liberation of an intermediate compact group $H_N\subset G_N\subset U_N$ is the intermediate compact quantum group $H_N\subset G_N\subset U_N$ given by
$$G_N^*=<G_N,H_N^*>$$
with the generation operation being taken in a topological sense, as an operation for the closed subgroups of the free unitary quantum group $U_N^+$.
\end{definition}

This definition is something that we already met before, in a more general setting. As a main result regarding this operation, we have:

\begin{theorem}
The half-liberations of the uniform easy groups, namely
$$O_N,U_N\ ,\ H_N^2=H_N,H_N^4,H_N^6,\ldots,H_N^\infty=K_N$$
coincide with their usual half-liberations, taken in the easy sense, namely
$$O_N^*,U_N^*\ ,\ H_N^{2*}=H_N^*,H_N^{4*},H_N^{6*},\ldots,H_N^{\infty*}=K_N^*$$
obtained by liberating, and then by imposing the relations $abc=cba$.
\end{theorem}

\begin{proof}
This is something standard. First of all, it follows from \cite{tw2} that the easy compact groups $S_N\subset G_N\subset U_N$ satisfying the uniformity assumption $G_{N-1}=G_N\cap U_{N-1}^+$ are precisely those in the statement, with the usual convention for reflections, namely:
$$H_N^s=\mathbb Z_s\wr S_N$$

In order to compute the half-liberations in our sense, we use the fact that the operations $<\,,>$ and $\cap$ are ``dual'' to each other via Tannakian duality $G\leftrightarrow C$, as follows:
$$C_{<G,H>}=C_G\cap C_H\quad,\quad C_{G\cap H}=<C_G,C_H>$$

With standard easy quantum group notations, if we denote by $D$ the category of partitions for $G_N$, and by $G_N^\times$ the easy half-liberation of $G_N$, we have then:
\begin{eqnarray*}
C_{G_N^*}
&=&C_{G_N}\cap C_{H_N^*}\\
&=&span(D)\cap span(P_{even}^*)\\
&=&span(D\cap P_{even}^*)\\
&=&span(D\cap NC_{even},\slash\hskip-1.6mm|\hskip-1.6mm\backslash)\\
&=&C_{G_N^\times}
\end{eqnarray*}

Here all the equalities are well-known and standard. Thus $G_N^*=G_N^\times$, as claimed.
\end{proof}

Summarizing, we have so far a notion of half-liberation for the intermediate compact groups $H_N\subset G_N\subset U_N$, which works well in the easy case. Next, we have:

\begin{definition}
The half-liberation of an intermediate compact group 
$$T_N\subset G_N\subset U_N$$
is the intermediate compact quantum group 
$$T_N^*\subset G_N^*\subset U_N$$
given by the following formula,
$$G_N^*=<G_N,T_N^*>$$
with the generation operation being taken as usual in a topological sense.
\end{definition}

Our first task is to verify that this more general notion is compatible with the one that we already have. This is something non-trivial, and we have indeed:

\begin{theorem}
For an intermediate compact group 
$$H_N\subset G_N\subset U_N$$
its ``soft'' and ``hard'' half-liberations, from Definitions 12.26 and 12.28, coincide.
\end{theorem}

\begin{proof}
We must prove that for any intermediate compact group $H_N\subset G_N\subset U_N$, we have the following equality, between closed subgroups of $U_N^+$:
$$<G_N,H_N^*>=<G_N,T_N^*>$$

It is enough to solve the problem for the smallest possible group under consideration, namely $G_N=H_N$. Thus, we are led into the following question:
$$H_N^*=<H_N,T_N^*>$$

Now let us denote by $G\leftrightarrow C$ the standard Tannakian correspondence, from chapter 2, with as usual $C=(Hom(u^{\otimes k},u^{\otimes l}))$, let also $G\to PG$ be the projective version construction, and let us denote the construction in \cite{bdu} as follows: 
$$(G\subset U_N)\to([G]\subset O_N^*)$$

We have then, by using a number of standard facts:
\begin{eqnarray*}
PH_N^*=P<H_N,T_N^*>
&\iff&PK_N=P<H_N,T_N^*>\\
&\iff&C_{PK_N}=C_{P<H_N,T_N^*>}\\
&\iff&C_{PK_N}=C_{PH_N}\cap C_{PT_N^*}\\
&\iff&C_{PK_N}=C_{PH_N}\cap C_{P\mathbb T_N}\\
&\iff&C_{PK_N}=C_{P<H_N,\mathbb T_N>}\\
&\iff&C_{PK_N}=C_{PK_N}
\end{eqnarray*}

Thus the projective versions coincide, and so the affine lifts must coincide as well.
\end{proof}

As a comment here, the above proof is not the only one. Since $H_N,H_N^*$ are both easy, coming from $P_{even},P_{even}^*$, our question $H_N^*=<H_N,T_N^*>$ becomes:
$$span(P_{even}^*)=span(P_{even})\cap C_{T_N^*}$$

But this can be proved by standard combinatorics, based on the standard fact that the half-classical combinatorics comes from the infinite symmetric group $S_\infty$.

\bigskip

Summarizing, we have now a notion of half-liberation for the intermediate compact groups $T_N\subset G_N\subset U_N$, which works well in the easy case. Next, we have:

\begin{proposition}
For any compact group $T_N\subset G_N\subset U_N$ we have the formula
$$G_N^*=[G_N^\circ]$$
where $G_N^\circ=<G_N,\mathbb T_N>$, and where $E_N\to[E_N]$ is the construction in \cite{bdu}.
\end{proposition}

\begin{proof}
This can be proved by using the same method as for Theorem 12.29. With the notations from there, we have the following computation:
\begin{eqnarray*}
PG_N^*=P[G_N^\circ]
&\iff&PG_N^*=PG_N^\circ\\
&\iff&P<G_N^*,T_N^*>=P<G_N,\mathbb T_N>\\
&\iff&C_{P<G_N^*,T_N^*>}=C_{P<G_N,\mathbb T_N>}\\
&\iff&C_{PG_N}\cap C_{PT_N^*}=C_{PG_N}\cap C_{P\mathbb T_N^*}\\
&\iff&C_{PG_N}\cap C_{PT_N^*}=C_{PG_N}\cap C_{PT_N^*}
\end{eqnarray*}

Thus the projective versions coincide, and so the affine lifts must coincide as well.
\end{proof}

In the unitary case, the situation is similar, and we have:

\begin{theorem}
For any compact group $T_N\subset G_N\subset U_N$ we have the formula
$$G_N^*=[[G_N^\circ]]$$
where $G_N^\circ=<G_N,\mathbb T_N>$, and where $E_N\to[[E_N]]$ is the construction in \cite{bb2}.
\end{theorem}

\begin{proof}
The computation here is identical with the one in the proof of Proposition 12.30, with technical ingredients coming this time from \cite{bb2}.
\end{proof}

As a further theme of discussion, let us discuss now the non-easy extension of the notion of orientability. Things are quite tricky here, and we must start as follows:

\index{classical version}
\index{discrete version}
\index{real version}
\index{free version}
\index{smooth version}
\index{unitary version}

\begin{definition}
Associated to any closed subgroup $G_N\subset U_N^+$ are its classical, discrete and real versions, given by
$$G_N^c=G_N\cap U_N$$
$$G_N^d=G_N\cap K_N^+$$
$$G_N^r=G_N\cap O_N^+$$
as well as its free, smooth and unitary versions, given by
$$G_N^f=<G_N,H_N^+>$$
$$G_N^s=<G_N,O_N>$$
$$G_N^u=<G_N,K_N>$$
where $<\,,>$ is the usual, non-easy topological generation operation.
\end{definition}

Observe the difference, and notational clash, with some of the notions used before. To be more precise, as explained in Part I, it is believed that we should have a formula of type $\{\,,\}=<\,,>$, but this is not clear at all, and the problem comes from this.

\bigskip

A second issue comes when composing the above operations, and more specifically those involving the generation operation, once again due to the conjectural status of the formula $\{\,,\}=<\,,>$. Due to this fact, instead of formulating a result here, we have to formulate a second definition, complementary to Definition 12.32, as follows:

\begin{definition}
Associated to any closed subgroup $G_N\subset U_N^+$ are the mixes of its classical, discrete and real versions, given by
$$G_N^{cd}=G_N\cap K_N$$
$$G_N^{cr}=G_N\cap O_N^+$$
$$G_N^{dr}=G_N\cap H_N^+$$
as well as the mixes of its free, smooth and unitary versions, given by
$$G_N^{fs}=<G_N,O_N^+>$$
$$G_N^{fu}=<G_N,K_N^+>$$
$$G_N^{us}=<G_N,U_N>$$
where $<\,,>$ is the usual, non-easy topological generation operation.
\end{definition}

Now back to our orientation questions, the slicing and bi-orientability conditions lead us again into $\{\,,\}$ vs. $<\,,>$ troubles, and are therefore rather to be ignored. The orientability conditions, however, have the following analogue:

\index{orientability}

\begin{definition}
A closed subgroup $G_N\subset U_N^+$ is called ``oriented'' if
$$G_N=<G_N^{cd},G_N^{cr},G_N^{dr}>$$
$$G_N=G_N^{fs}\cap G_N^{fu}\cap G_N^{su}$$
and ``weakly oriented'' if the following conditions hold,
$$G_N=<G_N^c,G_N^d,G_N^r>$$
$$G_N=G_N^f\cap G_N^s\cap G_N^u$$
where the various versions are those in Definition 12.32 and Definition 12.33.
\end{definition}

With these notions, our claim is that some classification results are possible:

\bigskip

(1) In the classical case, we believe that the uniform, half-homogeneous, oriented groups are those that we know, with some bistochastic versions excluded. This is of course something quite heavy, well beyond easiness, with the potential tools available for proving such things coming from advanced finite group theory and Lie algebra theory.  Our uniformity axiom could play a key role here, when combined with \cite{sto}, in order to exclude all the exceptional objects which might appear on the way.

\bigskip

(2) In the free case, under similar assumptions, we believe that the solutions should be those that we know, once again with some bistochastic versions excluded. This is something heavy, too, related to a well-known conjecture, namely $<G_N,S_N^+>=\{\bar{G}_N,S_N^+\}$. Indeed, assuming that we would have such a formula, and perhaps some more formulae of the same type as well, we could in principle work out our way inside the cube, from the edge and face projections to $G_N$ itself, and in this process $G_N$ would become easy. This would be the straightforward strategy here.

\bigskip

(3) In the group dual case, the orientability axiom simplifies, because the group duals are discrete in our sense. We believe that the uniform, twistable, oriented group duals should appear as combinations of certain abelian groups, which appear in the classical case, with duals of varieties of real reflection groups, which appear in the real case. This is probably the easiest question in the present series, and the most reasonable one, to start with. However, there are no concrete results so far, in this direction.

\bigskip

Finally, let us discuss the notion of Fourier liberation, which conjecturally solves some of the problems raised in the above. We first have the following standard notion:

\begin{proposition}
Given a closed subgroup $G\subset U_N^+$ and a matrix $Q\in U_N$, we let $T_Q\subset G$ be the diagonal torus of $G$, with fundamental representation spinned by $Q$:
$$C(T_Q)=C(G)\Big/\left<(QuQ^*)_{ij}=0\Big|\forall i\neq j\right>$$
This torus is then a group dual, $T_Q=\widehat{\Lambda}_Q$, where $\Lambda_Q=<g_1,\ldots,g_N>$ is the discrete group generated by the elements $g_i=(QuQ^*)_{ii}$, which are unitaries inside $C(T_Q)$.
\end{proposition}

\begin{proof}
This follows indeed from definitions, because, as said in the statement, $T_Q$ is by definition a diagonal torus. Equivalently, since $v=QuQ^*$ is a unitary corepresentation, its diagonal entries $g_i=v_{ii}$, when regarded inside $C(T_Q)$, are unitaries, and satisfy:
$$\Delta(g_i)=g_i\otimes g_i$$

Thus $C(T_Q)$ is a group algebra, and more specifically we have $C(T_Q)=C^*(\Lambda_Q)$, where $\Lambda_Q=<g_1,\ldots,g_N>$ is the group in the statement, and this gives the result.
\end{proof}

Summarizing, associated to any closed subgroup $G\subset U_N^+$ is a whole family of tori, indexed by the unitaries $U\in U_N$. As a first result regarding these tori, we have:

\begin{theorem}
For the quantum permutation group $S_N^+$, the discrete group quotient $F_N\to\Lambda_Q$ with $Q\in U_N$ comes from the following relations:
$$\begin{cases}
g_i=1&{\rm if}\ \sum_lQ_{il}\neq 0\\
g_ig_j=1&{\rm if}\ \sum_lQ_{il}Q_{jl}\neq 0\\ 
g_ig_jg_k=1&{\rm if}\ \sum_lQ_{il}Q_{jl}Q_{kl}\neq 0
\end{cases}$$
Also, given a decomposition $N=N_1+\ldots+N_k$, for the matrix $Q=diag(F_{N_1},\ldots,F_{N_k})$, where $F_N=\frac{1}{\sqrt{N}}(\xi^{ij})_{ij}$ with $\xi=e^{2\pi i/N}$ is the Fourier matrix, we obtain
$$\Lambda_Q=\mathbb Z_{N_1}*\ldots*\mathbb Z_{N_k}$$
with dual embedded into $S_N^+$ in a standard way.
\end{theorem}

\begin{proof}
This can be proved by a direct computation, as follows:

\medskip

(1) Fix a unitary matrix $Q\in U_N$, and consider the following quantities:
$$\begin{cases}
c_i=\sum_lQ_{il}\\
c_{ij}=\sum_lQ_{il}Q_{jl}\\
d_{ijk}=\sum_l\bar{Q}_{il}\bar{Q}_{jl}Q_{kl}
\end{cases}$$

We write $w=QvQ^*$, where $v$ is the fundamental corepresentation of $C(S_N^+)$. Assume $X\simeq\{1,\ldots,N\}$, and let $\alpha$ be the coaction of $C(S_N^+)$ on $C(X)$. Let us set:
$$\varphi_i=\sum_l\bar{Q}_{il}\delta_l\in C(X)$$

Also, let $g_i=(QvQ^*)_{ii}\in C^*(\Lambda_Q)$. If $\beta$ is the restriction of $\alpha$ to $C^*(\Lambda_Q)$, then:
$$\beta(\varphi_i)=\varphi_i\otimes g_i$$

(2) Now recall that $C(X)$ is the universal $C^*$-algebra generated by elements $\delta_1,\ldots,\delta_N$ which are pairwise orthogonal projections. Writing these conditions in terms of the linearly independent elements $\varphi_i$ by means of the formulae $\delta_i=\sum_lQ_{il}\varphi_l$, we find that the universal relations for $C(X)$ in terms of the elements $\varphi_i$ are as follows:
$$\begin{cases}
\sum_ic_i\varphi_i=1\\
\varphi_i^*=\sum_jc_{ij}\varphi_j\\
\varphi_i\varphi_j=\sum_kd_{ijk}\varphi_k
\end{cases}$$

(3) Let $\widetilde{\Lambda}_Q$ be the group in the statement. Since $\beta$ preserves these relations, we get:
$$\begin{cases}
c_i(g_i-1)=0\\
c_{ij}(g_ig_j-1)=0\\
d_{ijk}(g_ig_j-g_k)=0
\end{cases}$$

We conclude from this that $\Lambda_Q$ is a quotient of $\widetilde{\Lambda}_Q$. On the other hand, it is immediate that we have a coaction map as follows:
$$C(X)\to C(X)\otimes C^*(\widetilde{\Lambda}_Q)$$

Thus $C(\widetilde{\Lambda}_Q)$ is a quotient of $C(S_N^+)$. Since $w$ is the fundamental corepresentation of $S_N^+$ with respect to the basis $\{\varphi_i\}$, it follows that the generator $w_{ii}$ is sent to $\widetilde{g}_i\in\widetilde{\Lambda}_Q$, while $w_{ij}$ is sent to zero. We conclude that $\widetilde{\Lambda}_Q$ is a quotient of $\Lambda_Q$. Since the above quotient maps send generators on generators, we conclude that $\Lambda_Q=\widetilde{\Lambda}_Q$, as desired.

\medskip

(4) We apply the result found in (3), with the $N$-element set $X$ there being:
$$X=\mathbb Z_{N_1}\sqcup\ldots\sqcup\mathbb Z_{N_k}$$

With this choice, we have $c_i=\delta_{i0}$ for any $i$. Also, we have $c_{ij}=0$, unless $i,j,k$ belong to the same block to $Q$, in which case $c_{ij}=\delta_{i+j,0}$, and also $d_{ijk} =0$, unless $i,j,k$ belong to the same block of $Q$, in which case $d_{ijk}=\delta_{i+j,k}$. We conclude from this that $\Lambda_Q$ is the free product of $k$ groups which have generating relations as follows:
$$g_ig_j=g_{i+j}\quad,\quad g_i^{-1}=g_{-i}$$

But this shows that our group is $\Lambda_Q=\mathbb Z_{N_1}*\ldots*\mathbb Z_{N_k}$, as stated.
\end{proof}

In connection with our liberation questions for the subgroups $G\subset S_N$, all this is quite interesting, and suggests formulating the following definition:

\begin{definition}
Consider a closed subgroup $G\subset U_N^+$.
\begin{enumerate}
\item Its standard tori $T_F$, with $F=F_{N_1}\otimes\ldots\otimes F_{N_k}$, and $N=N_1+\ldots+N_k$ being regarded as a partition, are called Fourier tori.

\item In the case where we have $G_N=<G_N^c,(T_F)_F>$, we say that $G_N$ appears as a Fourier liberation of its classical version $G_N^c$.
\end{enumerate}
\end{definition}

The conjecture is then that the easy quantum groups should appear as Fourier liberations. With respect to the basic examples, the situation in the free case is as follows:

\bigskip

(1) $O_N^+,U_N^+$ are diagonal liberations, so they are Fourier liberations as well.

\bigskip

(2) $B_N^+,C_N^+$ are Fourier liberations too, with this being standard.

\bigskip

(3) $S_N^+$ is a Fourier liberation too, being generated by its tori \cite{bcf}, \cite{chi}.

\bigskip

(4) $H_N^+,K_N^+$ remain to be investigated, by using the general theory in \cite{rw3}.

\bigskip

Finally, let us mention that the notion of Fourier liberation is something specific to the easy case. Indeed, for the general compact quantum groups, this will not work.

\section*{12e. Exercises} 

Things have been difficult in this chapter, and there have been open questions all around the place. As exercise, that we believe to be central to all this, we have:

\begin{exercise}
Prove that the easy quantum groups appear via Fourier liberation.
\end{exercise}

This is a quite difficult exercise, and ideally we would need here a global proof, based on the abstract notion of easiness. But some case-by-case verifications would be extremely useful as well, the point being that each such verification requires lots of work.

\part{Super-easiness}

\ \vskip50mm

\begin{center}
{\em Maybe one day we'll be united

And our love won't be divided

Maybe one day we'll be united

And our love won't be divided}
\end{center}

\chapter{Schur-Weyl twists}

\section*{13a. Ad-hoc twisting}

In this final part of the present book we discuss how our previous results, from Parts I, II, III, can be used in order to say more things about the closed subgroups $G\subset U_N^+$. Normally we would like all these subgroups to be easy, but as we know well, they aren't. However, as we will soon discover, with explicit examples, some of them are in fact not very far from being easy. In view of this, our goal with be very simple, as follows:

\begin{goal}
Extend the easiness theory, by not deviating much from its original spirit, into a super-easiness theory, covering as many examples as possible. 
\end{goal}

This looks quite reasonable, and the key words in all this are of course ``by not deviating much''. That is, easiness theory is something useful and practical, and we would like of course its super-easiness extension to be useful and practical too, avoiding too much generality, which can only lead into weak, useless abstractions. Let us mention too that the name ``super-easiness'' is not exactly a joke as it seems, because the first such extension, constructed in \cite{bsk}, crucially used a ``super-space'' idea, leading to this name. But more on this in chapter 14 below, where we will discuss the constructions in \cite{bsk}.

\bigskip

Getting started now, ``by not deviating much'' means to keep things Tannakian. So, the idea will be very simple, namely that of modifying the operation $\pi\to T_\pi$, and constructing new subgroups $G\subset U_N^+$ in this way. This operation was given by:
$$T_\pi(e_{i_1}\otimes\ldots\otimes e_{i_k})=\sum_{j_1,\ldots,j_l}\delta_\pi\begin{pmatrix}i_1&\ldots&i_k\\ j_1&\ldots&j_l\end{pmatrix}e_{j_1}\otimes\ldots\otimes e_{j_l}$$

So, assume that by whatever magical trick, we manage to ``twist'' this operation $\pi\to T_\pi$, into a new operation $\pi\to\dot{T}_\pi$, still satisfying the categorial conditions, or at least a suitable modification of these categorical conditions, as for any category of partitions $D\subset P$ to produce a Tannakian category $C$, via the usual formula, namely:
$$C_{kl}=span\left(\dot{T}_\pi\Big|\pi\in D(k,l)\right)$$

Then, this will be a win, because we can define afterwards for any $N\in\mathbb N$ a certain closed subgroup $\dot{G}_N\subset U_N^+$, via Tannakian duality, having by definition $C$ as Tannakian category. We can then call this new subgroup $\dot{G}_N\subset U_N^+$ super-easy, add it to the usual easy group $G_N\subset U_N^+$ associated to $D$, and we will have in this way a nice extension of the easy quantum group formalism, covering far more objects than before.

\bigskip

All this looks very nice, but in practice, in the classical case already, our problems are not exactly obvious, and we are led in this way to the following question:

\begin{question}
How far can we go with super-easiness, in the classical case? Can we fully cover the Lie types ABCD? What about EFG? What about the complex reflection groups $H_N^{sd}$? What about the exceptional complex reflection groups?
\end{question}

However, all this looks a bit scary, more or less requiring us to become top-level experts in Lie groups or reflection groups, so we will keep such things for later.

\bigskip

As another idea, dual to the above one, we can look into discrete group duals. We know indeed, and this since chapter 1, that the Tannakian categories of group duals appear via the combinatorics of the group, and so are not very far from easiness. Thus, we have a good idea here, and we can formulate our second question as follows:

\begin{question}
Shall we orient our super-easiness theory towards group duals, for instance with, as a first objective, including the dual of the free group $F_N$?
\end{question}

As a first observation, although not exactly following the scheme $T_\pi\to\dot{T}_\pi$ evoked above, this looks quite reasonable, when compared for instance with Question 13.2. However, there is a downside to this, coming precisely from the ``reasonability'' of our question. Frankly, with such things we would not reach to any new mathematics, the dual of  $F_N$ being a sort of ``mathematical wheel'', that we would not like to reinvent, and with this in mind, Question 13.2 looks much more interesting, and so, more reasonable.

\bigskip

As a third attempt now to formulate a starting question, which is something natural as well, and certainly in tune with the spirit of abstract mathematics, we have:

\begin{question}
What about easiness naturally producing super-easiness, via various product operations, such as the usual product $\times$, or more complicated products?
\end{question}

To be more precise here, wanting for instance $G\times H$ to be super-easy when $G,H$ are easy certainly looks like a reasonable mathematical idea. However, a bit like it was the case with Question 13.3, this looks a bit boring, because the same pure mathematics leading to such ideas tell us to look for ``irreducible'' objects first.

\bigskip

As a fourth attempt now, coming by correcting a bit Question 13.4, we have:

\begin{question}
What about easiness naturally producing super-easiness, via various twisting operations?
\end{question}

This certainly looks more reasonable, because unlike the product operations from Question 13.4, bringing weakness, the twisting operations usually preserve the ``irreducibility'' of the quantum group. That is, if we start with a quality, indecomposable quantum group, then the twist will be as well a quality, indecomposable quantum group. But, as a downside, this latter question promises to be as technical as Question 13.2.

\bigskip

So, this is the situation, nothing both interesting and easily doable, on the spot, and shall we look for a fifth idea, or ask the cat. Hope you're with me for asking the cat, all this being a bit tiring. And cat, Paul Adrien Maurice by his name, answers:

\begin{cat}
Shut up and compute. You have $SU_2,SO_3$. Anticommutation twists too. Don't call yourself physicist until you tried these.
\end{cat}

Thanks cat, this sounds reasonable, let's not forget indeed about $SU_2,SO_3$, which are the alpha and omega of physics. And not forget either about commutation replaced by anticommutation, that's one good learning from basic quantum mechanics too.

\bigskip

We will do all this in the present Part IV of this book, the plan being:

\begin{plan}
In order to develop super-easiness theory, we must:
\begin{enumerate}
\item Cover the anticommutation twists $G'\subset U_N^+$ of the easy Lie groups $G\subset U_N$.

\item Cover $SU_2$, and more generally the symplectic groups $Sp_N\subset U_N$.

\item Cover $SO_3$, and more generally $S_Z^+$, with $Z$ finite quantum space.

\item And finally, merge these $3$ extensions into a super-easiness theory.
\end{enumerate}
\end{plan}

So, this will be the plan for the present Part IV of this book. In this chapter we will do (1), following \cite{bbc} and related papers, which will turn to be something quite simple. Then in chapter 14 and chapter 15 we will do (2) and (3), following respectively \cite{bsk}, and \cite{wa2} and related papers. And finally at the end of chapter 15 we will attempt to solve (4), which will lead us into a lot of thinking, and into some open problems too.

\bigskip

Getting started now, we need to twist the easy Lie groups $G\subset U_N$, and why not the more general easy quantum groups $G\subset U_N^+$. However, this does not look like a trivial question, and this even for $G=U_N$ itself, because we must keep some of the commutation relations $ab=ba$ between basic coordinates, while changing some other into anticommutation, $ab=-ba$. And things become even more complicated when looking for instance at $U_N^*$, with each relation of type $abc=cba$ waiting to be studied, and then either kept as such, or replaced by its anticommutation counterpart, $abc=-cba$.

\bigskip

Fortunately, there is a clever answer to this, providing us from doing too many computations. The idea will be that of twisting first the simplest objects that we have, namely the associated spheres. We haven't talked so far in this book about the spheres associated to our main unitary groups, but never too late, and here is their definition:

\index{quantum sphere}
\index{free sphere}

\begin{definition}
We have quantum spheres, with inclusions between them,
$$\xymatrix@R=13mm@C=10mm{
S^{N-1}_{\mathbb R,+}\ar[r]&\mathbb TS^{N-1}_{\mathbb R,+}\ar[r]&S^{N-1}_{\mathbb C,+}\\
S^{N-1}_{\mathbb R,*}\ar[r]\ar[u]&\mathbb TS^{N-1}_{\mathbb R,*}\ar[r]\ar[u]&S^{N-1}_{\mathbb C,*}\ar[u]\\
S^{N-1}_\mathbb R\ar[r]\ar[u]&\mathbb TS^{N-1}_\mathbb R\ar[r]\ar[u]&S^{N-1}_\mathbb C\ar[u]}$$
with the free complex sphere, on top right, being given by
$$C(S^{N-1}_{\mathbb C,+})
=C^*\left(x_1,\ldots,x_N\Big|\sum_ix_ix_i^*=\sum_ix_i^*x_i=1\right)$$
and with the other spheres being obtained as subspaces, in the obvious way.
\end{definition}

There are many things that can be said about these spheres, the idea being that the usual correspondences $S^{N-1}_\mathbb R\leftrightarrow O_N$ and $S^{N-1}_\mathbb C\leftrightarrow U_N$ from classical geometry, that you surely know well, extend to all the above spheres. In particular, with a suitable notion of quantum isometry group, the quantum isometry groups of our spheres are as follows:
$$\xymatrix@R=13mm@C=13mm{
O_N^+\ar[r]&\mathbb TO_N^+\ar[r]&U_N^+\\
O_N^*\ar[r]\ar[u]&\mathbb TO_N^*\ar[r]\ar[u]&U_N^*\ar[u]\\
O_N\ar[r]\ar[u]&\mathbb TO_N\ar[r]\ar[u]&U_N\ar[u]}$$

In what follows we will not really need this, and for more on such spheres, we refer to \cite{bgo} and related papers. However, based on this philosophy, we can try first to twist the spheres, which looks like an easy task, because we will be dealing here with single indices instead of double indices, and then, using this, twist afterwards the quantum groups as well. So, this will be our plan, and getting started now, we first have:

\index{twisting}
\index{twisted commutation}
\index{twisted half-commutation}
\index{twisted sphere}

\begin{theorem}
We have quantum spheres as follows, obtained via the twisted commutation relations $ab=\pm ba$, and twisted half-commutation relations $abc=\pm cba$,
$$\xymatrix@R=13mm@C=10mm{
S^{N-1}_{\mathbb R,+}\ar[r]&\mathbb TS^{N-1}_{\mathbb R,+}\ar[r]&S^{N-1}_{\mathbb C,+}\\
S^{N-1,\prime}_{\mathbb R,*}\ar[r]\ar[u]&\mathbb TS^{N-1,\prime}_{\mathbb R,*}\ar[r]\ar[u]&S^{N-1,\prime}_{\mathbb C,*}\ar[u]\\
S^{N-1,\prime}_\mathbb R\ar[r]\ar[u]&\mathbb TS^{N-1,\prime}_\mathbb R\ar[r]\ar[u]&S^{N-1,\prime}_\mathbb C\ar[u]}$$
with the precise signs being as follows:
\begin{enumerate}
\item The signs on the bottom correspond to anticommutation of distinct coordinates, and their adjoints. That is, with $z_i=x_i,x_i^*$ and $\varepsilon_{ij}=1-\delta_{ij}$, the formula is:
$$z_iz_j=(-1)^{\varepsilon_{ij}}z_jz_i$$

\item The signs in the middle come from functoriality, as for the spheres in the middle to contain those on the bottom. That is, the formula is:
$$z_iz_jz_k=(-1)^{\varepsilon_{ij}+\varepsilon_{jk}+\varepsilon_{ik}}z_kz_jz_i$$
\end{enumerate}
\end{theorem}

\begin{proof}
This is something elementary, from \cite{ba1}, the idea being as follows:

\medskip

(1) Here there is nothing to prove, because we can define the spheres on the bottom by the following formulae, with $z_i=x_i,x_i^*$ and $\varepsilon_{ij}=1-\delta_{ij}$ being as above:
$$C(S^{N-1,\prime}_\mathbb R)=C(S^{N-1}_{\mathbb R,+})\Big/\Big<x_ix_j=(-1)^{\varepsilon_{ij}}x_jx_i\Big>$$
$$C(S^{N-1,\prime}_\mathbb C)=C(S^{N-1}_{\mathbb C,+})\Big/\Big<z_iz_j=(-1)^{\varepsilon_{ij}}z_jz_i\Big>$$

(2) Here our claim is that, if we want to construct half-classical twisted spheres, via relations of type $abc=\pm cba$ between the coordinates $x_i$ and their adjoints $x_i^*$, as for these spheres to contain the twisted spheres constructed in (1), the only possible choice for these relations is as follows, with $z_i=x_i,x_i^*$ and $\varepsilon_{ij}=1-\delta_{ij}$ being as above:
$$z_iz_jz_k=(-1)^{\varepsilon_{ij}+\varepsilon_{jk}+\varepsilon_{ik}}z_kz_jz_i$$

But this is something clear, coming from the following computation, inside of the quotient algebras corresponding to the twisted spheres constructed in (1):
\begin{eqnarray*}
z_iz_jz_k
&=&(-1)^{\varepsilon_{ij}}z_jz_iz_k\\
&=&(-1)^{\varepsilon_{ij}+\varepsilon_{ik}}z_jz_kz_i\\
&=&(-1)^{\varepsilon_{ij}+\varepsilon_{jk}+\varepsilon_{ik}}z_kz_jz_i
\end{eqnarray*}

Thus, we are led to the conclusion in the statement, the spheres being given by:
$$C(S^{N-1,\prime}_{\mathbb R,*})=C(S^{N-1}_{\mathbb R,+})\Big/\Big<x_ix_jx_k=(-1)^{\varepsilon_{ij}+\varepsilon_{jk}+\varepsilon_{ik}}x_kx_jx_i\Big>$$
$$C(S^{N-1,\prime}_{\mathbb C,*})=C(S^{N-1}_{\mathbb C,+})\Big/\Big<z_iz_jz_k=(-1)^{\varepsilon_{ij}+\varepsilon_{jk}+\varepsilon_{ik}}z_kz_jz_i\Big>$$

Thus, we have constructed our spheres, and embeddings, as desired.
\end{proof}

Let us twist now the unitary quantum groups $U$. We would like these to act on the corresponding spheres, $U\curvearrowright S$. Thus, we would like to have morphisms, as follows:
$$\Phi(x_i)=\sum_jx_j\otimes u_{ji}$$

But this leads, via Theorem 13.9, to the following result:

\index{twisted orthogonal group}
\index{twisted unitary group}
\index{twisted rotation group}

\begin{theorem}
We have twisted orthogonal and unitary groups, as follows,
$$\xymatrix@R=15mm@C=15mm{
O_N^+\ar[r]&U_N^+\\
O_N'\ar[r]\ar[u]&U_N'\ar[u]}$$
defined via the following relations, with the convention $\alpha=a,a^*$ and $\beta=b,b^*$:
$$\alpha\beta=\begin{cases}
-\beta\alpha&{\rm for}\ a,b\in\{u_{ij}\}\ {\rm distinct,\ on\ the\ same\ row\ or\ column\ of\ }u\\
\beta\alpha&{\rm otherwise}
\end{cases}$$
These quantum groups act on the corresponding twisted real and complex spheres.
\end{theorem}

\begin{proof}
This is something routine, the idea being as follows:

\medskip

(1) Let us first discuss the construction of the quantum group $O_N'$. We must prove that the algebra $C(O_N')$ obtained from $C(O_N^+)$ via the relations in the statement has a comultiplication $\Delta$, a counit $\varepsilon$, and an antipode $S$. Regarding $\Delta$, let us set:
$$U_{ij}=\sum_ku_{ik}\otimes u_{kj}$$

(2) For $j\neq k$ we have the following computation:
\begin{eqnarray*}
U_{ij}U_{ik}
&=&\sum_{s\neq t}u_{is}u_{it}\otimes u_{sj}u_{tk}+\sum_su_{is}u_{is}\otimes u_{sj}u_{sk}\\
&=&\sum_{s\neq t}-u_{it}u_{is}\otimes u_{tk}u_{sj}+\sum_su_{is}u_{is}\otimes(-u_{sk}u_{sj})\\
&=&-U_{ik}U_{ij}
\end{eqnarray*}

Also, for $i\neq k,j\neq l$ we have the following computation:
\begin{eqnarray*}
U_{ij}U_{kl}
&=&\sum_{s\neq t}u_{is}u_{kt}\otimes u_{sj}u_{tl}+\sum_su_{is}u_{ks}\otimes u_{sj}u_{sl}\\
&=&\sum_{s\neq t}u_{kt}u_{is}\otimes u_{tl}u_{sj}+\sum_s(-u_{ks}u_{is})\otimes(-u_{sl}u_{sj})\\
&=&U_{kl}U_{ij}
\end{eqnarray*}

Thus, we can define a comultiplication map for $C(O_N')$, by setting:
$$\Delta(u_{ij})=U_{ij}$$

(3) Regarding now the counit $\varepsilon$ and the antipode $S$, things are clear here, by using the same method, and with no computations needed, the formulae to be satisfied being trivially satisfied. We conclude that $O_N'$ is a compact quantum group, and the proof for $U_N'$ is similar, by adding $*$ exponents everywhere in the above computations. 

\medskip

(4) Finally, the last assertion is clear too, by doing some elementary computations, of the same type as above, and with the remark that the converse holds too, in the sense that if we want a quantum group $U\subset U_N^+$ to be defined by relations of type $ab=\pm ba$, and to have an action $U\curvearrowright S$ on the corresponding twisted sphere, we are led to the relations in the statement. We refer to \cite{ba1} for further details on all this.
\end{proof}

In order to discuss now the half-classical case, given three coordinates $a,b,c\in\{u_{ij}\}$, let us set $span(a,b,c)=(r,c)$, where $r,c\in\{1,2,3\}$ are the number of rows and columns spanned by $a,b,c$. In other words, if we write $a=u_{ij},b=u_{kl},c=u_{pq}$ then $r=\#\{i,k,p\}$ and $l=\#\{j,l,q\}$. With this convention, we have the following result:

\index{twisted half-classical group}

\begin{theorem}
We have intermediate quantum groups as follows,
$$\xymatrix@R=12mm@C=12mm{
O_N^+\ar[r]&\mathbb TO_N^+\ar[r]&U_N^+\\
O_N^{*\prime}\ar[r]\ar[u]&\mathbb TO_N^{*\prime}\ar[r]\ar[u]&U_N^{*\prime}\ar[u]\\
O_N'\ar[r]\ar[u]&\mathbb TO_N'\ar[r]\ar[u]&U_N'\ar[u]}$$
defined via the following relations, with $\alpha=a,a^*$, $\beta=b,b^*$ and $\gamma=c,c^*$,
$$\alpha\beta\gamma=\begin{cases}
-\gamma\beta\alpha&{\rm for}\ a,b,c\in\{u_{ij}\}\ {\rm with}\ span(a,b,c)=(\leq 2,3)\ {\rm or}\ (3,\leq 2)\\
\gamma\beta\alpha&{\rm otherwise}
\end{cases}$$
which act on the corresponding twisted half-classical real and complex spheres.
\end{theorem}

\begin{proof}
We use the same method as for Theorem 13.10, but with the combinatorics being now more complicated. Observe first that the rules for the various commutation and anticommutation signs in the statement can be summarized as follows:
$$\begin{matrix}
r\backslash c&1&2&3\\
1&+&+&-\\
2&+&+&-\\
3&-&-&+
\end{matrix}$$

Let us first prove the result for $O_N^{*\prime}$. We must construct here morphisms $\Delta,\varepsilon,S$, and the proof, similar to the proof of Theorem 13.10, goes as follows:

\medskip

(1) We first construct $\Delta$. For this purpose, we must prove that $U_{ij}=\sum_ku_{ik}\otimes u_{kj}$ satisfy the relations in the statement. We have the following computation:
\begin{eqnarray*}
U_{ia}U_{jb}U_{kc}
&=&\sum_{xyz}u_{ix}u_{jy}u_{kz}\otimes u_{xa}u_{yb}u_{zc}\\
&=&\sum_{xyz}\pm u_{kz}u_{jy}u_{ix}\otimes\pm u_{zc}u_{yb}u_{xa}\\
&=&\pm U_{kc}U_{jb}U_{ia}
\end{eqnarray*}

We must show that, when examining the precise two $\pm$ signs in the middle formula, their product produces the correct $\pm$ sign at the end. But the point is that both these signs depend only on $s=span(x,y,z)$, and for $s=1,2,3$ respectively, we have:

\medskip

-- For a $(3,1)$ span we obtain $+-$, $+-$, $-+$, so a product $-$ as needed.

\smallskip

-- For a $(2,1)$ span we obtain $++$, $++$, $--$, so a product $+$ as needed.

\smallskip

-- For a $(3,3)$ span we obtain $--$, $--$, $++$, so a product $+$ as needed.

\smallskip

-- For a $(3,2)$ span we obtain $+-$, $+-$, $-+$, so a product $-$ as needed.

\smallskip

-- For a $(2,2)$ span we obtain $++$, $++$, $--$, so a product $+$ as needed.

\medskip

Together with the fact that our problem is invariant under $(r,c)\to(c,r)$, and with the fact that for a $(1,1)$ span there is nothing to prove, this finishes the proof for $\Delta$.

\medskip

(2) The construction of the counit, via the formula $\varepsilon(u_{ij})=\delta_{ij}$, requires the Kronecker symbols $\delta_{ij}$ to commute/anticommute according to the above table. Equivalently, we must prove that the situation $\delta_{ij}\delta_{kl}\delta_{pq}=1$ can appear only in a case where the above table indicates ``+''. But this is clear, because $\delta_{ij}\delta_{kl}\delta_{pq}=1$ implies $r=c$.

\medskip

(3) Finally, the construction of the antipode, via the formula $S(u_{ij})=u_{ji}$, is clear too, because this requires the choice of our $\pm$ signs to be invariant under transposition, and this is true, the above table being symmetric. 

\medskip

(4) We conclude that $O_N^{*\prime}$ is indeed a compact quantum group, and the proof for $U_N^{*\prime}$ is similar, by adding $*$ exponents everywhere in the above. Finally, the last assertion is clear too, exactly as in the proof of Theorem 13.10. We refer to \cite{ba1} for details.
\end{proof}

The above results can be summarized as follows:

\begin{theorem}
We have quantum groups as follows, obtained via the twisted commutation relations $ab=\pm ba$, and twisted half-commutation relations $abc=\pm cba$,
$$\xymatrix@R=12mm@C=12mm{
O_N^+\ar[r]&\mathbb TO_N^+\ar[r]&U_N^+\\
O_N^{*\prime}\ar[r]\ar[u]&\mathbb TO_N^{*\prime}\ar[r]\ar[u]&U_N^{*\prime}\ar[u]\\
O_N'\ar[r]\ar[u]&\mathbb TO_N'\ar[r]\ar[u]&U_N'\ar[u]}$$
with the various signs coming as follows:
\begin{enumerate}
\item The signs for $O_N'$ correspond to anticommutation of distinct entries on rows and columns, and commutation otherwise, with this coming from $O_N'\curvearrowright S^{N-1,\prime}_\mathbb R$.

\item The signs for $O_N^{*\prime},U_N',U_N^{*\prime}$ come as well from the signs for $S^{N-1,\prime}_\mathbb R$, either via the requirement $O_N'\subset U$, or via the requirement $U\curvearrowright S$. 
\end{enumerate}
\end{theorem}

\begin{proof}
This is a summary of Theorem 13.10 and Theorem 13.11, along with a few supplementary facts, coming from the proofs of these results.
\end{proof}

\section*{13b. Schur-Weyl twists}

Let us review now the above construction of the twists, which was something quite ad-hoc, and replace this by something more conceptual. Let us start with:

\begin{proposition}
The intermediate easy quantum groups 
$$H_N\subset G\subset U_N^+$$
come via Tannakian duality from the intermediate categories of partitions
$$P_{even}\supset D\supset\mathcal {NC}_2$$
with $P_{even}(k,l)\subset P(k,l)$ being the category of partitions whose blocks have even size.
\end{proposition}

\begin{proof}
This is something coming from the general easiness theory. Indeed, the easy quantum groups appear as certain intermediate compact quantum groups, as follows:
$$S_N\subset G\subset U_N^+$$

To be more precise, such a quantum group is easy when the corresponding Tannakian category comes from an intermediate category of partitions, as follows:
$$P\supset D\supset\mathcal {NC}_2$$

Now since this correspondence makes correspond $H_N\leftrightarrow P_{even}$, once again as explained in chapter 2, we are led to the conclusion in the statement.
\end{proof}

The idea now will be that the twisting operation $G\to\bar{G}$, in the easy case, can be implemented, via Tannakian duality as usual, via a signature operation on $P_{even}$. Given a partition $\tau\in P(k,l)$, let us call ``switch'' the operation which consists in switching two neighbors, belonging to different blocks, in the upper row, or in the lower row. Also, we use the standard embedding $S_k\subset P_2(k,k)$, via the pairings having only up-to-down strings. With these conventions, we have the following result, from \cite{ba1}:

\index{signature of partitions}
\index{number of crossings}
\index{signature map}
\index{nocrossing form}

\begin{theorem}
There is a signature map $\varepsilon:P_{even}\to\{-1,1\}$, given by 
$$\varepsilon(\tau)=(-1)^c$$
where $c$ is the number of switches needed to make $\tau$ noncrossing. In addition:
\begin{enumerate}
\item For $\tau\in S_k$, this is the usual signature.

\item For $\tau\in P_2$ we have $(-1)^c$, where $c$ is the number of crossings.

\item For $\tau\leq\pi\in NC_{even}$, the signature is $1$.
\end{enumerate}
\end{theorem}

\begin{proof}
In order to show that the signature map $\varepsilon:P_{even}\to\{-1,1\}$ in the statement, given by $\varepsilon(\tau)=(-1)^c$, is well-defined, we must prove that the number $c$ in the statement is well-defined modulo 2. It is enough to perform the verification for the noncrossing partitions. More precisely, given $\tau,\tau'\in NC_{even}$ having the same block structure, we must prove that the number of switches $c$ required for the passage $\tau\to\tau'$ is even.

\medskip

In order to do so, observe that any partition $\tau\in P(k,l)$ can be put in ``standard form'', by ordering its blocks according to the appearence of the first leg in each block, counting clockwise from top left, and then by performing the switches as for block 1 to be at left, then for block 2 to be at left, and so on. Here the required switches are also uniquely determined, by the order coming from counting clockwise from top left. 

\medskip

Here is an example of such an algorithmic switching operation, with block 1 being first put at left, by using two switches, then with block 2 left unchanged, and then with block 3 being put at left as well, but at right of blocks 1 and 2, with one switch:
$$
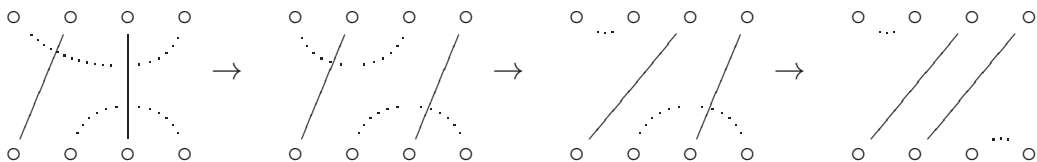
\xymatrix@R=3mm@C=3mm{\circ\ar@/_/@{.}[drr]&\circ\ar@{-}[dddl]&\circ\ar@{-}[ddd]&\circ\\
&&\ar@/_/@{.}[ur]&\\
&&\ar@/^/@{.}[dr]&\\
\circ&\circ\ar@/^/@{.}[ur]&\circ&\circ}
\xymatrix@R=4mm@C=1mm{&\\\to\\&\\& }
\xymatrix@R=3mm@C=3mm{\circ\ar@/_/@{.}[dr]&\circ\ar@{-}[dddl]&\circ&\circ\ar@{-}[dddl]\\
&\ar@/_/@{.}[ur]&&\\
&&\ar@/^/@{.}[dr]&\\
\circ&\circ\ar@/^/@{.}[ur]&\circ&\circ}
\xymatrix@R=4mm@C=1mm{&\\\to\\&\\&}
\xymatrix@R=3mm@C=3mm{\circ\ar@/_/@{.}[r]&\circ&\circ\ar@{-}[dddll]&\circ\ar@{-}[dddl]\\
&&&\\
&&\ar@/^/@{.}[dr]&\\
\circ&\circ\ar@/^/@{.}[ur]&\circ&\circ}
\xymatrix@R=4mm@C=1mm{&\\\to\\&\\& }
\xymatrix@R=3mm@C=3mm{\circ\ar@/_/@{.}[r]&\circ&\circ\ar@{-}[dddll]&\circ\ar@{-}[dddll]\\
&&&\\
&&&\\
\circ&\circ&\circ\ar@/^/@{.}[r]&\circ}$$

The point now is that, under the assumption $\tau\in NC_{even}(k,l)$, each of the moves required for putting a leg at left, and hence for putting a whole block at left, requires an even number of switches. Thus, putting $\tau$ is standard form requires an even number of switches. Now given $\tau,\tau'\in NC_{even}$ having the same block structure, the standard form coincides, so the number of switches $c$ required for the passage $\tau\to\tau'$ is indeed even.

\medskip

Regarding now the remaining assertions, these are all elementary:

\medskip

(1) For $\tau\in S_k$ the standard form is $\tau'=id$, and the passage $\tau\to id$ comes by composing with a number of transpositions, which gives the signature. 

\medskip

(2) For a general $\tau\in P_2$, the standard form is of type $\tau'=|\ldots|^{\cup\ldots\cup}_{\cap\ldots\cap}$, and the passage $\tau\to\tau'$ requires $c$ mod 2 switches, where $c$ is the number of crossings. 

\medskip

(3) Assuming that $\tau\in P_{even}$ comes from $\pi\in NC_{even}$ by merging a certain number of blocks, we can prove that the signature is 1 by proceeding by recurrence.
\end{proof}

With the above result in hand, we can now formulate:

\index{twisted linear map}
\index{twisted Kronecker symbol}

\begin{definition}
Associated to any partition $\pi\in P_{even}(k,l)$ is the linear map
$$T_\pi':(\mathbb C^N)^{\otimes k}\to(\mathbb C^N)^{\otimes l}$$
given by the following formula, with $e_1,\ldots,e_N$ being the standard basis of $\mathbb C^N$,
$$T_\pi'(e_{i_1}\otimes\ldots\otimes e_{i_k})=\sum_{j_1\ldots j_l}\delta'_\pi\begin{pmatrix}i_1&\ldots&i_k\\ j_1&\ldots&j_l\end{pmatrix}e_{j_1}\otimes\ldots\otimes e_{j_l}$$
and where $\delta'_\pi\in\{-1,0,1\}$ is $\delta'_\pi=\varepsilon(\tau)$ if $\tau\geq\pi$, and $\delta'_\pi=0$ otherwise, with $\tau=\ker(^i_j)$.
\end{definition}

In other words, what we are doing here is to add signatures to the usual formula of $T_\pi$. Indeed, observe that the usual formula for $T_\pi$ can be written as folllows:
$$T_\pi(e_{i_1}\otimes\ldots\otimes e_{i_k})=\sum_{j:\ker(^i_j)\geq\pi}e_{j_1}\otimes\ldots\otimes e_{j_l}$$

Now by inserting signs, coming from the signature map $\varepsilon:P_{even}\to\{\pm1\}$, we are led to the following formula, which coincides with the one given above:
$$T_\pi'(e_{i_1}\otimes\ldots\otimes e_{i_k})=\sum_{\tau\geq\pi}\varepsilon(\tau)\sum_{j:\ker(^i_j)=\tau}e_{j_1}\otimes\ldots\otimes e_{j_l}$$

We will be back later to this analogy, with more details on what can be done with it. For the moment, we must first prove a key categorical result, as follows:

\begin{proposition}
The assignement $\pi\to T_\pi'$ is categorical, in the sense that
$$T_\pi'\otimes T_\sigma'=T_{[\pi\sigma]}'\quad,\quad 
T_\pi'T_\sigma'=N^{c(\pi,\sigma)}T_{[^\sigma_\pi]}'\quad,\quad
(T_\pi')^*=T_{\pi^*}'$$
where $c(\pi,\sigma)$ are certain positive integers.
\end{proposition}

\begin{proof}
We have to go back to the proof from the untwisted case, from chapter 2, and insert signs. We have to check three conditions, as follows:

\medskip

\underline{1. Concatenation}. In the untwisted case, this was based on the following formula:
$$\delta_\pi\begin{pmatrix}i_1\ldots i_p\\ j_1\ldots j_q\end{pmatrix}
\delta_\sigma\begin{pmatrix}k_1\ldots k_r\\ l_1\ldots l_s\end{pmatrix}
=\delta_{[\pi\sigma]}\begin{pmatrix}i_1\ldots i_p&k_1\ldots k_r\\ j_1\ldots j_q&l_1\ldots l_s\end{pmatrix}$$

In the twisted case, it is enough to check the following formula:
$$\varepsilon\left(\ker\begin{pmatrix}i_1\ldots i_p\\ j_1\ldots j_q\end{pmatrix}\right)
\varepsilon\left(\ker\begin{pmatrix}k_1\ldots k_r\\ l_1\ldots l_s\end{pmatrix}\right)=
\varepsilon\left(\ker\begin{pmatrix}i_1\ldots i_p&k_1\ldots k_r\\ j_1\ldots j_q&l_1\ldots l_s\end{pmatrix}\right)$$

Let us denote by $\tau,\nu$ the partitions on the left, so that the partition on the right is of the form $\rho\leq[\tau\nu]$. Now by switching to the noncrossing form, $\tau\to\tau'$ and $\nu\to\nu'$, the partition on the right transforms into $\rho\to\rho'\leq[\tau'\nu']$. Now since the partition $[\tau'\nu']$ is noncrossing, we can use Theorem 13.14 (3), and we obtain the result.

\medskip

\underline{2. Composition}. In the untwisted case, this was based on the following formula:
$$\sum_{j_1\ldots j_q}\delta_\pi\begin{pmatrix}i_1\ldots i_p\\ j_1\ldots j_q\end{pmatrix}
\delta_\sigma\begin{pmatrix}j_1\ldots j_q\\ k_1\ldots k_r\end{pmatrix}
=N^{c(\pi,\sigma)}\delta_{[^\pi_\sigma]}\begin{pmatrix}i_1\ldots i_p\\ k_1\ldots k_r\end{pmatrix}$$

In order to prove now the result in the twisted case, it is enough to check that the signs match. More precisely, we must establish the following formula:
$$\varepsilon\left(\ker\begin{pmatrix}i_1\ldots i_p\\ j_1\ldots j_q\end{pmatrix}\right)
\varepsilon\left(\ker\begin{pmatrix}j_1\ldots j_q\\ k_1\ldots k_r\end{pmatrix}\right)
=\varepsilon\left(\ker\begin{pmatrix}i_1\ldots i_p\\ k_1\ldots k_r\end{pmatrix}\right)$$

Let $\tau,\nu$ be the partitions on the left, so that the partition on the right is of the form $\rho\leq[^\tau_\nu]$. Our claim is that we can jointly switch $\tau,\nu$ to the noncrossing form. Indeed, we can first switch as for $\ker(j_1\ldots j_q)$ to become noncrossing, and then switch the upper legs of $\tau$, and the lower legs of $\nu$, as for both these partitions to become noncrossing. Now observe that when switching in this way to the noncrossing form, $\tau\to\tau'$ and $\nu\to\nu'$, the partition on the right transforms into $\rho\to\rho'\leq[^{\tau'}_{\nu'}]$. Now since the partition $[^{\tau'}_{\nu'}]$ is noncrossing, we can apply Theorem 13.14 (3), and we obtain the result.

\medskip

\underline{3. Involution}. Here we must prove the following formula:
$$\delta_\pi'\begin{pmatrix}i_1\ldots i_p\\ j_1\ldots j_q\end{pmatrix}=\delta_{\pi^*}'\begin{pmatrix}j_1\ldots j_q\\ i_1\ldots i_p\end{pmatrix}$$

But this is clear from the definition of $\delta_\pi'$, and we are done.
\end{proof}

As a conclusion, our twisted construction $\pi\to T_\pi'$ has all the needed properties for producing quantum groups, via Tannakian duality, and we can now formulate:

\begin{theorem}
Given a category of partitions $D\subset P_{even}$, the construction
$$Hom(u^{\otimes k},u^{\otimes l})=span\left(T_\pi'\Big|\pi\in D(k,l)\right)$$
produces via Tannakian duality a quantum group $G_N'\subset U_N^+$, for any $N\in\mathbb N$.
\end{theorem}

\begin{proof}
This follows indeed from the Tannakian results from chapter 2, exactly as in the easy case, by using this time Proposition 13.16 as technical ingredient. To be more precise, Proposition 13.16 shows that the linear spaces on the right form a Tannakian category, and so the results in chapter 2 apply, and give the result.
\end{proof}

We can unify the easy quantum groups, or at least the examples coming from categories $D\subset P_{even}$, with the quantum groups constructed above, as follows:

\index{Schur-Weyl twist}
\index{twisted quantum group}
\index{q-easy quantum group}
\index{quizzy quantum group}

\begin{definition}
A closed subgroup $G\subset U_N^+$ is called $q$-easy, or quizzy, with deformation parameter $q=\pm1$, when its tensor category appears as follows,
$$Hom(u^{\otimes k},u^{\otimes l})=span\left(\dot{T}_\pi\Big|\pi\in D(k,l)\right)$$
for a certain category of partitions $D\subset P_{even}$, where, for $q=1,-1$:
$$\dot{T}=T,T'$$
The Schur-Weyl twist of $G$ is the quizzy quantum group $G'\subset U_N^+$ obtained via $q\to-q$.
\end{definition}

We first have to check that when applying the Schur-Weyl twisting to the basic unitary quantum groups, we obtain the previous ad-hoc twists. This is indeed the case:

\begin{theorem}
The twisted unitary quantum groups introduced before,
$$\xymatrix@R=12mm@C=12mm{
O_N^+\ar[r]&\mathbb TO_N^+\ar[r]&U_N^+\\
O_N^{*\prime}\ar[r]\ar[u]&\mathbb TO_N^{*\prime}\ar[r]\ar[u]&U_N^{*\prime}\ar[u]\\
O_N'\ar[r]\ar[u]&\mathbb TO_N'\ar[r]\ar[u]&U_N'\ar[u]}$$
appear as Schur-Weyl twists of the basic unitary quantum groups.
\end{theorem}

\begin{proof}
This is something routine, in several steps, as follows:

\medskip

(1) The basic crossing, $\ker\binom{ij}{ji}$ with $i\neq j$, comes from the transposition $\tau\in S_2$, so its signature is $-1$. As for its degenerated version $\ker\binom{ii}{ii}$, this is noncrossing, so here the signature is $1$. We conclude that the linear map associated to the basic crossing is:
$$T'_{\slash\!\!\!\backslash}(e_i\otimes e_j)
=\begin{cases}
-e_j\otimes e_i&{\rm for}\ i\neq j\\
e_j\otimes e_i&{\rm otherwise}
\end{cases}$$

Regarding now the half-classical crossing, namely $\ker\binom{ijk}{kji}$ with $i,j,k$ distinct, our claim is that the signature is once again $-1$. Indeed, this follows by examining the signatures of the various degenerations of this half-classical crossing, and more specifically from the following signature computations, obtained by counting the crossings, in the first case, by switching twice as to put the partition in noncrossing form, in the next 3 cases, and by observing that the partition is noncrossing, in the last case:
$$
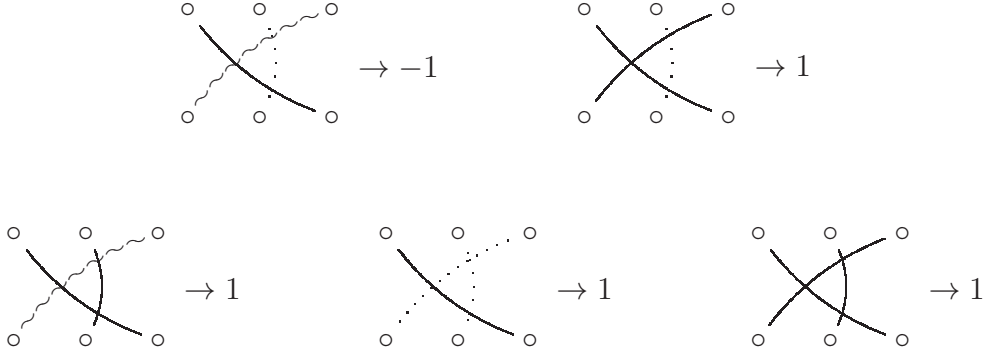
\xymatrix@R=10mm@C=5mm{\circ\ar@/_/@{-}[drr]&\circ\ar@/^/@{.}[d]&\circ\ar@/_/@{~}[dll]\\
\circ&\circ&\circ}
\xymatrix@R=4mm@C=1mm{&\\\to -1\\&\\& }\qquad\quad
\xymatrix@R=10mm@C=5mm{\circ\ar@/_/@{-}[drr]&\circ\ar@/^/@{.}[d]&\circ\ar@/_/@{-}[dll]\\
\circ&\circ&\circ}
\xymatrix@R=4mm@C=1mm{&\\\to 1\\&\\& }$$
$$\xymatrix@R=10mm@C=5mm{\circ\ar@/_/@{-}[drr]&\circ\ar@/^/@{-}[d]&\circ\ar@/_/@{~}[dll]\\
\circ&\circ&\circ}
\xymatrix@R=4mm@C=1mm{&\\\to 1\\&\\& }\qquad\quad
\xymatrix@R=10mm@C=5mm{\circ\ar@/_/@{-}[drr]&\circ\ar@/^/@{.}[d]&\circ\ar@/_/@{.}[dll]\\
\circ&\circ&\circ}
\xymatrix@R=4mm@C=1mm{&\\\to 1\\&\\& }\qquad\quad
\xymatrix@R=10mm@C=5mm{\circ\ar@/_/@{-}[drr]&\circ\ar@/^/@{-}[d]&\circ\ar@/_/@{-}[dll]\\
\circ&\circ&\circ}
\xymatrix@R=4mm@C=1mm{&\\\to 1\\&\\& }$$

Thus, we are led to the following formula, for the half-classical crossing:
$$T'_{\slash\hskip-1.6mm\backslash\hskip-1.1mm|\hskip0.5mm}(e_i\otimes e_j\otimes e_k)
=\begin{cases}
-e_k\otimes e_j\otimes e_i&{\rm for}\ i,j,k\ {\rm distinct}\\
e_k\otimes e_j\otimes e_i&{\rm otherwise}
\end{cases}$$

(2) Our claim now if that for an orthogonal quantum group $G$, the following holds, with the quantum group $O_N'$ being the one in Theorem 13.10:
$$T'_{\slash\!\!\!\backslash}\in End(u^{\otimes 2})\iff G\subset O_N'$$

Indeed, by using the formula of $T'_{\slash\!\!\!\backslash}$ found in (1), we obtain:
\begin{eqnarray*}
(T'_{\slash\!\!\!\backslash}\otimes1)u^{\otimes 2}(e_i\otimes e_j\otimes1)
&=&\sum_ke_k\otimes e_k\otimes u_{ki}u_{kj}\\
&-&\sum_{k\neq l}e_l\otimes e_k\otimes u_{ki}u_{lj}
\end{eqnarray*}

On the other hand, we have as well the following formula:
\begin{eqnarray*}
u^{\otimes 2}(T'_{\slash\!\!\!\backslash}\otimes1)(e_i\otimes e_j\otimes1)
&=&\begin{cases}
\sum_{kl}e_l\otimes e_k\otimes u_{li}u_{ki}&{\rm if}\ i=j\\
-\sum_{kl}e_l\otimes e_k\otimes u_{lj}u_{ki}&{\rm if}\ i\neq j
\end{cases}
\end{eqnarray*}

For $i=j$ the conditions are $u_{ki}^2=u_{ki}^2$ for any $k$, and $u_{ki}u_{li}=-u_{li}u_{ki}$ for any $k\neq l$. For $i\neq j$ the conditions are $u_{ki}u_{kj}=-u_{kj}u_{ki}$ for any $k$, and $u_{ki}u_{lj}=u_{lj}u_{ki}$ for any $k\neq l$. Thus we have exactly the relations between the coordinates of $O_N'$, and we are done.

\medskip

(3) Our claim now if that for an orthogonal quantum group $G$, the following holds, with the quantum group $O_N^{*\prime}$ being the one in Theorem 13.11:
$$T'_{\slash\hskip-1.6mm\backslash\hskip-1.1mm|\hskip0.5mm}\in End(u^{\otimes 3})\iff G\subset O_N^{*\prime}$$

Indeed, by using the formula of $T'_{\slash\hskip-1.6mm\backslash\hskip-1.1mm|\hskip0.5mm}$ found in (1), we obtain:
\begin{eqnarray*}
(T'_{\slash\hskip-1.6mm\backslash\hskip-1.1mm|\hskip0.5mm}\otimes1)u^{\otimes 2}(e_i\otimes e_j\otimes e_k\otimes1)
&=&\sum_{abc\ not\ distinct}e_c\otimes e_b\otimes e_a\otimes u_{ai}u_{bj}u_{ck}\\
&-&\sum_{a,b,c\ distinct}e_c\otimes e_b\otimes e_a\otimes u_{ai}u_{bj}u_{ck}
\end{eqnarray*}

On the other hand, we have as well the following formula:
\begin{eqnarray*}
&&u^{\otimes 2}(T'_{\slash\hskip-1.6mm\backslash\hskip-1.1mm|\hskip0.5mm}\otimes1)(e_i\otimes e_j\otimes e_k\otimes1)\\
&&=\begin{cases}
\sum_{abc}e_c\otimes e_b\otimes e_a\otimes u_{ck}u_{bj}u_{ai}&{\rm for}\ i,j,k\ {\rm not\ distinct}\\
-\sum_{abc}e_c\otimes e_b\otimes e_a\otimes u_{ck}u_{bj}u_{ai}&{\rm for}\ i,j,k\ {\rm distinct}
\end{cases}
\end{eqnarray*}

For $i,j,k$ not distinct the conditions are $u_{ai}u_{bj}u_{ck}=u_{ck}u_{bj}u_{ai}$ for $a,b,c$ not distinct, and $u_{ai}u_{bj}u_{ck}=-u_{ck}u_{bj}u_{ai}$ for $a,b,c$ distinct. For $i,j,k$ distinct the conditions are $u_{ai}u_{bj}u_{ck}=-u_{ck}u_{bj}u_{ai}$ for $a,b,c$ not distinct, and $u_{ai}u_{bj}u_{ck}=u_{ck}u_{bj}u_{ai}$ for $a,b,c$ distinct. Thus we have the relations between the coordinates of $O_N^{*\prime}$, as desired.

\medskip

(4) Now with the above in hand, we obtain that the Schur-Weyl twists of $O_N,O_N^*$ are indeed the quantum groups $O_N',O_N^{*\prime}$ from Theorem 13.10 and Theorem 13.11. 

\medskip

(5) The proof in the unitary case is similar, by adding signs in the above computations (2,3), the conclusion being that the Schur-Weyl twists of $U_N,U_N^*$ are $U_N',U_N^{*\prime}$. 
\end{proof}

\section*{13c. Reflection groups}

Let us clarify now the relation between the maps $T_\pi,T'_\pi$. By using the formulae from the proof of Theorem 13.19, we obtain the following formulae:
$$T'_{\slash\!\!\!\backslash}=-T_{\slash\!\!\!\backslash}+2T_{\ker(^{aa}_{aa})}$$
$$T'_{\slash\hskip-1.6mm\backslash\hskip-1.1mm|\hskip0.5mm}=-\bar{T}_{\slash\hskip-1.6mm\backslash\hskip-1.1mm|\hskip0.5mm}+2T_{\ker(^{aab}_{baa})}+2T_{\ker(^{aba}_{aba})}+2T_{\ker(^{baa}_{aab})}-4T_{\ker(^{aaa}_{aaa})}$$

In general, the answer comes from the M\"obius inversion formula. We recall that the M\"obius function of any lattice, and in particular of $P_{even}$, is given by:
$$\mu(\sigma,\pi)=\begin{cases}
1&{\rm if}\ \sigma=\pi\\
-\sum_{\sigma\leq\tau<\pi}\mu(\sigma,\tau)&{\rm if}\ \sigma<\pi\\
0&{\rm if}\ \sigma\not\leq\pi
\end{cases}$$

With this notation, we have the following useful result:

\index{M\"obius formula}

\begin{proposition}
For any partition $\pi\in P_{even}$ we have the formula
$$T'_\pi=\sum_{\tau\leq\pi}\alpha_\tau T_\tau$$
where $\alpha_\sigma=\sum_{\sigma\leq\tau\leq\pi}\varepsilon(\tau)\mu(\sigma,\tau)$, with $\mu$ being the M\"obius function of $P_{even}$.
\end{proposition}

\begin{proof}
The linear combinations $T=\sum_{\tau\leq\pi}\alpha_\tau T_\tau$ acts on tensors as follows:
\begin{eqnarray*}
T(e_{i_1}\otimes\ldots\otimes e_{i_k})
&=&\sum_{\tau\leq\pi}\alpha_\tau T_\tau(e_{i_1}\otimes\ldots\otimes e_{i_k})\\
&=&\sum_{\tau\leq\pi}\alpha_\tau\sum_{\sigma\leq\tau}\sum_{j:\ker(^i_j)=\sigma}e_{j_1}\otimes\ldots\otimes e_{j_l}\\
&=&\sum_{\sigma\leq\pi}\left(\sum_{\sigma\leq\tau\leq\pi}\alpha_\tau\right)\sum_{j:\ker(^i_j)=\sigma}e_{j_1}\otimes\ldots\otimes e_{j_l}
\end{eqnarray*}

Thus, in order to have $T'_\pi=\sum_{\tau\leq\pi}\alpha_\tau T_\tau$, we must have $\varepsilon(\sigma)=\sum_{\sigma\leq\tau\leq\pi}\alpha_\tau$, for any $\sigma\leq\pi$. But this problem can be solved by using the M\"obius inversion formula, and we obtain the numbers $\alpha_\sigma=\sum_{\sigma\leq\tau\leq\pi}\varepsilon(\tau)\mu(\sigma,\tau)$ in the statement.
\end{proof}

With the above results in hand, let us get now to the question of twisting the quantum reflection groups. It is convenient to include in our discussion two more quantum groups, coming from \cite{rw3} and denoted $H_N^{[\infty]},K_N^{[\infty]}$, constructed as follows:

\begin{proposition}
We have quantum groups $H_N^{[\infty]},K_N^{[\infty]}$ as follows, constructed by using the relations $\alpha\beta\gamma=0$ for any $a\neq c$ on the same row or column of $u$:
$$\xymatrix@R=15mm@C=17mm{
K_N\ar[r]&K_N^*\ar[r]&K_N^{[\infty]}\ar[r]&K_N^+\\
H_N\ar[r]\ar[u]&H_N^*\ar[r]\ar[u]&H_N^{[\infty]}\ar[r]\ar[u]&H_N^+\ar[u]}$$
These quantum groups are both easy, with the corresponding categories of partitions, denoted $P_{even}^{[\infty]}\subset P_{even}$ and $\mathcal P_{even}^{[\infty]}\subset\mathcal P_{even}$, being generated by $\eta=\ker(^{iij}_{jii})$.
\end{proposition}

\begin{proof}
This is routine, by using the fact that the relations $\alpha\beta\gamma=0$ in the statement are equivalent to the condition $\eta\in End(u^{\otimes k})$, with $|k|=3$. For details here, and for more on these two quantum groups, which are very interesting objects, and that we have actually already met in the above, we refer to the paper of Raum-Weber \cite{rw3}.
\end{proof}

In order to discuss now the Schur-Weyl twisting of the various quantum reflection groups that we have, we will need the following technical result:

\begin{proposition}
We have the following equalities,
\begin{eqnarray*}
P_{even}^*&=&\left\{\pi\in P_{even}\Big|\varepsilon(\tau)=1,\forall\tau\leq\pi,|\tau|=2\right\}
\\
P_{even}^{[\infty]}&=&\left\{\pi\in P_{even}\Big|\sigma\in P_{even}^*,\forall\sigma\subset\pi\right\}\\
P_{even}^{[\infty]}&=&\left\{\pi\in P_{even}\Big|\varepsilon(\tau)=1,\forall\tau\leq\pi\right\}
\end{eqnarray*}
where $\varepsilon:P_{even}\to\{\pm1\}$ is the signature of even permutations.
\end{proposition}

\begin{proof}
This is routine combinatorics, from \cite{ba1}, \cite{rw3}, the idea being as follows:

\medskip

(1) Given $\pi\in P_{even}$, we have $\tau\leq\pi,|\tau|=2$ precisely when $\tau=\pi^\beta$ is the partition obtained from $\pi$ by merging all the legs of a certain subpartition $\beta\subset\pi$, and by merging as well all the other blocks. Now observe that $\pi^\beta$ does not depend on $\pi$, but only on $\beta$, and that the number of switches required for making $\pi^\beta$ noncrossing is $c=N_\bullet-N_\circ$ modulo 2, where $N_\bullet/N_\circ$ is the number of black/white legs of $\beta$, when labelling the legs of $\pi$ counterclockwise $\circ\bullet\circ\bullet\ldots$ Thus $\varepsilon(\pi^\beta)=1$ holds precisely when $\beta\in\pi$ has the same number of black and white legs, and this gives the result.

\medskip

(2) This simply follows from the equality $P_{even}^{[\infty]}=<\eta>$ coming from Proposition 13.21, by computing $<\eta>$, and for the complete proof here we refer to \cite{rw3}.

\medskip

(3) We use the fact, also from \cite{rw3}, that the relations $g_ig_ig_j=g_jg_ig_i$ are trivially satisfied for real reflections. Thus, we have:
$$P_{even}^{[\infty]}(k,l)=\left\{\ker\begin{pmatrix}i_1&\ldots&i_k\\ j_1&\ldots&j_l\end{pmatrix}\Big|g_{i_1}\ldots g_{i_k}=g_{j_1}\ldots g_{j_l}\ {\rm inside}\ \mathbb Z_2^{*N}\right\}$$

We conclude that $P_{even}^{[\infty]}$ appears from $NC_{even}$ by ``inflating blocks'', in the sense that each $\pi\in P_{even}^{[\infty]}$ can be transformed into a partition $\pi'\in NC_{even}$ by deleting pairs of consecutive legs, belonging to the same block. Now since this operation leaves invariant modulo 2 the number $c\in\mathbb N$ of switches in the definition of the signature, it leaves invariant the signature $\varepsilon=(-1)^c$ itself, and we obtain the inclusion ``$\subset$'' in the statement. Conversely, given $\pi\in P_{even}$ satisfying $\varepsilon(\tau)=1$, $\forall\tau\leq\pi$, our claim is that:
$$\rho\leq\sigma\subset\pi,|\rho|=2\implies\varepsilon(\rho)=1$$

Indeed, let us denote by $\alpha,\beta$ the two blocks of $\rho$, and by $\gamma$ the remaining blocks of $\pi$, merged altogether. We know that the partitions $\tau_1=(\alpha\wedge\gamma,\beta)$, $\tau_2=(\beta\wedge\gamma,\alpha)$, $\tau_3=(\alpha,\beta,\gamma)$ are all even. On the other hand, putting these partitions in noncrossing form requires respectively $s+t,s'+t,s+s'+t$ switches, where $t$ is the number of switches needed for putting $\rho=(\alpha,\beta)$ in noncrossing form. Thus $t$ is even, and we are done. With the above claim in hand, we conclude, by using the second equality in the statement, that we have $\sigma\in P_{even}^*$. Thus we have $\pi\in P_{even}^{[\infty]}$, which ends the proof of ``$\supset$''.
\end{proof}

With the above result in hand, we can now prove:

\begin{theorem}
The basic quantum reflection groups, namely
$$\xymatrix@R=13mm@C=13mm{
H_N^+\ar[r]&\mathbb TH_N^+\ar[r]&K_N^+\\
H_N^*\ar[r]\ar[u]&\mathbb TH_N^*\ar[r]\ar[u]&K_N^*\ar[u]\\
H_N\ar[r]\ar[u]&\mathbb TH_N\ar[r]\ar[u]&K_N\ar[u]}$$
equal their own Schur-Weyl twists.
\end{theorem}

\begin{proof}
This result basically comes from the results that we have:

\medskip

(1) In the real case, the verifications are as follows:

\medskip

-- $H_N^+$. We know from Theorem 13.14 that for $\pi\in NC_{even}$ we have $T'_\pi=T_\pi$, and since we are in the situation $D\subset NC_{even}$, the definitions of $G,G'$ coincide.

\medskip

-- $H_N^{[\infty]}$. Here we can use the same argument as in (1), based this time on the description of $P_{even}^{[\infty]}$ involving the signatures found in Proposition 13.22.

\medskip

-- $H_N^*$. We have $H_N^*=H_N^{[\infty]}\cap O_N^*$, so $H_N^{*\prime}\subset H_N^{[\infty]}$ is the subgroup obtained via the defining relations for $O_N^{*\prime}$. But all the $abc=-cba$ relations defining $H_N^{*\prime}$ are automatic, of type $0=0$, and it follows that $H_N^{*\prime}\subset H_N^{[\infty]}$ is the subgroup obtained via the relations $abc=cba$, for any $a,b,c\in\{u_{ij}\}$. Thus we have $H_N^{*\prime}=H_N^{[\infty]}\cap O_N^*=H_N^*$, as claimed.

\medskip

-- $H_N$. We have $H_N=H_N^*\cap O_N$, and by functoriality, $H'_N=H_N^{*\prime}\cap O_N'=H_N^*\cap O_N'$. But this latter intersection is easily seen to be equal to $H_N$, as claimed.

\medskip

(2) In the complex case the proof is similar, and we refer here to \cite{ba1}.
\end{proof}

In the orthogonal case, we can say more about all this, and we have the following result, fully covering all the easy quantum groups classified in \cite{rw3}:

\begin{proposition}
The twists of the easy quantum groups $H_N\subset G\subset O_N^+$ are as follows:
\begin{enumerate}
\item For $G=O_N,O_N^*$ we obtain $G'=O_N',O_N^{*\prime}$.

\item For $G\neq O_N,O_N^*$ we have $G=G'$.
\end{enumerate}
\end{proposition}

\begin{proof}
We use the classification result in \cite{rw3}, explained in chapter 10. We have to examine the 3 cases left, namely $G=O_N^+,H_N^{\diamond r},H_N^\Gamma$, and the proof goes as follows:

\medskip

(1) Let $G=O_N^+$. We know from the above that for $\pi\in NC_{even}$ we have $T'_\pi=T_\pi$, and since we are in the situation $D\subset NC_{even}$, the definitions of $G,G'$ coincide.

\medskip

(2) Let $G=H_N^{\diamond r}$. We know that the generating partition is:
$$\pi_r=\ker\begin{pmatrix}1&\ldots&r&r&\ldots&1\\1&\ldots&r&r&\ldots&1\end{pmatrix}$$

By symmetry, putting this partition in noncrossing form requires the same number of upper switches and lower switches, and so requires an even number of total switches. Thus $\pi_r$ is even, and the same argument shows in fact that all its subpartitions are even as well. It follows that we have $T_{\pi_r}=T'_{\pi_r}$, and this gives the result.

\medskip

(3) Let $G=H_N^\Gamma$. We denote by $P_{even}^{[\infty]}\subset D\subset P_{even}$ the corresponding category of partitions. According to the description of $P_{even}^{[\infty]}$ worked out in \cite{rw3}, and discussed in the above, this category contains the following type of partition:
$$
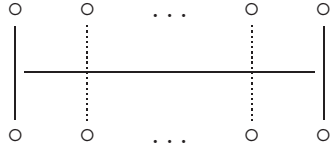
\xymatrix@R=5mm@C=5mm{
\circ\ar@{-}[dd]&\circ\ar@{.}[dd]&\ldots&\circ\ar@{.}[dd]&\circ\ar@{-}[dd]\\
\ar@{-}[rrrr]&&&&\\
\circ&\circ&\ldots&\circ&\circ}$$

The point now is that, by ``capping'' with such partitions, we can merge any pair of blocks of $\pi\in D$, by staying inside $D$. Thus, $D$ has the following property:
$$\tau\leq\pi\in D\implies\tau\in D$$

We deduce from this that $T'_\pi$ is an intertwiner for $G$, and so $G\subset G'$. By symmetry we must have $G'\subset G$ as well, and this finishes the proof.
\end{proof}

As a conclusion to all this, we have:

\begin{theorem}
The easy quantum groups $H_N\subset G\subset O_N^+$ and their twists are
$$
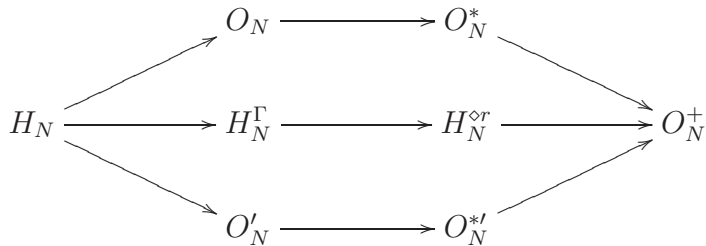
\xymatrix@R=7mm@C=20mm{
&O_N\ar[r]&O_N^*\ar[rd]\\
H_N\ar[r]\ar[ur]\ar[rd]&H_N^\Gamma\ar[r]&H_N^{\diamond r}\ar[r]&O_N^+\\
&O_N'\ar[r]&O_N^{*\prime}\ar[ru]}$$
and the set formed by these quantum groups is stable by intersections.
\end{theorem}

\begin{proof}
The first assertion follows indeed from \cite{rw3}, and from our twisting results. Regarding now the intersection assertion, we have the following intersection diagram:
$$
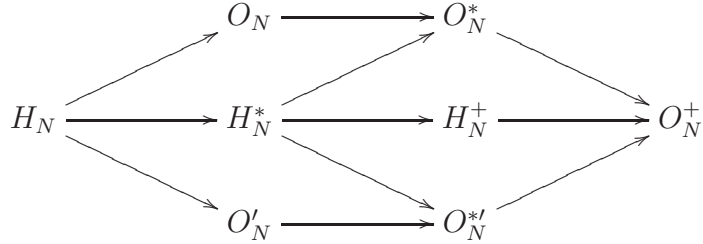
\xymatrix@R=7mm@C=20mm{
&O_N\ar[r]&O_N^*\ar[rd]\\
H_N\ar[r]\ar[ur]\ar[rd]&H_N^*\ar[r]\ar[ur]\ar[dr]&H_N^+\ar[r]&O_N^+\\
&O_N'\ar[r]&O_N^{*\prime}\ar[ru]}$$

More precisely, this diagram has the property that any intersection $G\cap H$ appears on the diagram, as the biggest quantum group contained in both $G,H$. But with this diagram in hand, our claim is that the assertion follows. Indeed, the intersections between the quantum groups $O_N^\times$ are their twists are all on this diagram, and hence on the diagram in the statement as well. Regarding now the intersections of an easy quantum group $H_N\subset G\subset H_N^+$ with the twists $O_N',O_N^{*\prime}$, we can use again the above diagram. Indeed, from $H_N^+\cap O_N^{*\prime}=H_N^*$ we deduce that both $G\cap O_N'$ and $G\cap O_N^{*\prime}$ appear as certain intermediate easy quantum groups $H_N\subset K^\times\subset H_N^*$, and we are done.
\end{proof}

\section*{13d. Integration theory}

Let us discuss now integration questions, over our twisted quantum groups. The result here, valid for any Schur-Weyl twist in our sense, is as follows:

\index{twisted integration}
\index{twisted Weingarten formula}

\begin{theorem}
We have the Weingarten type formula
$$\int_{\bar{G}}u_{i_1j_1}^{e_1}\ldots u_{i_kj_k}^{e_k}=\sum_{\pi,\sigma\in D(k)}\delta'_\pi(i_1\ldots i_k)\delta'_\sigma(j_1\ldots j_k)W_{kN}(\pi,\sigma)$$
where $W_{kN}=G_{kN}^{-1}$, with $G_{kN}(\pi,\sigma)=N^{|\pi\vee\sigma|}$, for $\pi,\sigma\in D(k)$.
\end{theorem}

\begin{proof}
This follows exactly as in the untwisted case, the idea being that the signs will cancel. Let us recall indeed from Definition 13.15 and the comments afterwards that the twisted vectors $\xi'_\pi$ associated to the partitions $\pi\in P_{even}(k)$ are as follows: 
$$\xi'_\pi=\sum_{\tau\geq\pi}\varepsilon(\tau)\sum_{i:\ker(i)=\tau}e_{i_1}\otimes\ldots\otimes e_{i_k}$$

Thus, the Gram matrix of these vectors is given by:
\begin{eqnarray*}
<\xi_\pi',\xi_\sigma'>
&=&\sum_{\tau\geq\pi\vee\sigma}\varepsilon(\tau)^2\left|\left\{(i_1,\ldots,i_k)\Big|\ker i=\tau\right\}\right|\\
&=&\sum_{\tau\geq\pi\vee\sigma}\left|\left\{(i_1,\ldots,i_k)\Big|\ker i=\tau\right\}\right|\\
&=&N^{|\pi\vee\sigma|}
\end{eqnarray*}

Thus the Gram matrix is the same as in the untwisted case, and so the Weingarten matrix is the same as well as in the untwisted case, and this gives the result.
\end{proof}

Summarizing, the integration problematics is not very interesting compared to the one in the classical case, with the difference simply coming from some signs.

\bigskip

Finally, for the sake of completness, let us record as well a number of interesting results, for the most in relation with noncommutative geometry, which is the place where twisting is mostly needed. In relation with tori, we have the following result:

\begin{theorem}
The diagonal tori of the twisted quantum groups are
$$\xymatrix@R=13.5mm@C=13.5mm{
T_N^+\ar[r]&\mathbb TT_N^+\ar[r]&\mathbb T_N^+\\
T_N^*\ar[r]\ar[u]&\mathbb TT_N^*\ar[r]\ar[u]&\mathbb T_N^*\ar[u]\\
T_N\ar[r]\ar[u]&\mathbb TT_N\ar[r]\ar[u]&\mathbb T_N\ar[u]}$$
exactly as in the untwisted case.
\end{theorem}

\begin{proof}
This is clear for the quantum reflection groups, which are not twistable, and for the quantum unitary groups this is elementary as well, coming from definitions.
\end{proof}

Regarding now the twisted spheres, we first have the following result:

\begin{theorem}
The twisted spheres have the following properties:
\begin{enumerate}
\item They have affine actions of the twisted unitary quantum groups.

\item They have unique invariant Haar functionals, which are ergodic.

\item Their Haar functionals are given by Weingarten type formulae.

\item They appear, via the GNS construction, as first row spaces.
\end{enumerate}
\end{theorem}

\begin{proof}
The proofs here are similar to those from the untwisted case, via a routine computation, by adding signs where needed, and with the main technical ingredient, namely the Weingarten formula, being available from Theorem 13.26. See \cite{ba1}.
\end{proof}

We have as well the following result, whose proof is more delicate:

\index{twisted isometry groups}

\begin{theorem}
We have the quantum isometry group formula 
$$U'=G^+(S')$$
in all the $9$ main twisted cases.
\end{theorem}

\begin{proof}
The proofs here are similar to those from the untwisted case, via a routine computation, by adding signs where needed, which amounts in replacing the usual commutators $[a,b]=ab-ba$ by twisted commutators, given by:
$$[[a,b]]=ab+ba$$

There is one subtle point, however, coming from the fact that the linear independence of various products of coordinates of length 1,2,3, which was something clear in the untwisted case, is now a non-trivial question. But this can be solved via a technical application of the Weingarten formula, from Theorem 13.26. 
\end{proof}

We refer to \cite{bbc}, \cite{bgo} and related papers, for more on all this.

\section*{13e. Exercises}

Things have been quite straightforward in this chapter, but as a non-trivial exercise here, going well beyond what was done in the above, we have:

\begin{exercise}
Further extend the quizziness theory developed in this chapter, as to cover some more examples of twists, which are not quizzy in our sense.
\end{exercise}

We will be actually back to such things, in the remainder of this book.

\chapter{Symplectic groups}

\section*{14a. Super-space}

In this chapter we discuss an alternative idea for reaching to super-easiness, by using the approach in \cite{bsk}, motivated by the symplectic groups, $Sp_N\subset U_N$ with $N\in2\mathbb N$, whose first and main particular case, appearing at $N=2$, is a very familiar group, namely:
$$Sp_2=SU_2$$ 

The starting point is the analogy, that we know well since the beginning of this book, between the representation theory of $O_N^+$ and $SU_2$. For our purposes, we will first need a functional analytic approach to $SU_2$. This can be done as follows:

\index{special unitary group}
\index{super-identity}

\begin{theorem}
The algebra of continuous functions on $SU_2$ appears as
$$C(SU_2)=C^*\left((u_{ij})_{i,j=1,2}\Big|u=F\bar{u}F^{-1}={\rm unitary}\right)$$
where $F$ is the following matrix,
$$F=\begin{pmatrix}0&1\\ -1&0\end{pmatrix}$$
called super-identity matrix. 
\end{theorem}

\begin{proof}
This can be done in several steps, as follows:

\medskip

(1) Let us first compute $SU_2$. Consider an arbitrary $2\times2$ complex matrix:
$$U=\begin{pmatrix}a&b\\c&d\end{pmatrix}$$

Assuming $\det U=1$, the unitarity condition $U^{-1}=U^*$ reads:
$$\begin{pmatrix}d&-b\\-c&a\end{pmatrix}
=\begin{pmatrix}\bar{a}&\bar{c}\\\bar{b}&\bar{d}\end{pmatrix}$$

Thus we must have $d=\bar{a}$, $c=-\bar{b}$, and we obtain the following formula:
$$SU_2=\left\{\begin{pmatrix}a&b\\-\bar{b}&\bar{a}\end{pmatrix}\Big|\ |a|^2+|b|^2=1\right\}$$

(2) With the above formula in hand, the fundamental corepresentation of $SU_2$ is:
$$u=\begin{pmatrix}a&b\\-\bar{b}&\bar{a}\end{pmatrix}$$

Now observe that we have the following equality:
$$\begin{pmatrix}a&b\\ -\bar{b}&\bar{a}\end{pmatrix}
\begin{pmatrix}0&1\\ -1&0\end{pmatrix}
=\begin{pmatrix}-b&a\\-\bar{a}&-\bar{b}\end{pmatrix}
=\begin{pmatrix}0&1\\ -1&0\end{pmatrix}
\begin{pmatrix}\bar{a}&\bar{b}\\ -b&a\end{pmatrix}$$

Thus, with $F$ being as in the statement, we have $uF=F\bar{u}$, and so:
$$u=F\bar{u}F^{-1}$$

We conclude that, if $A$ is the universal algebra in the statement, we have:
$$A\to C(SU_2)$$

(3) Conversely now, let us compute the universal algebra $A$ in the statement. For this purpose, let us write its fundamental corepresentation as follows:
$$u=\begin{pmatrix}a&b\\c&d\end{pmatrix}$$

We have $uF=F\bar{u}$, with these quantities being respectively given by:
$$uF=\begin{pmatrix}a&b\\c&d\end{pmatrix}
\begin{pmatrix}0&1\\ -1&0\end{pmatrix}
=\begin{pmatrix}-b&a\\-d&c\end{pmatrix}$$
$$F\bar{u}=\begin{pmatrix}0&1\\ -1&0\end{pmatrix}
\begin{pmatrix}a^*&b^*\\c^*&d^*\end{pmatrix}
=\begin{pmatrix}c^*&d^*\\-a^*&-b^*\end{pmatrix}$$

Thus we must have $d=a^*$, $c=-b^*$, and we obtain the following formula:
$$u=\begin{pmatrix}a&b\\-b^*&a^*\end{pmatrix}$$

We also know that this matrix must be unitary, and we have:
$$uu^*=\begin{pmatrix}a&b\\-b^*&a^*\end{pmatrix}
\begin{pmatrix}a^*&-b\\b^*&a\end{pmatrix}
=\begin{pmatrix}aa^*+bb^*&ba-ab\\a^*b^*-b^*a^*&a^*a+b^*b\end{pmatrix}$$
$$u^*u=\begin{pmatrix}a^*&-b\\b^*&a\end{pmatrix}
\begin{pmatrix}a&b\\-b^*&a^*\end{pmatrix}
=\begin{pmatrix}a^*a+bb^*&a^*b-ba^*\\b^*a-ab^*&aa^*+b^*b\end{pmatrix}$$

Thus, the unitarity equations for $u$ are as follows:
$$aa^*=a^*a=1-bb^*=1-b^*b$$
$$ab=ba,a^*b=ba^*,ab^*=a^*b,a^*b^*=b^*a^*$$

It follows that $a,b,a^*,b^*$ commute, so our algebra is commutative. Now since this algebra is commutative, the involution $*$ becomes the usual conjugation $-$, and so:
$$u=\begin{pmatrix}a&b\\-\bar{b}&\bar{a}\end{pmatrix}$$

But this tells us that we have $A=C(X)$ with $X\subset SU_2$, and so we have a quotient map $C(SU_2)\to A$, which is inverse to the map constructed in (2), as desired.
\end{proof}

Now with the above result in hand, we can see right away the relation with $O_N^+$, and more specifically with $O_2^+$. Indeed, this latter quantum group appears as follows:
$$C(O_2^+)=C^*\left((u_{ij})_{i,j=1,2}\Big|u=\bar{u}={\rm unitary}\right)$$

Thus, $SU_2$ appears from $O_2^+$ by replacing the identity with the super-identity, or perhaps vice versa. In practice now, the idea for unifying is quite clear, namely that of looking at quantum groups appearing via relations as follows: 
$$u=F\bar{u}F^{-1}={\rm unitary}$$

In order to clarify what exact matrices $F\in GL_N(\mathbb C)$ we can use, we must do some computations. Following \cite{ba1}, \cite{bsk}, \cite{bdv}, we first have the following result:

\begin{proposition}
Given a closed subgroup $G\subset U_N^+$, with irreducible fundamental corepresentation $u=(u_{ij})$, this corepresentation is self-adjoint, $u\sim\bar{u}$, precisely when 
$$u=F\bar{u}F^{-1}$$
for some $F\in U_N$, satisfying $F\bar{F}=\pm 1$. Moreover, when $N$ is odd we must have $F\bar{F}=1$. 
\end{proposition}

\begin{proof}
Since $u$ is self-adjoint, $u\sim\bar{u}$, we must have $u=F\bar{u}F^{-1}$, for a certain matrix $F\in GL_N(\mathbb C)$. We obtain from this, by using our assumption that $u$ is irreducible:
\begin{eqnarray*}
u=F\bar{u}F^{-1}
&\implies&\bar{u}=\bar{F}u\bar{F}^{-1}\\
&\implies&u=(F\bar{F})u(F\bar{F})^{-1}\\
&\implies&F\bar{F}=c1\\
&\implies&\bar{F}F=\bar{c}1\\
&\implies&c\in\mathbb R
\end{eqnarray*}

Now by rescaling we can assume $c=\pm1$, so we have proved so far that:
$$F\bar{F}=\pm 1$$

In order to establish now the formula $FF^*=1$, we can proceed as follows:
\begin{eqnarray*}
(id\otimes S)u=u^*
&\implies&(id\otimes S)\bar{u}=u^t\\
&\implies&(id\otimes S)(F\bar{u}F^{-1})=Fu^tF^{-1}\\
&\implies&u^*=Fu^tF^{-1}\\
&\implies&u=(F^*)^{-1}\bar{u}F^*\\
&\implies&\bar{u}=F^*u(F^*)^{-1}\\
&\implies&\bar{u}=F^*F\bar{u}F^{-1}(F^*)^{-1}\\
&\implies&FF^*=d1
\end{eqnarray*}

We have $FF^*>0$, so $d>0$. On the other hand, from $F\bar{F}=\pm 1$, $FF^*=d1$ we get:
$$|\det F|^2=\det(F\bar{F})=(\pm1)^N$$
$$|\det F|^2=\det(FF^*)=d^N$$

Since $d>0$ we obtain from this $d=1$, and so $FF^*=1$ as claimed. We obtain as well that when $N$ is odd the sign must be 1, and so $F\bar{F}=1$, as claimed.
\end{proof}

Once again following Bichon-De Rijdt-Vaes \cite{bdv}, where these questions were first studied, up to an orthogonal base change we can assume that our matrix is as follows, where $N=2p+q$ and $\varepsilon=\pm 1$, with the $1_q$ block at right disappearing if $\varepsilon=-1$:
$$F=\begin{pmatrix}
0&1\ \ \ \\
\varepsilon 1&0_{(0)}\\
&&\ddots\\
&&&0&1\ \ \ \\
&&&\varepsilon 1&0_{(p)}\\
&&&&&1_{(1)}\\
&&&&&&\ddots\\
&&&&&&&1_{(q)}
\end{pmatrix}$$

To be more precise, in the case $\varepsilon=1$, the super-identity is the following matrix:
$$F=\begin{pmatrix}
0&1\ \ \ \\
1&0_{(1)}\\
&&\ddots\\
&&&0&1\ \ \ \\
&&&1&0_{(p)}\\
&&&&&1_{(1)}\\
&&&&&&\ddots\\
&&&&&&&1_{(q)}
\end{pmatrix}$$

In the case $\varepsilon=-1$ now, the diagonal terms vanish, and the super-identity is:
$$F=\begin{pmatrix}
0&1\ \ \ \\
-1&0_{(1)}\\
&&\ddots\\
&&&0&1\ \ \ \\
&&&-1&0_{(p)}
\end{pmatrix}$$

We are therefore led into the following definition, from \cite{bsk}:

\index{super-space}
\index{super-identity}

\begin{definition}
The ``super-space'' $\mathbb C^N_F$ is the usual space $\mathbb C^N$, with its standard basis $\{e_1,\ldots,e_N\}$, with a chosen sign $\varepsilon=\pm 1$, and a chosen involution on the set of indices,
$$i\to\bar{i}$$
with $F$ being the ``super-identity'' matrix, $F_{ij}=\delta_{i\bar{j}}$ for $i\leq j$ and $F_{ij}=\varepsilon\delta_{i\bar{j}}$ for $i\geq j$.
\end{definition}

In what follows we will usually assume that $F$ is the explicit matrix appearing above. Indeed, up to a permutation of the indices, we have a decomposition $n=2p+q$ such that the involution is, in standard permutation notation:
$$(12)\ldots (2p-1,2p)(2p+1)\ldots (q)$$

Let us construct now some basic compact quantum groups, in our ``super'' setting. Once again following \cite{bsk}, let us formulate:

\index{super-orthogonal group}
\index{super-orthogonal quantum group}

\begin{definition}
Associated to the super-space $\mathbb C^N_F$ are $O_F,O_F^+$, given by
$$O_F=\left\{U\in U_N\Big|U=F\bar{U}F^{-1}\right\}$$
$$C(O_F^+)=C^*\left((u_{ij})_{i,j=1,\ldots,n}\Big|u=F\bar{u}F^{-1}={\rm unitary}\right)$$
called super-orthogonal group, and super-orthogonal quantum group.
\end{definition}

As explained in \cite{bsk}, it it possible to considerably extend this list, and we will be back to this, but for our purposes now, this is what we need for the moment. Indeed, as we will see next, Definition 14.3 and Definition 14.4 are all that we need, for including $SU_2$ and the other symplectic groups $Sp_N$ into a generalized easiness theory.

\bigskip

We have indeed the following result, from \cite{bsk}, making the connection with our unification problem for $O_N^+$ and $SU_2$, and more or less solving it:

\index{symplectic group}
\index{quantum symplectic group}
\index{free symplectic group}

\begin{theorem}
The basic orthogonal groups and quantum groups are as follows:
\begin{enumerate}
\item At $\varepsilon=-1$ we have $O_F=Sp_N$ and $O_F^+=Sp_N^+$.

\item At $\varepsilon=-1$ and $N=2$ we have $O_F=O_F^+=SU_2$.

\item At $\varepsilon=1$ we have $O_F=O_N$ and $O_F^+=O_N^+$.
\end{enumerate}
\end{theorem}

\begin{proof}
These results are all elementary, as follows:

\medskip

(1) At $\varepsilon=-1$ this follows from definitions, because the symplectic group $Sp_N\subset U_N$ is by definition the following group:
$$Sp_N=\left\{U\in U_N\Big|U=F\bar{U}F^{-1}\right\}$$

(2) Still at $\varepsilon=-1$, the equation $U=F\bar{U}F^{-1}$ tells us that the symplectic matrices $U\in Sp_N$ are exactly the unitaries $U\in U_N$ which are patterned as follows:
$$U=\begin{pmatrix}
a&b&\ldots\\
-\bar{b}&\bar{a}\\
\vdots&&\ddots
\end{pmatrix}$$

In particular, the symplectic matrices at $N=2$ are as follows:
$$U=\begin{pmatrix}
a&b\\
-\bar{b}&\bar{a}
\end{pmatrix}$$

Thus we have $Sp_2=U_2$, and the formula $Sp_2^+=Sp_2$ is elementary as well.

\medskip

(3) At $\varepsilon=1$ now, consider the root of unity $\rho=e^{\pi i/4}$, and set:
$$J=\frac{1}{\sqrt{2}}\begin{pmatrix}\rho&\rho^7\\ \rho^3&\rho^5\end{pmatrix}$$

This matrix $J$ is then unitary, and we have:
$$J\begin{pmatrix}0&1\\1&0\end{pmatrix}J^t=1$$

Thus the following matrix is unitary as well, and satisfies $KFK^t=1$:
$$K=\begin{pmatrix}J^{(1)}\\&\ddots\\&&J^{(p)}\\&&&1_q\end{pmatrix}$$

Thus in terms of the matrix $V=KUK^*$ we have:
$$U=F\bar{U}F^{-1}={\rm unitary}
\quad\iff\quad V=\bar{V}={\rm unitary}$$

We obtain in this way an isomorphism $O_F^+=O_N^+$ as in the statement, and by passing to classical versions, we obtain as well $O_F=O_N$, as desired.
\end{proof}

Summarizing, we have so far a good idea for defining super-easiness, by using the super-space $\mathbb C_F^N$ instead of the usual space $\mathbb C^N$. There is of course still a lot of work to be done, in order to reach to such a theory, at the representation theory level. Following \cite{bsk} and subsequent papers, we will do this next, directly in a more general setting.

\section*{14b. Formal twists}

In what follows, our idea will be that of replacing the Kronecker symbols $\delta_\pi\in\{0,1\}$ by more general quantities $\bar{\delta}_\pi\in\mathbb T\cup\{0\}$. Our motivation comes from the symplectic group $Sp_N$, which is covered by such a formalism, with $\bar{\delta}_\pi\in\{-1,0,1\}$, and by a number of more technical examples, of ``quantum'' nature, which suggest using $\bar{\delta}_\pi\in\mathbb T\cup\{0\}$. Let us first work out the needed basic algebra. We first have:

\index{generalized Kronecker function}

\begin{definition}
A generalized Kronecker function on a category of partitions $D\subset P$ is a collection of numbers $\bar{\delta}_\pi(^i_j)\in\mathbb T\cup\{0\}$, with $\pi\in D$, such that the formula
$$\bar{T}_\pi(e_{i_1}\otimes\ldots\otimes e_{i_k})=\sum_{j_1\ldots j_l}\bar{\delta}_\pi\begin{pmatrix}i_1&\ldots&i_k\\ j_1&\ldots&j_l\end{pmatrix}
e_{j_1}\otimes\ldots\otimes e_{j_l}$$
defines a correspondence $\pi\to\bar{T}_\pi$ which is categorical, in the sense that we have
$$\bar{T}_{[\pi\sigma]}=\bar{T}_\pi\otimes\bar{T}_\sigma\quad,\quad\bar{T}_{[^\sigma_\pi]}=N^{c(^\sigma_\pi)}\bar{T}_\pi\bar{T}_\sigma\quad,\quad \bar{T}_{\pi^*}=\bar{T}_\pi^*$$
as well as $\bar{T}_{|^{\hskip-1.1mm\circ}_{\hskip-1.1mm\circ}}=\bar{T}_{|^{\hskip-1.1mm\bullet}_{\hskip-1.1mm\bullet}}=id$, and $\bar{T}_{\!\!{\ }_\circ\hskip-1.1mm\cap_{\!\!\bullet}}=\bar{T}_{\!\!{\ }_\bullet\hskip-1.1mm\cap_{\!\!\circ}}=(1\to\sum_ie_i\otimes e_i)$.
\end{definition}

To be more precise here, $N\in\mathbb N$ is given, and if we assume $\pi\in D(k,l)$, then the above numbers $\bar{\delta}_\pi(^i_j)\in\mathbb T\cup\{0\}$ must be defined for any multi-indices $i,j$ having respective lengths $k,l$, and taking their individual indices from the set $\{1,\ldots,N\}$. In more concrete terms now, we have the following description:

\begin{proposition}
$\bar{\delta}:D\to\mathbb T\cup\{0\}$ is a generalized Kronecker function when
$$\bar{\delta}_\pi\begin{pmatrix}i_1&\ldots&i_p\\j_1&\ldots&j_q\end{pmatrix}\bar{\delta}_\sigma\begin{pmatrix}k_1&\ldots&k_r\\l_1&\ldots&l_s\end{pmatrix}
=\bar{\delta}_{[\pi\sigma]}\begin{pmatrix}i_1&\ldots&i_p&k_1&\ldots&k_r\\j_1&\ldots&j_q&l_1&\ldots&l_s\end{pmatrix}$$
$$\sum_{j_1\ldots j_q}\bar{\delta}_\sigma\begin{pmatrix}i_1&\ldots&i_p\\j_1&\ldots&j_q\end{pmatrix}
\bar{\delta}_\pi\begin{pmatrix}j_1&\ldots&j_q\\k_1&\ldots&k_r\end{pmatrix}
=N^{c(_\pi^\sigma)}\bar{\delta}_{[^\sigma_\pi]}\begin{pmatrix}i_1&\ldots&i_p\\k_1&\ldots&k_r\end{pmatrix}$$
$$\bar{\delta}_\pi\begin{pmatrix}i_1&\ldots&i_p\\ j_1&\ldots& j_q\end{pmatrix}
=\overline{\bar{\delta}_{\pi^*}\begin{pmatrix}j_1&\ldots& j_q\\ i_1&\ldots&i_p\end{pmatrix}}$$
and when $\bar{\delta}_\pi=\delta_\pi$ for $\pi=|^{\hskip-1.3mm\circ}_{\hskip-1.3mm\circ},|^{\hskip-1.3mm\bullet}_{\hskip-1.3mm\bullet},{\ }_\circ\hskip-1.25mm\cap_{\!\!\bullet},{\ }_\bullet\hskip-1.25mm\cap_{\!\!\circ}$.
\end{proposition}

\begin{proof}
The proof here is routine, as in \cite{bsp}. We include it, in view of some further use, later on. The concatenation axiom follows from the following computation:
\begin{eqnarray*}
&&(\bar{T}_\pi\otimes\bar{T}_\sigma)(e_{i_1}\otimes\ldots\otimes e_{i_p}\otimes e_{k_1}\otimes\ldots\otimes e_{k_r})\\
&=&\sum_{j_1\ldots j_q}\sum_{l_1\ldots l_s}\bar{\delta}_\pi\begin{pmatrix}i_1&\ldots&i_p\\j_1&\ldots&j_q\end{pmatrix}\bar{\delta}_\sigma\begin{pmatrix}k_1&\ldots&k_r\\l_1&\ldots&l_s\end{pmatrix}e_{j_1}\otimes\ldots\otimes e_{j_q}\otimes e_{l_1}\otimes\ldots\otimes e_{l_s}\\
&=&\sum_{j_1\ldots j_q}\sum_{l_1\ldots l_s}\bar{\delta}_{[\pi\sigma]}\begin{pmatrix}i_1&\ldots&i_p&k_1&\ldots&k_r\\j_1&\ldots&j_q&l_1&\ldots&l_s\end{pmatrix}e_{j_1}\otimes\ldots\otimes e_{j_q}\otimes e_{l_1}\otimes\ldots\otimes e_{l_s}\\
&=&\bar{T}_{[\pi\sigma]}(e_{i_1}\otimes\ldots\otimes e_{i_p}\otimes e_{k_1}\otimes\ldots\otimes e_{k_r})
\end{eqnarray*}

The composition axiom follows from the following computation:
\begin{eqnarray*}
&&\bar{T}_\pi\bar{T}_\sigma(e_{i_1}\otimes\ldots\otimes e_{i_p})\\
&=&\sum_{j_1\ldots j_q}\bar{\delta}_\sigma\begin{pmatrix}i_1&\ldots&i_p\\j_1&\ldots&j_q\end{pmatrix}
\sum_{k_1\ldots k_r}\bar{\delta}_\pi\begin{pmatrix}j_1&\ldots&j_q\\k_1&\ldots&k_r\end{pmatrix}e_{k_1}\otimes\ldots\otimes e_{k_r}\\
&=&\sum_{k_1\ldots k_r}N^{c(_\pi^\sigma)}\bar{\delta}_{[^\sigma_\pi]}\begin{pmatrix}i_1&\ldots&i_p\\k_1&\ldots&k_r\end{pmatrix}e_{k_1}\otimes\ldots\otimes e_{k_r}\\
&=&N^{c(_\pi^\sigma)}\bar{T}_{[^\sigma_\pi]}(e_{i_1}\otimes\ldots\otimes e_{i_p})
\end{eqnarray*}

The involution axiom follows from the following computation:
\begin{eqnarray*}
&&\bar{T}_\pi^*(e_{j_1}\otimes\ldots\otimes e_{j_q})\\
&=&\sum_{i_1\ldots i_p}<\bar{T}_\pi^*(e_{j_1}\otimes\ldots\otimes e_{j_q}),e_{i_1}\otimes\ldots\otimes e_{i_p}>e_{i_1}\otimes\ldots\otimes e_{i_p}\\
&=&\sum_{i_1\ldots i_p}\bar{\delta}_\pi\begin{pmatrix}i_1&\ldots&i_p\\ j_1&\ldots& j_q\end{pmatrix}e_{i_1}\otimes\ldots\otimes e_{i_p}\\
&=&\bar{T}_{\pi^*}(e_{j_1}\otimes\ldots\otimes e_{j_q})
\end{eqnarray*}

Finally, the identity and duality axioms follow from the last conditions in the statement. As for the converse, this follows as well from the above computations.
\end{proof}

At the quantum group level now, we are led to:

\index{generalized easy}

\begin{definition}
Given a category of partitions $D\subset P$, and a generalized Kronecker function $\bar{\delta}:D\to\mathbb T\cup\{0\}$, the formula $C=span(\bar{T}_\pi|\pi\in D)$, where
$$\bar{T}_\pi(e_{i_1}\otimes\ldots\otimes e_{i_k})=\sum_{j_1\ldots j_l}\bar{\delta}_\pi\begin{pmatrix}i_1&\ldots&i_k\\ j_1&\ldots&j_l\end{pmatrix}
e_{j_1}\otimes\ldots\otimes e_{j_l}$$
defines a closed subgroup $\bar{G}\subset U_N^+$. We call such quantum groups generalized easy.
\end{definition}

As basic examples, we have of course the usual easy quantum groups $S_N\subset G\subset U_N^+$. We will see in what follows that there are many other interesting examples. In order to compute such quantum groups, we will use the following result:

\begin{proposition}
Given a Kronecker function $\bar{\delta}:D\to\mathbb T\cup\{0\}$, and assuming $D=<\pi_1,\ldots,\pi_r>$ with $\pi_i\in D(k_i,l_i)$, the associated quantum group is given by
$$C(\bar{G})=C(U_N^+)\Big/\left<\bar{T}_{\pi_i}\in Hom(u^{\otimes k_i},u^{\otimes l_i})\Big|i=1,\ldots,r\right>$$
with the usual rules for the exponents, namely $u^\circ=u,u^\bullet=\bar{u}$ and multiplicativity.
\end{proposition}

\begin{proof}
This follows indeed from Proposition 14.7, and from the fact that, according to Definition 14.6, the correspondence $\pi\to\bar{T}_\pi$ is categorical.
\end{proof}

More concretely now, we can try to compute the quantum groups associated to some basic subcategories $E\subset D$. We first have the following result:
 
\begin{proposition}
Consider a generalized Kronecker function $\bar{\delta}:D\to\mathbb T\cup\{0\}$.
\begin{enumerate}
\item We have $\bar{\delta}=\delta$ on the subcategory ${\mathcal NC}_2\subset D$.

\item If $D={\mathcal NC}_2$, the associated quantum group is $U_N^+$.
\end{enumerate}
\end{proposition}

\begin{proof}
Observe first that we have indeed ${\mathcal NC}_2\subset D$, because $D$ must contain the identity and duality partitions $|^{\hskip-1.3mm\circ}_{\hskip-1.3mm\circ},|^{\hskip-1.3mm\bullet}_{\hskip-1.3mm\bullet},{\ }_\circ\hskip-1.25mm\cap_{\!\!\bullet},{\ }_\bullet\hskip-1.25mm\cap_{\!\!\circ}$, which generate ${\mathcal NC}_2$.

\medskip

(1) This follows indeed from the fact, from Proposition 14.7, that we have $\bar{\delta}_\pi=\delta_\pi$ for the standard generators $\pi=|^{\hskip-1.3mm\circ}_{\hskip-1.3mm\circ},|^{\hskip-1.3mm\bullet}_{\hskip-1.3mm\bullet},{\ }_\circ\hskip-1.25mm\cap_{\!\!\bullet},{\ }_\bullet\hskip-1.25mm\cap_{\!\!\circ}$ of the category ${\mathcal NC}_2$. 

\medskip

(2) This follows from (1).
\end{proof}

In relation now with the standard cube, given a category of partitions $D\subset P$, we can associate to it the following categories of partitions:
$$
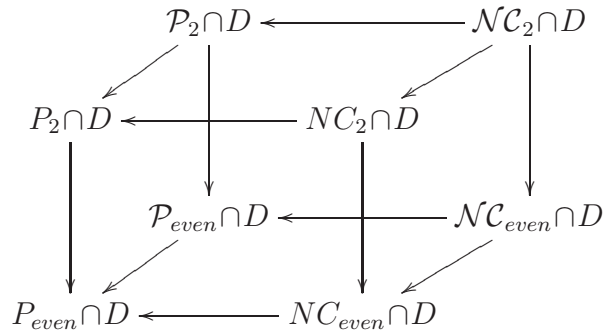
\xymatrix@R=7mm@C=0mm{
&\mathcal P_2\!\cap\!D\ar[dl]\ar[dd]&&\mathcal{NC}_2\!\cap\!D\ar[ll]\ar[dl]\ar[dd]\\
P_2\!\cap\!D\ar[dd]&&NC_2\!\cap\!D\ar[ll]\ar[dd]\\
&\mathcal P_{even}\!\cap\!D\ar[dl]&&\mathcal{NC}_{even}\!\cap\!D\ar[ll]\ar[dl]\\
P_{even}\!\cap\!D&&NC_{even}\!\cap\!D\ar[ll]}$$

Thus, we have 8 basic examples of generalized easy quantum groups, as follows:

\begin{definition}
Given a generalized Kronecker function $\bar{\delta}:D\to\mathbb T\cup\{0\}$, we let
$$
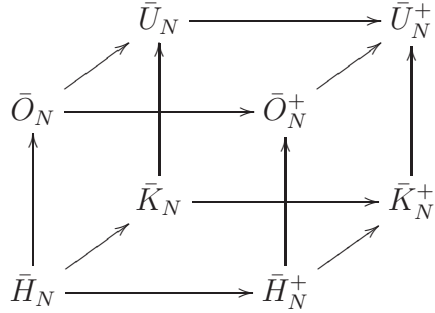
\xymatrix@R=5mm@C=8mm{
&\bar{U}_N\ar[rr]&&\bar{U}_N^+\\
\bar{O}_N\ar[rr]\ar[ur]&&\bar{O}_N^+\ar[ur]\\
&\bar{K}_N\ar[uu]\ar[rr]&&\bar{K}_N^+\ar[uu]\\
\bar{H}_N\ar[uu]\ar[rr]\ar[ur]&&\bar{H}_N^+\ar[uu]\ar[ur]}$$
be the generalized easy quantum groups associated to the above $8$ categories.
\end{definition}

We already know from Proposition 14.10 that we have $\bar{U}_N^+=U_N^+$. We can extend this observation, by analyzing the other quantum groups as well, and we obtain:

\begin{theorem}
The basic quantum unitary and reflection groups are as follows:
\begin{enumerate}
\item We have $\bar{U}_N^+=U_N^+$.

\item If $NC_2\subset D$, then $\bar{O}_N^+\subset U_N^+$ appears from the following intertwiner:
$$\bar{T}_\pi(e_i)=\sum_j\bar{\delta}_\pi\binom{i}{j}e_j\quad:\quad\pi=|^{\hskip-1.3mm\circ}_{\hskip-1.3mm\bullet}$$

\item If $\mathcal P_2\subset D$, then $\bar{U}_N\subset U_N^+$ appears from the following intertwiners: 
$$\bar{T}_\pi(e_i\otimes e_j)=\sum_{kl}\bar{\delta}_\pi\begin{pmatrix}i&j\\ k&l\end{pmatrix}e_k\otimes e_l\quad:\quad\pi=\slash^{\hskip-1mm\circ}_{\hskip-2.6mm\circ}\hskip-2.7mm\backslash^{\hskip-2.6mm\circ}_{\hskip-1mm\circ}\ ,\ \slash^{\hskip-1mm\circ}_{\hskip-2.6mm\circ}\hskip-2.7mm\backslash^{\hskip-2.6mm\bullet}_{\hskip-1mm\bullet}$$

\item If $\mathcal{NC}_{even}\subset D$, then $\bar{K}_N^+\subset U_N^+$ appears from the following intertwiners:
$$\bar{T}_\pi(1)=\sum_i\bar{\delta}_\pi(i\;i\;i\;i)\quad:\quad\pi={\ }_{\hskip-1mm\circ}\hskip-2.2mm\sqcap_{\hskip-1mm\circ}\hskip-2.3mm\sqcap_{\hskip-1mm\bullet}\hskip-1.4mm\sqcap_{\hskip-1mm\bullet}\ ,\,{\ }_{\hskip-1mm\circ}\hskip-2.2mm\sqcap_{\hskip-1mm\bullet}\hskip-2.3mm\sqcap_{\hskip-1mm\circ}\hskip-1.4mm\sqcap_{\hskip-1mm\bullet}\ ,\,{\ }_{\hskip-1mm\circ}\hskip-2.2mm\sqcap_{\hskip-1mm\bullet}\hskip-2.3mm\sqcap_{\hskip-1mm\bullet}\hskip-1.4mm\sqcap_{\hskip-1mm\circ}
\ ,\, {\ }_{\hskip-1mm\bullet}\hskip-2.2mm\sqcap_{\hskip-1mm\circ}\hskip-2.3mm\sqcap_{\hskip-1mm\circ}\hskip-1.4mm\sqcap_{\hskip-1mm\bullet}\ ,\,{\ }_{\hskip-1mm\bullet}\hskip-2.2mm\sqcap_{\hskip-1mm\circ}\hskip-2.3mm\sqcap_{\hskip-1mm\bullet}\hskip-1.4mm\sqcap_{\hskip-1mm\circ}\ ,\,{\ }_{\hskip-1mm\bullet}\hskip-2.2mm\sqcap_{\hskip-1mm\bullet}\hskip-2.3mm\sqcap_{\hskip-1mm\circ}\hskip-1.4mm\sqcap_{\hskip-1mm\circ}$$

\item If $P_2\subset D$ then $\bar{O}_N\subset U_N^+$ appears from the following intertwiners:
$$\bar{T}_\pi\quad:\quad\pi=|^{\hskip-1.3mm\circ}_{\hskip-1.3mm\bullet}\ ,\ \slash^{\hskip-1mm\circ}_{\hskip-2.6mm\circ}\hskip-2.7mm\backslash^{\hskip-2.6mm\circ}_{\hskip-1mm\circ}$$ 

\item If $NC_{even}\subset D$ then $\bar{H}_N^+\subset U_N^+$ appears from the following intertwiners:
$$\bar{T}_\pi\quad:\quad\pi=|^{\hskip-1.3mm\circ}_{\hskip-1.3mm\bullet}\ ,\ {\ }_{\hskip-1mm\circ}\hskip-2.2mm\sqcap_{\hskip-1mm\circ}\hskip-2.3mm\sqcap_{\hskip-1mm\circ}\hskip-1.4mm\sqcap_{\hskip-1mm\circ}$$

\item If $\mathcal P_{even}\subset D$ then $\bar{K}_N\subset U_N^+$ appears from the following intertwiners:
$$\bar{T}_\pi\quad:\quad\pi=\slash^{\hskip-1mm\circ}_{\hskip-2.6mm\circ}\hskip-2.7mm\backslash^{\hskip-2.6mm\circ}_{\hskip-1mm\circ}\ ,\ \slash^{\hskip-1mm\circ}_{\hskip-2.6mm\circ}\hskip-2.7mm\backslash^{\hskip-2.6mm\bullet}_{\hskip-1mm\bullet}\ ,\ {\ }_{\hskip-1mm\circ}\hskip-2.2mm\sqcap_{\hskip-1mm\bullet}\hskip-2.3mm\sqcap_{\hskip-1mm\circ}\hskip-1.4mm\sqcap_{\hskip-1mm\bullet}$$

\item If $P_{even}\subset D$ then $\bar{H}_N\subset U_N^+$ appears from the following intertwiners:
$$\bar{T}_\pi\quad:\quad\pi=|^{\hskip-1.3mm\circ}_{\hskip-1.3mm\bullet}\ ,\ \slash^{\hskip-1mm\circ}_{\hskip-2.6mm\circ}\hskip-2.7mm\backslash^{\hskip-2.6mm\circ}_{\hskip-1mm\circ}\ ,\ {\ }_{\hskip-1mm\circ}\hskip-2.2mm\sqcap_{\hskip-1mm\circ}\hskip-2.3mm\sqcap_{\hskip-1mm\circ}\hskip-1.4mm\sqcap_{\hskip-1mm\circ}$$
\end{enumerate}
In addition, the diagram formed by these quantum groups is a generation diagram.
\end{theorem}

\begin{proof}
All these assertions are elementary, and follow from some well-known presentation results for the categories in question. To be more precise:

\medskip

(1) This is something that we already know, from Proposition 14.10.

\medskip

(2,3) These results follow from the following well-known formulae:
$$NC_2=<|^{\hskip-1.3mm\circ}_{\hskip-1.3mm\bullet}>\quad,\quad
\mathcal P_2=<\slash^{\hskip-1mm\circ}_{\hskip-2.6mm\circ}\hskip-2.7mm\backslash^{\hskip-2.6mm\circ}_{\hskip-1mm\circ}\ ,\ \slash^{\hskip-1mm\circ}_{\hskip-2.6mm\circ}\hskip-2.7mm\backslash^{\hskip-2.6mm\bullet}_{\hskip-1mm\bullet}>$$

(4) It is well-known that, in the uncolored setting, we have $NC_{even}=<\sqcap\hskip-1.9mm\sqcap\hskip-1.9mm\sqcap>$. By restricting the attention to the category $\mathcal{NC}_{even}\subset NC_{even}$ of the partitions which are matching, it follows that this subcategory is generating by the 6 matching colorings of $\sqcap\hskip-1.6mm\sqcap\hskip-1.6mm\sqcap$. Thus, we have the following presentation formula, which gives the result:
$$\mathcal{NC}_{even}=<{\ }_{\hskip-1mm\circ}\hskip-2.2mm\sqcap_{\hskip-1mm\circ}\hskip-2.3mm\sqcap_{\hskip-1mm\bullet}\hskip-1.4mm\sqcap_{\hskip-1mm\bullet}\ ,\,{\ }_{\hskip-1mm\circ}\hskip-2.2mm\sqcap_{\hskip-1mm\bullet}\hskip-2.3mm\sqcap_{\hskip-1mm\circ}\hskip-1.4mm\sqcap_{\hskip-1mm\bullet}\ ,\,{\ }_{\hskip-1mm\circ}\hskip-2.2mm\sqcap_{\hskip-1mm\bullet}\hskip-2.3mm\sqcap_{\hskip-1mm\bullet}\hskip-1.4mm\sqcap_{\hskip-1mm\circ}\ ,\, {\ }_{\hskip-1mm\bullet}\hskip-2.2mm\sqcap_{\hskip-1mm\circ}\hskip-2.3mm\sqcap_{\hskip-1mm\circ}\hskip-1.4mm\sqcap_{\hskip-1mm\bullet}\ ,\,{\ }_{\hskip-1mm\bullet}\hskip-2.2mm\sqcap_{\hskip-1mm\circ}\hskip-2.3mm\sqcap_{\hskip-1mm\bullet}\hskip-1.4mm\sqcap_{\hskip-1mm\circ}\ ,\,{\ }_{\hskip-1mm\bullet}\hskip-2.2mm\sqcap_{\hskip-1mm\bullet}\hskip-2.3mm\sqcap_{\hskip-1mm\circ}\hskip-1.4mm\sqcap_{\hskip-1mm\circ}>$$

(5,6) These follow from the following presentations results, which are well-known:
$$P_2=<|^{\hskip-1.3mm\circ}_{\hskip-1.3mm\bullet}\ ,\ \slash^{\hskip-1mm\circ}_{\hskip-2.6mm\circ}\hskip-2.7mm\backslash^{\hskip-2.6mm\circ}_{\hskip-1mm\circ}>\quad,\quad 
NC_{even}=<|^{\hskip-1.3mm\circ}_{\hskip-1.3mm\bullet}\ ,\ {\ }_{\hskip-1mm\circ}\hskip-2.2mm\sqcap_{\hskip-1mm\circ}\hskip-2.3mm\sqcap_{\hskip-1mm\circ}\hskip-1.4mm\sqcap_{\hskip-1mm\circ}>$$

(7,8) Indeed, by proceeding as in the proof of (4), we conclude, starting from the above presentation results, that we have the following presentation results as well:
$$\mathcal P_{even}=<\slash^{\hskip-1mm\circ}_{\hskip-2.6mm\circ}\hskip-2.7mm\backslash^{\hskip-2.6mm\circ}_{\hskip-1mm\circ}\ ,\ \slash^{\hskip-1mm\circ}_{\hskip-2.6mm\circ}\hskip-2.7mm\backslash^{\hskip-2.6mm\bullet}_{\hskip-1mm\bullet}\ ,\ {\ }_{\hskip-1mm\circ}\hskip-2.2mm\sqcap_{\hskip-1mm\bullet}\hskip-2.3mm\sqcap_{\hskip-1mm\circ}\hskip-1.4mm\sqcap_{\hskip-1mm\bullet}>\quad,\quad 
P_{even}=<|^{\hskip-1.3mm\circ}_{\hskip-1.3mm\bullet}\ ,\ \slash^{\hskip-1mm\circ}_{\hskip-2.6mm\circ}\hskip-2.7mm\backslash^{\hskip-2.6mm\circ}_{\hskip-1mm\circ}\ ,\ {\ }_{\hskip-1mm\circ}\hskip-2.2mm\sqcap_{\hskip-1mm\circ}\hskip-2.3mm\sqcap_{\hskip-1mm\circ}\hskip-1.4mm\sqcap_{\hskip-1mm\circ}>$$

Finally, the last assertion is clear, because the diagram of the categories producing our quantum groups is an intersection diagram.
\end{proof}

The basic example of a generalized Kronecker function is the usual Kronecker function $\delta:P\to\{0,1\}$. Besides being defined on the whole $P$, and taking only positive values, this function $\delta_\pi(^i_j)$ has the remarkable property of depending only on $\ker(^i_j)\in P$. This suggests formulating the following definition:

\index{pure Kronecker function}

\begin{definition}
A generalized Kronecker function $\bar{\delta}:D\to\mathbb T\cup\{0\}$ is called pure when it is given by a formula of type $\bar{\delta}_\pi(^i_j)=\varphi(\pi,[^i_j])$, where
$$\begin{bmatrix}i_1&\ldots&i_k\\ j_1&\ldots&j_l\end{bmatrix}
=\ker\begin{pmatrix}i_1&\ldots&i_k\\ j_1&\ldots&j_l\end{pmatrix}$$
and where $\varphi:D\times P\to\mathbb T\cup\{0\}$ is a certain function.
\end{definition}

As already mentioned, the usual Kronecker function is pure. In fact, we have the following formula, where $\leq$ is the order relation on $P$ obtained by merging blocks:
$$\delta_\pi\begin{pmatrix}i_1&\ldots&i_k\\ j_1&\ldots&j_l\end{pmatrix}=
\begin{cases}
1&{\rm if}\ \begin{bmatrix}i_1&\ldots&i_k\\ j_1&\ldots&j_l\end{bmatrix}\leq\pi\\
0&{\rm otherwise}
\end{cases}$$

Another interesting example is the ``twisted'' version of $\delta$. Consider indeed the standard embeddings $S_\infty\subset P_2\subset P_{even}$, with the convention that the permutations act vertically, from top to bottom. We have then the following result, from \cite{ba1}:

\begin{proposition}
The signature of the permutations $S_\infty\to\{\pm1\}$ extends into a signature map $\varepsilon_2:P_{even}\to\{\pm1\}$, given by $\varepsilon_2(\pi)=(-1)^{c(\pi)}$, 
where $c(\pi)$ is the number of switches needed for putting $\pi$ in noncrossing form, and the formula
$$\bar{\delta}_\pi\begin{pmatrix}i_1&\ldots&i_k\\ j_1&\ldots&j_l\end{pmatrix}=\delta_\pi\begin{pmatrix}i_1&\ldots&i_k\\ j_1&\ldots&j_l\end{pmatrix}\varepsilon_2\begin{bmatrix}i_1&\ldots&i_k\\ j_1&\ldots&j_l\end{bmatrix}$$
defines a generalized Kronecker map on $P_{even}$, which is pure.
\end{proposition}

\begin{proof}
The idea indeed is that the number $c(\pi)$ in the statement is well-defined modulo 2, and so we have a signature map as above, extending the usual signature of the permutations. The proof of the categorical axioms is routine as well, see \cite{ba1}.
\end{proof}

At the quantum group level now, the result, also from \cite{ba1}, is as follows:

\begin{proposition}
The basic quantum unitary and reflection groups associated to the generalized Kronecker map constructed above are as follows,
$$
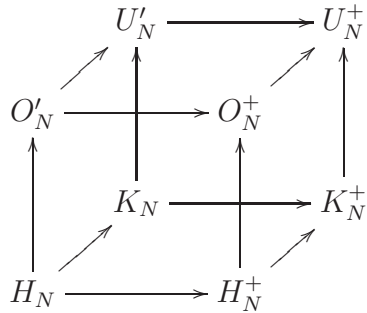
\xymatrix@R=5mm@C=5mm{
&U_N'\ar[rr]&&U_N^+\\
O_N'\ar[rr]\ar[ur]&&O_N^+\ar[ur]\\
&K_N\ar[uu]\ar[rr]&&K_N^+\ar[uu]\\
H_N\ar[uu]\ar[rr]\ar[ur]&&H_N^+\ar[uu]\ar[ur]}$$
with prime denoting $q=-1$ twists, obtained by stating that the standard coordinates and their adjoints anticommute on rows and columns, and commute otherwise.
\end{proposition}

\begin{proof}
The signature map $\varepsilon_2$ being trivial on $NC_{even}$, the free quantum groups, namely $O_N^+,U_N^+,H_N^+,K_N^+$, are not twistable. Regarding now $\bar{O}_N,\bar{U}_N$, the point here is that the linear map associated to the basic crossing is as follows:
$$\bar{T}_{\slash\!\!\!\backslash}(e_i\otimes e_j)
=\begin{cases}
-e_j\otimes e_i&{\rm for}\ i\neq j\\
e_i\otimes e_i&{\rm for}\ i=j
\end{cases}$$

Thus, the corresponding quantum groups $\bar{O}_N,\bar{U}_N$ are obtained from $O_N,U_N$ by replacing the commutation relations $ab=ba$ between coordinates by relations of type $ab=\pm ba$, and so appear as $q$-deformations of $O_N,U_N$, with deformation parameter $q=-1$. As for the remaining quantum groups, $H_N,K_N$, these are not twistable. The proof here is quite tricky, using some combinatorics in order to express the twisted maps $T_\pi$ in terms of the untwisted ones, via a M\"obius inversion type formula. See \cite{ba1}, \cite{ba2}.
\end{proof}

Now back to our considerations, regarding the pure Kronecker functions, the above construction suggests looking into functions of the following type:
$$\bar{\delta}_\pi\begin{pmatrix}i_1&\ldots&i_k\\ j_1&\ldots&j_l\end{pmatrix}=\delta_\pi\begin{pmatrix}i_1&\ldots&i_k\\ j_1&\ldots&j_l\end{pmatrix}\varepsilon\begin{bmatrix}i_1&\ldots&i_k\\ j_1&\ldots&j_l\end{bmatrix}$$

In order to define such functions, we only need to define $\varepsilon$ on the partitions of type $\nu=\ker(^i_j)$ which satisfy $\delta_\pi(^i_j)=0$, which means $\nu\leq\pi$. Thus, we are in need of:

\begin{definition}
The completion of a category of partitions $D$ is given by:
$$\bar{D}=\left\{\nu\Big|\exists\pi\in D,\nu\leq\pi\right\}$$
In other words, $\bar{D}$ is obtained from $D$ by allowing the joining of blocks. 
\end{definition}

Observe that $\bar{D}$ is a category of partitions as well. This follows indeed from definitions. As a basic example, the completion of $D=P_2$ is the category $\bar{D}=P_{even}$.

\bigskip

We have the following axioms for the generalized signatures:

\begin{definition}
A pre-signature on a complete category $\bar{D}$ is a map $\varepsilon:\bar{D}\to\mathbb T$ which is trivial on $\bar{D}\cap NC$, and which satisfies the following conditions:
\begin{enumerate}
\item $\varepsilon(\rho)=\varepsilon(\pi)\varepsilon(\sigma)$, with $\rho\leq[\pi\sigma]$ being obtained by joining left and right blocks.

\item $\varepsilon(\rho)=\varepsilon(\pi)\varepsilon(\sigma)$, with $\rho\leq[^\sigma_\pi]$ being obtained by joining up and down blocks.
\end{enumerate}
In the case where  $\varepsilon(\pi^*)=\bar{\varepsilon}(\pi)$ for any $\pi$, we call $\varepsilon$ a signature.
\end{definition}

As basic examples, we have the trivial signature $\varepsilon_1:P\to\{1\}$, as well as the ``standard'' signature $\varepsilon_2:P_{even}\to\{\pm1\}$. We will see later on that these two maps are particular cases of a pre-signature construction which works at any $s\in\mathbb N$.

\bigskip

Regarding our various axioms for the signatures, the only point which can be probably improved is our assumption that $\varepsilon$ must be trivial on $\bar{D}\cap NC$. This is actually quite a strong assumption, that we will not fully need, but which is verified for all the examples that we have. We believe that this might be actually automatic.

\bigskip

The interest in the signatures comes from the following result:

\begin{proposition}
Given a category of partitions $D$, and a signature on its completion $\varepsilon:\bar{D}\to\mathbb T$, the construction
$$\bar{\delta}_\pi\begin{pmatrix}i_1&\ldots&i_k\\ j_1&\ldots&j_l\end{pmatrix}=\delta_\pi\begin{pmatrix}i_1&\ldots&i_k\\ j_1&\ldots&j_l\end{pmatrix}\varepsilon\begin{bmatrix}i_1&\ldots&i_k\\ j_1&\ldots&j_l\end{bmatrix}$$
produces a generalized Kronecker function, which is pure.
\end{proposition}

\begin{proof}
Our first claim is that a map $\varepsilon:\bar{D}\to\mathbb T$ which is trivial on $\bar{D}\cap NC$ is a signature precisely when it satisfies the following conditions, for any choice of the multi-indices, such that the corresponding kernels belong to $\bar{D}$:
$$\varepsilon\begin{bmatrix}i_1&\ldots&i_p\\ j_1&\ldots&j_q\end{bmatrix}\varepsilon\begin{bmatrix}k_1&\ldots&k_r\\ l_1&\ldots&l_s\end{bmatrix}=
\varepsilon\begin{bmatrix}i_1&\ldots&i_p&k_1&\ldots&k_r\\ j_1&\ldots &j_q&l_1&\ldots&l_s\end{bmatrix}$$
$$\varepsilon\begin{bmatrix}i_1&\ldots&i_p\\ j_1&\ldots&j_q\end{bmatrix}
\varepsilon\begin{bmatrix}j_1&\ldots&j_q\\ k_1&\ldots&k_r\end{bmatrix}
=\varepsilon\begin{bmatrix}i_1&\ldots&i_p\\ k_1&\ldots&k_r\end{bmatrix}$$
$$\varepsilon\begin{bmatrix}i_1&\ldots&i_p\\ j_1&\ldots&j_q\end{bmatrix}=\bar{\varepsilon}\begin{bmatrix}j_1&\ldots&j_q\\ i_1&\ldots&i_p\end{bmatrix}$$

Indeed, these conditions are equivalent to the axioms that we have:

\medskip

(1) With $\{i,j\}\cap\{k,l\}=\emptyset$, our first condition reads $\varepsilon([\pi\sigma])=\varepsilon(\pi)\varepsilon(\sigma)$, where $\pi=[^i_j]$ and $\sigma=[^k_l]$. In the general case now, where $\{i,j\}\cap\{k,l\}$ is not necessarily empty, the partition $\rho=[^{i\,k}_{j\,l}]$ satisfies $\rho\leq[\pi\sigma]$, and is in fact obtained from $[\pi\sigma]$ by joining certain left and right blocks. Moreover, by choosing suitable indices $i,j,k,l$, we see that each such joining is allowed, and we conclude that our first axiom for $\varepsilon$ is equivalent to the condition $\varepsilon(\rho)=\varepsilon(\pi)\varepsilon(\sigma)$, for any such partition $\rho$, as stated.

\medskip

(2) The proof here is similar to the proof of (1), with the horizontal concatenation operation replaced by the vertical concatenation operation, and with the remark that the middle indices $j$ won't interfere, these indices being connected by definition.

\medskip

(3) This is trivial, because we can set $\pi=[^i_j]$, and we obtain the result.

\medskip

Now with this claim in hand, the result follows from our previous result and its proof, because the conditions (1-3), together with our assumption that $\varepsilon$ is trivial on the noncrossing partitions, allow us to insert the signature terms, in the formulae there.
\end{proof}

Summarizing, the signatures produce natural examples of pure Kronecker functions. At the theoretical level, one interesting question is whether any pure Kronecker function appears from a signature, in the above sense. We do not know.

\section*{14c. Super structures}

We recall that we have the following notion, from \cite{bsk}:

\begin{definition}
A super-structure on $\mathbb C^N$ is a linear map $J:\mathbb C^N\to\mathbb C^N$ satisfying $JJ^*=1,J\bar{J}=\pm 1$. Such a map is called normalized when it is of the form 
$$J(e_i)=w_ie_{\tau(i)}$$
with $\tau\in S_N$ being an involution, $\tau^2=id$, and with $w_i\in\{\pm1\}$ being certain signs, which are either trivial, $w_i=1$ for any $i$, or satisfy $w_iw_{\tau(i)}=-1$ for any $i$.
\end{definition}

In relation now with our considerations, our claim is that such a super-structure naturally produces a generalized Kronecker function, and that the associated generalized easy quantum groups are $O_N^{J+}$ and its unitary and discrete versions from \cite{bsk}. Let us first discuss the case $J\bar{J}=1$. Here we have the following result:

\begin{proposition} 
Associated to the map $J(e_i)=e_{\tau(i)}$, with $\tau\in S_N$ being an involution, is a generalized Kronecker function $\delta_\pi^\tau\in\{0,1\}$, defined on $P_{even}$, given by 
$$\delta^\tau_\pi=\prod_{\beta\in\pi}\delta^\tau_\beta$$
and whose values on the one-block partitions $\beta=1^k_l$ are given by $\delta_\beta^\tau(^i_j)=1$ precisely when the multi-index $(^i_j)$ is of the following form,
$$\begin{pmatrix}
\tau^{x_1}(p)&\tau^{\bar{x}_2}(p)&\tau^{x_3}(p)&\tau^{\bar{x}_4}(p)&\ldots\ldots\,\\
\tau^{y_1}(p)&\tau^{\bar{y}_2}(p)&\tau^{y_3}(p)&\tau^{\bar{y}_4}(p)&\ldots\ldots\, 
\end{pmatrix}$$
for a certain $p\in\{1,\ldots,N\}$, where $k=(x_1\ldots x_k)$ and $l=(y_1\ldots y_l)$ are the writings of $k,l$ using $\circ,\bullet$ symbols, with the conventions $\bar{\circ}=\bullet,\bar{\bullet}=\circ$, $\tau^\circ=id,\tau^\bullet=\tau$.
\end{proposition}

\begin{proof}
We have to check our various categorical conditions, for the function $\delta_\pi^\tau$ from the statement. The proof goes as follows:

\medskip

(1) Regarding the tensor product, by using $\delta_\pi^\tau=\prod_{\beta\in\pi}\delta_\beta^\tau$, along with the fact that the blocks of $[\pi\sigma]$ are the blocks of $\pi$ or of $\sigma$ we obtain, as desired:
\begin{eqnarray*}
\delta_\pi^\tau\begin{pmatrix}i_1&\ldots&i_p\\ j_1&\ldots&j_q\end{pmatrix}\delta_\sigma^\tau\begin{pmatrix}k_1&\ldots&k_r\\ l_1&\ldots&l_s\end{pmatrix}
&=&\prod_{\beta\in\pi}\delta_\beta^\tau\begin{pmatrix}i_1&\ldots&i_p\!{\ }_{|\beta}\\ j_1&\ldots&j_q\!{\ }_{|\beta}\end{pmatrix}\prod_{\gamma\in\sigma}\delta_\gamma^\tau\begin{pmatrix}k_1&\ldots&k_r\!{\ }_{|\gamma}\\ l_1&\ldots&l_s\!{\ }_{|\gamma}\end{pmatrix}\\
&=&\prod_{\rho\in[\pi\sigma]}\delta_\rho^\tau\begin{pmatrix}i_1&\ldots&i_p&\!\!k_1&\ldots&k_r\!{\ }_{|\rho}\\ j_1&\ldots&j_q&\!\!l_1&\ldots&l_s\!{\ }_{|\rho}\end{pmatrix}\\
&=&\delta_{[\pi\sigma]}^\tau\begin{pmatrix}i_1&\ldots&i_p&\!\!k_1&\ldots&k_r\\ j_1&\ldots&j_q&\!\!l_1&\ldots&l_s\end{pmatrix}
\end{eqnarray*}

(2) Regarding now the composition, we have to verify here that:
$$\sum_{j_1\ldots j_q}\delta^\tau_\sigma\begin{pmatrix}i_1&\ldots&i_p\\j_1&\ldots&j_q\end{pmatrix}
\delta^\tau_\pi\begin{pmatrix}j_1&\ldots&j_q\\k_1&\ldots&k_r\end{pmatrix}
=N^{c(_\pi^\sigma)}\delta^\tau_{[^\sigma_\pi]}\begin{pmatrix}i_1&\ldots&i_p\\k_1&\ldots&k_r\end{pmatrix}$$

In this formula, the sum on the left equals the cardinality of the following set:
$$S=\left\{j_1,\ldots,j_q\Big|\delta^\tau_\sigma\begin{pmatrix}i_1&\ldots&i_p\\j_1&\ldots&j_q\end{pmatrix}=1,\delta^\tau_\pi\begin{pmatrix}j_1&\ldots&j_q\\k_1&\ldots&k_r\end{pmatrix}=1\right\}$$

In order to compute the cardinality of this set, we have two cases, as follows:

\medskip

\underline{Case $c(_\pi^\sigma)=0$}. In this case there are no middle components when concatenating $\pi,\sigma$, the above multi-indices $j\in S$ are uniquely determined by $\pi,\sigma,i,k$, and so we have $\#S\in\{0,1\}$. Moreover, since these $j$ multi-indices are irrelevant with respect to the question of finding the precise cardinality $\#S\in\{0,1\}$, we have:
$$\#S=\delta^\tau_{[^\sigma_\pi]}\begin{pmatrix}i_1&\ldots&i_p\\k_1&\ldots&k_r\end{pmatrix}$$

Thus, we obtain the formula that we wanted to prove.

\medskip

\underline{Case $c(_\pi^\sigma)>0$}. In this case there are middle components when concatenating $\pi,\sigma$, but we can isolate them when computing $\#S$, by using our condition $\delta_\pi^\tau=\prod_{\beta\in\pi}\delta_\beta^\tau$. Thus, in order to prove the formula in this case, we just have to count the number of free multi-indices $j$, and prove that their number equals $N^{c(^\sigma_\pi)}$. But the free indices for each component must come from a single index $p\in\{1,\ldots,N\}$, via the formula for $\delta_{1^k_l}^\tau$ from the statement, and so we obtain as multiplicity the number $N^{c(^\sigma_\pi)}$, as claimed.

\medskip

(3) Regarding the involution axiom, here we must check that we have:
$$\delta^\tau_\pi\begin{pmatrix}i_1&\ldots&i_p\\ j_1&\ldots& j_q\end{pmatrix}
=\overline{\delta^\tau_{\pi^*}\begin{pmatrix}j_1&\ldots& j_q\\ i_1&\ldots&i_p\end{pmatrix}}$$

But this is clear, because we can use here the condition $\delta_\pi^\tau=\prod_{\beta\in\pi}\delta_\beta^\tau$, together with the formula of $\delta_{1^k_l}^\tau$, which is invariant under upside-down turning.

\medskip

(4) The identity axiom follows from the following computation:
\begin{eqnarray*}
\delta_{\hskip-2mm{\ }_{{\ }_\circ}\hskip-1.45mm|\hskip-3.85mm{\ }^{{\ }^\circ}}^\tau\binom{i}{j}=1
&\iff&\exists p,\binom{i}{j}=\binom{\tau^\circ(p)}{\tau^\circ(p)}\\
&\iff&\exists p,\binom{i}{j}=\binom{p}{p}\\
&\iff&i=j
\end{eqnarray*}

(5) The duality axiom follows from a similar computation, as follows:
\begin{eqnarray*}
\delta_{\!\!{\ }_\circ\hskip-1.1mm\cap\hskip-2.3mm{\ }_\bullet}^\tau(ij)=1
&\iff&\exists p,(i,j)=(\tau^\circ(p),\tau^{\bar{\bullet}}(p))\\
&\iff&\exists p,(i,j)=(\tau^\circ(p),\tau^\circ(p))\\
&\iff&\exists p,(i,j)=(p,p)\\
&\iff&i=j
\end{eqnarray*}

Thus, we have indeed a generalized Kronecker function, as claimed.
\end{proof}

At the quantum group level now, we first have the following result:

\begin{proposition}
The basic quantum unitary and reflection groups associated to the generalized Kronecker map constructed above are as follows,
$$\xymatrix@R=6mm@C=0mm{
&U_N\ar[rr]&&U_N^+\\
O_N^J\ar[rr]\ar[ur]&&O_N^{J+}\ar[ur]\\
&K_N\ar[uu]\ar[rr]&&K_N^+\ar[uu]\\
(K_p\wr\mathbb Z_2)\times K_q\ar[uu]\ar[rr]\ar[ur]&&(K_p^+\wr_*\mathbb Z_2)\,\hat{*}\,K_q^+\ar[uu]\ar[ur]}$$
where $O_N^{J+}$ is the quantum group constructed before, $O_N^J=\{U\in U_N|U=J\bar{U}J^{-1}\}$ is its classical version, and where $\wr_*$ is a free wreath product, and $\hat{*}$ is a dual free product.
\end{proposition}

\begin{proof}
We denote as usual by $\bar{G}_N^\times$ the analogues of the basic easy quantum groups $G_N^\times$, The computation of these quantum groups goes as follows:

\medskip

(1) We know from the above that we have $\bar{U}_N^+=U_N^+$. 

\medskip

(2) In order to identify the quantum group in the statement with $O_N^{J+}$, we have to find inside our category the conditions which imply $u=J\bar{u}J^{-1}$. There are two ways in doing so. First, we can use the identity partition, the computation being as follows: 
\begin{eqnarray*}
\delta_{\hskip-2mm{\ }_{{\ }_\circ}\hskip-1.45mm|\hskip-3.85mm{\ }^{{\ }^\bullet}}^\tau\binom{i}{j}=1
&\iff&\exists p,\binom{i}{j}=\binom{\tau^{\bullet(p)}}{\tau^\circ(p)}\\
&\iff&\exists p,\binom{i}{j}=\binom{\tau(p)}{p}\\
&\iff&i=\tau(j)
\end{eqnarray*}

We can equally use the duality partition, as follows:
\begin{eqnarray*}
\delta^\tau_{\!\!{\ }_\circ\hskip-1.1mm\cap\hskip-2.3mm{\ }_\circ}(ij)=1
&\iff&\exists p,(i,j)=(\tau^\circ(p),\tau^{\bar{\circ}}(p))\\
&\iff&\exists p,(i,j)=(\tau^\circ(p),\tau^\bullet(p))\\
&\iff&\exists p,(i,j)=(p,\tau(p))\\
&\iff&j=\tau(i)
\end{eqnarray*}

Summarizing, in both cases we have reached to the conclusion that we must have $u=J\bar{u}J^{-1}$, and this shows that we have $\bar{O}_N^+=O_N^{J+}$, as claimed.

\medskip

(3) In order to compute $\bar{U}_N$, we must impose to $\bar{U}_N^+=U_N^+$ the intertwining conditions coming from the basic crossing $\slash\hskip-2.1mm\backslash$\,, colored with its possible 4 matching colorings. By using the condition $\delta^\tau_\pi=\prod_{\beta\in\pi}\delta^\tau_\beta$, along with the identity axiom for $\delta^\tau$, and its conjugate, we conclude that for all the 4 matching versions of $\pi=\slash\hskip-2.1mm\backslash$, we have $\bar{T}_\pi=T_\pi$. Now since we have $T_\pi(e_i\otimes e_j)=e_j\otimes e_i$, we obtain in this way the usual commutation relations between the variables $\{u_{ij},u_{ij}^*\}$. Thus $\bar{U}_N$ is classical, and so we have $\bar{U}_N=U_N$.

\medskip

(4) Regarding the analogue of $O_N$, we can compute it by intersecting, as follows: 
$$\bar{O}_N=\bar{U}_N\cap\bar{O}_N^+=U_N\cap O_N^{J+}=O_N^J$$

(5) Let us compute now $\bar{K}_N$. We know that this appears as a subgroup of $\bar{U}_N=U_N$, via the relations coming from the intertwiner $\bar{T}_\pi$, where $\pi=\sqcap\hskip-1.7mm\sqcap\hskip-1.7mm\sqcap$, with any of its 6 matching colorings. The best here is to choose an alternating coloring, as follows:
$$\pi={\ }_{\hskip-1mm\circ}\hskip-2.2mm\sqcap_{\hskip-1mm\bullet}\hskip-2.3mm\sqcap_{\hskip-1mm\circ}\hskip-1.4mm\sqcap_{\hskip-1mm\bullet}$$

By Frobenius duality, we can say as well that $\bar{K}_N\subset U_N$ appears via the relations coming from the intertwiner $\bar{T}_\sigma$, where $\sigma$ is the rotated version of $\pi$, given by:
$$\xymatrix@R=2mm@C=3mm{\\ \\ \sigma\ \ =\\ \\}\ \ \ 
\xymatrix@R=2mm@C=3mm{
\circ\ar@/_/@{-}[dr]&&\bullet\\
&\ar@/_/@{-}[ur]\ar@{-}[dd]\\
\\
&\ar@/^/@{-}[dr]&&\\
\circ\ar@/^/@{-}[ur]&&\bullet}$$

According to our conventions, the linear map associated to this partition is independent of $\tau$, given by the following formula:
$$\bar{T}_\sigma(e_i\otimes e_j)=\delta_{ij}e_i\otimes e_i$$

By using this formula, we obtain the following relations:
\begin{eqnarray*}
\bar{T}_\sigma(u\otimes\bar{u})(e_i\otimes e_j)=\sum_{ab}\delta_{ab}u_{ai}u_{bj}^*\otimes e_a\otimes e_b\\
(u\otimes\bar{u})\bar{T}_\sigma(e_i\otimes e_j)=\sum_{ab}\delta_{ij}u_{ai}u_{bj}^*\otimes e_a\otimes e_b
\end{eqnarray*}

Thus, the condition $\bar{T}_\sigma\in End(u\otimes\bar{u})$, which defines $\bar{K}_N\subset U_N$, is equivalent to:
$$\delta_{ab}u_{ai}u_{bj}^*=\delta_{ij}u_{ai}u_{bj}^*\quad,\quad\forall a,b,i,j$$

But this latter condition holds trivially when $a=b,i=j$ or when $a\neq b,i\neq j$. As for the remaining two cases, namely $a=b,i\neq j$ and $a\neq b,i=j$, here this condition tells us that the distinct matrix entries of any group element $U\in\bar{K}_N$ must have products 0, on each row and each column. We conclude that we have $\bar{K}_N=K_N$, as claimed.

\medskip

(6) Let us compute now $\bar{H}_N$. This group appears as an intersection, as follows:
\begin{eqnarray*}
\bar{H}_N
&=&\bar{K}_N\cap\bar{O}_N\\
&=&K_N\cap O_N^J\\
&=&\left\{U\in K_N|U=J\bar{U}J^{-1}\right\}
\end{eqnarray*}

Now observe that, since the matrices $UJ,J\bar{U}$ are obtained from $U,\bar{U}$ by interchanging the rows and columns $i,i+1$, with $i=1,3,\ldots,2p-1$, the relation $UJ=J\bar{U}$ reads:
$$U=\begin{pmatrix}
a&b&\ldots&v\\
\bar{b}&\bar{a}&\ldots&\bar{v}\\
\ldots&\ldots&\ldots&\ldots\\
w&\bar{w}&\ldots&X
\end{pmatrix}$$

To be more precise, here $a,b,\ldots\in\mathbb C$, $v,\ldots\in\mathbb C^q$ are row vectors, $w,\ldots\in\mathbb C^q$ are column vectors, and $X\in M_q(\mathbb C)$. In the case $U\in K_N$, as only one entry in each row/column may be non-zero, the vectors $v,\ldots$ and $w,\ldots$ must vanish, we must have $X\in K_q$, and all $(^a_{\bar{b}}{\ }^b_{\bar{a}})$ blocks must be either $(^0_0{\ }^0_0)$ or the following form, with $z\in\mathbb T$:
$$\begin{pmatrix}z&0\\ 0&\bar{z}\end{pmatrix}\quad,\quad
\begin{pmatrix}0&z\\ \bar{z}&0\end{pmatrix}$$

Thus, we obtain in this way an identification $\bar{H}_N=(K_p\wr\mathbb Z_2)\times K_q$, as claimed. 

\medskip

(7) In order to compute now $\bar{K}_N^+$, we use the fact that this appears inside $\bar{U}_N^+=U_N^+$ by imposing the relations associated to the partition $\sqcap\hskip-1.7mm\sqcap\hskip-1.7mm\sqcap$, with its six possible matching decorations for the legs, namely $\circ\circ\bullet\,\bullet$, $\circ\bullet\circ\,\bullet$, $\circ\bullet\bullet\,\circ$, $\bullet\circ\circ\,\bullet$, $\bullet\circ\bullet\,\circ$, $\bullet\bullet\circ\,\circ$. By using now Frobenius duality, as in the proof of (5), we are led to the relations coming from the following 6 partitions, obtained by rotating:
$$
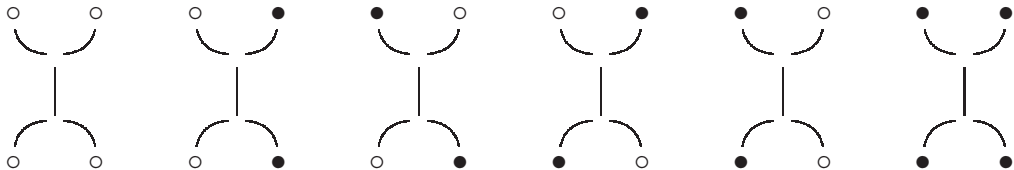
\xymatrix@R=2mm@C=2mm{
\circ\ar@/_/@{-}[dr]&&\circ\\
&\ar@/_/@{-}[ur]\ar@{-}[dd]\\
\\
&\ar@/^/@{-}[dr]&&\\
\circ\ar@/^/@{-}[ur]&&\circ}\quad 
\xymatrix@R=2mm@C=2mm{
\circ\ar@/_/@{-}[dr]&&\bullet\\
&\ar@/_/@{-}[ur]\ar@{-}[dd]\\
\\
&\ar@/^/@{-}[dr]&&\\
\circ\ar@/^/@{-}[ur]&&\bullet}\quad 
\xymatrix@R=2mm@C=2mm{
\bullet\ar@/_/@{-}[dr]&&\circ\\
&\ar@/_/@{-}[ur]\ar@{-}[dd]\\
\\
&\ar@/^/@{-}[dr]&&\\
\circ\ar@/^/@{-}[ur]&&\bullet}\quad 
\xymatrix@R=2mm@C=2mm{
\circ\ar@/_/@{-}[dr]&&\bullet\\
&\ar@/_/@{-}[ur]\ar@{-}[dd]\\
\\
&\ar@/^/@{-}[dr]&&\\
\bullet\ar@/^/@{-}[ur]&&\circ}\quad 
\xymatrix@R=2mm@C=2mm{
\bullet\ar@/_/@{-}[dr]&&\circ\\
&\ar@/_/@{-}[ur]\ar@{-}[dd]\\
\\
&\ar@/^/@{-}[dr]&&\\
\bullet\ar@/^/@{-}[ur]&&\circ}\quad 
\xymatrix@R=2mm@C=2mm{
\bullet\ar@/_/@{-}[dr]&&\bullet\\
&\ar@/_/@{-}[ur]\ar@{-}[dd]\\
\\
&\ar@/^/@{-}[dr]&&\\
\bullet\ar@/^/@{-}[ur]&&\bullet}$$

We already know, from the proof of (5), that the conditions coming from the $2^{nd}$ partition simply state that we must have $ab^*=0$, for any two distinct entries $a,b\in\{u_{ij}\}$ on the same row or the same column. Regarding the $1^{st}$, $5^{th}$, $6^{th}$ partitions, the situation here is similar, with the relations being respectively $ab=0,a^*b=0,a^*b^*=0$, with $a,b\in\{u_{ij}\}$ being as above. As for the remaining relations, coming from the $3^{rd}$, $4^{th}$ partitions, these tell us that the entries of $u_{ij}$ must be normal. Summing up, we have reached to the definition of $K_N^+$, and so we have $\bar{K}_N^+=K_N^+$, as claimed.

\medskip

(8) Finally, the quantum group $\bar{H}_N^+$ appears as an intersection, as follows:
$$\bar{H}_N^+=\bar{K}_N^+\cap\bar{O}_N^+=K_N^+\cap O_N^{J+}$$

In order to compute this intersection, we use the interpretation of the commutation relation $uJ=J\bar{u}$ given in the proof of (6), with all the complex conjugates there replaced of course by adjoints. By reasoning as there, we first conclude that the vectors $v,\ldots$ and $w,\ldots$ must vanish, and so that $X$ must be the fundamental corepresentation of $K_q^+$, and that we are in a dual free product situation. The study of the $2\times2$ blocks is similar, and this gives a free wreath product decomposition $K_p^+\wr_*\mathbb Z_2$ for the first component. Thus, we obtain $\bar{H}_N^+=(K_p^+\wr_*\mathbb Z_2)\,\hat{*}\,K_q^+$, as claimed.
\end{proof}

We can improve the above result by using the following observation, from \cite{bsk}:

\begin{proposition}
With the fundamental corepresentation $V=CUC^*$, where
$$C=\begin{pmatrix}
\Gamma_{(1)}\\
&&\ddots\\
&&&\Gamma_{(p)}\\
&&&&&1_{(1)}\\
&&&&&&\ddots\\
&&&&&&&1_{(q)}
\end{pmatrix}\quad:\quad\Gamma=\frac{1}{\sqrt{2}}\begin{pmatrix}\rho&\rho^7\\ \rho^3&\rho^5\end{pmatrix}\ ,\ \rho=e^{\pi i/4}$$
the relations defining $O_N^{J+}$ become $V=\bar{V}=$ unitary, and so we have $O_N^{J+}\sim O_N^+$.
\end{proposition}

\begin{proof}
Observe that the above matrix $\Gamma$ is unitary, and that we have $\Gamma(^0_1{\ }^1_0)\Gamma^t=1$. Thus the matrix $C$ is unitary as well, and satisfies $CJC^t=1$. It follows that in terms of $V=CUC^*$ the relations $U=J\bar{U}J^{-1}=$ unitary defining $O_N^{J+}$ simply read $V=\bar{V}=$ unitary, so we obtain an isomorphism $O_N^{J+}\simeq O_N^+$, as in the statement. See \cite{bsk}.
\end{proof}

We can now formulate an improved version of Proposition 14.21, as follows:

\begin{theorem}
The basic quantum unitary and reflection groups associated to the generalized Kronecker map constructed above are as follows,
$$
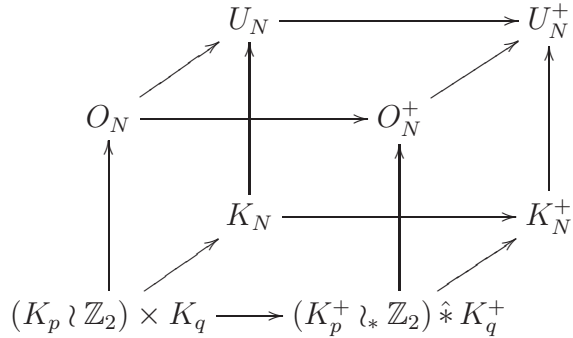
\xymatrix@R=6mm@C=0mm{
&U_N\ar[rr]&&U_N^+\\
O_N\ar[rr]\ar[ur]&&O_N^+\ar[ur]\\
&K_N\ar[uu]\ar[rr]&&K_N^+\ar[uu]\\
(K_p\wr\mathbb Z_2)\times K_q\ar[uu]\ar[rr]\ar[ur]&&(K_p^+\wr_*\mathbb Z_2)\,\hat{*}\,K_q^+\ar[uu]\ar[ur]}$$
where the quantum groups in the upper part of the diagram are taken with respect to the fundamental corepresentation $V=CUC^*$ constructed  above.
\end{theorem}

\begin{proof}
The computation of these quantum groups $\dot{G}_N^\times$ goes as follows:

\medskip

(1) We know from the above that we have $\bar{U}_N^+=U_N^+$. Since the matrix $C$ is unitary, it follows that we have $\dot{U}_N^+=U_N^+$ as well.

\medskip

(2) According to our results, we also have $\bar{O}_N^+=O_N^{J+}$, and so when changing the fundamental corepresentation, we obtain $\dot{O}_N^+=O_N^+$. 

\medskip

(3) We also know that the subgroups $\bar{U}_N\subset\bar{U}_N^+$ and $\bar{O}_N\subset\bar{O}_N^+$ simply appear by taking the classical version. Thus, the subgroups $\dot{U}_N\subset\dot{U}_N^+$ and $\dot{O}_N\subset\dot{O}_N^+$ appear as well by taking the classical version, and so $\dot{U}_N=U_N$, $\dot{O}_N=O_N$.

\medskip

(4) Finally, the quantum groups on the bottom are those from Proposition 14.21, with the remark of course that the 4 vertical embeddings are not the standard ones.
\end{proof}

\section*{14d. Symplectic groups}

We discuss now the case $J\bar{J}=-1$, which covers the group $SU_2$, and more generally the symplectic groups $Sp_N\subset U_N$. In this case, we have: 

\begin{proposition} 
Associated to a negative super-structure, $J(e_i)=w_ie_{\tau(i)}$, is the generalized Kronecker function $\delta_\pi^J\in\{-1,0,1\}$, defined on $P_{even}$, given by 
$$\delta^J_\pi\begin{pmatrix}i_1&\ldots&i_k\\ j_1&\ldots&j_l\end{pmatrix}=\delta^\tau_\pi\begin{pmatrix}i_1&\ldots&i_k\\ j_1&\ldots&j_l\end{pmatrix}w^J_\pi\begin{pmatrix}i_1&\ldots&i_k\\ j_1&\ldots&j_l\end{pmatrix}$$
where $\delta^\tau_\pi\in\{0,1\}$ is the generalized Kronecker function constructed in section 5 above, in the positive case, and where the sign on the right is given by
$$w^J_\pi\begin{pmatrix}i_1&\ldots&i_k\\ j_1&\ldots&j_l\end{pmatrix}=\prod_{r=1}^kw_{i_r}^{x_r-r}\prod_{s=1}^lw_{j_s}^{y_s-s}$$
where $k=(x_1\ldots x_k)$ and $l=(y_1\ldots y_l)$ are the writings of the colored integers $k,l$ by using $0,1$ symbols, with the conventions $\circ=0$, $\bullet=1$. 
\end{proposition}

\begin{proof}
In order to check our various categorical conditions, we just have to insert $w^J$ signs in the proof for positive super-structures, a bit in the same way as we did before, when talking twisting in chapter 13, for the signatures. This can be done as follows:

\medskip

(1) Regarding the tensor product, here we can use the following formula:
$$w^J_\pi\begin{pmatrix}i_1&\ldots&i_p\\j_1&\ldots&j_q\end{pmatrix}w^J_\sigma\begin{pmatrix}k_1&\ldots&k_r\\l_1&\ldots&l_s\end{pmatrix}
=w^J_{[\pi\sigma]}\begin{pmatrix}i_1&\ldots&i_p&k_1&\ldots&k_r\\j_1&\ldots&j_q&l_1&\ldots&l_s\end{pmatrix}$$

(2) For the composition, we can use the following formula:
$$w^J_\pi\begin{pmatrix}i_1&\ldots&i_p\\j_1&\ldots&j_q\end{pmatrix}
w^J_\sigma\begin{pmatrix}j_1&\ldots&j_q\\k_1&\ldots&k_r\end{pmatrix}
=w^J_{[^\pi_\sigma]}\begin{pmatrix}i_1&\ldots&i_p\\k_1&\ldots&k_r\end{pmatrix}$$

(3) For the involution axiom, we can use the following formula:
$$w^J_\pi\begin{pmatrix}i_1&\ldots&i_p\\ j_1&\ldots& j_q\end{pmatrix}
=w^J_{\pi^*}\begin{pmatrix}j_1&\ldots& j_q\\ i_1&\ldots&i_p\end{pmatrix}$$

(4) The identity axiom follows from the following computation:
$$w_{\hskip-2mm{\ }_{{\ }_\circ}\hskip-1.45mm|\hskip-3.85mm{\ }^{{\ }^\circ}}^J\binom{p}{p}=w_p^{0-1}w_p^{0-1}=w_p^2=1$$

(5) As for the duality axiom, this follows from a similar computation, namely:
$$w^J_{\!\!{\ }_\circ\hskip-1.1mm\cap\hskip-2.3mm{\ }_\bullet}(p,p)=w_p^{0-1}w_p^{1-2}=w_p^2=1$$

Thus, we have indeed a generalized Kronecker function, as claimed.
\end{proof}

At the quantum group level now, we first have:

\begin{proposition}
The quantum group associated to the category of pairings $P_2$, via the above generalized Kronecker function, is the quantum group $O_N^{J+}=Sp_N^+$.
\end{proposition}

\begin{proof}
In order to prove the result, we have to find inside our category the conditions which imply $u=J\bar{u}J^{-1}$. For the identity partition, we recall that we have: 
\begin{eqnarray*}
\delta_{\hskip-2mm{\ }_{{\ }_\circ}\hskip-1.45mm|\hskip-3.85mm{\ }^{{\ }^\bullet}}^\tau\binom{i}{j}=1
&\iff&\exists p,\binom{i}{j}=\binom{\tau^{\bullet(p)}}{\tau^\circ(p)}\\
&\iff&\exists p,\binom{i}{j}=\binom{\tau(p)}{p}\\
&\iff&i=\tau(j)
\end{eqnarray*}

In the case of negative structures, we have to add the following sign:
$$w^J_{\hskip-2mm{\ }_{{\ }_\circ}\hskip-1.45mm|\hskip-3.85mm{\ }^{{\ }^\bullet}}\binom{\tau(p)}{p}=w_{\tau(p)}^{1-1}w_p^{0-1}=1\cdot w_p=w_p$$

We can equally use the duality partition, where we have:
\begin{eqnarray*}
\delta^\tau_{\!\!{\ }_\circ\hskip-1.1mm\cap\hskip-2.3mm{\ }_\circ}(ij)=1
&\iff&\exists p,(i,j)=(\tau^\circ(p),\tau^{\bar{\circ}}(p))\\
&\iff&\exists p,(i,j)=(\tau^\circ(p),\tau^\bullet(p))\\
&\iff&j=\tau(i)
\end{eqnarray*}

In the case of negative structures, we have to add the following sign:
$$w^J_{\!\!{\ }_\circ\hskip-1.1mm\cap\hskip-2.3mm{\ }_\circ}(p,\tau(p))=w_p^{0-1}w_{\tau(p)}^{0-2}=w_p\cdot1=w_p$$

Summarizing, in both cases we have $u=J\bar{u}J^{-1}$, and this gives the result. 
\end{proof}

More generally, we have the following result, improving some previous findings:

\begin{theorem}
The basic quantum unitary and reflection groups associated to the generalized Kronecker map constructed above are as follows,
$$
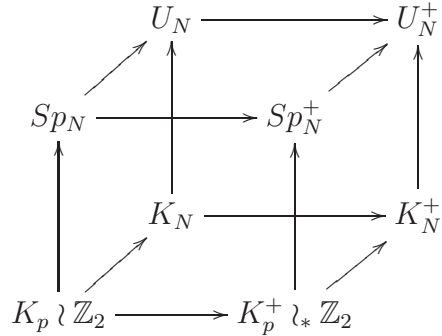
\xymatrix@R=6mm@C=3mm{
&U_N\ar[rr]&&U_N^+\\
Sp_N\ar[rr]\ar[ur]&&Sp_N^+\ar[ur]\\
&K_N\ar[uu]\ar[rr]&&K_N^+\ar[uu]\\
K_p\wr\mathbb Z_2\ar[uu]\ar[rr]\ar[ur]&&K_p^+\wr_*\mathbb Z_2\ar[uu]\ar[ur]}$$
with respect to the usual fundamental corepresentation.
\end{theorem}

\begin{proof}
We follow the proof of from the positive case, with the usual convention that $\bar{G}_K^\times$ denote the various quantum groups to be computed.

\medskip

(1) We know from the above that we have $\bar{U}_N^+=U_N^+$.

\medskip

(2) We also know from the above that we have $\bar{O}_N^+=Sp_N^+$.

\medskip

(3) Regarding now $\bar{U}_N$, here we must impose to the standard coordinates of $\bar{U}_N^+=U_N^+$ the relations coming from the basic crossing $\slash\hskip-2.1mm\backslash$\,, colored with its possible 4 matching colorings. But for any such matching coloring, the sign to be added is:
$$w^J_{\slash\hskip-1.6mm\backslash}\begin{pmatrix}i&j\\ j&i\end{pmatrix}=w_i^2\cdot w_j^2=1\cdot1=1$$

Thus we have $\bar{T}_\pi(e_i\otimes e_j)=e_j\otimes e_i$, and we obtain in this way usual commutation relations between the variables $\{u_{ij},u_{ij}^*\}$. Thus $\bar{U}_N$ is classical, and so $\bar{U}_N=U_N$.

\medskip

(4) The quantum group $\bar{O}_N$ appears as an intersection, as follows:
$$\bar{O}_N=\bar{U}_N\cap\bar{O}_N^+=U_N\cap Sp_N^+=Sp_N$$

(5) Regarding $\bar{K}_N$, we can follow again the method from the proof from the positive case. Indeed, in the main computation there, we have to add the following sign:
$$w^J_{\,|\hskip-0.6mm-\hskip-0.6mm|}\begin{pmatrix}i&i\\ i&i\end{pmatrix}=w_i^2=1$$

(6) The quantum group $\bar{H}_N$ appears as an intersection, as follows:
$$\bar{H}_N=\bar{K}_N\cap\bar{O}_N^+=K_N\cap Sp_N$$

In order to compute this intersection, we recall that for a $N\times N$ scalar matrix $U$, the relation $J\bar{U}=UJ$ reads:
$$U=\begin{pmatrix}a&b&\ldots\\
-\bar{b}&\bar{a}&\ldots\\
\ldots&\ldots&\ldots
\end{pmatrix}$$

Here $a,b,\ldots\in\mathbb C$. In the case $U\in K_N$, the $(^{\ a}_{-\bar{b}}{\ }^b_{\bar{a}})$ blocks can be either $(^0_0{\ }^0_0)$ or of the following form, with $z\in\mathbb T$:
$$\begin{pmatrix}z&0\\ 0&z\end{pmatrix}\quad
\begin{pmatrix}0&z\\ -\bar{z}&0\end{pmatrix}$$

By using once again the orthogonality condition, we must have exactly one nonzero block on each row and column, so we obtain $\bar{H}_N=K_p\wr\mathbb Z_2$, as claimed.

\medskip

(7) The proof of $\bar{K}_N^+=K_N^+$ follows by using the same argument as in the proof from the positive case, because in each case, the sign to be added is 1.

\medskip

(8) Finally, the quantum group $\bar{H}_N^+$ appears as an intersection, as follows:
$$\bar{H}_N^+=\bar{K}_N^+\cap\bar{O}_N^+=K_N^+\cap Sp_N^+$$

In order to compute this intersection, we use the interpretation of the commutation relation $uJ=J\bar{u}$ given in the proof of (6), with all the complex conjugates there replaced of course by adjoints. By reasoning as there, we conclude that we have a free wreath product decomposition $\bar{H}_N^+=K_p^+\wr_*\mathbb Z_2$, as claimed.
\end{proof}

\section*{14e. Exercises} 

We are now into difficult axiomatization questions in relation with super-easiness, and in relation with the above, we can only recommend working some more, and we have:

\begin{exercise}
Further build on the notion of super-easiness axiomatized in the above, notably by clarifying  the twists of the symplectic group $Sp_N$.
\end{exercise}

To be more precise here, in order for some general and satisfactory super-easiness theory to cover both the material from chapter 13, and from this chapter, we must solve first this exercise, and develop our improved theory afterwards. And with the remark that, even when doing so, things will not be over, because we still have some interesting quantum groups to be discussed, not covered by this. More on this in the next chapter.

\chapter{Quantum symmetries}

\section*{15a. Quantum symmetries}

We discuss here, following \cite{ba1}, \cite{wa2}, the quantum symmetry groups $S_Z^+$ of the finite quantum spaces $Z$, generalizing the quantum group $S_N^+$, coming from $Z=\{1,\ldots,N\}$, as well as the group $SO_3$, coming from $Z=M_2$. As with the quantum groups from the previous 2 chapters, these are waiting for a unification with the easy quantum groups.

\bigskip

In order to get started, we must talk about finite quantum spaces. In view of the general $C^*$-algebra theory explained in chapter 2, we have the following definition:

\index{quantum space}
\index{finite quantum space}
\index{counting measure}

\begin{definition}
A finite quantum space $Z$ is the abstract dual of a finite dimensional $C^*$-algebra $B$, according to the following formula:
$$C(Z)=B$$
The formal number of elements of such a space is $|Z|=\dim B$. By decomposing the algebra $B$, we have a formula of the following type:
$$C(Z)=M_{n_1}(\mathbb C)\oplus\ldots\oplus M_{n_k}(\mathbb C)$$
With $n_1=\ldots=n_k=1$ we obtain in this way the space $Z=\{1,\ldots,k\}$. Also, when $k=1$ the equation is $C(Z)=M_n(\mathbb C)$, and the solution will be denoted $Z=M_n$.
\end{definition}

In order to do some mathematics on such spaces, the very first observation is that we can talk about the formal number of points of such a space, as follows:
$$|Z|=\dim B$$

Alternatively, by decomposing the algebra $B$ as a sum of matrix algebras, as in Definition 15.1, we have the following formula for the formal number of points:
$$|Z|=n_1^2+\ldots+n_k^2$$

Pictorially, this suggests representing $Z$ as a set of $|Z|$ points in the plane, arranged in squares having sides $n_1,\ldots,n_k$, coming from the matrix blocks of $B$, as follows:
$$\begin{matrix}
\circ&\circ&\circ\\
\circ&\circ&\circ&&\ldots&&\circ&\circ\\
\circ&\circ&\circ&&&&\circ&\circ
\end{matrix}$$

As a second piece of mathematics, we can talk about counting measures, as follows:

\begin{definition}
Given a finite quantum space $Z$, we construct the functional
$$tr:C(Z)\to B(l^2(Z))\to\mathbb C$$
obtained by applying the regular representation, and the normalized matrix trace, and we call it integration with respect to the normalized counting measure on $Z$.
\end{definition}

To be more precise, consider the algebra $B=C(Z)$, which is by definition finite dimensional. We can make act $B$ on itself, by left multiplication:
$$\pi:B\to\mathcal L(B)\quad,\quad 
a\to(b\to ab)$$

The target of $\pi$ being a matrix algebra, $\mathcal L(B)\simeq M_N(\mathbb C)$ with $N=\dim B$, we can further compose with the normalized matrix trace, and we obtain $tr$:
$$tr=\frac{1}{N}\,Tr\circ\pi$$

As basic examples, for both $Z=\{1,\ldots,N\}$ and $Z=M_N$ we obtain the usual trace. In general, with $C(Z)=M_{n_1}(\mathbb C)\oplus\ldots\oplus M_{n_k}(\mathbb C)$, the weights of $tr$ are:
$$c_i=\frac{n_i^2}{\sum_in_i^2}$$

Pictorially, this suggests fine-tuning our previous picture of $Z$, by adding to each point the unnormalized trace of the corresponding element of $B$, as follows:
$$\begin{matrix}
\bullet_{n_1}&\circ_0&\circ_0\\
\circ_0&\bullet_{n_1}&\circ_0&&\ldots&&\bullet_{n_k}&\circ_0\\
\circ_0&\circ_0&\bullet_{n_1}&&&&\circ_0&\bullet_{n_k}
\end{matrix}$$

Here we have represented the points on the diagonals with solid circles, since they are of different nature from the off-diagonal ones, the attached numbers being nonzero. However, this picture is not complete either, and we can do better, as follows:

\begin{definition}
Given a finite quantum space $Z$, coming via a formula of type
$$C(Z)=M_{n_1}(\mathbb C)\oplus\ldots\oplus M_{n_k}(\mathbb C)$$
we use the following equivalent conventions for drawing $Z$:
\begin{enumerate}
\item Triple indices. We represent $Z$ as a set of $N=|Z|$ points, with each point being decorated with a triple index $ija$, coming from the standard basis $\{e_{ij}^a\}\subset B$. 

\item Double indices. As before, but by ignoring the index $a$, with the convention that $i,j$ belong to various indexing sets, one for each of the matrix blocks of $B$.

\item Single indices. As before, but with each point being now decorated with a single index, playing the role of the previous triple indices $ija$, or double indices $ij$.
\end{enumerate}
\end{definition}

All the above conventions are useful, and in practice, we will be mostly using the single index convention from (3). As an illustration, consider the space $Z=\{1,\ldots,k\}$. According to our single index convention, we can represent this space as a set of $k$ points, decorated by some indices, which must be chosen different. But the obvious choice for these $k$ different indices is $1,\ldots,k$, and we are led to the following picture:
$$\bullet_1\quad\bullet_2\quad\ldots\quad\bullet_k$$

As another illustration, consider the space $Z=M_n$. Here the picture is as follows, using double indices, which can be regarded as well as being single indices:
$$\begin{matrix}
\bullet_{11}&\circ_{12}&\circ_{13}\\
\circ_{21}&\bullet_{22}&\circ_{23}\\
\circ_{31}&\circ_{32}&\bullet_{33}
\end{matrix}$$

\smallskip

As yet another illustration, for the space $Z=M_3\sqcup M_2$, which appears by definition from the algebra $B=M_3(\mathbb C)\oplus M_2(\mathbb C)$, we are in need of triple indices, which can be of course regarded as single indices, in order to label all the points, and the picture is:
$$\begin{matrix}
\bullet_{111}&\circ_{121}&\circ_{131}\\
\circ_{211}&\bullet_{221}&\circ_{231}&&\ &&\bullet_{112}&\circ_{122}\\
\circ_{311}&\circ_{321}&\bullet_{331}&&&&\circ_{212}&\bullet_{222}
\end{matrix}$$

\smallskip

Let us study now the quantum group actions $G\curvearrowright Z$. If we denote by $\mu,\eta$ the multiplication and unit map of the algebra $C(Z)$, we have the following result:

\index{coaction}

\begin{proposition}
Consider a linear map $\Phi:C(Z)\to C(Z)\otimes C(G)$, written as
$$\Phi(e_i)=\sum_je_j\otimes u_{ji}$$
with $\{e_i\}$ being a linear space basis of $C(Z)$, chosen orthonormal with respect to $tr$.
\begin{enumerate}
\item $\Phi$ is a linear space coaction $\iff$ $u$ is a corepresentation.

\item $\Phi$ is multiplicative $\iff$ $\mu\in Hom(u^{\otimes 2},u)$.

\item $\Phi$ is unital $\iff$ $\eta\in Hom(1,u)$.

\item $\Phi$ leaves invariant $tr$ $\iff$ $\eta\in Hom(1,u^*)$.

\item If these conditions hold, $\Phi$ is involutive $\iff$ $u$ is unitary.
\end{enumerate}
\end{proposition}

\begin{proof}
This is similar to the proof for $S_N^+$ from chapter 2, as follows:

\medskip

(1) There are two axioms to be processed here, and we have indeed:
$$(id\otimes\Delta)\Phi=(\Phi\otimes id)\Phi
\iff\Delta(u_{ji})=\sum_ku_{jk}\otimes u_{ki}$$
$$(id\otimes\varepsilon)\Phi=id
\iff\varepsilon(u_{ji})=\delta_{ji}$$

(2) By using $\Phi(e_i)=u(e_i\otimes 1)$ we have the following identities, which give the result:
$$\Phi(e_ie_k)
=u(\mu\otimes id)(e_i\otimes e_k\otimes 1)$$
$$\Phi(e_i)\Phi(e_k)
=(\mu\otimes id)u^{\otimes 2}(e_i\otimes e_k\otimes 1)$$

(3) From $\Phi(e_i)=u(e_i\otimes1)$ we obtain by linearity, as desired:
$$\Phi(1)=u(1\otimes1)$$

(4) This follows from the following computation, by applying the involution:
\begin{eqnarray*}
(tr\otimes id)\Phi(e_i)=tr(e_i)1
&\iff&\sum_jtr(e_j)u_{ji}=tr(e_i)1\\
&\iff&\sum_ju_{ji}^*1_j=1_i\\
&\iff&(u^*1)_i=1_i\\
&\iff&u^*1=1
\end{eqnarray*}

(5) Assuming that (1-4) are satisfied, and that $\Phi$ is involutive, we have:
\begin{eqnarray*}
(u^*u)_{ik}
&=&\sum_lu_{li}^*u_{lk}\\
&=&\sum_{jl}tr(e_j^*e_l)u_{ji}^*u_{lk}\\
&=&(tr\otimes id)\sum_{jl}e_j^*e_l\otimes u_{ji}^*u_{lk}\\
&=&(tr\otimes id)(\Phi(e_i)^*\Phi(e_k))\\
&=&(tr\otimes id)\Phi(e_i^*e_k)\\
&=&tr(e_i^*e_k)1\\
&=&\delta_{ik}
\end{eqnarray*}

Thus $u^*u=1$, and since we know from (1) that $u$ is a corepresentation, it follows that $u$ is unitary. The proof of the converse is standard too, by using a similar computation.
\end{proof}

Following now \cite{ba1}, \cite{wa2}, we have the following result, extending the basic theory of $S_N^+$ from chapter 2 to the present finite quantum space setting:

\index{twisted quantum permutation}
\index{quantum automorphism group}
\index{finite quantum space}

\begin{theorem}
Given a finite quantum space $Z$, there is a universal compact quantum group $S_Z^+$ acting on $Z$, and leaving the counting measure invariant. We have
$$C(S_Z^+)=C(U_N^+)\Big/\Big<\mu\in Hom(u^{\otimes2},u),\eta\in Fix(u)\Big>$$
where $N=|Z|$, and where $\mu,\eta$ are the multiplication and unit maps of the algebra $C(Z)$. For the classical space $Z=\{1,\ldots,N\}$ we have $S_Z^+=S_N^+$. 
\end{theorem}

\begin{proof}
Here the first two assertions follow from Proposition 15.4, by using the standard fact that the complex conjugate of a corepresentation is a corepresentation too. As for the last assertion, regarding $S_N^+$, this follows from the results in chapter 2.
\end{proof}

The above result is quite conceptual, and we will see some applications in a moment. However, for many concrete questions, nothing beats multimatrix bases and indices. So, following the original paper of Wang \cite{wa2}, let us discuss this. We first have:

\begin{definition}
Given a finite quantum space $Z$, we let $\{e_i\}$ be the standard basis of $B=C(Z)$, so that the multiplication, involution and unit of $B$ are given by
$$e_ie_j=e_{ij}\quad,\quad 
e_i^*=e_{\bar{i}}\quad,\quad
1=\sum_{i=\bar{i}}e_i$$
where $(i,j)\to ij$ is the standard partially defined multiplication on the indices, with the convention $e_\emptyset=0$, and where $i\to\bar{i}$ is the standard involution on the indices.
\end{definition}

To be more precise, let $\{e_{ab}^r\}\subset B$ be the multimatrix basis. We set $i=(abr)$, and with this convention, the multiplication, coming from $e_{ab}^re_{cd}^p=\delta_{rp}\delta_{bc}e_{ad}^r$, is given by:
$$(abr)(cdp)=\begin{cases}
(adr)&{\rm if}\ b=c,\ r=p\\
\emptyset&{\rm otherwise}
\end{cases}$$

As for the involution, coming from $(e_{ab}^r)^*=e_{ba}^r$, this is given by:
$$\overline{(a,b,r)}=(b,a,r)$$

Finally, the unit formula comes from the following formula for the unit $1\in B$:
$$1=\sum_{ar}e_{aa}^r$$

Regarding now the generalized quantum permutation groups $S_Z^+$, the construction in Theorem 15.5 reformulates as follows, by using the above formalism:

\index{generalized magic}

\begin{proposition}
Given a finite quantum space $Z$, with basis $\{e_i\}\subset C(Z)$ as above, the algebra $C(S_Z^+)$ is generated by variables $u_{ij}$ with the following relations,
$$\sum_{ij=p}u_{ik}u_{jl}=u_{p,kl}\quad,\quad\sum_{kl=p}u_{ik}u_{jl}=u_{ij,p}$$
$$\sum_{i=\bar{i}}u_{ij}=\delta_{j\bar{j}}\quad,\quad\sum_{j=\bar{j}}u_{ij}=\delta_{i\bar{i}}$$
$$u_{ij}^*=u_{\bar{i}\hskip0.3mm\bar{j}}$$
with the fundamental corepresentation being the matrix $u=(u_{ij})$. We call a matrix $u=(u_{ij})$ satisfying the above relations ``generalized magic''.
\end{proposition}

\begin{proof}
We recall from Theorem 15.5 that the algebra $C(S_Z^+)$ appears as follows, where $N=|Z|$, and where $\mu,\eta$ are the multiplication and unit maps of $C(Z)$:
$$C(S_Z^+)=C(U_N^+)\Big/\Big<\mu\in Hom(u^{\otimes2},u),\eta\in Fix(u)\Big>$$

But the relations $\mu\in Hom(u^{\otimes2},u)$ and $\eta\in Fix(u)$ produce the 1st and 4th relations in the statement, then the biunitarity of $u$ gives the 5th relation, and finally the 2nd and 3rd relations follow from the 1st and 4th relations, by using the antipode.
\end{proof}

As an illustration, consider the case $Z=\{1,\ldots,N\}$. Here the index multiplication is $ii=i$ and $ij=\emptyset$ for $i\neq j$, and the involution is $\bar{i}=i$. Thus, our relations are as follows, corresponding to the standard magic conditions on a matrix $u=(u_{ij})$:
$$u_{ik}u_{il}=\delta_{kl}u_{ik}\quad,\quad u_{ik}u_{jk}=\delta_{ij}u_{ik}$$
$$\sum_iu_{ij}=1\quad,\quad\sum_ju_{ij}=1$$
$$u_{ij}^*=u_{ij}$$

As a second illustration now, which is something new, we have:

\begin{theorem}
For the space $Z=M_2$, coming via $C(Z)=M_2(\mathbb C)$, we have 
$$S_Z^+=SO_3$$
with the action $SO_3\curvearrowright M_2(\mathbb C)$ being the standard one, coming from $SU_2\to SO_3$.
\end{theorem}

\begin{proof}
This is something quite tricky, the idea being as follows:

\medskip

(1) First, we have an action by conjugation $SU_2\curvearrowright M_2(\mathbb C)$, and this action produces, via the canonical quotient map $SU_2\to SO_3$, an action as follows:
$$SO_3\curvearrowright M_2(\mathbb C)$$

(2) Then, it is routine to check, by using computations like those from the proof of $S_N^+=S_N$ at $N\leq3$, from chapter 2, that any action $G\curvearrowright M_2(\mathbb C)$ must come from a classical group. Thus the action $SO_3\curvearrowright M_2(\mathbb C)$ is universal, as claimed.

\medskip

(3) This was for the idea, and we will actually come back to this in a moment, in a more general setting, and with a new proof, complete this time. 
\end{proof}

As a conclusion so far, the quantum symmetry groups $S_Z^+$ unify the previous quantum permutation groups $S_N^+$, which are the most basic easy quantum groups, with the group $SO_3$. And this is of course very good news, in view our super-easiness purposes, because we have in this a way a potential method for including $SO_3$ in our easiness theory.

\section*{15b. Representation theory} 

Let us develop now some basic theory for the quantum symmetry groups $S_Z^+$, and their closed subgroups $G\subset S_Z^+$. We have here the following key result, from \cite{ba1}:

\index{quantum automorphism group}
\index{quantum symmetry group}
\index{Temperley-Lieb algebra}
\index{Marchenko-Pastur law}
\index{fusion rules}

\begin{theorem}
The quantum groups $S_Z^+$ have the following properties: 
\begin{enumerate}
\item The associated Tannakian categories are $TL_N$, with $N=|Z|$.

\item The main character follows the Marchenko-Pastur law $\pi_1$, when $|Z|\geq4$.

\item The fusion rules for $S_Z^+$ with $|Z|\geq4$ are the same as for $SO_3$.
\end{enumerate}
\end{theorem}

\begin{proof}
This result is from \cite{ba1}, the idea being as follows:

\medskip

(1) Let us pick our orthogonal basis $\{e_i\}$ as in Definition 15.6, so that we have, for a certain involution $i\to\bar{i}$ on the index set, the following formula:
$$e_i^*=e_{\bar{i}}$$

With this convention, we have the following computation:
\begin{eqnarray*}
\Phi(e_i)=\sum_je_j\otimes u_{ji}
&\implies&\Phi(e_i)^*=\sum_je_j^*\otimes u_{ji}^*\\
&\implies&\Phi(e_{\bar{i}})=\sum_je_{\bar{j}}\otimes u_{ji}^*\\
&\implies&\Phi(e_i)=\sum_je_j\otimes u_{\bar{i}\hskip0.3mm\bar{j}}^*
\end{eqnarray*}

Thus $u_{ji}^*=u_{\bar{i}\hskip0.3mm\bar{j}}$, so $u\sim\bar{u}$. Now with this result in hand, the proof goes as for the proof for $S_N^+$, from chapter 2. To be more precise, the result follows from the fact that the multiplication and unit of any complex algebra, and in particular of the algebra $C(Z)$ that we are interested in here, can be modelled by the following two diagrams:
$$m=|\cup|\qquad,\qquad u=\cap$$

Indeed, this is certainly true algebrically, and well-known, with as an illustration here, the associativity formula $m(m\otimes id)=(id\otimes m)m$ being checked as follows:
$$\begin{matrix}
|&\cup&|\ |&|\\
|&&\cup&|
\end{matrix}
\ \ =\ \ \begin{matrix}
|&|\ |&\cup&|\\
|&\cup&&|
\end{matrix}$$

As in what regards the $*$-structure, things here are fine too, because our choice for the trace from Definition 15.2 leads to the following formula regarding the adjoints, corresponding to $mm^*=N$, and so to the basic Temperley-Lieb calculus rule $\bigcirc=N$:
$$\mu\mu^*=N\cdot id$$

We conclude that the Tannakian category associated to $S_Z^+$ is, as claimed:
\begin{eqnarray*}
C
&=&<\mu,\eta>\\
&=&<m,u>\\
&=&<|\cup|,\cap>\\
&=&TL_N
\end{eqnarray*}

(2) The proof here is exactly as for $S_N^+$, by using moments. To be more precise, according to (1) these moments are the Catalan numbers, which are the moments of $\pi_1$.

\medskip

(3) Once again same proof as for $S_N^+$, by using the fact that the moments of $\chi$ are the Catalan numbers, which naturally leads to the Clebsch-Gordan rules.
\end{proof}

We can merge and reformulate our main results so far in the following way:

\index{projective version}

\begin{theorem}
The quantun groups $S_Z^+$ have the following properties: 
\begin{enumerate}
\item For $Z=\{1,\ldots,N\}$ we have $S_Z^+=S_N^+$.

\item For the space $Z=M_N$ we have $S_Z^+=PO_N^+=PU_N^+$.

\item In particular, for the space $Z=M_2$ we have $S_Z^+=SO_3$.

\item The fusion rules for $S_Z^+$ with $|Z|\geq4$ are independent of $Z$.

\item Thus, the fusion rules for $S_Z^+$ with $|Z|\geq4$ are the same as for $SO_3$.
\end{enumerate}
\end{theorem}

\begin{proof}
This is basically a compact form of what has been said above, with a new result added, and with some technicalities left aside, the idea being as follows:

\medskip

(1) This is something that we know from Theorem 15.5.

\medskip

(2) We recall from chapter 2 that we have $PO_N^+=PU_N^+$. Consider the standard vector space action of the free unitary group $U_N^+$, and its adjoint action:
$$U_N^+\curvearrowright\mathbb C^N\quad,\quad 
PU_N^+\curvearrowright M_N(\mathbb C)$$

By universality of $S_{M_N}^+$, we must have inclusions as follows:
$$PO_N^+\subset PU_N^+\subset S_{M_N}^+$$

On the other hand, the main character of $O_N^+$ with $N\geq2$ being semicircular, the main character of $PO_N^+$ must be Marchenko-Pastur. Thus the inclusion $PO_N^+\subset S_{M_N}^+$ has the property that it keeps fixed the law of main character, and by Peter-Weyl theory we conclude that this inclusion must be an isomorphism, as desired.

\medskip

(3) This is something that we know from Theorem 15.9, and that can be deduced as well from (2), by using the formula $PO_2^+=SO_3$, which is something elementary. Alternatively, this follows without computations from (4) below, because the inclusion of quantum groups $SO_3\subset S_{M_2}^+$ has the property that it preserves the fusion rules.

\medskip

(4) This is something that we know from Theorem 15.9.

\medskip

(5) This follows from (3,4), as already pointed out in Theorem 15.9.
\end{proof}

As an application of our extended formalism, the Cayley theorem for the finite quantum groups holds in the $S_Z^+$ setting. We have indeed the following result:

\index{Cayley theorem}
\index{no-Cayley theorem}
\index{Cayley embedding}

\begin{theorem}
Any finite quantum group $G$ has a Cayley embedding, as follows:
$$G\subset S_G^+$$
However, there are finite quantum groups which are not quantum permutation groups.
\end{theorem}

\begin{proof}
There are two statements here, the idea being as follows:

\medskip

(1) We have an action $G\curvearrowright G$, which leaves invariant the Haar measure. Now since the counting measure is left and right invariant, so is the Haar measure. We conclude that $G\curvearrowright G$ leaves invariant the counting measure, and so $G\subset S_G^+$, as claimed.

\medskip

(2) Regarding the second assertion, this is something non-trivial, the simplest counterexample being a certain quantum group $G$ appearing as a split abelian extension associated to the factorization $S_4=\mathbb Z_4S_3$, having cardinality $|G|=24$. 
\end{proof}

Finally, some interesting phenomena appear in the ``homogeneous'' case, where our quantum space is of the form $Z=M_K\times \{1\ldots,L\}$. Here we first have:

\index{wreath product}

\begin{proposition}
The classical symmetry group of $Z=M_K\times \{1\ldots,L\}$ is
$$S_Z=PU_K\wr S_L$$
with on the right a wreath product, equal by definition to $PU_K^L\rtimes S_L$.
\end{proposition}

\begin{proof}
The fact that we have an inclusion $PU_K\wr S_L\subset S_Z$ is standard, and this follows as well by taking the classical version of the inclusion $PU_K^+\wr_*S_L^+\subset S_Z^+$, established below. As for the fact that this inclusion $PU_K\wr S_L\subset S_Z$ is an isomorphism, this can be proved by picking an arbitrary element $g\in S_Z$, and decomposing it.
\end{proof}

Quite surprisingly, the quantum analogue of the above result fails:

\index{free wreath product}

\begin{theorem}
The quantum symmetry group of $Z=M_K\times \{1\ldots,L\}$ satisfies:
$$PU_K^+\wr_*S_L^+\subset S_Z^+$$
However, this inclusion is not an isomorphism at $K,L\geq2$.
\end{theorem}

\begin{proof}
We have several assertions to be proved, the idea being as follows:

\medskip

(1) The fact that we have $PU_K^+\wr_*S_L^+\subset S_Z^+$ is well-known and routine, by checking the fact that the matrix $w_{ija,klb}=u_{ij,kl}^{(a)}v_{ab}$ is a generalized magic unitary.

\medskip

(2) The inclusion $PU_K^+\wr_*S_L^+\subset S_Z^+$ is not an isomorphism, for instance by using \cite{twa}, along with the fact that $\pi_1\boxtimes\pi_1\neq\pi_1$ where $\pi_1$ is the Marchenko-Pastur law.
\end{proof}

\section*{15c. Twisting results} 

Now going towards $S_4^+$, let us start with the following definition, from \cite{bb2}:

\index{twisting}
\index{twisted determinant}
\index{anticommutation}

\begin{definition}
We let $SO_3'\subset O_3'$ be the subgroup coming from the relation
$$\sum_{\sigma\in S_3}u_{1\sigma(1)}u_{2\sigma(2)}u_{3\sigma(3)}=1$$
called twisted determinant one condition.
\end{definition}

Normally, we should prove here that $C(SO_3')$ is indeed a Woronowicz algebra. This is of course possible, by doing some computations, but we will not need to do these computations, because this follows from the following result, from \cite{bb2}:

\index{Klein group}

\begin{theorem}
We have an isomorphism of compact quantum groups
$$S_4^+=SO_3'$$
given by the Fourier transform over the Klein group $K=\mathbb Z_2\times\mathbb Z_2$.
\end{theorem}

\begin{proof}
Consider the following matrix, coming from the action of $SO_3'$ on $\mathbb C^4$:
$$u^+=\begin{pmatrix}1&0\\0&u\end{pmatrix}$$

We apply to this matrix the Fourier transform over the Klein group $K=\mathbb Z_2\times\mathbb Z_2$: 
$$v=
\frac{1}{4}
\begin{pmatrix}
1&1&1&1\\
1&-1&-1&1\\
1&-1&1&-1\\
1&1&-1&-1
\end{pmatrix}
\begin{pmatrix}
1&0&0&0\\
0&u_{11}&u_{12}&u_{13}\\
0&u_{21}&u_{22}&u_{23}\\
0&u_{31}&u_{32}&u_{33}
\end{pmatrix}
\begin{pmatrix}
1&1&1&1\\
1&-1&-1&1\\
1&-1&1&-1\\
1&1&-1&-1
\end{pmatrix}$$

This matrix is then magic, and vice versa, so the Fourier transform over $K$ converts the relations in Definition 15.14 into the magic relations. But this gives the result.
\end{proof}

There are many more things that can be said here, and we have:

\begin{theorem}
The quantum group $S_4^+=SO_3'$ has the following properties:
\begin{enumerate}
\item It appears as a cocycle twist of $SO_3$.

\item Its fusion rules are the same as for $SO_3$.

\item Its subgroups are basically twists of the subgroups of $SO_3$.
\end{enumerate}
\end{theorem}

\begin{proof}
These are more advanced results, from \cite{bb2}, the idea being as follows:

\medskip

(1) This follows by suitably reformulating the definition of $SO_3'$ given above in purely algebraic terms, using cocycles, and for details here, we refer to \cite{bb2}. In what concerns us, we will actually discuss a generalization of this, right next, following \cite{bbs}.

\medskip

(2) This is something that we know well, via numerous proofs, and we can add to our trophy list one more proof, coming from (1), via standard cocycle twisting theory.

\medskip

(3) The idea here is that the subgroups $G\subset SO_3$ are subject to an ADE classification result, and the subgroups of $SO_3'$ are basically twists of these, $G'\subset SO_3'$. 
\end{proof}

An interesting extension of the $S_4^+=SO_3'$ result comes by looking at the general case $N=n^2$, with $n\in\mathbb N$. We will prove that we have a twisting result, as follows: 
$$PO_n^+=(S_N^+)^\sigma$$

In order to explain this material, from \cite{bbs}, which is quite technical, requiring good algebraic knowledge, let us begin with some generalities. We first have:

\begin{proposition}
Given a finite group $G$, the algebra $C(S_{\widehat{G}}^+)$ is isomorphic to the abstract algebra presented by generators $x_{gh}$ with $g,h\in G$, with the following relations:
$$x_{1g}=x_{g1}=\delta_{1g}\quad,\quad
x_{s,gh}=\sum_{t\in G}x_{st^{-1},g}x_{th}\quad,\quad 
x_{gh,s}=\sum_{t\in G}x_{gt^{-1}}x_{h,ts}$$
The comultiplication, counit and antipode are given by the formulae
$$\Delta(x_{gh})=\sum_{s\in G}x_{gs}\otimes x_{sh}\quad,\quad 
\varepsilon(x_{gh})=\delta_{gh}\quad,\quad
S(x_{gh})=x_{h^{-1}g^{-1}}$$
on the standard generators $x_{gh}$.
\end{proposition}

\begin{proof}
This follows indeed from a direct verification, based either on Theorem 15.5, or on its equivalent formulation from Proposition 15.7. See \cite{bbs}.
\end{proof}

Let us discuss now the twisted version of the above result. Consider a 2-cocycle on $G$, which is by definition a map $\sigma:G\times G\to\mathbb C^*$ satisfying:
$$\sigma_{gh,s}\sigma_{gh}=\sigma_{g,hs}\sigma_{hs}\quad,\quad 
\sigma_{g1}=\sigma_{1g}=1$$

Given such a cocycle, we can construct the associated twisted group algebra $C(\widehat{G}_\sigma)$, as being the vector space $C(\widehat{G})=C^*(G)$, with product $e_ge_h=\sigma_{gh}e_{gh}$. We have:

\begin{proposition}
The algebra $C(S_{\widehat{G}_\sigma}^+)$ is isomorphic to the abstract algebra presented by generators $x_{gh}$ with $g,h\in G$, with the relations $x_{1g}=x_{g1}=\delta_{1g}$ and:
$$\sigma_{gh}x_{s,gh}=\sum_{t\in G}\sigma_{st^{-1},t}x_{st^{-1},g}x_{th}\quad,\quad
\sigma_{gh}^{-1}x_{gh,s}=\sum_{t\in G}\sigma_{t^{-1},ts}^{-1}x_{gt^{-1}}x_{h,ts}$$
The comultiplication, counit and antipode are given by the formulae
$$\Delta(x_{gh})=\sum_{s\in G}x_{gs}\otimes x_{sh}\quad,\quad 
\varepsilon(x_{gh})=\delta_{gh}\quad,\quad
S(x_{gh})=\sigma_{h^{-1}h}\sigma_{g^{-1}g}^{-1}x_{h^{-1}g^{-1}}$$
on the standard generators $x_{gh}$.
\end{proposition}

\begin{proof}
Once again, this follows from a direct verification, explained in \cite{bbs}.
\end{proof}

\index{twisting}
\index{cocycle twisting}

In what follows, we will prove that the quantum groups $S_{\widehat{G}}^+$ and $S_{\widehat{G}_\sigma}^+$ are related by a cocycle twisting operation. Let $A$ be a Hopf algebra. We recall that a left 2-cocycle is a convolution invertible linear map
$\sigma:A\otimes A\to\mathbb C$ satisfying:
$$\sigma_{x_1y_1}\sigma_{x_2y_2,z}=\sigma_{y_1z_1}\sigma_{x,y_2z_2}\quad,\quad 
\sigma_{x1}=\sigma_{1x}=\varepsilon(x)$$

Note that $\sigma$ is a left 2-cocycle if and only if $\sigma^{-1}$, the convolution inverse of $\sigma$, is a right 2-cocycle, in the sense that we have:
$$\sigma^{-1}_{x_1y_1,z}\sigma^{-1}_{x_1y_2}=\sigma^{-1}_{x,y_1z_1}\sigma^{-1}_{y_2z_2}\quad,\quad 
\sigma^{-1}_{x1}=\sigma^{-1}_{1x}=\varepsilon(x)$$

Given a left 2-cocycle $\sigma$ on $A$, one can form the 2-cocycle twist $A^\sigma$ as follows. As a coalgebra, $A^\sigma=A$, and an element $x\in A$, when considered in $A^\sigma$, is denoted $[x]$. The product in $A^\sigma$ is then defined, in Sweedler notation, by: 
$$[x][y]=\sum\sigma_{x_1y_1}\sigma^{-1}_{x_3y_3}[x_2y_2]$$

We can now state and prove a main theorem from \cite{bbs}, as follows:

\begin{theorem}
If $G$ is a finite group and $\sigma$ is a $2$-cocycle on $G$, the Hopf algebras
$$C(S_{\widehat{G}}^+)\quad,\quad C(S_{\widehat{G}_\sigma}^+)$$
are $2$-cocycle twists of each other, in the above sense.
\end{theorem}

\begin{proof}
In order to prove this result, we use the following Hopf algebra map: 
$$\pi:C(S_{\widehat{G}}^+)\to C(\widehat{G})\quad,\quad
x_{gh}\to\delta_{gh}e_g$$

Our 2-cocycle $\sigma:G\times G\to\mathbb C^*$ can be extended by linearity into a linear map as follows, which is a left and right 2-cocycle in the above sense:
$$\sigma:C(\widehat{G})\otimes C(\widehat{G})\to\mathbb C$$

Consider now the following composition:
$$\alpha=\sigma(\pi \otimes \pi):C(S_{\widehat{G}}^+)\otimes C(S_{\widehat{G}}^+)\to C(\widehat{G})\otimes C(\widehat{G})\to\mathbb C$$

Then $\alpha$ is a left and right 2-cocycle, because it is induced by a cocycle on a group algebra, and so is its convolution inverse $\alpha^{-1}$. Thus we can construct the twisted algebra $C(S_{\widehat{G}}^+)^{\alpha^{-1}}$, and inside this algebra we have the following computation:
$$[x_{gh}][x_{rs}]
=\alpha^{-1}(x_g,x_r)\alpha(x_h,x_s)[x_{gh}x_{rs}]
=\sigma_{gr}^{-1}\sigma_{hs}[x_{gh}x_{rs}]$$

By using this, we obtain the following formula:
$$\sum_{t\in G}\sigma_{st^{-1},t}[x_{st^{-1},g}][x_{th}]
=\sum_{t\in G}\sigma_{st^{-1},t}\sigma_{st^{-1},t}^{-1}\sigma_{gh}[x_{st^{-1},g}x_{th}]
=\sigma_{gh}[x_{s,gh}]$$

Similarly, we have the following formula:
$$\sum_{t\in G}\sigma_{t^{-1},ts}^{-1}[x_{g,t^{-1}}][x_{h,ts}]=\sigma_{gh}^{-1}[x_{gh,s}]$$

We deduce from this that there exists a Hopf algebra map, as follows:
$$\Phi:C(S_{\widehat{G}_\sigma}^+)\to C(S_{\widehat{G}}^+)^{\alpha^{-1}}\quad,\quad 
x_{gh}\to [x_{g,h}]$$

This map is clearly surjective, and is injective as well, by a standard fusion semiring argument, because both Hopf algebras have the same fusion semiring.
\end{proof}

Summarizing, we have proved our main twisting result. Our purpose in what follows will be that of working out versions and particular cases of it. We first have: 

\begin{proposition}
If $G$ is a finite group and $\sigma$ is a $2$-cocycle on $G$, then
$$\Phi(x_{g_1h_1}\ldots x_{g_mh_m})=\Omega(g_1,\ldots,g_m)^{-1}\Omega(h_1,\ldots,h_m)x_{g_1h_1}\ldots x_{g_mh_m}$$
with the coefficients on the right being given by the formula
$$\Omega(g_1,\ldots,g_m)=\prod_{k=1}^{m-1}\sigma_{g_1\ldots g_k,g_{k+1}}$$
is a coalgebra isomorphism $C(S_{\widehat{G}_\sigma}^+)\to C(S_{\widehat{G}}^+)$, commuting with the Haar integrals.
\end{proposition}

\begin{proof}
This is indeed just a technical reformulation of Theorem 15.19.
\end{proof}

Let us discuss now some concrete applications of the general results established above. Consider the group $G=\mathbb Z_n^2$, let $w=e^{2\pi i/n}$, and consider the following cocycle: 
$$\sigma:G\times G\to\mathbb C^*\quad,\quad 
\sigma_{(ij)(kl)}=w^{jk}$$ 

In order to understand what is the formula that we obtain, we must do some computations. Let $E_{ij}$ with $i,j \in\mathbb Z_n$ be the standard basis of $M_n(\mathbb C)$. We first have:

\begin{proposition}
The linear map given by
$$\psi(e_{(i,j)})=\sum_{k=0}^{n-1}{w}^{ki}E_{k,k+j}$$
defines an isomorphism of algebras $\psi:C(\widehat{G}_\sigma)\simeq M_n(\mathbb C)$. 
\end{proposition}

\begin{proof}
Consider indeed the following linear map:
$$\psi'(E_{ij})=\frac{1}{n}\sum_{k=0}^{n-1}{w}^{-ik}e_{(k,j-i)}$$
 
It is routine to check that both $\psi,\psi'$ are morphisms of algebras, and that these maps are inverse to each other. In particular, $\psi$ is an isomorphism of algebras, as stated.
\end{proof}

Next in line, we have the following result:

\begin{proposition}
The algebra map given by
$$\varphi(u_{ij}u_{kl}) = \frac{1}{n}\sum_{a,b=0}^{n-1}{w}^{ai-bj}x_{(a,k-i),(b,l-j)}$$
defines a Hopf algebra isomorphism $\varphi:C(S_{M_n}^+)\simeq C(S_{\widehat{G}_\sigma}^+)$.
\end{proposition}

\begin{proof}
Consider the universal coactions on the two algebras in the statement:
\begin{eqnarray*}
\alpha:M_n(\mathbb C)&\to&M_n({\mathbb C})\otimes C(S_{M_n}^+)\\
\beta:C(\widehat{G}_\sigma)&\to&C(\widehat{G}_\sigma)\otimes C(S_{\widehat{G}_\sigma}^+)
 \end{eqnarray*}
 
In terms of the standard bases, these coactions are given by:
\begin{eqnarray*}
\alpha(E_{ij})&=&\sum_{kl}E_{kl}\otimes u_{ki}u_{lj}\\
\beta(e_{(i,j)})&=&\sum_{kl} e_{(k,l)}\otimes x_{(k,l),(i,j)}
\end{eqnarray*}

We use now the identification $C(\widehat{G}_\sigma)\simeq M_n(\mathbb C)$ from Proposition 15.21. This identification produces a coaction map, as follows:
$$\gamma:M_n(\mathbb C)\to M_n(\mathbb C)\otimes C(S_{\widehat{G}_\sigma}^+)$$

Now observe that this map is given by the following formula:
 $$\gamma(E_{ij})=\frac{1}{n}\sum_{ab}E_{ab}\otimes\sum_{kr}w^{ar-ik} x_{(r,b-a),(k,j-i)}$$

By comparing with the formula of $\alpha$, we obtain the isomorphism in the statement.
\end{proof}

We will need one more result of this type, as follows:

\begin{proposition}
The algebra map given by
$$\rho(x_{(a,b),(i,j)})=\frac{1}{n^2}\sum_{klrs}w^{ki+lj-ra-sb}p_{(r,s),(k,l)}$$
defines a Hopf algebra isomorphism $\rho:C(S_{\widehat{G}}^+)\simeq C(S_G^+)$.
\end{proposition}

\begin{proof}
We have a Fourier transform isomorphism, as follows:
$$C(\widehat{G})\simeq C(G)$$

Thus the algebras in the statement are indeed isomorphic.
\end{proof}

As a conclusion to all this, we have the following result, from \cite{bbs}:

\begin{theorem}
Let $n\geq 2$ and $w=e^{2\pi i/n}$. Then
$$\Theta(u_{ij}u_{kl})=\frac{1}{n}\sum_{ab=0}^{n-1}w^{-a(k-i)+b(l-j)}p_{ia,jb}$$
defines a coalgebra isomorphism $C(PO_n^+)\to C(S_{n^2}^+)$ commuting with the Haar integrals.
\end{theorem}
 
\begin{proof}
We know from Theorem 15.10 that we have identifications as follows, where the projective version of $(A,u)$ is the pair $(PA,v)$, with $PA=<v_{ij}>$ and $v=u\otimes\bar{u}$:
$$PO_n^+=PU_n^+=S_{M_n}^+$$

With this in hand, the result follows from Theorem 15.19 and Proposition 15.21, by combining them with the various isomorphisms established above.
\end{proof}

\section*{15d. Quantum reflections} 

Let us start with the following straightforward extension of the usual notion of finite graph, from \cite{fpi}, obtained by using a finite quantum space as set of vertices:

\index{quantum graph}
\index{finite quantum graph}

\begin{definition}
We call ``finite quantum graph'' a pair of type
$$X=(Z,d)$$
with $Z$ being a finite quantum space, and with $d\in M_N(\mathbb C)$ being a matrix.
\end{definition}

This is of course something quite general. In the case $Z=\{1,\ldots,N\}$ for instance, what we have here is a directed graph, with the edges $i\to j$ colored by complex numbers $d_{ij}\in\mathbb C$, and with self-edges $i\to i$ allowed too, again colored by numbers $d_{ii}\in\mathbb C$. In the general case, however, where $Z$ is arbitrary, the need for extra conditions of type $d=d^*$, or $d_{ii}=0$, or $d\in M_N(\mathbb R)$, or $d\in M_N(0,1)$ and so on, is not very natural, as we will soon discover, and it is best to use Definition 15.25 as such, with no restrictions on $d$.

\bigskip

In general, a quantum graph can be represented as a colored oriented graph on $\{1,\ldots,N\}$, where $N=|Z|$, with the vertices being decorated by single indices $i$, and with the colors being complex numbers, namely the entries of $d$. This is similar to the formalism from before, but there is a discussion here in what regards the exact choice of the colors, which are usually irrelevant in connection with our symmetry problematics, and so can be true colors instead of complex numbers. More on this later.

\bigskip

With the above notion in hand, we have the following definition, also from \cite{fpi}:

\begin{definition}
The quantum automorphism group of $X=(Z,d)$ is the subgroup 
$$G^+(X)\subset S_Z^+$$
obtained via the relation $du=ud$, where $u=(u_{ij})$ is the fundamental corepresentation.
\end{definition}

Again, this is something very natural, coming as a continuation of the constructions for usual graphs. We refer to \cite{fpi}, \cite{twa} for more on this notion, and for a number of advanced computations, in relation with free wreath products. At an elementary level, a first problem is that of working out the basics of the correspondence $X\to G^+(X)$. There are several things to be done here, namely simplices, complementation, color independence, multi-simplices, and with a few twists, all this basically extends well. 

\bigskip

Let us start with the simplices. As we will soon discover, things are quite tricky here, leading us in particular to the conclusion that the simplex based on an arbitrary finite quantum space $Z$ is not a usual graph, with $d\in M_N(0,1)$ where $N=|Z|$, but rather a sort of ``signed graph'', with $d\in M_N(-1,0,1)$. Let us start our study with:

\index{empty graph}
\index{simplex}

\begin{theorem}
Given a finite quantum space $Z$, we have
$$G^+(Z_{empty})=G^+(Z_{full})=S_Z^+$$
where $Z_{empty}$ is the empty graph on the vertex set $Z$, coming from the matrix $d=0$, and where $Z_{full}$ is the simplex on the vertex set $Z$, coming from the matrix 
$$d=NP_1-1_N$$
where $N=|Z|$, and where $P_1$ is the orthogonal projection on the unit $1\in C(Z)$.
\end{theorem}

\begin{proof}
This is something quite tricky, the idea being as follows:

\medskip

(1) First of all, the formula $G^+(Z_{empty})=S_Z^+$ is clear from definitions, because the commutation of $u$ with the matrix $d=0$ is automatic.

\medskip

(2) Regarding $G^+(Z_{full})=S_Z^+$, let us first discuss the classical case, $Z=\{1,\ldots,N\}$. Here the simplex $Z_{full}$ is the graph having having edges between any two vertices, whose adjacency matrix is $d=\mathbb I_N-1_N$, where $\mathbb I_N$ is the all-1 matrix. The commutation of $u$ with $1_N$ being automatic, and the commutation with $\mathbb I_N$ being automatic too, $u$ being bistochastic, we have $[u,d]=0$, and so $G^+(Z_{full})=S_Z^+$ in this case, as stated.

\medskip

(3) In the general case, we know from Theorem 15.5 that we have $\eta\in Fix(u)$, with $\eta:\mathbb C\to C(Z)$ being the unit map. Thus we have $P_1\in End(u)$, and so the condition $[u,P_1]=0$ is automatic. Together with the fact that in the classical case we have $\mathbb I_N=NP_1$, this suggests to define the adjacency matrix of the simplex as being $d=NP_1-1_N$, and with this definition, we have indeed $G^+(Z_{full})=S_Z^+$, as claimed.
\end{proof}

Let us study now the simplices $Z_{full}$ found in Theorem 15.27. In the classical case, $Z=\{1,\ldots,N\}$, what we have is of course the usual simplex. However, in the general case things are more mysterious, the first result here being as follows:

\begin{proposition}
The adjacency matrix of the simplex $Z_{full}$, given by definition by $d=NP_1-1_N$, is a matrix $d\in M_N(-1,0,1)$, which can be computed as follows:
\begin{enumerate}
\item In single index notation, $d_{ij}=\delta_{i\bar{i}}\delta_{j\bar{j}}-\delta_{ij}$.

\item In double index notation, $d_{ab,cd}=\delta_{ab}\delta_{cd}-\delta_{ac}\delta_{bd}$.

\item In triple index notation, $d_{abp,cdq}=\delta_{ab}\delta_{cd}-\delta_{ac}\delta_{bd}\delta_{pq}$.
\end{enumerate}
\end{proposition}

\begin{proof}
According to our single index conventions, from Definition 15.3, the adjacency matrix of the simplex is the one in the statement, namely:
\begin{eqnarray*}
d_{ij}
&=&(NP_1-1_N)_{ij}\\
&=&\bar{1}_i1_j-\delta_{ij}\\
&=&\delta_{i\bar{i}}\delta_{j\bar{j}}-\delta_{ij}
\end{eqnarray*}

In double index notation now, with $i=(ab)$ and $j=(cd)$, and $a,b,c,d$ being usual matrix indices, each thought to be attached to the corresponding matrix block of $C(Z)$, the formula that we obtain in the second one in the statement, namely:
\begin{eqnarray*}
d_{ab,cd}
&=&\delta_{ab,ba}\delta_{cd,dc}-\delta_{ab,cd}\\
&=&\delta_{ab}\delta_{cd}-\delta_{ac}\delta_{bd}
\end{eqnarray*}

Finally, in standard triple index notation, $i=(abp)$ and $j=(cdq)$, with $a,b,c,d$ being now usual numeric matrix indices, ranging in $1,2,3,\ldots\,$, and with $p,q$ standing for corresponding blocks of the algebra $C(Z)$, the formula that we obtain is:
\begin{eqnarray*}
d_{abp,cdq}
&=&\delta_{abp,bap}\delta_{cdq,dcq}-\delta_{abp,cdq}\\
&=&\delta_{ab}\delta_{cd}-\delta_{ac}\delta_{bd}\delta_{pq}
\end{eqnarray*}

Thus, we are led to the conclusions in the statement.
\end{proof}

At the level of examples, for $Z=\{1,\ldots,N\}$ the best is to use the above formula (1). The involution on the index set is $\bar{i}=i$, and we obtain, as we should:
$$d_{ij}=1-\delta_{ij}$$

As a more interesting example now, for the quantum space $Z=M_n$, coming by definition via the formula $C(Z)=M_n(\mathbb C)$, the situation is as follows:

\begin{proposition}
The simplex $Z_{full}$ with $Z=M_n$ is as follows: 
\begin{enumerate}
\item The vertices are $n^2$ points in the plane, arranged in square form.

\item Usual edges, worth $1$, are drawn between distinct points on the diagonal.

\item In addition, each off-diagonal point comes with a self-edge, worth $-1$.
\end{enumerate}
\end{proposition}

\begin{proof}
Here the most convenient is to use the double index formula from Proposition 15.28 (2), which tells us that $d$ is as follows, with indices $a,b,c,d\in\{1,\ldots,n\}$:
$$d_{ab,cd}=\delta_{ab}\delta_{cd}-\delta_{ac}\delta_{bd}$$

This quantity can be $-1,0,1$, and the study goes as follows:

\medskip

-- Case $d_{ab,cd}=1$. This can only happen when $\delta_{ab}\delta_{cd}=1$ and $\delta_{ac}\delta_{bd}=0$, corresponding to a formula of type $d_{aa,cc}=0$, with $a\neq c$, and so to the edges in (2).

\medskip

-- Case $d_{ab,cd}=-1$. This can only happen when $\delta_{ab}\delta_{cd}=0$ and $\delta_{ac}\delta_{bd}=1$, corresponding to a formula of type $d_{ab,ab}=0$, with $a\neq b$, and so to the self-edges in (3).
\end{proof}

The above result is quite interesting, and as an illustration, here is the pictorial representation of the simplex $Z_{full}$ on the vertex set $Z=M_3$, with the convention that the solid arrows are worth $-1$, and the dashed arrows are worth 1:
$$\xymatrix@R=15pt@C=15pt
{\bullet\ar@{--}[dr]\ar@/^/@{--}[ddrr]&\circ^\circlearrowright&\circ^\circlearrowright\\
\circ^\circlearrowright&\bullet\ar@{--}[dr]&\circ^\circlearrowright\\
\circ^\circlearrowright&\circ^\circlearrowright&\bullet}$$

More generally, we can in fact compute $Z_{full}$ for any finite quantum space $Z$, with the result here, which will be our final saying on the subject, being as follows:

\begin{theorem}
Consider a finite quantum space $Z$, and write it as follows, according to the decomposition formula $C(Z)=M_{n_1}(\mathbb C)\oplus\ldots\oplus M_{n_k}(\mathbb C)$ for its function algebra:
$$Z=M_{n_1}\sqcup\ldots\sqcup M_{n_k}$$
The simplex $Z_{full}$ is then the classical simplex formed by the points lying on the diagonals of $M_{n_1},\ldots,M_{n_k}$, with self-edges added, each worth $-1$, at the non-diagonal points.
\end{theorem}

\begin{proof}
The study here is quite similar to the one from the proof of Proposition 15.29, but by using this time the triple index formula from Proposition 15.28 (3), namely:
$$d_{abp,cdq}=\delta_{ab}\delta_{cd}-\delta_{ac}\delta_{bd}\delta_{pq}$$

Indeed, this quantity can be $-1,0,1$, and the $1$ case appears precisely as follows, leading to the classical simplex mentioned in the statement:
$$d_{aap,ccq}=1\quad,\quad\forall ap\neq cq$$

As for the remaining $-1$ case, this appears precisely as follows, leading this time to the self-edges worth $-1$, also mentioned in the statement:
$$d_{abp,abp}=1\quad,\quad\forall a\neq b$$

Thus, we are led to the conclusion in the statement.
\end{proof}

As an illustration, here is the simplex on the vertex set $Z=M_3\sqcup M_2$, with again the convention that the solid arrows are worth $-1$, and the dashed arrows are worth 1:
$$
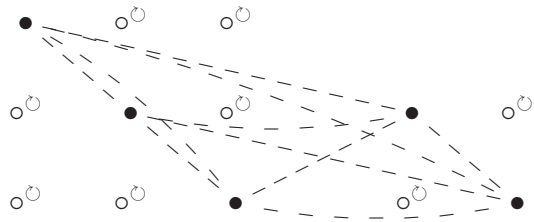
\xymatrix@R=18pt@C=20pt
{\bullet\ar@{--}[dr]\ar@/^/@{--}[ddrr]
\ar@{--}[drrrr]\ar@/^/@{--}[ddrrrrr]
&\circ^\circlearrowright&\circ^\circlearrowright\\
\circ^\circlearrowright&\bullet\ar@{--}[dr]\ar@/_/@{--}[rrr]\ar@{--}[drrrr]
&\circ^\circlearrowright&&\bullet\ar@{--}[dr]&\circ^\circlearrowright\\
\circ^\circlearrowright&\circ^\circlearrowright&\bullet\ar@{--}[urr]\ar@/_/@{--}[rrr]
&&\circ^\circlearrowright&\bullet}$$

\medskip

Long story short, we know what the simplex $Z_{full}$ is, and we have the formula $G^+(Z_{empty})=G^+(Z_{full})=S_Z^+$, exactly as in the $Z=\{1,\ldots,N\}$ case. Now with the above results in hand, we can talk as well about complementation, as follows:

\begin{theorem}
For any finite quantum graph $X$ we have the formula
$$G^+(X)=G^+(X^c)$$
where $X\to X^c$ is the complementation operation, given by $d_X+d_{X^c}=d_{Z_{full}}$.
\end{theorem}

\begin{proof}
This follows from Theorem 15.27, and more specifically from the following commutation relation, which is automatic, as explained there:
$$[u,d_{Z_{full}}]=0$$

Let us mention too that, in what concerns the pictorial representation of $X^c$, this can be deduced from what we have Theorem 15.30, in the obvious way.
\end{proof}

With this technology in hand, we can talk about twisted quantum reflections. The idea here, from \cite{ba2}, will be that the twisted analogues of the quantum reflection groups $H_N^{s+}\subset S_{sN}^+$ will be the quantum automorphism groups $S_{Z\to Y}^+$ of the fibrations of finite quantum spaces $Z\to Y$, which correspond by definition to the Markov inclusions of finite dimensional $C^*$-algebras $C(Y)\subset C(Z)$. In order to discuss this, let us start with:

\index{Markov inclusion}
\index{Markov fibration}
\index{fibration}
\index{quantum space fibration}

\begin{definition}
A fibration of finite quantum spaces $Z\to Y$ corresponds to an inclusion of finite dimensional $C^*$-algebras
$$C(Y)\subset C(Z)$$
which is Markov, in the sense that it commutes with the canonical traces.
\end{definition}

Here the commutation condition with the canonical traces means that the composition $C(Y)\subset C(Z)\to\mathbb C$ should equal the canonical trace $C(Y)\to\mathbb C$. At the level of the corresponding quantum spaces, this means that the quotient map $Z\to Y$ must commute with the corresponding counting measures, and this is where our term ``fibration'' comes from. In order to talk now about the quantum symmetry groups $S_{Z\to Y}^+$, we will need:

\index{Jones projection}

\begin{proposition}
Given a fibration $Z\to Y$, a closed subgroup $G\subset S_Z^+$ leaves invariant $Y$ precisely when its magic unitary $u=(u_{ij})$ satisfies the condition
$$e\in End(u)$$
where $e:C(Z)\to C(Z)$ is the Jones projection, onto the subalgebra $C(Y)\subset C(Z)$.
\end{proposition}

\begin{proof}
This is something that we know well, in the commutative case, where $Z$ is a usual finite set, and the proof in general is similar.
\end{proof}

We can now talk about twisted quantum reflection groups, as follows:

\index{twisted reflection group}
\index{twisted quantum reflection group}

\begin{theorem}
Any fibration of finite quantum spaces $Z\to Y$ has a quantum symmetry group, which is the biggest acting on $Z$ by leaving $Y$ invariant:
$$S_{Z\to Y}^+\subset S_Z^+$$
At the level of algebras of functions, this quantum group $S_{Z\to Y}^+$ is obtained as follows, with $e:C(Z)\to C(Y)$ being the Jones projection:
$$C(S_{Z\to Y}^+)=C(S_Z^+)\Big/\Big<e\in End(u)\Big>$$
We call these quantum groups $S_{Z\to Y}^+$ twisted quantum reflection groups. 
\end{theorem}

\begin{proof}
This follows indeed from Proposition 15.33.
\end{proof}

As a basic example, let us discuss the commutative case. Here we have:

\begin{proposition}
In the commutative case, the fibration $Z\to Y$ must be of the following special form, with $N,s$ being certain integers,
$$\{1,\ldots,N\}\times\{1,\ldots,s\}\to\{1,\ldots,N\}\quad,\quad (i,a)\to i$$
and we obtain the quantum reflection groups studied in chapter 6,
$$(S_{Z\to Y}^+\subset S_Z^+)\ = \ (H_N^{s+}\subset S_{sN}^+)$$
via some standard identifications.
\end{proposition}

\begin{proof}
In the commutative case our fibration must be a usual fibration of finite spaces, $\{1,\ldots,M\}\to\{1,\ldots,N\}$, commuting with the counting measures. But this shows that our fibration must be of the following special form, with $N,s\in\mathbb N$:
$$\{1,\ldots,N\}\times\{1,\ldots,s\}\to\{1,\ldots,N\}\quad,\quad (i,a)\to i$$

Regarding now the quantum symmetry group, we have the following formula for it, with $e:\mathbb C^N\otimes\mathbb C^s\to\mathbb C^N$ being the Jones projection for the inclusion $\mathbb C^N\subset\mathbb C^N\otimes\mathbb C^s$:
$$C(S_{Z\to Y}^+)=C(S_{sN}^+)\Big/\Big<e\in End(u)\Big>$$

On the other hand, recall that the quantum reflection group $H_N^{s+}\subset S_{sN}^+$ appears via the condition that the corresponding magic matrix must be sudoku:
$$u=\begin{pmatrix}
a^0&a^1&\ldots&a^{s-1}\\
a^{s-1}&a^0&\ldots&a^{s-2}\\
\vdots&\vdots&&\vdots\\
a^1&a^2&\ldots&a^0
\end{pmatrix}$$

But, as explained in \cite{ba2}, this is the same as saying that the quantum group $H_N^{s+}\subset S_{sN}^+$ appears as the symmetry group of the multi-simplex associated to the fibration $\{1,\ldots,N\}\times\{1,\ldots,s\}\to\{1,\ldots,N\}$, so we have an identification as follows:
$$(S_{Z\to Y}^+\subset S_Z^+)\ = \ (H_N^{s+}\subset S_{sN}^+)$$

Thus, we are led to the conclusions in the statement.
\end{proof}

Observe that in Proposition 15.35 the fibration $Z\to Y$ is ``trivial'', in the sense that it is of the following special form:
$$Y\times T\to Y\quad,\quad (i,a)\to i$$

However, in the general quantum case, there are many interesting fibrations $Z\to Y$ which are not trivial, and in what follows we will not make any assumption on our fibrations, and use Definition 15.32 and Theorem 15.34 as stated.

\bigskip

Following \cite{ba2}, we will prove now that the Tannakian category of $S_{Z\to Y}^+$, which is by definition a generalization of $S_Z^+$, is the Fuss-Catalan category, which is a generalization of the Temperley-Lieb category, introduced by Bisch and Jones in \cite{bjo}.

\bigskip

In order to do so, let us first reformulate Theorem 15.34 in a more convenient way, in purely functional analytic terms, and also as a self-contained statement, as follows: 

\begin{theorem}
Any Markov inclusion of finite dimensional algebras $D\subset B$ has a quantum symmetry group $S_{D\subset B}^+$. The corresponding Woronowicz algebra is generated by the coefficients of a biunitary matrix $v=(v_{ij})$ subject to the conditions
$$m\in Hom(v^{\otimes 2},v)\quad,\quad 
u\in Hom(1,v)\quad,\quad 
e\in End(v)$$
where $m:B\otimes B\to B$ is the multiplication, $u:\mathbb C\to B$ is the unit and $e:B\to B$ is the projection onto $D$, with respect to the scalar product $<x,y>=tr(xy^*)$.
\end{theorem}

\begin{proof}
This is a reformulation of Theorem 15.34, with several modifications made. Indeed, by using the algebras $D=C(Y)$, $B=C(Z)$ instead of the quantum spaces $Y,Z$ used there, and also by calling the fundamental corepresentation $v=(v_{ij})$, in order to avoid confusion with the unit $u:\mathbb C\to B$, the formula in Theorem 15.28 reads:
$$C(S_{D\subset B}^+)=C(S_B^+)\Big/\Big<e\in End(v)\Big>$$

Also, we know from Theorem 15.5 that we have the following formula, again by using $B$ instead of $Z$, and by calling the fundamental corepresentation $v=(v_{ij})$:
$$C(S_B^+)=C(U_N^+)\Big/\Big<m\in Hom(v^{\otimes2},v),u\in Fix(v)\Big>$$

Thus, we are led to the conclusion in the statement.
\end{proof}

Let us first discuss in detail the Temperley-Lieb algebra, as a continuation of the material above. In the present context, we have the following definition:

\index{Temperley-Lieb algebra}

\begin{definition}
The $\mathbb N$-algebra $TL^2$ of index $\delta >0$ is defined as follows:
\begin{enumerate}
\item The space $TL^2(m,n)$ consists of linear combinations
of noncrossing pairings between $2m$ points and $2n$ points:
$${{TL^2}}(m,n)=
\left\{\sum\,\alpha\,\,
\begin{matrix}\cdots\cdots&\leftarrow&2m\,\,\ {\rm points}\\
\mathfrak{W}&\leftarrow&m+n\ \ {\rm strings}\\
\cdot\,\cdot&\leftarrow&2n\,\,{\rm points}
\end{matrix}\right\}$$

\item The operations $\circ$, $\otimes$, $*$ are induced by the vertical and horizontal concatenation and the upside-down turning of diagrams:
$$A\circ B=\binom{B}{A}\quad,\quad A\otimes B=AB\quad,\quad A^*=\forall$$

\item With the rule $\bigcirc=\delta$, erasing a circle is the same as multiplying by $\delta$.
\end{enumerate}
\end{definition}

Our first task will be that of finding a suitable presentation for this algebra. Consider the following two elements $u\in TL^2(0,1)$ and $m\in TL^2(2,1)$:
$$u=\delta^{-\frac{1}{2}}\,\cap
\quad,\quad
m=\delta^{\frac{1}{2}}\,|\cup|$$ 

With this convention, we have the following result:

\begin{theorem}
The following relations are a presentation of $TL^2$ by the above rescaled diagrams $u\in TL^2(0,1)$ and $m\in TL^2(2,1)$:
\begin{enumerate}
\item $mm^*=\delta^2$.

\item $u^*u=1$.

\item $m(m\otimes1)=m(1\otimes m)$.

\item $m(1\otimes u)=m(u\otimes1)=1$.

\item $(m\otimes1)(1\otimes m^*)=(1\otimes m)(m^*\otimes1)=m^*m$.
\end{enumerate}
\end{theorem}

\begin{proof}
This is something well-known, and elementary, obtained by drawing diagrams, and for details here, we refer for instance to \cite{bjo}.
\end{proof}

In more concrete terms, the above result says that $u,m$ satisfy the above relations, which is something clear, and that if $C$ is a $\mathbb N$-algebra and $v\in C(0,1)$ and $n\in C(2,1)$ satisfy the same relations then there exists a $\mathbb N$-algebra morphism as follows:
$$TL^2\to C\quad,\quad 
u\to v\quad,\quad 
m\to n$$

Now let $B$ be a finite dimensional $C^*$-algebra, with its canonical trace. We have a scalar product $<x,y>=tr(xy^*)$ on $B$, so $B$ is an object in the category of finite dimensional Hilbert spaces. Consider the unit $u$ and the multiplication $m$ of $B$:
$$u\in\mathbb NB(0,1)\quad,\quad 
m\in\mathbb NB(2,1)$$

The relations in Theorem 15.38 are then satisfied, and one can deduce from this that in this case, the category of representations of $S_B^+$ is the completion of $TL^2$, as we already know. Getting now to Fuss-Catalan algebras, we have here:

\index{Fuss-Catalan algebra}
\index{bicolored pairings}

\begin{definition}
A Fuss-Catalan diagram is a planar diagram formed by an upper row of $4m$ points, a lower row of $4n$ points, both colored
$$\circ\bullet\bullet\circ\circ\bullet\bullet\ldots$$
and by $2m+2n$ noncrossing strings joining these $4m+4n$ points, with the rule that the points which are joined must have the same color.
\end{definition}

Fix $\beta>0$ and $\omega>0$. The $\mathbb N$-algebra $FC$ is defined as follows. The spaces $FC(m,n)$ consist of linear combinations of Fuss-Catalan diagrams:
$$FC(m,n)=\left\{\sum\,\alpha\,\,
\begin{matrix}\circ\bullet\bullet\circ\circ\bullet\bullet\circ\ldots\ldots&\leftarrow&4m\,\,\ {\rm colored\ points}\\
\ &\ &m+n {\rm\ black\ strings}\\
\mathfrak{W}\mbox{\ \ }&\leftarrow&{\rm and}\\
\ & \ & m+n{\rm\ white\ strings}\\
\circ\bullet\bullet\circ\circ\bullet\bullet\circ\ldots\ldots&\leftarrow&4n\,\,{\rm colored\ points}\end{matrix}\right\}$$

As before with the Temperley-Lieb algebra, the operations $\circ$, $\otimes$, $*$ are induced by vertical and horizontal concatenation and upside-down turning of diagrams, but this time with the rule that erasing a black/white circle is the same as multiplying by $\beta$/$\omega$:
$$A\circ B=\binom{B}{A}\quad,\quad 
A\otimes B=AB\quad,\quad 
A^*=\forall$$
$$\mbox{\tiny{black}}\rightarrow\bigcirc=\beta
\qquad,\qquad
\mbox{\tiny{white}}\rightarrow\bigcirc=\omega$$

Let $\delta=\beta\omega$. We have the following bicolored analogues of the elements $u,m$:
$$u=\delta^{-\frac{1}{2}}\,\,\bigcap\;\!\!\!\!\!\!\!\!\cap
\quad,\quad
m=\delta^{\frac{1}{2}}\,\,||\bigcup\;\!\!\!\!\!\!\!\!\cup\ ||$$

Consider also the black and white Jones projections, namely:
$$e=\omega^{-1}\,\,|\ \begin{matrix}\cup\cr\cap\end{matrix}\ |
\quad,\quad 
f=\beta^{-1}\,\,|||\ \begin{matrix}\cup\cr\cap\end{matrix}\ |||$$

For simplifying writing we identify $x$ and $x\otimes 1$. We have the following result:

\begin{theorem}
The following relations, with $f=\beta^{-2}(1\otimes me)m^*$, are a presentation of $FC$ by $m\in FC(2,1)$, $u\in FC(0,1)$ and $e\in FC(1)$:
\begin{enumerate}
\item The relations in Theorem 15.38, with $\delta =\beta\omega$.

\item $e=e^2=e^*$, $f=f^*$ and $(1\otimes f)f=f(1\otimes f)$.

\item $eu=u$.

\item $mem^*=m(1\otimes e)m^*=\beta^2$.

\item $mm(e\otimes e\otimes e)=emm(e\otimes1\otimes e)$.
\end{enumerate}
\end{theorem}

\begin{proof}
This is indeed something quite routine.
\end{proof}

Getting back now to the inclusions $D\subset B$, we have the following result:

\begin{theorem}
Given a Markov inclusion $D\subset B$, we have
$$<m,u,e>=FC$$
as an equality of $\mathbb N$-algebras.
\end{theorem}

\begin{proof}
It is routine to check that the linear maps $m,u,e$ associated to an inclusion $D\subset B$ as in the statement satisfy the relations (1-5) in Theorem 15.40. Thus, we obtain a certain $\mathbb N$-algebra surjective morphism, as follows:
$$J:FC\to<m,u,e>$$

But it is routine to prove that this morphism $J$ is faithful on $\Delta FC$, and then by Frobenius reciprocity faithfulness has to hold on the whole $FC$.
\end{proof}

Getting back now to quantum groups, we have:

\index{planar algebra}
\index{Fuss-Catalan algebra}
\index{twisted reflection group}

\begin{theorem}
Given a Markov inclusion $D\subset B$, the category of representations of its quantum symmetry group $S_{D\subset B}^+$ is the completion of $FC$.
\end{theorem}

\begin{proof}
Since $S_{D\subset B}^+$ comes by definition from the relations corresponding to $m,u,e$, its tensor category of corepresentations is the completion of the tensor category $<m,u,e>$. Thus Theorem 15.41 applies, and gives an isomorphism $<m,u,e>\simeq FC$.
\end{proof}

In terms of finite quantum spaces and quantum graphs, the conclusion is that the quantum automorphism groups $S_{Z\to Y}^+$ of the Markov fibrations $Z\to Y$, which can be thought of as being the ``twisted versions'' of the quantum reflection groups $H_N^{s+}$, correspond to the Fuss-Catalan algebras. We refer here to \cite{ba2} and related papers. 

\section*{15e. Exercises} 

Welcome to non-trivial mathematics, as mathematicians are supposed to like, and we are now in position of formulating a difficult exercise, for you reader, namely:

\begin{exercise}
Come up with a nice and general super-easiness theory, covering the various quantum groups from the previous 2 chapters, and from the present chapter, namely: twists, symplectic groups, and general quantum permutations and reflections.
\end{exercise}

Be said in passing, even if you solve this difficult exercise, things will be not over yet. Indeed, we can still talk about exceptional Lie groups, or about Drinfeld-Jimbo deformations, with parameter of your choice, real, or even better, root of unity. 

\bigskip

Cat says he most likes the perspective of looking into exceptional Lie groups, and that there should be some good questions here, in relation with liberation and twisting. So, unless you solve this specific problem, don't expect much recognition from the felines.

\chapter{Easy geometry}

\section*{16a. Easy geometries}

We have kept the best for the end. All the theory developed in this book, be that easy or super-easy, concerns quantum groups. But this is just the tip of the iceberg, because these quantum groups belong to noncommutative geometries, that we can study too. 

\bigskip

As an example here, $O_N$ certainly belongs to the real geometry, that of $\mathbb R^N$, and $U_N$ certainly belongs to the complex geometry, that of $\mathbb C^N$, But this suggests that $O_N^+$ should belong to a certain ``free real geometry'', that of $\mathbb R^N_+$, whatever this beast means, then $U_N^+$ should belong to a certain ``free complex geometry'', that of $\mathbb C^N_+$, and so on.

\bigskip

In practice now, in order to get started, we must axiomatize the abstract notion of ``noncommutative geometry''. And things are quite tricky here, because there is no hope of having quantum spaces of type $\mathbb R^N_+$ or $\mathbb C^N_+$, simply because the coordinates on these, in the $C^*$-algebra sense, would be unbounded. Thus, we must find something else.

\bigskip

An idea here would be that of restricting the attention to the spheres. That is, replacing $\mathbb R^N,\mathbb C^N$ by the corresponding spheres $S^{N-1}_\mathbb R,S^{N-1}_\mathbb C$, then replacing as well $\mathbb R^N_+,\mathbb C^N_+$ by the corresponding spheres $S^{N-1}_{\mathbb R,+},S^{N-1}_{\mathbb C,+}$, that we met in chapter 13, given by:
$$C(S^{N-1}_{\mathbb R,+})
=C^*\left(x_1,\ldots,x_N\Big|x_i=x_i^*,\sum_ix_i^2=1\right)$$
$$C(S^{N-1}_{\mathbb C,+})
=C^*\left(x_1,\ldots,x_N\Big|\sum_ix_ix_i^*=\sum_ix_i^*x_i=1\right)$$

This idea seems to work fine, but when getting to more complicated spheres $S$, for instance those coming from the various intermediate liberations $U_N\subset U\subset U_N^+$, there are some troubles with the correspondence $S\leftrightarrow U$, and also with the associated torus $T\subset S$, and reflection group $K\subset U$, with the overall problem being that $(S,T,U,K)$ does not seem to always satisfy the natural conditions that this quadruplet satisfies, in the usual cases, namely classical real, classical complex, free real and free complex.

\bigskip

In short, we are stuck, and we must ask the cat. And cat says:

\begin{cat}
Noncommutative geometries are reserved to felines. But you can have a taste of them by looking at the quadruplets $(S,T,U,K)$ consisting of a quantum sphere, torus, unitary group and reflection group, with correspondences between them
$$
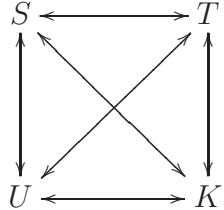
\xymatrix@R=54pt@C=54pt{
S\ar[r]\ar[d]\ar[dr]&T\ar[l]\ar[d]\ar[dl]\\
U\ar[u]\ar[ur]\ar[r]&K\ar[l]\ar[ul]\ar[u]
}$$
similar to those from the classical real and classical complex cases. And when the unitary group $U$ is easy or super-easy, you can call your geometry easy, or super-easy.
\end{cat}

Thanks cat, this looks quite clever, indeed. Getting to work now, let us first check that the classical real and classical complex geometries are indeed noncommutative geometries, in the human sense suggested by cat, and do as well the verification for the free real and free complex geometries. Fortunately all this works well, and we have:

\begin{theorem}
We have basic quadruplets $(S,T,U,K)$, as follows:
\begin{enumerate}
\item A classical real and a classical complex quadruplet, as follows:
$$
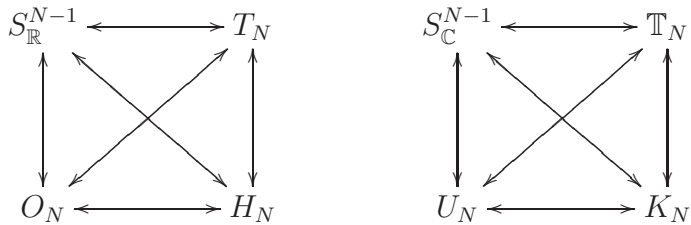
\xymatrix@R=50pt@C=50pt{
S^{N-1}_\mathbb R\ar[r]\ar[d]\ar[dr]&T_N\ar[l]\ar[d]\ar[dl]\\
O_N\ar[u]\ar[ur]\ar[r]&H_N\ar[l]\ar[ul]\ar[u]}
\qquad\qquad 
\xymatrix@R=50pt@C=50pt{
S^{N-1}_\mathbb C\ar[r]\ar[d]\ar[dr]&\mathbb T_N\ar[l]\ar[d]\ar[dl]\\
U_N\ar[u]\ar[ur]\ar[r]&K_N\ar[l]\ar[ul]\ar[u]}$$

\item A free real and a free complex quadruplet, as follows:
$$
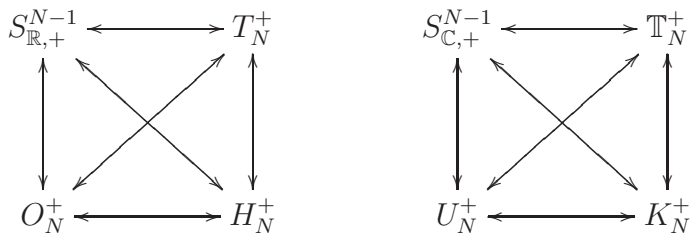
\xymatrix@R=50pt@C=50pt{
S^{N-1}_{\mathbb R,+}\ar[r]\ar[d]\ar[dr]&T_N^+\ar[l]\ar[d]\ar[dl]\\
O_N^+\ar[u]\ar[ur]\ar[r]&H_N^+\ar[l]\ar[ul]\ar[u]}
\qquad\qquad
\xymatrix@R=50pt@C=50pt{
S^{N-1}_{\mathbb C,+}\ar[r]\ar[d]\ar[dr]&\mathbb T_N^+\ar[l]\ar[d]\ar[dl]\\
U_N^+\ar[u]\ar[ur]\ar[r]&K_N^+\ar[l]\ar[ul]\ar[u]}$$
\end{enumerate}
Moreover, for all these quadruplets, the unitary quantum group $U$ is easy.
\end{theorem}

\begin{proof}
Here the various objects appearing in the above diagrams are objects that we know well, constructed at various places, in this book, and the last assertion, regarding easiness, is something that we know well too. As for the fact that, in each of the 4 cases under investigation, we have indeed a full set of 12 correspondences between these objects, this is something quite routine, and for details here, we refer to \cite{ba3}.
\end{proof}

With the advice Cat 16.1 in mind, we can reformulate Theorem 16.2 as follows:

\begin{theorem}
We have $4$ basic noncommutative geometries,
$$\xymatrix@R=52pt@C=52pt{
\mathbb R^N_+\ar[r]&\mathbb C^N_+\\
\mathbb R^N\ar[u]\ar[r]&\mathbb C^N\ar[u]
}$$
called classical real and complex, and free real and complex.
\end{theorem}

\begin{proof}
This is indeed a reformulation of Theorem 16.2, with the convention that a quadruplet as there corresponds to a noncommutative geometry, that we can denote and call as we like, and with the best notations and terminology being as above.
\end{proof}

So far so good, nothing spectacular, we just have here confirmation of a fact that we know well, and this since chapter 2 of the present book, namely that $O_N,U_N$ have free analogues $O_N^+,U_N^+$, and so everything classical is supposed to have a free counterpart.

\bigskip

Getting now to more technical aspects, as already mentioned before receiving the advice Cat 16.1, when trying to construct noncommutative geometries by starting with more complicated unitary quantum groups $U$, things can be quite tricky. Clarifying all this was in fact an open problem, all over the 10s, with several papers written on the subject. In order to discuss the solution, which came in the late 10s, let us start with more details regarding our axioms for noncommutative geometries. As explained in \cite{ba3}, these axioms, formally replacing the vague indications Cat 16.1, are as follows:

\begin{definition}
A quadruplet $(S,T,U,K)$ is said to produce a noncommutative geometry when one can pass from each object to all the other objects, as follows,
$$\begin{matrix}
S&=&S_{<O_N,T>}&=&S_U&=&S_{<O_N,K>}\\
\\
S\cap\mathbb T_N^+&=&T&=&U\cap\mathbb T_N^+&=&K\cap\mathbb T_N^+\\
\\
G^+(S)&=&<O_N,T>&=&U&=&<O_N,K>\\
\\
K^+(S)&=&K^+(T)&=&U\cap K_N^+&=&K
\end{matrix}$$
with the usual convention that all this is up to the equivalence relation.
\end{definition}

There axioms can look a bit complicated, at a first glance, but they are in fact very simple and natural, inspired from what happens in the classical case, and with a look at the free case too. To be more precise, what we have above are all sorts of objects and operations that we know well, along with the operation $U\to S_U$, which consists in taking the first row space, the operation $X\to G^+(X)$, which consists in taking the quantum isometry group, and finally the operation $X\to K^+(X)$, which consists in taking the reflection isometry group, $K^+(X)=G^+(X)\cap K_N^+$. As for the fact that these axioms are indeed satisfied in the 4 main cases of interest, as claimed in Theorem 16.2, this is something well-known in the classical case, and is routine to check in the free case.

\bigskip

In the easy case now, in order to reach to concrete results, it is convenient to formulate an independent definition, using the easy generation operation $\{\,,\}$ instead of the usual generation operation $<\,,>$. This definition, which is more or less a particular case of Definition 16.4, modulo the usual issues with the differences between $\{\,,\}$ and $<\,,>$, that we have met several times, and are very familiar with, is as follows:

\begin{definition}
A quadruplet $(S,T,U,K)$ is said to produce an easy geometry when $U,K$ are easy, and one can pass from each object to all the other objects, as follows,
$$\begin{matrix}
S&=&S_{\{O_N,K^+(T)\}}&=&S_U&=&S_{\{O_N,K\}}\\
\\
S\cap\mathbb T_N^+&=&T&=&U\cap\mathbb T_N^+&=&K\cap\mathbb T_N^+\\
\\
G^+(S)&=&\{O_N,K^+(T)\}&=&U&=&\{O_N,K\}\\
\\
K^+(S)&=&K^+(T)&=&U\cap K_N^+&=&K
\end{matrix}$$
with the usual convention that all this is up to the equivalence relation.
\end{definition}

Getting now into classification results, the idea is to focus on the quantum group content of the above definition. Indeed, we know that both the quantum groups $U,K$ are easy, and that the following easy generation formula must be satisfied:
$$U=\{O_N,K\}$$

Combinatorially, this leads to the following statement:

\index{easy geometry}
\index{category of partitions}
\index{category of pairings}
\index{easy generation}

\begin{proposition}
An easy geometry is uniquely determined by a pair $(D,E)$ of categories of partitions, which must be as follows,
$$\mathcal{NC}_2\subset D\subset P_2\quad,\quad 
\mathcal{NC}_{even}\subset E\subset P_{even}$$
and which are subject to the following intersection and generation conditions,
$$D=E\cap P_2\quad,\quad 
E=<D,\mathcal{NC}_{even}>$$
and to the usual axioms for the associated quadruplet $(S,T,U,K)$, where $U,K$ are respectively the easy quantum groups associated to the categories $D,E$.
\end{proposition}

\begin{proof}
This comes from the following conditions, with the first one being the one mentioned above, and with the second one being part of our general axioms: 
$$U=\{O_N,K\}\quad,\quad 
K=U\cap K_N^+$$

Indeed, $U,K$ must be easy, coming from certain categories of partitions $D,E$. It is clear that $D,E$ must appear as intermediate categories, as in the statement, and the fact that the intersection and generation conditions must be satisfied follows from:
\begin{eqnarray*}
U=\{O_N,K\}&\iff&D=E\cap P_2\\
K=U\cap K_N^+&\iff&E=<D,\mathcal{NC}_{even}>
\end{eqnarray*}

Thus, we are led to the conclusion in the statement.
\end{proof}

In order to discuss now classification results, we would need some technical results regarding the intermediate easy quantum groups as follows:
$$O_N\subset U\subset U_N^+\quad,\quad 
H_N\subset K\subset K_N^+$$

But we do have such results, as explained in Part II and Part III of the present book, and by using this, we are led to the following classification result:

\begin{theorem}
Under strong combinatorial axioms, of easiness and uniformity type, we have only $9$ noncommutative geometries, namely:
$$
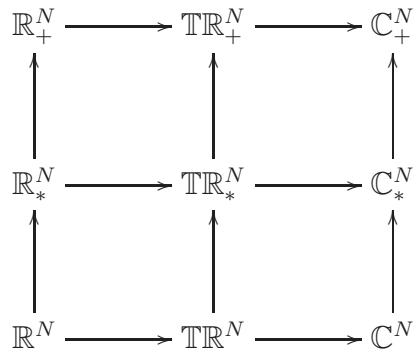
\xymatrix@R=40pt@C=40pt{
\mathbb R^N_+\ar[r]&\mathbb T\mathbb R^N_+\ar[r]&\mathbb C^N_+\\
\mathbb R^N_*\ar[u]\ar[r]&\mathbb T\mathbb R^N_*\ar[u]\ar[r]&\mathbb C^N_*\ar[u]\\
\mathbb R^N\ar[u]\ar[r]&\mathbb T\mathbb R^N\ar[u]\ar[r]&\mathbb C^N\ar[u]
}$$
Moreover, under even stronger combinatorial axioms, including a slicing condition, the $4$ basic geometries, those at the corners, are the only ones.
\end{theorem}

\begin{proof}
This is something quite technical, for which we refer to \cite{ba3}, but that we can basically understand with our quantum group knowledge, the idea being as follows:

\medskip

(1) To start with, the geometries in the statement correspond to the main examples of intermediate quantum groups $O_N\subset U\subset U_N^+$, that we met in Part II, namely:
$$\xymatrix@R=13mm@C=13mm{
O_N^+\ar[r]&\mathbb TO_N^+\ar[r]&U_N^+\\
O_N^*\ar[r]\ar[u]&\mathbb TO_N^*\ar[r]\ar[u]&U_N^*\ar[u]\\
O_N\ar[r]\ar[u]&\mathbb TO_N\ar[r]\ar[u]&U_N\ar[u]}$$

(2) As for the corresponding reflection groups, these correspond to the main examples of intermediate quantum groups $H_N\subset U\subset K_N^+$, that we met in Part III, namely:
$$\xymatrix@R=13mm@C=13mm{
H_N^+\ar[r]&\mathbb TH_N^+\ar[r]&K_N^+\\
H_N^*\ar[r]\ar[u]&\mathbb TH_N^*\ar[r]\ar[u]&K_N^*\ar[u]\\
H_N\ar[r]\ar[u]&\mathbb TH_N\ar[r]\ar[u]&K_N\ar[u]}$$

(3) With these conventions made, telling us who the quantum groups $U,K$ are, in each of the 9 cases under investigation, we can complete our quadruplets with objects $S,T$, by using either of the formulae involving them from Definition 16.5, and the verification of the axioms from Definition 16.5 is straightforward, in each of these 9 cases.

\medskip

(4) Finally, the classification assertions are more technical, whose proofs are basically based on the study of the correspondence $U\leftrightarrow K$ coming from Definition 16.5. To be more precise, we know from Proposition 16.6 that the unitary group $U$ of our easy geometry must come from a category of pairings $D\subset P_2$ satisfying the following condition:
$$D=<D,\mathcal{NC}_{even}>\cap  P_2$$

(5) But this equation can be solved by using the classification results discussed in Part II and Part III of the present book, and we are led to the conclusions in the statement, with the uniformity axiom there being something that we know well, from chapter 6, and with the slicing axiom being something that we are familiar with too, from chapter 7.

\medskip

(6) So, this was for the idea, and in practice now, all this needs a massive amount of routine verifications, at each single step, and all this is explained in \cite{ba3}. 
\end{proof}

\section*{16b. Free manifolds}

With the above discussion done, we are left with the question of picking up the potentially truly interesting geometries, from the list of 9 geometries from Theorem 16.7, and then getting to work, and developing these geometries. Leaving aside the classical geometries, which are not our business, in this book, and leaving aside as well the ``hybrid'' geometries, those appearing on the middle vertical, which do not look that interesting, we are left with 4 geometries, namely the half-classical and free ones, as follows:
$$\xymatrix@R=52pt@C=52pt{
\mathbb R^N_+\ar[r]&\mathbb C^N_+\\
\mathbb R^N_*\ar[u]\ar[r]&\mathbb C^N_*\ar[u]
}$$

These 4 geometries are all interesting, and it is possible to say as few common things about them, but previous experience from the present book with the corresponding unitary groups $U$ suggests that we should split our study, on one hand discussing the free geometries, with whatever ``free tools'' that we can find for studying them, and on the other hand discussing the half-classical geometries, which are quite close to the classical geometries, with tools inspired from classical geometry. So, this will be our plan, with freeness coming first, and leaving the half-classical geometries for later.

\bigskip

Getting started now, with free geometry, we would like to enlarge our collection of free manifolds, which for the moment consists of the 4 basic objects, $S,T,U,K$, plus of course of the various other free quantum groups investigated before in this book, which obviously belong to free geometry too. In practice, leaving the quantum groups aside, this leads us into the question of unifying the spheres $S$ with the unitary groups $U$.

\bigskip

So, this will be our first task, finding a suitable collection of ``free homogeneous spaces'', generalizing at the same time the free spheres $S$, and the free unitary groups $U$. This can be done at several levels of generality, and central here is the construction of the free spaces of partial isometries, which can be done in fact for any easy quantum group. In order to explain this, let us start with the classical case. We have here:

\begin{definition}
Associated to any integers $L\leq M,N$ are the spaces
$$O_{MN}^L=\left\{T:E\to F\ {\rm isometry}\Big|E\subset\mathbb R^N,F\subset\mathbb R^M,\dim_\mathbb RE=L\right\}$$
$$U_{MN}^L=\left\{T:E\to F\ {\rm isometry}\Big|E\subset\mathbb C^N,F\subset\mathbb C^M,\dim_\mathbb CE=L\right\}$$
where the notion of isometry is with respect to the usual real/complex scalar products.
\end{definition}

As a first observation, at $L=M=N$ we obtain the groups $O_N,U_N$:
$$O_{NN}^N=O_N\quad,\quad 
U_{NN}^N=U_N$$ 

Another interesting specialization is $L=M=1$. Here the elements of $O_{1N}^1$ are the isometries $T:E\to\mathbb R$, with $E\subset\mathbb R^N$ one-dimensional. But such an isometry is uniquely determined by $T^{-1}(1)\in\mathbb R^N$, which must belong to $S^{N-1}_\mathbb R$. Thus, we have $O_{1N}^1=S^{N-1}_\mathbb R$. Similarly, in the complex case we have $U_{1N}^1=S^{N-1}_\mathbb C$, and so our results here are:
$$O_{1N}^1=S^{N-1}_\mathbb R\quad,\quad 
U_{1N}^1=S^{N-1}_\mathbb C$$

Yet another interesting specialization is $L=N=1$. Here the elements of $O_{1N}^1$ are the isometries $T:\mathbb R\to F$, with $F\subset\mathbb R^M$ one-dimensional. But such an isometry is uniquely determined by $T(1)\in\mathbb R^M$, which must belong to $S^{M-1}_\mathbb R$. Thus, we have $O_{M1}^1=S^{M-1}_\mathbb R$. Similarly, in the complex case we have $U_{M1}^1=S^{M-1}_\mathbb C$, and so our results here are:
$$O_{M1}^1=S^{M-1}_\mathbb R\quad,\quad
U_{M1}^1=S^{M-1}_\mathbb C$$

In general, the most convenient is to view the elements of $O_{MN}^L,U_{MN}^L$ as rectangular matrices, and to use matrix calculus for their study. We have indeed:

\begin{proposition}
We have identifications of compact spaces
$$O_{MN}^L\simeq\left\{U\in M_{M\times N}(\mathbb R)\Big|UU^t={\rm projection\ of\ trace}\ L\right\}$$
$$U_{MN}^L\simeq\left\{U\in M_{M\times N}(\mathbb C)\Big|UU^*={\rm projection\ of\ trace}\ L\right\}$$
with each partial isometry being identified with the corresponding rectangular matrix.
\end{proposition}

\begin{proof}
We can indeed identify the partial isometries $T:E\to F$ with their corresponding extensions $U:\mathbb R^N\to\mathbb R^M$, $U:\mathbb C^N\to\mathbb C^M$, obtained by setting $U_{E^\perp}=0$. Then, we can identify these latter maps $U$ with the corresponding rectangular matrices.
\end{proof}

As an illustration, at $L=M=N$ we recover in this way the usual matrix description of $O_N,U_N$. Also, at $L=M=1$ we obtain the usual description of $S^{N-1}_\mathbb R,S^{N-1}_\mathbb C$, as row spaces over the corresponding groups $O_N,U_N$. Finally, at $L=N=1$ we obtain the usual description of $S^{N-1}_\mathbb R,S^{N-1}_\mathbb C$, as column spaces over the corresponding groups $O_N,U_N$. 

\bigskip

Now back to the general case, observe that the isometries $T:E\to F$, or rather their extensions $U:\mathbb K^N\to\mathbb K^M$, with $\mathbb K=\mathbb R,\mathbb C$, obtained by setting $U_{E^\perp}=0$, can be composed with the isometries of $\mathbb K^M,\mathbb K^N$, according to the following scheme:
$$
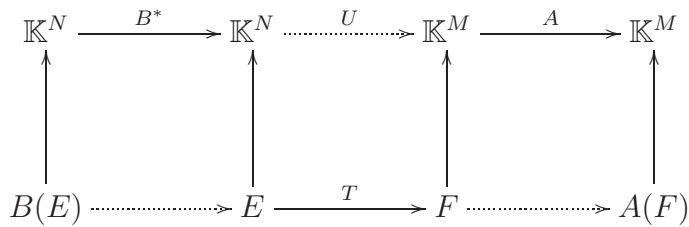
\xymatrix@R=17mm@C=17mm{
\mathbb K^N\ar[r]^{B^*}&\mathbb K^N\ar@.[r]^U&\mathbb K^M\ar[r]^A&\mathbb K^M\\
B(E)\ar@.[r]\ar[u]&E\ar[r]^T\ar[u]&F\ar@.[r]\ar[u]&A(F)\ar[u]
}$$

With the identifications in Proposition 16.9 made, the precise statement here is:

\begin{proposition}
We have action maps as follows, which are both transitive,
$$O_M\times O_N\curvearrowright O_{MN}^L\quad,\quad 
(A,B)U=AUB^t$$
$$U_M\times U_N\curvearrowright U_{MN}^L\quad,\quad 
(A,B)U=AUB^*$$
whose stabilizers are respectively $O_L\times O_{M-L}\times O_{N-L}$ and $U_L\times U_{M-L}\times U_{N-L}$.
\end{proposition}

\begin{proof}
We have indeed action maps as in the statement, which are transitive. Let us compute now the stabilizer $G$ of the following point:
$$U=\begin{pmatrix}1&0\\0&0\end{pmatrix}$$

Since $(A,B)\in G$ satisfy $AU=UB$, their components must be of the following form:
$$A=\begin{pmatrix}x&*\\0&a\end{pmatrix}\quad,\quad 
B=\begin{pmatrix}x&0\\ *&b\end{pmatrix}$$

Now since $A,B$ are unitaries, these matrices follow to be block-diagonal, and so:
$$G=\left\{(A,B)\Big|A=\begin{pmatrix}x&0\\0&a\end{pmatrix},B=\begin{pmatrix}x&0\\ 0&b\end{pmatrix}\right\}$$

The stabilizer of $U$ is parametrized by triples $(x,a,b)$ belonging to $O_L\times O_{M-L}\times O_{N-L}$ and $U_L\times U_{M-L}\times U_{N-L}$, and we are led to the conclusion in the statement.
\end{proof}

Finally, let us work out the quotient space description of $O_{MN}^L,U_{MN}^L$. We have here:

\begin{theorem}
We have isomorphisms of homogeneous spaces as follows,
\begin{eqnarray*}
O_{MN}^L&=&(O_M\times O_N)/(O_L\times O_{M-L}\times O_{N-L})\\
U_{MN}^L&=&(U_M\times U_N)/(U_L\times U_{M-L}\times U_{N-L})
\end{eqnarray*}
with the quotient maps being given by $(A,B)\to AUB^*$, where $U=(^1_0{\ }^0_0)$.
\end{theorem}

\begin{proof}
This is just a reformulation of Proposition 16.10, by taking into account the fact that the fixed point used in the proof there was $U=(^1_0{\ }^0_0)$.
\end{proof}

Once again, the basic examples here come from the cases $L=M=N$ and $L=M=1$. At $L=M=N$ the quotient spaces at right are respectively:
$$O_N\quad,\quad U_N$$

At $L=M=1$  the quotient spaces at right are respectively:
$$O_N/O_{N-1}\quad,\quad U_N/U_{N-1}$$

In fact, in the general $L=M$ case we obtain the following spaces:
$$O_{MN}^M=O_N/O_{N-M}
\quad,\quad 
U_{MN}^M=U_N/U_{N-M}$$

Similarly, the examples coming from the cases $L=M=N$ and $L=N=1$ are particular cases of the general $L=N$ case, where we obtain the following spaces:
$$O_{MN}^N=O_N/O_{M-N}
\quad,\quad 
U_{MN}^N=U_N/U_{M-N}$$

Summarizing, we have here some basic homogeneous spaces, unifying the spheres with the rotation groups. The point now is that we can liberate these spaces, as follows:

\begin{definition}
Associated to any integers $L\leq M,N$ are the algebras
\begin{eqnarray*}
C(O_{MN}^{L+})&=&C^*\left((u_{ij})_{i=1,\ldots,M,j=1,\ldots,N}\Big|u=\bar{u},uu^t={\rm projection\ of\ trace}\ L\right)\\
C(U_{MN}^{L+})&=&C^*\left((u_{ij})_{i=1,\ldots,M,j=1,\ldots,N}\Big|uu^*,\bar{u}u^t={\rm projections\ of\ trace}\ L\right)
\end{eqnarray*}
with the trace being by definition the sum of the diagonal entries.
\end{definition}

Observe that the above universal algebras are indeed well-defined, as it was previously  the case for the free spheres, and this due to the trace conditions, which read: 
$$\sum_{ij}u_{ij}u_{ij}^*
=\sum_{ij}u_{ij}^*u_{ij}
=L$$

We have inclusions between the various spaces constructed so far, as follows:
$$\xymatrix@R=15mm@C=15mm{
O_{MN}^{L+}\ar[r]&U_{MN}^{L+}\\
O_{MN}^L\ar[r]\ar[u]&U_{MN}^L\ar[u]}$$

At the level of basic examples now, at $L=M=1$ and at $L=N=1$ we obtain the following diagrams, showing that our formalism covers indeed the free spheres:
$$\xymatrix@R=15mm@C=15mm{
S^{N-1}_{\mathbb R,+}\ar[r]&S^{N-1}_{\mathbb C,+}\\
S^{N-1}_\mathbb R\ar[r]\ar[u]&S^{N-1}_\mathbb C\ar[u]}
\qquad\qquad 
\xymatrix@R=15mm@C=15mm{
S^{M-1}_{\mathbb R,+}\ar[r]&S^{M-1}_{\mathbb C,+}\\
S^{M-1}_\mathbb R\ar[r]\ar[u]&S^{M-1}_\mathbb C\ar[u]}$$

We have as well the following result, in relation with the free rotation groups:

\begin{proposition}
At $L=M=N$ we obtain the diagram
$$\xymatrix@R=15mm@C=15mm{
O_N^+\ar[r]&U_N^+\\
O_N\ar[r]\ar[u]&U_N\ar[u]}$$
consisting of the groups $O_N,U_N$, and their liberations.
\end{proposition}

\begin{proof}
We recall that the various quantum groups in the statement are constructed as follows, with the symbol $\times$ standing once again for ``commutative'' and ``free'':
\begin{eqnarray*}
C(O_N^\times)&=&C^*_\times\left((u_{ij})_{i,j=1,\ldots,N}\Big|u=\bar{u},uu^t=u^tu=1\right)\\
C(U_N^\times)&=&C^*_\times\left((u_{ij})_{i,j=1,\ldots,N}\Big|uu^*=u^*u=1,\bar{u}u^t=u^t\bar{u}=1\right)
\end{eqnarray*}

On the other hand, according to Proposition 16.9 and to Definition 16.12, we have the following presentation results:
\begin{eqnarray*}
C(O_{NN}^{N\times})&=&C^*_\times\left((u_{ij})_{i,j=1,\ldots,N}\Big|u=\bar{u},uu^t={\rm projection\ of\ trace}\ N\right)\\
C(U_{NN}^{N\times})&=&C^*_\times\left((u_{ij})_{i,j=1,\ldots,N}\Big|uu^*,\bar{u}u^t={\rm projections\ of\ trace}\ N\right)
\end{eqnarray*}

We use now the standard fact that if $p=aa^*$ is a projection then $q=a^*a$ is a projection too. We use as well the following formulae:
$$Tr(uu^*)=Tr(u^t\bar{u})\quad,\quad 
Tr(\bar{u}u^t)=Tr(u^*u)$$

We therefore obtain the following formulae:
\begin{eqnarray*}
C(O_{NN}^{N\times})&=&C^*_\times\left((u_{ij})_{i,j=1,\ldots,N}\Big|u=\bar{u},\ uu^t,u^tu={\rm projections\ of\ trace}\ N\right)\\
C(U_{NN}^{N\times})&=&C^*_\times\left((u_{ij})_{i,j=1,\ldots,N}\Big|uu^*,u^*u,\bar{u}u^t,u^t\bar{u}={\rm projections\ of\ trace}\ N\right)
\end{eqnarray*}

Now observe that, in tensor product notation, the conditions at right are all of the form $(tr\otimes id)p=1$. Thus, $p$ must be follows, for the above conditions:
$$p=uu^*,u^*u,\bar{u}u^t,u^t\bar{u}$$

We therefore obtain that, for any faithful state $\varphi$, we have $(tr\otimes\varphi)(1-p)=0$. It follows from this that the following projections must be all equal to the identity:
$$p=uu^*,u^*u,\bar{u}u^t,u^t\bar{u}$$

But this leads to the conclusion in the statement.
\end{proof}

Regarding now the homogeneous space structure of $O_{MN}^{L\times},U_{MN}^{L\times}$, the situation here is a bit more complicated in the free case than in the classical case, due to a number of algebraic and analytic issues. We first have the following result:

\begin{proposition}
The spaces $U_{MN}^{L\times}$ have the following properties:
\begin{enumerate}
\item We have an action $U_M^\times\times U_N^\times\curvearrowright U_{MN}^{L\times}$, given by $u_{ij}\to\sum_{kl}u_{kl}\otimes a_{ki}\otimes b_{lj}^*$.

\item We have a map $U_M^\times\times U_N^\times\to U_{MN}^{L\times}$, given by $u_{ij}\to\sum_{r\leq L}a_{ri}\otimes b_{rj}^*$.
\end{enumerate}
Similar results hold for the spaces $O_{MN}^{L\times}$, with all the $*$ exponents removed.
\end{proposition}

\begin{proof}
In the classical case, consider the following action and quotient maps:
$$U_M\times U_N\curvearrowright U_{MN}^L\quad,\quad 
U_M\times U_N\to U_{MN}^L$$

The transposes of these two maps are as follows, where $J=(^1_0{\ }^0_0)$:
\begin{eqnarray*}
\varphi&\to&((U,A,B)\to\varphi(AUB^*))\\
\varphi&\to&((A,B)\to\varphi(AJB^*))
\end{eqnarray*}

But with $\varphi=u_{ij}$ we obtain precisely the formulae in the statement. The proof in the orthogonal case is similar. Regarding now the free case, the proof goes as follows:

\medskip

(1) Assuming $uu^*u=u$, let us set:
$$U_{ij}=\sum_{kl}u_{kl}\otimes a_{ki}\otimes b_{lj}^*$$

We have then the following computation:
\begin{eqnarray*}
(UU^*U)_{ij}
&=&\sum_{pq}\sum_{klmnst}u_{kl}u_{mn}^*u_{st}\otimes a_{ki}a_{mq}^*a_{sq}\otimes b_{lp}^*b_{np}b_{tj}^*\\
&=&\sum_{klmt}u_{kl}u_{ml}^*u_{mt}\otimes a_{ki}\otimes b_{tj}^*\\
&=&\sum_{kt}u_{kt}\otimes a_{ki}\otimes b_{tj}^*\\
&=&U_{ij}
\end{eqnarray*}

Also, assuming that we have $\sum_{ij}u_{ij}u_{ij}^*=L$, we obtain:
\begin{eqnarray*}
\sum_{ij}U_{ij}U_{ij}^*
&=&\sum_{ij}\sum_{klst}u_{kl}u_{st}^*\otimes a_{ki}a_{si}^*\otimes b_{lj}^*b_{tj}\\
&=&\sum_{kl}u_{kl}u_{kl}^*\otimes1\otimes1\\
&=&L
\end{eqnarray*}

(2) Assuming $uu^*u=u$, let us set:
$$V_{ij}=\sum_{r\leq L}a_{ri}\otimes b_{rj}^*$$

We have then the following computation:
\begin{eqnarray*}
(VV^*V)_{ij}
&=&\sum_{pq}\sum_{x,y,z\leq L}a_{xi}a_{yq}^*a_{zq}\otimes b_{xp}^*b_{yp}b_{zj}^*\\
&=&\sum_{x\leq L}a_{xi}\otimes b_{xj}^*\\
&=&V_{ij}
\end{eqnarray*}

Also, assuming that we have $\sum_{ij}u_{ij}u_{ij}^*=L$, we obtain:
\begin{eqnarray*}
\sum_{ij}V_{ij}V_{ij}^*
&=&\sum_{ij}\sum_{r,s\leq L}a_{ri}a_{si}^*\otimes b_{rj}^*b_{sj}\\
&=&\sum_{l\leq L}1\\
&=&L
\end{eqnarray*}

By removing all the $*$ exponents, we obtain as well the orthogonal results.
\end{proof}

Let us examine now the relation between the above maps. In the classical case, given a quotient space $X=G/H$, the associated action and quotient maps are given by:
$$\begin{cases}
a:X\times G\to X&:\quad (Hg,h)\to Hgh\\
p:G\to X&:\quad g\to Hg
\end{cases}$$

Thus we have $a(p(g),h)=p(gh)$. In our context, a similar result holds: 

\begin{theorem}
With $G=G_M\times G_N$ and $X=G_{MN}^L$, where $G_N=O_N^\times,U_N^\times$, we have
$$\xymatrix@R=15mm@C=30mm{
G\times G\ar[r]^m\ar[d]_{p\times id}&G\ar[d]^p\\
X\times G\ar[r]^a&X
}$$
where $a,p$ are the action map and the map constructed in Proposition 16.14.
\end{theorem}

\begin{proof}
At the level of the associated algebras of functions, we must prove that the following diagram commutes, where $\Phi,\alpha$ are morphisms of algebras induced by $a,p$:
$$\xymatrix@R=15mm@C=25mm{
C(X)\ar[r]^\Phi\ar[d]_\alpha&C(X\times G)\ar[d]^{\alpha\otimes id}\\
C(G)\ar[r]^\Delta&C(G\times G)
}$$

When going right, and then down, the composition is as follows:
\begin{eqnarray*}
(\alpha\otimes id)\Phi(u_{ij})
&=&(\alpha\otimes id)\sum_{kl}u_{kl}\otimes a_{ki}\otimes b_{lj}^*\\
&=&\sum_{kl}\sum_{r\leq L}a_{rk}\otimes b_{rl}^*\otimes a_{ki}\otimes b_{lj}^*
\end{eqnarray*}

On the other hand, when going down, and then right, the composition is as follows, where $F_{23}$ is the flip between the second and the third components:
\begin{eqnarray*}
\Delta\pi(u_{ij})
&=&F_{23}(\Delta\otimes\Delta)\sum_{r\leq L}a_{ri}\otimes b_{rj}^*\\
&=&F_{23}\left(\sum_{r\leq L}\sum_{kl}a_{rk}\otimes a_{ki}\otimes b_{rl}^*\otimes b_{lj}^*\right)
\end{eqnarray*}

Thus the above diagram commutes indeed, and this gives the result.
\end{proof}

Let us discuss now some discrete extensions of the above constructions. We have:

\begin{definition}
Associated to any partial permutation, $\sigma:I\simeq J$ with $I\subset\{1,\ldots,N\}$ and $J\subset\{1,\ldots,M\}$, is the real/complex partial isometry
$$T_\sigma:span\left(e_i\Big|i\in I\right)\to span\left(e_j\Big|j\in J\right)$$
given on the standard basis elements by $T_\sigma(e_i)=e_{\sigma(i)}$.
\end{definition}

Let $S_{MN}^L$ be the set of partial permutations $\sigma:I\simeq J$ as above, with range $I\subset\{1,\ldots,N\}$ and target $J\subset\{1,\ldots,M\}$, and with $L=|I|=|J|$. We have:

\begin{proposition}
The space of partial permutations signed by elements of $\mathbb Z_s$,
$$H_{MN}^{sL}=\left\{T(e_i)=w_ie_{\sigma(i)}\Big|\sigma\in S_{MN}^L,w_i\in\mathbb Z_s\right\}$$
is isomorphic to the quotient space 
$$(H_M^s\times H_N^s)/(H_L^s\times H_{M-L}^s\times H_{N-L}^s)$$
via a standard isomorphism.
\end{proposition}

\begin{proof}
This follows by adapting the computations in the proof of Proposition 16.10 and Theorem 16.11. Indeed, we have an action map as follows, which is transitive:
$$H_M^s\times H_N^s\to H_{MN}^{sL}\quad,\quad 
(A,B)U=AUB^*$$

Consider now the following point:
$$U=\begin{pmatrix}1&0\\0&0\end{pmatrix}$$

The stabilizer of this point follows to be the following group:
$$H_L^s\times H_{M-L}^s\times H_{N-L}^s$$

To be more precise, this group is embedded via:
$$(x,a,b)\to\left[\begin{pmatrix}x&0\\0&a\end{pmatrix},\begin{pmatrix}x&0\\0&b\end{pmatrix}\right]$$

But this gives the result.
\end{proof}

In the free case now, the idea is similar, by using inspiration from the construction of the quantum group $H_N^{s+}=\mathbb Z_s\wr_*S_N^+$ in \cite{bb+}. The result here is as follows:

\begin{proposition}
The compact quantum space $H_{MN}^{sL+}$ associated to the algebra
$$C(H_{MN}^{sL+})=C(U_{MN}^{L+})\Big/\left<u_{ij}u_{ij}^*=u_{ij}^*u_{ij}=p_{ij}={\rm projections},u_{ij}^s=p_{ij}\right>$$
has an action map, and is the target of a quotient map, as in Theorem 16.15.
\end{proposition}

\begin{proof}
We must show that if the variables $u_{ij}$ satisfy the relations in the statement, then these relations are satisfied as well for the following variables: 
$$U_{ij}=\sum_{kl}u_{kl}\otimes a_{ki}\otimes b_{lj}^*\quad,\quad 
V_{ij}=\sum_{r\leq L}a_{ri}\otimes b_{rj}^*$$

We use the fact that the standard coordinates $a_{ij},b_{ij}$ on the quantum groups $H_M^{s+},H_N^{s+}$ satisfy the following relations, for any $x\neq y$ on the same row or column of $a,b$:
$$xy=xy^*=0$$
 
We obtain, by using these relations, the following formula:
$$U_{ij}U_{ij}^*
=\sum_{klmn}u_{kl}u_{mn}^*\otimes a_{ki}a_{mi}^*\otimes b_{lj}^*b_{mj}
=\sum_{kl}u_{kl}u_{kl}^*\otimes a_{ki}a_{ki}^*\otimes b_{lj}^*b_{lj}$$

On the other hand, we have as well the following formula:
$$V_{ij}V_{ij}^*
=\sum_{r,t\leq L}a_{ri}a_{ti}^*\otimes b_{rj}^*b_{tj}
=\sum_{r\leq L}a_{ri}a_{ri}^*\otimes b_{rj}^*b_{rj}$$

In terms of the projections $x_{ij}=a_{ij}a_{ij}^*$, $y_{ij}=b_{ij}b_{ij}^*$, $p_{ij}=u_{ij}u_{ij}^*$, we have:
$$U_{ij}U_{ij}^*=\sum_{kl}p_{kl}\otimes x_{ki}\otimes y_{lj}\quad,\quad 
V_{ij}V_{ij}^*=\sum_{r\leq L}x_{ri}\otimes y_{rj}$$

By repeating the computation, we conclude that these elements are projections. Also, a similar computation shows that $U_{ij}^*U_{ij},V_{ij}^*V_{ij}$ are given by the same formulae. Finally, once again by using the relations of type $xy=xy^*=0$, we have:
$$U_{ij}^s
=\sum_{k_rl_r}u_{k_1l_1}\ldots u_{k_sl_s}\otimes a_{k_1i}\ldots a_{k_si}\otimes b_{l_1j}^*\ldots b_{l_sj}^*
=\sum_{kl}u_{kl}^s\otimes a_{ki}^s\otimes(b_{lj}^*)^s$$

On the other hand, we have as well the following formula:
$$V_{ij}^s
=\sum_{r_l\leq L}a_{r_1i}\ldots a_{r_si}\otimes b_{r_1j}^*\ldots b_{r_sj}^*
=\sum_{r\leq L}a_{ri}^s\otimes(b_{rj}^*)^s$$

Thus the conditions of type $u_{ij}^s=p_{ij}$ are satisfied as well, and we are done.
\end{proof}

Let us discuss now the general case. We have the following result:

\begin{proposition}
The various spaces $G_{MN}^L$ constructed so far appear by imposing to the standard coordinates of $U_{MN}^{L+}$ the relations
$$\sum_{i_1\ldots i_s}\sum_{j_1\ldots j_s}\delta_\pi(i)\delta_\sigma(j)u_{i_1j_1}^{e_1}\ldots u_{i_sj_s}^{e_s}=L^{|\pi\vee\sigma|}$$
with $s=(e_1,\ldots,e_s)$ ranging over all the colored integers, and with $\pi,\sigma\in D(0,s)$.
\end{proposition}

\begin{proof}
According to the various constructions above, the relations defining the quantum space $G_{MN}^L$ can be written as follows, with $\sigma$ ranging over a family of generators, with no upper legs, of the corresponding category of partitions $D$:
$$\sum_{j_1\ldots j_s}\delta_\sigma(j)u_{i_1j_1}^{e_1}\ldots u_{i_sj_s}^{e_s}=\delta_\sigma(i)$$

We therefore obtain the relations in the statement, as follows:
\begin{eqnarray*}
\sum_{i_1\ldots i_s}\sum_{j_1\ldots j_s}\delta_\pi(i)\delta_\sigma(j)u_{i_1j_1}^{e_1}\ldots u_{i_sj_s}^{e_s}
&=&\sum_{i_1\ldots i_s}\delta_\pi(i)\sum_{j_1\ldots j_s}\delta_\sigma(j)u_{i_1j_1}^{e_1}\ldots u_{i_sj_s}^{e_s}\\
&=&\sum_{i_1\ldots i_s}\delta_\pi(i)\delta_\sigma(i)\\
&=&L^{|\pi\vee\sigma|}
\end{eqnarray*}

As for the converse, this follows by using the relations in the statement, by keeping $\pi$ fixed, and by making $\sigma$ vary over all the partitions in the category.
\end{proof}

In the general case now, where $G=(G_N)$ is an arbitary uniform easy quantum group, we can construct spaces $G_{MN}^L$ by using the above relations, and we have:

\begin{theorem}
The spaces $G_{MN}^L\subset U_{MN}^{L+}$ constructed by imposing the relations 
$$\sum_{i_1\ldots i_s}\sum_{j_1\ldots j_s}\delta_\pi(i)\delta_\sigma(j)u_{i_1j_1}^{e_1}\ldots u_{i_sj_s}^{e_s}=L^{|\pi\vee\sigma|}$$
with $\pi,\sigma$ ranging over all the partitions in the associated category, having no upper legs, are subject to an action map/quotient map diagram, as in Theorem 16.15.
\end{theorem}

\begin{proof}
We proceed as in the proof of Proposition 16.18. We must prove that, if the variables $u_{ij}$ satisfy the relations in the statement, then so do the following variables:
$$U_{ij}=\sum_{kl}u_{kl}\otimes a_{ki}\otimes b_{lj}^*\quad,\quad 
V_{ij}=\sum_{r\leq L}a_{ri}\otimes b_{rj}^*$$

Regarding the variables $U_{ij}$, the computation here goes as follows:
\begin{eqnarray*}
&&\sum_{i_1\ldots i_s}\sum_{j_1\ldots j_s}\delta_\pi(i)\delta_\sigma(j)U_{i_1j_1}^{e_1}\ldots U_{i_sj_s}^{e_s}\\
&=&\sum_{i_1\ldots i_s}\sum_{j_1\ldots j_s}\sum_{k_1\ldots k_s}\sum_{l_1\ldots l_s}u_{k_1l_1}^{e_1}\ldots u_{k_sl_s}^{e_s}\otimes \delta_\pi(i)\delta_\sigma(j)a_{k_1i_1}^{e_1}\ldots a_{k_si_s}^{e_s}\otimes(b_{l_sj_s}^{e_s}\ldots b_{l_1j_1}^{e_1})^*\\
&=&\sum_{k_1\ldots k_s}\sum_{l_1\ldots l_s}\delta_\pi(k)\delta_\sigma(l)u_{k_1l_1}^{e_1}\ldots u_{k_sl_s}^{e_s}\\
&=&L^{|\pi\vee\sigma|}
\end{eqnarray*}

For the variables $V_{ij}$ the proof is similar, as follows:
\begin{eqnarray*}
&&\sum_{i_1\ldots i_s}\sum_{j_1\ldots j_s}\delta_\pi(i)\delta_\sigma(j)V_{i_1j_1}^{e_1}\ldots V_{i_sj_s}^{e_s}\\
&=&\sum_{i_1\ldots i_s}\sum_{j_1\ldots j_s}\sum_{l_1,\ldots,l_s\leq L}\delta_\pi(i)\delta_\sigma(j)a_{l_1i_1}^{e_1}\ldots a_{l_si_s}^{e_s}\otimes(b_{l_sj_s}^{e_s}\ldots b_{l_1j_1}^{e_1})^*\\
&=&\sum_{l_1,\ldots,l_s\leq L}\delta_\pi(l)\delta_\sigma(l)\\
&=&L^{|\pi\vee\sigma|}
\end{eqnarray*}

Thus we have constructed an action map, and a quotient map, as in Proposition 16.18, and the commutation of the diagram in Theorem 16.15 is then trivial.
\end{proof}

Summarizing, and getting back now to our general free geometry motivations, leaving the free tori $T$ aside, which are quite special, dually being of ``classical'' nature, we have enlarged our collection of free manifolds $\{S,U,K\}$ to something far more general, consisting of the spaces $\{G_{MN}^L\}$ constructed above. This is of course just the tip of the iceberg, and it is possible to say far more things about this, first with a detailed study of these spaces $G_{MN}^L$, from an algebraic and analytic viewpoint, based on a Weingarten integration formula for them, and then with various generalizations of this formalism, and some axiomatization work as well. For more on all this, we refer to \cite{ba3}.

\section*{16c. Projective geometry}

All the above is quite interesting, but even more interesting is what happens in relation with projective versions. Let us go back to the diagram of 9 main geometries, namely:
$$\xymatrix@R=40pt@C=42pt{
\mathbb R^N_+\ar[r]&\mathbb T\mathbb R^N_+\ar[r]&\mathbb C^N_+\\
\mathbb R^N_*\ar[u]\ar[r]&\mathbb T\mathbb R^N_*\ar[u]\ar[r]&\mathbb C^N_*\ar[u]\\
\mathbb R^N\ar[u]\ar[r]&\mathbb T\mathbb R^N\ar[u]\ar[r]&\mathbb C^N\ar[u]
}$$

These geometries are by definition affine, and our claim is that some drastic simplifications appear when looking at the corresponding projective geometries:
$$\xymatrix@R=34pt@C=33pt{
P^{N-1}_+\ar[r]&P^{N-1}_+\ar[r]&P^{N-1}_+\\
P^{N-1}_\mathbb C\ar[r]\ar[u]&P^{N-1}_\mathbb C\ar[r]\ar[u]&P^{N-1}_\mathbb C\ar[u]\\
P^{N-1}_\mathbb R\ar[r]\ar[u]&P^{N-1}_\mathbb R\ar[r]\ar[u]&P^{N-1}_\mathbb C\ar[u]}$$

Thus, we are led to the conclusion that, under certain combinatorial axioms, there should be only 3 projective geometries, namely the real, complex and free ones:
$$P^{N-1}_\mathbb R\subset P^{N-1}_\mathbb C\subset P^{N-1}_+$$

And isn't this beautiful, what we have here is some sort of ``threefold way'', which looks very conceptual. In order to discuss this, let us start with:

\begin{proposition}
We have presentation results as follows,
\begin{eqnarray*}
C(P^{N-1}_\mathbb R)&=&C^*_{comm}\left((p_{ij})_{i,j=1,\ldots,N}\Big|p=\bar{p}=p^t=p^2,Tr(p)=1\right)\\
C(P^{N-1}_\mathbb C)&=&C^*_{comm}\left((p_{ij})_{i,j=1,\ldots,N}\Big|p=p^*=p^2,Tr(p)=1\right)
\end{eqnarray*}
for the algebras of continuous functions on the real and complex projective spaces.
\end{proposition}

\begin{proof}
We use the fact that the projective spaces $P^{N-1}_\mathbb R,P^{N-1}_\mathbb C$ can be respectively identified with the spaces of rank one projections in $M_N(\mathbb R),M_N(\mathbb C)$. With this picture in mind, it is clear that we have arrows $\leftarrow$. In order to construct now arrows $\to$, consider the universal algebras on the right, $A_R,A_C$. These algebras being both commutative, by the Gelfand theorem we can write, with $X_R,X_C$ being certain compact spaces:
$$A_R=C(X_R)\quad,\quad 
A_C=C(X_C)$$

Now by using the coordinate functions $p_{ij}$, we conclude that $X_R,X_C$ are certain spaces of rank one projections in $M_N(\mathbb R),M_N(\mathbb C)$. In other words, we have embeddings:
$$X_R\subset P^{N-1}_\mathbb R\quad,\quad 
X_C\subset P^{N-1}_\mathbb C$$

By transposing we obtain arrows $\to$, as desired.
\end{proof}

The above result suggests the following definition:

\begin{definition}
Associated to any $N\in\mathbb N$ is the following universal algebra,
$$C(P^{N-1}_+)=C^*\left((p_{ij})_{i,j=1,\ldots,N}\Big|p=p^*=p^2,Tr(p)=1\right)$$
whose abstract spectrum is called ``free projective space''.
\end{definition}

Observe that, according to our presentation results for the real and complex projective spaces $P^{N-1}_\mathbb R$ and $P^{N-1}_\mathbb C$, we have embeddings of compact quantum spaces, as follows:
$$P^{N-1}_\mathbb R\subset P^{N-1}_\mathbb C\subset P^{N-1}_+$$

Our first goal will be that of explaining why, in analogy with the uniqueness of the quantum group $PO_N^+=PU_N^+$, the free projective space $P^{N-1}_+$ is unique, and scalarless. 

\bigskip

Let us first discuss the relation with the various noncommutative spheres. Given a closed subset $X\subset S^{N-1}_{\mathbb R,+}$, its projective version is by definition the quotient space $X\to PX$ determined by the fact that $C(PX)\subset C(X)$ is the subalgebra generated by the variables $p_{ij}=x_ix_j$. With this convention, we have the following result:

\begin{theorem}
The projective versions of the main $9$ spheres are
$$\xymatrix@R=34pt@C=33pt{
P^{N-1}_+\ar[r]&P^{N-1}_+\ar[r]&P^{N-1}_+\\
P^{N-1}_\mathbb C\ar[r]\ar[u]&P^{N-1}_\mathbb C\ar[r]\ar[u]&P^{N-1}_\mathbb C\ar[u]\\
P^{N-1}_\mathbb R\ar[r]\ar[u]&P^{N-1}_\mathbb R\ar[r]\ar[u]&P^{N-1}_\mathbb C\ar[u]}$$
involving only the $3$ projective spaces $P^{N-1}_\mathbb R\subset P^{N-1}_\mathbb C\subset P^{N-1}_+$.
\end{theorem}

\begin{proof}
This is something quite elementary, the idea being that in the free case, the results follow in analogy with the result $PO_N^+=PU_N^+$ that we already know, then in the half-classical case the results follow in analogy with the result $PO_N^*=PU_N$, that we already know too, and finally in the hybrid cases the results are clear. See \cite{ba3}.
\end{proof}

Getting purely projective now, we can axiomatize our spaces, as follows:

\begin{definition}
A monomial projective space is a closed subset $P\subset P^{N-1}_+$ obtained via relations of type
$$p_{i_1i_2}\ldots p_{i_{k-1}i_k}=p_{i_{\sigma(1)}i_{\sigma(2)}}\ldots p_{i_{\sigma(k-1)}i_{\sigma(k)}},\ \forall (i_1,\ldots,i_k)\in\{1,\ldots,N\}^k$$
with $\sigma$ ranging over a certain subset of the infinite symmetric group
$$S_\infty=\bigcup_{k\in2\mathbb N}S_k$$
which is stable under the operation $\sigma\to|\sigma|$.
\end{definition}

Here the stability under the operation $\sigma\to|\sigma|$ means that if the above relation associated to $\sigma$ holds, then the following relation, associated to $|\sigma|$, must hold as well:
$$p_{i_0i_1}\ldots p_{i_ki_{k+1}}=p_{i_0i_{\sigma(1)}}p_{i_{\sigma(2)}i_{\sigma(3)}}\ldots p_{i_{\sigma(k-2)}i_{\sigma(k-1)}}p_{i_{\sigma(k)}i_{k+1}}$$

As an illustration, the basic projective spaces are all monomial:

\begin{proposition}
The $3$ projective spaces are all monomial, with the permutations
$$
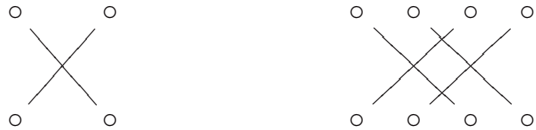
\xymatrix@R=10mm@C=8mm{\circ\ar@{-}[dr]&\circ\ar@{-}[dl]\\\circ&\circ}\qquad\ \qquad\ \qquad 
\xymatrix@R=10mm@C=3mm{\circ\ar@{-}[drr]&\circ\ar@{-}[drr]&\circ\ar@{-}[dll]&\circ\ar@{-}[dll]\\\circ&\circ&\circ&\circ}$$
producing respectively the spaces $P^{N-1}_\mathbb R,P^{N-1}_\mathbb C$, and with no relation needed for $P^{N-1}_+$.
\end{proposition}

\begin{proof}
We must divide the algebra $C(P^{N-1}_+)$ by the relations associated to the diagrams in the statement, as well as those associated to their shifted versions, given by:
$$
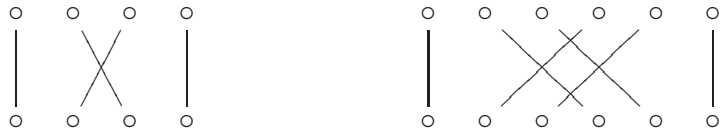
\xymatrix@R=10mm@C=3mm{\circ\ar@{-}[d]&\circ\ar@{-}[dr]&\circ\ar@{-}[dl]&\circ\ar@{-}[d]\\\circ&\circ&\circ&\circ}
\qquad\ \qquad\ \qquad 
\xymatrix@R=10mm@C=3mm{\circ\ar@{-}[d]&\circ\ar@{-}[drr]&\circ\ar@{-}[drr]&\circ\ar@{-}[dll]&\circ\ar@{-}[dll]&\circ\ar@{-}[d]\\\circ&\circ&\circ&\circ&\circ&\circ}$$ 

(1) The basic crossing, and its shifted version, produce the following relations: 
$$p_{ab}=p_{ba}\quad,\quad 
p_{ab}p_{cd}=p_{ac}p_{bd}$$

Now by using these relations several times, we obtain the following formula:
$$p_{ab}p_{cd}
=p_{ac}p_{bd}
=p_{ca}p_{db}
=p_{cd}p_{ab}$$

Thus, the space produced by the basic crossing is classical, $P\subset P^{N-1}_\mathbb C$. By using one more time the relations $p_{ab}=p_{ba}$ we conclude that we have $P=P^{N-1}_\mathbb R$, as claimed.

\medskip

(2) The fattened crossing, and its shifted version, produce the following relations:
$$p_{ab}p_{cd}=p_{cd}p_{ab}$$
$$p_{ab}p_{cd}p_{ef}=p_{ad}p_{eb}p_{cf}$$

The first relations tell us that the projective space must be classical, $P\subset P^{N-1}_\mathbb C$. Now observe that with $p_{ij}=z_i\bar{z}_j$, the second relations read:
$$z_a\bar{z}_bz_c\bar{z}_dz_e\bar{z}_f=z_a\bar{z}_dz_e\bar{z}_bz_c\bar{z}_f$$

Since these relations are automatic, we have $P=P^{N-1}_\mathbb C$, and we are done.
\end{proof}

We can now formulate our classification result, as follows:

\begin{theorem}
The basic projective spaces, namely 
$$P^{N-1}_\mathbb R\subset P^{N-1}_\mathbb C\subset P^{N-1}_+$$
are the only monomial ones.
\end{theorem}

\begin{proof}
We follow the proof from the affine case. Let $\mathcal R_\sigma$ be the collection of relations associated to a permutation $\sigma\in S_k$ with $k\in 2\mathbb N$, as in Definition 16.24. We fix a monomial projective space $P\subset P^{N-1}_+$, and we associate to it subsets $G_k\subset S_k$, as follows:
$$G_k=\begin{cases}
\{\sigma\in S_k|\mathcal R_\sigma\ {\rm hold\ over\ }P\}&(k\ {\rm even})\\
\{\sigma\in S_k|\mathcal R_{|\sigma}\ {\rm hold\ over\ }P\}&(k\ {\rm odd})
\end{cases}$$

As in the affine case, we obtain in this way a filtered group $G=(G_k)$, which is stable under removing outer strings, and under removing neighboring strings.  Thus the computations from the affine case apply, and show that we have only 3 possible situations, corresponding to the 3 projective spaces in Proposition 16.25. See \cite{ba3}.
\end{proof}

Let us discuss now similar results for the projective quantum groups. We have:

\begin{definition}
A projective category of pairings is a collection of subsets 
$$NC_2(2k,2l)\subset E(k,l)\subset P_2(2k,2l)$$
stable under the usual categorical operations, and satisfying $\sigma\in E\implies |\sigma|\in E$.
\end{definition}

As basic examples for this notion, we have the following projective categories of pairings, where $P_2^*$ is the category of matching pairings:
$$NC_2\subset P_2^*\subset P_2$$

This follows indeed from definitions. Now with the above notion in hand, we can formulate the following projective analogue of the notion of easiness:

\begin{definition}
An intermediate compact quantum group $PO_N\subset H\subset PO_N^+$ is called projectively easy when its Tannakian category
$$span(NC_2(2k,2l))\subset Hom(v^{\otimes k},v^{\otimes l})\subset span(P_2(2k,2l))$$
comes via via the following formula, using the standard $\pi\to T_\pi$ construction,
$$Hom(v^{\otimes k},v^{\otimes l})=span(E(k,l))$$
for a certain projective category of pairings $E=(E(k,l))$.
\end{definition}

Observe that, given an easy quantum group $O_N\subset G\subset O_N^+$, its projective version $PO_N\subset PG\subset PO_N^+$ is projectively easy in our sense. In particular the basic projective quantum groups $PO_N\subset PU_N\subset PO_N^+$ are all projectively easy in our sense, coming from the categories $NC_2\subset P_2^*\subset P_2$. We have in fact the following general result:

\begin{theorem}
We have a bijective correspondence between the affine and projective categories of partitions, given by the operation
$$G\to PG$$ 
at the level of the corresponding affine and projective easy quantum groups.
\end{theorem}

\begin{proof}
The construction of correspondence $D\to E$ is clear, simply by setting:
$$E(k,l)=D(2k,2l)$$

Indeed, due to the axioms in Definition 16.27, the conditions in Definition 16.28 are satisfied. Conversely, given $E=(E(k,l))$ as in Definition 16.28, we can set:
$$D(k,l)=\begin{cases}
E(k,l)&(k,l\ {\rm even})\\
\{\sigma:|\sigma\in E(k+1,l+1)\}&(k,l\ {\rm odd})
\end{cases}$$

Our claim is that $D=(D(k,l))$ is a category of partitions. Indeed:

\medskip

(1) The composition action is clear. Indeed, when looking at the numbers of legs involved, in the even case this is clear, and in the odd case, this follows from:
$$|\sigma,|\sigma'\in E
\implies|^\sigma_\tau\in E
\implies{\ }^\sigma_\tau\in D$$

(2) For the tensor product axiom, we have 4 cases to be investigated, depending on the parity of the number of legs of $\sigma,\tau$, as follows:

\medskip

-- The even/even case is clear. 

\medskip

-- The odd/even case follows from the following computation:
$$|\sigma,\tau\in E
\implies|\sigma\tau\in E
\implies\sigma\tau\in D$$

-- Regarding now the even/odd case, this can be solved as follows:
\begin{eqnarray*}
\sigma,|\tau\in E
&\implies&|\sigma|,|\tau\in E\\
&\implies&|\sigma||\tau\in E\\
&\implies&|\sigma\tau\in E\\
&\implies&\sigma\tau\in D
\end{eqnarray*}

-- As for the remaining odd/odd case, here the computation is as follows:
\begin{eqnarray*}
|\sigma,|\tau\in E
&\implies&||\sigma|,|\tau\in E\\
&\implies&||\sigma||\tau\in E\\
&\implies&\sigma\tau\in E\\
&\implies&\sigma\tau\in D
\end{eqnarray*}

(3) Finally, the conjugation axiom is clear from definitions. It is also clear that both compositions $D\to E\to D$ and $E\to D\to E$ are the identities, as claimed. As for the quantum group assertion, this is clear as well from definitions.
\end{proof}

In order to develop now free projective geometry, a first piece of work is that of developing a theory of free Grassmannians, free flag manifolds, and free Stiefel manifolds, based on the affine theory of the spaces of quantum partial isometries, discussed before. To be more precise, the definition of the free Grassmannians is straightforward, as follows, and the definition of the free flag manifolds and free Stiefel manifolds is very similar:
$$C(Gr_{LN}^+)=C^*\left((p_{ij})_{i,j=1,\ldots,N}\Big|p=p^*=p^2,Tr(p)=L\right)$$

Most of the arguments from the affine case carry over in the projective setting, and with solid and useful affine results to rely upon being available from what we said before. For more on all this, we refer to \cite{bgo} and related papers, and to the book \cite{ba3}.

\section*{16d. Matrix models}

As a last topic of discussion, again following \cite{ba3}, \cite{bb2}, \cite{bb3}, we would like to talk about matrix models for our manifolds. In the quantum group case, we have:

\begin{definition}
A matrix model for a closed subgroup $G\subset U_N^+$ is a morphism
$$\pi:C(G)\to M_K(C(T))$$
where $T$ is a compact space, and $K\geq1$ is an integer.
\end{definition}

More generally, we can model in this way the standard coordinates $x_i\in C(X)$ of the various algebraic manifolds $X\subset S^{N-1}_{\mathbb C,+}$. It is then elementary to show that, under the technical assumption $X^c\neq\emptyset$, there exists a universal $K\times K$ model for the algebra $C(X)$, which factorizes as follows, with $X^{(K)}\subset X$ being a certain algebraic submanifold: 
$$\pi_K:C(X)\to C(X^{(K)})\subset M_K(C(T_K))$$

To be more precise, the universal $K\times K$ model space $T_K$ appears by imposing to the complex $K\times K$ matrices the relations defining $X$, and the algebra $C(X^{(K)})$ is then by definition the image of $\pi_K$. In relation with this, we can set as well:
$$X^{(\infty)}=\bigcup_{K\in\mathbb N}X^{(K)}$$

We are led in this way to a filtration of $X$, as follows:
$$X^c= X^{(1)}\subset X^{(2)}\subset X^{(3)}\subset\ldots\ldots\subset X^{(\infty)}\subset X$$

It is possible to say a few non-trivial things about these manifolds $X^{(K)}$, by using algebraic and functional analytic techniques. In the quantum group case, far more things can be of course said. We refer to \cite{ba3}, \cite{bb2}, \cite{bb3} for a discussion here.

\bigskip

Generally speaking, the matrix models are an excellent topic of research, making the link with random matrix theory, and with serious analysis in general. For more on all this, philosophy and some concrete results as well, we refer to \cite{ba3}, \cite{bb2}, \cite{bb3}.

\section*{16e. Exercises} 

Congratulations for having read this book, and no exercises here, for celebrating this accumulated knowledge, a basic beer will do. This being said, in case you really solved all our exercises so far, and are looking for some more, here is one for you:

\begin{exercise}
Further develop free geometry, including the free Laplacian, and free harmonic functions, and then look into free PDE, and free physics. Is that true physics, at very small scales, quarks or below, what you get?
\end{exercise}

In the hope that you will like this exercise, which is something quite subtle. And please but please, do not take a bureaucratic approach to it, by labelling it ``job for physicists''. There is no such thing as a physicist knowing quantum physics, with Einstein himself being the best example, and help from anyone, including matematicians like you, is needed.

\baselineskip=14pt

\printindex

\end{document}